\documentclass{amsbook}
\usepackage[letterpaper,hmarginratio=1:1]{geometry}
\numberwithin{equation}{section}
\usepackage{chngcntr}
\counterwithout{section}{chapter}

\numberwithin{figure}{chapter} 

\usepackage{amsmath, amsthm, amsfonts, amsbsy, thmtools, amssymb,bm,stmaryrd,xfrac,mathrsfs,imakeidx}
\usepackage[sans]{dsfont}
\usepackage{thm-restate}
\usepackage{hyperref}
\usepackage[mathscr]{euscript}
\usepackage{cleveref}
 \usepackage{subcaption}

\usepackage{tikz}
\usetikzlibrary{shapes.geometric, arrows}

\DeclareMathOperator{\Id}{Id}
\DeclareMathOperator{\Real}{Re}
\DeclareMathOperator{\Imaginary}{Im}

\DeclareMathOperator{\Ai}{Ai}

\DeclareMathOperator{\TW}{TW}

\DeclareMathOperator{\TV}{TV}

\DeclareMathOperator{\no}{no}
\DeclareMathOperator{\ea}{ea}
\DeclareMathOperator{\so}{so}
\DeclareMathOperator{\we}{we}
\DeclareMathOperator{\Adm}{Adm}

\DeclareMathOperator{\ext}{ext}
\DeclareMathOperator{\dist}{dist}
\DeclareMathOperator{\supp}{supp}
\DeclareMathOperator{\semci}{sc}
\DeclareMathOperator{\Tra}{Tra}
\DeclareMathOperator{\PV}{PV}
\DeclareMathOperator{\fin}{fin}
\DeclareMathOperator{\Tr}{Tr}
\DeclareMathOperator{\diag}{diag}

\DeclareMathOperator{\eig}{eig}
\DeclareMathOperator{\dL}{L}
\DeclareMathOperator{\emp}{emp}
\DeclareMathOperator{\inv}{inv}

\DeclareMathOperator{\wind}{wind}

\newcommand{\cE}{{\mathcal E}}

\newcommand{ \sfx }{{\mathsf x}}
\newcommand{ \sfy }{{\mathsf y}}
\newcommand{\sfX}{{\mathsf X}}

\newcommand{ \sft }{{\mathsf t}}

\newcommand{\sfz}{{\mathsf z}}

\newcommand{ \sfs }{{\mathsf s}}

\newcommand{\sfT}{{\mathsf T}}

\newcommand{\fa}{{\mathfrak a}}

\newcommand{\fb}{{\mathfrak b}}
\newcommand{\fc}{{\mathfrak c}}

\newcommand{\fd}{{\mathfrak d}}

\newcommand{\ft}{{\mathfrak t}}

\newcommand{\fK}{{\mathfrak K}}

\newcommand{\fM}{{\mathfrak M}}

\newcommand{\fR}{{\mathfrak R}}

\newcommand{\bmv}{{\bm{v}}}

\newcommand{\rL}{{\rm L}}

\newcommand{\ri}{\mathrm{i}}

\newcommand{\bE}{\mathbb{E}}

\newcommand{\bP}{\mathbb{P}}

\newcommand{\bR}{{\mathbb R}}

\newcommand{\al}{\alpha}
\newcommand{\la}{\lambda}

\newtheorem{thm}{Theorem}[section]
\newtheorem{thr}[thm]{Theorem}
\newtheorem{prop}[thm]{Proposition}
\newtheorem{lem}[thm]{Lemma}
\newtheorem{cor}[thm]{Corollary}

\theoremstyle{remark}
\newtheorem{rem}[thm]{Remark}
\theoremstyle{definition}

\theoremstyle{definition}
\newtheorem{definition}[thm]{Definition}
\newtheorem{example}[thm]{Example}
\newtheorem{assumption}[thm]{Assumption}

\title{Strong Characterization for the Airy Line Ensemble}
\author{Amol Aggarwal$^{1, 2}$}
\author{Jiaoyang Huang$^3$}
\address[1]{Department of Mathematics, Columbia University}
\address[2]{Clay Mathematics Institute}
\address[3]{Department of Statistics and Data Science, University of Pennsylvania}

\makeindex 

\begin{document}

\begin{abstract}
	
	In this paper we show that a Brownian Gibbsian line ensemble whose top curve approximates a parabola must be given by the parabolic Airy line ensemble. More specifically, we prove that if $\bm{\mathcal{L}} = (\mathcal{L}_1, \mathcal{L}_2, \ldots )$ is a line ensemble satisfying the Brownian Gibbs property, such that for any $\varepsilon > 0$ there exists a constant $\mathfrak{K} (\varepsilon) > 0$ with
	\begin{flalign*} 
		\mathbb{P} \Big[ \big| \mathcal{L}_1 (t) + 2^{-1/2} t^2 \big| \le \varepsilon t^2 + \mathfrak{K} (\varepsilon) \Big] \ge 1 - \varepsilon, \qquad \text{for all $t \in \mathbb{R}$},
	\end{flalign*} 
	
	\noindent then $\bm{\mathcal{L}}$ is the parabolic Airy line ensemble, up to an independent affine shift. Specializing this result to the case when $\bm{\mathcal{L}} (t) + 2^{-1/2} t^2$ is translation-invariant confirms a prediction of Okounkov and Sheffield from 2006 and Corwin--Hammond from 2014.
	
\end{abstract}

\maketitle

\setcounter{tocdepth}{1}

\tableofcontents

\chapter{Results and Preliminaries}	

\label{BRIDGES0}

\section{Introduction}

\subsection{Preface}

\begin{figure}
	\center
	\includegraphics[scale=.6, trim=0 0.5cm 0 0.5cm, clip]{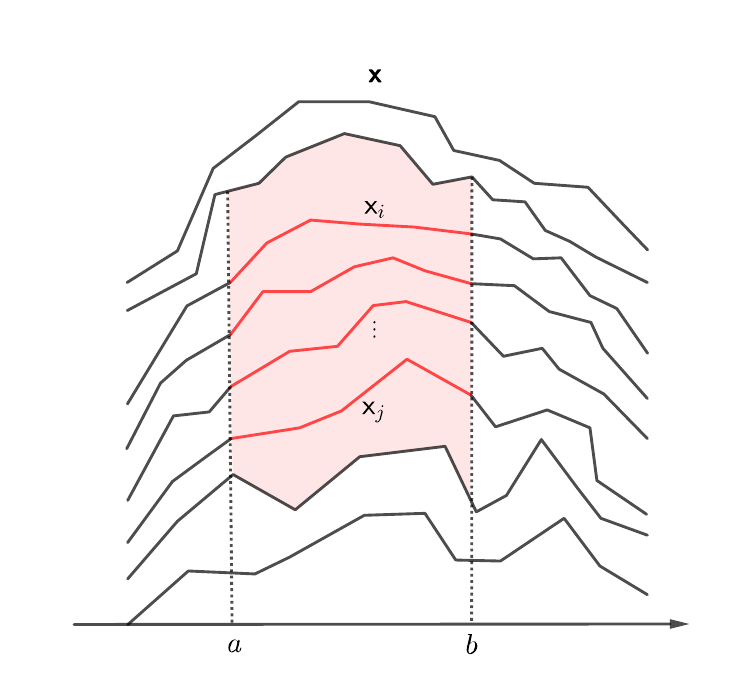}
	
	\caption{Depicted above is an example of Brownian Gibbsian line ensemble, where the red curves can be resampled in the shaded region. }
	\label{f:Line_Ensemble}
\end{figure}

A fundamental question in probability theory and mathematical physics concerns the classification of Gibbs measures for statistical mechanical systems. In this paper we analyze such questions for \emph{Brownian Gibbsian line ensembles}, which are infinite sequences of random functions (or curves) $\bm{\mathsf{x}} = (\mathsf{x}_1, \mathsf{x}_2, \ldots )$, with each $\mathsf{x}_j : \mathbb{R} \rightarrow \mathbb{R}$ continuous, that satisfy the \emph{Brownian Gibbs property}. The latter imposes two constraints. The first is that the $\mathsf{x}_j$ are ordered, meaning that $\mathsf{x}_1 > \mathsf{x}_2 > \cdots $ almost surely. The second is a resampling condition indicating that $\bm{\mathsf{x}}$ behaves as a family of two-sided Brownian motions conditioned to never intersect. More specifically, for any integers $1 \le i \le j$ and real numbers $a < b$, upon conditioning on $\mathsf{x}_k (s)$ with either $k \notin [i, j]$ or $s \notin [a, b]$, the law of the remaining $(\mathsf{x}_i, \mathsf{x}_{i+1}, \ldots , \mathsf{x}_j)$ on $[a, b]$ is given by standard Brownian bridges (whose starting and ending points are determined by the conditioning) conditioned to not intersect, stay below $\mathsf{x}_{i-1}$, and stay above $\mathsf{x}_{j+1}$ (that is, to satisfy $\mathsf{x}_{i-1} > \mathsf{x}_i > \cdots > \mathsf{x}_{j+1}$, where $\mathsf{x}_0 = \infty$). See \Cref{f:Line_Ensemble} for a depiciton. 

A prominent example of a Brownian Gibbsian line ensemble is the \emph{parabolic Airy line ensemble}, introduced by Prah\"{o}fer--Spohn \cite{SIDP} as the scaling limit for the multi-layer polynuclear growth (PNG) model; it can also be viewed as the edge limit for $n$ non-intersecting Brownian bridges, sometimes called the Brownian watermelon. These models are exactly solvable, or integrable, through the framework of determinantal point processes. In \cite{SIDP} (as well as in the prior works of Mac{\^e}do \cite{UPCSE} and Forrester--Nagao--Honner \cite{COUS}, and in the subsequent work of Johansson \cite{DPGDP}), the multi-point correlation functions of the Airy line ensemble were computed in terms of Airy functions. These calculations in particular implied that its curves decay parabolically, but become jointly translation-invariant after simultaneously shifting them by a parabola. The top curve of the Airy line ensemble is known as the Airy$_2$ process, whose one-point marginal is the Tracy--Widom disitribution governing fluctuations for the largest eigenvalue of a Gaussian Unitary Ensemble (GUE) random matrix \cite{LDK}. By combining these integrable inputs with a probabilistic analysis, Corwin--Hammond \cite{PLE} realized the parabolic Airy line ensemble as a family of continuous functions satisfying the Brownian Gibbs property (that is, as a Brownian Gibbsian line ensemble); the Airy line ensemble (incorporating the above parabolic shift) was later shown by Corwin--Sun \cite{ELE} to be ergodic under translations.

Over the past two decades, the Airy line ensemble has become a central object in random surfaces and stochastic growth models. In particular, it has long been understood that many random surfaces exhibit boundary-induced phase transitions, in that they can admit sharp interfaces separating faceted regions (where the surface is almost deterministically flat) from rough ones (where it appears more random). For Ising crystals, this phenomenon dates back to the Wulff construction (see the books of Dobrushin--Koteck\'{y}--Shlosman \cite{C} and Cerf \cite{CPM}) and, for other surfaces, to the work of Jockusch--Propp--Shor \cite{RTT} (who studied random domino tilings of the Aztec diamond).

At this interface, also called the arctic boundary or facet edge, the level lines of the random surface height function are believed to exhibit $n^{1/3}$ fluctuations on domains of diameter $n$; this $1 / 3$ exponent is closely related to the Pokrovsky--Talapov law that predicts the behavior of facet transitions in two-dimensional crystals \cite{GSSPD,TTDIC}. Upon rescaling by $n^{1/3}$, it is further believed that these level lines converge to the Airy line ensemble in the large $n$ limit. While this prediction remains unproven in general, it has been established for various solvable models, starting with the Brownian watermelon in \cite{DPGDP} and random plane partitions by Okounkov--Reshetikhin in \cite{CFPALG,okounkov2007random}. We refer to the survey of Johansson \cite{EFLS} for an exposition and extensive list of further references, as well as to the work of Ferrari--Shlosman \cite{PM} for additional predictions in this direction. 

The relation to stochastic growth models is that their height fluctuations (under wedge initial data) should converge to the Airy$_2$ process, a ubiquitous feature of systems in the Kardar--Parisi--Zhang (KPZ) universality class \cite{DSGI}; see the surveys of Corwin \cite{EUC} and Quastel \cite{IT}. One explanation for this is that, at least in some cases, these models can be exactly mapped to the facet edge of a corresponding random surface, also sometimes called a Gibbsian line ensemble (as the level lines of the surface height function form a line ensemble satisfying a Gibbs property). This idea was initially applied in \cite{RTT}, which used the shuffling algorithm introduced by Elkies--Kuperberg--Larsen--Propp \cite{AMT} to map the discrete-time totally asymmetric simple exclusion process (TASEP) to the arctic boundary for a random domino tiling. Following the framework of Rost \cite{NBPP}, \cite{RTT} showed a hydrodynamical limit for this TASEP, yielding the limit shape for the arctic boundary.

Such correspondences have since been more fruitfully used in reverse, to show that the Airy fluctuations for random surfaces imply those for stochastic growth models. This was first applied to analyze determinantal systems, such as TASEPs \cite{SFRM,TBP} through random tilings, and PNG models \cite{NPRTM,SIDP,DPGDP} (generalizations of longest increasing subsequences of random permutations, studied by Baik--Deift--Johansson \cite{DLLIS}) and Brownian last passage percolation (by Baryshnikov \cite{GQ}, O'Connell--Yor \cite{AT,PTRW}, and Warren \cite{II}) through Brownian watermelons. See  \cite{EODLE,RMDP,RGM} for surveys on these earlier papers. Later work of Hammond \cite{RLE,MCPWP,PRPW,ERDP} used the associated line ensembles to provide a detailed probabilistic analysis of the on-scale polymer geometry for Brownian last passage percolation. More recent papers of Matetski--Quastel--Remenik \cite{TFP} and Dauvergne--Ortmann--Vir\'{a}g \cite{TL} analyzed the full space-time scaling limit for TASEPs and last passage models, under arbitrary initial data. The latter in particular showed how the above correspondences with random surfaces (for them, the Brownian watermelon) could be used to describe this limit entirely in terms of the Airy line ensemble, further solidifying its role in the KPZ universality class. There is now a vast literature utilizing Gibbsian line ensembles to elucidate the probabilistic structure behind KPZ models. For examples just in the last several years, we refer to the papers (and references therein) \cite{ITFP,AFP,ACPIC,DLE} that used line ensembles to prove Brownian comparison results for KPZ models; \cite{ECT,LCDP,UTEL,EL} that used them to analyze the fine structure of the continuum directed polymer; \cite{FGPS,DOL,LGC,FGSDP,GNL} that used them to examine the fractal behavior of the directed landscape; and \cite{SHDLP,ETFP,NTM} that used them to study exceptional times, and related applications, for the KPZ fixed point.

The reasons for the effectiveness of random surface models, in understanding convergence to Airy statistics, can be viewed as twofold. The first reason is algebraic; if the model is integrable, then its solvable underpinnings often become more visible when one examines the random surface as a whole, as opposed to only its arctic boundary. Indeed, the former combinatorially corresponds to a Gelfand--Tsetlin pattern, which contains significantly more structure than the latter, which corresponds to its first (or last) column. This structure enables the introduction of natural $2 + 1$ dimensional dynamics on these random surface models, which project precisely to many of the $1 + 1$ dimensional growth systems in the KPZ universality class. See the works of Borodin--Ferrari \cite{AGRS}, O'Connell \cite{DPQL}, and Borodin--Corwin \cite{P} for examples of this perspective.

The second (which is more relevant to the impetus of this paper) is probabilistic and relates to the Gibbs property satisfied by random surfaces defined by local Boltzmann weights. Although the microscopic Gibbs property behind such a model might depend on the details of its definition, the general intuition is that this Gibbs property should converge to the Brownian one around a facet edge. Indeed, in such regions, the random surface becomes more flat, so its level lines become more sparse and separated. Hence, any local interactions between them should be asymptotically lost, making these level lines behave as random walks that do not intersect. Taking their scaling limit, one then expects to find an infinite family of non-intersecting Brownian bridges. 

Facilitated by the extensive array of methodology to show convergence to Brownian bridges, the above heuristic has been justified for wide classes of random surfaces, both solvable and not. In particular, assuming certain tightness and curvature conditions for their topmost curve, it has been proven that any limit point for the edge of such Gibbsian line ensembles must be Brownian ones; see the works of Dimitrov--Wu \cite{CRWB,TLE}, Barraquand--Corwin--Dimitrov \cite{STEL}, and Serio \cite{TDLE}. Ideas of this nature had earlier been used to prove qualitative results (such as local Brownian continuity for the height function) for integrable, but non-determinantal, models in the KPZ universality class. These include the KPZ equation and O'Connell--Yor polymer \cite{LE}; asymmetric simple exclusion process and stochastic six-vertex model \cite{TFSVM}; and log-gamma polymer \cite{TDLEE,STEL} (where the more involved associated Gibbsian line ensembles arose from works of O'Connell--Warren \cite{MSE} and Nica \cite{IDL}; Borodin--Bufetov--Wheeler \cite{SVMP}; and Corwin--O'Connell--Sepp\"{a}l\"{a}inen--Zygouras \cite{TCF} and Johnston--O'Connell \cite{SLNP}, respectively).

The above frameworks are well-suited to the qualitative task of showing that any edge limit for a random surface must be a Brownian Gibbsian line ensemble (up to tightness and curvature constraints for the extreme level line). However, they do not address the quantitative task of pinning the limit down as the parabolic Airy line ensemble. Therefore, a basic question that arises is if there is an axiomatic characterization of the Airy line ensemble or, more specifically, some practical criterion for when a Brownian Gibbsian line ensemble must be the parabolic Airy one. The purpose of this paper is to establish such a criterion. \\

This criterion can be stated as follows (see \Cref{l0} below). Let $\bm{\mathcal{L}} = (\mathcal{L}_1, \mathcal{L}_2, \ldots )$ denote a Brownian Gibbsian line ensemble. Suppose for any $\varepsilon > 0$ that there is a constant $\mathfrak{K} (\varepsilon) > 0$ with
\begin{flalign}
	\label{l1t} 
	\mathbb{P} \Big[ \big| \mathcal{L}_1 (t) + 2^{-1/2} t^2 \big| \le \varepsilon t^2 + \mathfrak{K}(\varepsilon) \Big] \ge 1 - \varepsilon, \qquad \text{for each $t \in \mathbb{R}$}.
\end{flalign}

\noindent Then $\bm{\mathcal{L}}$ is the parabolic Airy line ensemble, up to an independent affine shift; see \Cref{lr} below. 

Informally, \eqref{l1t} states that the top curve $\mathcal{L}_1$ of $\bm{\mathcal{L}}$ likely satisfies $\mathcal{L}_1 (t) = -\big( 2^{-1/2} + o(1) \big) t^2$ (with the constant $\mathfrak{K}(\varepsilon)$ in \eqref{l1t} being used to correct this approximation for small $t$). Let us mention that some type of quadratic decay for $\mathcal{L}_1$ must be imposed for the above characterization to hold. For instance, the Airy line ensembles with wanderers introduced by Adler--Ferrari--van Moerbeke \cite{PWUC} form examples of Brownian Gibbsian line ensembles for which $\mathcal{L}_1$ only decays linearly. 

Observe that \eqref{l1t} incorporates the scenario when the parabolically shifted line ensemble $\bm{\mathcal{L}} (t) + 2^{-1/2} t^2$ is translation-invariant in $t$. In this case, it was predicted\footnote{It was unpublished at the time but has since appeared in various forms in print, such as \cite[Conjecture 3.2]{PLE} and \cite[Conjecture 1.7]{ELE}.}  by Okounkov and Sheffield in 2006 that $\bm{\mathcal{L}}$ is given by a parabolic Airy line ensemble, up to an independent overall constant shift. This is in the spirit of classifications for translation-invariant Gibbs measures of discrete random surfaces, proven by Sheffield \cite{RS} (but is also of a distinct nature, since here the base space of the  line ensemble is not discrete, and also since here translation-invariance holds in only one coordinate, not both).

Our result \Cref{lr} quickly implies this prediction (see \Cref{linvariant} below), and further generalizes upon it in two ways. First, our assumption \eqref{l1t} only constrains the top curve of the ensemble, instead of imposing that all of its curves be jointly translation-invariant. The notion that sufficient information on the top curve could determine the entire line ensemble also appeared in the work of Dimitrov--Matetski \cite{CLE,CHLE}, though the control they required was quite siginificant, namely, knowledge of its full law (all of its finite-dimensional marginals). Those results in particular implied that $\bm{\mathcal{L}}$ is a parabolic Airy line ensemble, if one happened to know in advance that $\mathcal{L}_1$ were an Airy$_2$ process. Prior to our work, the latter seemed to be quite an involved task, though had been done by Quastel--Sarkar \cite{CEPEFP} and Vir\'{a}g \cite{HTL} for some special Gibbsian line ensembles, such as the KPZ one (as we explain below, our results directly imply an alternative proof of this KPZ result; see \Cref{converge0}).

Second, \eqref{l1t} only requires the limiting trajectory of $\mathcal{L}_1 (t)$ to approximate a parabola, as opposed to stipulating it to be exactly translation-invariant upon a parabolic shift. In agreement with the terminology from \cite[Section 10.4]{RS}, one might therefore refer to \Cref{lr} as a \emph{strong characterization} for the Airy line ensemble. Strong characterizations for Gibbs measures of random surface models, with power law correlation decay, appear to be quite rare in the literature (outside of the fairly distant setting of random lozenge tilings \cite{ULTS}). \\

Before continuing, let us briefly comment on two potential applications that our characterization may lead to in the future (for which both of the above-mentioned improvements would seem to be quite useful). The first concerns stochastic growth models; many such systems proven to be in the KPZ universality class are not fully solvable in the sense of being determinantal, but instead satisfy a Yang--Baxter equation. These include the stochastic six-vertex model \cite{SVMR,SSVM} and its degenerations (which encompass the KPZ equation and ASEP); certain random polymers \cite{AT,SDPBC,SLP,RWDR}; $q$-deformations of the TASEP \cite{ERTD,P} and PNG model \cite{DPG}; and various other systems. For all of these models, it is known that the one-point marginals of their height functions under wedge initial data converge to the Tracy--Widom GUE distribution; however, their full convergence to the Airy$_2$ process is still open for most of them, except for the ASEP, KPZ equation, and O'Connell--Yor polymer \cite{CEPEFP,HTL}. Using the Yang--Baxter equation alone, it is possible to map the height functions for all of the above models to the arctic boundary of an associated Gibbsian line ensemble; this was first done for the stochastic six-vertex model in \cite{SVMP,TFSVM}, and later systematized to other models through the bijectivization framework\footnote{The use of this framework (and its special case called stochasticization \cite{SSE}) to produce more general families of ``colored line ensembles'' from stochastic models with a Yang--Baxter equation are elaborated and extended upon by Aggarwal--Borodin \cite{UU}. The later work \cite{SLCA} of Aggarwal--Corwin--Hegde probabilistically analyzed these colored line ensembles (partly using the characterization developed in this paper) to derive the scaling limit of the colored ASEP and stochastic six-vertex models.} of Bufetov, Mucciconi, and Petrov \cite{FSF,RFSVM}. One-point convergence results for these models verify the tightness and curvature assumptions for the top curves of these ensembles, which might enable one to extend the frameworks developed in \cite{CRWB,STEL,TLE,TDLE} to show that they converge to Brownian Gibbsian line ensembles. Our characterization \Cref{lr} would then apply, proving their convergence to the Airy line ensemble, and hence of their top curves (tracking the height function of the associated stochastic growth model) to the Airy$_2$ process. In \Cref{EquationPolymer} below, we provide the very quick implementation of this idea for two examples (where the qualitative framework has already been set up), namely, the KPZ equation (\Cref{converge0}) and log-gamma polymer (\Cref{converge1}).

The second potential use of our characterization result is towards proving convergence of edge statistics for random surfaces to the Airy line ensemble. At the moment, there seem to be few (if any) natural examples of non-determinantal random surface models for which this statement has been proven.\footnote{Even one-point convergence statements seem to be rare in this context, but see the recent work of Ayyer--Chhita--Johansson \cite{FMTP} for such a result at the edge of the domain-wall ice model.} As mentioned previously, there exists a fairly robust framework \cite{CRWB,STEL,TLE,TDLE} for proving convergence of edge limits of random surfaces to Brownian line ensembles, assuming certain tightness and curvature constraints for the extreme level line. An obstruction that remains is thus in verifying these constraints; they can be reformulated as a weak local law\footnote{This terminology was adopted from works of Erd\H{o}s--Schlein--Yau \cite{SLSSD,LSLCD}, showing local semicircle laws for random matrices.} at the facet edge for general random surfaces, meaning that their limit shape phenomena hold not only on global scales, but also on mesoscopic ones (of dimensions $n^{1/3} \times n^{2/3}$) near the edge. In the bulk of the liquid region, such local laws have been proven for random tilings \cite{ULTS} by an inductive application of the associated variational principle on progressively smaller scales. It is an enticing question to see if those ideas can be extended to the facet edge of general random surface models, which together with our characterization might lead to universality results for the Airy line ensemble. 

We now return to the characterization \Cref{lr} and proceed to describe some of the ideas behind its proof (see \Cref{ProofL} below for a more precise exposition).

\subsection{Proof Overview}

\label{LProof1}

In what follows, we denote the parabolic Airy line ensemble by $\bm{\mathcal{R}} = (\mathcal{R}_1, \mathcal{R}_2, \ldots )$.  To be consistent with previous works, it will be its rescaling $2^{-1/2} \cdot \bm{\mathcal{R}}$ that satisfies the Brownian Gibbs property (of variance $1$); we denote $\bm{\mathcal{S}} = 2^{-1/2} \cdot \bm{\mathcal{R}}$, so that $\mathcal{S}_j (t) = 2^{-1/2} \cdot \mathcal{R}_j (t)$ for each $(j, t) \in \mathbb{Z}_{\ge 1} \times \mathbb{R}$. We will show that a line ensemble $\bm{\mathcal{L}}$ satisfying \eqref{l1t} is equal to $\bm{\mathcal{S}}$, up to an affine shift. To this end, we will prove a sequence of results indicating that $\bm{\mathcal{L}}$ is close to $\bm{\mathcal{S}}$, in an increasingly fine sense. To explain this further, we first recall from the work of Soshnikov \cite{FNP} (see also that of Dauvergne--Vir\'{a}g \cite{BPLE}) that with high probability the paths in $\bm{\mathcal{S}}$ satisfy
\begin{flalign}
	\label{2rj2} 
	 \mathcal{S}_j (t) + 2^{-1/2} t^2  =  - 2^{-7/6} (3 \pi)^{2/3} j^{2/3} + \mathcal{O} (j^{o(1)-1/3}).
\end{flalign}

\begin{figure}
	\center
	\includegraphics[scale = .6]{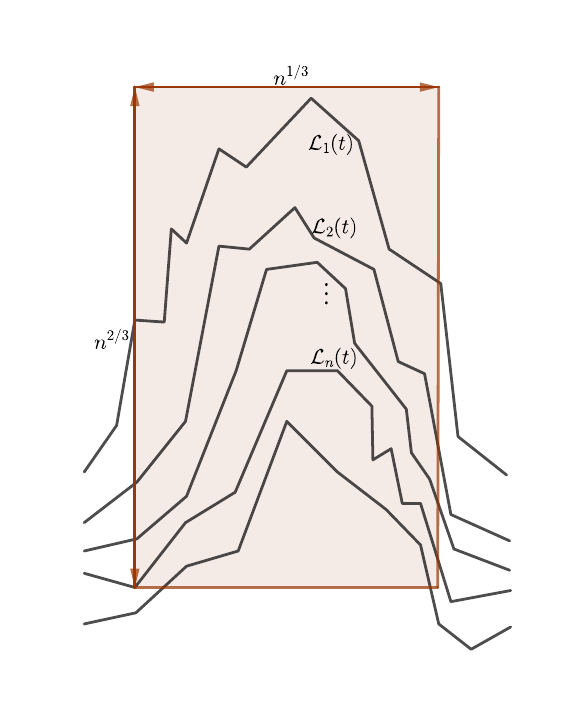}
	
	\caption{Shown above is an order $n^{1/3} \times n^{2/3}$ time-space $(t, x)$ scaling window to examine the curves $\mathcal{L}_k$ with $k$ of order $n$. 	
	\label{f:scaling}}
\end{figure}

\noindent In particular, given an integer $n \ge 1$, we likely have $\mathcal{S}_n (0)$ is of order $-  n^{2/3}$. More generally, if $t$ is of order $n^{1/3}$ and $j$ is of order $n$, then \eqref{2rj2} implies that $\mathcal{S}_j (t)$ is likely of order $- n^{2/3}$. For this reason, we will compare the top $n$ curves of $\bm{\mathcal{L}}$ to those of $\bm{\mathcal{S}}$ when the time $t$ and space $x$ parameters are of order $n^{1/3}$ and $n^{2/3}$, respectively; see \Cref{f:scaling} for a depiction.

On a more refined level, \eqref{2rj2} indicates that, while the ``deep'' paths (of high index) in $\bm{\mathcal{S}}$ are of large magnitude, they tightly concentrate around smooth, deterministic functions, in both the time (horizontal) and space (vertical) directions. In the time direction, \eqref{2rj2} implies that with high probability $\mathcal{S}_n (t)$ closely approximates a parabola of curvature $-2^{-1/2}$, namely,  
\begin{flalign}
	\label{rntft} 
	\big| \mathcal{S}_n (t) - \mathfrak{f}_n (t) \big| = o (1), \qquad \text{for $\mathfrak{f}_n (t) = -2^{-1/2} t^2 - 2^{-7/6} (3\pi)^{2/3} n^{2/3}$}.
\end{flalign}

\noindent In the space direction, \eqref{2rj2} implies for any $t$ that $\big( \mathcal{S}_{n+1} (t), \mathcal{S}_{n+2} (t), \ldots , \mathcal{S}_{2n} (t) \big)$, obtained from the paths of $\bm{\mathcal{S}}$ at time $t$ with indices in $\{ n+1, n+2, \ldots , 2n \}$, likely approximates a smooth, deterministic profile. More specifically, for any index $k \in [1, n]$, we likely have 
\begin{flalign}
	\label{rkngy}
	\big| \mathcal{S}_{k+n} (t) - n^{2/3} \cdot \mathfrak{g}_{tn^{-1/3}} (kn^{-1}) \big| = o(1), \quad \text{for $\mathfrak{g}_r (y) = - 2^{-7/6} (3\pi)^{2/3} (y+1)^{2/3} - 2^{-1/2} r^2$},
\end{flalign}

\noindent where we observe that the profile $\mathfrak{g}_r (y)$ is smooth in $y \in [0, 1]$. To establish \Cref{lr}, we will first prove that weaker variants of the bounds \eqref{2rj2}, \eqref{rntft}, and \eqref{rkngy} hold for $\bm{\mathcal L}$. In the first, we allow for a larger error; in the last two, we replace the deterministic functions $\mathfrak{f}_n$ and $\mathfrak{g}_r$ with unspecified, random functions $h_n$ and $\gamma_r$ (that likely satisfy some similar properties to $\mathfrak{f}_n$ and $\mathfrak{g}_r$, respectively). 

In particular, we will proceed by proving the following four, increasingly precise, statements. In what follows, $A > 1$ is an arbitrary constant; $n$ is an integer parameter that we view as tending to infinity; $t \in [-An^{1/3}, An^{1/3}]$ is a time parameter; $i$ and $j$ are indices with $1 \le i \le j \le An$; and all claims below hold with high probability.

\begin{enumerate}
	\item \emph{On-scale estimates}: The scaling in \eqref{2rj2} is valid for $\bm{\mathcal{L}}$, in two senses (\Cref{sclprobability}).	
	\begin{enumerate} 
		\item \emph{Path locations}: $-450 j^{2/3} \le \mathcal{L}_j (t) + 2^{-1/2} t^2 \le -j^{2/3} / 200$ likely holds, if $j \ge n / A$. 
		\item \emph{Gap upper bound}: $\big| \mathcal{L}_i (t) - \mathcal{L}_j (t) \big| \le \mathcal{O} (j^{2 /3} - i^{2/3}) + (\log n)^{25} i^{-1/3}$ likely holds.
	\end{enumerate} 
	\item \emph{Global law and regularity}: For $\bm{\mathcal{L}}$, \eqref{2rj2} likely holds but with a larger error $o(n^{2/3})$ (\Cref{p:globallaw2}), and \eqref{rkngy} likely holds but with an unknown, regular function $\gamma_r$ replacing $\mathfrak{g}_r$ (\Cref{p:closerho0}).   
	\begin{enumerate} 
		\item \emph{Global law}: $\big| \mathcal{L}_j (t) + 2^{1/2} t^2 + 2^{-7/6} (3\pi)^{2/3} j^{2/3} \big| = o(n^{2/3})$ likely holds.
		\item \emph{Spatial regularity}: There likely is an almost smooth (whose first $50$ derivatives are uniformly bounded in $n$), random function $\gamma_{tn^{-1/3}}: [0, 1] \rightarrow \mathbb{R}$, such that $\big| \mathcal{L}_{k+n} (t) - n^{2/3} \cdot \gamma_{tn^{-1/3}} ( k / n) \big| = o(1)$ whenever $1 \le k \le n$. 
	\end{enumerate} 
	
	\item \emph{Curvature approximation}: There likely is a random function $h_n : [-An^{1/3}, An^{1/3}] \rightarrow \mathbb{R}$, so that $\big| h_n'' (s) + 2^{-1/2} \big| = o(1)$ and $\big| \mathcal{L}_n (s) - h_n (s) \big| = o(1)$  for all $|s| \le An^{1/3}$ (\Cref{h0x2}).
	\item \emph{Airy statistics}: The ensemble $\bm{\mathcal{L}}$ has Airy statistics (\Cref{xdifferenceconverge} and \Cref{qsf}).
	\begin{enumerate} 
		\item \emph{Airy gaps}: The joint law of the gaps $\big( \mathcal{L}_1 (s) - \mathcal{L}_2 (s), \mathcal{L}_2 (s) - \mathcal{L}_3 (s), \ldots \big)$ coincides with that of the Airy point process $\big( \mathcal{S}_1 (0) - \mathcal{S}_2 (0), \mathcal{S}_2 (0) - \mathcal{S}_3 (0), \ldots \big)$. 
		\item \emph{Airy line ensemble}: Up to an affine shift, $\bm{\mathcal{L}}$ is a scaled parabolic Airy line ensemble $\bm{\mathcal{S}}$. 
	\end{enumerate}
\end{enumerate}

After providing a more detailed proof discussion and reviewing miscellaneous preliminary results in \Cref{BRIDGES0}, we will show the on-scale estimates in \Cref{GAPSCALE}. After proving several results about limit shapes for non-intersecting Brownian bridges in \Cref{EDGESHAPE} and couplings for them in \Cref{RectangleLCouple}, we will establish the global law and regularity for $\bm{\mathcal{L}}$ in \Cref{GlobalRegular}. We will prove the curvature approximation in the second half of \Cref{APPROXIMATECURVE}, and that $\bm{\mathcal{L}}$ has Airy statistics in \Cref{STATISTICSBRIDGES}. \\

Before proceeding, let us mention that it would be interesting to classify all Brownian Gibbsian line ensembles with a (finite) top curve. We suspect that there are no such ensembles whose top curve decays super-quadratically. Proving such a statement along the lines of our methodology would necessitate several changes, though, starting at least with a reworking of the on-scale estimates (whose scaling exponents would be different in this setting). As mentioned previously, there exist Brownian line ensembles whose top curves decay sub-quadratically, given by the Airy line ensembles with a finite number $k$ of wanderers. The top $k$ curves of these ensembles decay linearly, while the $(k+1)$-th one decays quadratically. It might be possible to extend our framework to show that Brownian Gibbsian line ensembles, whose $(k+1)$-th curve approximates a parabola, must be an Airy line ensemble with at most $k$ wanderers. In particular, it is plausible that the on-scale estimates; global law and regularity; and curvature approximation could be directly adaptable to that context, while the conclusion of Airy statistics would require some modification.

As shown by Dimitrov \cite{WLE}, there additionally exist Airy line ensembles with infinitely many wanderers. These are obtainable by taking a scaling limit of a last passage percolation model with two sets of parameters, as studied by Borodin--P\'{e}ch\'{e} \cite{KTSP}. One might (speculatively) hypothesize that any Brownian Gibbsian line ensemble must also be obtainable as suitable limits of the latter model. However, the full set of such line ensembles that would arise from such a limiting procedure, and their quantitative (or qualitative) properties, are not yet transparent to us. \\

Returning to the proof outline, we next describe the above four statements, and some ideas underlying their proofs, in greater detail. As we will see, an obstacle we will repeatedly face is the lack of control on the curves $\mathcal{L}_j$ of $\bm{\mathcal{L}}$. Even up until midway through the last (fourth) statement of the above overview, our estimates on the $\mathcal{L}_j$ will be quite poor, unable to forbid them from fluctuating more than the parabolic Airy line ensemble itself. On the other hand, over the past twenty-five years, Dyson Brownian motion (and equivalent familes of non-intersecting Brownian bridges, including Brownian watermelons as a special case) has been comprehensively understood, both from the perspective of exact solvability (starting with the works of Br\'{e}zin--Hikami \cite{USCG} and Johansson \cite{ULSDM}) and stochastic analysis (see the reviews of Guionnet \cite{LDSCRM} and Erd\H{o}s--Yau \cite{DRM}). 

A substantial portion of our analysis is therefore centered on devising a series of comparisons between the line ensemble $\bm{\mathcal{L}}$ and Dyson Brownian motion; this will enable us to transfer results about the latter (that are sometimes already available in the literature, which we will explain as they arise) to the former. These two systems initially appear to be quite different and, indeed, the first forms of our comparisons will be fairly coarse (though sufficient to prove the on-scale estimates, for example). However, as we continue to learn more about the line ensemble $\bm{\mathcal{L}}$, we will use the bounds obtained from previous comparisons to concoct new and improved ones, eventually reaching the level where we can compare exact Airy statistics. 

Now let us outline our proof in more detail. The main purpose of the below outline is to serve as a guide for readers examining in greater depth the arguments presented in the body of this paper; on occasion, they may wish to consult this outline to recall the overarching ideas and intuition underlying these arguments, to help them navigate through its lengthier details. As such, this outline may be skimmed or skipped on an initial reading, especially since it may get a bit involved at some points.

\subsection{On-Scale Estimates}

Before discussing our proof of the on-scale estimates, we first describe a coupling, called gap monotonicity, that will play an extensive role in many of our arguments.	

\subsubsection{Gap Monotonicity} 

\label{Gap0} 

\begin{figure}
	\center
\includegraphics[scale=.8]{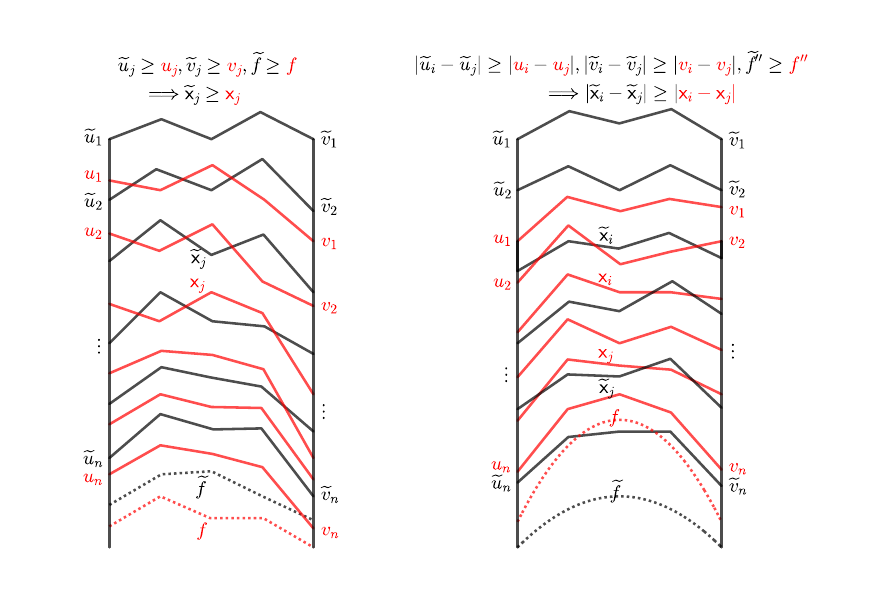}

\caption{Shown to the left is a depiction for height monotonicity.  Shown to the right is a depiction for gap monotonicity }
\label{f:comparison}
	\end{figure}

Monotone couplings have long been fundamental in the analysis of random surfaces. In the context of Brownian Gibbsian line ensembles, the most commonly used such coupling is called \emph{height monotonicity}, which indicates that non-intersecting Brownian bridges are increasing in their boundary data. More precisely, sample two families of $n$ non-intersecting Brownian bridges $\bm{\mathsf{x}} = (\mathsf{x}_1, \mathsf{x}_2, \ldots , \mathsf{x}_n)$ and $\widetilde{\bm{\mathsf{x}}} = (\widetilde{\mathsf{x}}_1, \widetilde{\mathsf{x}}_2, \ldots , \widetilde{\mathsf{x}}_n)$ starting at $n$-tuples $\bm{u}$ and $\widetilde{\bm{u}}$, respectively; ending at $n$-tuples $\bm{v}$ and $\bm{\widetilde{v}}$, respectively; and conditioned to stay above lower boundary curves $f$ and $\widetilde{f}$, respectively.\footnote{One can also constrain them to lie below upper boundaries, but we will not do so in this introductory exposition.} Assume that $\bm{u} \le \widetilde{\bm{u}}$, that $\bm{v} \le \widetilde{\bm{v}}$, and that $f \le \widetilde{f}$. Then, it is possible to couple $\bm{\mathsf{x}}$ and $\widetilde{\bm{\mathsf{x}}}$ such that each $\mathsf{x}_j \le \widetilde{\mathsf{x}}_j$. See the left side of \Cref{f:comparison} for a depiction.

In this paper we further require a different type of monotonicity that compares not the Brownian bridges themselves but the gaps between them. We refer to this as \emph{gap monotonicity}, stated as \Cref{monotonedifference} below, which indicates that the gaps between non-intersecting Brownian bridges are increasing in the gaps of their starting and ending data, and also in the convexity of their lower boundary curves. More precisely, assume that each $|u_i - u_j| \le |\widetilde{u}_i - \widetilde{u}_j|$ and $|v_i - v_j| \le |\widetilde{v}_i - \widetilde{v}_j|$, and also that  $f'' \le \widetilde{f}''$ (in the sense of distributions). Then, it is possible to couple $\bm{\mathsf{x}}$ and $\widetilde{\bm{\mathsf{x}}}$ such that each $|\mathsf{x}_i - \mathsf{x}_j| \le |\widetilde{\mathsf{x}}_i - \widetilde{\mathsf{x}}_j|$. See the right side of \Cref{f:comparison} for a depiction.

Perhaps the simplest proof of height monotonicity (see \cite{PLE}, for example) proceeds by first discretizing the Brownian bridges into non-intersecting Bernoulli random walks; coupling the latter under a local Markov chain (such as the Glauber dynamics) that preserves height orderings; running this chain until it mixes; and taking any limit point of the dynamics as a height monotone coupling. Such a proof cannot apply for gap monotonicity, as it is false in this discrete setup (see \Cref{discrete1}).

To prove gap monotonicity, we instead proceed by first ``semi-discretizing'' the Brownian bridges into Gaussian random walks that are continuous in space but discrete in time. They constitute $T \in \mathbb{Z}_{\ge 1}$ steps, which allows us to induct on $T$. To this end, we introduce a non-local Markov chain, which alternates between resampling the first step of all walks simultaneously and their remaining $T-1$ steps. Using the inductive hypothesis (replacing $T$ by either $2$ or $T-1$), we show that we can couple these dynamics so as to preserve gap orderings. By again running this chain until it mixes, this reduces proving semi-discrete gap monotonicity to verifying its base case $T = 2$, which is done directly, by induction on the number of paths. See \Cref{Coupling1} below for further details.

\subsubsection{Path Location Bounds}

\label{Location0}

The first aspect of the on-scale estimates, described in \Cref{LProof1}, states that for $t \in [-An^{1/3}, An^{1/3}]$ and $n / A \le j \le An$ we likely have 
\begin{flalign}
	\label{jljj} 
	-450 j^{2/3} \le \mathcal{L}_j (t) + 2^{-1/2} t^2 \le - \displaystyle\frac{j^{2/3}}{100}.
\end{flalign}  

\noindent This estimates the deep curves of $\bm{\mathcal{L}}$, only assuming the bound \eqref{l1t} on its top curve. While many of the previously mentioned works on Gibbsian line ensembles do show some control on these deep curves $\mathcal{L}_j$, their bounds are usually not optimal (a large power, if not exponential) in their dependence on the index $j$. For certain specific ensembles relating to last passage percolation models, the true dependence on $j$ (up to constants, as in \eqref{jljj}) was shown by Basu--Ganguly--Hammond--Hegde \cite{ISEGW}, by relating such estimates to the geometry of non-intersecting geodesics. 

\begin{figure}
	\center
\includegraphics[scale=.6]{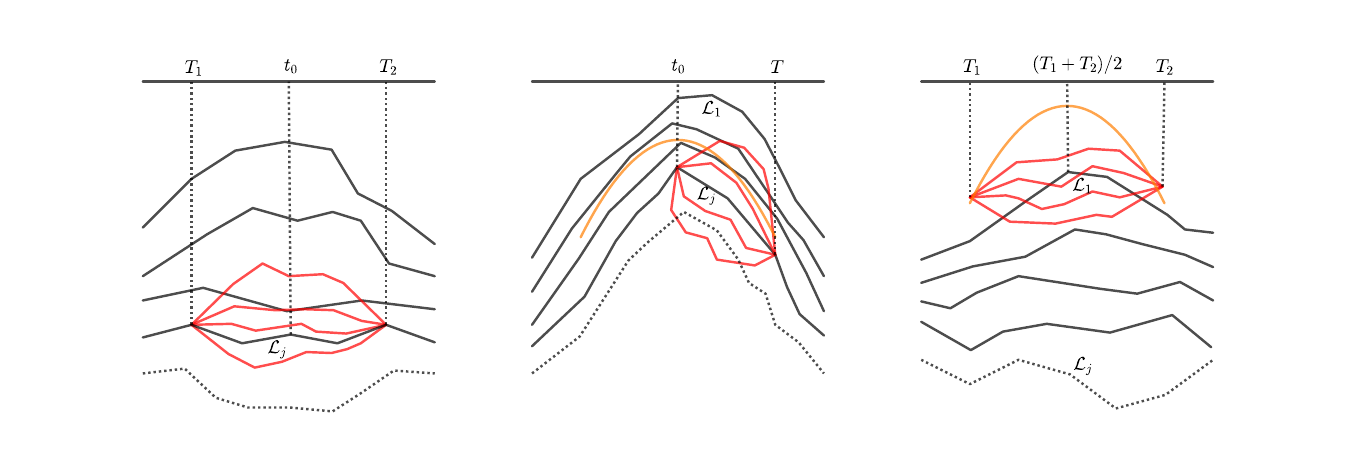}

\caption{Shown above are the three scenarios discussed in \Cref{Location0}, where the black curves are of $\bm{\mathcal{L}}$;  the red ones are the watermelons we eventually compare them to; and the orange one is the parabola that $\mathcal{L}_1$ should be close to, by \eqref{l1t}. On the left, $\mathcal{L}_j$ cannot be too low at time $t_0$ (even after pushing some curves in $\bm{\mathcal{L}}$ down to form the red watermelon). The curve $\mathcal{L}_1$ fails to approximate the orange parabola on the middle (where it is too high, even after pushing some curves in $\bm{\mathcal{L}}$ down to form the red watermelon) and on the right (where it is too low, even after pushing some curves in $\bm{\mathcal{L}}$ up to form the red watermelon).}
	
\label{f:path_location}
	\end{figure}

For general Brownian Gibbsian line ensembles satisfying \eqref{l1t}, this connection between $\bm{\mathcal{L}}$ and last passage percolation is lost, and so our proof instead uses only the Brownian Gibbs property. In particular, we show on $[-An^{1/3}, An^{1/3}]$ that $\mathcal{L}_j$ can neither be very high ($\mathcal{L}_j (t) + 2^{-1/2} t^2 > -j^{2/3} / 100$), nor very low ($\mathcal{L}_j (t) + 2^{-1/2} t^2 < -450 j^{2/3}$). This will be a quick consequence of combining the following three statements, where in all of them we assume that $\mathcal{L}_1 (t)$ is close to the parabola $-2^{-1/2} t^2$, as is likely implied by \eqref{l1t}. While the statements of, and reasoning behind, these claims in this exposition will be imprecise, their proper justification will be obtained by applying height monotonicity to compare $\bm{\mathcal{L}}$ to Brownian watermelons; see \Cref{f:path_location} for depictions and \Cref{LLocation} below for further details.

1. If $\mathcal{L}_j$ is very low at a point $t_0$, then it is likely low on a long interval (\Cref{probabilityeventll0}). Otherwise, there would exist two points $T_1 < t_0 < T_2$ not very far from $t_0$, such that $\mathcal{L}_j (T_1)$ and $\mathcal{L}_j (T_2)$ are much higher than $\mathcal{L}_j (t_0)$. Then resampling the top $j$ curves of $\bm{\mathcal{L}}$ on $[T_1, T_2]$, one finds that the conditional boundary data of these $j$ paths is too high to likely allow their bottom curve to drop to $\mathcal{L}_j (t_0)$ at time $t_0$, which is a contradiction. See the left side of \Cref{f:path_location}.

2. If $\mathcal{L}_j$ is very high at a point $t_0$, then it is likely low on a long interval to the right of $t_0$ (\Cref{probabilityeventclh0}). Otherwise, there would exist a point $T > t_0$ not very far from $t_0$, such that $\mathcal{L}_j (T)$ is not low. Resampling the top $j$ curves of $\bm{\mathcal{L}}$ on $[t_0, T]$, one finds that their conditional starting data at time $t_0$ is high enough (and their conditional ending data at time $T$ is not low enough to counteract them) to cause their top curve $\mathcal{L}_1$ to likely ``shoot'' far above the parabola $-2^{-1/2} t^2$ at some time $R \in [t_0, T]$, which contradicts \eqref{l1t}. See the middle of \Cref{f:path_location}.

3. The curve $\mathcal{L}_j$ is likely not low on any long interval $[T_1, T_2]$ (\Cref{probabilityeventl0}). Otherwise, resampling the top $j-1$ paths in $\bm{\mathcal{L}}$ on $[T_1, T_2]$ (and possibly moving their boundary data up a bit), one finds that the conditional lower boundary $\mathcal{L}_j$ is too low to affect them; as such, it can be removed. Since the interval $[T_1, T_2]$ is long, in the absence of a lower boundary, the top curve $\mathcal{L}_1$ of these $j-1$ paths will likely stay close to the line connecting $\mathcal{L}_1 (T_1) \approx -2^{-1/2} T_1^2$ to $\mathcal{L}_1 (T_2) \approx -2^{-1/2} T_2^2$. Thus, it cannot reach high enough to meet the parabola $-2^{-1/2} t^2$ at, say, the midpoint $(T_1+T_2) / 2$ of $[T_1, T_2]$, which again contradicts \eqref{l1t}. See the right side of \Cref{f:path_location}.

\subsubsection{Gap Upper Bound} 

\label{UpperGap} 

The second aspect of the on-scale estimates, described in \Cref{LProof1}, states for any $1 \le i \le j \le n$ and $t \in [-An^{1/3}, An^{1/3}]$ that we likely have 
\begin{flalign} 
	\label{ltiljji0} 
	\big| \mathcal{L}_i (t) - \mathcal{L}_j (t) \big| = \mathcal{O} (j^{2/3} - i^{2/3}) + (\log n)^{25} i^{-1/3}.
\end{flalign} 

\noindent To show this, we imagine $A \ge 1$ as large (but uniformly bounded); set $\mathsf{T} = An^{1/3}$; and resample the top $2n$ curves of $\bm{\mathcal{L}}$ on $[-2\mathsf{T}, 2\mathsf{T}]$, which become non-intersecting Brownian bridges starting at $\bm{u} = \big( \mathcal{L}_1 (-2\mathsf{T}), \ldots , \mathcal{L}_{2n} (-2\mathsf{T}) \big)$; ending at $\bm{v} = \big( \mathcal{L}_1 (2\mathsf{T}), \ldots , \mathcal{L}_{2n} (2\mathsf{T}) \big)$; and conditioned to stay above $\mathcal{L}_{2n+1}$. By gap monotonicity (recall \Cref{Gap0}), removing the lower boundary $\mathcal{L}_{2n+1}$ increases the gaps between the $\mathcal{L}_n$. So, it suffices to prove the gap upper bound \eqref{ltiljji0}, but for $2n$ non-intersecting Brownian bridges $\bm{\mathsf{x}} = (\mathsf{x}_1, \mathsf{x}_2, \ldots , \mathsf{x}_{2n})$ starting at $\bm{u}$ and ending at $\bm{v}$. 

At this point, we use a known fact relating the law of the non-intersecting Brownian bridges $\bm{\mathsf{x}}$ without upper and lower boundary, to that of certain random matrix spectra; it states the following. For any matrix $\bm{M}$, let $\eig (\bm{M})$ denote its ordered sequence of eigenvalues; also set $\bm{U}$ and $\bm{V}$ to be the $2n \times 2n$ diagonal matrices, with $\eig (\bm{U}) = \bm{u}$ and $\eig (\bm{V}) = \bm{v}$. Fix $t \in [-2 \mathsf{T}, 2 \mathsf{T}]$; setting $S = \mathsf{T} - t^2 / 4 \mathsf{T}$, the law of $\bm{\mathsf{x}}(t)$ is given by $\eig (\bm{A}_t + S^{1/2} \cdot \bm{G})$, where $\bm{G}$ is a $2n\times 2n$ GUE random matrix and $\bm{A}_t = ( 1 / 2 - t / 2\mathsf{T}) \cdot \bm{U} + ( 1 / 2 + t / 2 \mathsf{T}) \cdot \bm{W} \bm{V} \bm{W}^*$, for $\bm{W}$ a certain (not Haar distributed) unitary random matrix (see \Cref{t:lawt} below for the precise statement). Hence, the
\begin{flalign}
	\label{xta0} 
	\text{law of $\bm{\mathsf{x}}(t)$ is given by Dyson Brownian motion run for time $S$, with initial data $\eig (\bm{A}_t)$}.  
\end{flalign}

While Dyson Brownian motion is now well understood, effectively using \eqref{xta0} is typically complicated by the involved law of $\bm{A}_t$. However, in our setting, we will only require a bound on the norm of $\bm{A}_t$ (after subtracting a multiple of the identity from it), namely, that we likely have $\| \bm{A}_t - 2^{3/2} \mathsf{T}^2\cdot \Id  \| \le 450 n^{2/3}$; this quickly follows from the same bound on $\bm{U}$ and $\bm{V}$, which hold by the path location estimate \eqref{jljj}. Thus, by \eqref{xta0}, the law of $\bm{\mathsf{x}}(t)$ is given by Dyson Brownian motion, run for time $S$, on initial data supported on an interval of width $900 n^{2/3}$. 

For $t \in [-An^{1/3}, An^{1/3}] = [-\mathsf{T}, \mathsf{T}]$, we have $S = \mathsf{T} - t^2 / 4 \mathsf{T} \ge 3 \mathsf{T} / 4 > An^{1/3} / 2$. So, for $A$ large this amounts to running Dyson Brownian motion for a long time, on initial data supported on a bounded interval.\footnote{Under our normalization of Dyson Brownian motion, $n^{1/3}$ and $n^{2/3}$ are the natural scalings of time and space, respectively. So, we only take the prefactors $A / 2$ and $3000$ into account when using the words,  ``long'' and ``bounded.''} It is known in this context that the first $n$ particles equilibrate to have gaps likely satisfying \eqref{ltiljji0} (for example, this sort of statement was shown by Lee--Schnelli \cite{EUDM}; the slightly improved formulation we use appeared in \cite{ERGIM}), implying the gap upper bound. See \Cref{ProbabilityScale} below (which also includes some H\"{o}lder regularity bounds and improvements of the path location estimates, which we do not discuss here but will be useful later in the paper) for further details.

\subsection{Global Law and Regularity} 

\label{RegularGlobal0} 

The proofs of the global law and regularity are based on the notion that non-intersecting Brownian bridges without lower and upper boundaries are simpler to analyze than those with them; the relation \eqref{xta0} to Dyson Brownian motion already provides one manifestation of this phenomenon. To realize this idea, we will restrict the ensemble $\bm{\mathcal{L}}$ to a tall rectangle, giving rise to a family of non-intersecting Brownian bridges with a lower (and no upper) boundary; we will then implement two tasks. The first is to introduce a coupling that compares a family of non-intersecting Brownian bridges on a tall rectangle with lower boundary, to one with the same starting and ending data but without a boundary; we refer to it as the boundary removal coupling. The second is to prove versions of the global law and spatial regularity for non-intersecting Brownian bridges, without boundary, on a tall rectangle. For the global law, the latter will require regularity estimates at the edge of certain limit shapes; we explain this first. 

From this point (particularly in this \Cref{RegularGlobal0} and the next \Cref{Derivative2Approximate}), this proof outline will begin becoming more analytically involved.

\subsubsection{Limit Shapes Near the Edge} 	

\label{ShapeEdge0} 

It has been known since works of Guionnet and Zeitouni \cite{LDASI,FOAMI} (proving earlier predictions of Matytsin \cite{LLI}) that non-intersecting Brownian bridges, without upper and lower boundaries, exhibit a limit shape phenomenon\footnote{By \eqref{xta0}, this amounted to a result on Dyson Brownian motion, namely a large deviations principle for it.} in the following sense (see \Cref{rhot} below for a more precise statement, under a slightly different normalization). Fix real numbers $a < b$ (which will act as times) and $R > 0$ (which will parameterize the number of Brownian bridges). For each integer $n \ge 1$, let $\bm{u}^n = (u_1, u_2, \ldots , u_{Rn})$ and $\bm{v}^n = (v_1, v_2, \ldots , v_{Rn})$ denote $(Rn)$-tuples, such that $n^{-2/3} \cdot \bm{u}$ and $n^{-2/3} \cdot \bm{v}$ both converge to some given profiles. Then letting $\bm{\mathsf{x}}^n = (\mathsf{x}_1^n, \mathsf{x}_2^n, \ldots , \mathsf{x}_{Rn}^n)$ denote $Rn$ non-intersecting Brownian bridges on $[an^{1/3}, bn^{1/3}]$, starting at $\bm{u}$ and ending at $\bm{v}$, their rescaled trajectories\footnote{In fullest generality, it is technically only the associated height function that converges in this way, but in this introductory exposition we ignore that subtlety (which becomes irrelevant in the presence of the gap upper bound).} $n^{-2/3} \cdot \mathsf{x}_j (tn^{1/3})$ converge to a limit shape $G(t, jn^{-1})$, for each $(t, j) \in [a, b] \times [1, Rn]$. Some properties of this function $G: [a, b] \times [0, R] \rightarrow \mathbb{R}$ are known (see \Cref{LimitBridges} below for further details), for example that it satisfies a partial differential equation on the region where it is smooth, given by
\begin{flalign}
	\label{gyy4gt}
	\partial_y^2 G + (\partial_y G)^{-4} \cdot \partial_t^2 G = 0.
\end{flalign}

The assumption that $n^{-2/3} \cdot \bm{u}$ and $n^{-2/3} \cdot \bm{v}$ converge to given profiles already underscores the relevance of the on-scale estimates from \Cref{LProof1}. Setting $\bm{u} = \big( \mathcal{L}_1 (an^{1/3}), \ldots , \mathcal{L}_{Rn} (an^{1/3}) \big)$ and $\bm{v} = \big( \mathcal{L}_1 (bn^{1/3}), \ldots , \mathcal{L}_{Rn} (bn^{1/3}) \big)$, one can only hope for $n^{-2/3} \cdot \bm{u}$ and $n^{-2/3} \cdot \bm{v}$ to have (subsequential) limits if some form of \eqref{jljj} holds (perhaps with different constants). The conclusion above here is also reminiscent of the global law from \Cref{LProof1}; both provide deterministic approximations for the Brownian paths, up to error $o(n^{2/3})$. However, the deterministic approximation in the global law had an exact formula, but here it is given by the less transparent function $G$.

While the full limit shape $G$ is usually quite inexplicit indeed, we will show under certain conditions that it admits a ``universal behavior'' near the edge $y = 0$ (corresponding to the top curves in $\bm{\mathsf{x}}$). Specifically, for fixed $\mathfrak{t} \in (a, b)$, there exist constants $\mathfrak{a}, \mathfrak{b} \in \mathbb{R}$ and $\mathfrak{c} > 0$ such that 
\begin{flalign}
	\label{g0abc} 
	G(t, y) \approx \mathfrak{a} + \mathfrak{b} t - \mathfrak{c} t^2 - 2^{-4/3} (3\pi)^{2/3} \mathfrak{c}^{-1/3} y^{2/3}, \qquad \text{for $(t, y) \approx (\mathfrak{t}, 0)$}. 
\end{flalign}

\noindent It will be central for the approximation error in \eqref{g0abc} to remain uniformly small as the parameter $R$ grows, which we will ultimately take to be large (to later compare $\bm{\mathcal{L}}$ to a system of Brownian bridges without boundaries); see \Cref{p:limitprofile} below (where the $R$ here is $L^{3/2}$ there). 

The proof of \eqref{g0abc} is based on a purely deterministic analysis of the limit shape $G$ and the associated partial differential equation \eqref{gyy4gt}. Non-uniformly elliptic equations similar to \eqref{gyy4gt} (though different in that they were constrained to be Lipschitz, as they arose from limit shapes of dimer models) were analyzed from a real analytic perspective by De Silva--Savin \cite{MCFARS} and from a complex analytic one by Kenyon--Okounkov \cite{LSCE} and Astala--Duse--Prause--Zhong \cite{DMC}. In that setting, the last work \cite{DMC} proved a variant of \eqref{g0abc}, though they did not investigate the uniformity of that approximation in the size of the underlying domain.

In our context, this uniformity is in fact false in general. To witness it, we must impose hypotheses on the boundary data for $G$ (\Cref{integralmu0mu1} and \Cref{gapmu0mu1}), stipulating the existence of a constant $C > 1$ so that for each boundary point $s \in \{ a, b \}$ we have 
\begin{flalign}
	\label{cg2} 
	-C - Cy^{2/3} \le G(s, y) \le C - C^{-1} y^{2/3}, \quad \text{and} \quad \big| G(s, y) - G(s, y') \big| \le C \big|y^{2/3} - (y')^{2/3} \big|,
\end{flalign}

\noindent for each $y, y' \in [0, R]$. Observe that two bounds in \eqref{cg2} constitute continuum counterparts of the two on-scale estimates (path locations and gap upper bound) from \Cref{LProof1}.

The proof of \eqref{g0abc} first requires an \emph{a priori} estimate on how $\partial_y G(t,y)$ diverges for small $y$, namely, $\partial_y G (t, y) \sim -y^{-1/3}$, where the implicit constants are uniform in $R$ (\Cref{p:checka}). To verify the upper bound on $|\partial_y G|$, we show a continuous variant (\Cref{limitdifferencecompare}) of gap monotonicity, indicating that the $y$-derivatives of limit shapes are increasing in those of their boundary data. Since \eqref{cg2} upper bounds $|\partial_y G|$ on the boundary, we can use this to upper bound the $y$-derivative of $G$ by that of an explicit limit shape, which can directly be seen to have the $y^{-1/3}$ divergence.

To lower bound $|\partial_y G|$, we instead use the following property about limit shapes \cite{FOAMI}. Let $\varrho = -(\partial_y G)^{-1}$, fix $t \in (a, b)$, and set $\tau = (b-t)(t-a) / (b-a)$. Then $\varrho (t, \cdot)$ is the density for a measure $\nu_{\tau}$, given by the free convolution between some measure $\nu$ of total mass $R$ and the semicircle distribution of size $\tau$. This is specific to limit shapes for Brownian bridges without upper and lower boundaries, and it can be viewed as a continuum counterpart of \eqref{xta0}. Similarly to \eqref{xta0}, effectively using this fact is complicated by the fact that little is in general known about $\nu$. 

So, we develop a general estimate for such measures when $\tau \gtrsim 1$ stating that, if the first bound in \eqref{cg2} holds at $s=t$, then $\varrho \lesssim 1$ holds uniformly in $R = \nu(\mathbb{R})$ for $y \le 1$ (\Cref{p:densityest}). While the former bound in \eqref{cg2} was only stipulated to hold at $s \in \{a, b \}$, it can be shown (by a continuous variant of height monotonicity, \Cref{limitheightcompare}) to extend to $s \in [a, b]$. This verifies the assumption in the above free convolution result, yielding for $y \le 1$ that $\varrho \lesssim 1$, and so $|\partial_y G| \gtrsim 1$. Improving this bound to $|\partial_y G| \gtrsim y^{-1/3}$ requires further effort (involving elliptic regularity and another application of the continuum gap monotonicity). See \Cref{s:density} below for further details.

Given the above, to establish \eqref{g0abc}, we next use the fact \cite{LLI} that \eqref{gyy4gt} can be equivalently written as a complex Burgers equation for the complex slope $f = \partial_t G - \mathrm{i} \cdot (\partial_y G)^{-1}$, providing $f$ a holomorphic structure; such ideas were also fruitful in prior works \cite{LSCE,DMC} analyzing dimer limit shapes. In particular, defining the complex coordinate $z(t, y) = y - t \cdot f(t, y)$, this indicates that $f = F(z)$, for some holomorphic function $F$. We show that the previously mentioned bounds for $\partial_y G$ imply uniform derivative estimates for $F$ (\Cref{p:analytic}), enabling $F$ to be Taylor expanded. Translating this expansion for $F$ into one for $G$ eventually yields the approximation \eqref{g0abc}. See \Cref{s:shape1} below for further details.

\subsubsection{Boundary Removal Coupling} 

\label{CoupleBoundary2}

\begin{figure}
	\center
\includegraphics[scale=.6]{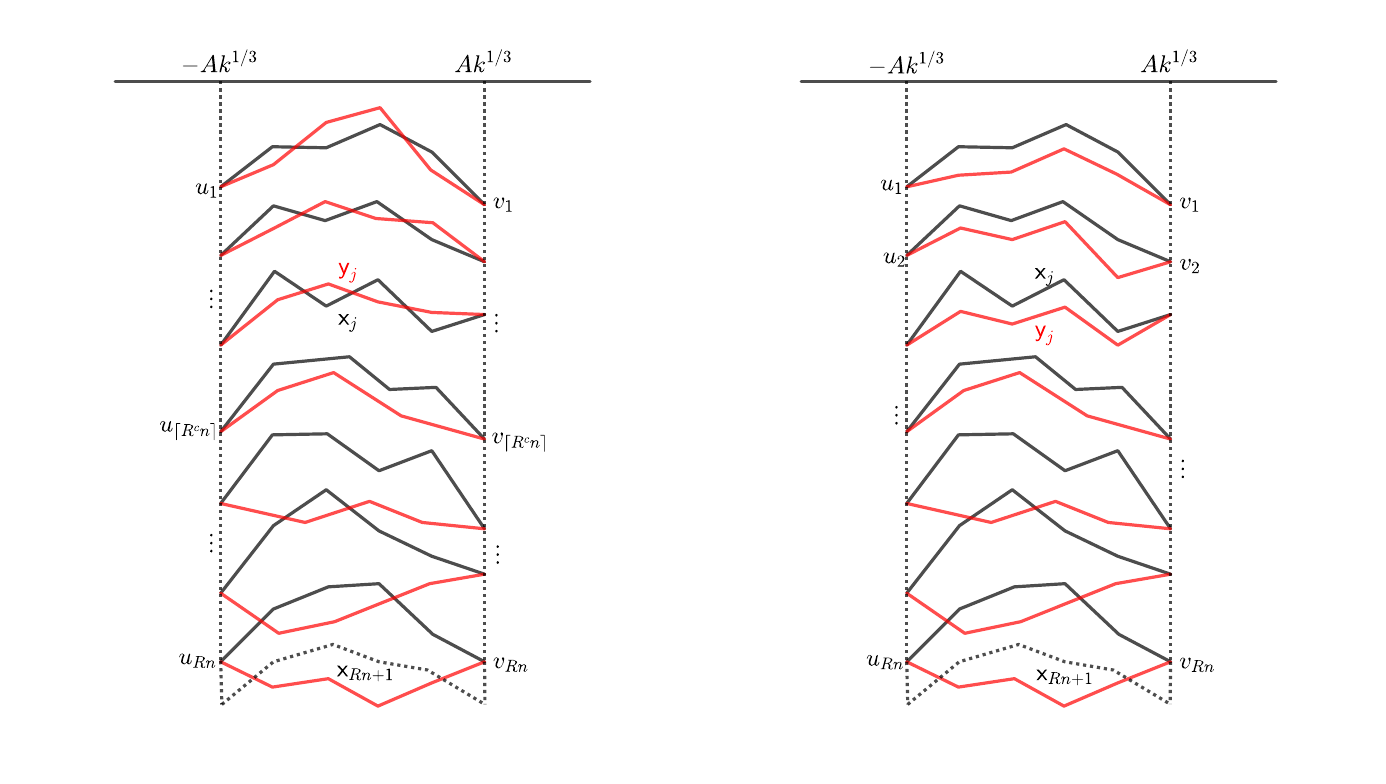}

\caption{Depicted above is the boundary removal coupling.}
\label{f:remove_boundary2}
	\end{figure}

The boundary removal coupling can be described as follows (see \Cref{c:finalcouple} below for the precise statement, where the $R$ here is $L^{3/2}$ there). Fix a bounded real number $A \ge 1$; let $\bm{\mathsf{x}} = (\mathsf{x}_1, \mathsf{x}_2, \ldots )$ denote a Brownian Gibbsian line ensemble likely satisfying path location bounds of the type \eqref{jljj}; and let $R > 1$ be a real number and $n \ge 1$ be an integer (that can now be arbitrarily large with respect to each other). Sample $Rn$ non-intersecting Brownian bridges $\bm{\mathsf{y}} = (\mathsf{y}_1, \mathsf{y}_2, \ldots , \mathsf{y}_{Rn})$ on $[-An^{1/3}, An^{1/3}]$, starting at $\bm{u} = \big( \mathsf{x}_1 (-An^{1/3}), \ldots , \mathsf{x}_{Rn} (-An^{1/3}) \big)$ and ending at $\bm{v} = \big( \mathsf{x}_1 (An^{1/3}), \ldots , \mathsf{x}_{Rn} (An^{1/3}) \big)$. As $\bm{\mathsf{y}}$ has no lower boundary, height monotonicity (recall \Cref{Gap0}) yields a coupling between $\bm{\mathsf{x}}$ and $\bm{\mathsf{y}}$ such that $\mathsf{x}_j \ge \mathsf{y}_j$ almost surely, for each $j \in [1, Rn]$; see the right side of \Cref{f:remove_boundary2} for a depiciton. We will show that there exists a coupling satisfying a bound in the reverse direction for small indices $j$, that is, for $c = 2^{-4000}$ we have with high probability that 
\begin{flalign} 
	\label{yjxjrcn} 
	\mathsf{y}_j \ge \mathsf{x}_j - R^{-c} n^{2/3}, \qquad \text{for each $j \in [1, R^c n]$.}
\end{flalign} 

\noindent See the left side of \Cref{f:remove_boundary2} for a depiction. Together, these couplings suggest that the top $R^c n$ curves in $\bm{\mathsf{x}}$ and $\bm{\mathsf{y}}$ should be close, for large $R$.

To exhibit the boundary removal coupling, we first reduce it to a ``preliminary coupling'' that introduces a lower boundary $f : [-An^{1/3}, An^{1/3}] \rightarrow \mathbb{R}$ for $\bm{\mathsf{y}}$. It essentially states the following (see \Cref{p:comparison2} below for the precise formulation). Assume that with high probability the path location estimates of the type \eqref{jljj} hold for $\bm{\mathsf{x}}$ and moreover that, (i) $\mathsf{x}_{Rn+1}$ is not too far above $f$, namely, $f \ge \mathsf{x}_{Rn+1} - (R^{\alpha} n)^{2/3}$ for some $\alpha < 1$ and, (ii) its paths $\mathsf{x}_j$ are regular, namely, they satisfy a H\"{o}lder bound that is governed by a parameter $\beta \in (0, 1)$ in a specific way (see \Cref{eventregularityimproved} below), where smaller $\beta$ implies improved regularity. Then, the preliminary coupling between $\bm{\mathsf{x}}$ and $\bm{\mathsf{y}}$ ensures for $c_0 = 2^{-3500}$ that $\mathsf{y}_j \ge \mathsf{x}_j - R^{c_0 (2\beta-7/8)} n^{2/3}$ likely holds,\footnote{Here, the constant $2\beta - 7 / 8$ in the exponent is not optimal. It can at least be improved to $2\beta-1^-$, but we do not know what the optimal constant should be.} for each $j \in [1, R^{c_0} n]$. 

For the original line ensemble $\bm{\mathsf{x}}$ in the boundary removal coupling, we will show that (ii) holds at $\beta = 3 / 8$ (\Cref{p:Lip}), so $c_0 (2\beta - 7 / 8) = -c_0 / 8 < -c$, yielding the negative exponent in \eqref{yjxjrcn}. This $\beta = 3 / 8$ result is established in \Cref{RegularImproved} below, and its proof amounts to an inductive series of comparisons between $\bm{\mathsf{x}}$ and Dyson Brownian motion, the latter of whose regularity properties are well understood (some of which are collected in \Cref{s:linelong} below). Since $\bm{\mathsf{y}}$ has no lower boundary, (i) cannot be literally be true as written, but we may view its bottom path $\mathsf{y}_{Rn}$ as a lower boundary for its remaining $Rn-1$ ones $(\mathsf{y}_1, \mathsf{y}_2, \ldots , \mathsf{y}_{Rn-1})$. A weak estimate on how much $\mathsf{x}_{Rn}$ and $\mathsf{y}_{Rn}$ can oscillate on the interval $[-An^{1/3}, An^{1/3}]$ (recall $\bm{\mathsf{x}}$ and $\bm{\mathsf{y}}$ share the same endpoints) will suffice to show that $|\mathsf{x}_{Rn} - \mathsf{y}_{Rn}| \le (R^{\alpha} n)^{2/3}$ holds with high probability for some $\alpha<1$. This enables us to deduce the boundary removal coupling as a consequence of the preliminary one, applied to $(\mathsf{y}_1, \mathsf{y}_2, \ldots , \mathsf{y}_{Rn-1})$ with lower boundary $f = \mathsf{y}_{Rn}$. See \Cref{RectangleCouple} below for further details. 

\begin{figure}
	\center
\includegraphics[scale=.8]{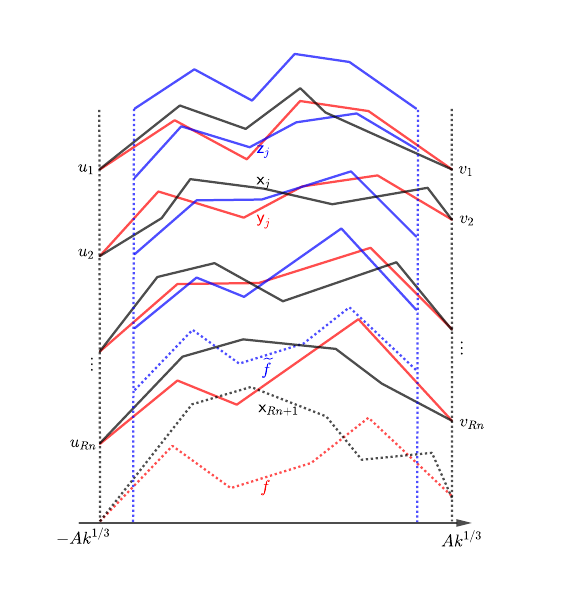}

\caption{Shown above is a depiction for the proof of the preliminary coupling described in \Cref{CoupleBoundary2}.}
\label{f:rescaling2}
	\end{figure}

It remains to prove that the preliminary coupling exists, which can heuristically be explained through a diffusive scaling. Fix some parameter $\vartheta \sim R^{2(\alpha-1)/3}$ and define the non-intersecting Brownian bridges $\bm{\mathsf{z}} = (\mathsf{z}_1, \mathsf{z}_2, \ldots , \mathsf{z}_{Rn})$ with lower boundary $\widetilde{f}$, by first diffusively ``shrinking'' $\bm{\mathsf{y}}$ (with lower boundary $f$) by factors of $1 - \vartheta$ in space and $(1-\vartheta)^2$ in time, and then by translating them up slightly. See \Cref{f:rescaling2} for a depiction.
This has two effects. First, it can be shown that the original lower boundary $f \sim -(Rn)^{2/3}$ is quite low, so shrinking it to $\widetilde{f}$ lifts it up considerably, by $\vartheta \cdot |f| \sim (R^{\alpha} n)^{2/3}$. Due to (i), this (upon the proper choice of constants) pulls this lower boundary above $\mathsf{x}_{Rn+1}$. Second, it changes the time domain of the bridges slightly, from $[-An^{1/3}, An^{1/3}]$ to $\big[ -A n^{1/3} (1-\vartheta)^2, A n^{1/3} (1-\vartheta)^2 \big]$. By the H\"{o}lder bound (ii), the paths in $\bm{\mathsf{x}}$ cannot increase much when passing from the former interval to the latter one. Hence, the starting and ending data for $\bm{\mathsf{z}}$ on $\big[ -An^{1/3} (1-\vartheta)^2, An^{1/3} (1-\vartheta)^2 \big]$ likely continues to dominate that of $\bm{\mathsf{x}}$. 

Height monotonicity then provides a coupling between $\bm{\mathsf{x}}$ and $\bm{\mathsf{z}}$ such that $\mathsf{z}_j \ge \mathsf{x}_j$ for each $j \in [1, Rn]$. Since the top curves of $\bm{\mathsf{y}}$ are only barely perturbed under the shrinking to $\bm{\mathsf{z}}$, this coupling lower bounds the top curves of $\bm{\mathsf{y}}$ by those of $\bm{\mathsf{x}}$. More specifically, we can deduce for some explicit constants $d = d(\alpha, \beta) > 0$ and $\Delta = \Delta (\alpha, \beta) > 0$ that $\mathsf{y}_j \ge \mathsf{x}_j - (R^{d(\alpha - \Delta)} n)^{2/3}$  for each $j \in [1, R^d n+1]$ (\Cref{p:comparison1}). Replacing $R$ by $R^d$, this effectively reduces the original exponent $\alpha$ appearing in the preliminary coupling to $\alpha - \Delta$. While this procedure might not immediately yield an exponent of $2\beta - 7 / 8$, it can be applied recursively. By repeating it many times, we can reduce $\alpha$ to just above the value $\alpha_0$ where $\Delta (\alpha_0, \beta) \ge 1/250$ (or is bounded below by any small positive constant). A calculation reveals that $\alpha_0 < 2\beta - 7 / 8$, which yields the exponent stated in the preliminary coupling. See \Cref{Couple0Proof} below for further details. 

\subsubsection{Global Law}

\begin{figure}
	\center
\includegraphics[scale=.5, trim=1cm 4cm 1cm 4cm, clip]{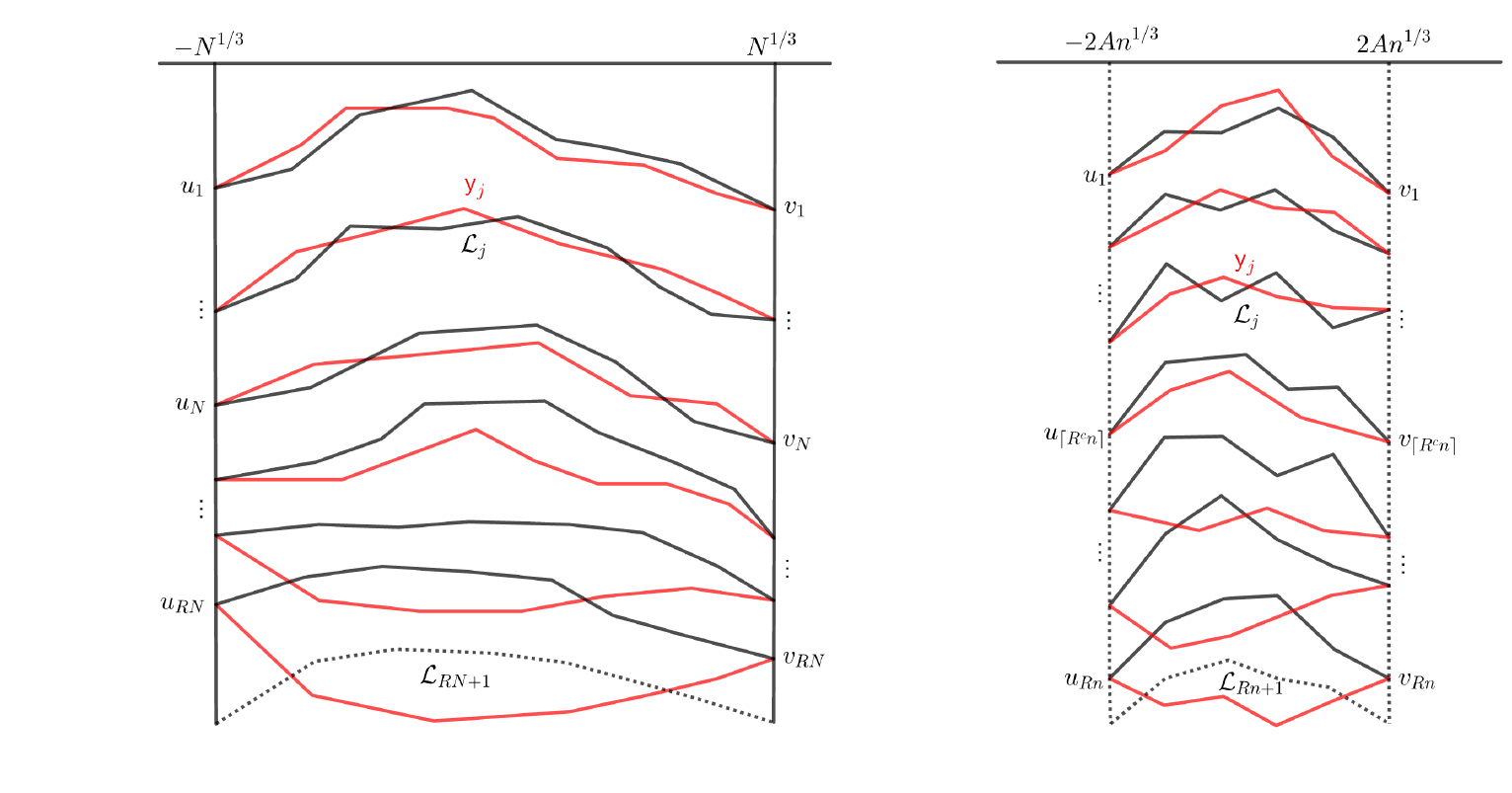}

\caption{Shown to the left is a depiction of how the boundary removal coupling is used to prove the global law, where $R$ is large but fixed in $n$. Shown to the right is a depiction of how it is used to prove the spatial regularity, where $R$ is growing much faster than $n$. In both cases, the boundary removal coupling is used to compare the original black curves $(\mathcal{L}_1, \mathcal{L}_2, \ldots , \mathcal{L}_{Rn})$, which have the lower boundary $\mathcal{L}_{Rn+1}$, to a family of red ones without any bottom boundary.}
\label{f:global_spatial}
	\end{figure}

The global law from \Cref{LProof1} states that, for any fixed $\delta > 0$,
\begin{flalign}
	\label{lnt22}
	\big| \mathcal{L}_n (t) + 2^{-1/2} t^2 + 2^{-7/6} (3\pi)^{2/3} n^{2/3} \big| < \delta n^{2/3},
\end{flalign}

\noindent likely holds. To establish it, we let $\theta \in (0, 1)$ and $R > 1$ be small and very large (relative to $\delta$) real numbers, respectively, both of which are fixed in $n$ (we also assume $R$ is much larger than $\theta^{-1}$). Letting $N = n / \theta^3$, we then apply the boundary removal coupling (recall \Cref{CoupleBoundary2}) to the top $RN$ paths of $\bm{\mathcal{L}}$. Sampling $RN$ non-intersecting Brownian bridges without boundaries $\bm{\mathsf{y}} = (\mathsf{y}_1, \mathsf{y}_2, \ldots , \mathsf{y}_{RN})$ on $[-N^{1/3}, N^{1/3}]$, starting at  $\bm{u} = \big( \mathcal{L}_1 (-N^{1/3}), \ldots , \mathcal{L}_{RN} (-N^{1/3}) \big)$ and ending at $\bm{v} = \big( \mathcal{L}_1 (N^{1/3}), \ldots , \mathcal{L}_{RN} (N^{1/3}) \big)$, this enables us to approximate the upper $N \ge n$ paths $\bm{\mathcal{L}}$ by those of $\bm{\mathsf{y}}$, up to an error of $R^{-c} N^{2/3} \ll \delta n^{2/3}$. See the left side of \Cref{f:global_spatial} for a depiction.

We next apply limit shape results from \cite{LDASI,FOAMI} to $\bm{\mathsf{y}}$ (recall \Cref{ShapeEdge0}). Although they were originally only stated for infinite sequences of Brownian bridge ensembles with starting and ending data that ``converge,'' a compactness argument can be used to show a finite variant of it (\Cref{pluconverge}). This provides a limit shape $G : [-1, 1] \times [0, 1] \rightarrow \mathbb{R}$ such that we likely have $\big| \mathsf{y}_j (tN^{1/3}) - N^{2/3} \cdot G(t, j / N) \big| \ll \delta n^{1/3}$. Our edge behavior result described in \Cref{ShapeEdge0}, whose hypothesis \eqref{cg2} can be verified by the two on-scale estimates from \Cref{LProof1}, then applies to this limit shape $G$ and yields (at this point, unknown) constants $(\mathfrak{a}, \mathfrak{b}, \mathfrak{c})$ so that 
\begin{flalign*} 
	G(s, y) \approx \mathfrak{a} + \mathfrak{b} s - \mathfrak{c} s^2 - 2^{-4/3} (3\pi)^{2/3} \mathfrak{c}^{-1/3} y^{2/3}, \qquad \text{for small $(s, y) \in [-\theta, \theta] \times [0, \theta]$}.
\end{flalign*} 

\noindent  Combining this at $(s, y) = (\theta t,  j / N)$ with the previous statements, we deduce 
\begin{flalign*} 
	\big| \mathcal{L}_j (tn^{1/3}) -n^{2/3} \cdot ( \theta^{-2} \mathfrak{a} + \theta^{-1} \mathfrak{b} t - \mathfrak{c} t^2) + 2^{-4/3} (3\pi)^{2/3} \mathfrak{c}^{-1/3} j^{2/3} \big| \ll \delta n^{2/3},
\end{flalign*} 

\noindent likely holds for $j \le n \le \theta N$. Matching this against the behavior \eqref{l1t} imposed on the top curve $j=1$, we obtain $(\mathfrak{a}, \mathfrak{b}, \mathfrak{c}) = (0, 0, 2^{-1/2})$, giving \eqref{lnt22}. See \Cref{ProofGlobal} below for further details.

\subsubsection{Spatial Regularity} 

\label{LRegular}

For any $t \in [-A, A]$, the spatial regularity from \Cref{LProof1} exhibits a random, almost smooth (say, with bounded first $50$ derivatives) function $\gamma_t : [0, 1] \rightarrow \mathbb{R}$ such that 
\begin{flalign}
\label{ljn1} 
	\big| \mathcal{L}_{j+n} (tn^{1/3}) - n^{2/3} \cdot \gamma_t (jn^{-1}) \big| = o(1), \qquad \text{likely holds for each $1 \le j \le n$}.
\end{flalign} 

\noindent This provides a much stronger approximation than the global law \eqref{lnt22}, at the expense of making the approximating function $\gamma_t$ less explicit. Its proof again makes use of the boundary removal coupling (recall \Cref{CoupleBoundary2}), but now takes $R \gg n^{2/c}$ (at $c = 2^{-4000}$) to grow much faster than $n$. Sampling $Rn$ non-intersecting Brownian bridges $\bm{\mathsf{y}} = (\mathsf{y}_1, \mathsf{y}_2, \ldots , \mathsf{y}_{Rn})$ on $[-2An^{1/3}, 2An^{1/3}]$, starting at $\bm{u} = \big( \mathcal{L}_1 (-2An^{1/3}), \ldots , \mathcal{L}_{Rn} (-2An^{1/3}) \big)$ and ending at $\bm{v} = \big( \mathcal{L}_1 (2An^{1/3}), \ldots , \mathcal{L}_{Rn} (2An^{1/3}) \big)$, this approximates the upper $R^c n \gg 2n$ paths of $\bm{\mathcal{L}}$ by those of $\bm{\mathsf{y}}$, up to an error of $R^{-c} n^{2/3} = o(1)$. See the right side of \Cref{f:global_spatial} for a depiction. It thus suffices prove spatial regularity for $\bm{\mathsf{y}}$.

  The benefit in this is that, since $\bm{\mathsf{y}}$ does not have an upper or lower boundary, it admits the interpretation \eqref{xta0} in terms of Dyson Brownian motion. Rigidity results of Huang--Landon \cite{MCTMGP} apply to the latter and imply that the $\mathsf{y}_j (t)$ closely concentrate, up to error $o(1)$, around the quantiles of a measure $\nu_{\tau}$, given by the free convolution between some measure $\nu$ of total mass $R$ and the semicircle distribution of size $\tau = A - t^2 / 4A \gtrsim 1$. The spatial regularity of $\bm{\mathsf{y}}$ then amounts to ensuring that the density for this measure $\nu_{\tau}$ is almost smooth, but this is once again complicated by the fact that little is known about $\nu$. So, we develop a general result about such measures $\nu_{\tau}$, closely related to the one described in \Cref{ShapeEdge0}. It states that, under certain conditions (which can in our context can be later verified by the two on-scale estimates described in \Cref{LProof1}), when $\tau \gtrsim 1$ the derivatives of the density for $\nu_{\tau}$ are uniformly bounded in $R = \nu(\mathbb{R})$ (\Cref{c:rhoderbound}). This confirms the spatial regularity for $\bm{\mathsf{y}}$ and thus for $\bm{\mathcal{L}}$. See \Cref{ProofRegular} below for further details.

\subsection{Curvature Approximation}

\label{Derivative2Approximate}

The curvature approximation from \Cref{LProof1} exhibits a random, twice-differentiable function $h_n : [-An^{1/3}, An^{1/3}] \rightarrow \mathbb{R}$ so that we likely have 
\begin{flalign} 
	\label{hn2} 
	\big| h_n'' (s) + 2^{-1/2} \big| = o(1), \quad \text{and} \quad  \big| h_n (s) - \mathcal{L}_n (s) \big| = o(1), \qquad  \text{for all $s \in [-An^{1/3}, An^{1/3}]$}.
\end{flalign} 
\noindent As for spatial regularity (recall \Cref{LRegular}), it provides a stronger approximation than the global law \eqref{lnt22}, though with a less explicit approximating function $h_n$. To establish it, we make use of a concentration bound for non-intersecting Brownian bridges with smooth boundary data, proven in \cite[Sections 4-7]{U} (based on carefully ``patching together'' concentration bounds for Dyson Brownian motion, inspired by ideas of Laslier--Toninelli \cite{TDMS}). That bound can be described as follows (see \Cref{gh} below for a more precise statement, under a slightly different normalization). 

Let $a < b$ be real numbers; let $k \ge 1$ be a large integer (which we view as tending to $\infty$); and let $\bm{\mathsf{x}} = (\mathsf{x}_1, \mathsf{x}_2, \ldots , \mathsf{x}_k)$ denote non-intersecting Brownian bridges on $[ak^{1/3}, bk^{1/3}]$, starting at $\bm{u} = (u_1, u_2, \ldots , u_k)$; ending at $\bm{v} = (v_1, v_2, \ldots , v_k)$; and conditioned to lie above and below functions $f : [ak^{1/3}, bk^{1/3}] \rightarrow \mathbb{R}$ and $g : [ak^{1/3}, bk^{1/3}] \rightarrow \mathbb{R}$, respectively. Assume that there is an almost smooth (with bounded first $50$ derivatives) solution $G : [a, b] \times [0, 1] \rightarrow \mathbb{R}$ to the limit shape partial differential equation \eqref{gyy4gt}, which is close to $(\bm{u}; \bm{v}; f; g)$ along the boundary, namely, for each $j \in [1, k]$ and $t \in [a, b]$ we have
\begin{flalign}
	\label{uvghk}
	\begin{aligned} 
	& \big| k^{2/3} \cdot G(a, jk^{-1}) - u_j \big| = o(1); \qquad \quad \big| k^{2/3} \cdot G(b, jk^{-1}) - v_j \big| = o(1); \\
	&  k^{2/3} \cdot G(t, 0) = g(tk^{1/3}); \qquad \qquad \qquad k^{2/3} \cdot G(t, 1) = f(tk^{1/3}).
	\end{aligned}
\end{flalign}

\noindent Then, $\big| \mathsf{x}_j (tk^{1/3}) - k^{2/3} \cdot G (t, jk^{-1}) \big| = o(1)$ holds for all $(t, j) \in [a, b] \times [1, k]$, with high probability. 

The condition \eqref{uvghk} can be viewed as a constraint on the boundary data $(\bm{u}; \bm{v}; f; g)$ for $\bm{\mathsf{x}}$, as it implies that they must approximate a smooth function $G$. In particular, the first two bounds in \eqref{uvghk} underscore the relevance of spatial regularity for $\bm{\mathcal{L}}$. Fixing $k = 2n / 3$ and letting $\bm{u} = \big( \mathcal{L}_{k+1} (-Ak^{1/3}), \ldots , \mathcal{L}_{2k} (-Ak^{1/3}) \big)$ and $\bm{v} = \big( \mathcal{L}_{k+1} (Ak^{1/3}), \ldots , \mathcal{L}_{2k} (Ak^{1/3}) \big)$, they can only hold if \eqref{ljn1} does (with the $n$ equal to $k$ here) for an almost smooth function $\gamma_t (\cdot) = G(t, \cdot)$, at $t \in \{ -A, A \}$.

\begin{figure}
	\center
\includegraphics[scale=2]{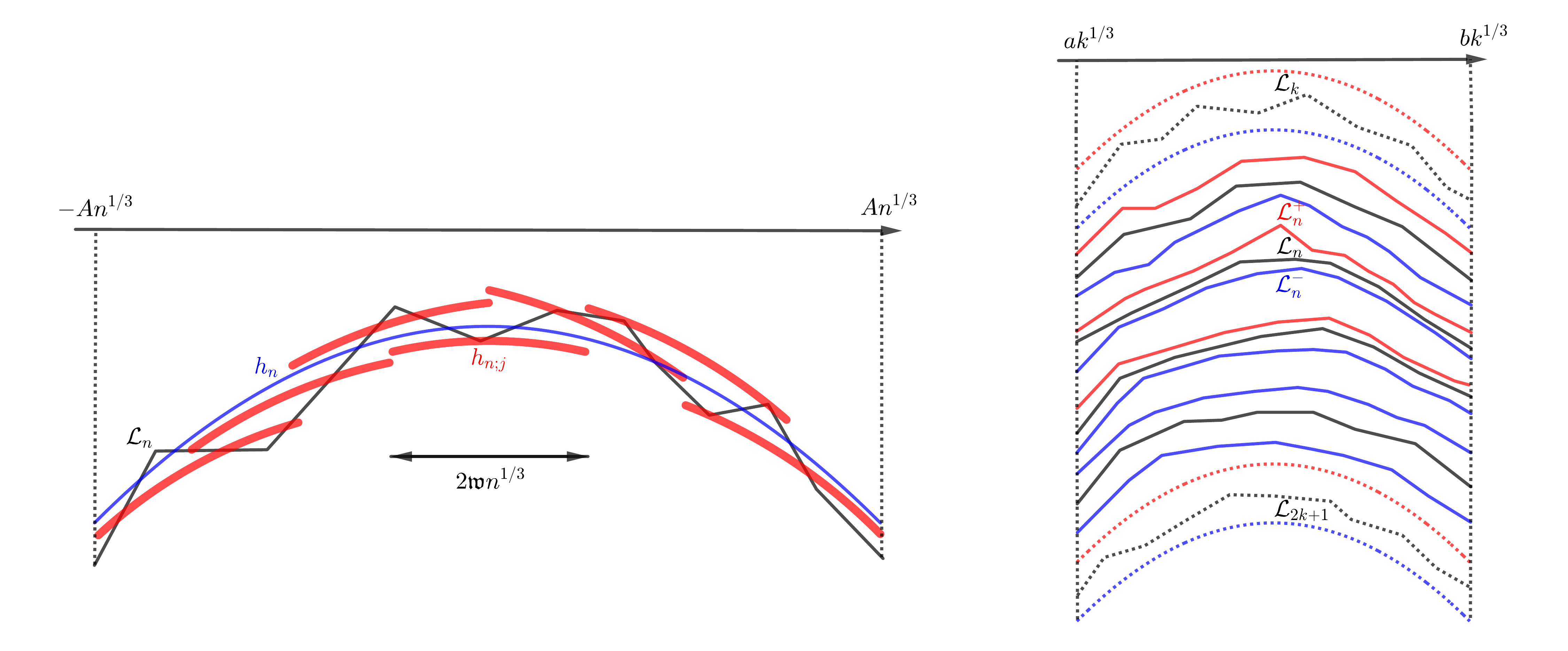}

\caption{Shown to the left is the gluing used to produce a global $h_n$ satisfying \eqref{hn2}. Shown to the right is the sandwiching of $\bm{\mathcal L}$ between $\bm{\mathcal L}^+$ and $\bm{\mathcal L}^-$.}
\label{f:curvature}
	\end{figure}

While this spatial regularity ensures that the starting and ending data $(\bm{u}; \bm{v})$ approximate almost smooth profiles, it makes no such guarantee for the upper and lower boundaries $(f; g)$. In our case, the latter are $(\mathcal{L}_k; \mathcal{L}_{2k+1})$, and we do not know how to directly show that they are close to smooth functions. To circumvent this issue, we let $\eta \in (0, 1)$ be a small parameter and subdivide $[-An^{1/3}, An^{1/3}]$ into subintervals $\{ I_j \}$ of length $2 \eta n^{1/3}$. For each $j \in [ 1, A / \eta]$, we will produce a ``local'' approximating function $h_{n;j}$ likely satisfying \eqref{hn2} on $I_j$ (\Cref{h0xlocal}), and then we will ``glue'' these local $h_{n;j}$ together (see the proof of \Cref{h0x2}) to form a global one satisfying \eqref{hn2} on $[-An^{1/3}, An^{1/3}]$. See the left side of \Cref{f:curvature} for a depiction. 
 This reduces us to proving a version of \eqref{hn2}, but only on short intervals of length $2 \eta n^{1/3}$. The value in this is that, on such thin domains, one does not expect the upper and lower boundaries $f = \mathcal{L}_k$ and $g = \mathcal{L}_{2k+1}$ to substantially affect most of the middle curves. 

To make this precise, we introduce two families $\bm{\mathcal{L}}^-$ and $\bm{\mathcal{L}}^+$ of non-intersecting Brownian bridges and sandwich $\bm{\mathcal{L}}$ between them. Their starting and ending data will nearly coincide with those of $\bm{\mathcal{L}}$, which are almost smooth by spatial regularity \eqref{ljn1}. Their (upper and lower) boundaries will also be smooth, and will lie slightly above and below those of $\bm{\mathcal{L}}$, respectively. See the right side of \Cref{f:curvature} for a depiction. This regularity will enable us to approximate (as in \eqref{uvghk}) the boundary data for $\bm{\mathcal{L}}^-$ and $\bm{\mathcal{L}}^+$ by almost smooth limit shapes $G^-$ and $G^+$ satisfying \eqref{gyy4gt}, respectively. By the global law \eqref{lnt22}, the boundary data for $G^-$ and $G^+$ approximate the function $\mathfrak{G} (t, y) = -2^{-1/2} t^2 - 2^{-7/6} (3\pi)^{2/3} y^{2/3}$. Using properties of solutions to the equation \eqref{gyy4gt}, we will further show (\Cref{fg1g2}) that these limit shapes essentially satisfy, (i) they are extremely close to each other in the middle, so $\big| G^+ (t,1 / 2) - G^- (t, 1 / 2) \big| = o(k^{-2/3})$ and, (ii) they approximate $\mathfrak{G}$ also in their derivatives, so $\partial_t^2 G^+ = \partial_t^2 \mathfrak{G} + o(1) = o(1) -2^{-1/2}$ and similarly $\partial_t^2 G^- = o(1) -2^{-1/2}$.  

 Thus the above concentration bound applies to $\bm{\mathcal{L}}^-$ and $\bm{\mathcal{L}}^+$, which with the sandwiching of $\bm{\mathcal{L}}$ between them gives $k^{2/3} \cdot G^- ( t, 1 / 2) - o(1) \le \mathcal{L}_n^- (tk^{1/3}) \le \mathcal{L}_n (tk^{1/3}) \le \mathcal{L}_n^+ (tk^{1/3}) = k^{2/3} \cdot G^+ ( t, 1 / 2) + o(1)$ with high probability. By (i), the left and right sides of this inequality are within $o(1)$ of each other, and so $\mathcal{L}_n (t k^{1/3}) = k^{2/3} \cdot G^+ (t, 1 / 2) + o(1)$. Then taking $h_{n;j} (t) = k^{2/3} \cdot G^+ (tk^{-1/3}, 1 / 2)$ yields the second statement in \eqref{hn2}. By (ii), we also have $h_{n;j}'' (t) = \partial_t^2 G^+ (tk^{-1/3}, 1 / 2 ) = o(1) - 2^{-1/2}$, confirming the first statement in \eqref{hn2}. See \Cref{DerivativePath2} and \Cref{Proof0Omega} below for further details, where the proper implementation of the above framework involves an induction on scales argument.

\subsection{Airy Statistics}

Although we now have the curvature approximation \eqref{hn2}, we are not yet able to directly compare $\bm{\mathcal{L}}$ to the scaled parabolic Airy line ensemble $\bm{\mathcal{S}}$. Indeed, the bound there $|h_n'' + 2^{-1/2}| = o(1)$ on the approximating function $h_n$ still in principle allows it to have large oscillations, of sizes up to  $o(n^{2/3})$ on the time interval $[-n^{1/3}, n^{1/3}]$; these can already dominate the fluctuations of $\bm{\mathcal{S}}$. We instead first pin down a more robust family of statistics for $\bm{\mathcal{L}}$ (as opposed to its entire law), given by its gaps $\big(\mathcal{L}_1 (t) - \mathcal{L}_2 (t), \mathcal{L}_2 (t) - \mathcal{L}_3 (t), \ldots \big)$ at a given time $t$.

\subsubsection{Airy Gaps} 

The first aspect of Airy statistics, as described in \Cref{LProof1}, states that for any fixed $t \in \mathbb{R}$ the law of $\big( \mathcal{L}_1 (t) - \mathcal{L}_2 (t), \mathcal{L}_2 (t) - \mathcal{L}_3 (t), \ldots \big)$ coincides with that of $\big( \mathcal{S}_1 (0) - \mathcal{S}_2 (0), \mathcal{S}_2 (0) - \mathcal{S}_3 (0), \ldots \big)$. To establish this, we show that the former stochastically is lower bounded by the latter (\Cref{xjxj11}), and also is stochastically upper bounded by it (\Cref{xjxj12}). The proofs of both use gap monotonicity (recall \Cref{Gap0}) in different ways. 

\begin{figure}
	\center
\includegraphics[scale=1, trim=1cm 0.5cm 0 0, clip]{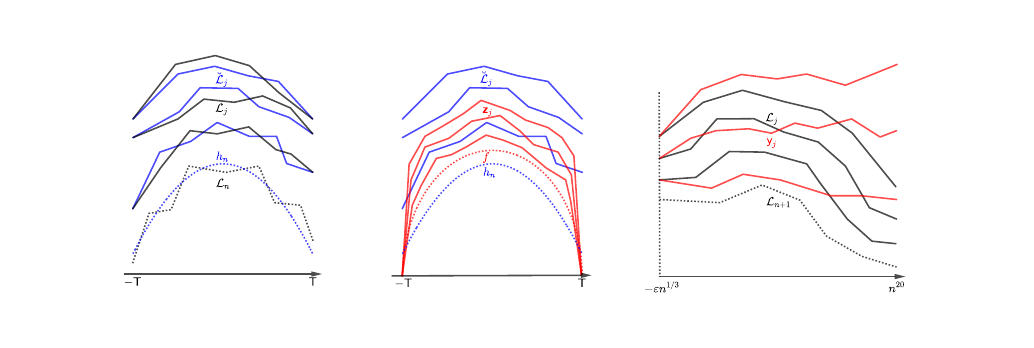}

\caption{Shown to the left and middle is a depiction for the proof of the Airy gap lower bound for $\bm{\mathcal L}$. Shown to the right is a depiction for the proof of the corresponding upper bound.}
\label{f:Airy_Gaps}
	\end{figure}

To prove the lower bound, we fix a large integer $n \gg 1$ and make use of the curvature approximation \eqref{hn2}, recalling the function $h_n$ appearing there. Let $A \ge 1$ be a large real number, bounded in $n$, and set $\mathsf{T} = An^{1/3}$. Sample $n-1$ non-intersecting Brownian bridges $\breve{\bm{\mathcal{L}}} = (\breve{\mathcal{L}}_1, \breve{\mathcal{L}}_2, \ldots , \breve{\mathcal{L}}_{n-1})$ on $[-\mathsf{T}, \mathsf{T}]$, starting at $\bm{u} = \big( \mathcal{L}_1 (-\mathsf{T}), \ldots , \mathcal{L}_{n-1} (-\mathsf{T}) \big)$; ending at $\bm{v} = \big( \mathcal{L}_1 (\mathsf{T}), \ldots , \mathcal{L}_{n-1}  (\mathsf{T}) \big)$; and conditioned to lie above $h_n$. Since $(\mathcal{L}_1, \mathcal{L}_2, \ldots , \mathcal{L}_{n-1})$ start at $\bm{u}$, end at $\bm{v}$, and are conditioned to lie above $\mathcal{L}_n$, \eqref{hn2} with height monotonicity (recall \Cref{Gap0}) yields a coupling between $\bm{\mathcal{L}}$ and $\breve{\bm{\mathcal{L}}}$ such that $\mathcal{L}_j = \breve{\mathcal{L}}_j + o(1)$ for each $j \in [1, n-1]$. See the left side of \Cref{f:Airy_Gaps} for a depiction. Thus, we must lower bound the gaps of $\breve{\bm{\mathcal{L}}}$.

To this end, let $\bm{\mathsf{z}} = (\mathsf{z}_1, \mathsf{z}_2, \ldots , \mathsf{z}_{n-1})$ denote non-intersecting Brownian bridges on $[-\mathsf{T}, \mathsf{T}]$, starting and ending at $\bm{0}_{n-1} = (0, 0, \ldots , 0)$, and conditioned to lie above a stretched semicircle $f (s) = 2^{1/2} \sigma \mathsf{T} (\mathsf{T}^2 - s^2)^{1/2}$, for some real number $\sigma = 1 + o(1)$ near but slightly larger than $1$. The gaps of the starting and ending data for $\bm{\mathsf{z}}$ are then smaller than those of $\bm{u}$ and $\bm{v}$, and the lower boundary $h_n$ is more convex than $f$, since $f'' \le -2^{-1/2} \sigma \le -2^{-1/2} - o(1) \le h_n''$. Hence, gap monotonicity applies and implies that the gaps of $\breve{\bm{\mathcal{L}}}$ stochastically dominate those of $\bm{\mathsf{z}}$. We further show that top curves $\big( \mathsf{z}_1, \mathsf{z}_2, \ldots )$ in $\bm{\mathsf{z}}$ converge to an Airy line ensemble (\Cref{f1x2converge}), by again using height monotonicity, now to compare $\bm{\mathsf{z}}$ to the top $n-1$ curves of a Brownian watermelon (with about $A^3 n$ paths, the $n$-th of which is known to concentrate tightly around the semicircle $f$). Therefore, $\big( \mathsf{z}_1 (t) - \mathsf{z}_2 (t), \mathsf{z}_2 (t) - \mathsf{z}_3 (t), \ldots \big)$ converges to $\big( \mathcal{S}_1 (0) - \mathcal{S}_2 (0), \mathcal{S}_2 (0) - \mathcal{S}_3 (0), \ldots \big)$. Combining this with the above comparisons yields the Airy gap lower bound for $\bm{\mathcal{L}}$. See the middle side of \Cref{f:Airy_Gaps} for a depiction.

To prove the upper bound, we instead rely on the global law \eqref{lnt22} (as opposed to the curvature approximation \eqref{hn2}). Again let $n \gg 1$ be a large integer, and now fix a small real number $\varepsilon \in (0, 1)$, independent of $n$. Let $\bm{\mathsf{y}} = (\mathsf{y}_1, \mathsf{y}_2, \ldots , \mathsf{y}_n)$ denote Dyson Brownian motion, starting at time $-\varepsilon n^{1/3}$, with initial data $\bm{u}' = \big( \mathcal{L}_1 (-\varepsilon n^{1/3}), \ldots , \mathcal{L}_n (-\varepsilon n^{1/3}) \big)$. Conditioning on the locations $\bm{v}' = \big( \mathsf{y}_1 (n^{20}), \ldots , \mathsf{y}_n (n^{20}) \big)$ of $\bm{\mathsf{y}}$ at time $n^{20}$, the law of $\bm{\mathsf{y}}$ on $[-\varepsilon n^{1/3}, n^{20}]$ is then given by $n$ non-intersecting Brownian bridges, starting at $\bm{u}'$ and ending at $\bm{v}'$, without boundaries. One can verify (under a few mild modifications to the above setup that we do not detail here) that the gaps $\big| \mathsf{y}_i (n^{20}) - \mathsf{y}_j (n^{20}) \big|$ of $\bm{\mathsf{y}}$ after being run for such a long time $n^{20} + \varepsilon n^{1/3}$ are likely very large, and in particular greater than those $\big| \mathcal{L}_i (n^{20}) - \mathcal{L}_j (n^{20}) \big|$ of $\bm{\mathcal{L}}$ allowed by the gap upper bound \eqref{ltiljji0}. As $\bm{\mathsf{y}}$ has no lower boundary, gap monotonicity thus applies and implies that the gaps of $\bm{\mathcal{L}}$ are stochastically dominated by those of $\bm{\mathsf{y}}$. See the right side of \Cref{f:Airy_Gaps} for a depiction.

Results by Capitaine--P\'{e}ch\'{e} \cite{FESD}, on edge statistics of Dyson Brownian motion under general initial data, can then be used on $\bm{\mathsf{y}}$. They indicate that $\varsigma \cdot \big( \mathsf{y}_1 (t), \mathsf{y}_2 (t), \ldots \big)$ converges (after recentering) to the Airy point process, where the rescaling factor $\varsigma$ admits an explicit formula in terms of the initial data $\bm{u}'$. Using the approximation for $\bm{u}'$ provided by the global law \eqref{lnt22}, we show that $\varsigma \approx 2^{1/2}$ (\Cref{sigmanux12}). It follows that $\big( \mathsf{y}_1 (t) - \mathsf{y}_2 (t), \mathsf{y}_2 (t) - \mathsf{y}_3 (t), \ldots \big)$ converges to $\big( \mathcal{S}_1 (0) - \mathcal{S}_2 (0), \mathcal{S}_2 (0) - \mathcal{S}_3 (0), \ldots \big)$, which with the above comparison between $\bm{\mathcal{L}}$ and $\bm{\mathsf{y}}$ yields the Airy gap upper bound for $\bm{\mathcal{L}}$. See \Cref{ProofDifference} below for further details.

\subsubsection{Airy Line Ensemble}

By \eqref{2rj2}, the fact that the gaps of $\bm{\mathcal{L}}$ at any fixed time $t \in \mathbb{R}$ are given by those of an Airy point process implies strong concentration bounds for $\bm{\mathcal{L}} (t)$, up to an overall (random) shift. In particular, fix large integers $N \gg n \gg 1$ and denote the $N$-tuples $\bm{u} = (u_1, \ldots , u_N) = \big( \mathcal{L}_1 (-n^{1/3}), \ldots , \mathcal{L}_N (-n^{-1/3}) \big)$ and $\bm{v} = (v_1, \ldots , v_N) = \big( \mathcal{L}_1 (n^{1/3}),  \ldots , \mathcal{L}_N (n^{1/3}) \big)$. Then, \eqref{2rj2} yields random variables $u, v \in \mathbb{R}$ (we may take $u = u_N + 2^{-1/2} n^{2/3} + 2^{-7/6} (3\pi)^{2/3} N^{2/3}$ and $v = v_N + 2^{-1/2} n^{2/3} + 2^{-7/6} (3 \pi)^{2/3} N^{2/3}$) such that for $j \in [1, N]$ we have with high probability that
\begin{flalign}
	\label{ujuvjv} 
	\begin{aligned}
	& u_j = u  -2^{-1/2} n^{2/3} -2^{-7/6} (3\pi)^{2/3} j^{2/3} + \mathcal{O} (j^{o(1)-1/3}); \\
	& v_j = v  -2^{-1/2} n^{2/3} -2^{-7/6} (3\pi)^{2/3} j^{-2/3} + \mathcal{O} (j^{o(1)-1/3}).
	\end{aligned}
\end{flalign}

\noindent Condition on $\bm{\mathcal{L}} (-n^{1/3})$ and $\bm{\mathcal{L}} (n^{1/3})$ (thus fixing $u$ and $v$), so $(\mathcal{L}_1, \mathcal{L}_2, \ldots , \mathcal{L}_N)$ are $N$ non-intersecting Brownian bridges starting at $\bm{u}$, ending at $\bm{v}$, and conditioned to lie above $\mathcal{L}_{N+1}$. Further restrict to the (likely) event that \eqref{ujuvjv} holds. By subtracting an affine shift\footnote{This is ultimately what gives the residual independent affine shift in the characterization \Cref{lr} for $\bm{\mathcal{L}}$.} from the $\mathcal{L}_j$ (given in terms of $u$ and $v$ by $t \cdot \xi + \zeta = t \cdot (v - u) / 2n^{1/3}  + (u + v ) / 2$), we may assume $u = v = 0$ in \eqref{ujuvjv}. 

To show that $\bm{\mathcal{L}}$ is a scaled parabolic Airy line ensemble, it then suffices to establish the following more general statement (\Cref{uvxrconverge}). Sample $N$ non-intersecting Brownian bridges $\bm{\mathsf{x}} = (\mathsf{x}_1, \mathsf{x}_2, \ldots , \mathsf{x}_N)$ on $[-n^{1/3}, n^{1/3}]$, starting at $\bm{u}$ and ending at $\bm{v}$ satisfying \eqref{ujuvjv} with $(u, v) = (0, 0)$, and conditioned to lie above a (not too irregular) lower boundary curve $f$. Then $(\mathsf{x}_1, \mathsf{x}_2, \ldots)$ converges to $\bm{\mathcal{S}}$, as $n$ tends to $\infty$. To prove this, we use \eqref{ujuvjv} to sandwich $\bm{\mathsf{x}}$ between two rescaled parabolic Airy line ensembles with approximately equal curvatures. This sort of idea was also fruitful in analyzing edge statistics for random tilings \cite{ESTP}, once a concentration bound for the associated paths around explicit parabolas, as strong as \eqref{rntft}, was proven in the time direction.

	\begin{figure}
	\center
\includegraphics[scale=.6, trim=0cm 1cm 0 0, clip]{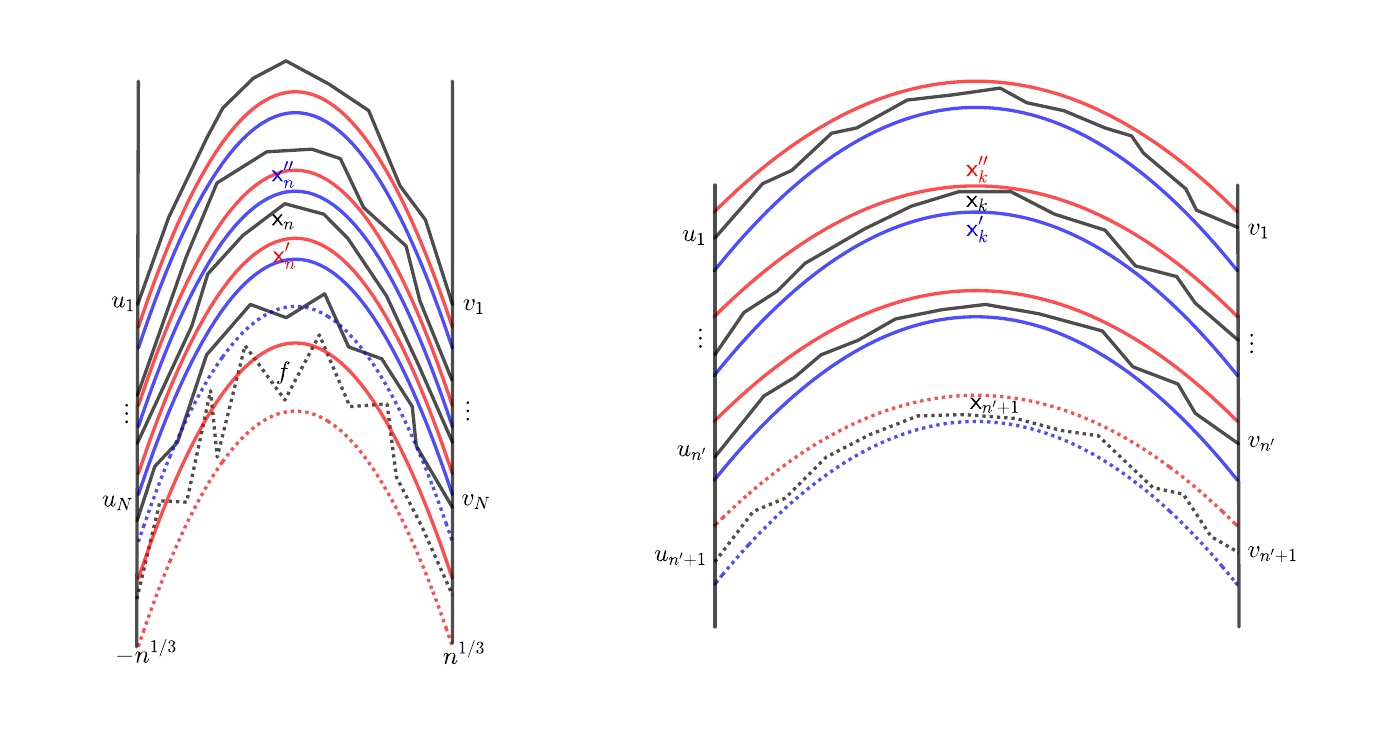}

\caption{Shown to the left is the coupling on a short interval, used to show that the curves in $\bm{\mathsf{x}}$ are close to parabolas (where the top blue curve is at $\infty$ and thus not depicted). Shown to the right is the coupling on a long interval, used to show convergence to the Airy line ensemble.}
\label{f:Airy_Line3}
	\end{figure}

In our context, we instead have the concentration bound \eqref{ujuvjv} in the spatial direction, so we must apply this idea in two ways. The first uses \eqref{ujuvjv} to sandwich $\bm{\mathsf{x}}$ between two rescaled parabolic Airy line ensembles on a tall rectangle; see the left side of \Cref{f:Airy_Line3} for a depiction. 
This enables us to closely approximate $\mathsf{x}_n$ by a parabola in the time direction, verifying that \eqref{rntft} holds for $\bm{\mathsf{x}}$ with that explicit choice of $\mathfrak{f}_n$ (\Cref{uvxrn10}). The second is to use this near-parabolicity to sandwich $\bm{\mathsf{x}}$ between two rescaled parabolic Airy line ensembles on a long interval (as in \cite{ESTP}), to conclude its convergence to the Airy line ensemble. See the right side of \Cref{f:Airy_Line3} for a depiction, and \Cref{ProofL0} below for further details.

\subsection{Guide for the Paper}

Here we provide some guiding commentary on the remainder of this paper. When examining a chapter, readers may wish to consult the information below to help them decide where to focus their attention, depending on their interests. 

\emph{Chapter \ref{BRIDGES0}}: Section \ref{Line} states our main characterization, \Cref{lr}. Section \ref{ProofL} gives the structure of its proof, by precisely stating the results enumerated in \Cref{LProof1} and directly establishing \Cref{lr} assuming them. The proofs of these results constitute the remaining chapters of the paper, which can largely (with a few exceptions) be read independently from each other. Section \ref{Estimates} collects (mainly known) results. Those that we will use most frequently are the scaling invariances (\Cref{linear} and \Cref{scale}) and height monotonicity (\Cref{monotoneheight} and \Cref{uvv}), followed by concentration estimates for Dyson Brownian motion (\Cref{concentrationequation}) and Brownian watermelons (\Cref{estimatexj}). The remaining results in \Cref{Estimates} can be skimmed (or skipped) on an initial reading, and returned to as they are mentioned later in the paper. 

\emph{Chapter \ref{GAPSCALE}}: Section \ref{Coupling1} establishes gap monotonicity; the ideas behind its (sort of nonstandard) proof will play little role after this section, but the statement itself will be used extensively throughout this work. The remainder of this chapter, \Cref{LLocation} and \Cref{ProbabilityScale}, establishes the on-scale estimates. The proofs there are based on the Brownian Gibbs property, so they may feel natural to experts in such arguments.

\emph{Chapter \ref{EDGESHAPE}}: This chapter constitutes an entirely deterministic analysis of limit shapes for non-intersecting Brownian bridges. Section \ref{LimitBridges} collects (mainly known) results about these limit shapes. Section \ref{s:density} proves \emph{a priori} on these limit shapes under boundary data satisfying \eqref{cg2}. Many arguments there are broadly parallel to (and are continuum variants of) ones based on the Brownian Gibbs property, though they often also use as input continuum regularity estimates that do not seem transparent at the discrete level (such as in the proof of \Cref{p:densityrg}). Section \ref{s:shape1} establishes the approximation \eqref{g0abc}; a heuristic is provided in \Cref{Profile0} (which underlies the computations in \Cref{FunctionEdge}), assuming some uniform derivative estimates for the limit shapes. The latter uniform estimates constitute the core of this section and are shown in \Cref{SetCharacteristic}. 

\emph{Chapter \ref{RectangleLCouple}}: Section \ref{s:linelong} provides several concentration results for non-intersecting Brownian bridges, whose proofs mainly consist of combining results from previous chapters; the results there can be skimmed (or skipped) initially, and returned to as they are mentioned later. Section \ref{RectangleCouple} precisely states the boundary removal coupling (described in \Cref{CoupleBoundary2}). The primary effort in its proof is in exhibiting an initial coupling (\Cref{p:comparison2}) and verifying certain H\"{o}lder estimates (\Cref{p:Lip}), both of which are deferred to later sections in this chapter. So, the primary new content in \Cref{RectangleCouple} is its definitions and intermediate statements (which may get a bit involved). Section \ref{Couple0Proof} exhibits the initial coupling; the heuristic underlying its proof is provided in \Cref{Compare00}. Section \ref{RegularImproved} verifies the H\"{o}lder estimates, where the definitions and statements in \Cref{RegularImprovedProof} play a central role (see also \Cref{proofd00} for an explanation for one of the exponents appearing there). 

\emph{Chapter \ref{GlobalRegular}}: Section \ref{ProofRegular} and \Cref{ProofGlobal} show spatial regularity and the global law for $\bm{\mathcal{L}}$, respectively. The proofs in this chapter consist of direct (but slightly tedious) combinations of results already established in the previous chapters.

\emph{Chapter \ref{APPROXIMATECURVE}}: Section \ref{DerivativePath2} establishes the curvature approximation. The core of this section constitutes the inductive definitions and statements in \Cref{LOCALHJ}. They are shown in the following \Cref{Proof0Omega}, whose main aspect is the proof of \Cref{omega4}. Numerous error parameters and rescalings appear in this chapter, and the reader may wish to keep \Cref{delta0omegaklk} and \Cref{x4} in mind while navigating them.

\emph{Chapter \ref{STATISTICSBRIDGES}}: This chapter concludes Airy statistics for $\bm{\mathcal{L}}$, where the main novelty is in showing that its gaps coincide with those of the Airy line ensemble; this is done in \Cref{ProofDifference}. Section \ref{ProofL0}, which is closer to prior work (such as \cite[Section 3.5]{ESTP}), uses this to confirm Airy line ensemble statistics. As in \Cref{LLocation} and \Cref{ProbabilityScale}, the proofs in this chapter are based on the Brownian Gibbs property, and so they may again feel natural to experts in such arguments. 

\subsection{Notation}

\label{FNotation} 

In what follows, we set $\overline{\mathbb{R}} = \mathbb{R} \cup \{ -\infty, \infty \}$;\index{R@$\overline{\mathbb{R}}$} its topology is specified by prescribing the standard one on $\mathbb{R}$ and imposing the $\{ -\infty \}$ and $\{ \infty \}$ are each both open and closed sets. We also let $\mathbb{H} = \{ z \in \mathbb{C} : \Imaginary z > 0\}$ denote the upper half-plane, $\overline{\mathbb{H}}$ denote its closure, $\mathbb{H}^- = \{ z \in \mathbb{C} : \Imaginary z < 0\}$ \index{H@$\mathbb{H}^-$}denote the lower half-plane; and $\overline{\mathbb{H}^-}$ denote its closure. Moreover, for any subset $I \subseteq \mathbb{R}$ and measurable functions $f, g : I \rightarrow \overline{\mathbb{R}}$ we write $f < g$ (equivalently, $g > f$) if $f(t) < g(t)$ for each $t \in I$; we similarly write $f \le g$ (equivalently, $g \ge f$) if $f(t) \le g(t)$ for each $t \in I$. For any sets $A_0 \subseteq A$ and function $f : A \rightarrow \mathbb{C}$, let $f |_{A_0}$ denote the restriction of $f$ to $A_0$. In what follows, for any topological space $I$, we let $\mathcal{C}(I)$ \index{C@$\mathcal{C}(I)$}denote the space of real-valued, continuous functions $f : I \rightarrow \mathbb{R}$. Given an integer $d \ge 1$ and a subset $U \subset \mathbb{R}^d$, a function $f : U \rightarrow \mathbb{C}$ is called real analytic if, for every point $z_0$ in the interior of $U$, it admits a power series expansion that converges absolutely in a neighborhood of $z_0$.{\index{R@real analytic}} 

For any integer $d \ge 1$ and $d$-tuple $\gamma = (\gamma_1, \gamma_2, \ldots , \gamma_d) \in \mathbb{Z}_{\ge 0}^d$, define $|\gamma| = \sum_{j=1}^d \gamma_j$ and $\partial_{\gamma} = \prod_{j=1}^d \partial_j^{\gamma_j}$, where we have abbreviated the differential operator $\partial_j = \partial_{x_j} = \partial / \partial x_j $ for each $j \in [1, d]$. For any integer $k \ge 0$ and open subset $\mathfrak{R} \subseteq \mathbb{R}^d$, let $\mathcal{C}^k (\mathfrak{R})$\index{C@$\mathcal{C}^k (\mathfrak{R})$, $\mathcal{C}^k (\overline{\mathfrak{R}})$} denote the set of $f \in \mathcal{C} (\overline{\mathfrak{R}})$ such that $\partial_{\gamma} f$ is continuous on $\mathfrak{R}$, for each $\gamma \in \mathbb{Z}_{\ge 0}^d$ with $|\gamma| \le k$. Further let $\mathcal{C}^k (\overline{\mathfrak{R}})$ denote the set of functions $f \in \mathcal{C}^k (\mathfrak{R})$ such that $\partial_{\gamma} f$ extends continuously to $\overline{\mathfrak{R}}$, for each $\gamma \in \mathbb{Z}_{\ge 0}^d$ with $|\gamma| \le k$. For any function $f \in \mathcal{C} (\mathfrak{R})$ and integer $k \in \mathbb{Z}_{\ge 0}$, we further define the (semi)norms $\| f \|_0 = \| f \|_{0; \mathfrak{R}}$,\index{F@$\lVert f \rVert_0$} $[f] = [f]_{k; 0; \mathfrak{R}}$,\index{F@$[f]_k$} and $\| f \|_{\mathcal{C}^k (\overline{\mathfrak{R}})} = \| f \|_{\mathcal{C}^k (\mathfrak{R})} = \| f \|_k = \| f \|_{k; 0; \mathfrak{R}}$\index{F@$\lVert f \rVert_{\mathcal{C}^k (\mathfrak{R})}$} on these spaces by setting
\begin{flalign}\label{e:norms}
	\| f \|_0 = \displaystyle\sup_{z \in \mathfrak{R}} \big| f(z) \big|; \qquad [f]_k = \displaystyle\max_{\substack{\gamma \in \mathbb{Z}_{\ge 0}^d \\ |\gamma| = k}} \big\| \partial_{\gamma} f \|_0; \qquad \| f \|_{\mathcal{C}^k (\overline{\mathfrak{R}})} = \displaystyle\sum_{j=0}^k [f]_j.
\end{flalign}

\noindent The norms $\| f \|_0$ and $\| f \|_{\mathcal{C}^k (\overline{\mathfrak{R}})}$ extend by continuity to the closures of the spaces $\mathcal{C} (\overline{\mathfrak{R}})$ and $\mathcal{C}^k (\overline{\mathfrak{R}})$ (with respect to these norms), respectively.

For any real numbers $a, b \in \mathbb{R}$ with $a \le b$, we set $\llbracket a, b \rrbracket = [a, b] \cap \mathbb{Z}$.\index{A@$\llbracket a, b \rrbracket$} We also let $a \vee b = \max \{ a, b \}$\index{A@$a \vee b$} and $a \wedge b = \min \{ a, b \}$\index{A@$a \wedge b$}. For any integer $k \ge 1$, we denote the entries of any $k$-tuple $\bm{y} \in \mathbb{C}^k$ by $\bm{y} = (y_1 ,y_2, \ldots , y_k)$, unless we specify the indexing otherwise. For any $k$-tuples $\bm{x}, \bm{y} \in \mathbb{R}^k$, we write $\bm{x} < \bm{y}$ (equivalently, $\bm{y} > \bm{x}$) if $x_j < y_j$ for each $j \in \llbracket 1, k \rrbracket$; we similarly write $\bm{x} \le \bm{y}$ (equivalently, $\bm{y} \ge \bm{x}$) if $x_j \le y_j$ for each $j \in \llbracket 1, k \rrbracket$. We also let $\mathbb{W}_k = \{ \bm{y} \in \mathbb{R}^k : y_1 > y_2 > \cdots >  y_k \}$\index{W@$\mathbb{W}_k$} and let $\overline{\mathbb{W}}_k$ denote the closure of $\mathbb{W}_k$. Further let $\bm{0}_k = (0, 0, \ldots , 0) \in \overline{\mathbb{W}}_k$, \index{0@$\bm{0}_k$}where $0$ appears with multiplicity $k$.

For any integer $k \ge 1$ and subset $\mathfrak{S} \subseteq \mathbb{R}^k$, we let $\partial \mathfrak{S}$ denote the boundary of $\mathfrak{S}$; for any point $z \in \mathbb{R}^k$, we also let $\dist (z, \mathfrak{S}) = \inf_{s \in \mathfrak{S}} |z-s|$.\index{D@$\dist$} For any complex numbers $a, b \in \mathbb{C}$, and vector $\bm{x} \in \mathbb{C}^k$, we set $a \bm{x} + b = (ax_1 + b, ax_2 + b, \ldots , ax_k + b) \in \mathbb{C}^k$. For any interval $I \subset \mathbb{R}^k$ and set $\mathfrak{S}$ of vectors $\mathfrak{S} \subset \mathbb{R}^k$ or of functions $\mathfrak{S} \in \mathcal{C} (I)$, we similarly denote $a \cdot \mathfrak{S} + b = \{ as + b \}_{s \in \mathfrak{S}}$. For any additional such set $\mathfrak{S}'$, denote $\mathfrak{S} + \mathfrak{S}' = \{ s + s' : s \in \mathfrak{S}, s' \in \mathfrak{S}' \}$. 

Let $\mathscr{P}_{\fin} = \mathscr{P}_{\fin} (\mathbb{R})$\index{P@$\mathscr{P}_{\fin}$, $\mathscr{P}$, $\mathscr{P}_0$} denote the set of nonnegative measures $\mu$ on $\mathbb{R}$ with finite total mass, $\mu (\mathbb{R}) < \infty$. Further let $\mathscr{P} = \mathscr{P} (\mathbb{R}) \subset \mathscr{P}_{\fin}$ denote the set of probability measures on $\mathbb{R}$, and let $\mathscr{P}_0 = \mathscr{P}_0 (\mathbb{R}) \subset \mathscr{P}$ denote the set of probability measures that are compactly supported; the support of any measure $\nu \in \mathscr{P}$ is denoted by $\supp \nu$. We say that a probability measure $\mu \in \mathscr{P}$ has density $\varrho$ (with respect to Lebesgue measure) if $\varrho : \mathbb{R} \rightarrow \mathbb{R}$ is a measurable function satisfying $\mu (dx) = \varrho(x) dx$. For any real number $x \in \mathbb{R}$, we let $\delta_x \in \mathscr{P}_0$ denote the delta function at $x$.\index{0@$\delta_x$} For any sequence $\bm{a} = (a_1, a_2, \ldots , a_n) \in \overline{\mathbb{W}}_n$, we denote its \emph{empirical measure} $\emp (\bm{a}) \in \mathscr{P}$ by\index{E@$\emp$}
\begin{flalign}
	\label{aemp} 
	\emp (\bm{a}) = \displaystyle\frac{1}{n} \displaystyle\sum_{j=1}^n \delta_{a_j}.
\end{flalign}

\noindent We denote the complement of any event $\mathscr{E}$ by $\mathscr{E}^{\complement}$.\index{E@$\mathscr{E}^{\complement}$}

Throughout, given some integer parameter $n \ge 1$ and event $\mathscr{E}_n$ depending on $n$, we will often make statements of the following form. There exists a constant $c > 0$, independent of $n$ (but possibly dependent on other parameters), such that $\mathbb{P}[\mathscr{E}_n^{\complement}] \le f(c, n)$ holds for an explicit function $f : \mathbb{R}_{> 0} \times \mathbb{Z}_{\ge 0} \rightarrow \mathbb{R}_{\ge 0}$, which is non-decreasing in $c$ and satisfies $\lim_{n \rightarrow \infty} f (c, n) = 0$ and $\lim_{c \rightarrow 0} f (c, n) > 1$ (an example is $f(c, n) = c^{-1} e^{-c (\log n)^2}$). When proving such statements we will often implicitly (and without comment) assume that $n \ge N_0$ is sufficiently large. Indeed, suppose there exist $N_0 \ge 1$ and $c_0 > 0$ such that $\mathbb{P} \big[ \mathscr{E}_n^{\complement} \big] \le f(c_0, n)$ holds whenever $n \ge N_0$. Since $\lim_{c \rightarrow 0} f(c, n) > 1$, there exists a constant $c_1 > 0$ such that for $n \le N_0$ we have $f(c_1, n) \ge 1$, in which case $\mathbb{P} \big[\mathscr{E}_n^{\complement} \big] \le 1 \le f(c_1, n)$ continues to hold. Thus, taking $c = \min \{ c_0, c_1 \}$ guarantees that $\mathbb{P} \big[ \mathscr{E}_n^{\complement} \big] \le f(c, n)$ holds for all $n \ge 1$.

\subsection*{Acknowledgements} 

The research of Amol Aggarwal was partially supported by a Packard Fellowship for Science and Engineering, a Clay Research Fellowship, by NSF grant DMS-1926686, and by the IAS School of Mathematics. The research of Jiaoyang Huang was supported by NSF grant DMS-2054835 and DMS-2331096, and the Sloan Research Fellowship.  The authors heartily thank Alan Hammond for extensive detailed and valuable comments throughout this paper. They further thank Alexei Borodin and Ivan Corwin for providing very helpful suggestions, and also Evgeni Dimitrov, Milind Hedge, Andrei Okounkov, Fabio Toninelli, Horng-Tzer Yau, and Lingfu Zhang for insightful conversations. They are also grateful to the anonymous referees, for providing detailed comments that improved the presentation of this paper. They wish to acknowledge the NSF grant DMS-1928930, which supported their participation in the Fall 2021 semester program at MSRI in Berkeley, California titled, “Universality and Integrability in Random Matrix Theory and Interacting Particle Systems.''

	\section{Results}
	
	\label{Line} 
	
	\subsection{Brownian Gibbs Property and the Airy Line Ensemble}
	
	\label{EnsemblesCurve} 
	
	In this section we introduce notation for non-intersecting Brownian bridges. Let $\Sigma \subseteq \mathbb{Z}_{\ge 1}$ and $I \subseteq \mathbb{R}$ denote intervals. Let $\mathfrak{X} = \mathfrak{X} (\Sigma; I)$ denote the set of continuous functions $f : \Sigma \times I \rightarrow \mathbb{R}$, whose topology is determined by uniform convergence on compact subsets of $\Sigma \times I$; we denote the associated Borel $\sigma$-algebra by $\mathscr{C} = \mathscr{C} (\Sigma \times I)$. Since $\Sigma$ is discrete, any such $f$ can be interpreted as an element of $\Sigma \times \mathcal{C}(I)$.\index{S@$\Sigma \times \mathcal{C}(I)$} An \emph{$\Sigma \times I$ indexed line ensemble} is a $\mathfrak{X}$-valued random variable $\bm{\mathsf{x}} \in \Sigma \times \mathcal{C} (I)$ defined on a probability space $(\Omega, \mathscr{B}, \mathbb{P})$ that is $(\mathscr{B}, \mathscr{C})$-measurable. We will frequently denote such a line ensemble by $\bm{\mathsf{x}} = (\mathsf{x}_j)_{j \in \Sigma}$, where $\mathsf{x}_j : I \rightarrow \mathbb{R}$ is a (random) continuous function for each $j \in \Sigma$; in this case, we also set $\bm{\mathsf{x}} (t) = \big( \mathsf{x}_j (t) \big)_{j \in \Sigma}$ for each $t \in I$. 
	
	We next provide notation for the probability measure of $n$ non-intersecting Brownian bridges with given starting and ending points, and for given upper and lower boundaries.
	
	\begin{figure}
	\center
\includegraphics[scale=.8]{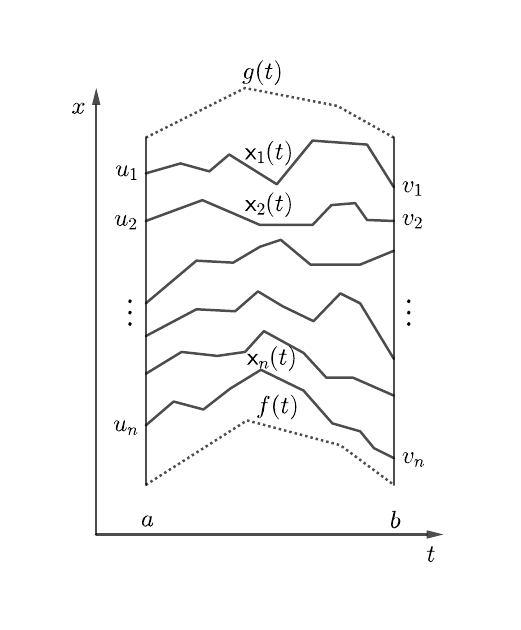}

\caption{Depicted above is a sample from $\mathsf{Q}_{f; g}^{\bm{u}; \bm{v}} (\sigma)$.}
\label{f:bridge1}
	\end{figure}
	
	\begin{definition}
		
		\label{qxyfg}
		
		Fix an integer $n \ge 1$; a real number $\sigma > 0$; two $n$-tuples $\bm{u}, \bm{v} \in \mathbb{W}_n$; an interval $[a, b] \subseteq \mathbb{R}$; and continuous functions $f, g : [a, b] \rightarrow \overline{\mathbb{R}}$ with $f < g$, such that $f(a) < u_n < u_1 < g(a)$ and $f(b) < v_n < v_1 < g(b)$. Let $\mathsf{Q}_{f; g}^{\bm{u}; \bm{v}} (\sigma)$ denote the law of the $\llbracket 1, n \rrbracket \times [a, b]$ indexed line ensemble $\bm{\mathsf{x}} = (\mathsf{x}_1, \mathsf{x}_2, \ldots , \mathsf{x}_n) \in \llbracket 1, n \rrbracket \times \mathcal{C} \big( [a, b] \big)$, given by $n$ independent Brownian motions of variance $\sigma$ on the time interval $t \in [a, b]$, conditioned on satisfying the following three properties.
		\begin{enumerate} 
			\item The $\mathsf{x}_j$ do not intersect, that is, $\bm{\mathsf{x}} (t) \in \mathbb{W}_n$ for each $t \in (a, b)$.
			\item The $\mathsf{x}_j$ start at $u_j$ and end at $v_j$, that is, $\mathsf{x}_j (a) = u_j$ and $\mathsf{x}_j (b) = v_j$ for each $j \in \llbracket 1, n \rrbracket$.
			\item The $\mathsf{x}_j$ are bounded below by $f$ and above by $g$, that is, $f < \mathsf{x}_j < g$ for each $j \in \llbracket 1, n \rrbracket$. 
		\end{enumerate} 
				
		\noindent The above definition can be extended to allow for $\bm{u}, \bm{v} \in \overline{\mathbb{W}}_n$ through a limiting procedure. Specifically, suppose that either $\bm{u} \in \overline{\mathbb{W}}_n \setminus \mathbb{W}_n$ or $\bm{v} \in \overline{\mathbb{W}}_n \setminus \mathbb{W}_n$. Then, for any $\varepsilon > 0$, define $\bm{u}^{\varepsilon}, \bm{v}^{\varepsilon} \in \mathbb{W}_n$ by setting $u_j^{\varepsilon} = u_j - j \varepsilon$ and $v_j^{\varepsilon} = v_j - j \varepsilon$, for each $j \in \llbracket 1, n \rrbracket$. Also define $f^{\varepsilon}, g^{\varepsilon} : [a, b] \rightarrow \overline{\mathbb{R}}$ by setting $f^{\varepsilon} (s) = f(s) - (n+1) \varepsilon$ and $g^{\varepsilon} (s) = g(s)$, for each $s \in [a, b]$. It is directly verified\footnote{See \Cref{boundaryconvergeensemble} below.} that, as $\varepsilon$ tends to $0$, the measure $\mathsf{Q}_{f^{\varepsilon}; g^{\varepsilon}}^{\bm{u}^{\varepsilon}; \bm{v}^{\varepsilon}}$ converges to a limit, which we denote by $\mathsf{Q}_{f;g}^{\bm{u}; \bm{v}}$. 
			
		We refer to $\bm{u}$ as \emph{starting data} for $\bm{\mathsf{x}}$, and to $\bm{v}$ as its \emph{ending data}. We also refer to $f$ as the \emph{lower boundary} for $\bm{\mathsf{x}}$, and to $g$ as its \emph{upper boundary}.  When discussing $\mathsf{Q}_{f;g}^{\bm{u};\bm{v}} (\sigma)$, it is always assumed that $f(a) \le u_n \le u_1 \le g(a)$ and $f(b) \le v_n \le v_1 \le g(b)$, and that $f < g$, even when not stated explicitly. See \Cref{f:bridge1}. 
		
		If $g = \infty$, then we abbreviate $\mathsf{Q}_{f; g}^{\bm{u}; \bm{v}} (\sigma) = \mathsf{Q}_f^{\bm{u}; \bm{v}} (\sigma)$;\index{Q@$\mathsf{Q}_{f;g}^{\bm{u}; \bm{v}}$} if additionally $f = -\infty$, we set $\mathsf{Q}_{f; \infty}^{\bm{u}; \bm{v}} (\sigma) = \mathsf{Q}_{f}^{\bm{u}; \bm{v}} (\sigma) = \mathsf{Q}^{\bm{u}; \bm{v}} (\sigma)$. If $\sigma = 1$, then we omit the parameter $\sigma$ from the notation, writing $\mathsf{Q}_{f; g}^{\bm{u}; \bm{v}} = \mathsf{Q}_{f; g}^{\bm{u}; \bm{v}} (1)$, $\mathsf{Q}_f^{\bm{u}; \bm{v}} = \mathsf{Q}_f^{\bm{u}; \bm{v}} (1)$, and $\mathsf{Q}^{\bm{u}; \bm{v}} = \mathsf{Q}^{\bm{u}; \bm{v}} (1)$.  
		
	\end{definition}

	We next describe a resampling property from \cite{PLE}. 
	
	\begin{definition} 
		
	\label{property} 
	
	Fix intervals $\Sigma \subseteq \mathbb{Z}_{\ge 1}$ and $I \subseteq \mathbb{R}$, as well as a $\Sigma \times I$ indexed line ensemble $\bm{\mathsf{x}} = (\mathsf{x}_j)_{j \in \Sigma}$. For any integers $1 \le j \le k$ such that $\llbracket j, k \rrbracket \subseteq \Sigma$, define the $\llbracket j, k \rrbracket \times I$ indexed line ensemble $\bm{\mathsf{x}}_{\llbracket j, k \rrbracket} = (\mathsf{x}_j, \mathsf{x}_{j+1}, \ldots , \mathsf{x}_k) \in \llbracket j, k \rrbracket \times \mathcal{C}(I)$. For any intervals $\Sigma' \subseteq \Sigma$ and $I' \subseteq I$, further let $\mathcal{F}_{\ext} ( \Sigma' \times I')$\index{F@$\mathcal{F}_{\ext} (\Sigma' \times I')$} denote the $\sigma$-algebra generated by the $\big( \mathsf{x}_j (t) \big)$, over all $j \notin \Sigma'$ or $t \notin I'$. 
	
	Fix a real number $\sigma > 0$. We say that a $\Sigma \times I$ indexed line ensemble $\bm{\mathsf{x}}$ has the \emph{Brownian Gibbs property} of variance $\sigma$ if we almost surely have $\mathsf{x}_1 (t) > \mathsf{x}_2 (t) > \cdots $, for each real number $t \in I$, and the following holds, for any bounded intervals $\llbracket k_1, k_2 \rrbracket \subseteq \Sigma$ and $(a, b) \subseteq I$. The law of $\big( \mathsf{x}_j (t) \big)$, over $(j, t) \in \llbracket k_1, k_2 \rrbracket \times [a, b]$, conditional on $\mathcal{F}_{\ext} \big( \llbracket k_1, k_2 \rrbracket \times (a, b) \big)$ is given by the non-intersecting Brownian bridge measure $\mathsf{Q}_{f, g}^{\bm{u}; \bm{v}} (\sigma)$. Here, the entrance and exit data $\bm{u}, \bm{v} \in \mathbb{R}^{k_2 - k_1 + 1}$ are given by $\bm{u} = \big( \mathsf{x}_{k_1} (a), \mathsf{x}_{k_1 + 1} (a), \ldots , \mathsf{x}_{k_2} (a) \big)$ and $\bm{v} = \big( \mathsf{x}_{k_1} (b), \mathsf{x}_{k_1+1} (b), \ldots , \mathsf{x}_{k_2} (b) \big)$, and the boundary data $f, g : [a, b] \rightarrow \mathbb{R}$ are given by $f = \mathsf{x}_{k_1-1} |_{[a, b]}$, and $g = \mathsf{x}_{k_2+1} |_{[a, b]}$ (setting $\mathsf{x}_j = \infty$ if $j < \min \Sigma$ and $\mathsf{x}_j = -\infty$ if $j > \max \Sigma$). If $\sigma = 1$, we omit it from the notation, saying that $\bm{\mathsf{x}}$ satisfies the Brownian Gibbs property. 
	
	\end{definition}

	We next require some notation on edge statistics. 
	
	\begin{definition} 
		
		\label{kernellimit} 
		
		For any $s, t, x, y \in \mathbb{R}$, the \emph{extended Airy kernel} $\mathcal{K} : \mathbb{R}^4 \rightarrow \mathbb{R}$ is given by
		\begin{flalign*}
			\mathcal{K} (s, x; t, y) = \displaystyle\int_0^{\infty} e^{u (t - s)} \Ai (x + u) \Ai (y + u) \mathrm{d}u, \qquad & \text{if $s \ge t$}; \\
			\mathcal{K} (s, x; t, y) = - \displaystyle\int_{-\infty}^0 e^{u (t - s)} \Ai (x + u) \Ai (y + u) \mathrm{d}u, \qquad & \text{if $s < t$},
		\end{flalign*}
		
		\noindent where we recall that the Airy function $\Ai: \mathbb{R} \rightarrow \mathbb{R}$ is given by
		\begin{flalign*}
			\Ai (x) = \displaystyle\frac{1}{\pi} \displaystyle\int_{0}^{\infty} \cos \Big( \displaystyle\frac{z^3}{3} + xz \Big) \mathrm{d}z.
		\end{flalign*}
		
	\end{definition} 
	
	From this, we define the Airy line ensemble. 
	
	\begin{definition}
		
		\label{ensemblewalks}
		
		The \emph{(stationary) Airy line ensemble} $\bm{\mathcal{A}} = (\mathcal{A}_1, \mathcal{A}_2, \ldots ) \in \mathbb{Z}_{\ge 1} \times \mathcal{C} (\mathbb{R})$\index{A@$\bm{\mathcal{A}}$; Airy line ensemble} is an infinite collection of random continuous curves $\mathcal{A}_i: \mathbb{R} \rightarrow \mathbb{R}$, ordered as $\mathcal{A}_1 (t) > \mathcal{A}_2 (t) > \cdots$ for each $t \in \mathbb{R}$, such that 
		\begin{flalign}
			\label{probabilityaxjtj}
		 		\mathrm{d} \mathbb{P} \Bigg[ \bigcap_{j = 1}^m \big\{ (t_j, y_j) \in \bm{\mathcal{A}} \big\} \Bigg] = \det \big[ \mathcal{K} (t_i, y_i; t_j, y_j) \big]_{1 \le i, j \le m} \displaystyle\prod_{j = 1}^m d y_j,
		\end{flalign}
		
		\noindent for any $(t_1, y_1), (t_2, y_2), \ldots , (t_m, y_m) \in \mathbb{R}^2$. Here, we have written $(t, y) \in \bm{\mathcal{A}}$ if there exists some integer $k \ge 1$ such that $\mathcal{A}_k (t) = y$. The existence of such an ensemble was shown as \cite[Theorem 3.1]{PLE} (and the uniqueness follows from the prescription \eqref{probabilityaxjtj} of its multi-point distributions; see \cite[Lemma 3.1]{CLE}). 
		
		We abbreviate the \emph{parabolic Airy line ensemble} $\bm{\mathcal{R}} = \big( \mathcal{A}_1 (t) - t^2, \mathcal{A}_2 (t) - t^2, \ldots \big) \in \mathbb{Z}_{\ge 1} \times \mathcal{C}(\mathbb{R})$, which may be viewed as a function $\mathcal{R}: \mathbb{Z}_{\ge 1} \times \mathbb{R} \rightarrow \mathbb{R}$ by setting $\mathcal{R} (i, t) = \mathcal{R}_i (t) = \mathcal{A}_i (t) - t^2$.\index{R@$\bm{\mathcal{R}}$; parabolic Airy line ensemble} Further define the \emph{scaled parabolic Airy line ensemble} $\bm{\mathcal{S}} = ( \mathcal{S}_1, \mathcal{S}_2, \ldots ) \in \mathbb{Z}_{\ge 1} \times \mathcal{C}(\mathbb{R})$ by setting $\bm{\mathcal{S}} = 2^{-1/2} \cdot \bm{\mathcal{R}}$, that is, $\mathcal{S}_i (t) = 2^{-1/2} \cdot \mathcal{R}_i (t)$ for each $(i, t) \in \mathbb{Z}_{\ge 1} \times \mathbb{R}$.\index{S@$\bm{\mathcal{S}}$ scaled parabolic Airy line ensemble}
		
	\end{definition}

	The following lemma from \cite{PLE} states that the scaled parabolic Airy line ensemble satisfies the Brownian Gibbs property (after rescaling by $2^{-1/2}$). 
	
	\begin{lem}[{\cite[Theorem 3.1]{PLE}}]
		
		\label{propertya} 
		
		The ensemble $\bm{\mathcal{S}}$ has the Brownian Gibbs property.
		
	\end{lem}

	For any real number $s \in \mathbb{R}$, we let $\Theta_s : \mathcal{C} (\mathbb{R}) \rightarrow \mathcal{C}(\mathbb{R})$ denote the translation operator acting on any function $f \in \mathcal{C} (\mathbb{R})$ by setting $\Theta_s f (x) = f (x+s)$, for each $x \in \mathbb{R}$. This operator also acts on $\mathbb{Z}_{\ge 1} \times \mathbb{R}$ indexed line ensembles $\bm{\mathsf{x}} = (\mathsf{x}_1, \mathsf{x}_2, \ldots )$ by setting $\Theta_s \bm{\mathsf{x}} = (\Theta_s \mathsf{x}_1, \Theta_s \mathsf{x}_2, \ldots )$. As such, it also acts on measurable sets in the Borel $\sigma$-algebra $\mathscr{C} = \mathscr{C} (\mathbb{Z}_{\ge 1} \times \mathbb{R})$. An $\mathbb{Z}_{\ge 1} \times \mathbb{R}$ indexed line ensemble $\bm{\mathsf{x}}$ is called \emph{translation-invariant} if the law of $\bm{\mathsf{x}}$ is equal to that of $\Theta_s \bm{\mathsf{x}}$, for each $s \in \mathbb{R}$. 
	
	We further say that an event $\mathscr{F}$ is translation-invariant if $\Theta_s \mathscr{F} = \mathscr{F}$, for any $s \in \mathbb{R}$. For any real number $\sigma > 0$, we let $\Tra (\sigma)$ denote the set of probability measures $\mu$ associated with a $\mathbb{Z}_{\ge 1} \times \mathbb{R}$ indexed line ensemble $\bm{\mathsf{x}} = (\mathsf{x}_1, \mathsf{x}_2, \ldots )$ satisfying the Brownian Gibbs property, such that the ensemble $\big( \mathsf{x}_1 (t) + \sigma t^2, \mathsf{x}_2 (t) + \sigma t^2, \ldots \big) \in \mathbb{Z}_{\ge 1} \times \mathcal{C} (\mathbb{R})$ is translation-invariant. We call a measure $\mu \in \Tra (\sigma)$ \emph{extremal} if, for any real number $p \in (0, 1)$ and measures $\mu_1, \mu_2 \in \Tra (\sigma)$ such that $\mu = p \mu_1 + (1 - p) \mu_2$, we have $\mu_1 = \mu = \mu_2$.

	\begin{lem}[{\cite[Proposition 1.13]{ELE}}]
		
		\label{translationa} 
		
		The law of $\bm{\mathcal{S}}$ is in $\Tra (2^{-1/2})$ and is extremal. 
		
	\end{lem}

	\subsection{Line Ensembles With Parabolic Decay} 
	
	\label{Ensemble2} 
	
	In this section we state our results, which constitute characterizations for line ensembles satisfying the Brownian Gibbs property and certain growth conditions. The latter conditions are explained through the following definition and assumption, which describe the family of line ensembles we will analyze. The definition introduces the event on which a given point in the top curve of the ensemble is between two parabolas of approximately equal leading coefficients (chosen to be $-2^{-1/2}$, to agree with the behavior of $\bm{\mathcal{S}}$); the assumption states that this event is likely.

	\begin{definition}
		
		\label{eventparabola} 
		
		Let $\Sigma \subseteq \mathbb{Z}_{\ge 1}$ denote an interval with $1 \in \Sigma$; let $I \subseteq \mathbb{R}$ denote an interval (not necessarily bounded); and let $\bm{\mathsf{x}} = (\mathsf{x}_s)_{s \in \Sigma} \in \Sigma \times \mathcal{C} (I)$ denote a $\Sigma \times I$ indexed line ensemble. For any real numbers $\varepsilon > 0$, $C > 1$ and $t \in I$, define the event $\textbf{PAR}_{\varepsilon} (t; C) = \textbf{PAR}_{\varepsilon}^{\bm{\mathsf{x}}} (t; C)$ by\index{P@$\textbf{PAR}$}
		\begin{flalign}
			\label{eventtec} 
			\textbf{PAR}_{\varepsilon} (t; C) = \big\{ -(2^{-1/2} + \varepsilon) t^2 - C \le \mathsf{x}_1 (t) \le -(2^{-1/2} - \varepsilon) t^2 + C \big\}.  
		\end{flalign}
	\end{definition}

	\begin{assumption}
	
		\label{l0} 
		
		Let $\bm{\mathcal{L}} = (\mathcal{L}_1, \mathcal{L}_2, \ldots ) \in \mathbb{Z}_{\ge 1} \times \mathcal{C} (\mathbb{R})$ denote a $\mathbb{Z}_{\ge 1} \times \mathbb{R}$ indexed line ensemble satisfying the Brownian Gibbs property. Assume that there exists a function\footnote{Whenever adopting this assumption, we will view $\mathfrak{K}$ as fixed. In particular, underlying constants might depend on $\mathfrak{K}$, even when this dependence is not stated explicitly.} $\mathfrak{K} : \mathbb{R}_{> 0} \rightarrow \mathbb{R}_{>0}$ such that, for each $\varepsilon> 0$, we have $\mathbb{P} \big[ \textbf{PAR}_{\varepsilon}^{\bm{\mathcal{L}}} (t; \mathfrak{C}) \big] \ge 1 - \varepsilon$, for any real numbers $t \in \mathbb{R}$ and $\mathfrak{C} \ge \mathfrak{K} (\varepsilon)$.
		
	\end{assumption}

	The following assumption classifies those line ensembles satisfying \eqref{l0} as a combination of scaled parabolic Airy line ensembles; it will be established in \Cref{Prooflr} below.
	
	\begin{thr}
		
		\label{lr} 
		
		Adopt \Cref{l0}. There exist two random variables $\mathfrak{l}, \mathfrak{c} \in \mathbb{R}$, and a scaled parabolic Airy line ensemble $\bm{\mathcal{S}} = (\mathcal{S}_1, \mathcal{S}_2, \ldots ) \in \mathbb{Z}_{\ge 1} \times \mathcal{C} (\mathbb{R})$ (as in \Cref{ensemblewalks}) independent from them, such that $\mathcal{L}_j (t) = \mathcal{S}_j (t) + \mathfrak{l} t + \mathfrak{c}$, for each $(j, t) \in \mathbb{Z}_{\ge 1} \times \mathbb{R}$. 
		
	\end{thr}
	
	Let us discuss a few consequences of \Cref{lr}. To do so, it will be useful to rescale the parabolic Airy line ensemble. For any real number $\sigma > 0$, define the $\mathbb{Z}_{\ge 1} \times \mathbb{R}$ indexed line ensemble $\bm{\mathcal{S}}^{(\sigma)} = \big( \mathcal{S}_1^{(\sigma)}, \mathcal{S}_2^{(\sigma)}, \ldots \big) \in \mathbb{Z}_{\ge 1} \times \mathcal{C}(\mathbb{R})$ by setting\index{S@$\bm{\mathcal{S}}^{(\sigma)}$; rescaled Airy line ensemble}
	\begin{flalign}
		\label{sigmar}
		\mathcal{S}_j^{(\sigma)} (t) = \sigma^{-1} \cdot \mathcal{S}_j (\sigma^2 t), \quad \text{for each $(j, t) \in \mathbb{Z}_{\ge 1} \times \mathbb{R}$}.
	\end{flalign}
	
	\begin{rem}
		
		\label{sigmascale} 
		
		Since, for any real number $\sigma > 0$, the law of any Brownian bridge $B(t)$ is equal to that of $\sigma^{-1} \cdot B (\sigma^2 t)$ (see \Cref{scale} below), and since $\bm{\mathcal{S}}$ satisfies the Brownian Gibbs property (recall \Cref{propertya}), $\bm{\mathcal{S}}^{(\sigma)}$ does as well for any $\sigma > 0$.
		
	\end{rem}

	From \Cref{lr}, we can quickly derive the following corollary classifying line ensembles with more specific rates of decay; it will be established in \Cref{Proofsql1} below.

	\begin{cor}
		
		\label{constantk1} 
		
		Fix a real number $\sigma > 0$ and a $\mathbb{Z}_{\ge 1} \times \mathbb{R}$ indexed line ensemble $\bm{\mathcal{L}} = (\mathcal{L}_1, \mathcal{L}_2, \ldots ) \in \mathbb{Z}_{\ge 1} \times \mathcal{C} (\mathbb{R})$ satisfying the Brownian Gibbs property; set $q = 2^{1/6} \sigma^{1/3}$.
		
		\begin{enumerate} 
			\item Assume for any real number $\varepsilon > 0$ that there exists a constant $C = C(\varepsilon) > 1$ such that 
		\begin{flalign}
			\label{sigmat2}
			\mathbb{P} \big[- \sigma (1 + \varepsilon) t^2 - C \le \mathcal{L}_1 (t) \le -\sigma (1- \varepsilon) t^2 + C \big] \ge 1 - \varepsilon,
		\end{flalign}
	
		\noindent for any real number $t \in \mathbb{R}$. Then there exist two random variables $\mathfrak{l}, \mathfrak{c} \in \mathbb{R}$ and a rescaled parabolic Airy line ensemble $\bm{\mathcal{S}}^{(q)} = \big(\mathcal{S}_1^{(q)}, \mathcal{S}_2^{(q)}, \ldots \big) \in \mathbb{Z}_{\ge 1} \times \mathcal{C} (\mathbb{R})$ (as in \eqref{sigmar}) independent from them, such that $\mathcal{L}_j (t) = \mathcal{S}_j^{(q)} (t) + \mathfrak{l} t + \mathfrak{c}$, for each $(j, t) \in \mathbb{Z}_{\ge 1} \times \mathbb{R}$. 
		
		\item \label{constantk20} Further fix a real number $\ell \in \mathbb{R}$. Assume for any real number $\varepsilon > 0$ that there exists a constant $C = C(\varepsilon) > 1$ such that 
		\begin{flalign}
			\label{sigmalinear} 
			\mathbb{P} \big[ -\sigma t^2 + \ell t - \varepsilon |t| - C \le \mathcal{L}_1 (t) \le -\sigma t^2 + \ell t + \varepsilon |t| + C \big] \ge 1 - \varepsilon,
		\end{flalign}
		
		\noindent for any real number $t \in \mathbb{R}$. Then there exists a random variable $\mathfrak{c} \in \mathbb{R}$ and a rescaled parabolic Airy line ensemble $\bm{\mathcal{S}}^{(q)} = \big(\mathcal{S}_1^{(q)}, \mathcal{S}_2^{(q)}, \ldots \big) \in \mathbb{Z}_{\ge 1} \times \mathcal{C} (\mathbb{R})$ (as in \eqref{sigmar}) independent from $\mathfrak{c}$, such that $\mathcal{L}_j (t) = \mathcal{S}_j^{(q)} (t) + \ell t + \mathfrak{c}$, for each $(j, t) \in \mathbb{Z}_{\ge 1} \times \mathbb{R}$. 
		
		\end{enumerate} 
		
	\end{cor}

	From \Cref{constantk1}, we can quickly establish the following result characterizing extremal line ensembles satisfying the Brownian Gibbs property; it will also be established in \Cref{Proofsql1} below.

	\begin{cor}
		
		\label{linvariant} 
		
		Let $\bm{\mathcal{L}} = (\mathcal{L}_1, \mathcal{L}_2, \ldots ) \in \mathbb{Z}_{\ge 1} \times \mathcal{C} (\mathbb{R})$ be a $\mathbb{Z}_{\ge 1} \times \mathbb{R}$ indexed line ensemble; denote its associated probability measure by $\mu$. If $\mu \in \Tra (2^{-1/2})$ and $\mu$ is extremal, then there exists a (deterministic) constant $c \in \mathbb{R}$ such that $\mathcal{L}_j (t) = \mathcal{S}_j (t) + c$, for each $(j, t) \in \mathbb{Z}_{\ge 1} \times \mathbb{R}$. 
		
	\end{cor}
	
	\section{Proof of Characterization} 
	
	\label{ProofL} 
	
	In this section we establish \Cref{lr} assuming several statements that will be established later, which consist of two types of results. The first provides various properties of the line ensemble $\bm{\mathcal{L}}$ (defined on the infinite line $\mathbb{R}$) satisfying \Cref{l0}; they are given in \Cref{Estimatex} and \Cref{Global}. The second analyzes the asymptotic behaviors of families of non-intersecting Brownian bridges on finite intervals; they are given in \Cref{EventDerivativeSmooth} and \Cref{BridgesCurve}. We then establish \Cref{lr} in \Cref{Prooflr} and \Cref{Proof4a}; we conclude by establishing \Cref{constantk1} and \Cref{linvariant} as consequences of \Cref{lr} in \Cref{Proofsql1}.

	\subsection{On-Scale Events}
	
	\label{Estimatex}
	
	In this section we state two results indicating a coarse similarity between any line ensemble $\bm{\mathcal{L}}$ satisfying \Cref{l0} and the scaled parabolic Airy line ensemble $\bm{\mathcal{S}}$ of \Cref{ensemblewalks}. The first will imply that the top curve of $\bm{\mathcal{L}}$ is close to (within $o(n^{2/3})$ of) a parabola along a long interval (of length growing faster than $n^{2/3}$). The second will bound the locations of and gaps between (and also the H\"{o}lder regularity of) the paths in $\bm{\mathcal{L}}$, showing that they are of the same order as those in $\bm{\mathcal{S}}$. 
	
	To make these notions precise, we begin with the following two definitions. The first prescribes a certain mesh $\mathfrak{T}_k$, and also the event on which $\textbf{PAR}$ (from \Cref{eventparabola}) holds at each point on an interval. The second prescribes the event on which the $k$-th curve of a line ensemble is below the parabola $-2^{-1/2} t^2$, within by a given amount ($b$ and $B$, which could be negative). Throughout this section, we let $\bm{\mathsf{x}} = (\mathsf{x}_1, \mathsf{x}_2, \ldots ) \in \mathbb{Z}_{\ge 1} \times \mathcal{C}(\mathbb{R})$ denote a $\mathbb{Z}_{\ge 1} \times \mathbb{R}$ indexed line ensemble satisfying the Brownian Gibbs property.

	\begin{definition} 
		
		\label{parabolat} 
		
		For any integer $k \ge 1$; real numbers $\alpha \in (0, 1)$, $\varepsilon, C > 0$, and $A \ge 1$; and subset $\mathcal{T} \subset \mathbb{R}$ (which will often be a continuous interval or a finite mesh), define the event $\textbf{PAR}_{\varepsilon} (\mathcal{T}; C) = \textbf{PAR}_{\varepsilon}^{\bm{\mathsf{x}}} (\mathcal{T}; C)$ by
		\begin{flalign*}
			\textbf{PAR}_{\varepsilon} (\mathcal{T}; C) = \bigcap_{t \in \mathcal{T}} \textbf{PAR}_{\varepsilon}^{\bm{\mathsf{x}}} (t; C).
		\end{flalign*} 
		
		\noindent Further define the set $\mathfrak{T}_k (\alpha; A) \subset \mathbb{R}$ by\index{T@$\mathfrak{T}_k (\alpha; A)$}
		\begin{flalign*}
		 \mathfrak{T}_k (\alpha; A) = \big\{ x \in [-Ak^{1/3}, Ak^{1/3}] : x\in (\alpha k^{1/3}) \cdot \mathbb{Z} \big\}.
		\end{flalign*}
		
	\end{definition} 

	\begin{definition} 
		
		\label{eventsregular1} 
		
		\noindent For any integer $k \ge 1$ and real numbers $t, b, B \in \mathbb{R}$ with $B \ge b$; define the \emph{medium position event} \index{M@$\textbf{MED}$; medium position event} $\textbf{MED}_k (t; b; B) = \textbf{MED}_k^{\bm{\mathsf{x}}} (t; b; B)$ by
		\begin{flalign*}
			& \quad \textbf{MED}_k (t; b; B) = \big\{ b \le - \mathsf{x}_k (t) - 2^{-1/2} t^2 \le B \big\}.
		\end{flalign*} 
		
		\noindent For any subset $\mathcal{T} \subseteq \mathbb{R}$, define the event $\textbf{MED}_k (\mathcal{T}; b; B) = \textbf{MED}_k^{\bm{\mathsf{x}}} (\mathcal{T}; b; B)$ by
		\begin{flalign*}
			\textbf{MED}_k ( \mathcal{T}; b; B ) = \bigcap_{t \in \mathcal{T}} \textbf{MED}_k (t; b; B).
		\end{flalign*}
		
		\noindent For each of these events, if $k = 1$, then we abbreviate the \emph{top curve events} $\textbf{TOP} (t; b)= \textbf{MED}_1 (t; b; b)$ and $\textbf{TOP} (\mathcal{T}; b) = \textbf{MED}_1 (\mathcal{T}; b; b)$.\index{T@$\textbf{TOP}$; top curve event} Observe in particular that $\textbf{TOP} (t; \varepsilon t^2 + C) = \textbf{PAR}_{\varepsilon} (t; C)$ for any real numbers $t \in \mathbb{R}$ and $\varepsilon, C > 0$.

	\end{definition}

	The following proposition states that, if $\textbf{PAR}_{\varepsilon}$ (with very small $\varepsilon$) holds at a sufficiently fine mesh of points, then the top curve $\mathsf{x}_1$ of $\bm{\mathsf{x}}$ is likely within $o(n^{2/3})$ of the parabola $2^{-1/2} t^2$, at every point on an interval of length larger than $n^{1/3}$. It is established in \Cref{Proofx1} below.

	\begin{prop} 
		
		\label{x1small0} 
		
		Fix real numbers $\varepsilon, \delta \in ( 0, 1/4)$ and $A \ge 1$, and denote $\vartheta = 1000 A^2 \varepsilon$. There exists a constant $C = C(\varepsilon, \delta, A) > 1$ such that, for any integer $k \ge C$, we have
		\begin{flalign*}
			\mathbb{P} \Big[ \textbf{\emph{PAR}}_{\varepsilon} \big( \mathfrak{T}_k (\varepsilon; 15A); \varepsilon k^{2/3} \big) \cap \textbf{\emph{TOP}} \big( [-10Ak^{1/3}, 10Ak^{1/3}]; \vartheta k^{2/3} \big)^{\complement} \Big] \le \delta.
		\end{flalign*}
		
	\end{prop} 
	
	The following corollary applies \Cref{x1small0} to a line ensemble $\bm{\mathcal{L}}$ satifsying \Cref{l0}. 
	
	\begin{cor} 
		
		\label{x1lsmall}  
		
		Adopt \Cref{l0} and fix real numbers $B, \delta, \vartheta > 0$. There exists a constant $C = C(B, \delta, \vartheta) > 1$ such that, for $n \ge C$, we have    
		\begin{flalign*}
			\mathbb{P} \Big[ \textbf{\emph{TOP}}^{\bm{\mathcal{L}}} \big( [-Bn^{1/3}, Bn^{1/3}]; \vartheta n^{2/3} \big) \Big] \ge 1 - \delta.
		\end{flalign*}
		
	\end{cor} 
	
	\begin{proof} 
		
		Since $\textbf{TOP}^{\bm{\mathcal{L}}} \big( [-Bn^{1/3}, Bn^{1/3}]; \vartheta n^{2/3} \big) \subseteq \textbf{TOP}^{\bm{\mathcal{L}}} \big( [-B'n^{1/3}, B'n^{1/3}]; \vartheta n^{2/3} \big)$ whenever $B \ge B'$, we may assume that $B \ge 10$. Define the real numbers $A = B / 10$ and $\varepsilon = \vartheta / 1000 A^2$. Then \Cref{l0} implies for sufficiently large $n$ that 
		\begin{flalign*}
			\displaystyle\inf_{t \in \mathbb{R}} \mathbb{P} \big[ \textbf{PAR}_{\varepsilon}^{\bm{\mathcal{L}}} ( t; \varepsilon n^{2/3}) \big] \geq 1 - \displaystyle\frac{\varepsilon \delta}{90A}.
		\end{flalign*}
		
		\noindent This, a union bound, and the fact that $\big|\mathfrak{T}_n (\varepsilon; 15A) \big| \le 45A \varepsilon^{-1}$, together imply that
		\begin{flalign}
			\label{l1} 
			\mathbb{P} \Big[ \textbf{PAR}_{\varepsilon}^{\bm{\mathcal{L}}} \big( \mathfrak{T}_n (\varepsilon; 15 A); \varepsilon n^{2/3} \big) \Big] \ge 1 - \displaystyle\frac{\delta}{2}.
		\end{flalign}
		
		\noindent  It follows that, for sufficiently large $n$, we have 
		\begin{flalign*}
			\mathbb{P} \Big[ \textbf{TOP}^{\bm{\mathcal{L}}}  \big( [ & -B n^{1/3}, B n^{1/3}]; \vartheta n^{2/3} \big) \Big]  = \mathbb{P} \Big[ \textbf{TOP}^{\bm{\mathcal{L}}} \big( [-10A n^{1/3}, 10A n^{1/3}]; \vartheta n^{2/3} \big) \Big] \ge 1 - \delta,
		\end{flalign*}			
		
		\noindent where in the first statement we used the definition of $A$; in the second, we applied \eqref{l1}, \Cref{x1small0} (with the $\delta$ there equal to $\delta/2$ here), and a union bound. This yields the corollary.
	\end{proof}

	We next define two additional events. The first states an upper bound for the gaps between the paths in $\bm{\mathsf{x}}$, indicating that they are comparable to those of the Airy line ensemble (in which the distance between the $i$-th and $j$-th curves is of order about $|j^{2/3} - i^{2/3}|$, which can be deduced from \Cref{kdeltad} below). The second provides a H\"{o}lder type estimate for the paths in $\bm{\mathsf{x}}$ (that is fairly weak in comparison to the one that holds for the Airy line ensemble).

	\begin{definition} 
		\label{gap} 
		
		For any integer $k \ge 1$; real number $B \in \mathbb{R}$; and subset $\mathcal{T} \subset \mathbb{R}$, define the \emph{gap event} $\textbf{GAP}_k (\mathcal{T}; B) = \textbf{GAP}_k^{\bm{\mathsf{x}}} (\mathcal{T}; B)$ by\index{G@$\textbf{GAP}$; gap event}
		\begin{align}\label{e:defgap}
			\textbf{GAP}_k \big( \mathcal{T}; B) = \bigcap_{t \in \mathcal{T}} \bigcap_{1\leq i < j\leq k} \big\{ \mathsf{x}_i (t) - \mathsf{x}_j (t) \le B (j^{2/3} - i^{2/3}) + (\log k)^{25} i^{-1/3}  \big\}.
		\end{align}
		
	\end{definition}

	\begin{definition} 
		\label{eventtsregular2} 
		
		For any integers $k, n \ge 1$; real number $B \ge 0$; and subset $\mathcal{T} \subset \mathbb{R}$, define the \emph{H\"{o}lder regular event} $\textbf{REG}_k (\mathcal{T}; B; n) = \textbf{REG}_k^{\bm{\mathsf{x}}} (\mathcal{T}; B; n)$\index{R@$\textbf{REG}$; H\"{o}lder regular event} by
		\begin{flalign} 
			\label{eventsregular0}
			\textbf{REG}_k \big( \mathcal{T}; B; n) = \bigcap_{t, t+s \in \mathcal{T}} \Big\{ \big| \mathsf{x}_k (t+s) - \mathsf{x}_k (t) \big| \le 4 \big( n |s| \big)^{1/2} + B k^{1/3} |s| + k^{-25} \Big\}.
		\end{flalign}
		
	\end{definition}

	We next define an event that is formed from intersecting the ones above; it prescribes when the gaps and locations of $\bm{\mathsf{x}}$ are ``on-scale'' with respect to (that is, within constant factors of) those in the parabolic Airy line ensemble (in addition to imposing the H\"{o}lder type regularity of \Cref{eventtsregular2}). In what follows, if one examines curves $\mathsf{x}_k$ with $k$ of order $n$, then the relevant scales of the time $t$ and space $\mathsf{x}$ parameters are $n^{1/3}$ and $n^{2/3}$, respectively; see \Cref{f:scaling}.
	
	\begin{definition} 
		
		\label{eventscl} 
		
		For any integer $n \ge 1$ and real numbers $A, B, R > 0$, define the \emph{on-scale event} $\textbf{SCL}_n (A; B; R) = \textbf{SCL}_n^{\bm{\mathsf{x}}} (A; B; R)$ by\index{S@$\textbf{SCL}$; on-scale event}
		\begin{flalign*}
			\textbf{SCL}_n (A; B; R) & = \bigcap_{k= \lceil n/B \rceil}^{\lfloor Bn \rfloor } \textbf{MED}_k \Big( [-3An^{1/3}, 3An^{1/3}]; \displaystyle\frac{k^{2/3}}{100}; 450 k^{2/3} \Big) \\
			& \qquad \cap \bigcap_{k= \lceil n/B \rceil}^{\lfloor Bn \rfloor} \textbf{REG}_k \big( [-An^{1/3}, An^{1/3}]; 1000AB; Bn \big) \\
			& \qquad  \cap \textbf{GAP}_n \big( [-An^{1/3}, An^{1/3}]; R).
		\end{flalign*}
		
	\end{definition} 
	
	The next theorem indicates that, if the top curve of $\bm{\mathsf{x}}$ is close to a parabola on a long interval, then the on-scale event likely holds on another long (but slightly shorter) interval. It is proven in \Cref{ProofScale} below.
	
	\begin{thr}[On-scale estimates]
		
		\label{sclprobability} 
		
		For any real numbers $A, B \ge 2$, there exist constants $c = c(A, B) > 0$, $C_1 = C_1 (B) > 1$, and $C_2 = C_2 (A, B) > 1$ such that the following holds. For any real number $R \ge C_2$, we have
		\begin{flalign*}
			\mathbb{P} \Big[ \textbf{\emph{TOP}} \big( [-C_2 n^{1/3}, C_2 n^{1/3}]; C_1^{-1} n^{2/3} \big) \cap \textbf{\emph{SCL}}_n (A; B; R)^{\complement} \Big] \le c^{-1} e^{-c (\log n)^2}. 
		\end{flalign*}
		
	\end{thr}

	\begin{rem} 
		
		In \Cref{sclprobability}, and in many other places throughout this article, we will show that unfavorable events hold with a probability of the form $Ce^{-c(\log n)^2}$. The precise form of this estimate is not so significant, but it is important that it decays superpolynomially in $n$, as this will enable us to apply union bounds over very fine meshes. 
		
	\end{rem}

	\subsection{Global Law and Regular Profile Events} 
	
	\label{Global}
	
	Recall that \Cref{sclprobability} indicates when paths in a line ensemble are within constant factors of those $\mathcal{S}_j (t)$ in the scaled parabolic Airy line ensemble $\bm{\mathcal{S}}$. Those of the latter are known to concentrate around a deterministic profile. Throughout, we define the function $\mathfrak{G} : \mathbb{R} \times \mathbb{R}_{\ge 0} \rightarrow \mathbb{R}$ by setting
	\begin{flalign}
		\label{gfunction0} 
		\mathfrak{G} (t, x) = -2^{-1/2} t^2 - 2^{-7/6} (3 \pi)^{2/3} x^{2/3}, \qquad \text{for any $(t, x) \in \mathbb{R} \times \mathbb{R}_{\ge 0}$}. 
	\end{flalign} 

	\noindent Then we have (see \Cref{kdeltad} or \Cref{airy} below) with high probability that
	\begin{flalign}
		\label{rjt} 
		 \mathcal{S}_j (t) \approx \mathfrak{G} (t, j) + \mathcal{O} (j^{o(1) - 1/3}) = n^{2/3} \cdot \mathfrak{G} (tn^{-1/3}, jn^{-1}) + \mathcal{O} (j^{o(1) - 1/3}), 
	\end{flalign}

	\noindent for any integer $n \ge 1$. We first state a result indicating that the curves $\mathcal{L}_j (t)$ of an ensemble $\bm{\mathcal{L}}$ under \Cref{l0}  satisfy the bound $\mathcal{L}_j (t) = n^{2/3} \cdot \mathfrak{G} (tn^{-1/3}, jn^{-1}) + o(n^{2/3})$. This might be viewed as a global law, or limit shape, for the line ensemble $\bm{\mathcal{L}}$. It is weaker than \eqref{rjt} but improves on the $\textbf{MED}$ event part appearing in $\textbf{SCL}$ (recall \Cref{eventscl}) arising in \Cref{sclprobability}.
	 
	We begin with the following definition for the event on which the global law holds.

	\begin{definition} 
		
		\label{functiong} 
	   
		 Fix an infinite sequence $\bm{\mathsf{x}} = (\mathsf{x}_1, \mathsf{x}_2, \ldots ) \in \mathbb{Z}_{\ge 1} \times \mathcal{C}(\mathbb{R})$ of continuous functions. For any integer $n \ge 1$, and real numbers $\delta, B > 0$, define the \emph{global law event} $\textbf{GBL}_n (\delta; B) = \textbf{GBL}_n^{\bm{\mathsf{x}}} (\delta; B)$\index{G@$\textbf{GBL}$; global law event} by (recalling \eqref{gfunction0})
		\begin{flalign*}
			\textbf{GBL}_n^{\bm{\mathsf{x}}} (\delta; B) = \bigcap_{|t| \le Bn^{1/3}} \bigcap_{j = 1}^{\lfloor Bn \rfloor} \Big\{ \big| \mathsf{x}_j (t) - n^{2/3} \cdot \mathfrak{G} (t n^{-1/3}, jn^{-1}) \big| \le \delta n^{2/3} \Big\}.
		\end{flalign*}
	\end{definition}

	The following theorem, to be established in \Cref{ProofxGlobal} below, states that the global law event likely holds for the ensemble $\bm{\mathcal{L}}$ satisfying \Cref{l0}. 
	
	\begin{thr}[Global law]
		
		\label{p:globallaw2}
	
		Adopt \Cref{l0}, and fix real numbers $B > 1$ and $\delta > 0$. There exists a constant $C = C(B, \delta) > 0$ such that, for $n \ge C$, we have     
		\begin{flalign*}
			\mathbb{P} \big[ \textbf{\emph{GBL}}_n^{\bm{\mathcal{L}}} (\delta; B) \big] \ge 1 - \delta.
		\end{flalign*}

	\end{thr}

	We next state results indicating that the locations of the paths in a line ensemble $\bm{\mathcal{L}}$ satisfying \Cref{l0} approximate a ``regular profile.'' The following definition makes that notion more precise. The first bound in \eqref{xjtgammat} below states that $\bm{x}$ at time $t$ approximates the profile $\gamma_t$; the second states that $\gamma_t$ is regular.

	\begin{definition}		
		
		\label{eventpfl} 
		
		Fix real numbers $a < b$, and let $\bm{x} = (x_1, x_2, \ldots , x_n) \in \llbracket 1, n \rrbracket \times \mathcal{C} \big( [a, b] \big)$ denote a sequence of functions. For any real numbers $\delta, B > 0$ and $t \in [a, b]$, we define the \emph{regular profile} event\index{P@$\textbf{PFL}$; regular profile event} $\textbf{PFL}^{\bm{x}} (t; \delta; B)$ to be that on which there exists a continuous function $\gamma_t : [0, 1] \rightarrow \mathbb{R}$ such that 
		\begin{flalign}
			\label{xjtgammat} 
			\displaystyle\max_{j \in \llbracket 1, n \rrbracket} \big| x_j (t) - \gamma_t (jn^{-1}) \big| \le \delta, \qquad \text{and} \qquad \big\| \gamma_t - \gamma_t (0) \big\|_{\mathcal{C}^{50}} \le B.
		\end{flalign}
		
	\end{definition}

	\begin{rem} 
		
		\label{pfl0event} 
		
		Observe that those $\bm{x}$ satisfying $\textbf{PFL}^{\bm{x}} (t; \delta; B)$ define a closed (and thus measurable) subset of $\llbracket 1, n \rrbracket \times \mathcal{C} \big( [a, b] \big)$. Indeed, let $\mathcal{Y}$ denote the set of continuous functions $\gamma_t : [0, 1] \rightarrow \mathbb{R}$ satisfying the second statement of \eqref{xjtgammat}. Then, the function $\max_{j \in \llbracket 1 ,n \rrbracket} \big| x_j (t) - \gamma_t (jn^{-1}) \big|$ is continuous in $\bm{x} \in \llbracket 1, n \rrbracket \times \mathcal{C} \big( [a, b] \big)$, uniformly in $\bm{x} \in \llbracket 1, n \rrbracket \times \mathcal{C} \big( [a, b] \big)$ and $\gamma_t \in \mathcal{Y}$. Hence, the function mapping $\bm{x} \in \llbracket 1, n \rrbracket \times \mathcal{C} \big( [a, b] \big)$ to $\inf_{\gamma_t \in \mathcal{Y}} \max_{j \in \llbracket 1, n \rrbracket} \big| x_j (t) - \gamma_t (jn^{-1}) \big|$ is continuous. Since $\textbf{PFL}^{\bm{x}} (t; \delta; B)$ is the event on which this function\footnote{The infimum over $\gamma_t \in \mathcal{Y}$ in it is attained, as a quick consequence of the Arzel\`{a}--Ascoli theorem.} is at most $\delta$, it is closed. 
		
	\end{rem}

	We will show through the following theorem that the $\{ n+1, n+2, \ldots , 2n \}$-th curves of the line ensemble $\bm{\mathcal{L}}$ satisfying \Cref{l0} satisfy the regular profile event with high probability, after rescaling and restricting to an intersection of $\textbf{TOP}$ events (recall \Cref{eventsregular1}). It will be established in \Cref{ProofRegular0} below. 	      
	
	\begin{thr}[Spatial regularity]
		
		\label{p:closerho0}
		
		Adopt \Cref{l0}, and fix a real number $A > 1$. There exist constants $\omega \in (0, 1/2)$, $c = c(A) \in (0, 1)$, and $C = C(A) > 1$ such that the following holds. Let $n \ge 1$ be an integer, and define the ensemble $\bm{l} = \bm{l}^{(n)} = (l_1, l_2, \ldots , l_n) \in \llbracket 1, n \rrbracket \times \mathcal{C} (\mathbb{R})$ by setting $l_j (t) = n^{-2/3} \cdot \mathcal{L}_{j+n} (tn^{1/3})$, for each $(j, t) \in \llbracket 1, n \rrbracket \times \mathbb{R}$. Then, we have 
		\begin{flalign*}
			\mathbb{P} \Bigg[ \bigcup_{|t| \le A} \textbf{\emph{PFL}}^{\bm{l}} \big( t;  n^{-1} (\log n)^7; C \big)^{\complement} \cap \bigcap_{k=1}^{\lfloor \omega^{-2} \rfloor} \textbf{\emph{TOP}}^{\bm{\mathcal{L}}} \big( [-Cn^{k\omega/3}, C n^{k\omega/3} &]; cn^{2k\omega/3}) \Bigg] \le Cn^{-50}.
		\end{flalign*}

	\end{thr}

	\subsection{Second Derivative Estimates for Paths} 
	
	\label{EventDerivativeSmooth} 
	
	In this section we state \Cref{h0x2} below. It indicates that, under certain conditions, the paths in a family of non-intersecting Brownian bridges are close to (random) curves with nearly constant second derivatives. The following assumption more precisely prescribes these conditions on the bridges, which will take place on the time interval $[-\xi n^{1/3}, \xi n^{1/3}]$ for some bounded $\xi > 0$. The first condition \eqref{gnuv} below states that the boundary data (consisting of the entrance and exit data, as well as the lower and upper boundaries) are approximated by a function $G$, which is a shift of the $\mathfrak{G}$ from \eqref{gfunction0}. The second \eqref{pflx2} states that the bridges likely satisfy a regular profile event at any fixed time $s \in [-\xi n^{1/3}, \xi n^{1/3}]$.

	\begin{assumption} 
		
		\label{derivativegxpfl}
		
		Define the function $G: \mathbb{R} \times \mathbb{R}_{\ge 0}$ by, for any $(t, x) \in \mathbb{R}_{\ge 0} \times \mathbb{R}$, setting 
		\begin{flalign}
			\label{g2121} 
			G(t,x) = \mathfrak{G}(t, x+1) = - 2^{-1/2} t^2 - 2^{-7/6} (3 \pi)^{2/3} (x+1)^{2/3}. 
		\end{flalign}
		
		\noindent Fix real numbers $\delta \in (0, 1/2)$ and $B > 1$; let $\xi \in (B^{-1}, B)$ be a real number; set $\mathfrak{R} \in [-\xi, \xi] \times [0, 1]$; and assume that $B>1$ is sufficiently large so that   
		\begin{flalign}
			\label{g50}
			\| G \|_{\mathcal{C}^{50} (\mathfrak{R})} \le B. 
		\end{flalign}
	
		\noindent Further let $\mathsf{T} = \xi n^{1/3}$. Fix  two $n$-tuples $\bm{u}, \bm{v} \in \mathbb{W}_n$ and two functions $f, g : [-\mathsf{T}, \mathsf{T}] \rightarrow \mathbb{R}$, such that $f < g$ and $f(-\mathsf{T}) < u_n < u_1 < g(-\mathsf{T})$ and $f(\mathsf{T}) < v_n < v_1 < g(\mathsf{T})$. Suppose that 
		\begin{flalign}
			\label{gnuv} 
			\begin{aligned} 
				& \displaystyle\max_{j \in \llbracket 1, n \rrbracket} \big| n^{-2/3} u_j - G (-\xi, jn^{-1}) \big| < \delta; \qquad \qquad \displaystyle\max_{j \in \llbracket 1, n \rrbracket} \big| n^{-2/3} v_j - G(\xi, jn^{-1}) \big| < \delta; \\
				& \displaystyle\sup_{s \in [-\mathsf{T}, \mathsf{T}]} \big| n^{-2/3} f (s) - G (s n^{-1/3}, 1) \big| < \delta; \qquad \displaystyle\sup_{s \in [-\mathsf{T}, \mathsf{T}]} \big| n^{-2/3} g(s) - G (s n^{-1/3}, 0) \big| < \delta.
			\end{aligned}
		\end{flalign}
		
		\noindent Let $\bm{\mathsf{x}} = (\mathsf{x}_1, \mathsf{x}_2, \ldots , \mathsf{x}_n) \in \llbracket 1, n \rrbracket \times \mathcal{C} \big( [-\mathsf{T}, \mathsf{T}] \big)$ be a family of $n$ non-intersecting Brownian bridges sampled from the measure $\mathsf{Q}_{f; g}^{\bm{u}; \bm{v}}$. Further define the rescaled family of non-intersecting Brownian bridges $\bm{x} = (x_1, x_2, \ldots , x_n) \in \llbracket 1, n \rrbracket \times \mathcal{C} \big( [-\xi, \xi] \big)$ by setting
		\begin{flalign}
			\label{xjsn} 
			x_j (s) = n^{-2/3} \cdot \mathsf{x}_j (n^{1/3} s), \qquad \text{for each $(j, s) \in \llbracket 1, n \rrbracket \times [-\xi, \xi]$}.
		\end{flalign}
		
		\noindent Assume for each real number $s \in [-\xi, \xi]$ that  
		\begin{flalign}
			\label{pflx2} 
			\mathbb{P} \big[ \textbf{PFL}^{\bm{x}} (s; n^{-19/20}; B) \big] \ge 1 - n^{-20}.
		\end{flalign}
		
	\end{assumption} 
	
	Observe that the event in \eqref{pflx2}  depends not only on the boundary data, but also on the random bridges in $\bm{\mathsf{x}}$ themselves. It imposes that we somehow knew ``in advance'' that these bridges likely have some regularity. In our eventual context, this knowledge will come from \Cref{p:closerho0}.
	
\begin{figure}
\centering
\begin{subfigure}{.45\textwidth}
  \centering
  \includegraphics[width=1\linewidth]{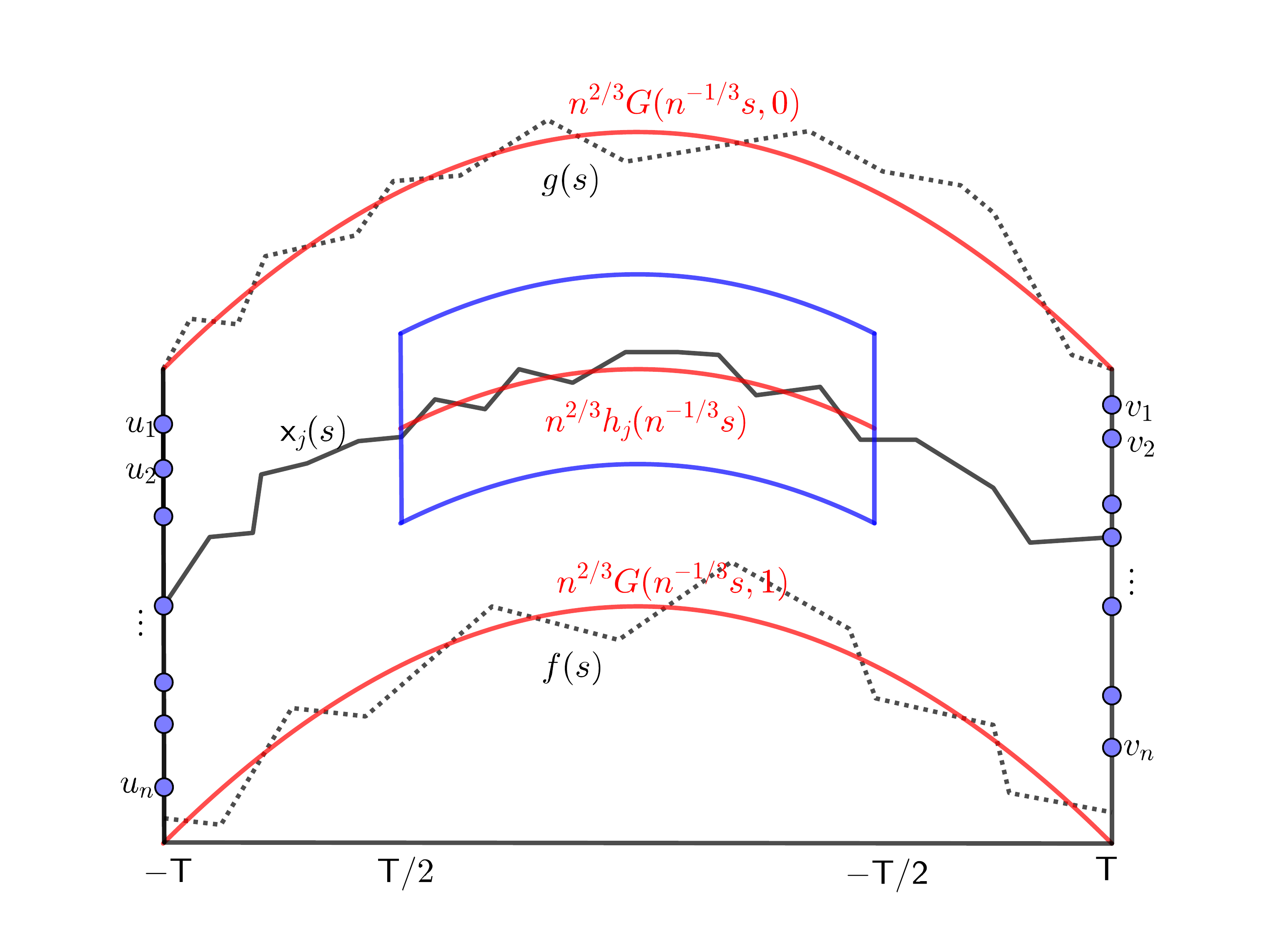}
\end{subfigure}%
\begin{subfigure}{.55\textwidth}
  \centering
  \includegraphics[width=1\linewidth]{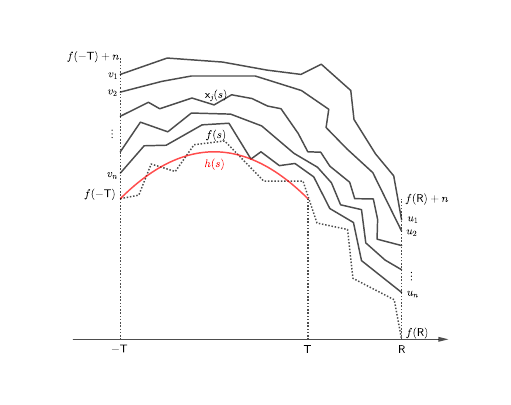}
\end{subfigure}
\caption{Shown to the left is a depiction for the conclusion of \Cref{h0x2}.  Shown to the right is a depiction for \Cref{xfdelta2} (where $\mathsf{T}$ and $\mathsf{R}$ are not drawn to scale).}
\label{f:assumption}
\end{figure}

	The following theorem, to be established in \Cref{Proofhj} below, states under \Cref{derivativegxpfl} that the paths in $\bm{\mathsf{x}}$ are near curves that have nearly constant second derivative $-2^{-1/2}$ in $t$. See the left side of \Cref{f:assumption}.

	\begin{thr}[Curvature approximation]
		
		\label{h0x2}
		
		Adopting \Cref{derivativegxpfl}, there exist constants $c = c(B) > $ and $C = C(B) > 1$ such that the following holds with probability at least $1 - C n^{-10}$ whenever $\delta < c$. For each integer $j \in \llbracket n/3, 2n/3 \rrbracket$, there exists a (random) twice-differentiable function $h_j : [ -\xi/2, \xi/2] \rightarrow \mathbb{R}$ with 
		\begin{flalign}
			\label{hjxj2} 
			\displaystyle\sup_{|s| \le \xi/2} \big| \partial_s^2 h_j (s) + 2^{-1/2} \big| \le \delta^{1/8} + (\log n)^{-1/4}, \qquad \text{and} \qquad \| h_j \|_{\mathcal{C}^1} \le C,
		\end{flalign}
		
		\noindent such that 
		\begin{flalign}
			\label{hjxj1}
			\displaystyle\sup_{|s| \le \mathsf{T}/2} \big| \mathsf{x}_j (s) - n^{2/3} \cdot h_j (s n^{-1/3}) \big| < n^{-1/5}.
		\end{flalign}
	\end{thr}

	\subsection{Brownian Bridges Above a Curve}
	
	\label{BridgesCurve}

	We will frequently make use of the process obtained by examining the Airy line ensemble $\bm{\mathcal{A}}$ at a given time, called the Airy point process.
	
	\begin{definition}
		
		\label{a0}
		
		Let $\bm{\mathcal{A}}$ denote the Airy line ensemble, as in \Cref{ensemblewalks}. The random infinite sequence\index{A@$\bm{\mathfrak{a}}$; Airy point process} $\bm{\mathfrak{a}} = (\mathfrak{a}_1, \mathfrak{a}_2, \ldots ) = \big( \mathcal{A}_1 (0), \mathcal{A}_2 (0), \ldots \big)$ of decreasing real numbers is the \emph{Airy point process}. 
		
	\end{definition} 

	In this section we state two results. The first indicates that, under certain conditions, the gaps between paths at a single time for an ensemble of non-intersecting Brownian bridges converge to those of the Airy point process. We begin by describing these conditions more precisely through the following assumption. The bound \eqref{fsfs} imposes that the lower bounary $f$ is close to a function with nearly constant second derivative; the bound \eqref{probabilityfnt} imposes that the entrance data for the ensemble satisfies a global law (analogous to \Cref{functiong}).   
	
	\begin{assumption}
		
		\label{xfdelta2}

		Let $n \ge 1$ be an integer; $\bm{\delta} = (\delta_1, \delta_2, \ldots ) \subset ( 0, 1/4)$ be a non-increasing sequence of real numbers satisfying $\lim_{k \rightarrow \infty} \delta_k = 0$ and $\delta_k \ge k^{-1/10}$ for each integer $k > 2^{20}$; and $\mathsf{T} = \mathsf{T}_n$ be a real number such that $\delta_n^{-1} n^{1/3} \le \mathsf{T} \le n^{1/2}$. Further set $\mathsf{R} = \mathsf{R}_n = n^{20} \ge \mathsf{T}$, and fix a function $f = f_n : [-\mathsf{T}, \mathsf{R}] \rightarrow \mathbb{R}$. Assume that there exists a function $h = h_n : [-\mathsf{T}, \mathsf{T}] \rightarrow \mathbb{R}$ satisfying
		\begin{flalign}
			\label{fsfs}
			\sup_{s \in [-\mathsf{T}, \mathsf{T}]} \big| \partial_s^2 h (s) + 2^{-1/2} \big| \le \delta_n; \qquad \displaystyle\sup_{s \in [-\mathsf{T}, \mathsf{T}]} \big| f (s) - h (s) \big| < \delta_n.  
		\end{flalign} 
		\noindent Let $\bm{u} = \bm{u}^n \in \mathbb{W}_n$ and $\bm{v} = \bm{v}^n \in \mathbb{W}_n$ be sequences such that 
		\begin{flalign}
			\label{fvvf} 
			 f(\mathsf{R}) \le v_n \le v_1 \le f(\mathsf{R}) + n.
		\end{flalign} 
		
		\noindent Sample non-intersecting Brownian bridges $\bm{\mathsf{x}} =\bm{\mathsf{x}}^n= (\mathsf{x}_1, \mathsf{x}_2, \ldots , \mathsf{x}_n) \in \llbracket 1, n \rrbracket \times \mathcal{C} \big( [-\mathsf{T}, \mathsf{R}] \big)$ under $\mathsf{Q}_{f}^{\bm{u}; \bm{v}}$. See the right side of \Cref{f:assumption}.
		
		Moreover, for any $t \in [-\mathsf{T}, \mathsf{R}]$, define the event
		\begin{flalign}
			\label{tfn}
			\mathscr{F} (t) = \mathscr{F}_n (t) = \bigg\{ \displaystyle\max_{j \in \llbracket 1, n \rrbracket} \big| \mathsf{x}_j (t) - 2^{-7/6} (3 \pi)^{2/3} (n^{2/3} - j^2) - f_n (t) \big| \le \delta_n n^{2/3} \bigg\}.
		\end{flalign}
		
		\noindent Assume that 
		\begin{flalign}
			\label{probabilityfnt} 
			\mathbb{P} \big[ \mathscr{F}_n (t) \big] \ge 1 - \delta_n, \qquad \text{for each real number $t \in [-n^{1/3}, n^{1/3}]$}.
		\end{flalign}      
	\end{assumption} 

	\noindent Similarly to \eqref{pflx2}, \eqref{probabilityfnt} imposes that we knew in advance that the curves in $\bm{\mathsf{x}}^n$ likely approximate a specific deterministic function at intermediate times $t$ in the domain. In our eventual context, this knowledge will come from \Cref{p:globallaw2}.
	
	The following theorem indicates that, under these two assumptions, the gaps between the bridges in $\bm{\mathsf{x}}^n$ converge to those of the Airy point process $\bm{\mathfrak{a}}$ (from \Cref{a0}). It is established in \Cref{DifferenceUpper} below.

	\begin{thr}[Airy gaps]
		
		\label{xdifferenceconverge}
		
		Adopt \Cref{xfdelta2}, and fix an integer $k \ge 1$ and a real number $t \in \mathbb{R}$. As $n$ tends to $\infty$, the $k$-tuple of gaps 
		\begin{flalign*} 
			2^{1/2} \cdot \big( \mathsf{x}_1 (t) - \mathsf{x}_2 (t), \mathsf{x}_2 (t)- \mathsf{x}_3 (t), \ldots , \mathsf{x}_k (t) - \mathsf{x}_{k+1} (t) \big),
		\end{flalign*}
		
		\noindent converges in law to that  $(\mathfrak{a}_1 - \mathfrak{a}_2, \mathfrak{a}_2 - \mathfrak{a}_3, \ldots , \mathfrak{a}_k - \mathfrak{a}_{k+1})$ of the Airy point process $\bm{\mathfrak{a}}$.
	\end{thr}
	
	The second result of this section indicates that, if the infinite ensemble $\bm{\mathcal{L}}$ satisfying \Cref{l0} has the property that the gaps between its paths at any time converges to those of the Airy point process, then $\bm{\mathcal{L}}$ must be a scaled parabolic Airy line ensemble, up to a (possibly random) affine shift. It is established in \Cref{ProofL00} below.

	\begin{prop} 
		
		\label{qsf}
		
		Adopt \Cref{l0}. Further assume for any integer $k \ge 1$ and real number $t \in \mathbb{R}$ that the $k$-tuple of gaps $2^{1/2} \cdot \big( \mathcal{L}_1 (t) - \mathcal{L}_2 (t), \mathcal{L}_2 (t) - \mathcal{L}_3 (t), \ldots , \mathcal{L}_k (t) - \mathcal{L}_{k+1} (t) \big)$ has the same law as that $(\mathfrak{a}_1 - \mathfrak{a}_2, \mathfrak{a}_2 - \mathfrak{a}_3, \ldots , \mathfrak{a}_k - \mathfrak{a}_{k+1})$ of the Airy point process $\bm{\mathfrak{a}}$. Then, there exist two random variables $\mathfrak{l}, \mathfrak{c} \in \mathbb{R}$, and a scaled parabolic Airy line ensemble $\bm{\mathcal{S}} = (\mathcal{S}_1, \mathcal{S}_2, \ldots ) \in \mathbb{Z}_{\ge 1} \times \mathcal{C}(\mathbb{R})$ independent from them, such that $\mathcal{L}_j (t) = \mathcal{S}_j (t) + \mathfrak{l} t + \mathfrak{c}$ for each $(j, t) \in \mathbb{Z}_{\ge 1} \times \mathbb{R}$. 
		
	\end{prop} 
	
	Let us mention that, although \Cref{qsf} as stated is a result about the line ensemble $\bm{\mathcal{L}}$ on the infinite line $\mathbb{R}$, it will quickly be reduced to one about line ensembles on finite intervals (see \Cref{uvxrconverge} below), which is our reason for including it here.

	\subsection{Proof of \Cref{lr}} 
	
	\label{Prooflr} 
	
	In this section we use the previous results to establish \Cref{lr}. We begin with the following lemma that will enable us to verify \Cref{xfdelta2} of \Cref{xdifferenceconverge}. In what follows, we recall $\mathcal{F}_{\ext}$ from \Cref{property}. It is proven at the end of this section.

	\begin{lem} 
		
	\label{fhl} 
	
	Adopt \Cref{l0}. For any real number $\delta \in (0, 1/4)$, there exists a constant $N_0 = N_0 (\delta) > 1$ so that the following holds. Let $n \ge N_0$ be any integer; set $\mathsf{T} = \delta^{-1} n^{1/3}$ and $\mathsf{R} = n^{20}$. There exists an event $\mathscr{A} = \mathscr{A}_n (\delta)$ with $\mathbb{P}[\mathscr{A}] \ge 1 - \delta$, measurable with respect to $\mathcal{F}_{\ext} = \mathcal{F}_{\ext} \big( \llbracket 1, n \rrbracket \times (-\mathsf{T}, \mathsf{R}) \big)$, such that, conditional on $\mathcal{F}_{\ext}$ and restricting to $\mathscr{A}$, the following three statements hold. 

	\begin{enumerate}
	\item We have $\mathcal{L}_1 (\mathsf{R}) \le \mathcal{L}_{n+1} (\mathsf{R}) + n$.
	\item For each $t \in [-n^{1/3}, n^{1/3}]$, we have  
	\begin{flalign*}
		 \mathbb{P}  \Bigg[ \displaystyle\max_{j \in \llbracket 1, n \rrbracket} \big| \mathcal{L}_j (t) + 2^{-1/2} t^2 + 2^{-7/6} (3 \pi)^{2/3} j^{2/3} \big| \le \displaystyle\frac{\delta n^{2/3}}{2}  \Bigg] \ge 1 - \delta.
	\end{flalign*}
		\item There exists a twice-differentiable function $h = h_n : [-\mathsf{T}, \mathsf{T}] \rightarrow \mathbb{R}$ such that 
	\begin{flalign}
		\label{th} 
		\begin{aligned}
			& \sup_{|t| \le \mathsf{T}} \big| \mathcal{L}_{n+1} (t) - h (t) \big| < n^{-1/6}; \qquad \displaystyle\sup_{|t| \le \mathsf{T}} \big| \partial_t^2 h(t) + 2^{-1/2} \big| \le \delta; \\
			& \qquad \qquad \displaystyle\sup_{|t| \le \mathsf{T}} \big| h(t) + 2^{-1/2} t^2 + 2^{-7/6} (3 \pi)^{2/3} n^{2/3} \big| \le \displaystyle\frac{\delta n^{2/3}}{2}.
		\end{aligned} 
	\end{flalign}	
\end{enumerate}  
	
	\end{lem} 
	
	Given \Cref{fhl}, we can quickly establish \Cref{lr} using the Airy gaps \Cref{xdifferenceconverge} and \Cref{qsf}.

	\begin{proof}[Proof of \Cref{lr}]

		By \Cref{fhl}, there is a non-increasing sequence $\bm{\delta} = (\delta_1, \delta_2, \ldots )$ of real numbers with $\lim_{j \rightarrow \infty} \delta_j = 0$ and $\delta_j \ge j^{-1/10}$ for each integer $j > 2^{20}$, such that the events $\mathscr{A}_n = \mathscr{A}_n (\delta_n)$ satisfying the properties listed in \Cref{fhl} (with each appearance of $\delta$ there replaced by $\delta_n$ here) exist. Set $\mathsf{T}_n = \delta_n^{-1} n^{1/3}$ and $\mathsf{R}_n = n^{20}$, and condition on $\mathcal{F}_{\ext} \big( \llbracket 1, n \rrbracket \times (-\mathsf{T}_n, \mathsf{R}_n) \big)$.
		
		Set $\bm{u}^n = \big( \mathcal{L}_1 (-\mathsf{T}_n), \mathcal{L}_2 (-\mathsf{T}_n), \ldots , \mathcal{L}_n (-\mathsf{T}_n) \big)$ and $\bm{v}^n = \big( \mathcal{L}_1 (\mathsf{R}_n), \mathcal{L}_2 (\mathsf{R}_n), \ldots , \mathcal{L}_n (\mathsf{R}_n) \big)$; define $f_n : [-\mathsf{T}_n, \mathsf{T}_n] \rightarrow \mathbb{R}$ by setting $f_n (s) = \mathcal{L}_{n+1} (s)$ for each $s \in [-\mathsf{T}_n, \mathsf{T}_n]$; and sample $\bm{\mathsf{x}}^n = (\mathsf{x}_1^n, \mathsf{x}_2^n, \ldots , \mathsf{x}_n^n) \in \llbracket 1, n \rrbracket \times \mathcal{C} \big( [-\mathsf{T}_n, \mathsf{T}_n] \big)$ under the measure $\mathsf{Q}_{f;g}^{\bm{u}^n; \bm{v}^n}$. Then, for any real numbers $r_1, r_2, \ldots , r_k \in \mathbb{R}$, we have by the Brownian Gibbs property that 
		\begin{flalign}
			\label{ljlj1t} 
			\mathbb{P} \bigg[ \bigcap_{j=1}^k \Big\{ 2^{-1/2} \big(\mathcal{L}_j (t) - \mathcal{L}_{j+1}(t) \big) \ge r_j \Big\} \bigg| \mathscr{A}_n \bigg] = \mathbb{P} \bigg[ \bigcap_{j=1}^k \Big\{ 2^{-1/2} \big( \mathsf{x}_j^n (t) - \mathsf{x}_{j+1}^n (t) \big) \ge r_j \Big\} \bigg| \mathscr{A}_n \bigg].
		\end{flalign}

		Next, observe upon restricting to the event $\mathscr{A}_n$ that \eqref{fsfs} and \eqref{fvvf} are verified by the first property in \Cref{fhl}, together with the first two bounds in \eqref{th}; moreover, \eqref{probabilityfnt} is verified by the last bound in \eqref{th}, together with the second property in \Cref{fhl}. Thus, we have by the Airy gaps \Cref{xdifferenceconverge} that   
		\begin{flalign}
			\label{xr2} 
			\displaystyle\lim_{n \rightarrow \infty} \mathbb{P} \bigg[ \bigcap_{j=1}^k \Big\{ 2^{-1/2} \big( \mathsf{x}_j^n (t) - \mathsf{x}_{j+1}^n (t) \big) \ge r_j \Big\} \bigg| \mathscr{A}_n \bigg] = \mathbb{P} \bigg[ \bigcap_{j=1}^k \{ \mathfrak{a}_j - \mathfrak{a}_{j+1} \ge r_j \} \bigg].
		\end{flalign} 
		
		\noindent Thus,
			\begin{flalign*} 
			\mathbb{P} \bigg[ \bigcap_{j=1}^k \Big\{ 2^{-1/2} \big( \mathcal{L}_j (t) - \mathcal{L}_{j+1} (t) \big) \ge r_j \Big\} \bigg] & = \displaystyle\lim_{n \rightarrow \infty}  \mathbb{P} \bigg[ \bigcap_{j=1}^k \Big\{ 2^{-1/2} \big( \mathsf{x}_j^n (t) - \mathsf{x}_{j+1}^n (t) \big) \ge r_j \Big\} \bigg| \mathscr{A}_n \bigg] \\ 
			& = \mathbb{P} \bigg[ \bigcap_{j=1}^k \{ \mathfrak{a}_j - \mathfrak{a}_{j+1} \ge r_j \} \bigg].
		\end{flalign*} 
	
		\noindent Here, the first statement follows from the limit as $n$ tends to $\infty$ of \eqref{ljlj1t}, and the facts that $\mathbb{P}[\mathscr{A}_n] \ge 1 - \delta_n$ (by \Cref{fhl}) and $\lim_{n \rightarrow \infty} \delta_n = 0$; the second follows from \eqref{xr2}. Since the $(r_j)$ were arbitrary, we deduce that the law of $2^{1/2} \cdot \big( \mathcal{L}_1 (t) - \mathcal{L}_2 (t), \mathcal{L}_2 (t) - \mathcal{L}_3 (t), \ldots , \mathcal{L}_k (t) - \mathcal{L}_{k+1} (t) \big)$ coincides with that of the gaps $(\mathfrak{a}_1 - \mathfrak{a}_2, \mathfrak{a}_2 - \mathfrak{a}_3, \ldots , \mathfrak{a}_k - \mathfrak{a}_{k+1})$ for the Airy point process, for any integer $k \ge 1$ and real number $t \in \mathbb{R}$. Thus, the theorem follows from \Cref{qsf}. 
	\end{proof}

	We now establish \Cref{fhl}; we adopt the notation and assumptions of that lemma in what follows. We will define $\mathscr{A}$ as the intersection $\mathscr{A} = \bigcap_{j=1}^3 \mathscr{A}^{(j)}$ of three events $\mathscr{A}^{(j)}$, measurable with respect to $\mathcal{F}_{\ext} = \mathcal{F}_{\ext} \big( \llbracket 1, n \rrbracket \times (-\mathsf{T}, \mathsf{R}) \big)$, that essentially correspond to the three parts of \Cref{fhl}. Let $\mathfrak{c}_1 > 0$, $\mathfrak{C}_1 > 1$, and $\mathfrak{C}_2 > 1$ denote the constants $c$, $C_1$, and $C_2$ from \Cref{sclprobability} at $(A, B, R) = (2, 2, \mathfrak{C}_2)$, respectively. We first define the event
	\begin{flalign}
		\label{a1} 
		\mathscr{A}^{(1)} = \textbf{GAP}_{n^{30}} (\mathsf{R}; \mathfrak{C}_2).
	\end{flalign}

	\noindent We next define $\mathscr{A}^{(2)} = \mathscr{A}_1^{(2)} \cap \mathscr{A}_2^{(2)}$ to be the event measurable with respect to $\mathcal{F}_{\ext}$, where $\mathscr{A}_1^{(2)}$ and $\mathscr{A}_2^{(2)}$ are given by
	\begin{flalign}
		\label{a2} 
		\begin{aligned} 
		\mathscr{A}_1^{(2)} & = \Big\{ \mathbb{P} \big[ \textbf{GBL}_n^{\bm{\mathcal{L}}} ( \delta^2; \delta^{-1} ) \big| \mathcal{F}_{\ext} \big] \ge 1 - \delta \Big\} \\
		 \mathscr{A}_2^{(2)} &=  \Bigg\{ \displaystyle\sup_{|t| \le \mathsf{T}} \big| \mathcal{L}_{n+1} (t) + 2^{-1/2} t^2 + 2^{-7/6} (3 \pi)^{2/3} n^{2/3} \big| < \displaystyle\frac{\delta n^{2/3}}{4} \bigg\},
		\end{aligned} 
	\end{flalign}
	
	\noindent where the probability is conditional on $\mathcal{F}_{\ext}$. Further let $\mathscr{A}^{(3)}$ denote the event that is measurable (as it defines an open set) with respect to $\mathcal{F}_{\ext}$, on which there exists a twice-differentiable function $h = h_n : [-\mathsf{T}, \mathsf{T}] \rightarrow \mathbb{R}$ satisfying the first two bounds in \eqref{th}.

	 The following lemmas say that each of the $\mathscr{A}^{(j)}$ is likely; we establish the former in this section and the latter in \Cref{Proof4a} below. In what follows, when we write that a statement holds, ``for $n$ sufficiently large,'' we means that there exists a constant $N_0 = N_0 (\delta)$ such that the statement holds whenever $n \ge N_0$. 
	 
	\begin{lem} 
		
		\label{a123} 
		
		For $n$ sufficiently large, we have $\mathbb{P} \big[ \mathscr{A}^{(1)} \cap \mathscr{A}^{(2)} \big] \ge 1 - \delta/2$.   
	\end{lem}

	\begin{lem} 
		
	\label{a4} 
	
	For $n$ sufficiently large, we have $\mathbb{P} \big[ \mathscr{A}^{(3)} \big] \ge 1- \delta/2$. 

	\end{lem} 

	Given \Cref{a123} and \Cref{a4}, we can quickly establish \Cref{fhl}. 
	
	\begin{proof}[Proof of \Cref{fhl}]
		
		Set $\mathscr{A} = \mathscr{A}_n (\delta) = \mathscr{A}^{(1)} \cap \mathscr{A}^{(2)} \cap \mathscr{A}^{(3)}$. By \Cref{a123}, \Cref{a4}, and a union bound, we have $\mathbb{P} [\mathscr{A}] \ge 1 - \delta$ for sufficiently large $n$. Since each of the $\mathscr{A}^{(j)}$ are measurable with respect to $\mathcal{F}_{\ext}$ by their definitions, it suffices to verify that the three properties listed in \Cref{fhl} hold on $\mathscr{A}$. To confirm that the first does, observe from the fact that $\mathscr{A} \subseteq \mathscr{A}^{(1)}$, \eqref{a1}, and the definition of the event $\textbf{GAP}$ from \Cref{gap}, that on $\mathscr{A}$ we have
		\begin{flalign*}
			\mathcal{L}_1 (\mathsf{R}) - \mathcal{L}_n (\mathsf{R}) \le \mathfrak{C}_2 \big( (n+1)^{2/3} + (\log n^{30})^{25} \big) \le n,
		\end{flalign*}					
		
		\noindent for sufficiently large $n$. That the second does follows from \Cref{functiong}, \eqref{gfunction0}, and \eqref{a2}, which define $\textbf{GBL}$, $\mathfrak{G}$, and $\mathscr{A}^{(2)}$, respectively (and the fact that $\delta^2 \le \delta/2$). To confirm that the third does, observe that the first two bounds in \eqref{th} hold by the definition of $\mathscr{A}^{(3)}$; the third bound in \eqref{th} holds by the last part of the definition \eqref{a2} of $\mathscr{A}^{(2)}$ together with the first bound in \eqref{th}. This establishes the lemma.
	\end{proof} 
	
	Let us now establish \Cref{a123}. 
	
	\begin{proof}[Proof of \Cref{a123}] 
		
		By a union bound, it suffices to show that 
		\begin{flalign} 
			\label{a1a2a3}
			\mathbb{P} \big[ \mathscr{A}^{(1)} \big] \ge 1 - \displaystyle\frac{\delta}{4}; \qquad \mathbb{P} \big[ \mathscr{A}^{(2)} \big] \ge 1 - \displaystyle\frac{\delta}{4}. 
		\end{flalign} 
		
		\noindent By \Cref{x1lsmall} (with the $(n, B, \vartheta, \delta)$ there equal to $( n^{30}, \mathfrak{C}_2, \mathfrak{C}_1^{-1}, \delta/8)$ here), we have  
		\begin{flalign}
			\label{1c2}
			\begin{aligned} 
			\mathbb{P} \Big[ \textbf{TOP}^{\bm{\mathcal{L}}}  \big( [ & -\mathfrak{C}_2 n^{10}, \mathfrak{C}_2 n^{10}]; \mathfrak{C}_1^{-1} n^{20} \big) \Big] \ge 1 - \displaystyle\frac{\delta}{8},
			\end{aligned} 
		\end{flalign}			
	
		\noindent for sufficiently large $n$. Hence, 
		\begin{flalign}
			\label{a11} 
			\mathbb{P} \big[ \mathscr{A}^{(1)} \big] \ge \mathbb{P} \big[ \textbf{SCL}_{n^{30}}^{\bm{\mathcal{L}}} (2; 2; \mathfrak{C}_2) \big] \ge 1 - \displaystyle\frac{\delta}{8} - \mathfrak{c}_1^{-1} e^{-\mathfrak{c}_1 (\log n)^2}.
		\end{flalign}
	
		\noindent Here, in the first bound we used the fact that $\textbf{SCL}_{n^{30}}^{\bm{\mathcal{L}}} (2; 2; \mathfrak{C}_2) \subseteq \textbf{GAP}_{n^{30}}^{\bm{\mathcal{L}}} \big( [-2n^{10}, 2n^{10}]; \mathfrak{C}_2 \big) \subseteq \textbf{GAP}_{n^{30}}^{\bm{\mathcal{L}}} (n^{20}; \mathfrak{C}_2) = \mathscr{A}^{(1)}$ (by \Cref{eventscl}, \eqref{a1}, and the fact that $\mathsf{R} = n^{20} \in [-2n^{20}, 2n^{20}]$); in the second, we applied \Cref{sclprobability} (with the $n$ there equal to $n^{30}$ here), \eqref{1c2}, and a union bound. The estimate \eqref{a11} then gives the first bound in \eqref{a1a2a3}, for sufficiently large $n$.
		
		To establish the second, observe by \eqref{a2} that it suffices to show (since $\delta \le 1/4$)  
		\begin{flalign}
			\label{a12a22}
			\mathbb{P} \big[ \mathscr{A}_1^{(2)} \big] \ge 1 - \delta^3; \qquad  \mathbb{P} \big[ \mathscr{A}_2^{(2)} \big] \ge 1 - \delta^4.
		\end{flalign} 
		
		\noindent To that end, first observe by \Cref{p:globallaw2} that, for sufficiently large $n$, we have 
	\begin{flalign*}
		\mathbb{P} \big[ \textbf{GBL}_{n+1}^{\bm{\mathcal{L}}} (\delta^4; \delta^{-1}) \big] \ge 1 - \delta^4.
	\end{flalign*} 
	
		\noindent Thus, the second bound in \eqref{a12a22} holds since $\textbf{GBL}_{n+1} (\delta^4; \delta^{-1}) \subseteq \mathscr{A}_2^{(2)}$ (by \eqref{a2} and the facts $(n+1)^{2/3} - n^{2/3} \le n^{-1/3} < \delta n^{2/3} / 60$ and $\delta^4 (n+1)^{2/3} < \delta n^{2/3} / 6$, which hold for sufficiently large $n$, as $\delta \le 1/2$). To show the first bound in \eqref{a12a22}, observe that 
		\begin{flalign*} 
			 \mathbb{E} \Big[ \mathbb{P} \big[ \textbf{GBL}_n^{\bm{\mathcal{L}}} (\delta^2; \delta^{-1})^{\complement} \big| \mathcal{F}_{\ext} \big] \Big] = \mathbb{P} \big[ \textbf{GBL}_n^{\bm{\mathcal{L}}} (\delta^2; \delta^{-1})^{\complement} \big] \le \mathbb{P}\big[ \textbf{GBL}_n^{\bm{\mathcal{L}}} (\delta^4; \delta^{-1})^{\complement} \big] \le \delta^4,
		\end{flalign*} 
	
	\noindent which implies by a Markov estimate that 
	\begin{flalign*}
		\mathbb{P} \Big[ \big(\mathscr{A}_1^{(2)} \big)^{\complement} \Big] = \mathbb{P} \Big[ \mathbb{P} \big[ \textbf{GBL}_n^{\bm{\mathcal{L}}} (\delta^2; \delta^{-1})^{\complement} \big| \mathcal{F}_{\ext} \big] > \delta \Big] \le \delta^3.
	\end{flalign*}

	\noindent This proves the first statement of \eqref{a12a22} and thus the lemma. 
	\end{proof}

	\subsection{Proof of \Cref{a4}} 
	
	\label{Proof4a}
	
	In this section we establish \Cref{a4}. This will follow from \Cref{h0x2}, after conditioning on an event on which the hypotheses in \Cref{derivativegxpfl} hold.

	To define the latter event, set $n_0 = \lfloor 2n/3 \rfloor$ and denote the ensemble $\bm{l} = \bm{l}^{(n_0)} = (l_1, l_2, \ldots , l_{n_0}) \in \llbracket 1, n_0 \rrbracket \times \mathcal{C} (\mathbb{R})$ by setting $l_j (t) = n_0^{-2/3} \cdot \mathcal{L}_{j+n_0} (tn_0^{1/3})$ for each $(j, t) \in \llbracket 1, n_0 \rrbracket \times \mathbb{R}$. Further let $\mathfrak{C}_3 > 1$ denote the constant $C(4\delta^{-1})$ from \Cref{p:closerho0}. We then define the event $\mathscr{A}^{(4)}$, measurable with respect to $\mathcal{G}_{\ext} = \mathcal{F}_{\ext} \big( \llbracket n_0 + 1, 2n_0 \rrbracket \times (-4\delta^{-1} n_0^{1/3}, 4\delta^{-1} n_0^{1/3}) \big)$, by setting $\mathscr{A}^{(4)} = \mathscr{A}_1^{(4)} \cap \mathscr{A}_2^{(4)}$, where
	\begin{flalign}
		\label{a4122} 
		\begin{aligned} 
		\mathscr{A}_1^{(4)} & = \Bigg\{ \mathbb{P} \bigg[ \bigcap_{|\delta t| \le 4} \textbf{PFL}^{\bm{l}} \big(t; n_0^{-1} (\log n_0)^7; \mathfrak{C}_3 \big) \Bigg| \mathcal{G}_{\ext} \bigg] \ge 1 - n_0^{-20} \Bigg\}; \\
		\mathscr{A}_2^{(4)} & = \bigcap_{\delta t \in \{ -4n_0^{1/3}, 4n_0^{1/3} \}} \bigcap_{j=1}^{n_0} \bigg\{ \big| \mathcal{L}_{j+n_0} (t) - n_0^{2/3} \cdot \mathfrak{G} (tn_0^{-1/3}, jn_0^{-1} + 1) \big| \le \delta^{20} n_0^{2/3} \bigg\} \\
		& \qquad \cap \bigcap_{|\delta t| \le 4 n_0^{1/3}} \bigg( \Big\{ \big| \mathcal{L}_{n_0} (t) - n_0^{2/3} \cdot \mathfrak{G} (tn_0^{-1/3}, 1) \big| \le \delta^{20} n_0^{2/3} \Big\} \\
		& \qquad \qquad \qquad \qquad \qquad \cap \Big\{ \big| \mathcal{L}_{2n_0+1} (t) - n_0^{2/3} \cdot \mathfrak{G} (tn_0^{-1/3}, 2)  \big| \le \delta^{20} n_0^{2/3} \Big\} \bigg).
		\end{aligned} 
	\end{flalign} 
	
	\noindent  recalling the function $\mathfrak{G}$ from \eqref{gfunction0}. The following lemma states that $\mathscr{A}^{(4)}$ is likely.
	\begin{lem} 
		
		\label{a5probability} 
		
		For sufficiently large $n$, we have $\mathbb{P} \big[ \mathscr{A}^{(4)} \big] \ge 1 - \delta/4$. 
	\end{lem} 

	\begin{proof} 
		
		By a union bound, it suffices to show that 
		\begin{flalign} 
		\label{a412} 
		\mathbb{P} \big[ \mathscr{A}_1^{(4)} \big] \ge 1 - \displaystyle\frac{\delta}{8}; \qquad \mathbb{P} \big[ \mathscr{A}_2^{(4)} \big] \ge 1 - \displaystyle\frac{\delta}{8}.
		\end{flalign} 
		
		\noindent To that end, first let $\mathfrak{c}_2 \in (0, 1)$ and $\omega > 0$ denote the constants $c(4\delta^{-1})$ and $\omega$ from \Cref{p:closerho0}, and recall that $\mathfrak{C}_3 > 1$ denotes the constant $C(4\delta^{-1})$ from that theorem. Define the event $\mathscr{A}_1^{(5)}$, which is measurable with respect to $\mathcal{G}_{\ext} = \mathcal{F}_{\ext} \big( \llbracket n_0+1, 2n_0 \rrbracket \times (-4\delta^{-1} n_0^{1/3}, 4\delta^{-1} n_0^{1/3}) \big)$, by
		\begin{flalign*}
		\mathscr{A}_1^{(5)} & = \bigcap_{k=1}^{\lfloor \omega^{-2} \rfloor} \mathscr{A}_1^{(5)} (k), \qquad \text{where} \qquad \mathscr{A}_1^{(5)} (k) = \textbf{TOP}^{\bm{\mathcal{L}}} \big( [-\mathfrak{C}_3 n_0^{k \omega/3}, \mathfrak{C}_3 n_0^{k\omega/3}]; \mathfrak{c}_2 n_0^{2k\omega/3} \big),
		\end{flalign*} 
	
		\noindent for any integer $k \ge 1$. By \Cref{x1lsmall} (with the $(n, B, \vartheta, \delta)$ there equal to $( n_0^{k\omega}, \mathfrak{C}_3, \mathfrak{c}_2, \omega^2 \delta/16)$ here), we have for sufficiently large $n$ that
		\begin{flalign*}
			\mathbb{P} \big[ \mathscr{A}_1^{(5)} (k) \big] \ge 1 - \displaystyle\frac{\omega^2 \delta}{16},
		\end{flalign*}
		
		\noindent and so a union bound yields 
		\begin{flalign}
			\label{a15}
			\mathbb{P} \big[ \mathscr{A}_1^{(5)} \big] \ge 1 - \displaystyle\sum_{k=1}^{\lfloor \omega^{-2} \rfloor} \Big( 1 - \mathbb{P} \big[ \mathscr{A}^{(5)} (k) \big] \Big) \ge 1 - \displaystyle\frac{\delta}{16}.
		\end{flalign}
		
		\noindent Now observe that \Cref{p:closerho0} implies the estimate
		\begin{flalign}
			\label{ea15} 
			\mathbb{P} \big[ \mathscr{E} \big| \mathscr{A}_1^{(5)} \big] \ge 1 - 2 \mathfrak{C}_3 n_0^{-50} \ge 1 - n_0^{-40}, \quad \text{where} \quad \mathscr{E} = \bigcap_{|t| \le 4 \delta^{-1} n_0^{1/3}} \textbf{PFL}^{\bm{l}} \big( t; n_0^{-1} (\log n_0)^7; \mathfrak{C}_3 \big),
		\end{flalign}
	
		\noindent and on the left side we conditioned on $\mathscr{A}_1^{(5)}$. Therefore, 
		\begin{flalign*}
			\mathbb{P} \Big[ \mathscr{A}_1^{(5)} \cap \big( \mathscr{A}_1^{(4)} \big)^{\complement} \Big] & = \mathbb{P} \bigg[ \mathscr{A}_1^{(5)} \cap \Big\{ \mathbb{P} \big[ \mathscr{E}^{\complement} \big| \mathcal{G}_{\ext} \big] \ge n_0^{-20} \Big\} \bigg] \\
			 & = \mathbb{P} \Big[ \mathbb{P} \big[ \mathscr{E}^{\complement} \cap \mathscr{A}_1^{(5)} \big| \mathcal{G}_{\ext} \big] \ge n_0^{-20} \Big] \le n_0^{20} \cdot \mathbb{E} \big[ \mathscr{E}^{\complement} \cap \mathscr{A}_1^{(5)} \big] \le n_0^{-20},
		\end{flalign*}
		
		\noindent where the first statement holds by the definition \eqref{a4122} of $\mathscr{A}_1^{(4)}$; the second since $\mathscr{A}_1^{(5)}$ is measurable with respect to $\mathcal{G}_{\ext}$; the third by a Markov estimate applied to the random variable $\mathbb{P} \big[ \mathscr{E}^{\complement} \cap \mathscr{A}_1^{(5)} \big]$; and the fourth by \eqref{ea15}. This, with \eqref{a15} and a union bound, together yield  
		\begin{flalign*} 
			\mathbb{P} \big[ \mathscr{A}_1^{(4)} \big] \ge \mathbb{P} \big[ \mathscr{A}_1^{(5)} \big] - n_0^{-20} \ge 1 - \displaystyle\frac{\delta}{16}  - n_0^{-20} \ge 1 - \displaystyle\frac{\delta}{8}, 
		\end{flalign*} 
	
		\noindent for sufficiently large $n$, which establishes the first bound in \eqref{a412}.
		
		To establish the second, we observe from \Cref{functiong} that $ \textbf{GBL}_{n_0}^{\bm{\mathcal{L}}} (\delta^{25}; 4\delta^{-1})\subseteq \mathscr{A}_2^{(4)} $ (where the containment holds for the last event defining $\mathscr{A}_2^{(4)}$ in \eqref{a4122}, since $\delta^{25} + \big| \mathfrak{G} (tn_0^{-1/3}, 2 + n_0^{-1}) - \mathfrak{G} (tn_0^{-1/3}, 2) \big| \le \delta^{25} + 3n_0^{-1} \le \delta^{20}$ for sufficiently large $n$). Hence, \Cref{p:globallaw2} implies that 
		\begin{flalign*} 
			\mathbb{P} \big[ \mathscr{A}_2^{(4)} \big] \ge \mathbb{P} \big[ \textbf{GBL}_{n_0}^{\bm{\mathcal{L}}} (\delta^{25}; 4\delta^{-1}) \big] \ge 1 - \delta^{25} \ge 1 - \displaystyle\frac{\delta}{8}, 
			\end{flalign*} 
		
		\noindent for sufficiently large $n$, verifying the second statement of \eqref{a412} and thus the lemma.	
	\end{proof} 

	Now we can establish \Cref{a4} using \Cref{h0x2}.

	\begin{proof}[Proof of \Cref{a4}] 
		
	Condition on $\mathcal{F}_{\ext} \big( \llbracket n_0 + 1, 2n_0 \rrbracket \times (-4\delta^{-1} n_0^{1/3}, 4\delta^{-1} n_0^{1/3}) \big)$ and restrict to the event $\mathscr{A}^{(4)}$. Setting $\mathfrak{R} = [-4\delta^{-1}, 4\delta^{-1}] \times [0, 1]$, fix $B' = B' (\delta)= \| G \|_{\mathcal{C}^{50} (\mathfrak{R})} + 4\delta^{-1}$. We then apply \Cref{h0x2}, with the $(n, B, \xi, \delta)$ there equal to $(n_0, B' + \mathfrak{C}_3, 4\delta^{-1}, \delta^{20})$ here,; the $(\bm{\mathsf{x}}; f, g)$ there equal to $(\bm{\mathcal{L}}_{\llbracket n_0+1, 2n_0 \rrbracket}; \mathcal{L}_{n_0}, \mathcal{L}_{2n_0+1})$ (restricted to the interval $[-4\delta^{-1} n_0^{1/3}, 4\delta^{-1} n_0^{1/3}]$) here; and the $\bm{u}$ and $\bm{v}$ there equal to 
	\begin{flalign*} 
		& \big( \mathcal{L}_{n_0+1} (-4\delta^{-1} n_0^{1/3}), \mathcal{L}_{n_0+2} (-4\delta^{-1} n_0^{1/3}), \ldots , \mathcal{L}_{2n_0} (-4\delta^{-1} n_0^{1/3}) \big); \\ 
		& \big( \mathcal{L}_{n_0 + 1} (4\delta^{-1} n_0^{1/3}), \mathcal{L}_{n_0 + 2} (4\delta^{-1} n_0^{1/3}), \ldots , \mathcal{L}_{2n_0} (4\delta^{-1} n_0^{1/3}) \big),
	\end{flalign*}

	\noindent here, respectively. To verify \Cref{derivativegxpfl}, observe that \eqref{g50} is satisfied since $B > B'$; \eqref{pflx2} is confirmed by the definition \eqref{a4122} of $\mathscr{A}_1^{(4)}$ (with the bound $n_0^{-1} (\log n_0)^{20} \le n_0^{-19/20}$); \eqref{gnuv} is confirmed by the definition \eqref{a4122} of $\mathscr{A}_2^{(4)}$.
	
	Hence, \Cref{h0x2} applies and yields a constant $C_1 = C_1 (\delta) > 1$ such that the following holds with probability at least $1 - C_1 n_0^{-10}$. There exists a (random) twice-differentiable function $\widetilde{h} : [-2\delta^{-1}, 2\delta^{-1}] \rightarrow \mathbb{R}$ such that
	\begin{flalign*}
		\displaystyle\sup_{|s| \le 2 \delta^{-1}} \big| \partial_s^2 \widetilde{h} (s) + 2^{-1/2} \big| \le \delta^2 + (\log n_0)^{-1/4}; \qquad \displaystyle\sup_{|s| \le 2\delta^{-1} n_0^{1/3}} \big| \mathcal{L}_{n+1} (s) - n_0^{2/3} \cdot \widetilde{h} (n_0^{-1/3} s) \big| \le n_0^{-1/5}.
	\end{flalign*} 

	\noindent Since $\delta^{-1} n^{1/3} \le 2 \delta^{-1} n_0^{1/3}$, defining $h: [-\delta^{-1} n^{1/3}, \delta^{-1} n^{1/3}] \rightarrow \mathbb{R}$ by setting $h(t) = n_0^{2/3} \cdot \widetilde{h} (n_0^{-1/3} t)$ for each $|t| \le \delta^{-1} n^{1/3}$, it follows that 
	\begin{flalign*}
		 \displaystyle\sup_{|t| \le \delta^{-1} n^{1/3}} \big| \partial_t^2 h(t) + 2^{-1/2} \big| \le \delta^2 + (\log n_0)^{-1/4} \le \delta; \quad \displaystyle\sup_{|t| \le \delta^{-1} n^{1/3}} \big| \mathcal{L}_{n+1} (t) - h(t) \big| \le n_0^{-1/5} < n^{-1/6},
	\end{flalign*} 

	\noindent for $n$ sufficiently large. In particular, $h$ satisfies the first two bounds in \eqref{th}, so $\mathscr{A}^{(3)}$ holds. Hence, 
	\begin{flalign*} 
		\mathbb{P} \big[ \mathscr{A}^{(3)} \big] \ge \mathbb{P} \big[ \mathscr{A}^{(4)} \big] - C_1 n_0^{-10} \ge 1 - \frac{\delta}{4} - C_1 n_0^{-10} \ge 1 - \displaystyle\frac{\delta}{2}, 
	\end{flalign*} 
	
	\noindent where in the second bound we applied \Cref{a5probability} and in the third we used that $n$ is sufficiently large. This establishes the lemma. 
	\end{proof}

		\subsection{Proofs of \Cref{constantk1} and \Cref{linvariant}}
	
	\label{Proofsql1}
	
	In this section we establish first \Cref{constantk1} and then \Cref{linvariant}, both of which are quick consequences of \Cref{lr}.
	
	\begin{proof}[Proof of \Cref{constantk1}]
		
		First observe by \eqref{sigmar} and \Cref{sigmascale} that we may assume that $\sigma = 2^{-1/2}$, so $q = 1$. Then, the first part of the corollary follows from \Cref{lr}. 
		
		To establish the second part, observe that we may also assume that $\ell = 0$, by subtracting $\ell t$ from the curves $\mathcal{L}_j (t)$ of $\bm{\mathcal{L}}$ and using the fact that that this affine transformation does not affect the Brownian Gibbs property (as it neither affects the law of a Brownian bridge nor the non-intersection property between curves; see \Cref{linear} for more details). Since \eqref{sigmalinear} implies \eqref{sigmat2} (taking the constant $C$ in the latter sufficiently large in comparison to that from the former), as $\sigma \varepsilon t^2 > \varepsilon |t|$ for $t$ sufficienly large, there exist random variables $\mathfrak{l}, \mathfrak{c} \in \mathbb{R}$ and a scaled parabolic Airy line ensemble $\bm{\mathcal{S}} = (\mathcal{S}_1, \mathcal{S}_2, \ldots )$ independent from $\mathfrak{l}$ and $\mathfrak{c}$, such that $\mathcal{L}_j (t) = \mathcal{S}_j (t) + \mathfrak{l} t + \mathfrak{c}$; we must show that $\mathfrak{l} = 0$ deterministically. 
		
		Fix a real number $\ell' > 0$. The translation-invariance of $\bm{\mathcal{S}} + 2^{-1/2} \cdot t^2$ implies for any $t > 0$ that
		\begin{flalign*}
			\mathbb{P} \big[ \mathcal{L}_1 (t) \le -2^{-1/2} t^2 + \ell' t \big] & = \mathbb{P} \big[ \mathcal{S}_1 (t) \le -2^{-1/2} t^2 + t (\ell' - \mathfrak{l}) - \mathfrak{c} \big]  = \mathbb{P} \big[ \mathcal{S}_1 (0) \le t (\ell' - \mathfrak{l}) - \mathfrak{c} \big].
		\end{flalign*} 
		
		\noindent Since $\ell' > 0$, the left side of this equality tends to $1$ as $t$ tends to $\infty$, by \eqref{sigmalinear}. Thus, 
		\begin{flalign*}
			\displaystyle\lim_{t \rightarrow \infty} \mathbb{P} \big[ \mathcal{S}_1 (0) \le t ( \displaystyle 2\ell' - \mathfrak{l} ) \big] \geq \displaystyle\lim_{t \rightarrow \infty} \mathbb{P} \big[ \mathcal{S}_1 (0) \le  t (\ell' - \mathfrak{l}) + \mathfrak{c} \big] = 1,
		\end{flalign*}
		
		\noindent where we used the fact that $\lim_{t \rightarrow \infty} \mathbb{P} [\ell' t >  \mathfrak{c}] = 1$. It follows that $\mathbb{P} [\mathfrak{l} \le 2\ell'] = 1$ for any $\ell' > 0$, so that $\mathfrak{l} \le 0$ almost surely. The proof that $\mathfrak{l} \ge 0$ is entirely analogous (by letting $t$ tend to $-\infty$, instead of to $\infty$, in the above). This shows $\mathfrak{l} = 0$, establishing the second statement of the corollary.
	\end{proof}

	\begin{proof}[Proof of \Cref{linvariant}] 
		
		Since $\mu \in \Tra (2^{-1/2})$, the ensemble $\bm{\mathcal{L}}$ satisfies \eqref{sigmalinear} at $(\sigma, q, \ell) = (2^{-1/2}, 1, 0)$. It follows that there exists a random variable $\mathfrak{c} \in \mathbb{R}$ and an independent scaled parabolic Airy line ensemble $\bm{\mathcal{S}} = (\mathcal{S}_1, \mathcal{S}_2, \ldots )$ such that $\mathcal{L}_j (t) = \mathcal{S}_j (t) + \mathfrak{c}$, for each $(j, t) \in \mathbb{Z}_{\ge 1} \times \mathbb{R}$. Since $\mu$ is extremal, this implies $\mathfrak{c}$ is some (deterministic) constant $c$, which establishes the corollary.
	\end{proof}

	\section{Miscellaneous Preliminaries} 
	
	\label{Estimates}
	
	In this section we collect various facts about non-intersecting Brownian bridges, free convolutions, and Dyson Brownian motion. These results are (essentially) known, though for completeness we include the proofs of those that we did not directly find in the literature in the appendix, \Cref{ProofBridgeSum}, below.
	
	\subsection{Strong Gibbs Property and Invariances}
	
	\label{PropertyPaths}
	
	In this section we review a more restrictive variant of the Brownian Gibbs property (referred to as the strong Brownian Gibbs property) and several transformations that leave non-intersecting Brownian bridge measures invariant; we begin with the former.
	
	\begin{definition} 
		
	\label{property2} 
	
	Fix subsets $I \subseteq  \mathbb{R}$ and $\Sigma \subseteq \mathbb{Z}_{\ge 1}$, and a $\Sigma \times I$ indexed line ensemble $\bm{\mathsf{x}} = (\mathsf{x}_s)_{s \in \Sigma} \in \Sigma \times \mathcal{C} (I)$. For any finite interval $\Sigma' \subseteq \Sigma$, a random variable $(\mathfrak{a}, \mathfrak{b}) \in I^2$ is called a \emph{$\Sigma'$-stopping domain} if, for any $a, b \in I$ with $a < b$, we have 
	\begin{flalign*} 
		\{ \mathfrak{a} < a, \mathfrak{b} > b \} \in \mathcal{F}_{\ext} \big( \Sigma' \times (a, b) \big).
	\end{flalign*} 
	
	\noindent For any $\Sigma'$-stopping domain $(\mathfrak{a}, \mathfrak{b})$, let $\mathcal{F}_{\ext} \big( \Sigma' \times (\mathfrak{a}, \mathfrak{b}) \big)$ denote the $\sigma$-algebra generated by all events $\mathscr{E}$ such that $\mathscr{E} \cap \{ \mathfrak{a} < a, \mathfrak{b} > b \} \in \mathcal{F}_{\ext} \big( \Sigma' \times (a, b) \big)$ holds for any real numbers $a < b$. 
	
	Let $\mathcal{C}^{\Sigma'}$ denote the set of $\big( |\Sigma'| + 2\big)$-tuples $(a, b, f_j)_{j \in \Sigma'}$ such that $a, b \in I$ with $a < b$ and $f_j \in \mathcal{C} \big( [a, b] \big)$ for each $j \in \Sigma'$. We consider $\mathcal{C}^{\Sigma'}$ under the topology with an open basis given by 
	\begin{flalign*}
	\bigg\{ (a, b, f_j)_{j \in \Sigma'} \in \mathcal{C}^{\Sigma'} : a < b, |a-A| + |b-B| < \varepsilon, \displaystyle\max_{j \in \Sigma'} \displaystyle\max_{a \vee A \le x \le b \wedge B} \big| f_j (x) - F_j (x) \big| < \varepsilon  \bigg\},
	\end{flalign*} 
	
	\noindent over $\varepsilon \in (0, 1)$ and $(A, B, F_j)_{j \in \Sigma} \in \mathcal{C}^{\Sigma'}$.
	
	A $\Sigma \times I$ indexed line ensemble $\bm{\mathsf{x}} = (\mathsf{x}_j)_{j \in \Sigma}$ is said to satisfy the \emph{strong Brownian Gibbs property} if, for any interval $\llbracket k_1, k_2 \rrbracket \subseteq \Sigma$; Borel measurable function $F : \mathcal{C}^{\Sigma'} \rightarrow \mathbb{R}$; and $\llbracket k_1, k_2 \rrbracket$-stopping domain $(\mathfrak{a}, \mathfrak{b})$, we have 
	\begin{flalign*}
		\mathbb{E} \Big[ F \big( \mathfrak{a}, \mathfrak{b}, \mathsf{x}_{k_1} |_{[\mathfrak{a}, \mathfrak{b}]}, \mathsf{x}_{k_1+1} |_{[\mathfrak{a}, \mathfrak{b}]}, \ldots , \mathsf{x}_{k_2} |_{[\mathfrak{a}, \mathfrak{b}]} & \big) \cdot \textbf{1}_{\mathfrak{a} < \mathfrak{b}} \Big| \mathcal{F}_{\ext} \big( \llbracket k_1, k_2 \rrbracket \times (\mathfrak{a}, \mathfrak{b}) \big) \Big] \\
		&  = \mathbb{E} \Big[ F (\mathfrak{a}, \mathfrak{b}, \mathsf{y}_{k_1}, \mathsf{y}_{k_1+1}, \ldots , \mathsf{y}_{k_2} ) \cdot \textbf{1}_{\mathfrak{a} < \mathfrak{b}} ]\Big],
	\end{flalign*}
	
	\noindent where the expectation on the right side is with respect to both $(\mathfrak{a}, \mathfrak{b})$ and non-intersecting Brownian bridges $\bm{\mathsf{y}} = (\mathsf{y}_{k_1}, \mathsf{y}_{k_1+1}, \ldots , \mathsf{y}_{k_2}) \in \llbracket k_1, k_2 \rrbracket \times \mathcal{C} \big( [\mathfrak{a}, \mathfrak{b}] \big)$ sampled according to the measure $\mathsf{Q}_{\mathsf{x}_{k_2+1}; \mathsf{x}_{k_1-1}}^{\bm{u}; \bm{v}}$ (setting $\mathsf{x}_{k_1-1} = \infty$ if $k_1 - 1 < \min \Sigma$ and $\mathsf{x}_{k_2+1} = -\infty$ if $k_2+1 > \max \Sigma$), whose entrance data is given by $\bm{u} = \big( \mathsf{x}_{k_1} (\mathfrak{a}), \mathsf{x}_{k_1+1} (\mathfrak{a}), \ldots , \mathsf{x}_{k_2} (\mathfrak{a}) \big)$ and exit data is given by $\bm{v} = \big(\mathsf{x}_{k_1} (\mathfrak{b}), \mathsf{x}_{k_1+1} (\mathfrak{b}), \ldots , \mathsf{x}_{k_2} (\mathfrak{b}) \big)$.  
	
	\end{definition}

	The following lemma indicates that line ensembles satisfying the Brownian Gibbs property also satisfy its above strong variant.
	
	\begin{lem}[{\cite[Lemma 2.5]{PLE}}]	
		
		\label{propertyproperty2}
		
		Fix intervals $\Sigma \subseteq \mathbb{Z}_{\ge 1}$ and $I \subseteq \mathbb{R}$. Any $\Sigma \times I$ indexed line ensemble satisfying the Brownian Gibbs property also satisfies the strong Brownian Gibbs property.
		
	\end{lem}

	Next we observe two invariance properties satisfied by non-intersecting Brownian bridges; the first is under affine transformations, and the second is under diffusive scaling.

	\begin{rem}
		
		\label{linear}
		
		Non-intersecting Brownian bridges satisfy the following invariance property under affine transformations. Adopt the notation of \Cref{qxyfg}, and fix real numbers $\alpha, \beta \in \mathbb{R}$. Define the $n$-tuples $\widetilde{\bm{u}}, \widetilde{\bm{v}} \in \mathbb{W}_n$ and functions $\widetilde{f}, \widetilde{g} : [a, b] \rightarrow \overline{\mathbb{R}}$ by setting 
		\begin{flalign*} 
			& \widetilde{u}_j = u_j + \alpha, \quad \text{and} \quad \widetilde{v}_j = v_j + (b-a) \beta + \alpha, \qquad \qquad \quad \text{for each $j \in \llbracket 0, n \rrbracket$}; \\
			& \widetilde{f}(t) = f(t) + t \beta + \alpha, \quad \text{and} \quad \widetilde{g}(t) = g(t) + t \beta + \alpha, \qquad \text{for each $t \in [a, b]$}.
		\end{flalign*} 
		
		\noindent Sampling $\widetilde{\bm{\mathsf{x}}} = (\widetilde{\mathsf{x}}_1, \widetilde{\mathsf{x}}_2, \ldots , \widetilde{\mathsf{x}}_n )$ under $\mathsf{Q}_{\tilde{f}, \tilde{g}}^{\tilde{\bm{u}}; \tilde{\bm{v}}}$, there is a coupling between $\widetilde{\bm{\mathsf{x}}}$ and $\bm{\mathsf{x}}$ such that $\widetilde{\mathsf{x}}_j (t) = \mathsf{x}_j (t) + (t-a) \beta + \alpha$ for each $t \in [a, b]$ and $j \in \llbracket 1, n \rrbracket$. 
		
		Indeed, this follows from the analogous affine invariance of a single Brownian bridge, together with the fact that affine transformations do not affect the non-intersecting property. More specifically, if $\big( \mathsf{x} (t) \big)$, for $t \in [a, b]$, is a Brownian bridge from some $u \in \mathbb{R}$ to some $v \in \mathbb{R}$ then $\big( x(t) + (t-a) \beta + \alpha \big)$ is a Brownian bridge from $u + \alpha$ to $v + (b-a) \beta + \alpha$, and any $\bm{\mathsf{y}} (t) \in \mathbb{W}_n$ (is non-intersecting) if and only if $\bm{\mathsf{y}}(t) + (t-a) \beta + \alpha \in \mathbb{W}_n$.
		
	\end{rem}

	\begin{rem}
		
		\label{scale}
		
		Non-intersecting Brownian bridges also satisfy the following invariance property under diffusive scaling. Again adopt the notation of \Cref{qxyfg}; assume that $(a, b) = (0, T)$, for some real number $T > 0$. Further fix a real number $\sigma > 0$, and set $\widetilde{T} = \sigma T$. Define the $n$-tuples $\widetilde{\bm{u}}, \widetilde{\bm{v}} \in \mathbb{W}_n$ and functions $\widetilde{f}, \widetilde{g} : [0, \widetilde{T}] \rightarrow \overline{\mathbb{R}}$ by setting 
		\begin{flalign*} 
			& \widetilde{u}_j = \sigma^{1/2} u_j, \quad \text{and} \quad \widetilde{v}_j = \sigma^{1/2} v_j, \qquad \qquad \qquad \qquad \qquad \text{for each $j \in \llbracket 0, n \rrbracket$}; \\
			& \widetilde{f} (t) = \sigma^{1/2} \cdot f(\sigma^{-1} t), \quad \text{and} \quad \widetilde{g} (t) = \sigma^{1/2} \cdot g(\sigma^{-1} t), \qquad \text{for each $t \in \big[ 0, \widetilde{T} \big]$}.
		\end{flalign*} 
		
		\noindent Sampling $\widetilde{\bm{\mathsf{x}}} = ( \widetilde{\mathsf{x}}_1, \widetilde{\mathsf{x}}_2, \ldots , \widetilde{\mathsf{x}}_n)$ under $\mathsf{Q}_{\tilde{f}, \tilde{g}}^{\tilde{\bm{u}}; \tilde{\bm{v}}}$, there is a coupling between $\widetilde{\bm{\mathsf{x}}}$ and $\bm{\mathsf{x}}$ such that $\widetilde{\mathsf{x}}_j (t) = \sigma^{1/2} \cdot \mathsf{x}_j (\sigma^{-1} t)$ for each $(j, t) \in \llbracket 1, n \rrbracket \times [0, \widetilde{T}]$. Similarly to in \Cref{linear}, this follows from the analogous scaling invariance of a single Brownian bridge.
		
	\end{rem}
	
	The following lemma is a (known) bound for the maximum of a Brownian bridge.
	
	\begin{lem}[{\cite[Chapter 4, Equation (3.40)]{MSC}}]
		
		\label{maximumx} 
		
		Fix a real number $T > 0$. Let $\mathsf{x}: [0, T] \rightarrow \mathbb{R}$ denote a Brownian bridge conditioned to start and end at $\mathsf{x} (0) = 0 = \mathsf{x} (T)$. For any real number $u > 0$, we have
		\begin{flalign*}
			\mathbb{P} \bigg[ \displaystyle\sup_{t \in [0, T]} \big| \mathsf{x} (t) \big| \ge u \bigg] = 2 e^{-u^2/2T}.
		\end{flalign*}
		
	\end{lem}

	We conclude this section with a Markov estimate; its quick proof is given in \Cref{Proofydelta}.
	
	\begin{lem} 
		
		\label{fg0g}
		
		Let $\Sigma' \subseteq \Sigma \subseteq \mathbb{Z}_{\ge 1}$ and $I' \subseteq I \subseteq \mathbb{R}$ be intervals; $a, b \in (0, 1]$ be real numbers; and $\bm{\mathsf{x}}$ be a $(\Sigma, I)$ indexed line ensemble. Let $\mathscr{E}$ denote an event measurable with respect to the $\sigma$-algebra generated by $\big( \mathsf{x}_j (t) \big)$ for $(j, t) \in \Sigma' \times I'$, such that $\mathbb{P}[\mathscr{E}] \ge 1-ab$. Then there exists an event $\mathscr{G}$, measurable with respect to $\mathcal{F}_{\ext} = \mathcal{F}_{\ext} (\Sigma' \times I')$, such that the following two statements hold. 
		\begin{enumerate} 
			\item We have $\mathbb{P} [\mathscr{G}] \ge 1 - b$. 
			\item Conditioning on $\mathcal{F}_{\ext}$, we have $\mathbb{P} [ \mathscr{E} | \mathcal{F}_{\ext}] \ge (1 - a) \cdot \textbf{\emph{1}}_{\mathcal{G}}$.
		\end{enumerate} 
	
	\end{lem}

	\subsection{Height Monotonicity, Concentration Bounds, and H\"{o}lder Estimates}

	\label{Coupling2}
	
	In this section we state monotone couplings for non-intersecting Brownian Gibbsian line ensembles, as well as concentration bounds and H\"{o}lder estimates they satisfy. The following lemma recalls a monotone coupling for non-intersecting Brownian bridges; we refer to it as \emph{height monotonicity}. 
	
	It is essentially due to \cite[Lemmas 2.6 and 2.7]{PLE}, though the latter was only stated if the endpoints $\bm{u}, \bm{v}, \widetilde{\bm{u}}, \widetilde{\bm{v}} \in \mathbb{W}_n$ are strictly ordered. We provide the quick extension of this proof to when $\bm{u}, \bm{v}, \widetilde{\bm{u}}, \widetilde{\bm{v}} \in \overline{\mathbb{W}}_n$ in \Cref{Proofydelta} below. 
	
	\begin{lem}
		
		\label{monotoneheight}
		
		Fix an integer $n \ge 1$; four $n$-tuples $\bm{u}, \widetilde{\bm{u}}, \bm{v}, \widetilde{\bm{v}} \in \overline{\mathbb{W}}_n$; an interval $[a, b] \in \mathbb{R}$; and continuous functions $f, \widetilde{f}, g, \widetilde{g}: [a, b] \rightarrow \overline{\mathbb{R}}$. Sample two families of non-intersecting Brownian bridges $\bm{\mathsf{x}} = (\mathsf{x}_1, \mathsf{x}_2, \ldots , \mathsf{x}_n) \in \llbracket 1, n \rrbracket \times \mathcal{C} \big( [a, b] \big)$ and $\widetilde{\bm{\mathsf{x}}} = (\widetilde{\mathsf{x}}_1, \widetilde{\mathsf{x}}_2, \ldots , \widetilde{\mathsf{x}}_n) \in \llbracket 1, n \rrbracket \times \mathcal{C} \big( [a, b] \big)$ from the measures $\mathsf{Q}_{f; g}^{\bm{u}; \bm{v}}$ and $\mathsf{Q}_{\tilde{f}; \tilde{g}}^{\tilde{\bm{u}}; \tilde{\bm{v}}}$, respectively. If 
		\begin{flalign}
			\label{fgxmonotone} 
			f \le \widetilde{f} ; \quad  g \le \widetilde{g}; \quad \bm{u} \le \widetilde{\bm{u}}; \quad \bm{v} \le \widetilde{\bm{v}},
		\end{flalign}
		
		\noindent  then there exists a coupling between $\bm{\mathsf{x}}$ and $\widetilde{\bm{\mathsf{x}}}$ so that $\mathsf{x}_j (t) \le \widetilde{\mathsf{x}}_j (t)$, for each $(j, t) \in \llbracket 1, n \rrbracket \times [a, b]$.
		
	\end{lem}
	
	We next state the following variant of the above coupling, due to \cite{U}, whose second part provides a linear bound on the difference between two families of non-intersecting Brownian bridges, which have the same starting data but different ending data.  
	
	\begin{lem}[{\cite[Lemma 2.4 and Remark 2.5]{U}}]
		
		\label{uvv}
		
		Fix an integer $n \ge 1$; a real number $B \ge 0$; a finite interval $[a, b] \subset \mathbb{R}$; four $n$-tuples $\bm{u}, \widetilde{\bm{u}}, \bm{v}, \widetilde{\bm{v}} \in \mathbb{W}_n$; and four continuous functions $f, \widetilde{f}, g, \widetilde{g} : [a, b] \rightarrow \overline{\mathbb{R}}$. Assume that  $f(a) < u_n < u_1 < g(a)$ and $f(b) < v_n< v_1 < g(b)$, and that $\widetilde{f}(a) < \widetilde{u}_n < \widetilde{u}_1 < \widetilde{g}(a)$ and $\widetilde{f}(b) < \widetilde{v}_n < \widetilde{v}_1 < \widetilde{g}(b)$. Also assume that 
		\begin{flalign*} 
			\max_{j \in \llbracket 1, n \rrbracket} \big\{ |u_j - \widetilde{u_j}|,  |v_j - \widetilde{v}_j| \big\} \le B; \qquad \sup_{t \in [a, b]} \Big\{ \big| f(t) - \widetilde{f}(t) \big|, \big| g(t) - \widetilde{g}(t) \big| \Big\} \le B
		\end{flalign*}
		
		\noindent Sample two families of non-intersecting Brownian bridges $\bm{\mathsf{x}} = (\mathsf{x}_1, \mathsf{x}_2, \ldots , \mathsf{x}_n) \in \llbracket 1, n \rrbracket \times \mathcal{C} \big( [a, b] \big)$ and $\widetilde{\bm{\mathsf{x}}} = (\widetilde{\mathsf{x}}_1, \widetilde{\mathsf{x}}_2, \ldots , \widetilde{\mathsf{x}}_n) \in \llbracket 1, n \rrbracket \times \mathcal{C} \big( [a, b] \big)$ from the measures $\mathsf{Q}_{f; g}^{\bm{u}; \bm{v}}$ and $\mathsf{Q}_{\tilde{f}; \tilde{g}}^{\tilde{\bm{u}}; \tilde{\bm{v}}}$, respectively. 
		
		\begin{enumerate} 
			
			\item  There is a coupling between $\bm{\mathsf{x}}$ and $\widetilde{\bm{\mathsf{x}}}$ so that $\big|\widetilde{\mathsf{x}}_j (t)- \mathsf{x}_j (t)\big| \leq B$ for each $(j, t) \in \llbracket 1, n \rrbracket \times [a, b]$. 
			
			\item Further assume that $\bm{u} = \widetilde{\bm{u}}$ and for each $t \in [a, b]$ that
			\begin{flalign*} 
				\big| f(t) - \widetilde{f} (t) \big| \le \frac{t-a}{b-a} \cdot B; \qquad \big| g(t) - \widetilde{g}(t) \big| \le \frac{t-a}{b-a}\cdot B.
			\end{flalign*} 
			
			\noindent Then, it is possible to couple $\bm{\mathsf{x}}$ and $\widetilde{\bm{\mathsf{x}}}$ such that 
			\begin{flalign*} 
				\big| \mathsf{x}_j (t) - \widetilde{\mathsf{x}}_j (t) \big| \le \frac{t-a}{b-a} \cdot B, \qquad \text{for each $(j, t) \in \llbracket 1, n \rrbracket \times [a, b]$}.
			\end{flalign*}

		\end{enumerate} 
	\end{lem}

	 We next recall the following H\"{o}lder estimate from \cite{BPLE} for non-intersecting Brownian bridges.
	
	\begin{lem}[{\cite[Proposition 3.5]{BPLE}}]
		
		\label{estimatexj2}

		There exist constants $c > 0$ and $C > 1$ such that the following holds. Let $n \ge 1$ be an integer, $B \ge 1$ be a real number, $[a, b] \subset \mathbb{R}$ be an interval, and $\bm{u}, \bm{v} \in \mathbb{W}_n$ be two $n$-tuples; set $\mathsf{T} = b - a$. Sampling $\bm{\mathsf{x}} = (\mathsf{x}_1, \mathsf{x}_2, \ldots , \mathsf{x}_n) \in \llbracket 1, n \rrbracket \times \mathcal{C} \big( [a, b] \big)$ under $\mathsf{Q}^{\bm{u}; \bm{v}}$, we have  
		\begin{flalign*}
			\mathbb{P} \Bigg[ \bigcup_{j=1}^n \bigcup_{a \le t < t+s \le b} \Big\{ \big| \mathsf{x}_j (t+s) - \mathsf{x}_j (t) - s \mathsf{T}^{-1} (v_j - u_j) \big| > B s^{1/2} \log (2s^{-1} \mathsf{T}) \Big\} \Bigg] \le C e^{Cn - cB^2}.
		\end{flalign*} 
	\end{lem}

	 We will also require the following variant of \Cref{estimatexj2} from \cite{U} that allows for a lower boundary $f$. Without at least some continuity constraints on $f$, H\"{o}lder bounds for the paths in $\bm{\mathsf{x}}$ cannot hold everywhere on $[a, b]$ (for example, if $\mathsf{x}_n (a) = f(a)$ and $f(a^+) > f(a)$, then $\mathsf{x}_n$ will necessarily be discontinuous at the endpoint $a$). However, the next lemma provides a H\"{o}lder estimate on these paths (in the absence of a continuity condition on $f$) away the boundaries $t \in \{ a, b \}$ of $[a, b]$.

	\begin{lem}[{\cite[Lemma 2.7(1)]{U}}]
		
		\label{estimatexj3}
		
		There exist constants $c > 0$ and $C > 1$ such that the following holds. Let $n$, $B$, $a$, $b$, $\mathsf{T}$, $\bm{u}$, and $\bm{v}$ be as in \Cref{estimatexj2}. Further let $A \ge 1$ be a real number and $f: [a, b] \rightarrow \mathbb{R}$ be a continuous function such that $f(a) \le u_n$ and $f(b) \le v_n$. Assume that $f(r) - u_n \le A \mathsf{T}^{1/2}$ and $f(r) - v_n \le A \mathsf{T}^{1/2}$, for each $r \in [a, b]$. Sampling $\bm{\mathsf{x}} = (\mathsf{x}_1, \mathsf{x}_2, \ldots , \mathsf{x}_n) \in \llbracket 1, n \rrbracket \times \mathcal{C} \big( [a, b] \big)$ under $\mathsf{Q}_f^{\bm{u}; \bm{v}}$, we have for any real number $0 < \kappa < \min \{ \mathsf{T} / 2, 1 \}$ that  
		\begin{flalign*}
			\mathbb{P} \Bigg[ \bigcup_{j=1}^n \bigcup_{a + \kappa \le t < t+s \le b - \kappa} \bigg\{ \big| \mathsf{x}_j (& t+s) - \mathsf{x}_j (t) - s \mathsf{T}^{-1} (v_j - u_j) \big| \\
			& > s^{1/2}  \big( B \log |2s^{-1} \mathsf{T}| + \kappa^{-1} \mathsf{T} (A+B) \big)^2 \bigg\} \Bigg] \le C e^{Cn - cB^2}.
		\end{flalign*} 
	
	\end{lem}

	We next state the following concentration bound for non-intersecting Brownian bridges from \cite{U}; it is analogous to ones that appear in the context of random tilings \cite[Theorem 21]{LSRT}. We first require the notion of a height function associated with a line ensemble.  
	
	\begin{definition} 
	
	\label{htw} 
	
	For any line ensemble $\bm{\mathsf{x}} = (\mathsf{x}_1, \mathsf{x}_2, \ldots , \mathsf{x}_n) \in \llbracket 1, n \rrbracket \times \mathcal{C}(\mathbb{R})$, we define the associated \emph{height function} $\mathsf{H} = \mathsf{H}^{\bm{\mathsf{x}}} : \mathbb{R}^2 \rightarrow \mathbb{R}$ by for any $(t, w) \in \mathbb{R}^2$ setting
		\begin{flalign*}
			\mathsf{H} (t, w) = \# \big\{ j \in \llbracket 1, n \rrbracket : \mathsf{x}_j (t) > w \big\}.
		\end{flalign*}
	
	\end{definition} 

	\begin{lem}[{\cite[Lemma A.1]{U}}]
		
		\label{concentrationbridge}
		
		Let $n \ge 1$ be an integer; $\mathsf{T} > 0$ and $r, B \ge 1$ be real numbers; $\bm{u}, \bm{v} \in \mathbb{W}_n$ be $n$-tuples; and $f, g : [0, \mathsf{T}] \rightarrow \overline{\mathbb{R}}$ be continuous functions with $f \le g$, such that $f(0) \le u_n \le u_1 \le g(0)$ and $f(\mathsf{T}) \le v_n \le v_1 \le g(\mathsf{T})$. Sample non-intersecting Brownian bridges $\bm{\mathsf{x}} \in \llbracket 1, n \rrbracket \times \mathcal{C} \big( [0, \mathsf{T}] \big)$ from the measure $\mathsf{Q}_{f; g}^{\bm{u}; \bm{v}}$. Fix real numbers $t \in [0, \mathsf{T}]$ and $w \in \big[ f(t), g(t) \big]$. Denoting the event $\mathscr{E} = \big\{ \mathsf{H} (t, w) \le B \big\}$, there exists a deterministic number $\mathfrak{Y} = \mathfrak{Y} (\bm{u}; \bm{v}; f; g; \mathsf{T}; t; w; B) \ge 0$ such that
		\begin{align}
			\label{e:concentration}
			\mathbb{P} \Big[ \big| \mathsf{H}(t, w) - \mathfrak{Y} \big| \ge rB^{1/2} \Big] \le 4 e^{-r^2/4} + 2 \cdot \mathbb{P} \big[\mathscr{E}^{\complement} \big].
		\end{align}
	
		\noindent In particular, setting $B = n$, we have $\mathbb{P} \big[ \big| \mathsf{H}(t, w) - \mathfrak{Y} \big| \ge rn^{1/2} \big] \le 4e^{-r^2/4}$. 
		
	\end{lem}

\subsection{Free Convolution With Semicircle Distributions}

	\label{TransformConvolution}
	
	In this section we recall various results concerning Stieltjes transforms and free convolutions with the semicircle distribution. Fix a measure $\mu \in \mathscr{P}_{\fin}$. We define the \emph{Stieltjes transform}\index{M@$m^{\mu}$; Stieltjes transform} of $\mu$ to be the function $m = m^{\mu} : \mathbb{H} \rightarrow \mathbb{H}$ by for any complex number $z \in \mathbb{H}$ setting
	\begin{flalign}
		\label{mz0} 
		m(z) = \displaystyle\int_{-\infty}^{\infty} \displaystyle\frac{\mu (dx)}{x-z}.
	\end{flalign}
	
	\noindent If $\mu$ has a density with respect to Lebesgue measure, that is, $\mu(d x) = \varrho (x) d x$ for some $\varrho \in L^1 (\mathbb{R})$, then $\varrho$ can be recovered from its Stieltjes transform by the identity \cite[Equation (8.14)]{FPRM}, 
	\begin{flalign}
		\label{mrho} 
		\pi^{-1} \displaystyle\lim_{y \rightarrow 0^+} \Imaginary m(x + \mathrm{i} y) = \varrho (x); \qquad \displaystyle\pi^{-1} \lim_{y \rightarrow 0^+} \Real m(x + \mathrm{i} y) = H \varrho (x),
	\end{flalign}
	
	\noindent for any $x \in \mathbb{R}$. In the latter, $Hf$ denotes the Hilbert transform of any function $f \in L^1 (\mathbb{R})$, given by 
	\begin{flalign*}
		Hf (x) = \pi^{-1} \cdot \PV \displaystyle\int_{-\infty}^{\infty} \displaystyle\frac{f(w) dw}{w-x},
	\end{flalign*}
	
	\noindent where $\PV$ denotes the Cauchy principal value (assuming the integral exists as a principal value).

	The \emph{semicircle distribution}\index{0@$\mu_{\semci}$, $\mu_{\semci}^{(t)}$; semicircle distribution}\index{0@$\varrho_{\semci}$, $\varrho_{\semci}^{(t)}$; semicircle density} is a measure $\mu_{\semci} \in \mathscr{P} (\mathbb{R})$ whose density $\varrho_{\semci} : \mathbb{R} \rightarrow \mathbb{R}_{\ge 0}$ with respect to the Lebesgue measure is given by 
	\begin{flalign}
		\label{rho1} 
		\varrho_{\semci} (x) = \displaystyle\frac{(4-x^2)^{1/2}}{2\pi} \cdot \textbf{1}_{x \in [-2, 2]}, \qquad \text{for all $x \in \mathbb{R}$}. 
	\end{flalign}
	
	\noindent For any real number $t > 0$, we  denote the rescaled semicircle density $\varrho_{\semci}^{(t)}$ and distribution $\mu_{\semci}^{(t)} \in \mathscr{P}$ by 
	\begin{flalign}
		\label{rhosct}
		\varrho_{\semci}^{(t)} (x) = t^{-1/2} \varrho_{\semci} (t^{-1/2} x); \qquad \mu_{\semci}^{(t)} = \varrho_{\semci}^{(t)} (x) d x.
	\end{flalign}

	\noindent We next discuss the free convolution of a probability measure $\mu \in \mathscr{P}$ with the (rescaled) semicircle distribution $\mu_{\semci}^{(t)}$. For any $t > 0$, denote the function $M = M^{\mu} = M^{t; \mu}: \mathbb{H} \rightarrow \mathbb{C}$ and the set $\Lambda_t = \Lambda_{t; \mu} \subseteq \mathbb{H}$ by 
	\begin{flalign}
		\label{mtlambdat} 
		M(z) = z - t m(z); \quad \Lambda_t = \Big\{ z \in \mathbb{H} : \Imaginary \big( z - tm (z) \big) > 0 \Big\} = \Bigg\{ z \in \mathbb{H} : \displaystyle\int_{-\infty}^{\infty} \displaystyle\frac{\mu(d x)}{|z-x|^2} < \displaystyle\frac{1}{t} \Bigg\}.
	\end{flalign}

	\begin{lem}[{\cite[Lemma 4]{FCSD}}] 
		
		\label{mz}

		 The function $M$ is a homeomorphism from $\overline{\Lambda}_t$ to $\overline{\mathbb{H}}$. Moreover, it is a holomorphic map from $\Lambda_t$ to $\mathbb{H}$ and a bijection from $\partial \Lambda_t$ to $\mathbb{R}$.
	
	\end{lem} 

	For any real number $t \ge 0$, define $m_t = m_t^{\mu} : \mathbb{H} \rightarrow \mathbb{H}$ as follows. First set $m_0 (z) = m(z)$; for any real number $t > 0$, define $m_t$ so that	
	\begin{flalign}
		\label{mt} 
		m_t \big( z - t m_0 (z) \big) = m_0 (z), \qquad \text{for any $z \in \Lambda_t$}.
	\end{flalign}
	
	\noindent Since by \Cref{mz} the function $M(z) = z - tm_0 (z)$ is a bijection from $\Lambda_t$ to $\mathbb{H}$, \eqref{mt} defines $m_t$ on $\mathbb{H}$. By \cite[Proposition 2]{FCSD}, $m_t$ is the Stieltjes transform of a measure $\mu_t \in \mathscr{P} (\mathbb{R})$. This measure is called the \emph{free convolution}\index{0@$\mu \boxplus \mu_{\semci}^{(t)}$; free convolution} between $\mu$ and $\mu_{\semci}^{(t)}$, and we often write $\mu_t = \mu \boxplus \mu_{\semci}^{(t)}$. By \cite[Corollary 2]{FCSD}, $\mu_t$ has a density $\varrho_t = \varrho_t^{\mu} : \mathbb{R} \rightarrow \mathbb{R}_{\ge 0}$ with respect to Lebesgue measure for $t > 0$.

	\begin{rem} 
		
		\label{mtscale} 
		
		While free convolutions are typically defined between probability measures, the relation \eqref{mt} also defines the free convolution of any measure $\mu\in \mathscr{P}_{\fin}$, satisfying $A = \mu(\mathbb{R}) < \infty$, with the rescaled semicircle distribution $\mu_{\semci}^{(t)}$. Indeed, define the probability measure $\widetilde{\mu} \in \mathscr{P}$ from $\mu$ by setting $\widetilde{\mu} (I)= A^{-1} \cdot \mu (A^{1/2} I)$, for any interval $I\subseteq \mathbb{R}$. Furthermore, for any real number $s \ge 0$, define the probability measure $\widetilde{\mu}_s = \widetilde{\mu} \boxplus \mu_{\semci}^{(s)}$, and denote its Stieltjes transform by $\widetilde{m}_s$. Then, define the free convolution $\mu_t = \mu \boxplus \mu_{\semci}^{(t)}$ and its Stieltjes transform $m_t = m_{t; \mu}$ by setting
		\begin{flalign*} 
			\mu_t (I) = A \cdot \widetilde{\mu}_t (A^{-1/2} I), \qquad \text{for any interval $I \subseteq \mathbb{R}$, so that} \qquad	m_t (z) = A^{1/2} \cdot \widetilde m_t ( A^{-1/2} z),
		\end{flalign*} 
		
		\noindent where the second equality follows from the first by \eqref{mz0}. Then, 
		\begin{align*}
			m_{t} \big( z-t m_{0}(z) \big) = A^{1/2} \cdot \widetilde{m}_t \Big( A^{-1/2} \big(z- t m_{0}(z) \big) \Big) & = A^{1/2} \cdot \widetilde{m}_t \big( A^{-1/2} z- t \widetilde{m}_{0}(A^{-1/2} z) \big) \\
			& = A^{1/2} \cdot \widetilde{m}_{0}(A^{-1/2} z)= m_{0}(z),
		\end{align*}
		
		\noindent so that \eqref{mt} continues to hold for $m_t$. In particular, \Cref{mz} also hold for $\mu$.
		
	\end{rem} 
	
	\begin{rem}
		
		\label{mtscalebeta}
		
		Let us describe a scaling invariance for time under free convolutions. Fix a measure $\mu \in \mathscr{P}_{\fin}$ with $\mu (\mathbb{R}) < \infty$, and let $m_s$ denote the Stieltjes transform of $\mu_s = \mu \boxplus \mu_{\semci}^{(s)}$, for any real number $s \ge 0$. Fix a real number $\beta > 0$, and define the measure $\widetilde{\mu} \in \mathscr{P}_{\fin}$ by setting $\widetilde{\mu} (I) = \mu (\beta^{1/2} \cdot I)$, for any interval $I \subseteq \mathbb{R}$. Denote the Stieltjes transform of $\widetilde{\mu}$ by $\widetilde{m} = m_{\tilde{\mu}}$, and let $\widetilde{m}_s$ denote the Stieltjes transform of $\widetilde{\mu}_s = \widetilde{\mu} \boxplus \mu_{\semci}^{(s)}$ for any $s > 0$. Then, observe for any real number $t \ge 0$ and complex number $z \in \mathbb{H}$ that 
		\begin{flalign} 
			\label{mtbeta} 
			\widetilde{m}_t (z) = \beta^{1/2} \cdot m_{\beta t} (\beta^{1/2} z).
		\end{flalign} 
		
		\noindent Indeed, this holds at $t = 0$ by \eqref{mz0}; for $t > 0$, we have
		\begin{flalign*}
			\widetilde{m}_t \big( z - t \widetilde{m}_0 (z) \big) = \widetilde{m}_0 (z) = \beta^{1/2} \cdot m_0 (\beta^{1/2} z) & = \beta^{1/2} \cdot m_{\beta t} \big( \beta^{1/2} z - \beta t m_0 (\beta^{1/2} z) \big) \\
			& = \beta^{1/2} \cdot m_{\beta t} \Big( \beta^{1/2} \big( z - t \widetilde{m}_0 (z) \big) \Big),
		\end{flalign*}
		
		\noindent which by \Cref{mz} (and \Cref{mtscale}, if $\mu$ is not a probability measure) implies \eqref{mtbeta}. The equality \eqref{mtbeta}, with the first statement of \eqref{mrho}, in particular implies that $\widetilde{\mu}_t (I) = \mu_{\beta t} (\beta^{1/2} \cdot I)$ for any interval $I \subseteq \mathbb{R}$.
	\end{rem}

	\subsection{Dyson Brownian Motion} 
	
	\label{LambdaEquation}

	In this section we recall properties about Dyson Brownian motion. Fix an integer $n \ge 1$ and a sequence $\bm{\lambda} (0) = \big( \lambda_1 (0), \lambda_2 (0), \ldots , \lambda_n (0) \big) \in \overline{\mathbb{W}}_n$. Define the sequence $\bm{\lambda} (t) = \big(\lambda_1 (t), \lambda_2 (t), \ldots , \lambda_n (t) \big) \in \overline{\mathbb{W}}_n$, for $t \ge 0$, to be the unique strong solution (see \cite[Proposition 4.3.5]{IRM} for its existence) to the stochastic differential equations 
	\begin{flalign}
		\label{lambdaequation}
		d \lambda_i (t) = d B_i (t)+ \displaystyle\sum_{\substack{1 \le j \le n \\ j \ne i}}	\displaystyle\frac{dt}{\lambda_i (t) - \lambda_j (t)} ,\quad 1\leq i\leq n.
	\end{flalign}
	
	\noindent The system \eqref{lambdaequation} is called \emph{Dyson Brownian motion} (with $\beta=2$),\index{D@Dyson Brownian motion} run for time $t$, with \emph{initial data} $\bm{\lambda} (0)$; the $\lambda_i$ are sometimes referred to as \emph{particles}.

	\begin{rem}
		
		\label{scalemotion}
		
		As in \Cref{scale}, Dyson Brownian motion admits the following invariance under diffusive scaling for any real number $\sigma > 0$. If $\bm{\lambda} (t) = \big( \lambda_1 (t), \lambda_2 (t), \ldots , \lambda_n (t) \big) \in \overline{\mathbb{W}}_n$ solves \eqref{lambdaequation} then, denoting $\widetilde{\lambda}_j (t) = \sigma^{1/2} \cdot \lambda_j (\sigma^{-1} t)$, the process $\widetilde{\bm{\lambda}} (t) = \big( \widetilde{\lambda}_1 (t), \widetilde{\lambda}_2 (t), \ldots , \widetilde{\lambda}_n (t) \big) \in \overline{\mathbb{W}}_n$ also solves \eqref{lambdaequation}. This again follows from the invariance of the Brownian motions $B_i$ under the same scaling.
	
	\end{rem}
	
	\begin{rem}
		
		\label{kscaling}
		
		To later analyze limit shapes, we will occasionally consider a scaled variant of \eqref{lambdaequation}. In particular, set $\widetilde{\lambda}_j (t) = n^{-1} \cdot \lambda_j (nt)$ for each $t > 0$ and $j \in \llbracket 1, n \rrbracket$, which amounts to scaling the time $t$ and space $x$ by $n^{-1}$. Then, the process $\widetilde{\bm{\lambda}} (t) = \big( \widetilde{\lambda}_1 (t), \widetilde{\lambda}_2 (t), \ldots , \widetilde{\lambda}_n (t) \big) \in \overline{\mathbb{W}}_n$ satisfies
		\begin{flalign}
		\label{1/nscaling}
		d\widetilde \lambda_i (t) =  \frac{d B_i (t)}{\sqrt n} + \frac{1}{n} \displaystyle\sum_{\substack{1 \le j \le n \\ j \ne i}}	\displaystyle\frac{dt}{\widetilde\lambda_i (t) - \widetilde\lambda_j (t)} ,\quad 1\leq i\leq n.
	\end{flalign}
	\end{rem}

	We next describe the relation between Dyson Brownian motion, random matrices, and non-intersecting Brownian bridges, to which end we require some additional terminology. A \emph{random matrix} is a matrix whose entries are random variables. The \emph{Gaussian Unitary Ensemble}\index{G@$\bm{G}_n$; Gaussian Unitary Ensemble (GUE)} is an $n \times n$ random Hermitian matrix $\bm{G} = \bm{G}_n$ with random complex entries $\{ w_{ij} \}$ (for $i, j \in \llbracket 1, n \rrbracket$) defined as follows. Its diagonal entries $\{ w_{jj} \}$ are standard real Gaussian random variables, and its upper-triangular entries $\{ w_{ij} \}_{i < j}$ are standard complex Gaussian random variables (that is, whose real and imaginary parts are independent Gaussian random variables, each of variance $1/2$); these entries are mutually independent, and the lower triangular entries $\{ w_{ij} \}_{i > j}$ are determined from the upper triangular ones by the Hermitian symmetry relation $w_{ij} = \overline{w_{ji}}$.
	
	The \emph{Hermitian Brownian motion} $\bm{G} (t) = \bm{G}_n (t)$ is a stochastic process (over $t \ge 0$) on $n\times n$ random matrices, whose entries $\big\{ w_{ij} (t) \big\}$ are defined as follows. Its diagonal entries $\big\{ w_{jj} (t) \big\}$ are Brownian motions of variance $1$, and its upper triangular entries $\big\{ w_{ij} (t) \big\}_{i < j}$ are standard complex Brownian motions (that is, whose real and imaginary parts are independent Brownian motions, each of variance $1/2$). These entries are again mutually independent, and the lower triangular entries $\big\{ w_{ij} (t) \big\}_{i > j}$ are determined from its upper triangular ones by symmetry, $w_{ij} (t) =\overline{ w_{ji} }(t)$. Observe that $\bm{G} (1)$ has the same law as a GUE matrix $\bm{G}$.
	
	 The following lemma from \cite{MERM} (stated as below in \cite{MCNP}) interprets Dyson Brownian motion in terms of sums of random matrices, and also in terms of non-intersecting Brownian motions conditioned to never intersect; we recall the definition of the latter in terms of Doob $h$-transforms from \cite[Section 6.2]{MCNP}.

	\begin{lem}[{\cite[Theorems 3 and 4]{MCNP}}]
		
		\label{lambdat}
		
		Fix an integer $n \ge 1$ and a sequence $\bm{\lambda} (0) \in \overline{\mathbb{W}}_n$. For any real number $t > 0$, let $\bm{\lambda} (t) \in \overline{\mathbb{W}}_n$ denote Dyson Brownian motion, run for time $t$, with initial data $\bm{\lambda} (0)$. Further let $\bm{A}$ denote an $n \times n$ diagonal matrix whose eigenvalues are given by $\bm{\lambda} (0)$, and let $\bm{G} (t) = \bm{G}_n (t)$ denote an $n \times n$ Hermitian Brownian motion. 
		
		\begin{enumerate}
			\item The law of the eigenvalues of $\bm{A} + \bm{G}(t)$ coincides with that of $\bm{\lambda} (t)$, jointly over $t \ge 0$. 
			\item Consider $n$ Brownian motions $\sfX = (\sfx_1, \sfx_2, \ldots , \sfx_n) \in \llbracket 1, n \rrbracket \times \mathcal{C} (\mathbb{R}_{\ge 0})$, with variances $1$ and starting data $\bm{\lambda} (0)$, conditioned to never intersect. Then, $\big( \sfX(t) \big)_{t \ge 0}= \big( \sfx_1(t), \mathsf{x}_2 (t), \ldots \big)$ has the same law as $\big( \bm{\lambda} (t) \big)_{t \ge 0}$. 
		\end{enumerate}
	\end{lem} 
	
	\begin{rem}
		
		\label{motionscale2}
		
		By the second part of \Cref{lambdat}, for any real number $\sigma > 0$, the paths $\sigma^{-2/3} \cdot \big( \sfx_j ( \sigma^{1/3} t) \big)$ are given by Brownian motions, with variances $\sigma^{-1}$, conditioned to never intersect.
		 
	\end{rem}

	\begin{rem}
		
		\label{tobridge}
		
		Given a real number $\mathsf{T} > 0$ and a Brownian bridge $B : [0, \mathsf{T}] \rightarrow \mathbb{R}$, conditioned to start at some $u \in \mathbb{R}$ and end at $0$ (namely, $B(0) = u$ and $B(\mathsf{T}) = 0$), observe that $W: \mathbb{R}_{> 0} \rightarrow \mathbb{R}$ defined by $W(t) = \mathsf{T}^{-1} (\mathsf{T} + t) \cdot B \big( \mathsf{T} t / (\mathsf{T} + t) \big)$ has the law of a Brownian motion starting at $u$ (that is, with $W(0) = u$). Thus, fixing $\bm{u} \in \overline{\mathbb{W}}_n$ and letting $\bm{\mathsf{y}} = (\mathsf{y}_1, \mathsf{y}_2, \ldots , \mathsf{y}_n) \in \llbracket 1, n \rrbracket \times \mathcal{C} \big( [0, \mathsf{T}] \big)$ denote non-intersecting Brownian bridges sampled under the measure $\mathsf{Q}^{\bm{u}; \bm{0}_n}$, then defining 
	\begin{align}
		\label{lambday} 
	\lambda_j(t)=\frac{\sfT+t}{\sfT} \cdot \sfy_j \left(\frac{\sfT t}{\sfT+t}\right),
	\end{align}

	\noindent for each $(j, t) \in \llbracket 1, n \rrbracket \times [0, \infty)$ the process $\bm{\lambda} (t) = \big( \lambda_1 (t), \lambda_2 (t), \ldots , \lambda_n (t) \big)$ defines Brownian motions starting from $\bm{u}$, conditioned to never intersect (as, under \eqref{lambday}, the paths in $\bm{\lambda} (t)$ do not intersect if and only if those in $\bm{\mathsf{y}} (t)$ do not intersect). By the second part of \Cref{lambdat}, this has the law of Dyson Brownian motion with initial data $\bm{u}$, run for time $t$. As such, we can view the latter as a special case of non-intersecting Brownian bridges. 
		
	\end{rem}

	\subsection{Estimates for Dyson Brownian Motion}
	
	\label{MotionEstimates}
	
	In this section we state concentration bounds and gap estimates for Dyson Brownian motion. We begin by recalling the concentration results from \cite{MCTMGP}, to which end we require the notion of a classical location with respect to a density.
	
	\begin{definition}
		
		\label{gammarho} 
		
		Let $\mu \in \mathscr{P}_{\fin}$ denote a measure of finite total mass $\mu (\mathbb{R}) = A$. For any integers $n \ge 1$ and $j \in \mathbb{Z}$, we define the \emph{classical location} (also called \emph{$n^{-1}$-quantiles}) with respect to $\mu$, $\gamma_j = \gamma_j^{\mu} = \gamma_{j; n} = \gamma_{j; n}^{\mu} \in \mathbb{R}$ by setting\index{0@$\gamma_{j;n}^{\mu}$; classical locations}
		\begin{flalign*}
			\gamma_j = \displaystyle\sup \Bigg\{ \gamma \in \mathbb{R} : \displaystyle\int_{\gamma}^\infty d \mu(x) \ge \displaystyle\frac{A(2j-1)}{2n} \Bigg\}, \qquad \text{if $j \in \llbracket 1, n \rrbracket$},
		\end{flalign*}  
		
		\noindent and also setting $\gamma_j = \infty$ if $j < 1$ and $\gamma_j = -\infty$ if $j > n$.
		
	\end{definition}
	
	The following lemma due to\footnote{In \cite{MCTMGP}, the probability on the right side of \eqref{lambdatprobability} was written to be $1 - C n^{-D}$ for any $D > 1$, but it can be seen from the proof (see that of \cite[Proposition 3.8]{MCTMGP}, where $\delta$ there is $5 / 4$ here) that it can be taken to be $1 - C e^{-(\log n)^2}$ instead.} \cite{MCTMGP} (together with the scale invariance \Cref{scalemotion}) provides a concentration, or \emph{rigidity}, estimate for the locations of bulk particles (namely, those sufficiently distant from the first and last) under Dyson Brownian motion around the classical locations of a free convolution measure.

	\begin{lem}[{\cite[Corollary 3.2]{MCTMGP}}]
		
		\label{concentrationequation}
		
		For any real number $D > 1$, there exists a constant $C = C (D) > 1$ such that the following holds. Fix an integer $n \ge 1$ and sequence $\bm{\lambda} (0) \in \overline{\mathbb{W}}_n$ with $-n^D \le \min \bm{\lambda} (0) \le \max \bm{\lambda} (0) \le n^D$. Denote the measure $\mu = n^{-1} \sum_{j=1}^n \delta_{\lambda_j (0) /n} \in \mathscr{P}$, and set $\mu_t = \mu \boxplus \mu_{\semci}^{(t)}$; also denote the classical locations $\gamma_j (t) = \gamma_{j; n}^{\mu_t} \in \mathbb{R}$. Letting $\bm{\lambda} (t) = \big( \lambda_1 (t), \lambda_2 (t), \ldots , \lambda_n (t) \big) \in \overline{\mathbb{W}}_n$ denote Dyson Brownian motion with initial data $\bm{\lambda} (0)$, we have 
		\begin{flalign}
			\label{lambdatprobability} 
			\mathbb{P} \Bigg[ \bigcap_{j = 1}^n \bigcap_{t \in [0, n^D]} \big\{ \gamma_{j+\lfloor (\log n)^6 \rfloor}(t)-n^{-D} \le n^{-1} \lambda_j (nt) \le \gamma_{j- \lfloor (\log n)^6 \rfloor}(t)+n^{-D} \big\} \Bigg] \ge 1 - C e^{ - (\log n)^2}. 
		\end{flalign}
	\end{lem}

	\begin{rem} 
		
		Let us briefly explain the differences between \cite[Corollary 3.2]{MCTMGP} and \Cref{concentrationequation}. The first is that \cite[Corollary 3.2]{MCTMGP} considers non-intersecting Brownian motions $\widetilde{\bm{\lambda}}(t)$ with variances $n^{-1}$, while \Cref{concentrationequation} considers them with variances $1$. The second is that \cite[Corollary 3.2]{MCTMGP} (with the $(\mathfrak{c}, \delta)$ there equal to $(2D+1, 1)$ here) shows with high probability that 
		\begin{flalign}
			\label{gamma2} 
			\widetilde{\gamma}_{j+ \lfloor (\log n)^6 \rfloor} (\widetilde{t}) - n^{-2D} \le \widetilde{\lambda}_j (\widetilde{t}) \le \widetilde{\gamma}_{j - \lfloor (\log n)^6 \rfloor} (\widetilde{t}) + n^{-2D},
		\end{flalign} 
	
		\noindent under the assumptions that $\widetilde{t} \le (\log n)^{-2}$ and $-a \le \min \bm{\lambda}(0) \le \max \bm{\lambda} (0) \le a$ for some fixed constant $a > 0$, where the $\widetilde{\gamma}_j$ are the classical locations of $\widetilde{\mu} \boxplus \mu_{\semci}^{(\tilde{t})}$, for $\widetilde{\mu} = n^{-1} \sum_{j=1}^n \delta_{\tilde{\lambda}_j (0)}$. On the other hand, \Cref{concentrationequation} allows $t \le n^D$ and $-n^D \le \min \bm{\lambda}(0) \le \max \bm{\lambda}(0) \le n^D$.
				
		Both of these are resolved by the rescalings given by \Cref{scalemotion} and \Cref{kscaling}. Indeed, adopting the assumptions of \Cref{concentrationequation} and setting $\widetilde{\lambda}_j (s) = n^{-D-1} \cdot \lambda_j (n^{2D+1} s)$, \cite[Corollary 3.2]{MCTMGP} yields \eqref{gamma2}. This, together with the fact that $\widetilde{\gamma}_j (t) = n^{-D} \cdot \gamma_j (t)$ (as a quick consequence of \Cref{mtscalebeta}), yields the event in \eqref{lambdatprobability} from \eqref{gamma2}.
		
		The third difference is that \cite[Corollary 3.2]{MCTMGP} states a probability bound of $1 - Cn^{-D}$ for any $D>1$, instead of $1 - e^{-c(\log n)^2}$. However, the proof of \cite[Corollary 3.2]{MCTMGP} implicitly shows this stronger probability bound. Indeed, the only places in the proof of \cite[Corollary 3.2]{MCTMGP} where a probability estimate appears are in \cite[Proposition 3.8]{MCTMGP} and \cite[Proposition 3.9]{MCTMGP}; the latter is already stated with the exponential estimate $1 - e^{-n}$. In the proof of the former, it is stated above \cite[Equation (3.35)]{MCTMGP} that the event $\Omega$ on which \cite[Proposition 3.9]{MCTMGP} holds has probability at least $1 - C |\mathcal{L}| \cdot \log n \cdot e^{-c (\log n)^{1+\delta}}$; here, $|\mathcal{L}|$ grows at most polynomially in $n$ (see \cite[Equations (3.2) and (3.16)]{MCTMGP}), and we take $\delta=1$. This yields a probability bound of $1 - C' e^{-c'(\log n)^2}$, as in \Cref{concentrationequation}.
	\end{rem}

	We next state a result bounding the gaps between the first particles under Dyson Brownian motion whose initial data is ``sufficiently small.'' 
	
	\begin{lem}[{\cite[Corollary 4.3]{ERGIM}}]
		
		\label{initialsmall2}
		
		For any real number $B > 1$, there exists a constant $c = c(B) > 0$ such that the following holds. Let $n \ge 1$ be an integer and $\bm{\lambda} = (\lambda_1, \lambda_2, \ldots , \lambda_n) \in \overline{\mathbb{W}}_n$ be a sequence of real numbers such that $\lambda_1 - \lambda_n < c n^{2/3}$. Letting $\bm{\lambda}(s) = \big( \lambda_1 (s), \lambda_2 (s), \ldots , \lambda_n (s) \big) \in \overline{\mathbb{W}}_n$ denote Dyson Brownian motion with initial data, run for time $s$. Then, 
		\begin{flalign*}
			\mathbb{P} \Bigg[ \bigcap_{t \in [1/B, B]} \bigcap_{1 \le j < k \le \lfloor n/2 \rfloor} \Big\{ \big| \lambda_j (tn^{1/3}) - \lambda_k (tn^{1/3}) \big| \le 25 t^{1/2} (k^{2/3}-j^{2/3} &) + (\log n)^{20} j^{-1/3} \Big\} \Bigg] \\
			& \qquad \ge 1- c^{-1} e^{-c (\log n)^2}.
		\end{flalign*}

	\end{lem}

	We next state the following result bounding the location of the last particle in Dyson Brownian motion, assuming its initial data is not too densely packed.      
	
	\begin{lem}[{\cite[Corollary 4.7]{ERGIM}}]
		
		\label{p:extreme}
		
		For any real numbers $B, D > 1$, there exist constants $c = c(B) > 1$, $C_1 = C_1 (B) > 1$ and $C_2 = C_2 (B, D) > 1$ such that the following holds. Let $k, n \ge 2$ be integers, and let $L \in [1, k^D]$ be a real number such that $n = L^{3/2} k$. Let $\bm{\lambda}(s) = \big( \lambda_1 (s), \lambda_2 (s), \ldots , \lambda_n (s) \big) \in \overline{\mathbb{W}}_n$ denote Dyson Brownian motion with initial data $\bm{\lambda}(0)$, run for time $s$. Suppose that, for some real number $M \ge 1$, we have
		\begin{align}\label{e:xi-xj2}
			\lambda_i(0)-\lambda_j(0)\geq \left(\frac{j-i}{BL^{3/4}k}-M \right)k^{2/3}, \qquad \text{for each $1 \le i \le j \le n$}.
		\end{align} 
		
		\noindent Then, for any $t \in [0, 1]$, we have 
		\begin{align}\label{e:eigbound}
			\begin{aligned}
				\mathbb{P} \Bigg[ \lambda_n (tk^{1/3}) \geq \lambda_n(0)-C_1 k^{2/3} \Big( tL^{3/4} \big| \log (2t^{-1}) \big|^2 + ( & Mt)^{1/2} L^{3/8}+ (tk^{-1})^{1/2} (\log n)^3 \Big)  \Bigg] \\
				& \qquad \qquad \qquad \ge 1 - C_2 e^{-c (\log n)^2}. 
			\end{aligned} 
		\end{align}
		
	\end{lem}

	\subsection{Edge Statistics of Dyson Brownian Motion}
	
	\label{MotionEdge} 
	
	In this section we state a result, essentially due to \cite{FESD}, on the edge statistics of Dyson Brownian motion (recall \Cref{LambdaEquation}). It requires the following assumption. 

	\begin{assumption}
		
		\label{yconvergedelta0}
		
		Fix real numbers $\delta_0 \in (0, 1)$ and $t > 0$; a sequence $\bm{\delta} = (\delta_1, \delta_2, \ldots )$ of real numbers such that $\lim_{n \rightarrow \infty} \delta_n = 0$; and a measure $\nu \in \mathscr{P}_0$ such that, for each integer $n \ge 1$,
		\begin{flalign}
			\label{yndelta1} 
			 \displaystyle\inf_{s \in \supp \nu} \displaystyle\lim_{\varepsilon \rightarrow 0} \displaystyle\int_{-\infty}^{\infty} \displaystyle\frac{\nu (dx)}{(s-x)^2 + \varepsilon^2} > t^{-1} + \delta_0.
		\end{flalign}
	
		\noindent By \eqref{yndelta1}, there exists a unique real solution $z_0 > \max (\supp \nu)$ to the equation 
		\begin{flalign}
			\label{z0tsigma} 
			\displaystyle\int_{-\infty}^{\infty} \displaystyle\frac{\nu (dy)}{(y-z_0)^2} = t^{-1}, \qquad \text{so set} \qquad \sigma = \sigma_{\nu; t} = \bigg( t^3 \displaystyle\int_{-\infty}^{\infty} \displaystyle\frac{\nu(dy)}{(z_0 - y)^3} \bigg)^{-1/3},
		\end{flalign} 
		
		\noindent as the integral on the left side of \eqref{z0tsigma} is continuous and decreasing in $z_0 > \max (\supp \nu)$.
		
		 For each integer $n \ge 1$, let $\bm{y} = \bm{y}^n = (y_1, y_2, \ldots , y_n) \in \overline{\mathbb{W}}_n$ be a sequence. Assume for any integer $n \ge 1$ and real numbers $a < b$ that we have the bounds 
		\begin{flalign}
			\label{yndelta} 
				& \Bigg| \displaystyle\frac{1}{n} \displaystyle\sum_{j = 1}^n \textbf{1}_{y_j /n \in [a, b]} - \displaystyle\int_a^b \nu (dx) \Bigg| \le \delta_n; \qquad  \displaystyle\max_{1 \le j \le n} \dist (n^{-1} y_j, \supp \nu) \le \delta_n.
		\end{flalign}
	
		\noindent For each integer $n \ge 1$, let $\bm{\lambda} = \bm{\lambda}^n \in \overline{\mathbb{W}}_n$ denote Dyson Brownian motion run for time $nt$, with initial data $\bm{y}^n$. 
		
	\end{assumption}

	Let us briefly explain \Cref{yconvergedelta0}. The second bound in \eqref{yndelta1} indicates that $\nu$ cannot decay too slowly near the edge of its support; the parameter $\sigma$ in \eqref{z0tsigma} will be a scaling factor; and \eqref{yndelta} imposes a rate of convergence to $\nu$ of the empirical measure of the initial data $\bm{y}$ for $\bm{\lambda}$. 
	
	The following lemma, which is a quick consequence (and a uniform variant) of \cite[Theorem 1.1]{FESD}, states that the gaps between the largest particles of $\bm{\lambda}$ converge to those of the Airy point process $\bm{\mathfrak{a}} = (\mathfrak{a}_1, \mathfrak{a}_2, \ldots )$ (recall \Cref{a0}). It is proven in \Cref{Proofydelta}.
	
	\begin{lem}
		
		\label{yconvergedelta}
		
		Adopt \Cref{yconvergedelta0}; fix an integer $k \ge 1$ and a $k$-tuple of real numbers $\bm{r} = (r_1, r_2, \ldots , r_k) \in \mathbb{R}_{\ge 0}^k$. For any $\varepsilon \in (0, 1)$, there exists a constant $N_0 = N_0 (\varepsilon, \delta_0, t, \bm{\delta}, \nu, k, \bm{r}) > 1$ (that does not depend on the $\bm{y}^n$) such that, for any integer $n \ge N_0$, we have 
		\begin{flalign}
			\label{convergelambda}
			\Bigg| \mathbb{P} \bigg[ \bigcap_{j=1}^k \big\{ \sigma n^{-1/3} (\lambda_j - \lambda_{j+1}) \ge r_j \big\} \bigg] - \mathbb{P} \bigg[ \bigcap_{j=1}^k \big\{ \mathfrak{a}_j - \mathfrak{a}_{j+1} \ge r_j \}\bigg] \Bigg| \le \varepsilon.
		\end{flalign}
		
	\end{lem}

	\subsection{Dyson Brownian Motion and Non-Intersecting Bridges}
	
	\label{MotionCurves} 
	
	In this section we recall results that relate non-intersecting Brownian bridges (with no upper or lower boundary) to Dyson Brownian motion. We first recall the following lemma giving a description for the law of the locations of these bridges at a single time; it is essentially due to \cite{LDASI,FOAMI} (see also the exposition in \cite[Section 2.1]{LDSI}), but we provide its short proof in \Cref{ProofBridgeSum1} below. In what follows, for any integer $k$ and $k$-tuple $\bm{a} = (a_1, a_2, \ldots , a_k) \in \mathbb{C}^k$, we let $\diag (\bm{a})$\index{D@$\diag$} denote the $k \times k$ diagonal matrix whose $(j, j)$ entry is $a_j$, for each $j \in \llbracket 1, k \rrbracket$. For any $n \times n$ Hermitian matrix $\bm{M}$, we also let $\eig (\bm{M}) \in \overline{\mathbb{W}}_n$\index{E@$\eig$} denote the $n$-tuple of eigenvalues of $\bm{M}$, ordered to be non-increasing; we additionally let $\bm{W}^*$ denote the conjugate transpose of any complex matrix $\bm{W}$.
	
	\begin{lem}
		
		\label{t:lawt}
		
		Let $n \ge 1$ be an integer and $\bm{u}, \bm{v} \in \overline{\mathbb{W}}_n$ be $n$-tuples. Define the $n \times n$ diagonal matrices $\bm{U} = \diag ( \bm{u})$ and $\bm{V} = \diag ( \bm{v})$, and let $\bm{G}$ denote an $n \times n$ GUE random matrix. Letting $\mathsf{T} > 0$ be a real number, and sample non-intersecting Brownian bridges $\bm{\mathsf{x}} = (\mathsf{x}_1, \mathsf{x}_2, \ldots , \mathsf{x}_n) \in \llbracket 1, n \rrbracket \times \mathcal{C} \big( [0, \mathsf{T}] \big)$ from the measure $\mathsf{Q}^{\bm{u}; \bm{v}}$. For any real number $t \in [0, \sfT]$, the $n$-tuple $ \bm{\mathsf{x}} (t) \in \overline{\mathbb{W}}_n$ has the same law as     
		\begin{align}\label{e:lawy}
			\eig \bigg( \bm{A} + \Big( \displaystyle\frac{t (\sfT-t)}{\sfT} \Big)^{1/2} \cdot \bm{G} \bigg), \quad \text{where} \quad \bm{A} = \frac{\sfT-t}{\sfT} \cdot \bm{U}+ \frac{t}{\sfT} \cdot \bm{W} \bm{V} \bm{W}^*.
		\end{align}
	
		\noindent Here, $\bm{W}$ is a random unitary matrix whose law is given by
		\begin{flalign}
			\label{wuv}
			\mathbb{P} [d \bm{W}] = Z^{-1} \exp \Big( \mathsf{T}^{-1} \Tr \bm{U} \bm{W} \bm{V} \bm{W}^* \Big) d \bm{W}, \quad Z = Z_n (\bm{U}, \bm{V}) = \displaystyle\int_{\mathrm{U}(n)} e^{- \mathsf{T}^{-1} \Tr \bm{U} \bm{W} \bm{V} \bm{W}^*} d \bm{W},
		\end{flalign}
		
		\noindent and $d \bm{W}$ denotes the Haar measure on the group $\mathrm{U}(n)$ of $n\times n$ unitary matrices.
	\end{lem}

	\begin{rem}
		
		\label{tuvwx}
		
		Adopting the notation of \Cref{t:lawt}, \Cref{t:lawt} indicates that the law of $ \bm{\mathsf{x}} (t)$ is given by Dyson Brownian motion with initial data $\eig (\bm{A})$, run for time $t ( 1 - t \mathsf{T}^{-1})$.
		
	\end{rem}

	The following corollary uses \Cref{t:lawt} with \Cref{initialsmall2} to bound the gaps between non-intersecting Brownian bridges, run for time much longer than the sizes of the supports of their starting and ending data; it is established in \Cref{ProofBridgeSum2} below.

	\begin{cor}
		\label{p:smallinitial}
		
		For any real numbers $A, B > 1$, there exist constants $c = c(A, B) > 0$, $C_1 = C_1 (B) > 1$, and $C_2 = C_2 (A, B) > 1$ such that the following holds. Let $n \ge 1$ be an integer; $T \in [C_1, AC_1]$ be a real number; and $\bm{u}, \bm{v} \in \overline{\mathbb{W}}_n$ be $n$-tuples with
		\begin{flalign}
			\label{buvb}
			-Bn^{2/3} \le \min \bm{u} \le \max \bm{u} \le Bn^{2/3}; \qquad -Bn^{2/3} \le \min \bm{v} \le \max \bm{v} \le Bn^{2/3}.
		\end{flalign}
		
		\noindent Sample non-intersecting Brownian bridges $\bm{\mathsf{x}} = (\mathsf{x}_1, \mathsf{x}_2, \ldots , \mathsf{x}_n) \in \llbracket 1, n \rrbracket \times \mathcal{C} \big( [0, Tn^{1/3}] \big)$ under the measure $\mathsf{Q}^{\bm{u}; \bm{v}}$. Then,
		\begin{flalign}
			\label{xjk2n13c}
			\begin{aligned} 
				\mathbb{P} \Bigg[ \bigcap_{t \in [T/4,3T/4]} \bigcap_{1 \le j \le k \le \lfloor n/2 \rfloor} \Big\{ \big| \mathsf{x}_j (tn^{1/3}) - \mathsf{x}_k (tn^{1/3}) \big| \le C_2 (& k^{2/3} - j^{2/3}) + (\log n)^{25} j^{-1/3} \Big\} \Bigg] \\
				& \qquad \qquad \quad \ge 1 - c^{-1} e^{-c(\log n)^2}.  
			\end{aligned} 
		\end{flalign}
		
	\end{cor}

	\subsection{Brownian Watermelon and Airy Line Ensemble Estimates}
	
	\label{PathsUV0}
	
	In this section we provide estimates for the locations of paths in the scaled parabolic Airy line ensemble and in an ensemble of non-intersecting Brownian bridges conditioned to start and end at $0$; the latter ensemble is sometimes referred to as a \emph{Brownian watermelon}. In what follows, for each real number $y \in [0, 1]$, let $\gamma_{\semci} (y)$ to be the classical location of the semicircle distribution,\index{0@$\gamma_{\semci}$, $\gamma_{\semci; n}$} defined to be  
	\begin{flalign} 
		\label{gammascy} 
		\text{unique} \quad \gamma \in [-2, 2] \quad \text{solving the equation} \quad (2\pi)^{-1} \int_{\gamma}^2 (4-x^2)^{1/2} dx = y.
	\end{flalign} 

	\noindent For any integers $n \ge 1$ and $j \in \mathbb{Z}$ we let $\gamma_{\semci; n} (j)$ be the classical location (recall \Cref{gammarho}) with respect to the semicircle distribution, given by
	\begin{flalign}
		\label{gammaj} 
		\gamma_{\semci; n} (j) = \gamma_{j; n}^{\mu_{\semci}} = \gamma_{\semci} \Big( \displaystyle\frac{2j-1}{2n} \Big), \quad \text{which satisfies} \quad \displaystyle\frac{1}{2 \pi} \displaystyle\int_{\gamma_{\semci; n} (j)}^2 (4-x^2)^{1/2} dx = \displaystyle\frac{2j-1}{2n}.
	\end{flalign}

	We begin with the following lemma which bounds the classical locations of the semicircle distribution $\gamma_{\semci}$ (recall \eqref{gammascy}) and their derivatives; we establish it in \Cref{ProofBridgeSum1} below. 
	
	\begin{lem} 
		
		\label{gammaderivative} 
		
		The following two statements hold.
		\begin{enumerate} 
			\item If $y \in [0, 1]$ then $2y^{2/3} \le 2 - \gamma_{\semci} (y) \le 8y^{2/3}$. 
			\item If $y \in [0, 1]$, then $-\gamma_{\semci}' (y) \ge 2^{-3/2} y^{-1/3}$. Moreover, if $y \in [0, 1/2]$, then $-\gamma_{\semci}' (y) \le \pi y^{-1/3}$. 
		\end{enumerate} 
		
	\end{lem}

	The next lemma from \cite{U} provides a concentration bound for paths in Brownian watermelons. 
	
	\begin{lem}[{\cite[Lemma 2.18]{U}}]
		
		\label{estimatexj} 
		
		For any real number $D > 1$, there exists a constant $C = C(D) > 1$ such that the following holds. Adopt the notation of \Cref{estimatexj2}; assume that $b-a \le n^D$; fix real numbers $u, v \in \mathbb{R}$; and assume that $\bm{u} = (u, u, \ldots , u) \in \overline{\mathbb{W}}_n$ and $\bm{v} = (v, v, \ldots , v) \in \overline{\mathbb{W}}_n$ (where $u$ and $v$ appear with multiplicity $n$).

		\begin{enumerate} 
		
		\item With probability at least $1 - C e^{-(\log n)^5}$, we have
		\begin{flalign*}
			\displaystyle\max_{j \in \llbracket 1, n \rrbracket}  \displaystyle\sup_{t \in [a, b]} \Bigg|  \mathsf{x}_j (t) - n^{1/2} \bigg( \displaystyle\frac{(b-t)(t-a)}{(b-a)} & \bigg)^{1/2} \cdot \gamma_{\semci; n} (j) - \displaystyle\frac{b - t}{b-a} \cdot u  - \displaystyle\frac{t-a}{b-a} \cdot v  \Bigg| \\ 
			& \le (\log n)^9 \cdot n^{-1/6} (b-a)^{1/2} \cdot \min \{  j, n-j+1 \}^{-1/3}.
		\end{flalign*} 
	
		\item With probability at least $1 - C e^{-(\log n)^5}$, we have
		\begin{flalign*}
			& \displaystyle\max_{j \in \llbracket 1, n \rrbracket} \displaystyle\sup_{t \in [a, b]} \Bigg( \bigg| \mathsf{x}_j (t) - \displaystyle\frac{b - t}{b-a} \cdot u - \displaystyle\frac{t-a}{b-a} \cdot v \bigg| -   (8n)^{1/2} \bigg( \displaystyle\frac{(b-t)(t-a)}{b-a} \bigg)^{1/2} \Bigg) \le  n^{-D}; \\
			& \displaystyle\min_{j \in \llbracket 1, n \rrbracket} \displaystyle\inf_{t \in [a, b]} \Bigg( \bigg| \mathsf{x}_j (t) - \displaystyle\frac{b - t}{b-a} \cdot u - \displaystyle\frac{t-a}{b-a} \cdot v \bigg| +  (8n)^{1/2} \bigg( \displaystyle\frac{(b-t)(t-a)}{b-a} \bigg)^{1/2} \Bigg) \ge - n^{-D}.
		\end{flalign*}
		\end{enumerate} 
	\end{lem}

	The following result from \cite{PLE} (upon applying the scale invariance of \Cref{scale}) indicating convergence of the top curves of the watermelon to the scaled parabolic Airy line ensemble. 
	
	\begin{lem}[{\cite[Theorem 3.1]{PLE}}]
		
		\label{convergea}
		
		Adopt the notation of \Cref{estimatexj}; assume $u = 0 = v$ and $(a, b) = (-Tn^{1/3}, Tn^{1/3})$; set $\sigma = T^{1/2}$; and define  
		\begin{flalign*} 
			\bm{\mathsf{X}}^n = ( \mathsf{X}_1^n, \mathsf{X}_2^n, \ldots , \mathsf{X}_n^n ) \in \llbracket 1, n \rrbracket \times \mathcal{C} \big( [-n^{1/3}, n^{1/3}] \big), \quad \text{where} \quad \mathsf{X}_j^n (t) = \sigma^{-1} \cdot \mathsf{x}_j^n (\sigma^2 t) - 2^{1/2} n^{2/3}. 
		\end{flalign*} 
		
		\noindent Then $\bm{\mathsf{X}}^n$ converges to $\bm{\mathcal{S}}$ on compact subsets of $\mathbb{Z}_{\ge 1} \times \mathbb{R}$, as $n$ tends to $\infty$.
	\end{lem} 
	
	The next lemma from \cite{BPLE} is a concentration bound for the $k$-th path of the rescaled parabolic Airy line ensemble, stating that it typically fluctuates by $\mathcal{O} (k^{-1/3})$ around a (deterministic) parabola. It was stated in \cite{BPLE} at $\sigma = 1$ and on the interval $s \in [0, t]$. That it also holds for arbitrary $\sigma > 0$ and on the interval $s \in [-t, t]$ follow from \eqref{sigmar} and the translation-invariance of $\bm{\mathcal{A}}$ (recall \Cref{translationa}), respectively.
	
	\begin{lem}[{\cite[Corollary 6.3]{BPLE}}]
		
		\label{kdeltad} 
		
		There exists a constant $c > 0$ such that the following holds. For any integer $k \ge 1$ and real numbers $t \ge 1$; $\sigma > 0$; and $u > c^{-1} \log (k+1)$, we have
		\begin{flalign*}
			\mathbb{P} \Bigg[ \displaystyle\sup_{s \in [-t, t]} \bigg| \mathcal{S}_k^{(\sigma)} (s) + 2^{-1/2} \sigma^3 s^2 + \sigma^{-1} 2^{-7/6} (3 \pi)^{2/3} k^{2/3} \bigg| \ge uk^{-1/3} \Bigg] \le c^{-1} t e^{-c \sigma u}.
		\end{flalign*}
		
	\end{lem}

	The next lemma provides upper and lower bounds for families of non-intersecting bridges whose $j$-th curve is of order $j^{2/3}$ (similarly to the parabolic Airy line ensemble, by \Cref{kdeltad}). It will be deduced through a comparison with the parabolic Airy line ensemble (as a consequence of \Cref{monotoneheight}, \Cref{estimatexj}, and \Cref{kdeltad}) in \Cref{ProofBridgeSum2} below.
	
	\begin{lem}\label{p:compareAiry}
		
		For any real numbers $A, B, d, D> 0$, there exist constants $c_1 = c_1 (d, D) > 1$ and $c_2 = c_2 (A, B) > 1$ such that the following holds. Fix an integer $n \ge 1$; a real number $M \in \mathbb{R}$; two $n$-tuples $\bm{u}, \bm{v} \in \overline{\mathbb{W}}_n$; an interval $[a, b] \in \mathbb{R}$ with $b-a \le n^D$; and a continuous function $f: [a, b] \rightarrow \overline{\mathbb{R}}$. Sample non-intersecting Brownian bridges $\bm{\mathsf{x}} = (\mathsf{x}_1, \mathsf{x}_2, \ldots , \mathsf{x}_n) \in \llbracket 1, n \rrbracket \times \mathcal{C} \big( [a, b] \big)$ from the measure $\mathsf{Q}_f^{\bm{u}; \bm{v}}$.  
		\begin{enumerate}
			\item Assume for each integer $j \in \llbracket 1, n \rrbracket$ and real number $t \in [a, b]$ that $\max \{ u_j, v_j \} \leq M - dj^{2/3}$ and $f(t)\leq M - d (n+1)^{2/3}$.  Then,     
			\begin{align}\label{e:xjupbound}
				\mathbb{P} \Bigg[ \bigcap_{j=1}^n \bigcap_{t \in [a, b]} \bigg\{ \sfx_j(t)\leq M+\frac{9\pi^2}{64d^3}(b-a)^2 - dj^{2/3} +2(\log n)^2 \bigg\} \Bigg] \ge 1 - c_1^{-1} e^{-c_1(\log n)^2}.
			\end{align}
			\item Assume that $b - a \le An^{1/3}$ and for each integer $j \in \llbracket 1, n \rrbracket $ that $\min \{ u_j, v_j \} \geq -B j^{2/3}-M$. Then, setting $A_0 = 2A^2 + B + 3$, we have
			\begin{align}\label{e:xjlowbound}
				\begin{aligned} 
					\mathbb{P} \Bigg[ \bigcap_{j=1}^n \bigcap_{t \in [a, b]} \bigg\{ \sfx_j(t) & \geq \frac{9\pi^2}{16 A_0^3} (t-a)(b-t) - M - 2 (\log n)^2 - A_0 j^{2/3}  \bigg\} \Bigg] \ge 1 - c_2^{-1} e^{-c_2 (\log n)^2}.
				\end{aligned} 
			\end{align}

		\end{enumerate}
	\end{lem}

\chapter{Gap Monotonicity and Likelihood of On-Scale Events}

\label{GAPSCALE}

	\section{Gap Monotonicity} 
	
	\label{Coupling1} 
	
	\subsection{Gap Couplings}
	
	\label{Coupling3}
	
	In this section we state monotone couplings for the gaps between the curves of non-intersecting Brownian Gibbsian line ensembles (that may have a lower boundary but no upper boundary). Throughout this section, for any integer $n \ge 1$ and real numbers $a < b$, we denote the entries of any $n$-tuple $\bm{w} \in \mathbb{R}^n$ by $\bm{w} = (w_1, w_2, \ldots , w_n)$ and of any line ensemble $\bm{\mathsf{y}} \in \llbracket 1, n \rrbracket \times \mathcal{C} \big( [a, b] \big)$ by $\bm{\mathsf{y}} = (\mathsf{y}_1 ,\mathsf{y}_2, \ldots , \mathsf{y}_n)$, unless stated otherwise.

	The next proposition states a variant of \Cref{monotoneheight} that provides monotone couplings for gaps $\mathsf{x}_j (t) - \mathsf{x}_{j+1} (t)$ between the curves in a line ensemble, instead of for the curves themselves. Instead of \eqref{fgxmonotone} we assume that the gaps between entries in $\bm{u}$ and $\bm{v}$ are bounded above by those in $\widetilde{\bm{u}}$ and $\widetilde{\bm{v}}$, respectively (see \eqref{uuvv}), and that $f$ is ``more concave'' than $\widetilde{f}$ (see \eqref{ff2}). We refer this result as \emph{gap monotonicity}; see the left side of \Cref{f:gap}. It is proven in \Cref{DifferenceDiscrete} below. In what follows we recall the measure $\mathsf{Q}$ prescribing non-intersecting Brownian bridges from \Cref{qxyfg}.
	
		\begin{prop}[Gap monotonicity]
		
		\label{monotonedifference}

		Fix an integer $n \ge 1$; four $n$-tuples $\bm{u}, \widetilde{\bm{u}}, \bm{v}, \widetilde{\bm{v}} \in \overline{\mathbb{W}}_n$; an interval $[a, b] \subset \mathbb{R}$; and continuous functions $f, \widetilde{f}: [a, b] \rightarrow \mathbb{R} \cup \{ -\infty \}$. Sample non-intersecting Brownian bridges $\bm{\mathsf{x}} (t)$ and $\widetilde{\bm{\mathsf{x}}} (t)$ from the measures $\mathsf{Q}_f^{\bm{u}; \bm{v}}$ and $\mathsf{Q}_{\tilde{f}}^{\tilde{\bm{u}}; \tilde{\bm{v}}}$, respectively. Assume
		\begin{flalign}
			\label{uuvv} 
			\begin{aligned} 
			& \qquad \qquad \quad 0 \le u_n - f(a) \le \widetilde{u}_n - \widetilde{f}(a); \quad \text{and} \quad 0 \le v_n - f(b) \le \widetilde{v}_n - \widetilde{f}(b); \\
			& u_j - u_{j+1} \le \widetilde{u}_j - \widetilde{u}_{j+1} \quad \text{and} \quad v_j - v_{j+1} \le \widetilde{v}_j - \widetilde{v}_{j+1}, \qquad \text{for each integer $j \in \llbracket 1, n-1 \rrbracket$}.
			\end{aligned} 
		\end{flalign}
	
		\noindent Moreover assume that the following two statements hold.
		\begin{enumerate} 
			\item Either we have $\widetilde{f} = -\infty$, or we have both $f > -\infty$ and $\widetilde{f} > -\infty$.
			\item The difference $\widetilde{f} - f$ is convex, that is, for any $s, t \in [a, b]$ and $r \in [0, 1]$ we have
		\begin{flalign}
			\label{ff2}
			& r \cdot f(s) -  f \big( rs + (1-r) t\big) + (1-r) \cdot f(t) \le r \cdot \widetilde{f} (s) - \widetilde{f} \big( rs + (1-r) t \big) + (1-r) \cdot \widetilde{f}(t).
		\end{flalign}	
	
		\end{enumerate}
		
		\noindent Then, there exists a coupling between $\bm{\mathsf{x}} (t)$ and $\widetilde{\bm{\mathsf{x}}} (t)$ such that $\mathsf{x}_n (t) - f(t) \le\widetilde{\mathsf{x}}_n (t) - \widetilde{f} (t)$ and $\mathsf{x}_j (t) - \mathsf{x}_{j+1} (t) \le\widetilde{\mathsf{x}}_j (t) -\widetilde{\mathsf{x}}_{j+1} (t)$, for each real number $t \in [a, b]$ and integer $j \in \llbracket 1, n-1 \rrbracket$.

	\end{prop}
	
					\begin{figure}
	\center
\includegraphics[scale = .65]{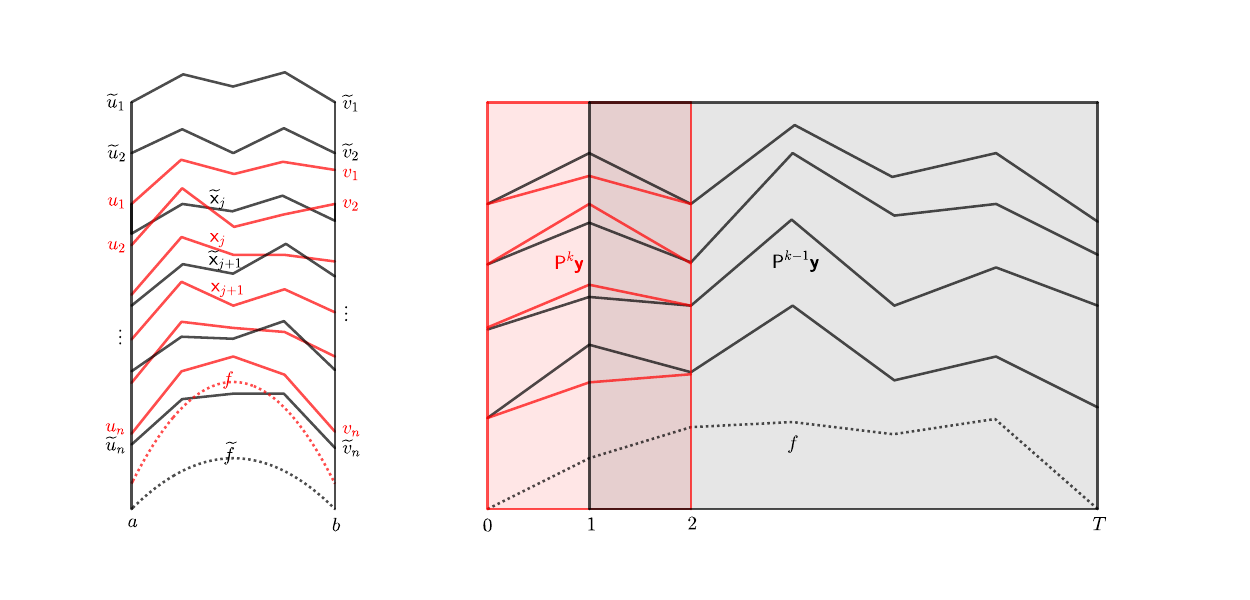}

\caption{Shown to the left is the gap monotonicity result, \Cref{monotonedifference}. Shown to the right are the alternating Markov dynamics from \Cref{dynamicalternating}, which alternate between resampling the Gaussian bridges in the red and gray boxes. }
\label{f:gap}
	\end{figure}

	\subsection{Semi-discrete Gap Monotonicity}
	
	\label{DifferenceDiscrete} 
	
	In this section we reduce \Cref{monotonedifference} to a semi-discrete analog of it, in which Brownian bridges are replaced by Gaussian ones. To explain this, for any integer $T \ge 1$, a \emph{($T$-step) Gaussian walk} starting at $u \in \mathbb{R}$ is a probability measure on $(T+1)$-tuples $\big( x(0), x(1), \ldots , x(T) \big) \in \mathbb{R}^{T+1}$ with $x(0) = u$ such that, for each $j \in \llbracket 1, T \rrbracket$, the jump $x(j) - x(j-1)$ is a centered Gaussian random variable of variance $1$. A \emph{($T$-step) Gaussian bridge} from $u$ to $v$ is a Gaussian walk starting at $u$, conditioned to end at $v$ (that is, $x(T) = v$). The following definition is similar to \Cref{qxyfg} and provides notation for non-intersecting Gaussian bridges.

	\begin{definition}
		
		\label{gxyfg}
		
		Fix integers $T, n \ge 1$; two $n$-tuples $\bm{u}, \bm{v} \in \overline{\mathbb{W}}_n$; and two functions $f, g: \llbracket 0, T \rrbracket \rightarrow \overline{\mathbb{R}}$ such that $f < g$, $f < \infty$, and $g > -\infty$. Let $\mathsf{G}_{f; g}^{\bm{u}; \bm{v}}$ denote the law on sequences $\bm{\mathsf{x}} (t) = \big( \mathsf{x}_1 (t), \mathsf{x}_2 (t), \ldots , \mathsf{x}_n (t) \big)$, with $t \in \llbracket 0, T \rrbracket$, given by $n$ independent $T$-step Gaussian bridges, conditioned on satisfying $\bm{\mathsf{x}} (t) \in \mathbb{W}_n$ for each $t \in \llbracket 0, T-1 \rrbracket$; $\mathsf{x}_j (0) = u_j$ and $\mathsf{x}_j (T) = v_j$ for each $j \in \llbracket 1, n \rrbracket$; and $f \le \mathsf{x}_j \le g$ for each $j \in \llbracket 1, n \rrbracket$. If $g = \infty$, then we abbreviate $\mathsf{G}_f^{\bm{u}; \bm{v}} = \mathsf{G}_{f; \infty}^{\bm{u}; \bm{v}}$. It is assumed here that $f(0) \le u_n \le u_1 \le g(0)$ and $f(T) \le v_n \le v_1 \le g(T)$, and that $f < g$, even when not stated explicitly.\index{G@$\mathsf{G}_{f;g}^{\bm{u}; \bm{v}}$}
		
	\end{definition}
	
	\begin{rem}
		\label{discretelinear}
		
		As in \Cref{linear}, non-intersecting Gaussian bridges satisfy the following useful invariance property under affine transformations. Adopt the notation of \Cref{gxyfg}, and fix real numbers $\alpha, \beta \in \mathbb{R}$. Define the $n$-tuples $\bm{u}', \bm{v}' \in \overline{\mathbb{W}}_n$ and functions $f', g' : \llbracket 0, n \rrbracket \rightarrow \overline{\mathbb{R}}$ by setting 
		\begin{flalign*} 
			& u_j' = u_j + \alpha, \quad \text{and} \quad v_j' = v_j + T \beta + \alpha, \qquad \qquad \qquad \quad \text{for each $j \in \llbracket 0, n \rrbracket$}; \\
			& f'(t) = f(t) + t \beta + \alpha, \quad \text{and} \quad g'(t) = g(t) + t \beta + \alpha, \qquad \text{for each $t \in \llbracket 0, T \rrbracket$}.
		\end{flalign*} 
		
		\noindent Sampling $\bm{\mathsf{x}}' = \big( \mathsf{x}_1' (t), \mathsf{x}_2' (t), \ldots , \mathsf{x}_n' (t) \big)$ under $\mathsf{G}_{f', g'}^{\bm{u}'; \bm{v}'}$, there is a coupling between $\bm{\mathsf{x}}'$ and $\bm{\mathsf{x}}$ such that $\mathsf{x}_j' (t) = \mathsf{x}_j (t) + \beta t + \alpha$ for each $t \in \llbracket 0, T \rrbracket$ and $j \in \llbracket 1, n \rrbracket$. 
		
		Indeed, this follows from the analogous affine invariance of a single Gaussian random bridge, together with the fact that affine transformations do not affect the non-intersecting property. More specifically, if $\big( x(t) \big)$ is a $T$-step Gaussian random walk from some $u \in \mathbb{R}$ to some $v \in \mathbb{R}$ then $\big( x(t) + t \beta + \alpha \big)$ is a $T$-step Gaussian random walk from $u + \alpha$ to $v + T \beta + \alpha$, and $\bm{\mathsf{x}} (t) \in \mathbb{W}_n$ if and only if $\bm{\mathsf{x}}(t) + t \beta + \alpha \in \mathbb{W}_n$ (for any $t \in \llbracket 0, T \rrbracket$ and $\alpha, \beta \in \mathbb{R}$).

	\end{rem}
	
	The following lemma from \cite{U} states a version of height monotonicity (the analog of \Cref{uvv}) for non-intersecting Gaussian bridges. 
	
	\begin{lem}[{\cite[Lemma B.6]{U}}]
		
		\label{fgk} 
		
		Fix integers $T, n \ge 1$; a real number $B \ge 0$; four $n$-tuples $\bm{u}, \widetilde{\bm{u}}, \bm{v}, \widetilde{\bm{v}} \in \overline{\mathbb{W}}_n$; and functions $f, \widetilde{f}, g, \widetilde{g}: \llbracket 0, T \rrbracket \rightarrow \overline{\mathbb{R}}$. Sample non-intersecting Gaussian bridges $\bm{\mathsf{x}} (t)$ and $\widetilde{\bm{\mathsf{x}}} (t)$ from the measures $\mathsf{G}_{f; g}^{\bm{u}; \bm{v}}$ and $\mathsf{G}_{\tilde{f}; \tilde{g}}^{\tilde{\bm{u}}; \tilde{\bm{v}}}$, respectively. If $u_j \le \widetilde{u}_j \le u_j + B$ and $v_j \le \widetilde{v}_j \le v_j + B$ for all $j \in \llbracket 1, n \rrbracket$, then the following two statements hold.

		\begin{enumerate} 
			\item If $f(t) \le \widetilde{f}(t) \le f(t) + B$ and $g(t) \le \widetilde{g}(t) \le g(t) + B$ for each $t \in \llbracket 0, T \rrbracket$, then there is a coupling between $\bm{\mathsf{x}}$ and $\widetilde{\bm{\mathsf{x}}}$ so that $\mathsf{x}_j (t) \le\widetilde{\mathsf{x}}_j (t) \le \mathsf{x}_j (t) + B$ for each $(t, j) \in \llbracket 0, T \rrbracket \times \llbracket 1, n \rrbracket$. 
			\item If $\bm{u} = \widetilde{\bm{u}}$, and $f(t) \le \widetilde{f} (t) \le f(t) + tT^{-1} B$ and $g(t) \le \widetilde{g}(t) \le g(t) + tT^{-1} B$ for each $t \in \llbracket 0, T \rrbracket$, then there is a coupling between $\bm{\mathsf{x}}$ and $\widetilde{\bm{\mathsf{x}}}$ so that $\mathsf{x}_j (t) \le\widetilde{\mathsf{x}}_j (t) \le \mathsf{x}_j (t) + tT^{-1} B$ for each $(t, j) \in \llbracket 0, T \rrbracket \times \llbracket 1, n \rrbracket$.
		\end{enumerate} 
		
	\end{lem}

	Stated next is an analog of \Cref{monotonedifference} for Gaussian bridges; its proof is in \Cref{2T} below.
	
	\begin{prop}
		
		\label{monotonedifferencediscrete}

		Fix integers $T, n \ge 1$; four $n$-tuples $\bm{u}, \widetilde{\bm{u}}, \bm{v}, \widetilde{\bm{v}} \in \overline{\mathbb{W}}_n$; and functions $f, \widetilde{f}: \llbracket 0, T\rrbracket \rightarrow \overline{\mathbb{R}}$. Sample non-intersecting Gaussian bridges $\bm{\mathsf{x}} (t)$ and $\widetilde{\bm{\mathsf{x}}} (t)$ from the measures $\mathsf{G}_f^{\bm{u}; \bm{v}}$ and $\mathsf{G}_{\tilde{f}}^{\tilde{\bm{u}}; \tilde{\bm{v}}}$, respectively. Assume that
		\begin{flalign}
			\label{y1} 
			\begin{aligned}
				& u_n - f(0) \le \widetilde{u}_n - \widetilde{f}(0), \quad \text{and} \quad v_n - f(T) \le \widetilde{v}_n - \widetilde{f}(T); \\
				& u_j - u_{j+1} \le \widetilde{u}_j - \widetilde{u}_{j+1} \quad \text{and} \quad v_j - v_{j+1} \le \widetilde{v}_j - \widetilde{v}_{j+1}, \qquad \qquad \text{for each $j \in \llbracket 1, n-1 \rrbracket$}.
			\end{aligned}
		\end{flalign}
			
		\noindent Moreover assume that we have $\widetilde{f} = -\infty$, or that we have $f > -\infty$, $\widetilde{f} > - \infty$, and
		\begin{flalign} 
		\label{y200} 
				f(t+1) - 2 f (t) + f(t-1) \le \widetilde{f} (t+1) - 2 \widetilde{f}(t) + \widetilde{f}(t-1), \qquad \text{for each $t \in \llbracket 1, T-1\rrbracket$}.
		\end{flalign}	
		
		\noindent Then, there exists a coupling between $\bm{\mathsf{x}} (t)$ and $\widetilde{\bm{\mathsf{x}}} (t)$ such that $\mathsf{x}_n (t) - f(t) \le\widetilde{\mathsf{x}}_n (t) - \widetilde{f} (t)$ and $\mathsf{x}_j (t) - \mathsf{x}_{j+1} (t) \le\widetilde{\mathsf{x}}_j (t) -\widetilde{\mathsf{x}}_{j+1} (t)$, for each $t \in \llbracket 0, T \rrbracket$ and $j \in \llbracket 1, n-1 \rrbracket$. 
		
	\end{prop}

	\begin{rem}
		\label{discrete1} 
		
		Unlike for \Cref{monotoneheight}, the fully discrete variant of \Cref{monotonedifference} obtained by replacing Gaussian bridges with Bernoulli random bridges, with jumps in $\{ -1 , 1 \}$, is false (which can eventually be attributed to the fact that the latter does not satisfy the affine invariance from \Cref{discretelinear}). Indeed, consider two pairs of non-intersecting Bernoulli random bridges $\bm{\mathsf{x}} = \big( \bm{\mathsf{x}}(0), \bm{\mathsf{x}}(1), \bm{\mathsf{x}}(2) \big)$ and $\widetilde{\bm{\mathsf{x}}} = \big(\widetilde{\mathsf{x}}(0),\widetilde{\mathsf{x}}(1),\widetilde{\mathsf{x}}(2) \big)$ on the interval $\llbracket 0, 2 \rrbracket$, both with infinite lower boundary $f = -\infty$; the first has starting points $(u_1, u_2) = (2, 0)$ and ending points $(v_1, v_2) = (4, 2)$, while the second has starting points $(\widetilde{u}_1, \widetilde{u}_2) = (3, 0)$ and ending points $(\widetilde{v}_1, \widetilde{v}_2) = (3, 0)$. 
		
		The analog of \Cref{monotonedifferencediscrete} would have suggested the existence of a coupling between $\bm{\mathsf{x}}$ and $\widetilde{\bm{\mathsf{x}}}$ such that $\mathsf{x}_1 (1) - \mathsf{x}_2 (1) \le\widetilde{\mathsf{x}}_1 (1) -\widetilde{\mathsf{x}}_2 (1)$. However, this is not possible. Indeed, the starting and ending data for $\bm{\mathsf{x}}$ deterministically imposes $\big( \mathsf{x}_1 (1), \mathsf{x}_2 (1) \big) = (3, 1)$, so that $\mathsf{x}_1 (1) - \mathsf{x}_2 (1) = 2$. On the other hand, for $\widetilde{\bm{\mathsf{x}}}$, we have $\big(\widetilde{\mathsf{x}}_1 (1),\widetilde{\mathsf{x}}_2 (1) \big) \in \big\{ (4, 1), (4, -1), (2, 1), (2, -1)\big\}$ each with probability $1 / 4$, so that $\widetilde{\mathsf{x}}_1 (1) -\widetilde{\mathsf{x}}_2 (1) = 1$ occurs with probability $1 / 2$. 
		
	\end{rem} 
	
	Given \Cref{monotonedifferencediscrete}, we can establish \Cref{monotonedifference}.

	\begin{proof}[Proof of \Cref{monotonedifference}]
		
		First assume that $\bm{u}, \bm{v}, \widetilde{\bm{u}}, \widetilde{\bm{v}} \in \mathbb{W}_n$. For each integer $T \ge 0$, define the $n$-tuples $\bm{u}^{(T)}, \bm{v}^{(T)} \in \mathbb{W}_n$ by setting $u_j^{(T)} = T^{1/2} u_j$ and $v_j^{(T)} = T^{1/2} v_j$, for each $j \in \llbracket 1, n \rrbracket$. Further define the functions $f^{(T)}, \widetilde{f}^{(T)} : \llbracket 0, T \rrbracket \rightarrow \overline{\mathbb{R}}$ by, for each $t \in \llbracket 0, T \rrbracket$, setting 
		\begin{flalign*}
			f^{(T)} (t) = T^{1/2} \cdot f \bigg( \frac{(T-t)a}{T} + \frac{tb}{T} \bigg); \qquad  \widetilde{f}^{(T)} (t) = T^{1/2} \cdot \widetilde{f} \bigg( \displaystyle\frac{(T-t)a}{T} + \displaystyle\frac{tb}{T} \bigg).
		\end{flalign*} 
		
		\noindent Sample the two families of non-intersecting Gaussian bridges $\bm{\mathsf{x}}^{(T)} = \big( \mathsf{x}_1^{(T)}, \mathsf{x}_2^{(T)}, \ldots , \mathsf{x}_n^{(T)} \big)$ and $\widetilde{\bm{\mathsf{x}}}^{(T)} = \big(\widetilde{\mathsf{x}}_1^{(T)},\widetilde{\mathsf{x}}_2^{(T)}, \ldots ,\widetilde{\mathsf{x}}_n^{(T)} \big)$ from the measures $\mathsf{G}_{f^{(T)}}^{\bm{u}^{(T)}; \bm{v}^{(T)}}$ and $\mathsf{G}_{\widetilde{f}^{(T)}}^{\widetilde{\bm{u}}^{(T)}; \widetilde{\bm{v}}^{(T)}}$, respectively. As $T$ tends to $\infty$, the joint laws over $t$ of $\big( T^{-1/2} \bm{\mathsf{x}}^{(T)} ( \lfloor \frac{(t-a)T}{b-a} \rfloor )\big)$ and $\big( T^{-1/2} \widetilde{\bm{\mathsf{x}}}^{(T)} (\lfloor \frac{(t-a)T}{b-a} \rfloor) \big)$ converge to those $\bm{\mathsf{x}} (t)$ and $\widetilde{\bm{\mathsf{x}}} (t)$, respectively. The latter follows from the convergence of an individual Gaussian bridge to a Brownian bridge,\footnote{Indeed, the joint marginals for $t \in \big\{ a + (b-a)i/T \big\}_{i \in \llbracket 0, T \rrbracket}$ of a Brownian bridge $B(t)$ on $[a, b]$ is given by a Gaussian random bridge.} together with the fact that non-intersection event we are conditioning has positive probability (as $\bm{u}, \bm{v}, \widetilde{\bm{u}}, \widetilde{\bm{v}} \in \mathbb{W}_n$).
			
		Next, by \Cref{monotonedifferencediscrete}, there exists a coupling between $\bm{\mathsf{x}}^{(T)}$ and $\widetilde{\bm{\mathsf{x}}}^{(T)}$ such that, for each $j \in \llbracket 1, n-1\rrbracket$ and $t \in \llbracket 0, T \rrbracket$, we have 
		\begin{flalign*} 
			& T^{-1/2} \big( \mathsf{x}_n^{(T)} (t) - f^{(T)} (t) \big) \le T^{-1/2} \big( \widetilde{\mathsf{x}}_n^{(T)} (t) - \widetilde{f}^{(T)} (t)\big); \\ 
			& T^{-1/2} \big( \mathsf{x}_j^{(T)} (t) - \mathsf{x}_{j+1}^{(T)} (t) \big) \le T^{-1/2} \big( \widetilde{\mathsf{x}}_j^{(T)} (t) - \widetilde{\mathsf{x}}_{j+1}^{(T)} (t) \big).
		\end{flalign*} 
		
		\noindent Taking any limit point of these couplings as $T$ tends to $\infty$ (this sequence of couplings is tight, since their marginals are) yields a coupling between $\bm{\mathsf{x}}$ and $\widetilde{\bm{\mathsf{x}}}$ such that $\mathsf{x}_n (t) - f(t) \le \widetilde{\mathsf{x}}_n (t) - \widetilde{f} (t)$ and $\mathsf{x}_j (t) - \mathsf{x}_{j+1} (t) \le\widetilde{\mathsf{x}}_j (t) -\widetilde{\mathsf{x}}_{j+1} (t)$, for each $t \in [a, b]$ and $j \in \llbracket 1, n-1 \rrbracket$.
		
		Now assume that at least one of $\bm{u}$, $\bm{v}$, $\widetilde{\bm{u}}$, or $\widetilde{\bm{v}}$ is in $\overline{\mathbb{W}}_n \setminus \mathbb{W}_n$.  Then, for any $\varepsilon > 0$, define $\bm{u}^{\varepsilon}, \bm{v}^{\varepsilon} \in \mathbb{W}_n$ by setting $u_j^{\varepsilon} = u_j - j \varepsilon$ and $v_j^{\varepsilon} = v_j - j \varepsilon$ for each $j \in \llbracket 1, n \rrbracket$; also define $f^{\varepsilon}, g^{\varepsilon} : [a, b] \rightarrow \mathbb{R}$ by setting $f^{\varepsilon} (s) = f(s) - (n+1) \varepsilon$ and $g^{\varepsilon} (s) = g(s)$, for each $s \in [a, b]$. Then $\bm{u}^{\varepsilon}, \bm{v}^{\varepsilon}, \widetilde{\bm{u}}^{\varepsilon}, \widetilde{\bm{v}}^{\varepsilon} \in \mathbb{W}_n$. Thus, sampling non-intersecting Brownian bridges $\bm{\mathsf{x}}^{\varepsilon}$ and $\widetilde{\bm{\mathsf{x}}}^{\varepsilon}$ from the measures $\mathsf{Q}_{f^{\varepsilon}}^{\bm{u}^{\varepsilon}; \bm{v}^{\varepsilon}}$ and $\mathsf{Q}_{\tilde{f}^{\varepsilon}}^{\tilde{\bm{u}}^{\varepsilon}; \tilde{\bm{v}}^{\varepsilon}}$, it follows from the above that there exists a coupling between $(\bm{\mathsf{x}}^{\varepsilon}, \widetilde{\bm{\mathsf{x}}}^{\varepsilon})$ such that $\mathsf{x}_n^{\varepsilon} (t) - f^{\varepsilon} (t) \le \widetilde{\mathsf{x}}_n^{\varepsilon} (t) - \widetilde{f}^{\varepsilon} (t)$ for each $t \in [a, b]$ and $\mathsf{x}_j^{\varepsilon} (t) - \mathsf{x}_{j+1}^{\varepsilon} (t) \le \widetilde{\mathsf{x}}_j^{\varepsilon} (t) - \widetilde{\mathsf{x}}_{j+1}^{\varepsilon} (t)$, for each $(j, t) \in \llbracket 1, n - 1 \rrbracket \times [a,b]$. Moreover, by (the $B = (n+1) \varepsilon$ case of) \Cref{uvv}, there exist couplings between $(\bm{\mathsf{x}}, \bm{\mathsf{x}}^{\varepsilon})$ and $(\widetilde{\bm{\mathsf{x}}}, \widetilde{\bm{\mathsf{x}}}^{\varepsilon})$ such that $\big| \mathsf{x}_j (t) - \mathsf{x}_j^{\varepsilon} (t) \big| \le (n+1) \varepsilon$ and $\big| \widetilde{\mathsf{x}}_j (t) - \widetilde{\mathsf{x}}_j^{\varepsilon} (t) \big| \le (n+1) \varepsilon$.
		
		Therefore, for any $\varepsilon > 0$, this induces a coupling between $(\bm{\mathsf{x}}, \widetilde{\bm{\mathsf{x}}})$ such that  
		\begin{flalign*}
			 \mathsf{x}_n (s) - f(s) \le \mathsf{x}_n^{\varepsilon} (s) - f^{\varepsilon} (s) + 2(n+1) \varepsilon \le \widetilde{\mathsf{x}}_n^{\varepsilon} (s) - \widetilde{f}^{\varepsilon} (s) + 2 (n+1) \varepsilon \le \widetilde{\mathsf{x}}_n (s) - \widetilde{f} (s) + 4(n+1) \varepsilon,
		\end{flalign*} 
	
		\noindent and 
		\begin{flalign*} 
			\mathsf{x}_j (s) - \mathsf{x}_{j+1} (s) & \le \mathsf{x}_j^{\varepsilon} (s) - \mathsf{x}_{j+1}^{\varepsilon} (s) + 2 (n+1) \varepsilon \\
			& \le \widetilde{\mathsf{x}}_j^{\varepsilon} (s) - \widetilde{\mathsf{x}}_{j+1}^{\varepsilon} (s) + 2(n+1) \varepsilon \le \widetilde{\mathsf{x}}_j (s) - \widetilde{\mathsf{x}}_{j+1} (s) + 4(n+1) \varepsilon,
		\end{flalign*}
		
		\noindent for each $s \in [a, b]$ and $j \in \llbracket 1, n-1 \rrbracket$. Taking any limit point of these couplings as $\varepsilon$ tends to $0$ (which form a tight sequence, since their marginals are tight) then yields the proposition, due to the convergence\footnote{See \Cref{boundaryconvergeensemble} below.} of $(\bm{\mathsf{x}}^{\varepsilon}, \widetilde{\bm{\mathsf{x}}}^{\varepsilon})$ to $(\bm{\mathsf{x}}, \widetilde{\bm{\mathsf{x}}})$ as $\varepsilon$ tends to $0$.
	\end{proof}

	\subsection{Reduction to the Case $T=2$}
	
	\label{2T} 
	
	The height monotonicity result \Cref{monotoneheight} was shown in \cite{PLE} by verifying that monotone couplings were preserved under certain local Markov (Glauber) dynamics. That proof using those dynamics does not seem to apply for gap monotonicity, but in this section we will use a (less local) Markov dynamic to establish \Cref{monotonedifferencediscrete}, assuming the following result stating it holds when $T = 2$. It will be established in \Cref{Proof2T} below.
	
	\begin{prop}
		
		\label{monotonedifferencediscrete2}
		
		If $T = 2$, then \Cref{monotonedifferencediscrete} holds. 
	\end{prop}
	
	The Markov dynamics we use are a semi-discrete analog of those introduced in \cite[Definition 4.5]{ESTP}; they are given by repeatedly alternating between resampling the Gaussian bridges on $t = 1$ (conditional on their values at $t \ne 1$) and on $t \in \llbracket 2, T-1]$ (conditional on their value at $t = 1$). See the right side of \Cref{f:gap}.

	\begin{definition}
		
		\label{dynamicalternating}
		
		Fix integers $T, n \ge 1$ and a function $f : \llbracket 0, T \rrbracket \rightarrow \overline{\mathbb{R}}$. For $t \in \llbracket 0, T \rrbracket$, let $\bm{\mathsf{y}} (t) = \big( \mathsf{y}_1 (t), \mathsf{y}_2 (t), \ldots , \mathsf{y}_n (t) \big) \in \overline{\mathbb{W}}_n$ be a family of $n$ non-intersecting paths of length $T + 1$. The \emph{alternating dynamics} is the discrete-time Markov chain\footnote{We may identify the state space of this Markov chain by $\mathbb{W}_n^{T-1}$, as $\mathsf{P}^k \bm{\mathsf{y}} (t)$ can be arbitrary elements of $\mathbb{W}_n$ for $t \in \llbracket 1, T-1 \rrbracket$ but must satisfy $\mathsf{P}^k \bm{\mathsf{y}} (t) = \bm{\mathsf{y}} (t)$ for $t \in \{ 0, T \}$.} whose state $\mathsf{P}^k \bm{\mathsf{y}} (t) = \big( \mathsf{P}^k \mathsf{y}_1 (t), \mathsf{P}^k \mathsf{y}_2 (t), \ldots \mathsf{P}^k \mathsf{y}_n (t) \big)$ at time $k \ge 0$ is determined as follows. If $k = 0$, set $\mathsf{P}^k \bm{\mathsf{y}} = \bm{\mathsf{y}}$. For $k \ge 1$, sample $\mathsf{P}^k \bm{\mathsf{y}}$ inductively as below; throughout, we set $\bm{\mathsf{y}}' = \mathsf{P}^{k-1} \bm{\mathsf{y}}$. 
		\begin{enumerate}
			\item If $k$ is odd, set $\mathsf{P}^k \mathsf{y}_j (t) = \mathsf{y}_j' (t)$ for each $t \in \llbracket 2, T \rrbracket$. For $t =1$, sample $\big( \mathsf{P}^k \mathsf{y}_j (t) \big)_{t \in \llbracket 0, 2 \rrbracket}$ as $2$-step non-intersecting Gaussian bridges under the measure $\mathsf{G}_{f |_{\llbracket 0, 2 \rrbracket}}^{\bm{\mathsf{y}}'(0); \bm{\mathsf{y}}'(2)}$.
			
			\item If $k$ is even, set $\mathsf{P}^k \mathsf{y}_j (t) = \mathsf{y}_j' (t)$ for each $t \in \llbracket 0, 1 \rrbracket$. For $t \in \llbracket 2, T \rrbracket$, sample $\big( \mathsf{P}^k \mathsf{y}_j (t) \big)_{t \in \llbracket 1, T \rrbracket}$ as $(T-1)$-step non-intersecting Gaussian bridges under the measure $\mathsf{G}_{f|_{\llbracket 1, T \rrbracket}}^{\bm{\mathsf{y}}'(1); \bm{\mathsf{y}}'(T)}$.\index{P@$\mathsf{P}$; alternating dynamics}
		\end{enumerate}
		
	\end{definition}
	
	\begin{rem}
		
		\label{gmeasurep}
		
		It follows from the Gibbs property (for non-intersecting Gaussian bridges) that $\mathsf{G}_f^{\bm{\mathsf{y}} (0); \bm{\mathsf{y}}(T)}$ is a stationary measure for the alternating dynamics.
		
	\end{rem}
	
	The following lemma states that the alternating dynamics converge to the measure $\mathsf{G}_f^{\bm{\mathsf{y}}(0), \bm{\mathsf{y}}(T)}$; its proof is given in \Cref{DynamicConverge} below as a consequence of a convergence theorem for Harris chains. In what follows, for any two probability measures $\nu_1$ and $\nu_2$, on a measurable space $\Omega$ with $\sigma$-algebra $\mathcal{F}$, we recall that the total variation distance between them is defined by
	\begin{flalign*} 
			d_{\TV} (\nu_1, \nu_2) = \displaystyle\sup_{A \in \mathcal{F}} \big| \nu_1 (A) - \nu_2 (A) \big|.
	\end{flalign*} 
	
	\begin{lem}
		
		\label{pkxtconverge} 
		
		Adopting the notation of \Cref{dynamicalternating}, the law of $\mathsf{P}^{2k} \bm{\mathsf{y}}$ converges as $k$ tends to $\infty$ to $\mathsf{G}_f^{\bm{\mathsf{y}}(0); \bm{\mathsf{y}}(T)}$, under the total variational distance norm.
		
	\end{lem} 
	
	Given \Cref{monotonedifferencediscrete2} and \Cref{pkxtconverge}, we can establish \Cref{monotonedifferencediscrete}.

	\begin{proof}[Proof of \Cref{monotonedifferencediscrete}]
		
		First observe that the proposition holds for $T \in \{ 1, 2 \}$. Indeed, if $T = 1$ then $\bm{\mathsf{x}} (t)$ and $\widetilde{\bm{\mathsf{x}}} (t)$ are (deterministically) fixed by $\bm{u}$, $\widetilde{\bm{u}}$, $\bm{v}$, and $\widetilde{\bm{v}}$, and \Cref{monotonedifferencediscrete2} indicates that the result holds for $T = 2$. Thus, let us verify it for $T > 2$ by induction on $T$. 
		
		Fix sequences of non-intersecting $T$-step walks $\bm{\mathsf{y}}(t) = \big( \mathsf{y}_1 (t), \mathsf{y}_2 (t), \ldots , \mathsf{y}_n (t) \big) \in \overline{\mathbb{W}}_n$ and $\widetilde{\bm{\mathsf{y}}}(t) = \big( \widetilde{\mathsf{y}}_1 (t), \widetilde{\mathsf{y}}_2 (t), \ldots , \widetilde{\mathsf{y}}_n (t) \big)$ such that for each $t \in \llbracket 0, T \rrbracket$ and $j \in \llbracket 1, n \rrbracket$ we have
		\begin{flalign}
			\label{y0} 
			\begin{aligned}
				& \mathsf{y}_j (0) = u_j; \qquad \mathsf{y}_j (T) = v_j; \qquad \widetilde{\mathsf{y}}_j (0) = \widetilde{u}_j; \qquad \widetilde{\mathsf{y}}_j (T) = \widetilde{v}_j; \\ 
				& \mathsf{y}_n (t) - f(t) \le \widetilde{\mathsf{y}}_n (t) - \widetilde{f} (t); \qquad \mathsf{y}_j (t) - \mathsf{y}_{j+1} (t) \le \widetilde{\mathsf{y}}_j (t) - \widetilde{\mathsf{y}}_{j+1} (t).
			\end{aligned} 
		\end{flalign} 
		
		\noindent Such $\bm{\mathsf{y}}$ and $\widetilde{\bm{\mathsf{y}}}$ are guaranteed to exist by \eqref{y1}.
				
		Applying the alternating dynamics $\mathsf{P}$ to $\bm{\mathsf{y}}$ and $\widetilde{\bm{\mathsf{y}}}$, we claim it is possible to couple $\mathsf{P}^k \bm{\mathsf{y}}$ and $\mathsf{P}^k \widetilde{\bm{\mathsf{y}}}$ in such a way that
		\begin{flalign} 
			\label{kpy} 
			\mathsf{P}^k \mathsf{y}_n (t) - f(t) \le \mathsf{P}^k \widetilde{\mathsf{y}}_n (t) - \widetilde{f} (t); \qquad  \mathsf{P}^k \mathsf{y}_j (t) - \mathsf{P}^k \mathsf{y}_{j+1} (t) \le \mathsf{P}^k \widetilde{\mathsf{y}}_j (t) - \mathsf{P}^k \widetilde{\mathsf{y}}_{j+1} (t),
		\end{flalign}
		
		\noindent for each $k \in \mathbb{Z}_{\ge 0}$, $t \in \llbracket 0, T \rrbracket$, and $j \in \llbracket 1, n \rrbracket$. This follows by induction on $k$. Indeed, the statement is true by \eqref{y0} at $k=0$, so let $k \ge 1$ and assume that 
		\begin{flalign} 
			\label{kpy1} 
			\mathsf{P}^{k-1} \mathsf{y}_n (t) - f(t) \le \mathsf{P}^{k-1} \widetilde{\mathsf{y}}_n (t) - \widetilde{f} (t); \qquad  \mathsf{P}^{k-1} \mathsf{y}_j (t) - \mathsf{P}^{k-1} \mathsf{y}_{j+1} (t) \le \mathsf{P}^{k-1} \widetilde{\mathsf{y}}_j (t) - \mathsf{P}^{k-1} \widetilde{\mathsf{y}}_{j+1} (t),
		\end{flalign}
	
		\noindent By \Cref{dynamicalternating}, if $k$ is odd, then $\mathsf{P}^k \bm{\mathsf{y}}$ is obtained from $\mathsf{P}^{k-1} \bm{\mathsf{y}}$ by sampling $\big( \mathsf{P}^k \mathsf{y}_j (t) \big)_{t \in \llbracket 0, 2 \rrbracket}$ according to $\mathsf{G}_{f|_{\llbracket 0, 2 \rrbracket}}^{\mathsf{P}^{k-1} \bm{\mathsf{y}}(0); \mathsf{P}^{k-1} \bm{\mathsf{y}}(2)}$ (and leaving $\mathsf{P}^{k-1} \mathsf{y}_j (t)$ fixed for $t \notin \llbracket 0, 2 \rrbracket$). If $k$ is even, then it is obtained by sampling $\big( \mathsf{P}^k \mathsf{y}_j (t) \big)_{t \in \llbracket 1, T \rrbracket}$ according to $\mathsf{G}_{f \in \llbracket 1, T \rrbracket}^{\mathsf{P}^{k-1} \bm{\mathsf{y}}(1); \mathsf{P}^{k-1} (T)}$ (and leaving $\mathsf{P}_{k-1} \mathsf{y}_j (t)$ fixed for $t \in \{ 0, 1 \}$).   
		
		If $k$ is even then, by \eqref{kpy1} and the inductive hypothesis in $T$, it is possible to couple $(\mathsf{P}^k \bm{\mathsf{y}}, \mathsf{P}^k \widetilde{\bm{\mathsf{y}}})$ under $\big( \mathsf{G}_{f|_{\llbracket 1, T \rrbracket}}^{\mathsf{P}^{k-1} \bm{\mathsf{y}} (1); \mathsf{P}^{k-1} \bm{\mathsf{y}} (T)}, \mathsf{G}_{\widetilde{f}|_{\llbracket 1, T \rrbracket}}^{\mathsf{P}^{k-1} \widetilde{\bm{\mathsf{y}}} (1); \mathsf{P}^{k-1} \widetilde{\bm{\mathsf{y}}} (T)} \big)$ (leaving all $\mathsf{y}_j (t)$ and $\widetilde{\mathsf{y}}_j (t)$ for $j \notin \llbracket 2, T\rrbracket$ fixed) such that \eqref{kpy} holds for all $t \in \llbracket 1, T \rrbracket$ almost surely; it also holds for $t = 0$, by our assumption \eqref{y1} on $\bm{u}$. The case when $k$ is odd is addressed entirely analogously, thereby verifying \eqref{kpy}.
		
		Now take any limit point, over even integers $k$ tending to $\infty$, of the coupling between $(\mathsf{P}^k \bm{\mathsf{y}}, \mathsf{P}^k \widetilde{\bm{\mathsf{y}}})$ guaranteeing \eqref{kpy}. Then applying \Cref{pkxtconverge} (to run the dynamics until they mix) gives the proposition.
	\end{proof}

\subsection{The Equal Boundary Case}	
	\label{Equaluv} 
	
	In this section we establish the following variant of \Cref{monotonedifferencediscrete2} that assumes that the endpoints of $\bm{\mathsf{x}}$ and $\widetilde{\bm{\mathsf{x}}}$ are equal, namely, $\bm{u} = \widetilde{\bm{u}}$ and $\bm{v} = \widetilde{\bm{v}}$. This variant further incorporates upper boundaries $g, \widetilde{g}$ to the non-intersecting Gaussian bridges $\bm{\mathsf{x}}, \widetilde{\bm{\mathsf{x}}}$ (in addition to the lower boundaries $f, \widetilde{f}$).
	
	\begin{prop}
		
		\label{2difference} 
		
		Fix an integer $n \ge 1$; two $n$-tuples $\bm{u}, \bm{v} \in \overline{\mathbb{W}}_n$; and four functions $f, \widetilde{f}, g, \widetilde{g} : \llbracket 0, 2 \rrbracket \rightarrow \overline{\mathbb{R}}$ with $f(1) \ge \widetilde{f}(1)$ and $g(1) \le \widetilde{g}(1)$. Sample non-intersecting Gaussian bridges $\bm{\mathsf{x}} (t)$ and $\widetilde{\bm{\mathsf{x}}} (t)$ from the measures $\mathsf{G}_{f;g}^{\bm{u}; \bm{v}}$ and $\mathsf{G}_{\tilde{f}; \tilde{g}}^{\bm{u}; \bm{v}}$, respectively. Then, there exists a coupling between $\bm{\mathsf{x}}$ and $\widetilde{\bm{\mathsf{x}}}$ such that, for each $j \in \llbracket 1, n-1 \rrbracket$,
		\begin{flalign}
			\label{xfguvequal}
			\mathsf{x}_n (1) - f(1) \le\widetilde{\mathsf{x}}_n (1) - \widetilde{f} (1); \quad g(1) - \mathsf{x}_1 (1) \le \widetilde{g}(1) -\widetilde{\mathsf{x}}_1 (1); \quad  \mathsf{x}_j (1) - \mathsf{x}_{j+1} (1) \le\widetilde{\mathsf{x}}_j (1) -\widetilde{\mathsf{x}}_{j+1} (1).
		\end{flalign}
	\end{prop}
	
	We will show \Cref{2difference} through the following lemma, which assumes that either $f(1) = \widetilde{f} (1)$ or $g(1) = \widetilde{g}(1)$, and gives a slightly stronger coupling. 
	
	\begin{lem}
		
		\label{3difference}
		
		Adopt the notation and assumptions of \Cref{2difference}. 
		\begin{enumerate} 
			\item  If $f(1) = \widetilde{f} (1)$, then there exists a coupling between $\bm{\mathsf{x}}$ and $\widetilde{\bm{\mathsf{x}}}$ such that \eqref{xfguvequal} holds, and thus that $\mathsf{x}_j (1) \le\widetilde{\mathsf{x}}_j (1)$ for each $j \in \llbracket 1, n \rrbracket$. 
			
			\item If $g(1) = \widetilde{g}(1)$, then there exists a coupling between $\bm{\mathsf{x}}$ and $\widetilde{\bm{\mathsf{x}}}$ such that \eqref{xfguvequal} holds, and thus that $\mathsf{x}_j (1) \ge\widetilde{\mathsf{x}}_j (1)$ for each $j \in \llbracket 1, n \rrbracket$.
			
		\end{enumerate} 
	\end{lem}

	Given \Cref{3difference}, we can quickly establish \Cref{2difference}. 
	
	\begin{proof}[Proof of \Cref{2difference}]

		Sample non-intersecting $2$-step Gaussian bridges $\widehat{\mathsf{x}} (t)$ from the measure $\mathsf{G}_{\tilde{f}; g}^{\bm{u}; \bm{v}}$ (so that it has lower boundary $\widetilde{f}$ and upper boundary $g$). Applying \Cref{3difference} twice yields couplings between $(\bm{\mathsf{x}}; \widehat{\bm{\mathsf{x}}})$ and $(\widehat{\bm{\mathsf{x}}}; \widetilde{\bm{\mathsf{x}}})$ such that 
		\begin{flalign*}
			\mathsf{x}_n (1) - f(1) \le \widehat{\mathsf{x}}_n (1) - \widetilde{f} (1); \quad g(1) - \mathsf{x}_1 (1) \le g(1) - \widehat{\mathsf{x}}_1 (1); \quad \mathsf{x}_j (1) - \mathsf{x}_{j+1} (1) \le \widehat{\mathsf{x}}_j (1) - \widehat{\mathsf{x}}_{j+1} (1); \\
			\widehat{\mathsf{x}}_n (1) - \widetilde{f}(1) \le\widetilde{\mathsf{x}}_n (1) - \widetilde{f} (1); \quad g(1) - \widehat{\mathsf{x}}_1 (1) \le \widetilde{g}(1) -\widetilde{\mathsf{x}}_1 (1); \quad \widehat{\mathsf{x}}_j (1) - \widehat{\mathsf{x}}_{j+1} (1) \le\widetilde{\mathsf{x}}_j (1) -\widetilde{\mathsf{x}}_{j+1} (1).
		\end{flalign*}
		
		\noindent Combining these couplings (first sampling $\widehat{\bm{\mathsf{x}}}$ conditional on $\bm{\mathsf{x}}$, and then sampling $\widetilde{\bm{\mathsf{x}}}$ conditional on $\bm{\mathsf{x}}$) yields one between $\bm{\mathsf{x}}$ and $\widetilde{\bm{\mathsf{x}}}$ such that \eqref{xfguvequal} holds.	
	\end{proof}

	Now we can establish \Cref{3difference}. 
	
	\begin{proof}[Proof of \Cref{3difference}]
		
		We only address the second case $g(1) = \widetilde{g}(1)$ the lemma, as its proof if $f(1) = \widetilde{f}(1)$ is entirely analogous; throughout, we set $f = f(1)$, $\widetilde{f} = \widetilde{f}(1)$, and $g = g(1) = \widetilde{g}(1)$. 
		
		We induct on $n \ge 1$. To verify the result if $n = 1$, observe since $f \ge \widetilde{f}$ that \Cref{fgk} yields a coupling between $\bm{\mathsf{x}}$ and $\widetilde{\bm{\mathsf{x}}}$ so that $\widetilde{\mathsf{x}}_1 (1) \le \mathsf{x}_1 (1) \le\widetilde{\mathsf{x}}_1 (1) + f - \widetilde{f}$; this confirms \eqref{xfguvequal} and the bound $\mathsf{x}_1 (1) \ge\widetilde{\mathsf{x}}_1 (1)$, establishing the lemma if $n = 1$.
		
		Next suppose $n > 1$. Let $\bm{\mathsf{y}} (t) = \big( \mathsf{y}_1 (t), \mathsf{y}_2 (t), \ldots , \mathsf{y}_n (t) \big)$ and $\widetilde{\bm{\mathsf{y}}} (t) = \big( \widetilde{\mathsf{y}}_1 (t), \widetilde{\mathsf{y}}_2 (t), \ldots , \widetilde{\mathsf{y}}_n (t) \big)$ be two families of non-intersecting $2$-step Gaussian random walks sampled under the measures $\mathsf{G}_{f; g}^{\bm{u}; \bm{v}}$ and $\mathsf{G}_{\tilde{f}; g}^{\bm{u}; \bm{v}}$, respectively (so that they have the same laws as $\bm{\mathsf{x}} (t)$ and $\widetilde{\bm{\mathsf{x}}} (t)$, respectively). By \Cref{fgk}, and the fact that $f \ge \widetilde{f}$, we may couple $\bm{\mathsf{y}}$ and $\widetilde{\bm{\mathsf{y}}}$ so that
		\begin{flalign}
			\label{yjyj}
			\mathsf{y}_j (1) \ge \widetilde{\mathsf{y}}_j (1), \quad \text{and} \quad \mathsf{y}_j (1) \le \widetilde{\mathsf{y}}_j (1) + f - \widetilde{f}, \qquad \text{for each $j \in \llbracket 1, n \rrbracket$}.
		\end{flalign}
		
		\noindent  Define $\widehat{f}, \breve{f}: \llbracket 0, 2 \rrbracket \rightarrow \mathbb{R}$ by setting $\widehat{f} (t) = f(t) = \breve{f}(t)$ for $t \in \{ 0, 2 \}$, and setting $\widehat{f} (1) = \mathsf{y}_n (1)$ and $\breve{f}(1) = \widetilde{\mathsf{y}}_n (1)$ for $t = 1$. Also let $\widehat{\bm{u}} = (u_1, u_2, \ldots , u_{n-1}) \in \mathbb{W}_{n-1}$ and $\widehat{\bm{v}} = (v_1, v_2, \ldots , v_{n-1}) \in \mathbb{W}_{n-1}$.

		Given $\bm{\mathsf{y}}$, we can sample $\bm{\mathsf{x}}$ by first fixing $\mathsf{x}_n (1) = \mathsf{y}_n (1)$, and then sampling the remaining points $\big( \mathsf{x}_1 (1), \mathsf{x}_2 (1), \ldots , \mathsf{x}_{n-1} (1) \big)$ according to the measure $\mathsf{G}_{\widehat{f}; g}^{\widehat{\bm{u}}; \widehat{\bm{v}}}$ (this is equivalent to first sampling the bottom point $\mathsf{x}_n (1)$ of $\bm{\mathsf{x}}$ according to its marginal, and then resampling the others conditional on $\mathsf{x}_n (1)$). Similarly, given $\widetilde{\bm{\mathsf{y}}}$, we can sample $\widetilde{\bm{\mathsf{x}}}$ by setting $\widetilde{x}_n (1) = \widetilde{\mathsf{y}}_n (1)$, and then resampling $\big(\widetilde{\mathsf{x}}_1 (t),\widetilde{\mathsf{x}}_2 (t), \ldots ,\widetilde{\mathsf{x}}_{n-1} (t) \big)$ according to $\mathsf{G}_{\breve{f}; g}^{\widehat{\bm{u}}; \widehat{\bm{v}}}$. See \Cref{f:sequential_sample}. 
		
				\begin{figure}
	\center
\includegraphics[scale = .65]{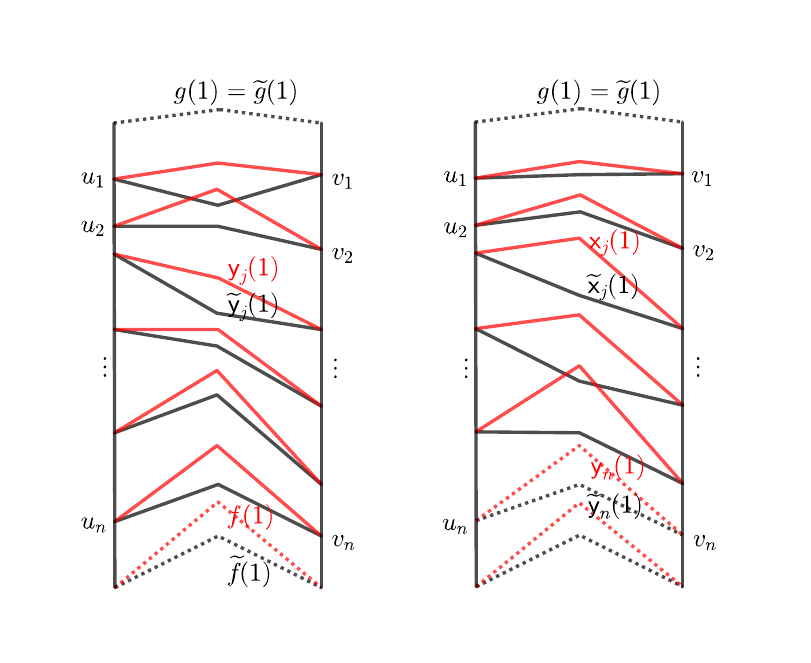}

\caption{In the proof of \Cref{3difference}, we first sample $\bm{\mathsf{y}} (t)$ and $\widetilde{\bm{\mathsf{y}}} (t)$, as shown on the left. Next, we fix $\mathsf{x}_n (1) = \mathsf{y}_n (1)$ and $\widetilde{\mathsf{x}}_n (1) = \widetilde{\mathsf{y}}_n (1)$, and sample the remaining $\big( \mathsf{x}_1 (1), \mathsf{x}_2 (1), \ldots , \mathsf{x}_{n-1} (1) \big)$ and $\big(\widetilde{\mathsf{x}}_1 (t),\widetilde{\mathsf{x}}_2 (t), \ldots ,\widetilde{\mathsf{x}}_{n-1} (t) \big)$, as shown on the right. }
\label{f:sequential_sample}
	\end{figure}
		
		Since \eqref{yjyj} gives $\mathsf{x}_n (1) = \mathsf{y}_n (1) \ge \widetilde{\mathsf{y}}_n (1) = \widetilde{\mathsf{x}}_n (1)$, the inductive hypothesis (and the fact that $g = \widetilde{g}$) yields a coupling between $\big(\mathsf{x}_1 (t), \mathsf{x}_2 (t), \ldots , \mathsf{x}_{n-1} (t) \big)$ and $\big(\widetilde{\mathsf{x}}_1 (t),\widetilde{\mathsf{x}}_2 (t), \ldots ,\widetilde{\mathsf{x}}_{n-1} (t) \big)$ so that
		\begin{flalign*}
			\mathsf{x}_j (1) - \mathsf{x}_{j+1} (1) \le\widetilde{\mathsf{x}}_j (1) -\widetilde{\mathsf{x}}_{j+1} (1); \qquad g - \mathsf{x}_1 (1) \le g -\widetilde{\mathsf{x}}_1 (1); \qquad \mathsf{x}_j (1) \ge\widetilde{\mathsf{x}}_j (1),
		\end{flalign*}
		
		\noindent for each $j \in \llbracket 1, n-1 \rrbracket$. By \eqref{yjyj} and the fact that $\mathsf{x}_n (1) = \mathsf{y}_n (1)$ and $\widetilde{x}_n (1) = \widetilde{\mathsf{y}}_n (1)$, we further have that $\mathsf{x}_n (1) \ge\widetilde{\mathsf{x}}_n (1)$ and $\mathsf{x}_n (1) - f = \mathsf{y}_n (1) - f \le \widetilde{\mathsf{y}}_n (1) - \widetilde{f} = \widetilde{\mathsf{x}}_n (1) - \widetilde{f}$. Thus, this coupling satisfies the required properties, which establishes the lemma.   
	\end{proof}

	\subsection{Proof of \Cref{monotonedifferencediscrete2}}
	
	\label{Proof2T}
	
	In this section we establish \Cref{monotonedifferencediscrete2}. We begin by reducing to the following case of it.
	
	\begin{lem}
		
		\label{4difference}
		
		If $u_n = v_n$, $\widetilde{u}_n = \widetilde{v}_n$, and $f(1) = \widetilde{f}(1)$, then \Cref{monotonedifferencediscrete2} holds.
	\end{lem} 
	
	Assuming \Cref{4difference}, we can quickly show \Cref{monotonedifferencediscrete2} holds in general.
	
	\begin{proof}[Proof of \Cref{monotonedifferencediscrete2}]
		
		We first reduce to the case when $u_n = v_n$ and $\widetilde{u}_n = \widetilde{v}_n$. Observe by using an affine shift to replace $\big( \mathsf{x}_j (t) \big)$ and $\big( f(t) \big)$ by 
		\begin{flalign*}
			\bigg( \mathsf{x}_j (t) - u_n + \displaystyle\frac{t}{2} (u_n - v_n) \bigg), \quad \text{and} \quad \bigg( f(t) - u_n + \displaystyle\frac{t}{2} (u_n - v_n) \bigg), \quad \text{respectively},
		\end{flalign*} 
		
		\noindent and $\big(\widetilde{\mathsf{x}}_j (t) \big)$ and $\big( \widetilde{f} (t) \big)$ with
		\begin{flalign*} 
			\bigg(\widetilde{\mathsf{x}}_j (t) - \widetilde{u}_n + \displaystyle\frac{t}{2} (\widetilde{u}_n - \widetilde{v}_n) \bigg), \quad \text{and} \quad \bigg( f(t) - \widetilde{u}_n + \displaystyle\frac{t}{2} (\widetilde{u}_n - \widetilde{v}_n) \bigg), \quad \text{respectively},
		\end{flalign*}
		
		\noindent we can assume by \Cref{discretelinear} (and the fact that such affine transformations do not affect the differences $\mathsf{x}_j (t) - \mathsf{x}_{j+1} (t)$, $\mathsf{x}_n (t) - f(t)$, $\widetilde{x}_j (t) -\widetilde{\mathsf{x}}_{j+1} (t)$, and $\widetilde{x}_n (t) - f(t)$) that $u_n = v_n = \widetilde{u}_n = \widetilde{v}_n$.
		
		Next, observe that $f(1) \ge \widetilde{f} (1)$, as repeated application of \eqref{y1} and \eqref{y200} yields
		\begin{flalign*}
			\widetilde{f}(2) - 2 \widetilde{f} (1) + \widetilde{f}(0) - u_n - v_n & \ge f(2) - u_n - 2f(1) + f(0) - v_n \\
			& \ge \widetilde{f}(2) - \widetilde{u}_n - 2f(1) + \widetilde{f}(0) - \widetilde{v}_n = \widetilde{f} (2) - 2 f (1) + \widetilde{f}(0) - u_n - v_n.
		\end{flalign*}

		To reduce to when $f(1) = \widetilde{f} (1)$, we follow the proof of \Cref{2difference} given \Cref{3difference}. Sample a family of $n$ non-intersecting $2$-step Gaussian bridges $\widehat{\bm{\mathsf{x}}}(t) = \big( \widehat{\mathsf{x}}_1 (t), \widehat{\mathsf{x}}_2 (t), \ldots , \widehat{\mathsf{x}}_n (t) \big)$ from the measure $\mathsf{G}_f^{\tilde{\bm{u}}; \tilde{\bm{v}}}$. Since \Cref{4difference} indicates that \Cref{monotonedifferencediscrete2} holds when $f(1) = \widetilde{f}(1)$, we obtain a coupling between $(\bm{\mathsf{x}}; \widehat{\bm{\mathsf{x}}})$ such that 
		\begin{flalign*}
			\mathsf{x}_j (t) - \mathsf{x}_{j+1} (t) \le \widehat{\mathsf{x}}_j (t) - \widehat{\mathsf{x}}_{j+1} (t); \qquad \mathsf{x}_n (t) - f(t) \le \widehat{\mathsf{x}}_n (t) - f(t).
		\end{flalign*} 
	
		\noindent Moreover, \Cref{2difference} provides a coupling between and $(\widehat{\bm{\mathsf{x}}}; \widetilde{\bm{\mathsf{x}}})$ such that 
		\begin{flalign*}
			& \widehat{\mathsf{x}}_j (t) - \widehat{\mathsf{x}}_{j+1} (t) \le\widetilde{\mathsf{x}}_j (t) -\widetilde{\mathsf{x}}_{j+1} (t); \qquad \widehat{\mathsf{x}}_n (t) - f(t) \le\widetilde{\mathsf{x}}_n (t) - \widetilde{f} (t).
		\end{flalign*} 
		
		\noindent So, combining these couplings yields one between $\bm{\mathsf{x}}$ and $\widetilde{\bm{\mathsf{x}}}$ satisfying the required properties. 
	\end{proof}

	Now let us establish \Cref{4difference}.
	
	\begin{proof}[Proof of \Cref{4difference}]
		
		We induct on the number 
		\begin{flalign*} 
			\ell = \ell (\bm{\mathsf{x}}; \widetilde{\bm{\mathsf{x}}})= \# \big\{ j \in \llbracket 1, n \rrbracket : u_j \ne \widetilde{u}_j \big\} + \# \big\{ j \in \llbracket 1, n \rrbracket : v_j \ne \widetilde{v}_j \big\} \in \llbracket 0, 2n-2 \rrbracket, 
		\end{flalign*} 
		
		\noindent of ``mismatches'' between the boundary data for $\bm{\mathsf{x}}$ and $\widetilde{\bm{\mathsf{x}}}$. The result is true for $\ell = 0 $ by \Cref{2difference}, so let us assume that $\ell \ge 1$ and prove the lemma assuming it holds for smaller $\ell$. 
		
		Let $k \le n$ be the smallest index such that $\mathsf{x}_j (0) =\widetilde{\mathsf{x}}_j (0)$ and $\mathsf{x}_j (2) =\widetilde{\mathsf{x}}_j (2)$, for each $j \in \llbracket k, n \rrbracket$. We may assume that $k > 1$, for otherwise the lemma follows from \Cref{2difference}. Then, by \eqref{y1}, we either have $\widetilde{\mathsf{x}}_{k-1} (0) > \mathsf{\mathsf{x}}_{k-1} (0)$ or $\widetilde{x}_{k-1} (2) > \mathsf{x}_{k-1} (2)$. The two cases are entirely analogous, so let us assume the latter holds and set $\Delta =\widetilde{\mathsf{x}}_{k-1} (2) - \mathsf{x}_{k-1} (2)$.

				\begin{figure}
	\center
\includegraphics[width=0.7\textwidth,trim=0 1cm 0 1cm, clip]{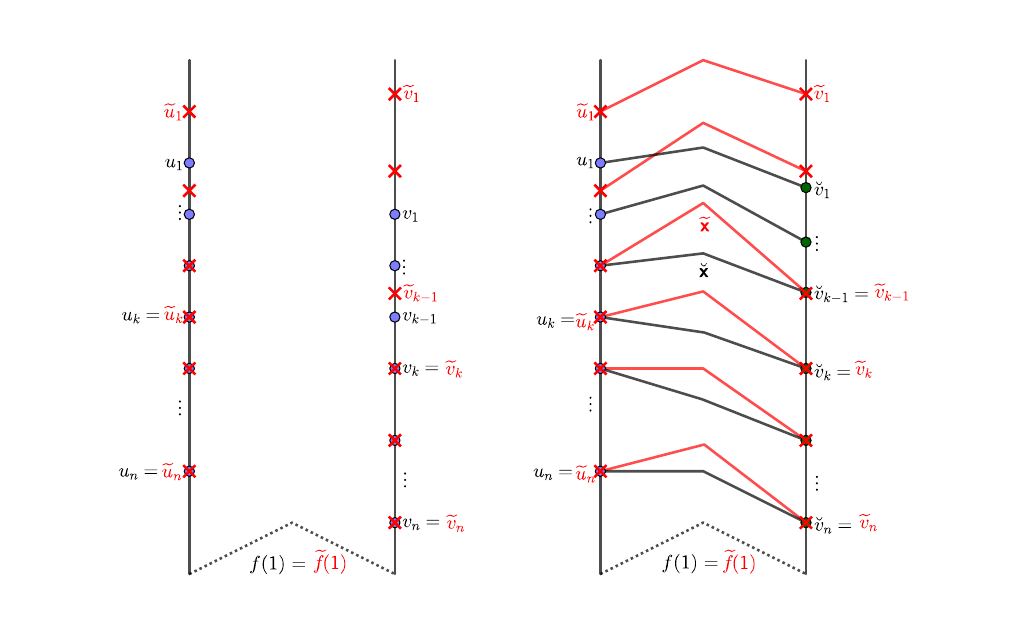}

\caption{Shown to the left are the boundary data $(\bm u, \bm v)$ and $(\widetilde{ \bm u}, \widetilde \bmv)$, satisfying $u_j =\widetilde u_j$ and $v_j=\widetilde v_j$ for $j \in \llbracket k, n \rrbracket$, and $\widetilde v_{k-1}\geq v_{k-1}$. Shown to the right is the new ending data $\breve{\bm{v}} = (\breve{v}_1, \breve{v}_2, \ldots , \breve{v}_n) =(v_1 + \Delta, v_2 + \Delta, \ldots , v_{k-1} + \Delta, v_k, v_{k+1}, \ldots , v_n)$, satisfying $\breve v_{k-1}=\widetilde v_{k-1}$, and the associated Gaussian bridges $\breve{\bm{\mathsf{x}}}$ (coupled with $\widetilde{\bm{\mathsf{x}}}$).}
\label{f:sequential_couple}
	\end{figure}

		Define the $n$-tuple $\breve{\bm{v}} = (\breve{v}_1, \breve{v}_2, \ldots , \breve{v}_n) =(v_1 + \Delta, v_2 + \Delta, \ldots , v_{k-1} + \Delta, v_k, v_{k+1}, \ldots , v_n) \in \mathbb{W}_n$, and sample the family $\breve{\bm{\mathsf{x}}} (t) = \big( \breve{\mathsf{x}}_1 (t), \breve{\mathsf{x}}_2 (t), \ldots , \breve{\mathsf{x}}_n (t) \big)$ of $n$ non-intersecting Gaussian bridges from the measure $\mathsf{G}_f^{\bm{u}; \breve{\bm{v}}}$; see \Cref{f:sequential_couple}.  Observe that $\ell (\breve{\mathsf{x}}; \widetilde{\bm{\mathsf{x}}}) < \ell (\bm{\mathsf{x}}; \widetilde{\bm{\mathsf{x}}}) = \ell$, since $\breve{v}_{j-1} - \breve{v}_j = v_{j-1} - v_j$ if $j \ne k$ and $\breve{v}_{k-1} - \breve{v}_k = v_{k-1} + \Delta - v_k = v_{k-1} + \Delta - \widetilde{v}_k = \widetilde{v}_{k-1} - \widetilde{v}_k$. Hence, the inductive hypothesis yields a coupling between $\breve{\bm{\mathsf{x}}}$ and $\widetilde{\bm{\mathsf{x}}}$ such that 
		\begin{flalign}
			\label{x2x}
			\breve{\mathsf{x}}_n (1) \le\widetilde{\mathsf{x}}_n (1), \quad \text{and} \quad \breve{\mathsf{x}}_j (1) - \breve{\mathsf{x}}_{j+1} (1) \le\widetilde{\mathsf{x}}_j (1) -\widetilde{\mathsf{x}}_{j+1} (1), \qquad \text{for each $j \in \llbracket 1, n-1 \rrbracket$}.
		\end{flalign}
		
		We claim that it is possible to couple $\bm{\mathsf{x}}$ and $\breve{\bm{\mathsf{x}}}$ in such a way that 
		\begin{flalign}
			\label{x1x}
			\mathsf{x}_n (1) \le \breve{\mathsf{x}}_n (1), \quad \text{and} \quad \mathsf{x}_j (1) - \mathsf{x}_{j+1} (1) \le \breve{\mathsf{x}}_j (1) - \breve{\mathsf{x}}_{j+1} (1), \qquad \text{for each $j \in \llbracket 1, n-1 \rrbracket$}.
		\end{flalign}
		
		\noindent Together with \eqref{x2x}, this would imply the existence of a coupling between $\bm{\mathsf{x}}$ and $\widetilde{\bm{\mathsf{x}}}$ satisfying the required properties. 
		
		It therefore remains to establish \eqref{x1x}, which proceeds  similarly to in the proof of \Cref{3difference}. Specifically, let $\bm{\mathsf{y}}(t) = \big( \mathsf{y}_1 (t), \mathsf{y}_2 (t), \ldots , \mathsf{y}_n (t) \big)$ and $\breve{\bm{\mathsf{y}}} = \big( \breve{\mathsf{y}}_1 (t), \breve{\mathsf{y}}_2 (t), \ldots , \breve{\mathsf{y}}_n (t) \big)$ be families of $n$ non-intersecting $2$-step Gaussian random walks, sampled under the measures $\mathsf{G}_f^{\bm{u}; \bm{v}}$ and $\mathsf{G}_f^{\bm{u}; \breve{\bm{v}}}$, respectively. By \Cref{fgk} (and the fact that $v_j \le \breve{v}_j \le v_j + \Delta$ for each $j \in \llbracket 1, n \rrbracket$), there is a coupling between $\bm{\mathsf{y}}$ and $\breve{\bm{\mathsf{y}}}$ such that
		\begin{flalign}
			\label{yj1yj1} 
			\mathsf{y}_j (1) \le \breve{\mathsf{y}}_j (1) \le \mathsf{y}_j (1) + \displaystyle\frac{\Delta}{2}, \qquad \text{for each $j \in \llbracket 1, n \rrbracket$}.
		\end{flalign}
		
		Define the starting points $\bm{u}' = (u_1, u_2, \ldots , u_{k-2}) \in \mathbb{W}_{k-2}$ and $\bm{u}'' = (u_k, u_{k+1}, \ldots , u_n) \in \mathbb{W}_{n-k+1}$, and define the ending points $\bm{v}', \breve{\bm{v}}' \in \mathbb{W}_{k-2}$ and $\bm{v}'', \breve{\bm{v}}'' \in \mathbb{W}_{n-k+1}$ similarly. Given $\bm{\mathsf{y}}$, we can sample $\bm{\mathsf{x}}$ by first fixing $\mathsf{x}_{k-1} (1) = \mathsf{y}_{k-1} (1)$, and then sampling $\bm{\mathsf{x}}' = \big( \mathsf{x}_1 (1), \mathsf{x}_2 (1), \ldots , \mathsf{x}_{k-2} (1) \big)$ and $\bm{\mathsf{x}}'' = \big( \mathsf{x}_k (1), \mathsf{x}_{k+1} (1), \ldots , \mathsf{x}_n (1) \big)$ from $\mathsf{G}_{\mathsf{x}_{k-1}(1)}^{\bm{u}'; \bm{v}'}$ and $\mathsf{G}_{f(1); \mathsf{x}_{k-1} (1)}^{\bm{u}''; \bm{v}''}$, respectively.\footnote{For any functions $h, g : \llbracket 0, 2 \rrbracket \rightarrow \mathbb{R}$, starting points $\bm{r}$, and ending points $\bm{w}$, we are implicitly setting $\mathsf{G}_{h(1); g(1)}^{\bm{r}; \bm{w}} = G_{h; g}^{\bm{r}; \bm{w}}$, as this measure only depends on $h$ and $g$ through $\big(h(1), g(1) \big)$.} Similarly, given $\breve{\bm{\mathsf{y}}}$, we can sample $\breve{\bm{\mathsf{x}}}$ by fixing $\breve{\mathsf{x}}_{k-1} (1) = \breve{\mathsf{y}}_{k-1} (1)$, and then sampling $\breve{\bm{\mathsf{x}}}' = \big( \breve{\mathsf{x}}_1 (1), \breve{\mathsf{x}}_2 (1), \ldots , \breve{\mathsf{x}}_{k-1} (1) \big)$ and $\breve{\bm{\mathsf{x}}}'' = \big( \breve{\mathsf{x}}_k (1), \breve{\mathsf{x}}_{k+1} (1), \ldots , \breve{\mathsf{x}}_n (1) \big)$ from $\mathsf{G}_{\breve{\mathsf{x}}_{k-1} (1)}^{\bm{u}'; \breve{\bm{v}}'}$ and $\mathsf{G}_{f(1); \breve{\mathsf{x}}_{k-1} (1)}^{\bm{u}''; \breve{\bm{v}}''}$, respectively. 
		
		By \eqref{yj1yj1} and the first part of \Cref{3difference}, it is possible to couple $\bm{\mathsf{x}}''$ and $\breve{\bm{\mathsf{x}}}''$ so that 
		\begin{flalign}
			\label{xjxj1}
			\mathsf{x}_n (1) \le \breve{\mathsf{x}}_n (1); \quad \text{and} \quad \mathsf{x}_j (1) - \mathsf{x}_{j+1} (1) \le \breve{\mathsf{x}}_j (1) - \breve{\mathsf{x}}_{j+1} (1), \qquad \text{for each $j \in \llbracket k-1, n-1 \rrbracket$}.
		\end{flalign}
		
		\noindent To couple $\bm{\mathsf{x}}'$ and $\breve{\bm{\mathsf{x}}}'$, observe that the starting data $\bm{u}'$ of these non-intersecting path ensembles coincide, and that their ending data $(\bm{v}'; \breve{\bm{v}}')$ coincide up to a shift, namely, $v_j = \breve{v}_j -\Delta$ for each $j \in \llbracket 1, k-1 \rrbracket$. Moreover, \eqref{yj1yj1} gives the bound $\breve{\mathsf{x}}_{k-1} (1) - \Delta / 2 \le \mathsf{x}_{k-1} (1)$. So, upon subtracting the linear function $t \Delta / 2$ from $\breve{\bm{\mathsf{x}}}'$ and using \Cref{discretelinear}, the ($g=\infty$ case of the) second part of \Cref{3difference} applies to yield a coupling between $\bm{\mathsf{x}}'$ and $\breve{\bm{\mathsf{x}}}'$ so that	
		\begin{flalign}
			\label{xjxj2} 
			\mathsf{x}_j (1) - \mathsf{x}_{j+1} (1) \le \breve{\mathsf{x}}_j (1) - \breve{\mathsf{x}}_{j+1} (1), \qquad \text{for each $j \in \llbracket 1, k- 2 \rrbracket$}.
		\end{flalign}
		
		\noindent By \eqref{xjxj1} and \eqref{xjxj2}, this couples $\bm{\mathsf{x}}$ and $\breve{\bm{\mathsf{x}}}$ in a way satisfying \eqref{x1x}, establishing the lemma.
	\end{proof}

	\section{Likelihood of Medium Position Events}
	\label{LLocation}

	In this section we establish that the $\textbf{TOP}$ events (recall \Cref{eventsregular1}) are likely upon restricting to the $\textbf{PAR}$ ones (see \Cref{x1small}), and also prove results indicating that the $\textbf{MED}$ events (recall \Cref{eventsregular1}) are likely upon restricting to the $\textbf{TOP}$ ones (see \Cref{ljt} below). The latter shows that the $\textbf{MED}$ part of the $\textbf{SCL}$ ones from \Cref{eventscl} is likely; the proof that the $\textbf{GAP}$ and $\textbf{REG}$ parts are also likely will appear in \Cref{ProbabilityScale} below. Throughout this section, we let $\bm{\mathsf{x}} = (\mathsf{x}_1, \mathsf{x}_2, \ldots ) \in \mathbb{Z}_{\ge 1} \times \mathcal{C}(\mathbb{R})$ denote a $\mathbb{Z}_{\ge 1} \times \mathbb{R}$ indexed line ensemble satisfying the Brownian Gibbs property. We also recall the set $\mathfrak{T}_k (\alpha; A)$ and the events $\textbf{PAR}$, $\textbf{MED}$, and $\textbf{TOP}$ from \Cref{parabolat} and \Cref{eventsregular1}.

	\subsection{Proof of \Cref{x1small0}} 
	
	\label{Proofx1}
	
\begin{figure}
	\center
\includegraphics[width=1\textwidth, trim={2cm 1cm 2cm 0cm},clip]{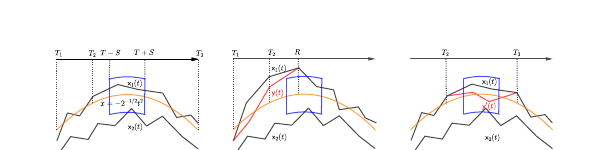}

\caption{Shown on the left is a depiction of \Cref{x1small}, indicating if $\mathsf x_1(t)$ and the parabola $-2^{-1/2}t^2$ are close at the three points $\{ T_1, T_2, T_3 \}$, then they are close on the entire interval $[T-S, T+S]$ (shown by the blue box). 
	Shown in the middle is a depiction that, if $\mathsf x_1(R)$ is too high, then so is $\mathsf x_1(T_2)\geq \mathsf y_1(T_2)$. Shown on the right is a depiction that, if $\mathsf x_1(t)$ and the parabola $-2^{-1/2}t^2$ are close at the two points $\{ T_2, T_3 \}$, then $\mathsf x_1(t)\geq \mathsf y'(t)$ cannot be too low for $t\in [T-S, T+S]$.}
\label{f:top_event}
	\end{figure}

	In this section we establish \Cref{x1small0}, which is a quick consequence of the next lemma, stating the following. Suppose that the top curve $\mathsf{x}_1 (t)$ of $\bm{\mathsf{x}}$ is close to the parabola $-2^{1/2} t^2$ at three points $T_1 = T - 5S$, $T_2 = T - 2S$, and $T_3 = T+3S$, for some parameter $S$ that is much smaller than another parameter $T$. Then $\mathsf{x}_1 (t)$ remains close, of distance much smaller than $T^2$, to this parabola on an interval between them; see the left side of \Cref{f:top_event}.

	\begin{lem} 
		
		\label{x1small} 
		
		Fix an integer $k \ge 1$ and two real numbers $\varepsilon \in ( 0, 1 / 4)$ and $A \ge 1$. Let $T \in [-10Ak^{1/3}, 10Ak^{1/3}]$ be a real number, and denote 
		\begin{flalign}\label{e:defC0}
			 \vartheta = 1000  A^2 \varepsilon; \qquad S = \varepsilon k^{1/3}; \qquad T_1 = T - 5S; \qquad T_2 = T - 2S; \qquad T_3 = T+3S.
		\end{flalign}
		
		\noindent We have 
		\begin{flalign*}
			\mathbb{P} \Big[ \textbf{\emph{PAR}}_{\varepsilon} \big( \{ T_1, T_2, T_3 \}; \varepsilon k^{2/3} \big) \cap \textbf{\emph{TOP}} \big( [T-S, T+S]; \vartheta k^{2/3} \big)^{\complement} \Big] \le 4 e^{-\varepsilon^3 k}. 
		\end{flalign*}
		
	\end{lem} 
	
	\begin{proof}[Proof of \Cref{x1small0}]

		 By \Cref{x1small}, at all points $T \in [-10Ak^{1/3}, 10Ak^{1/3}] \cap (S \cdot \mathbb{Z})$ (which would force $T_1, T_2, T_3 \in [-15Ak^{1/3}, 15Ak^{1/3}] \cap (S \cdot \mathbb{Z}) = \mathfrak{T}_k (\varepsilon; 15A)$), with a union bound,  we have  
		\begin{flalign*}
			\mathbb{P} \Big[ \textbf{PAR}_{\varepsilon} \big( \mathfrak{T}_k (\varepsilon; 15A); \varepsilon k^{2/3} \big) \cap \textbf{TOP} \big( [-10Ak^{1/3}, 10Ak^{1/3}]; \vartheta k^{2/3} \big)^{\complement} \Big] \le 120\varepsilon^{-1} A e^{-\varepsilon^3 k},
		\end{flalign*}
		
		\noindent where we also used the fact that the number of such $T$ in the union bound is at most  $30 \varepsilon^{-1} A$. This yields the proposition, by taking $k$ sufficiently large relative to $\varepsilon$, $\delta$, and $A$. 
	\end{proof}

	The proof of \Cref{x1small} will make use of the following events, which will also be used throughout this section.
	
	\begin{definition} 
		
		\label{hlevent} 

		For any integer $k \ge 1$ and real numbers $t, B \in \mathbb{R}$, define the \emph{low position event} $\textbf{LOW}_k (t; B) = \textbf{LOW}_k^{\bm{\mathsf{x}}} (t; B)$\index{L@$\textbf{LOW}$; low position event} and \emph{high position event} $\textbf{HIGH}_k (t; B) = \textbf{HIGH}_k^{\bm{\mathsf{x}}} (t; B)$\index{H@$\textbf{HIGH}$; high position event} by setting
		\begin{flalign*}
			& \textbf{LOW}_k (t; B) = \big\{ \mathsf{x}_k (t) \le -2^{-1/2} t^2-B  \big\}; \qquad  \textbf{HIGH}_k (t; B) = \big\{ \mathsf{x}_k (t) \ge -2^{-1/2} t^2-B  \big\}.
		\end{flalign*} 
	
		\noindent Similarly, for any subset $\mathcal{T} \subseteq \mathbb{R}$, define the events $\textbf{LOW}_k (\mathcal{T}; B) = \textbf{LOW}_k^{\bm{\mathsf{x}}} (\mathcal{T}; B)$ and $\textbf{HIGH}_k (\mathcal{T}; B) = \textbf{HIGH}_k^{\bm{\mathsf{x}}} (\mathcal{T}; B)$ by setting
		\begin{flalign*} 
			\textbf{LOW}_k (\mathcal{T}; B) = \bigcap_{t \in \mathcal{T}} \textbf{LOW}_k (t; B); \qquad  \textbf{HIGH}_k (\mathcal{T}; B) = \bigcap_{t \in \mathcal{T}} \textbf{HIGH}_k (t; B),
		\end{flalign*} 
		
	\end{definition} 
	
	\noindent Using the above notions, we can establish \Cref{x1small}. 
	
	\begin{proof}[Proof of \Cref{x1small}]
		
		By a union bound (and the facts that $\textbf{PAR}_{\varepsilon} \big( \{ T_1, T_2, T_3 \}; \varepsilon k^{2/3} \big) \subseteq \textbf{PAR}_{\varepsilon} \big( \{ T_1, T_2 \}; \varepsilon k^{2/3} \big)$ and $\textbf{PAR}_{\varepsilon} \big( \{ T_1, T_2, T_3 \}; \varepsilon k^{2/3} \big) \subseteq \textbf{PAR}_{\varepsilon} \big( \{ T_2, T_3 \}; \varepsilon k^{2/3} \big)$) it suffices to show 
		\begin{flalign}
			\label{eventphpl} 
			\begin{aligned} 
				& \mathbb{P} \Bigg[ \textbf{PAR}_{\varepsilon} \big( \{ T_1, T_2 \}; \varepsilon k^{2/3} \big) \cap \bigcup_{t \in [T-S,T+S]} \textbf{HIGH}_1 (t; -\vartheta k^{2/3}) \Bigg] \le 2 e^{-S^3}; \\
				& \mathbb{P} \Bigg[ \textbf{PAR}_{\varepsilon} \big( \{ T_2, T_3 \}; \varepsilon k^{2/3} \big) \cap \bigcup_{t \in [T-S, T+S]} \textbf{LOW}_1 (t; \vartheta k^{2/3}) \Bigg] \le 2 e^{-S^3}.
			\end{aligned}
		\end{flalign}
		
		We begin with verifying the first bound in \eqref{eventphpl}, to which end we condition on $\mathcal{F}_{\ext} \big( \{ 1 \} \times (T_1, T+S) \big)$ (recall \Cref{property}) and restrict to the event 
		\begin{flalign*} 
			\mathscr{E}_1 = \textbf{PAR}_{\varepsilon} (T_1; \varepsilon k^{2/3}) \cap \bigcup_{t \in [T-S, T+S]} \textbf{HIGH}_1 (t; -\vartheta k^{2/3}).
		\end{flalign*} 
	
		\noindent We will then show that $\mathsf{x}_1 (T_2)$ is likely larger than allowed by the event $\textbf{PAR}_{\varepsilon} (T_2; \varepsilon k^{2/3})$. Due to our restriction to $\mathscr{E}_1$, there exists some real number $R \in [T-S, T+S]$ such that $\mathsf{x}_1 (R) \ge \vartheta k^{2/3} -2^{-1/2} R^2$. Letting $R \in [T-S, T+S]$ be the largest such real number, we find that $(T_1, R)$ is a $\{ 1 \}$-stopping domain in the sense of \Cref{property2}. Thus, \Cref{propertyproperty2} implies that the law of $\mathsf{x}_1$ on $[T_1, R]$, conditional on $u = \mathsf{x}_1 (T_1)$, $v = \mathsf{x}_1 (R)$, and $f = \mathsf{x}_2 |_{[T_1, R]}$, is given by a Brownian bridge conditioned to start at $u$, end at $v$, and remain above $f$. See the middle of \Cref{f:top_event}. 
		
		Letting $\mathsf{y} : [T_1, R] \rightarrow \mathbb{R}$ denote a Brownian bridge conditioned to start at $u$ and end at $v$, \Cref{monotoneheight} yields a coupling between $\mathsf{x}_1$ and $\mathsf{y}$ such that $\mathsf{x}_1 (T_2) \ge \mathsf{y} (T_2)$. It follows that 
		\begin{flalign}
			\label{hxyp} 
			\begin{aligned}
			\mathbb{P} & \Bigg[ \textbf{PAR}_{\varepsilon} \big( \{ T_1, T_2 \}; \varepsilon k^{2/3} \big) \cap \bigcup_{t \in [T-S,T+S]} \textbf{HIGH}_1 (t; -\vartheta k^{2/3}) \Bigg] \\
			& \le \mathbb{P} \Big[ \big\{ \mathsf{x}_1 (T_2) \le \varepsilon k^{2/3} - (2^{-1/2} - \varepsilon) T_2^2 \big\} \cap \mathscr{E}_1  \Big] \le \mathbb{P} \Big[ \big\{ \mathsf{y} (T_2) \le \varepsilon k^{2/3} - (2^{-1/2} - \varepsilon) T_2^2 \big\} \cap \mathscr{E}_1 \Big].
			\end{aligned}
		\end{flalign} 
	
		\noindent Applying \Cref{maximumx} to, and using the affine invariance (\Cref{linear}) of, $\mathsf{y}$ yields 
		\begin{flalign}
			\label{y2t}
			\mathbb{P} \bigg[ \mathsf{y} (T_2) \le \displaystyle\frac{R - T_2}{R - T_1} \cdot u + \displaystyle\frac{T_2 - T_1}{R - T_1} \cdot v - 2a (R-T_1)^{1/2} \bigg] \le 2e^{-a^2},
		\end{flalign}
		
		\noindent for any real number $a > 0$. Let us show at $a = S^{3/2}$ that the quantity on the right side inside the probability in \eqref{y2t} is at least $\varepsilon k^{2/3} - (2^{-1/2} + \varepsilon) T_2^2$. This will follow from direct computation, using the fact that $v$ is ``lifted'' by $\vartheta k^{2/3}$, namely, that we have $v = \mathsf{x}_1 (R) \ge \vartheta k^{2/3} - 2^{-1/2} R^2$.
		
		To implement this, observe that 
		\begin{flalign}
			\label{rt1} 
			R - T_1 \le 6S; \qquad \frac{T_2 - T_1}{R - T_1} \ge \frac{1}{2},
		\end{flalign} 
		 
		 \noindent where the first holds since $T_1 = T-5S$ and $R \in [T-S, T+S]$, and the second holds since $T_2 - T_1 = 3S$ (and $R - T_1 \le 6S$). Therefore, 
		\begin{flalign}
			\label{rt2} 
			\displaystyle\frac{R - T_2}{R - T_1} \cdot u + \displaystyle\frac{T_2 - T_1}{R - T_1} \cdot v \ge \min \{ u, v - \vartheta k^{2/3} \} + \displaystyle\frac{\vartheta k^{2/3}}{2}.
		\end{flalign}
	
		\noindent To estimate the right side of this inequality, observe that 
		\begin{flalign}
			\label{rt3} 
			u \ge -\varepsilon k^{2/3} - (2^{-1/2} + \varepsilon) T_1^2 \ge -307 \varepsilon A^2 k^{2/3} - 2^{-1/2} T_2^2,
		\end{flalign} 
	
		\noindent where the first holds since $u = \mathsf{x}_1 (T_1)$ and we restricted to the event $\textbf{PAR}_{\varepsilon} (T_1; \varepsilon k^{2/3}) \subseteq \mathscr{E}_1$; the second holds since $T_1^2 = (T_2-3S)^2 \le T^2 + 81 \varepsilon A k^{2/3}$ (as $S = \varepsilon k^{1/3}$, $|T_2| \le |T| + 2S \le 12 Ak^{1/3}$, and $A \ge 1$) and $\varepsilon T_1^2 \le 225 A^2 \varepsilon k^{2/3}$ (as $|T_1| \le |T| + 5S \le 15Ak^{1/3}$). We also have 
		\begin{flalign}
			\label{rt4} 
			v - \vartheta k^{2/3} \ge  - 2^{-1/2} R^2 \ge -81 \varepsilon A k^{2/3} -2^{-1/2} T_2^2,
		\end{flalign} 
	
		\noindent where the first holds since $v = \mathsf{x}_1 (R) \ge \vartheta k^{2/3} - 2^{-1/2} R^2$ and the second holds since $|R|^2 \le \big(|T_2| + 3S \big)^2 \le T^2 + 81 A \varepsilon k^{2/3}$ (as $R \in [T_2 + S, T_2 + 3S]$, again with the facts that $S = \varepsilon k^{1/3}$, $|T_2| \le 12Ak^{1/3}$, and $A \ge 1$). Combining \eqref{rt1}, \eqref{rt2}, \eqref{rt3}, and \eqref{rt4} gives 
		\begin{flalign}
			\label{rt5}
			\displaystyle\frac{R-T_2}{R-T_1} \cdot u + \displaystyle\frac{T_2-T_1}{R-T_1} \cdot v - 2a(R-T_1)^{1/2} \ge \displaystyle\frac{\vartheta k^{2/3}}{2} - 307 \varepsilon A^2 k^{2/3} - 6a S^{1/2} - 2^{-1/2} T_2^2.
		\end{flalign}
		
		\noindent We will take $a = S^{3/2}$, so observe (by \eqref{e:defC0}, and the facts that $|T_2| \le 12Ak^{1/3}$ and $S = \varepsilon k^{1/3}$) that 
		\begin{flalign*} 
			\displaystyle\frac{\vartheta k^{2/3}}{2} = 500 A^2 \varepsilon k^{2/3} \ge 307 \varepsilon A k^{2/3} + 6S^2 + 2^{-1/2} \varepsilon T_2^2 + \varepsilon k^{2/3}.
		\end{flalign*} 
		
		\noindent This, with \eqref{rt5} and the $a = S^{3/2}$ case of \eqref{y2t}, gives $\mathbb{P} \big[ \{ \mathsf{y} (T_2) \le \varepsilon k^{2/3} - (2^{-1/2} - \varepsilon) T_2^2 \} \cap \mathscr{E}_1 \big] \le 2e^{-S^3}$, which by \eqref{hxyp} yields the first bound in \eqref{eventphpl}.
		
		Next we verify the second bound in \eqref{eventphpl}. Restrict to the event $\mathscr{E}_2 = \textbf{PAR}_{\varepsilon} \big( \{ T_2, T_3 \}; \varepsilon k^{2/3} \big)$. We then will show that $\mathsf{x}_1$ is likely larger than allowed by the event $\textbf{LOW}_1 \big( [T-S, T+S]; \varepsilon k^{2/3} \big)$; see the right side of \Cref{f:top_event}. To this end, conditional on $u' = \mathsf{x}_1 (T_2)$, $v' = \mathsf{x}_1 (T_3)$, and $f' = \mathsf{x}_2 |_{[T_2, T_3]}$, the law of $\mathsf{x}_1 |_{[T_2, T_3]}$ is given by a Brownian bridge conditioned to start at $u'$, end at $v'$, and remain above $f'$. Letting $\mathsf{y}' : [T_2, T_3] \rightarrow \mathbb{R}$ denote a Brownian bridge conditioned to start at $u'$ and end at $v'$, \Cref{monotoneheight} again yields a coupling between $\mathsf{x}_1$ and $\mathsf{y}'$ such that $\mathsf{x}_1 (t) \ge \mathsf{y}' (t)$, for each $t \in [T_2, T_3]$. It follows that 
		\begin{flalign}
			\label{pxy2}
			\begin{aligned}
			\mathbb{P} \Bigg[ \textbf{PAR}_{\varepsilon} \big( \{ T_2, T_3 \}; \varepsilon k^{2/3} & \big) \cap \bigcup_{t \in [T-S,T+S]} \textbf{LOW}_1 (t; \vartheta k^{2/3}) \Bigg] \\
			& \le \mathbb{P} \Bigg[ \bigcup_{t \in [T-S, T+S]} \big\{ \mathsf{x}_1 (t)  \le -2^{-1/2} t^2 - \vartheta k^{2/3} \big\} \cap \mathscr{E}_2 \Bigg] \\
			& \le \mathbb{P} \Bigg[ \bigcup_{t \in [T-S,T+S]} \big\{ \mathsf{y}' (t) \le -2^{-1/2} t^2 - \vartheta k^{2/3} \big\} \cap \mathscr{E}_2 \Bigg].
			\end{aligned} 
		\end{flalign}
		
		\noindent We once again use \Cref{maximumx} (and \Cref{linear}) to deduce for any real number $a > 0$ that 
		\begin{flalign}
			\label{yt2t3t2}
			\mathbb{P} \Bigg[ \bigcup_{t \in [T-S, T+S]} \bigg\{ \mathsf{y}' (t) \le \displaystyle\frac{T_3 - t}{T_3 - T_2} \cdot u + \displaystyle\frac{t - T_2}{T_3 - T_2} \cdot v - 2a (T_3-T_2)^{1/2} \bigg\} \Bigg] \le 2 e^{-a^2}.
		\end{flalign}
	
		\noindent Let us show at $a = S^{3/2}$ that the quantity on the right side of the probability in \eqref{yt2t3t2} is at least $-2^{-1/2} t^2 - \varepsilon k^{2/3}$; this will again follow from direct computation, using the fact that $\vartheta$ is large relative to $\varepsilon$.
		
		To implement this, observe for any $t \in [T-S, T+S]$ that   
		\begin{flalign}
			\label{rt6}
			\displaystyle\frac{T_3-t}{T_3 - T_2} \cdot u + \displaystyle\frac{t - T_2}{T_3 - T_2} \cdot v \ge \min \{ u, v \} & \ge -(2^{-1/2} + \varepsilon) \cdot \displaystyle\max \{ T_2^2, T_3^2 \} - \varepsilon k^{2/3},
		\end{flalign} 
	
		\noindent where the last bound holds since we restricted to $\mathscr{E}_2 = \textbf{PAR}_{\varepsilon} \big( \{ T_2, T_3 \}; \varepsilon k^{2/3} \big)$. We also have   
		\begin{flalign}
			\label{rt7} 
			-\max \{ T_2^2, T_3^2 \} \le - \big( |t| + 4S \big)^2 \le -t^2 - 104 A \varepsilon k^{2/3},
		\end{flalign}
	
		\noindent where the first statement holds since $t \in [T-S, T+S]$ and $(T_2, T_3) = (T-2S, T+3S)$, and the second holds since $|t| \le |T| + S \le 11Ak^{1/3}$ and $S = \varepsilon k^{1/3}$. Therefore,
		\begin{flalign*}
			\displaystyle\frac{T_3-t}{T_3-T_2} \cdot u + \displaystyle\frac{t-T_2}{T_3 - T_2} \cdot v - 2S^{3/2} (T_3 - T_2)^{1/2} & \ge -(2^{-1/2} + \varepsilon) t^2 - 111 A \varepsilon k^{2/3} \\
			& \ge -2^{-1/2} t^2 - 232 A^2 \varepsilon k^{2/3} \ge - 2^{-1/2} t^2 - \vartheta k^{2/3} ,
		\end{flalign*}
		
		\noindent where the first inequality follows from inserting \eqref{rt7} into \eqref{rt6} (and using the facts that $T_3 - T_2 = 5S$ and that $2^{-1/2} + \varepsilon \le 1$); the second follows from the fact that $|t| \le |T| + S \le 11Ak^{1/3}$; and the third follows from the definition \eqref{e:defC0} of $\vartheta$. Applying this in the $a = S^{3/2}$ case of \eqref{yt2t3t2} gives
		\begin{flalign*}
			\mathbb{P} \Bigg[ \bigcup_{t \in [T-S,T+S]} \big\{ \mathsf{y}' (t) \le -2^{-1/2} t^2 - \vartheta k^{2/3} \big\} \cap \mathscr{E}_2 \Bigg] \le 2 e^{-S^3},
		\end{flalign*}
		
		\noindent which by \eqref{pxy2} yields the second statement of \eqref{eventphpl} and thus the lemma. 
	\end{proof}

	\subsection{Likelihood of $\textbf{MED}$ Restricted to $\textbf{TOP}$}
	
	\label{ProofEventL} 
	
	In this section we state and establish \Cref{ljt}, which indicates the following.  If the top curve $\mathsf{x}_1 (t)$ of $\bm{\mathsf{x}}$ is close to $2^{-1/2} t^2$ on an interval with length of order $k^{2/3}$, then the distance between its $j$-th curve $\mathsf{x}_j(t)$ and this parabola is of order $j^{2/3}$, for each integer $j$ of order $k$.

	\begin{prop}
		
		\label{ljt}
		
		For any real numbers $A \ge 20$ and $B \ge 1$, there exists a constant $C = C(A,B) > 1$ such that 
		\begin{flalign*} 
			\mathbb{P} & \Bigg[ \bigcup_{j = \lfloor k/B \rfloor}^{\lfloor Bk \rfloor} \textbf{\emph{MED}}_j \bigg( [-Ak^{1/3}, Ak^{1/3}]; \frac{j^{2/3}}{100}; 450 j^{2/3} \bigg)^{\complement} \\
			& \qquad \qquad \cap \textbf{\emph{TOP}} \bigg( [-10ABk^{2/3}, 10ABk^{2/3}]; \displaystyle\frac{k^{2/3}}{200B} \bigg) \Bigg]  \le C e^{-(\log k)^2}. 
		\end{flalign*} 
		
	\end{prop}

	 The proof of \Cref{ljt} uses the following three lemmas (where we recall the $\textbf{LOW}$ and $\textbf{HIGH}$ events in them from \Cref{hlevent}). The first indicates that a line ensemble likely cannot remain low at every point of a long interval, if its top curve decays parabolically; it is shown in \Cref{Proofl0} below. The second indicates that, if $\mathsf{x}_k$ is very low at a given point, then there likely exists a long interval on which it is low at every point; it is shown in \Cref{Proofl0} below. The third indicates that, if $\mathsf{x}_k$ is too high at a given point and the top curve of the ensemble decays parabolically, then there again likely exists a long interval on which $\mathsf{x}_k$ is too low at every point; it is shown in \Cref{Prooflh0} below. These together show that $\mathsf{x}_k$ can neither be too low nor too high anywhere on the interval.

	\begin{lem} 
		
		\label{probabilityeventl0} 
		
		There exists a constant $C>1$ such that the following holds. Fix an integer $k \ge 1$ and real numbers $T_1, T_2 \in \mathbb{R}$ with $T_2 - T_1 = 8k^{1/3}$. Setting $T = (T_1 + T_2) / 2$, we have                                                                                                                 
		\begin{flalign*}                                                                                                                                                                                                                                                                                                                                                                          	                                                               
			\mathbb{P} \Big[ \textbf{\emph{LOW}}_k \big( [T_1, T_2]; 20 k^{2/3} \big) \cap \textbf{\emph{TOP}} \big( \{ T_1, T, T_2 \} &; k^{2/3} \big) \Big]  \le C e^{-(\log k)^3}.            
		\end{flalign*}                                                                                                                                                                                                 
		
	\end{lem}

	\begin{lem}                                                                                                                                                                                                    
		
		\label{probabilityeventll0}                                                                                                                                                                                    
		
		There exists a constant $C > 1$ such that the following holds. Let $k \ge 1$ be an integer and $T \in \mathbb{R}$ be a real number; set $S = 8k^{1/3}$. We have
		\begin{flalign*}
			\mathbb{P} \Big[ \textbf{\emph{HIGH}}_k & \big([T-S, T+S]; 450k^{2/3} \big)^{\complement} \cap \textbf{\emph{LOW}}_k \big( [T-3S, T-2S ]; 20k^{2/3} \big)^{\complement} \\
			& \quad \cap \textbf{\emph{LOW}}_k \big( [ T+2S, T+3S ]; 20 k^{2/3} \big)^{\complement} \Big] \le C e^{-(\log k)^3}.
		\end{flalign*}
	\end{lem}

	\begin{lem}
		
		\label{probabilityeventclh0}
		
		There exists a constant $C > 1$ such that the following holds. Let $k \ge 1$ be an integer and $T \in \mathbb{R}$ be a real number; set $S = 8k^{1/3}$. We have 
		\begin{flalign*}
			\mathbb{P} \bigg[  \textbf{\emph{LOW}}_k \Big( [T, T+S]; \displaystyle\frac{k^{2/3}}{100} & \Big)^{\complement} \cap  \textbf{\emph{LOW}}_k \big( [T+3S, T+4S]; 20k^{2/3} \big)^{\complement} \\
			&  \cap \textbf{\emph{TOP}} \Big( [T, T+4S]; \displaystyle\frac{k^{2/3}}{200} \Big) \bigg] \le C e^{-(\log k)^3}.
		\end{flalign*}
	\end{lem}
	
	Given the above three lemmas, we can establish \Cref{ljt}. 
	
	\begin{proof}[Proof of \Cref{ljt}]
		
		Throughout this proof, we will restrict to the event given by $\textbf{TOP} \big( [-10ABk^{1/3}, 10ABk^{1/3}]; k^{2/3}/200B \big)$, which we will abbreviate $\textbf{TOP}$. Fix an index $j \in \llbracket \lfloor k/B \rfloor, \lfloor Bk \rfloor \rrbracket$; denote $S = 8j^{1/3}$; and fix an interval $\mathcal{I} = [T_1, T_2] \subseteq [-ABk^{1/3}, ABk^{1/3}]$ of length $T_2 - T_1 = S$. It suffices to show for some constant $C>1$ that, off of an event of probability at most $C e^{-(\log j)^3}$, we have 
		\begin{flalign}
			\label{akl0} 
			\displaystyle\inf_{t \in \mathcal{I}} \big( \mathsf{x}_j (t) + 2^{-1/2} t^2 \big) \ge -450j^{2/3}; \qquad \displaystyle\inf_{|t| \le Ak^{1/3}} \big( \mathsf{x}_j (t) + 2^{-1/2} t^2 \big) \le -\displaystyle\frac{j^{2/3}}{100}.
		\end{flalign}
		
		\noindent Indeed, then the proposition would follow from taking a union bound over all (at most $Bk$) such indices $j$ and a set (of size at most $AB$) disjoint intervals $\mathcal{I}$ covering $[-Ak^{1/3}, Ak^{1/3}]$.
		
		To that end, since we have restricted to the event $\textbf{TOP} \subseteq \textbf{TOP} \big( [T_1-4S, T_2+4S]; j^{2/3} \big)$, we deduce the following by applying \Cref{probabilityeventl0} three times (with the $(T_1, T_2)$ there equal to $(T_1-3S, T_1-2S)$, $(T_1+2S, T_1+3S)$, and $(T_1+3S,T_1+4S)$ here). Off of an event of probability at most $C e^{-(\log j)^3}$, there exist real numbers $R_1 \in [T_1-3S, T_1-2S]$, $R_2 = [T_1+2S, T_1+3S]$, and $R_3 \in [T_1+3S, T_1+4S]$ such that 
		\begin{flalign}
			\label{r1r2r3} 
			\mathsf{x}_j (R) + 2^{-1/2} R^2 \ge - 20j^{2/3}, \qquad \text{for all $R \in \{ R_1, R_2, R_3 \}$}.
		\end{flalign} 
	
		In what follows, we restrict to the event on which \eqref{r1r2r3} holds. By \Cref{probabilityeventl0}, the existence of $R_1$ and $R_2$ satisfying \eqref{r1r2r3} yields that the first statement in \eqref{akl0} holds, off of an event of probability at most $Ce^{-(\log k)^3}$. The existence of $R_3$ satisfying \eqref{r1r2r3}, with our restriction to $\textbf{TOP} \subseteq \textbf{TOP} \big( [T_1, T_1+4S]; j^{2/3}/200 \big)$, yields by \Cref{probabilityeventclh0} that the second statement in \eqref{akl0} holds, off of a further event of probability at most $Ce^{-(\log k)^3}$. This establishes the proposition.
	\end{proof}

	\subsection{Avoiding Low Intervals}
	
	\label{Proofl0} 
	
	In this section we establish \Cref{probabilityeventl0}, according to the outline in \Cref{Location0}. The fact underlying its proof, which holds by \Cref{estimatexj}, is the following. In a family of $k-1$ non-intersecting Brownian bridges on some time interval $[-t, t]$, conditioned to start and end at position $-2^{-1/2} t^2$ (that is, a Brownian watermelon on $[-t,t]$, shifted down by $2^{-1/2} t^2$; recall \Cref{PathsUV0}), the top curve at time $0$ is located around $(2kt)^{1/2} - 2^{-1/2} t^2$; for $t = 4k^{1/3}$, this is much less than $0$. So by height monotonicity \Cref{monotoneheight}, in order to not contradict the $\textbf{PAR}$ event, the $k$-th curve $\mathsf{x}_k$ our Brownian line ensemble $\bm{\mathsf{x}}$ must ``push up'' the $(k-1)$-th one $\mathsf{x}_{k-1}$. Hence, $\mathsf{x}_k$ cannot be lower on all of $[-t, t]$ than the bottom curve of the above watermelon, whose minimum (again by \Cref{estimatexj}) is around $-(2kt)^{1/2} - 2^{-1/2} t^2 > -20k^{2/3}$.

		\begin{figure}
	\center
\includegraphics[scale = 1.2, trim=6cm 11cm 6cm 10cm, clip]{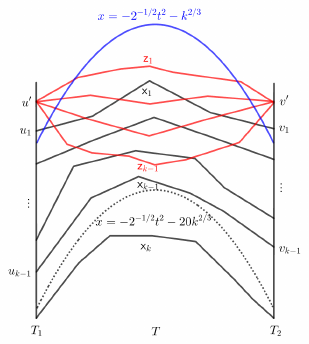}

\caption{Shown above is the coupling between $\bm{\mathsf{x}}$ and $\bm{\mathsf{z}}$ used to show \Cref{probabilityeventl0}. }
\label{f:Low_Top_event}
	\end{figure}

	\begin{proof}[Proof of \Cref{probabilityeventl0}] 
		
		By replacing $\bm{\mathsf{x}} (t)$ with $\bm{\mathsf{x}} (t - T) + 2^{-1/2} (T^2 - 2tT)$, which does not affect our assumptions (by \Cref{linear}, as it is an affine shift). Thus, we may assume throughout this proof that $T = 0$, and hence that $(T_1, T_2) = (-4k^{2/3}, 4k^{2/3})$. 
		
		For the remainder of this proof, we also condition on $\mathcal{F}_{\ext} \big( \llbracket 1, k-1 \rrbracket \times (T_1, T_2) \big)$ and restrict to the event $\mathscr{E} = \textbf{LOW}_k \big( [T_1, T_2]; 20 k^{2/3} \big) \cap \textbf{TOP} \big( \{ T_1, T_2 \}; k^{2/3} \big)$; we will then show that $\mathsf{x}_1 (0)$ is likely lower than allowed by the event $\textbf{TOP} (0; k^{2/3})$. In what follows, we define the $(k-1)$-tuples $\bm{u} = \bm{\mathsf{x}}_{\llbracket 1,k-1 \rrbracket} (T_1) \in \mathbb{W}_{k-1}$ and $\bm{v} = \bm{\mathsf{x}}_{\llbracket 1,k-1 \rrbracket} (T_2) \in \mathbb{W}_{k-1}$ , and the function $f = \mathsf{x}_k |_{[T_1, T_2]}$. Then, the law of $\big( \mathsf{x}_j (s) \big)$ over $(j, s) \in \llbracket 1, k-1 \rrbracket \times [T_1, T_2]$ is given by $\mathsf{Q}_f^{\bm{u}; \bm{v}}$. 
		
		Define $u' = k^{2/3} - 2^{-1/2} T_1^2 = k^{2/3} - 2^{-1/2} T_2^2$, and denote the $(k-1)$-tuple $\bm{u}' = (u', u', \ldots , u') \in \overline{\mathbb{W}}_{k-1}$ (where $u'$ has multiplicity $k-1$). Further define the function $f' : [T_1, T_2] \rightarrow \mathbb{R}$ by setting $f' (s) = -2^{-1/2} s^2 -20 k^{2/3}$ for $s \in [T_1, T_2]$. Then, sample two families of non-intersecting Brownian bridges $\bm{\mathsf{y}} = (\mathsf{y}_1, \mathsf{y}_2, \ldots , \mathsf{y}_{k-1}) \in \llbracket 1, k-1 \rrbracket \times \mathcal{C} \big( [T_1, T_2] \big)$ and $\bm{\mathsf{z}} = (\mathsf{z}_1, \mathsf{z}_2, \ldots , \mathsf{z}_{k-1} \big) \in \llbracket 1, k-1 \rrbracket \times \mathcal{C} \big( [T_1, T_2] \big)$ from the measures $\mathsf{Q}_{f'}^{\bm{u}'; \bm{u}'}$ and $\mathsf{Q}^{\bm{u}'; \bm{u}'}$, respectively; see \Cref{f:Low_Top_event} for the latter.
				
		Observe that $f' (s) = -2^{-1/2} s^2 -20 k^{2/3} \ge f (s)$, as we restricted to $\textbf{LOW}_k (20 k^{2/3}; T_1, T_2) \subseteq \mathscr{E}$. We also have $u' \ge \mathsf{x}_1 (T_1) \ge \mathsf{x}_j (T_1)$ and $v' \ge \mathsf{x}_1 (T_2) \ge \mathsf{x}_j (T_2)$ for each $j \in \llbracket 1, k-1 \rrbracket$, since we have restricted to $\textbf{TOP} \big( \{ T_1, T_2 \}; k^{2/3} \big) \subseteq \mathscr{E}$. Thus, \Cref{monotoneheight} yields a coupling between $\bm{\mathsf{x}}$ and $\bm{\mathsf{y}}$ with 
		\begin{flalign}
			\label{kxjyj} 
			\mathsf{x}_j (s) \le \mathsf{y}_j (s), \qquad \text{for each $(j, s) \in \llbracket 1, k-1 \rrbracket \times [T_1, T_2]$}.
		\end{flalign}
		
		Next by \Cref{estimatexj} (with the $(n; a, b; u, v; D)$ there equal to $(k-1; -4k^{1/3}, 4k^{1/3}; u', u'; 2)$ here, using the fact that $(t-T_1)(T_2 - t) \le 16k^{2/3}$ for $t \in [T_1, T_2] = [-4k^{1/3}, 4k^{1/3}]$) and the fact that $-2 \le \gamma_{\semci; k-1} (k-1) \le \gamma_{\semci; k-1} (1) \le 2$, there is a constant $C_1 > 1$ such that
		\begin{flalign}
			\label{zkt1t2} 
			\begin{aligned}
				\mathbb{P} \Bigg[ \displaystyle\bigcup_{t \in [T_1, T_2]} \big\{ \mathsf{z}_1 (t) \ge (5 - 2^{7/2}) k^{2/3} \big\} \Bigg] \le C_1 e^{-(\log k)^3},
			\end{aligned}
		\end{flalign} 
		
		\noindent and 
		\begin{flalign} 
			\label{zkt1t22}
			\begin{aligned}
				\mathbb{P} \Bigg[ \displaystyle\bigcap_{t \in [T_1, T_2]} \big\{ \mathsf{z}_{k-1} (t) \ge - (2^{7/2} + 3) k^{2/3} \big\} \Bigg] \ge 1 - C_1 e^{-(\log k)^3}.
			\end{aligned}
		\end{flalign}

		\noindent Now observe for any $t \in [T_1, T_2]$ that $-(2^{7/2} + 3) k^{2/3} \ge -2^{-1/2} t^2 - 20k^{2/3} = f'(t)$. Therefore, it follows from \eqref{zkt1t22} that $\bm{\mathsf{z}}$ remains above $f'$ with probability at least $1 - C_1 e^{-(\log k)^3}$. Since the law of $\bm{\mathsf{y}}$ is given by that of $\bm{\mathsf{z}}$, conditioned on $\mathsf{z}_j (s) \ge f' (s)$ for each $s \in [T_1, T_2]$, it follows that $\mathbb{P} \big[ \mathsf{y}_1 (0) \ge -k^{2/3} \big] \le (1 - C_1 e^{-(\log k)^3})^{-1} \cdot \mathbb{P} \big[ \mathsf{z}_1 (0) \ge -k^{2/3} \big] \le 2 \cdot \mathbb{P} \big[ \mathsf{z}_1 (0) \ge -k^{2/3} \big]$. Together with \eqref{kxjyj}, this gives
		\begin{flalign}
			\label{xk1zk1}
			\mathbb{P} \big[ \textbf{TOP} (0; k^{2/3}) \big]  \le 2 \cdot \mathbb{P} \big[ \mathsf{z}_1 (0) \ge - k^{2/3} \big].
		\end{flalign}
	
		\noindent Since $(5 - 2^{7/2}) k^{2/3} \leq -k^{2/3}$,  \eqref{zkt1t2} yields $\mathbb{P} \big[ \mathsf{z}_1 (0) \geq - k^{2/3} \big] \le C_1 e^{-(\log k)^3}$, which upon insertion into \eqref{xk1zk1} yields the lemma.
	\end{proof}

	\subsection{Low Interval From a Low Point}
	
	\label{Proofll0} 
	
	In this section we establish \Cref{probabilityeventll0}, following the outline in \Cref{Location0}.

	 	\begin{figure}
	\center
\includegraphics[width=0.5\textwidth, trim=0 0.5cm 0 0.5cm, clip]{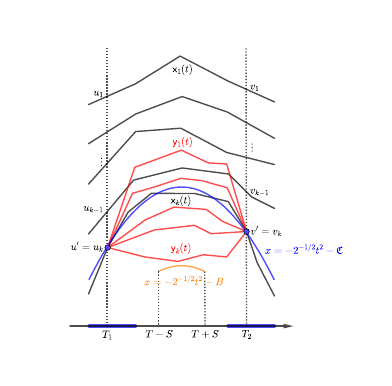}
\caption{\Cref{probabilityeventll0} states that, if $\mathsf x_k(t)$ fails to be consistently below the blue parabola $-2^{-1/2}t^2-20k^{2/3}$, both on the left blue interval  $[T-3S, T-2S]$ and the right one $[T+2S, T+3S]$, then it is likely above the orange parabola $-2^{-1/2}t^2-450k^{2/3}$ on the entire interval $[T-S, T+S]$. The proof proceeds by coupling $\bm {\mathsf x}$ with the red curves $\bm{\mathsf y}$, to likely satisfy ${\mathsf x}_k(t)\geq {\mathsf y}_k(t)\geq -2^{-1/2}t^2-450k^{2/3}$  on $[T-S, T+S]$.}
\label{f:Low_Low_event}
	\end{figure}

	\begin{proof}[Proof of \Cref{probabilityeventll0}]

		As in the proof of \Cref{probabilityeventl0}, we may assume $T = 0$. Throughout, we restrict to the event $\mathscr{E} = \textbf{LOW}_k ( [-3S, -2S]; 20k^{2/3} )^{\complement} \cap \textbf{LOW}_k ( [2S, 3S]; 20k^{2/3} )^{\complement}$. On this event, there exist times $T_1 \in [-3S, -2S ]$ and $T_2 \in [2S, 3S ]$ such that $\mathsf{x}_k (T_1) \ge -2^{-1/2} T_1^2-20k^{2/3}$ and $\mathsf{x}_k (T_2) \ge 2^{-1/2} T_2^2 -20k^{2/3}$. Assume $T_1$ is the smallest such time in $[-3S, -2S]$ and that $T_2$ is the largest such time in $[2S, 3S]$; then, $(T_1, T_2)$ is a $\llbracket 1, k \rrbracket$-stopping domain in the sense of \Cref{property2}. Hence \Cref{propertyproperty2} implies that the law of $\big( \mathsf{x}_j (s) \big)$ for $(j, s) \in \llbracket 1, k \rrbracket \times [T_1, T_2]$, conditional on the $k$-tuples $\bm{u} = \bm{\mathsf{x}}_{\llbracket 1, k \rrbracket} (T_1) \in \mathbb{W}_k$ and $\bm{v} = \mathsf{x}_{\llbracket 1, k \rrbracket} (T_2) \in \mathbb{W}_k$, and on the function $f = \mathsf{x}_{k+1} |_{[T_1, T_2]}$, is given by $\mathsf{Q}_f^{\bm{u}; \bm{v}}$. We will then show that $\mathsf{x}_k$ is likely higher than allowed by the complement of the event $\textbf{HIGH}_k \big( [-S,S]; 450k^{2/3} \big)$.
		
		To this end, set $u' = -2^{-1/2} T_1^2 - 20k^{2/3}$ and $v' = -2^{-1/2} T_2^2 - 20k^{2/3}$, and define the $k$-tuples $\bm{u}' = (u', u', \ldots , u') \in \overline{\mathbb{W}}_k$ and $\bm{v}' = (v', v', \ldots , v') \in \overline{\mathbb{W}}_k$ (where $u'$ and $v'$ both appear with multiplicity $k$). Then sample non-intersecting Brownian bridges $\bm{\mathsf{y}} = (\mathsf{y}_1, \mathsf{y}_2, \ldots , \mathsf{y}_k) \in \llbracket 1, k \rrbracket \times \mathcal{C} \big( [T_1, T_2] \big)$ from the measure $\mathsf{Q}^{\bm{u}'; \bm{v}'}$; see \Cref{f:Low_Low_event}. Since $\mathsf{x}_j (T_1) \ge \mathsf{x}_k (T_1) \ge -2^{-1/2} T_1^2 -20k^{2/3} = u' = \mathsf{y}_j (T_1)$ (and similarly $\mathsf{x}_j (T_2) \ge \mathsf{y}_j (T_2)$) for any $j \in \llbracket 1, k \rrbracket$, by \Cref{monotoneheight} we may couple $\bm{\mathsf{x}}$ with $\bm{\mathsf{y}}$ such that $\mathsf{x}_j (t) \ge \mathsf{y}_j (t)$, for each $(j, t) \in \llbracket 1, k \rrbracket \times [T_1, T_2]$. Hence, 
		\begin{flalign}
			\label{llclc}
			\begin{aligned}
				\mathbb{P} \Big[ \textbf{HIGH}_k \big( [-S, S]; 450k^{2/3}  \big)^{\complement} \cap \mathscr{E} \Big] & \le \mathbb{P} \Bigg[ \mathscr{E} \cap \bigcup_{t \in [-S, S]} \big\{ \mathsf{x}_k (t) \le -2^{-1/2} t^2 -450k^{2/3} \big\} \Bigg] \\ 
				& \le \mathbb{P} \Bigg[ \bigcup_{t \in [-S, S]} \big\{ \mathsf{y}_k (t) \le -2^{-1/2} t^2 -450k^{2/3} \big\} \Bigg].
			\end{aligned}
		\end{flalign}
		
		By the first part of \Cref{estimatexj} (with the $(n; a, b; u, v; D)$ there equal to $(k; T_1, T_2; u', v'; 2)$ here) and the fact that $\gamma_{\semci; k} (k) \ge -2$, there exists a constant $C_1 > 1$ such that 
		\begin{flalign}
			\label{ttstsyk}
			\begin{aligned} 
				\mathbb{P} & \bigg[ \bigcup_{t \in [-S, S]} \Big\{ \mathsf{y}_k (t) \le - 2^{-1/2} (tT_1 + tT_2 - T_1 T_2) - 20k^{2/3} - \big( 2k (T_2 - T_1) \big)^{1/2} \Big\} \bigg] \le C_1 e^{-(\log k)^3}.
			\end{aligned} 
		\end{flalign}
	
		\noindent Now observe for $t \in [-S, S]$ that 
		\begin{flalign*} 
			-2^{-1/2} (& tT_1 + tT_2 - T_1 T_2) - 20k^{2/3} - \big( 2k (T_2 - T_1) \big)^{1/2} \\ 
			&  \ge -2^{-1/2} t^2 - 2^{-1/2} (T_2-t)(t-T_1) - 30 k^{2/3} \ge -2^{-1/2} t^2 - 450k^{2/3},
		\end{flalign*} 
		
		\noindent where the first statement holds since $T_2 - T_1 \le 6S = 48k^{1/3}$, and the second holds since $2^{-1/2} (T_2-t)(t-T_1) \le 2^{-1/2} \cdot 9S^2 \le 410k^{2/3}$ (as $-3S \le T_1 \le t \le T_2 \le 3S$ and $S = 8k^{1/3}$). Inserting this into \eqref{ttstsyk} gives
		\begin{flalign*}
			\mathbb{P} \Bigg[ \bigcup_{t \in [-S,S]} \big\{ \mathsf{y}_k (t) \le -2^{-1/2} t^2 - 450k^{2/3} \big\} \Bigg] \le C_1 e^{-(\log k)^3},
		\end{flalign*}
		
		\noindent which, together with \eqref{llclc}, implies the lemma. 
	\end{proof}

	\subsection{Low Interval From $\textbf{TOP}$ and a High Point} 
	
	\label{Prooflh0} 
	
	In this section we establish \Cref{probabilityeventclh0}, following the outline in \Cref{Location0}.

	 	\begin{figure}
	\center
\includegraphics[scale=1, trim=0 0cm 0 0cm, clip]{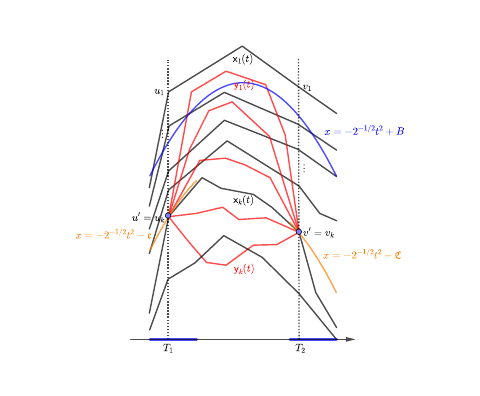}

\caption{Shown above is a depiction of \Cref{probabilityeventclh0}, indicating that if $\mathsf x_k$ is not (entirely) below $-2^{-1/2}t^2-k^{2/3} / 100$ (the left orange parabola) on the interval $[T-2S, T-S]$ (the left blue interval), and is also not below $-2^{-1/2}t^2-20k^{2/3}$ (the right orange parabola) on $[T+S, T+2S]$ (the right blue interval), then $\mathsf{x}_1$ is likely above $-2^{-1/2}t^{2}+k^{2/3}/200$ (the blue parabola) on $[T-2S, T+2S]$. The proof proceeds by coupling $\bm{\mathsf{x}}$ with the red curves $\bm{\mathsf{y}}$ to likely satisfy ${\mathsf x}_k(t)\geq {\mathsf y}_k(t) \ge  - 2^{-1/2} t^2+k^{2/3}/200$.}

\label{f:Low_High_event}

	\end{figure}

	\begin{proof}[Proof of \Cref{probabilityeventclh0}]
		
		As in the proof of \Cref{probabilityeventl0}, we may assume that $T = 0$. It will be convenient to set
		\begin{flalign}
			\label{cbeta} 
			\mathfrak{c} = \displaystyle\frac{k^{2/3}}{100}; \qquad B = \displaystyle\frac{k^{2/3}}{200}; \qquad \beta = \displaystyle\frac{1}{40000}.
		\end{flalign}
	
		\noindent Throughout, we restrict to the event $\mathscr{E} = \textbf{LOW}_k \big( [0, S]; \mathfrak{c} \big)^{\complement} \cap \textbf{LOW}_k \big( [3S, 4S]; 20k^{2/3} \big)^{\complement}$; we will then show that $\mathsf{x}_1$ is at some point likely larger than allowed by $\textbf{TOP} \big( [0, 4S]; B\big)$. On $\mathscr{E}$, there exist times $T_1 \in [0, S]$ and $T_2 \in [3S, 4S]$ such that $\mathsf{x}_k (T_1) \ge -2^{-1/2} T_1^2 - \mathfrak{c}$ and $\mathsf{x}_k (T_2) \ge -2^{-1/2} T_2^2 - 20k^{2/3}$, respectively. Assume $T_1$ is the smallest such time in $[0, S]$ and that $T_2$ is the largest such time in $[3S, 4S]$; then, $(T_1, T_2)$ is a $\llbracket 1, k \rrbracket$-stopping time in the sense of \Cref{property2}. Then \Cref{propertyproperty2} implies that the law of $\big(\mathsf{x}_j (s) \big)$ for $(j, t) \in \llbracket 1, k \rrbracket \times [T_1, T_2]$, conditional on the $k$-tuples $\bm{u} = \bm{\mathsf{x}}_{\llbracket 1, k \rrbracket} (T_1) \in \mathbb{W}_k$ and $\bm{v} = \bm{\mathsf{x}}_{\llbracket 1, k \rrbracket} (T_2) \in \mathbb{W}_k$, and the function $f = \mathsf{x}_{k+1} |_{[T_1, T_2]}$, is given by the non-intersecting Brownian bridge measure $\mathsf{Q}_f^{\bm{u}; \bm{v}}$.
		
		Set $u' = \mathsf{x}_k (T_1)$ and $v' = \mathsf{x}_k (T_2)$. Further define the $k$-tuples $\bm{u}' = (u', u', \ldots , u') \in \overline{\mathbb{W}}_k$ and $\bm{v}' = (v', v', \ldots , v') \in \overline{\mathbb{W}}_k$ (where $u'$ and $v'$ both appear with multiplicity $k$), and sample non-intersecting Brownian bridges $\bm{\mathsf{y}} = (\mathsf{y}_1, \mathsf{y}_2, \ldots , \mathsf{y}_k) \in \llbracket 1, k \rrbracket \times \mathcal{C} \big( [T_1, T_2] \big)$ according to the measure $\mathsf{Q}^{\bm{u}'; \bm{v}'}$; see \Cref{f:Low_High_event}. Since $\bm{u}' \le \bm{u}$ and $\bm{v}' \le \bm{v}$, \Cref{monotoneheight} indicates that we may couple $\bm{\mathsf{x}}$ and $\bm{\mathsf{y}}$ such that 
		\begin{flalign} 
			\label{xjtyjt2} 
			\mathsf{x}_j (t) \ge \mathsf{y}_j (t), \qquad \text{for each $(j, t) \in \llbracket 1, k \rrbracket \times [T_1, T_2]$}.
		\end{flalign} 
		
		By the first statement in \Cref{estimatexj} with $j=1$ (and $D = 2$), and the fact that $\gamma_{\semci; n} (1) \ge 2^{1/2}$ for sufficiently large $n$ (by \Cref{gammaderivative}), there exists a constant $C_1 > 1$ such that 
		\begin{flalign*}
			\mathbb{P} \Bigg[ \bigcap_{t \in [T_1, T_2]} \bigg\{ \mathsf{y}_1 (t) \ge \displaystyle\frac{T_2 - t}{T_2 - T_1} \cdot u' + \displaystyle\frac{t-T_1}{T_2-T_1} \cdot v' +  \Big(  2k \cdot \displaystyle\frac{(T_2-t)(t - T_1)}{T_2-T_1} & \Big)^{1/2} - (T_2 - T_1)^{1/2} \bigg\} \Bigg] \\
			& \qquad \ge 1 - C_1 e^{-(\log k)^3}.
		\end{flalign*}
		
		\noindent Since $T_2 - T_1 \le 4S = 32 k^{1/3}$, taking $t = R$ where $R = (1 - \beta) T_1 + \beta T_2$ (recalling from \eqref{cbeta} that $\beta = 1/3500$), we deduce that 
		\begin{flalign}
			\label{yjtalpha}
			\mathbb{P} \big[ \mathsf{y}_1 (R) \ge (1-\beta)u' + \beta v' + 8 \beta^{1/2} k^{2/3}  - 2^{5/2} k^{1/6} \big] \ge 1 - C_1 e^{-(\log k)^3}.
		\end{flalign}
		
		\noindent We have 
		\begin{flalign*} 
			(1-\beta & ) u' + \beta v' + 8 \beta^{1/2} k^{2/3} - 2^{5/2} k^{1/6}  \\
			& \ge 7 \beta^{1/2} k^{2/3} - 2^{-1/2} (1-\beta) T_1^2 - 2^{-1/2} \beta T_2^2 - 420 \beta k^{2/3} \\
			& = 7 \beta^{1/2} k^{2/3}  - 420 \beta k^{2/3}  - 2^{-1/2} \beta (1-\beta) (T_1 - T_2)^2 - 2^{-1/2} R^2 \\
			& \ge 7 \beta^{1/2} k^{2/3} - 1200 \beta k^{2/3} - 2^{-1/2} R^2 \ge B -2^{-1/2} R^2,
		\end{flalign*}
		
		\noindent where the first statement holds for sufficiently large $k$, as $u' \ge -2^{-1/2} T_1^2 - \mathfrak{c} \ge -2^{-1/2} T_1^2 - 400 \beta k^{2/3}$ (recall \eqref{cbeta}) and $v' \ge - 2^{-1/2} T_2^2 - 20k^{2/3}$; the second holds since $R = (1 - \beta) T_1 + \beta T_2$; the third holds since $|T_1 - T_2| \le 4S \le 32 k^{1/3}$; and the fourth holds since $\beta = 1/40000$ and $B = k^{2/3} / 200$. Together with \eqref{xjtyjt2} and \eqref{yjtalpha}, this yields 
		\begin{flalign*}
			\mathbb{P} \Big[ \textbf{TOP} \big( [0, 4S]; B \big) \cap \mathscr{E} \Big] \le \mathbb{P} \Big[ \big\{ \mathsf{x}_1 (R) \le B -2^{-1/2} R^2 \big\} \cap \mathscr{E} \Big] & \le \mathbb{P} \Big[ \big\{ \mathsf{y}_1 (R) \le B -2^{-1/2} R^2 \big\} \cap \mathscr{E} \Big] \\
			& \le C_1 e^{-(\log k)^3},
		\end{flalign*}
		
		\noindent which establishes the lemma.
	\end{proof}

	\section{Likelihood of On-Scale and Improved Medium Events} 
	
	\label{ProbabilityScale} 
	
	In this section we establish \Cref{sclprobability}, showing that the on-scale event $\textbf{SCL}$ (from \Cref{eventscl}) is likely upon restricting to the $\textbf{TOP}$ event. We also define an ``improved variant''  (see \Cref{eventimproved} below) of the $\textbf{MED}$ part of that event, which considerably extends the range of the index $k$ appearing there, and show it is likely (see \Cref{probabilityimproved} below). Throughout this section, we let $\bm{\mathsf{x}}$ denote a $\mathbb{Z}_{\ge 1} \times \mathbb{R}$ indexed line ensemble satisfying the Brownian Gibbs property and recall the notation of \Cref{Estimatex}.
	
	\subsection{Proof of \Cref{sclprobability}} 
	
	\label{ProofScale} 
	
	In this section we establish \Cref{sclprobability}. That the $\textbf{MED}$ part of that event is likely was shown as \Cref{ljt}, so we must show that the $\textbf{REG}$ and $\textbf{GAP}$ parts of that event are also likely. This is done through the first and second propositions below, which are established in \Cref{ProofEventr0} and \Cref{ProofEventGap}, respectively.
	
	\begin{prop} 
		
		\label{probabilityeventr0}

		For any real numbers $A, B \ge 1$, there exists a constant $C = C(A, B)> 1$ such that, for any integer $k \ge 1$, we have
		\begin{flalign*}
			\mathbb{P} \bigg[ \textbf{\emph{REG}}_k \big( [-Ak^{1/3}, Ak^{1/3} & ]; 3(2A+B); k \big)^{\complement}  \cap \textbf{\emph{TOP}}  \big( [-3Ak^{1/3}, 3Ak^{1/3}]; k^{2/3} \big) \\
			& \qquad \cap \textbf{\emph{MED}}_{k+1} \big( [-3Ak^{1/3}, 3Ak^{1/3}]; -Bk^{2/3}; Bk^{2/3} \big) \bigg] \ \le C e^{-(\log k)^2}.
		\end{flalign*} 
	\end{prop}

	\begin{prop}\label{p:gapupperbound}
		
		For any real number $B \ge 1$, there exist constants $c = c(B) > 0$, $A = A(B) > B \ge 1$, and $C = C(B) > 1$ such that, for any integer $k \ge 1$, we have 
		\begin{flalign*}
			\mathbb{P}  \Bigg[ \textbf{\emph{TOP}} \big( \{-Ak^{1/3}, Ak^{1/3} \}; Bk^{2/3} & \big) \cap \textbf{\emph{MED}}_k \big( \{ -Ak^{1/3}, Ak^{1/3} \}; -Bk^{2/3}; Bk^{2/3} \big) \\
			& \cap  \textbf{\emph{GAP}}_{\lfloor k/2 \rfloor} \bigg( \Big[ -\displaystyle\frac{Ak^{1/3}}{2}, \displaystyle\frac{Ak^{1/3}}{2} \Big]; C \bigg)^{\complement} \Bigg] \le c^{-1} e^{- c (\log k)^2}. 
		\end{flalign*} 
		
	\end{prop}

	Observe the slight difference in how the $\textbf{MED}$ and $\textbf{TOP}$ events are imposed in \Cref{probabilityeventr0} and \Cref{p:gapupperbound}. The former imposes them on an interval, which will be convenient its proof in \Cref{ProofEventr0} below. The latter imposes them at just two points, which will be useful to establish \Cref{probabilityimproved} below.

	\begin{proof}[Proof of \Cref{sclprobability}]
		
		Let $A_0 = A_0 (B)$ be the constant $A(B)$ from \Cref{p:gapupperbound}, and let $R$ denote the constant $C = C(B+900) > 0$ from \Cref{p:gapupperbound}. Set 
		\begin{flalign*} 
			A_1 = B (A_0 + A); \qquad C_1 = 200 B; \qquad C_2 = 30 A_1.
		\end{flalign*}  
	
		\noindent  Throughout this proof, we  restrict to the event $\textbf{TOP} \big( [-C_2 n^{1/3}, C_2 n^{1/3}]; C_1^{-1} n^{2/3} \big)$, which we abbreviate by $\textbf{TOP}$. By \Cref{eventscl}, it suffices to show that there exists a constant $c = c(A,B) > 0$ such that the following holds, for any fixed $k \in \llbracket \lfloor n/B \rfloor, \lfloor Bn \rfloor \rrbracket$. Off of an event of probability at most $c^{-1} e^{-c(\log n)^2}$, the events $\textbf{MED}_k = \textbf{MED}_k \big( [-3A_1 n^{1/3}, 3A_1 n^{1/3}]; k^{2/3} / 100; 450 k^{2/3} \big)$; $\textbf{REG}_k = \textbf{REG}_k \big( [-An^{1/3}, An^{1/3}]; 1000AB; Bn \big)$; and $\textbf{GAP}_n = \textbf{GAP}_n \big( [-An^{1/3}, An^{1/3}]; R \big)$ hold.
		
		Due to our restriction to $\textbf{TOP}$, and our choices of $(A_1, C_1, C_2)$, \Cref{ljt} (with the $(k, A)$ there equal to $(n, 3A_1)$ here) implies that $\textbf{MED}_k$, $\textbf{MED}_{k+1}$, and $\textbf{MED}_{2k}$ all hold off of an event with probability $c^{-1} e^{-c(\log n)^2}$. In particular, this verifies the first of the three events above. So, we further restrict to $\textbf{MED}_k \cap \textbf{MED}_{k+1} \cap \textbf{MED}_{2k}$ below. 
		
		By \Cref{probabilityeventr0} (with the $(A,B)$ there equal to $(AB, 900)$ here), and our restriction to $\textbf{MED}_{k+1} \cap \textbf{TOP}$, we have that $\textbf{REG}_k \big( [-AB k^{1/3}, AB k^{1/3}]; 3(2AB+900); k \big)$ holds off of an event of probability at most $c^{-1} e^{-c(\log n)^2}$. Since $\textbf{REG}_k \big( [-AB k^{1/3}, AB k^{1/3}]; 3(AB+900); k \big) \subseteq \textbf{REG}_k$ (as $An^{1/3} \le AB k^{1/3}$, $3(2AB+900) \le 1000AB$, and $k \le Bn$) this verifes the second event above. 
		
		By \Cref{p:gapupperbound} (with the $(k, A, B)$ there equal to $(2k, 2A, 900)$ here), with our restriction to $\textbf{MED}_k \cap \textbf{TOP}$, we have that $\textbf{GAP}_n$ holds off of an event of probability at most $c^{-1} e^{-c(\log n)^2}$. This verifies the third event above, confirming the theorem. 	
	\end{proof}

	\subsection{Likelihood of $\textbf{REG}$}
	
	\label{ProofEventr0}
	
	In this section we establish \Cref{probabilityeventr0}, which will be a quick consequence of the following variant of it addressing a single (pair of) time(s). 
	
	\begin{lem} 
		
		\label{probabilityeventr}

		There exists a constant $C > 1$ such that the following holds. For any real numbers $A, B \ge 1$ and $s, t \in \mathbb{R}$ with $-Ak^{1/3} \le t \le t + s \le Ak^{1/3}$, we have  
		\begin{flalign*}
			\mathbb{P} \bigg[ \Big\{ \big| \mathsf{x}_k (t+s) - \mathsf{x}_k (t) \big| & \ge 4 (ks)^{1/2} + 3 (A+B) k^{1/3} s + k^{-50} \Big\} \\
			& \qquad \quad \cap \textbf{\emph{MED}}_k \big( [t-k^{1/3}, t+s+k^{1/3}]; -Bk^{2/3}; B k^{2/3} \big) \bigg]  \le C e^{-(\log k)^3}.
		\end{flalign*} 
	\end{lem}

	\begin{proof}[Proof of \Cref{probabilityeventr0}] 
		
		Condition on $\mathcal{F}_{\exp} \big( \llbracket 1, k \rrbracket \times [-3Ak^{1/3}, 3Ak^{1/3}] \big)$ and restrict to the event $\mathscr{E} = \textbf{TOP} \big( \{ -3Ak^{1/3}, 3Ak^{1/3} \}; k^{2/3} \big)\cap \textbf{MED}_{k+1} \big( [ -3Ak^{1/3}, 3Ak^{1/3} ]; -Bk^{2/3}; Bk^{2/3} \big)$; observe that $\mathscr{E} \subseteq \textbf{MED}_k \big( \{ -3Ak^{1/3}, 3Ak^{1/3} \}; -Bk^{2/3}; Bk^{2/3} \big)$, since $\mathsf{x}_{k+1} \le \mathsf{x}_k \le \mathsf{x}_1$ and $B \ge 1$. Then, \Cref{estimatexj3} (where the $(n; A, B; a, b; \mathsf{T}; \kappa)$ there is $(k; Ak^{2/3}, k^4; -3Ak^{1/3}, 3Ak^{1/3}; 6Ak^{1/3}; 1 / 2)$ here) gives constants $c > 0$ and $C_1 > 1$ such that
		\begin{flalign}
			\label{tts2ak}
			\mathbb{P} \Bigg[ \bigcup_{\substack{|t| \le 2Ak^{1/3} \\ |t+s| \le 2Ak^{1/3}}} \bigg\{ \Big| \mathsf{x}_k (t + s) - \mathsf{x}_k (t) - s \cdot \displaystyle\frac{ \mathsf{x}_k (3Ak^{1/3}) - \mathsf{x}_k (-3Ak^{1/3})}{6Ak^{1/3}} \Big|  >  k^5 |s|^{1/3} \bigg\} \Bigg] \le C_1 e^{-ck^4}.
		\end{flalign}
		
		\noindent Since we have restricted to the event $\mathscr{E} \subseteq \textbf{MED}_k \big( \{ -3 Ak^{1/3}, 3Ak^{1/3} \}; -Bk^{2/3}; Bk^{2/3} \big)$, we have  
		\begin{flalign*} 
			\big| \mathsf{x}_k (3Ak^{1/3}) - \mathsf{x}_k (-3Ak^{1/3}) \big| \le 2B k^{2/3} \le 2k^2.
		\end{flalign*} 
		
		\noindent Inserting this into \eqref{tts2ak} (and using the fact that, for sufficiently large $k$, we have $2k^2 s \le k^5 s^{1/3}$, if $s \le 4Ak^{1/3}$) gives
		\begin{flalign}
			\label{tts2ak2}
			\mathbb{P} \Bigg[ \bigcup_{\substack{|t| \le 2Ak^{1/3} \\ |t+s| \le 2Ak^{1/3}}} \Big\{ \big| \mathsf{x}_k (t + s) - \mathsf{x}_k (t)  \big| > 2 k^5 |s|^{1/3} \Big\} \Bigg] \le C_1 e^{-ck^4}.
		\end{flalign}
		
		Next, define the set $\mathcal{T} = [-2Ak^{1/3}, 2Ak^{1/3}] \cap (k^{-200} \cdot\mathbb{Z})$. Applying \Cref{probabilityeventr} (with the $A$ there equal to $2A$ here) and a union bound over all $s, t \in \mathcal{T}$, we deduce the existence of a constant $C_2 > 1$ such that
		\begin{flalign}
			\label{ssts}
			\begin{aligned} 
				\mathbb{P} \Bigg[ & \bigcup_{t, t+s \in \mathcal{T}} \Big\{ \big| \mathsf{x}_k (t+s) - \mathsf{x}_k (t) \big| > 4 \big( k|s| \big)^{1/2} + 3 (2A+B) k^{1/3} |s| + k^{-50} \Big\} \\
				& \qquad \qquad \qquad \cap \textbf{MED}_k \big( [-3Ak^{1/3}, -3Ak^{1/3}]; -Bk^{2/3}; B k^{2/3} \big) \Bigg]  \le C_2 Ak^{201} e^{-(\log k)^3}.
			\end{aligned} 
		\end{flalign}
		
		\noindent Further observe for any real numbers $t, t+s \in [-2Ak^{1/3}, 2Ak^{1/3}]$ that, for sufficiently large $k$,
		\begin{flalign*}
			4| & ks|^{1/2} + 3(2A+B) k^{1/3} |s| + k^{-50} + 2k^{5} |t+s-t'-s'|^{1/3} + 2k^{5} |t-t'|^{1/3} \\
			& \le 4|ks|^{1/2} + 3(2A+B) k^{1/3} |s| + k^{-50} + 4 k^{-55} \le 4|ks|^{1/2} + 3(2A+B) k^{1/3} |s|+ k^{-25},
		\end{flalign*}
		
		\noindent where $t'$ and $t'+s'$ are the closest elements in $\mathcal{T}$ to $t'$ and $t'+s'$, respectively. Applying this, and using \eqref{ssts} to bound $\big| \mathsf{x}_k (t'+s') - \mathsf{x}_k (t') \big|$ and \eqref{tts2ak2} to bound $\big| \mathsf{x}_k (t+s) - \mathsf{x}_k (t'+s') \big|$ and $\big| \mathsf{x}_k (t') - \mathsf{x}_k (t) \big|$, yields the proposition.
	\end{proof}

	\begin{proof}[Proof of \Cref{probabilityeventr}]
		
		Throughout this proof, we set $r_1 = t - k^{1/3}$ and $r_2 = t + s + k^{1/3}$, and we also denote the event $\mathscr{E} = \textbf{MED}_k \big( \{ r_1, r_2 \}; -Bk^{2/3}; B k^{2/3}  \big)$. It suffices to show that 
		\begin{flalign}
			\label{xktsxktsk}
			\begin{aligned} 
				& \mathbb{P} \Big[ \big\{  \mathsf{x}_k (t+s) - \mathsf{x}_k (t) \le - 4(ks)^{1/2} - 3(A+B) k^{1/3} s- k^{-50} \big\} \cap \mathscr{E} \Big] \le C e^{-(\log k)^2}; \\
				& \mathbb{P} \Big[ \big\{  \mathsf{x}_k (t+s) - \mathsf{x}_k (t) \ge 4 (ks)^{1/2} + 3(A+B) k^{1/3} s + k^{-50} \big\} \cap \mathscr{E} \Big] \le C e^{-(\log k)^2}.
			\end{aligned} 
		\end{flalign}
		
		\noindent We only show the former bound in \eqref{xktsxktsk}, as the proof of the latter is entirely analogous (obtained by taking $r = r_1 = t - k^{1/3}$ below, instead of $r = r_2 = t + s + k^{1/3}$). To this end, set $r = t + s + k^{1/3}$; condition on $\mathcal{F}_{\ext} \big( \llbracket 1, k \rrbracket \times (t, r) \big)$; restrict to the event $\mathscr{E}$; and define the $k$-tuples $\bm{u} = \bm{\mathsf{x}}_{\llbracket 1, k \rrbracket} (t) \in \mathbb{W}_k$ and $\bm{v} = \bm{\mathsf{x}}_{\llbracket 1, k \rrbracket} (t) \in \mathbb{W}_k$, as well as the function $f = \mathsf{x}_{k+1} |_{[t, r]}$. Then the law of $\big( \mathsf{x}_j (t') \big)$, for $(j, t') \in \llbracket 1, k \rrbracket \times [t, r]$, is given by the non-intersecting Brownian bridge measure $\mathsf{Q}_f^{\bm{u}; \bm{v}}$. 
		
		Next set $u' = \mathsf{x}_k (t)$ and $v' = \mathsf{x}_k (r)$, and denote the $k$-tuples $\bm{u}' = (u', u', \ldots , u') \in \overline{\mathbb{W}}_k$ and $\bm{v}' = (v', v', \ldots , v') \in \overline{\mathbb{W}}_k$ (where $u'$ and $v'$ both appear with multiplicity $k$). Sample non-intersecting Brownian bridges $\bm{\mathsf{y}} = (\mathsf{y}_1, \mathsf{y}_2, \ldots , \mathsf{y}_k) \in \llbracket 1, k \rrbracket \times \mathcal{C} \big( [t, r] \big)$ from the measure $\mathsf{Q}^{\bm{u}; \bm{v}}$. Since $\bm{u} \ge \bm{u}'$ and $\bm{v} \ge \bm{v}'$, \Cref{monotoneheight} gives a coupling between $\bm{\mathsf{x}}$ and $\bm{\mathsf{y}}$ such that $\mathsf{x}_j (t') \ge \mathsf{y}_j (t')$, for each $(j, t') \in \llbracket 1, k \rrbracket \times [t, r]$. 
		
		From the second part of \Cref{estimatexj} (applied with the parameters $(n; a, b; u, v; t; D)$ there equal to $(k; t, r; u', v'; t + s; 50)$ here), there exists a constant $C > 1$ such that
		\begin{flalign*}
			\mathbb{P} \bigg[  \mathsf{y}_k (t+s) \le u' + \displaystyle\frac{s}{k^{1/3} + s} \cdot (v'-u') - 4(ks)^{1/2} - k^{-50} \bigg] \le C e^{-(\log k)^3},
		\end{flalign*}
		
		\noindent which together with the above coupling between $\bm{\mathsf{x}}$ and $\bm{\mathsf{y}}$ (with the facts that $u' = \mathsf{x}_k (t)$ and $v' = \mathsf{x}_k (r)$) yields
		\begin{flalign}
			\label{xkts2} 
			\mathbb{P} \bigg[ \mathsf{x}_k (t+s) - \mathsf{x}_k (t) \le -\displaystyle\frac{s}{k^{1/3} + s} \cdot \big| \mathsf{x}_k (r) - \mathsf{x}_k (t) \big| - 4(ks)^{1/2} - k^{-50} \bigg] \le C e^{-(\log k)^3}.
		\end{flalign}
		
		Since we are restricting to $\mathscr{E}$ and since $|t^2 - r^2| \le |r-t| \big( |t| + |r| \big) \le 3A k^{1/3} (k^{1/3} + s)$ (which holds since $r - t = k^{1/3} + s$ and $|t| + |r| \le 2Ak^{1/3} + k^{1/3} \le 3Ak^{1/3}$), we have  
		\begin{flalign*} 
			\big| \mathsf{x}_k (r) - \mathsf{x}_k (t) \big| \le 2^{-1/2} (t^2 - r^2) + 2 B k^{2/3} & \le 2 B k^{2/3} + 3 A k^{1/3} (k^{1/3} + s).
		\end{flalign*} 
		
		\noindent Inserting this into \eqref{xkts2}, we deduce the first bound in \eqref{xktsxktsk}. As mentioned previously, the proof of the second is very similar and thus omitted; this yields the lemma. 
	\end{proof} 
	
	\subsection{Likelihood of $\textbf{GAP}$}
	
	\label{ProofEventGap}
	
	In this section we establish \Cref{p:gapupperbound}, which will be a quick consequence of \Cref{monotonedifference} and \Cref{p:smallinitial}.

	\begin{proof}[Proof of \Cref{p:gapupperbound}]
		Recall the constant $C_1 (B) > 1$ from \Cref{p:smallinitial}, and let  $A =B \cdot C_1 (B) > B\geq 1$. We then restrict to the event $\mathscr{E} = \textbf{TOP} \big( \{ -Ak^{1/3}, Ak^{1/3} \}; Bk^{2/3} \big) \cap \textbf{MED}_k \big( \{ -Ak^{1/3}, Ak^{1/3} \}; -Bk^{2/3}; Bk^{2/3} \big)$. It suffices to show that for some constants $c = c(B) > 0$ and $C = C(B) > 1$ we have
		\begin{flalign}
			\label{xjxicb} 
			\mathbb{P} \Bigg[ \mathscr{E} \cap \bigcup_{|t| \le Ak^{1/3} / 2} \bigcup_{1 \le i < j \le \lfloor k/2 \rfloor} \big\{ \mathsf{x}_i (t) - \mathsf{x}_j (t) > C (j^{2/3} - i^{2/3}) + (\log k)^{25} i^{-1/3} \big\} \Bigg] \le c^{-1} e^{-c(\log k)^2}. 
		\end{flalign}
		
		To this end, define the $k$-tuples $\bm{u} = \bm{\mathsf{x}}_{\llbracket 1, k \rrbracket} (-Ak^{1/3}) \in \mathbb{W}_k$ and $\bm{v} = \bm{\mathsf{x}}_{\llbracket 1, k \rrbracket} (Ak^{1/3}) \in \mathbb{W}_k$. Then sample $k$ non-intersecting Brownian bridges $\bm{\mathsf{y}} = (\mathsf{y}_1, \mathsf{y}_2, \ldots , \mathsf{y}_k)$ from the measure $\mathsf{Q}^{\bm{u}; \bm{v}}$. Since the law of $\bm{\mathsf{x}}_{\llbracket 1, k \rrbracket}$ on $[-Ak^{1/3}, Ak^{1/3}]$ is given by $\mathsf{Q}_{\mathsf{x}_{k+1}}^{\bm{u}; \bm{v}}$, it follows from gap monotonicity \Cref{monotonedifference} that we may couple $\bm{\mathsf{x}}$ and $\bm{\mathsf{y}}$ such that 
		\begin{flalign}
			\label{xisyisdifference}
			\mathsf{x}_i (s) - \mathsf{x}_j (s) \le \mathsf{y}_i (s) - \mathsf{y}_j (s), \qquad \text{for any $1 \le i < j \le k$ and $s \in [-An^{1/3}, An^{1/3}]$}.
		\end{flalign}
		
		\noindent Further observe since we have restricted to $\mathscr{E}$ that $|u_1 + 2^{1/2} A^2 k^{2/3}| = \big| \mathsf{x}_1 (-Ak^{1/3}) + 2^{1/2} A^2 k^{2/3} \big| \le Bk^{2/3}$; by similar reasoning, we have $|u_k + 2^{1/2} A^2 k^{2/3}| \le Bk^{2/3}$, $|v_1 + 2^{1/2} A^2 k^{2/3}| \le Bk^{2/3}$, and $|v_k + 2^{1/2} A^2 k^{2/3}|\le Bk^{1/3}$. Hence, \Cref{p:smallinitial} (applied with the $\bm{\mathsf{x}}$ there equal to $\bm{\mathsf{y}}$ here, translated vertically by $2^{1/2} A^2 k^{2/3}$ and horizontally by $Ak^{1/3}$; the $A$ there equal to $2B$ here; and the $T$ there equal to $2A=2B\cdot C_1(B)$ here) yields constants $c = c(B) > 0$ and $C = C(B) > 1$ such that   
		\begin{flalign*}
			\mathbb{P} \Bigg[ \mathscr{E} \cap  \bigcup_{|t| \le Ak^{1/3} /2} \bigcup_{1\leq i < j\leq \lfloor k/2 \rfloor} \big\{ \mathsf{y}_i (tk^{1/3}) - \mathsf{y}_j (tk^{1/3})  > C (j^{2/3} - i^{2/3})  +  (& \log k)^{25} i^{-1/3} \big\} \Bigg] \\
			& \le  c^{-1} e^{-c (\log k)^2}. 
		\end{flalign*}
		
		\noindent This, together with \eqref{xisyisdifference} verifies \eqref{xjxicb} and thus the proposition.
	\end{proof}
	
	\subsection{Improved Medium Position Events} 
	
	\label{EventR} 
	
	Observe that the $\textbf{SCL}$ event from \Cref{eventscl} is the intersection of the medium position events $\textbf{MED}_k$ only for $k$ within a constant multiple of $n$. It will later be useful to have $\mathsf{x}_k$ be of order $-k^{2/3}$, for much larger values of $k$ (say, for $k \in [B^{-1} n, n^{100}]$). To this end, we define the following improvement of the medium position event.
	
	\begin{definition} 
		
		\label{eventimproved} 
		
		For any integer $n \ge 1$ and real numbers $A \ge 0$ and $B, C, R \ge 1$, define the \emph{improved medium position event} $\textbf{IMP}_n (A; B; C; R) = \textbf{IMP}_n^{\bm{\mathsf{x}}} (A; B; C; R)$ by setting\index{I@$\textbf{IMP}$; improved medium position event}
		\begin{flalign*}
			\textbf{IMP}_n (A; B; C; R) = \bigcap_{|t| \le An^{1/3}} \bigcap_{j = \lceil n/B \rceil}^{\lfloor Rn \rfloor} \big\{ C^{-1} n^{2/3} -Cj^{2/3} \le \mathsf{x}_j (t) \le C n^{2/3} -C^{-1} j^{2/3} \big\}. 
		\end{flalign*}
		
	\end{definition}

	What will later (in \Cref{EventProofB} below) be relevant for us is to have $\textbf{IMP}_n (A; B; C; R)$ hold when $R = n^D$ for some large (but uniformly bounded) $D > 1$. Observe, by \Cref{sclprobability}, that $\textbf{MED}_k$ is very likely if we restrict to the event $\textbf{TOP} \big( [-Ck^{1/3}, Ck^{1/3}]; C^{-1} k^{2/3} \big)$ for some sufficiently large constant $C > 1$. Thus $\textbf{IMP}$ would be very likely if we restricted to the intersection over, say $\log n$, many of these $\textbf{TOP}$ events (for example, for any $k$ equal to power of $2$ in $\llbracket n, n^{D+1} \rrbracket$); this would require us to take a union bound over $\log n$ many such events. Unfortunately, \Cref{l0} only indicates that for any given $k \in \llbracket B^{-1} n, n^{D+1} \rrbracket$, such a $\textbf{TOP}$ event  holds with probability $1 - \delta_k$ satisfying $\lim_{k \rightarrow \infty} \delta_k = 0$,  but without an effective rate. Thus, it is unclear if one can efficiently take a union bound over them. 
	
	The following proposition indicates that $\textbf{IMP}$ is very likely, upon restricting to only a uniformly bounded number (with respect to the index $k$) of $\textbf{TOP}$ and $\textbf{MED}$ events (for which a union bound can be taken). 
	
	\begin{prop}
		
		\label{probabilityimproved}
		
		For any real numbers $b \in ( 0, 1 / 4)$ and $A, B, D \ge 3$, there exist constants $c = c(A, b, B, D) > 0$, $C_1 = C_1 (B) > 1$, and $C_2 = C_2 (A, b, B) > 1$ such that 
		\begin{flalign}
			\label{probabilityestimateimproved} 
			\begin{aligned} 
				\mathbb{P} \Bigg[ & \textbf{\emph{IMP}}_n ( A; B; C_2; n^D)^{\complement}  \cap \bigcap_{|t| \le An^{1/3}} \big( \textbf{\emph{MED}}_{\lfloor n/4B \rfloor} (t; 2bn^{2/3}; Bn^{2/3}) \cap \textbf{\emph{TOP}} (t; bn^{2/3}) \big) \\
				&  \cap \bigcap_{t \in \{ -C_1 n^{10D}, C_1 n^{10D} \}} \big( \textbf{\emph{MED}}_{n^{30D}} (t; -Bn^{20D}; Bn^{20D}) \cap \textbf{\emph{TOP}} (t; Bn^{20D}) \big) \Bigg] \le c^{-1} e^{-c(\log n)^2}. 
			\end{aligned} 
		\end{flalign} 
		
	\end{prop}
	
	To show \eqref{probabilityestimateimproved}, we must show (on the $\textbf{TOP}$ and $\textbf{MED}$ events there) a lower and upper bound on $\mathsf{x}_k$, which amount to an upper and lower bound on $\mathsf{x}_1 - \mathsf{x}_k$, for any $k \in \llbracket B^{-1} n, n^{D+1} \rrbracket$. The upper bound on this difference will eventually follow from a $\textbf{GAP}$ event, which will be guaranteed by \Cref{p:gapupperbound}. To lower bound this difference, we will show that if $\mathsf{x}_1 - \mathsf{x}_k$ is too small, then a $\textbf{GAP}_k (t; \omega)$ event holds for a very small value of $\omega$ (equivalently, $\mathsf{x}_{\llbracket 1, k/2 \rrbracket} (t)$ has a very high density), which will contradict the $\textbf{MED}_{\lfloor n / 4B \rfloor}$ event in \eqref{probabilityestimateimproved}. 
	
	To make this precise, we begin with the following definition indicating when an $n$-tuple has a high density near its top entries, after applying Dyson Brownian motion for some time.
	
	\begin{definition}
		
		\label{eventpacked}
		
		Fix an integer $n \ge 1$; a bounded interval $(a, b) \subset \mathbb{R}_{\ge 0}$; and real numbers $\omega \ge 0$ and $\xi \in (0, 1)$. An $n$-tuple $\bm{u} = (u_1, u_2, \ldots , u_n) \in \mathbb{W}_n$ is called \emph{$\big( [a, b]; \omega; \xi \big)$-packed} if the following holds. Defining $\bm{\lambda} = (\lambda_1, \lambda_2, \ldots , \lambda_n) \in \llbracket 1, n \rrbracket \times \mathcal{C} (\mathbb{R}_{\ge 0})$ by letting $\bm{\lambda} (s)$ denote Dyson Brownian motion, run for time $s$, with initial data $\bm{\lambda} (0) = \bm{u}$, we have 
		\begin{flalign}
			\label{tlambdaj}
			\mathbb{P} \Bigg[ \bigcup_{s \in [a, b]} \bigcup_{1 \le j < k \le \lfloor n/2\rfloor}  \big\{ \lambda_j (sn^{1/3}) - \lambda_k (sn^{1/3}) > \omega (k^{2/3} - j^{2/3}) + (\log n)^{25} j^{-1/3} \big\} \Bigg] \le  \xi^{-1} e^{-\xi (\log n)^2}.
		\end{flalign}
		
		\noindent Given a line ensemble $\bm{\mathsf{x}} = (\mathsf{x}_1, \mathsf{x}_2, \ldots ) \in \mathbb{Z}_{\ge 1} \times \mathcal{C} (\mathbb{R})$ and real number $t \in \mathbb{R}$, let\index{P@$\textbf{PAC}$; packed event} $\textbf{PAC}_n \big(t; [a, b]; \omega; \xi \big)= \textbf{PAC}_n^{\bm{\mathsf{x}}} \big(t; [a, b]; \omega; \xi \big)$ denote the event that $\bm{\mathsf{x}}_{\llbracket 1, n \rrbracket} (tn^{1/3}) \in \mathbb{W}_n$ is $\big([a, b]; \omega; \xi\big)$-packed. 
		
	\end{definition}

	Since the probability on the left side of \eqref{tlambdaj} is continuous in the initial data $\bm{\lambda} (0) = \bm{u}$, the set of $\bm{\mathsf{x}}$ satisfying $\textbf{PAC}_n^{\bm{\mathsf{x}}} \big( t; [a, b]; \omega; \xi) $ defines a closed (and thus measurable) set.
	
	\begin{rem}
		
		\label{gappacked}
		
		Let $t \in \mathbb{R}$; $\omega, b \ge 0$; and $\xi \in (0, 1)$ be real numbers. Observe that we have $\textbf{PAC}_n \big(t; [0, b]; \omega, \xi \big) \subseteq \textbf{GAP}_{\lfloor n/2 \rfloor} (tn^{1/3}; \omega)$, for sufficiently large $n$. Indeed, $\bm{\mathsf{x}}_{\llbracket 1, n \rrbracket} (t)$ is $\big( [0,b]; \omega; \xi \big)$-packed on $\textbf{PAC}_n \big(t; [0, b]; \omega; \xi \big)$. Since the $s = 0 \in [0, b]$ case of the event in \eqref{tlambdaj} is deterministic for a given $\bm{u}$, it follows (for $n$ sufficiently large so that $\xi^{-1} e^{-\xi (\log n)^2} < 1$) that $\mathsf{x}_j (tn^{1/3}) - \mathsf{x}_k (tn^{1/3}) \le \omega (k^{2/3} - j^{2/3}) + (\log n)^{25} j^{-1/3}$ for each $1 \le j < k \le \lfloor n / 2 \rfloor$, and hence $\textbf{GAP}_{\lfloor n/2 \rfloor} (tn^{1/3}; \omega)$ holds.  
		
	\end{rem}

	The next lemma, which is a quick consequence of \Cref{initialsmall2}, indicates that an $n$-tuple $\bm{u} = (u_1, u_2, \ldots , u_n)$ with $u_1 - u_n$ sufficiently small is packed.
	
	\begin{lem} 
		
		\label{packedu1unsmall}
		
		For any real number $\omega > 0$, there exist real numbers $a = a(\omega) > 0$, $c = c(\omega) > 0$, and $\xi = \xi(\omega) \in (0, 1)$ such that the following holds. Let $n \ge 1$ be an integer and $\bm{u} = (u_1, u_2, \ldots , u_n) \in \mathbb{W}_n$ be an $n$-tuple with $u_1 - u_n < cn^{2/3}$. Then, $\bm{u}$ is $\big( [a, 5a]; \omega; \xi \big)$-packed.
	\end{lem} 
	
	\begin{proof}
		
		Let $C_0 > 1$ denote the constant $C$ from \Cref{initialsmall2}; fix a real number $B > 1$ such that $C_0 B^{-1/2} \le \omega$; let $c_0 = c_0 (B) > 1$ be the constant $c(5B)$ from \Cref{initialsmall2}; and set $a = (5B)^{-1}$ and $\xi = c = c_0$. If $u_1 - u_n < cn^{2/3}$, then defining $\bm{\lambda} = (\lambda_1, \lambda_2, \ldots , \lambda_n) \in \llbracket 1, n \rrbracket \times \mathcal{C} (\mathbb{R}_{\ge 0})$ by letting $\bm{\lambda}(s)$ denote Dyson Brownian motion with initial data $\bm{\lambda} (0) = \bm{u}$, \Cref{initialsmall2} yields 
		\begin{flalign*}
			\mathbb{P} \Bigg[ \bigcup_{t \in [a, 5a]} \bigcup_{1 \le j < k \le \lfloor n/2 \rfloor} \big\{ \lambda_j (tn^{1/3}) - \lambda_k (tn^{1/3}) \ge \omega (k^{2/3} - j^{2/3}) + ( \log n & )^{-20} j^{1/3} \big\} \Bigg] \\
			& \le \xi^{-1} e^{-\xi (\log n)^2},
		\end{flalign*} 
		
		\noindent where we used the facts that $[a, 5a] = \big[ (5B)^{-1}, B^{-1} \big] \subseteq \big[ (5B)^{-1}, 5B \big]$ and that $C_0 t^{1/2} \le C_0 B^{-1/2} \le \omega$ if $ t \le B^{-1}$; this verifies that $\bm{u}$ is $\big( [a, 5a]; \omega; \xi \big)$-packed.
	\end{proof}
	
	The following proposition indicates that, if the line ensemble $\bm{\mathsf{x}}$ is packed at some time $t_0$ (and some weak $\textbf{GAP}$ event holds), then $\bm{\mathsf{x}}$ is packed on certain other time intervals; see the left side of \Cref{f:pac}. We establish it in \Cref{ConfigurationPacked} below. 

	\begin{figure}
	\center
\includegraphics[scale = .75, trim=0 0.5cm 0 0.5cm, clip]{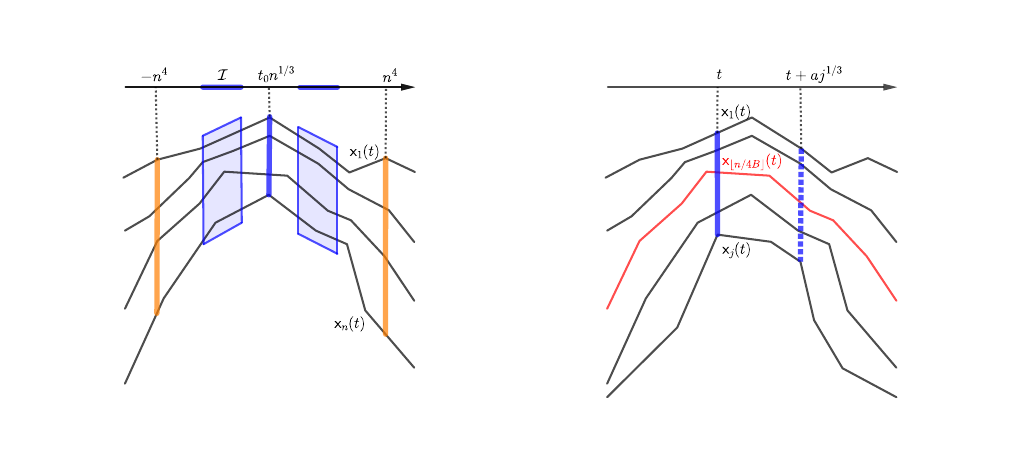}

\caption{Shown to the left is a depiction of \Cref{kgapc}, which states that if $\bm{\mathsf{x}}$ is packed at some time $t_0n^{1/3}$ (the blue line) and some weak $\textbf{GAP}$ event holds at times $\pm n^4$ (the orange lines), then $\bm{\mathsf{x}}$ is packed on certain other time intervals (the blue regions). Shown on the right is the setup for the proof of \Cref{probabilityimproved}.
}
\label{f:pac}
	\end{figure}
	
	\begin{prop}
		
		\label{kgapc} 
		
		For any real numbers $b \ge a \ge 0$ with $b \ge 2a$; $\omega > 0$; and $\xi \in (0, 1)$, there exists a constant $c = c (a, b, \omega, \xi) > 0$ such that the following holds. Let $n \ge 1$ be an integer and $t_0 \in [-n, n]$ be a real number. Abbreviate the event $\mathscr{A} = \textbf{\emph{PAC}}_n \big( t_0; [a,b]; \omega; \xi \big)$. Defining the interval $\mathcal{I}_1 = [t_0 - 2a, t_0 - a]$, we have for any $t \in \mathcal{I}_1$ that   
		\begin{flalign}
			\label{eventa10} 
			\begin{aligned} 
			\mathbb{P}\bigg[ & \mathscr{A} \cap \textbf{\emph{GAP}}_n ( -n^4; n)  \cap \textbf{\emph{PAC}}_n \Big( t ; [0, b-2a]; 2 \omega; \frac{\xi}{2} \Big)^{\complement} \bigg] \le c^{-1} e^{-c(\log n)^2}; \\ 
			\mathbb{P}\Big[ & \mathscr{A} \cap \textbf{\emph{GAP}}_n ( -n^4; n)  \cap \textbf{\emph{GAP}}_{\lfloor n/2 \rfloor} \big( n^{1/3} \cdot \mathcal{I}_1; 2\omega \big)^{\complement} \Big] \le c^{-1} e^{-c(\log n)^2}.
			\end{aligned} 
		\end{flalign}
	
		\noindent Similarly, defining the interval $\mathcal{I}_2 = [t_0 + a, t_0 + 2a]$, we have for any $t \in \mathcal{I}_2$ that 
			\begin{flalign}
			\label{eventa20} 
			\begin{aligned} 
			\mathbb{P}\bigg[ & \mathscr{A} \cap \textbf{\emph{GAP}}_n ( n^4 ; n )  \cap \textbf{\emph{PAC}}_n \Big( t ; [0, b-2a]; 2 \omega; \frac{\xi}{2} \Big)^{\complement} \bigg] \le c^{-1} e^{-c(\log n)^2}; \\ 
			\mathbb{P}\Big[ & \mathscr{A} \cap \textbf{\emph{GAP}}_n ( n^4 ; n )  \cap \textbf{\emph{GAP}}_{\lfloor n/2 \rfloor} \big( n^{1/3} \cdot \mathcal{I}_2; 2\omega \big)^{\complement} \Big] \le c^{-1} e^{-c(\log n)^2}.
			\end{aligned} 
		\end{flalign}
		
	\end{prop}

	Now we can establish \Cref{probabilityimproved}. 
	
	\begin{proof}[Proof of \Cref{probabilityimproved}]
		
		Let $C_1 \ge B \ge 2$ denote the constant $A(B)$ from \Cref{p:gapupperbound}, and define the events 
		\begin{flalign}
			\label{e1e2e3} 
			\begin{aligned}
				& \quad \mathscr{E}_1 = \textbf{MED}_{n^{30D}} \big( \{ -C_1 n^{10D}, C_1 n^{10D} \}; -Bn^{20D}; Bn^{20D} \big) \cap \textbf{TOP} \big( \{ - C_1 n^{10D}, C_1 n^{10D} \}; Bn^{20D} \big); \\
				& \mathscr{E}_2 = \textbf{TOP} \big( [-An^{1/3}, An^{1/3}]; bn^{2/3} \big); \quad \mathscr{E}_3 = \textbf{MED}_{\lfloor n/4B \rfloor} \big( [ -An^{1/3}, An^{1/3}]; 2bn^{2/3}; Bn^{2/3} \big).
			\end{aligned} 
		\end{flalign}
		
		\noindent By \Cref{eventimproved} and a union bound, it suffices to show that there exist constants $c = c(A, b, B, D) > 0$ and $C_2 = C_2 (A, b, B) > 1$ such that for each $j \in \llbracket B^{-1} n, n^{D+1} \rrbracket$ we have 
		\begin{flalign}
			\label{xk12} 
			\begin{aligned} 
				& \mathbb{P} \Bigg[ \bigcup_{|t| \le An^{1/3}} \big\{ \mathsf{x}_j (t) > C_2^{-1} n^{2/3} -C_2 j^{2/3} \big\} \cap \mathscr{E}_1 \cap \mathscr{E}_2  \Bigg] \le c^{-1} e^{-c (\log n)^2}; \\
				& \mathbb{P} \Bigg[ \bigcup_{|t| \le An^{1/3}} \big\{ \mathsf{x}_j (t) > C_2 n^{2/3} -C_2^{-1} j^{2/3} \big\} \cap \mathscr{E}_1 \cap \mathscr{E}_2 \cap \mathscr{E}_3 \Bigg] \le c^{-1} e^{-c (\log n)^2}; 
			\end{aligned}	
		\end{flalign}

		To this end, first observe by restricting to $\mathscr{E}_1$ that \Cref{p:gapupperbound} (with the $k$ there equal to $n^{30D}$ here) yields $c_1 = c_1 (B) > 0$ and $C_3 = C_3 (B) > 1$ such that   
		\begin{flalign}
			\label{probability4e} 
			\mathbb{P} [\mathscr{E}_4^{\complement} \cap \mathscr{E}_1 \big] \le c_1^{-1} e^{-c_1 (\log n)^2}, \qquad \text{where} \quad \mathscr{E}_4 = \textbf{GAP}_{\lfloor n^{30D} / 2 \rfloor} \big( [-n^{10D}, n^{10D}]; C_3 \big).
		\end{flalign}
		
		\noindent By \Cref{gap} of the $\textbf{GAP}$ event we have, for each $(k, t) \in \llbracket B^{-1} n, n^{D+1} \rrbracket \times [ -An^{1/3}, An^{1/3} ]$, that
		\begin{flalign}
			\label{x1txkt} 
			\mathsf{x}_1 (t) - \mathsf{x}_k (t) \le C_3 k^{2/3} + (\log n)^{25} \le 2C_3 k^{2/3}, \qquad \text{on the event $\mathscr{E}_4$}.
		\end{flalign}
		
		\noindent Next, on the event $\mathscr{E}_2$, we have for $t \in [-An^{1/3}, An^{1/3}]$ that
		\begin{flalign}
			\label{x1tn23} 
			-(A^2 + 1) n^{2/3} \le -b n^{2/3} -2^{-1/2} t^2\le \mathsf{x}_1 (t) \le bn^{2/3}-2^{-1/2} t^2  \le \displaystyle\frac{n^{2/3}}{4}, 
		\end{flalign}
		
		\noindent where in the first and last inequalities we used the fact that $b < 1 / 4$. Together with \eqref{x1txkt}, it follows for each $(k, t) \in \llbracket B^{-1} n, n^{D+1} \rrbracket \times [ -An^{1/3}, An^{1/3}]$ that 
		\begin{flalign*}
			\mathsf{x}_k (t) \ge - (A^2 + 1) n^{2/3} - 2C_3 k^{2/3} \ge -(2C_3 + A^2 + 1) B^{2/3} k^{2/3}, \qquad \text{on the event $\mathscr{E}_2 \cap \mathscr{E}_4$},
		\end{flalign*}  
		
		\noindent which together with \eqref{probability4e} establishes the first bound in \eqref{xk12} (by taking the $C_2$ there sufficiently large so that $C_2^{-1} n^{2/3} - C_2 k^{2/3} \le -(2C_3 + A^2  +1) B^{2/3} k^{2/3}$ for $k \ge B^{-1} n$). 
		
		To establish the second bound there, observe from the upper bound in \eqref{x1tn23} that it suffices to show that there exist constants $a = a(b) > 0$, $c = c(A, b, B, D) > 0$, and $C_2 = C_2 (A, b, B) > 1$ such that for each $(j, t) \in \llbracket B^{-1} n, n^{D+1} \rrbracket \times [1 -An^{1/3}, A n^{1/3} + 1]$ we have 
		\begin{flalign}
			\label{x1txjt1} 
			\mathbb{P} \bigg[ \bigcup_{s \in [t - 1, t + 1]} \big\{ \mathsf{x}_1 (s) - \mathsf{x}_j (s) < C_2^{-1} j^{2/3} \big\} \cap \mathscr{E}_1 \cap \mathscr{E}_2 \cap \mathscr{E}_3 \bigg] \le c^{-1} e^{-c (\log n)^2}.
		\end{flalign}
		
		\noindent Let us briefly outline how we will proceed.  First, we will apply \Cref{packedu1unsmall} to show that if $\mathsf{x}_1 (s) - \mathsf{x}_j (s) \le C_2^{-1} j^{2/3}$ for some $s \in [t - an^{1/3}, t + an^{1/3}]$, then $\bm{\mathsf{x}}_{\llbracket 1, j \rrbracket} (s)$ is packed. We will then use \Cref{kgapc} to deduce that $\bm{\mathsf{x}}_{\llbracket 1, j \rrbracket} (t + 2aj^{1/3})$ is packed for some constant $a > 0$, and then apply \Cref{kgapc} again (with \Cref{gappacked}) to deduce that $\bm{\mathsf{x}}_{\llbracket 1, j/2 \rrbracket} (t)$ has small gaps. The latter will contradict the $\textbf{TOP} \cap \textbf{MED}$ event defining $\mathscr{E}_2 \cap \mathscr{E}_3$ (recall \eqref{e1e2e3}). See the right side of \Cref{f:pac}.
		
		Now let us implement this in more detail. Define $C_2 = c(b)^{-1}$, where $c(b)$ is the constant from \Cref{packedu1unsmall} (with $\omega$ there taken to be $b$ here). We further fix some $(j, t) \in \llbracket B^{-1} n, n^{D+1} \rrbracket \times [-An^{1/3}, An^{1/3}]$ and define the event
		\begin{flalign}
			\label{evente5} 
			\mathscr{E}_5 = \mathscr{E}_5 (j, t) = \bigcup_{s \in [t-1, t + 1]} \big\{ \mathsf{x}_1 (s) - \mathsf{x}_j (s) < C_2^{-1} j^{2/3} \big\}.
		\end{flalign}
		
		\noindent We restrict to $\mathscr{E}_5$ in what follows, so that there exists some time $s_0 \in [t-1, t+1]$ such that $\mathsf{x}_1 (s_0) - \mathsf{x}_j (s_0) < C_2^{-1} j^{2/3}$; assume that $s_0$ is the smallest such time in $[t - an^{1/3}, t+an^{1/3}]$, so that $(s_0, n^4)$ is a $\llbracket 1, j \rrbracket$-stopping domain, in the sense of \Cref{property2}. By \Cref{packedu1unsmall}, there exist constants $a = a(b) > 0$ and $\xi = \xi (b) > 0$ such that $\mathscr{E}_5 \subseteq \textbf{PAC}_j \big( s_0; [2a, 10 a]; b; \xi \big)$. Hence, applying the first statement of \eqref{eventa20} in \Cref{kgapc} with the $(n; t_0; a, b; \omega; \xi)$ there equal to $(j; j^{-1/3} s_0; 2a, 10a; b; \xi)$ here yields a constant $c_2 = c_2 (A, b, B, D) > 0$ such that 
		\begin{flalign}
			\label{e5e4packed} 
			\mathbb{P} \bigg[ \mathscr{E}_5 \cap \mathscr{E}_4 \cap \textbf{PAC}_j \Big( tj^{-1/3} + 3a; [0, 6a]; 2b; \displaystyle\frac{\xi}{2} \Big)^{\complement} \bigg] \le c_2^{-1} e^{-c_2 (\log n)^2}.
		\end{flalign}
		
		\noindent Here, we have used the facts that $\mathscr{E}_4 \subseteq \textbf{GAP}_j (j^4; j )$ for $j \in \llbracket B^{-1} n, n^{D+1} \rrbracket$ (due to \eqref{probability4e}, the fact that $ j^4 \in [An^{1/3}, n^{10D}]$, the fact that $C_1 (k^{2/3} - i^{2/3}) + (30 D \log n)^{25} i^{-1/3} \le j (k^{2/3} - i^{2/3}) + (\log j)^{25} i^{-1/3}$ for any $1 \le i < k \le j$, and the fact that $n$ is sufficiently large), and the fact that $t + 3aj^{1/3} \in [s_0 + 2aj^{1/3}, s_0 + 4aj^{1/3}]$ for $s_0 \in [t-1, t+1]$ (as $n$ is sufficiently large).
		
		Applying the second statement of \eqref{eventa10} in \Cref{kgapc} with the parameters $(n; t_0; a, b; \omega; \xi)$ there equal to $( j; tj^{-1/3} + 3a; 2a, 6a; 2b; \xi / 2)$ here (now using the fact that $\mathscr{E}_4 \subseteq \textbf{GAP}_j ( -j^4; j)$) yields a constant $c_3 = c_3 (A, b, B, D) > 0$ such that 
		\begin{flalign*}
			\mathbb{P} \bigg[ \textbf{PAC}_j \Big( tj^{-1/3} + 3a; [0, 6a]; 2b; \displaystyle\frac{\xi}{2} \Big) \cap \mathscr{E}_4 \cap \textbf{GAP}_{\lfloor j/2 \rfloor} \big( [ t - aj^{1/3}, t + aj^{1/3}]; 4b \big)^{\complement} \bigg] \le c_3^{-1} e^{-c_3 (\log n)^2}.
		\end{flalign*}
		
		\noindent Together with \eqref{e5e4packed} and \eqref{probability4e}, and a union bound, this yields a constant $c_4 = c_4 (A, b, B, D) > 0$ such that 
		\begin{flalign}
			\label{x1txjt2} 
			\mathbb{P} \big[ \mathscr{E}_5 \cap \mathscr{E}_1\cap \textbf{GAP}_{\lfloor j/2 \rfloor} ( t; 4b \big)^{\complement} \big] \le c_4^{-1} e^{-c_4 (\log n)^2}. 
		\end{flalign}
		
		\noindent We then claim that $\mathscr{E}_2 \cap \mathscr{E}_3 \subseteq \textbf{GAP}_{\lfloor j/2 \rfloor} (t; 4b)^{\complement}$. To this end, first observe on $\textbf{GAP}_{\lfloor j/2 \rfloor} (t; 4b)$ that   
		\begin{flalign*}
			\mathsf{x}_1 (t) - \mathsf{x}_{\lfloor n / 4B \rfloor} (t) \le 4b \Big( \displaystyle\frac{n}{4B}\Big)^{2/3} + (\log n)^{25} \le bn^{2/3},
		\end{flalign*}
		
		\noindent where we used the facts that $B \ge 3$ and $n$ is sufficiently large. However, by \Cref{eventsregular1} and \eqref{e1e2e3}, we have $\mathsf{x}_1 (t) - \mathsf{x}_{\lfloor n/4B \rfloor} (t) \ge bn^{2/3}$ on $\mathscr{E}_2 \cap \mathscr{E}_3$, meaning that $\mathscr{E}_2 \cap \mathscr{E}_3 \subseteq \textbf{GAP}_{\lfloor j/2 \rfloor} (t; 4b)^{\complement}$. This, together with \eqref{x1txjt2}, establishes \eqref{x1txjt1} and thus the proposition.		
	\end{proof} 
	
	\subsection{Proof of \Cref{kgapc} } 
	
	\label{ConfigurationPacked}

 In this section we prove \Cref{kgapc}. We begin with the following lemma, which and shows that any $n$-tuple whose gaps are approximately upper bounded by a packed $n$-tuple is also packed. 
 
	\begin{lem}  
		
		\label{gap2} 
		
		Adopt the notation of \Cref{eventpacked}; fix a real number $\varsigma \ge 0$, and suppose that $\bm{u} \in \mathbb{W}_n$ is $\big( [a, b]; \omega; \xi \big)$-packed. Then any $n$-tuple $\widetilde{\bm{u}} \in \mathbb{W}_n$ such that $\widetilde{u}_j - \widetilde{u}_{j+1} \le  u_j - u_{j+1} + \varsigma / (4n^{4/3})$ for each integer $j \in \llbracket 1, n \rrbracket$ is $\big( [a, b]; \omega + \varsigma; \xi \big)$-packed. 
		
	\end{lem}

	Its proof uses the following two lemmas; the first and second are counterparts of height and gap monotonicity for Dyson Brownian motion, respectively. We omit their proofs, which are quick consequences of \Cref{uvv} and \Cref{monotonedifference}, with \Cref{tobridge} (the latter taken as $\mathsf{T}$ tends to $\infty$); see also \cite{ERGIM} for the second lemma below. 
	
	\begin{lem}
		
		\label{couplemotion}
		
		Let $n \ge 1$ be an integer; $\varsigma \ge 0$ be a real number; and $\bm{u}, \widetilde{\bm{u}} \in \mathbb{W}_n$ be $n$-tuples such that $\max_{j \in \llbracket 1, n \rrbracket} \big| u_j - \widetilde{u}_j \big| \le \varsigma$. Define $\bm{\lambda} = (\lambda_1, \lambda_2, \ldots , \lambda_n) \in \llbracket 1, n \rrbracket \times \mathcal{C} (\mathbb{R}_{\ge 0})$ and $\widetilde{\bm{\lambda}} = \big( \widetilde{\lambda}_1, \widetilde{\lambda}_2, \ldots , \widetilde{\lambda}_n \big) \in \llbracket 1, n \rrbracket \times \mathcal{C} (\mathbb{R}_{\ge 0})$ by letting $\bm{\lambda} (s)$ and $\widetilde{\bm{\lambda}} (s)$ denote Dyson Brownian motions, run for time $s$, with initial data $\bm{\lambda} (0) = \bm{u}$ and $\widetilde{\bm{\lambda}} (0) = \widetilde{\bm{u}}$, respectively. Then, there exists a coupling between $\bm{\lambda}$ and $\widetilde{\bm{\lambda}}$ such that $\big| \lambda_j (s) - \widetilde{\lambda}_j (s) \big| \le \varsigma$ for each $(j, s) \in \llbracket 1, n \rrbracket \times \mathbb{R}_{\ge 0}$.
	\end{lem}
	
	\begin{lem}[{\cite[Lemma 2.8(2)]{ERGIM}}]
		
		\label{gapmotion}
		
		Let $n \ge 1$ be an integer and $\bm{u}, \widetilde{\bm{u}} \in \mathbb{W}_n$ be $n$-tuples such that $u_j - u_{j+1} \le \widetilde{u}_j - \widetilde{u}_{j+1}$ for each $j \in \llbracket 1, n - 1 \rrbracket$. Define $\bm{\lambda} = (\lambda_1, \lambda_2, \ldots , \lambda_n) \in \llbracket 1, n \rrbracket \times \mathcal{C} (\mathbb{R}_{\ge 0})$ and $\widetilde{\bm{\lambda}} = \big( \widetilde{\lambda}_1, \widetilde{\lambda}_2, \ldots , \widetilde{\lambda}_n \big) \in \llbracket 1, n \rrbracket \times \mathcal{C} (\mathbb{R}_{\ge 0})$ by letting $\bm{\lambda} (s)$ and $\widetilde{\bm{\lambda}} (s)$ denote Dyson Brownian motions, run for time $s$, with initial data $\bm{\lambda} (0) = \bm{u}$ and $\widetilde{\bm{\lambda}} (0) = \widetilde{\bm{u}}$, respectively. Then, there exists a coupling between $\bm{\lambda}$ and $\widetilde{\bm{\lambda}}$ such that $\lambda_j (s) - \lambda_{j+1} (s) \le \widetilde{\lambda}_j (s) - \widetilde{\lambda}_{j+1} (s)$ for each $(j, s) \in \llbracket 1, n \rrbracket \times \mathbb{R}_{\ge 0}$.
	\end{lem}

	\begin{proof}[Proof of \Cref{gap2}]
		
		Define the $n$-tuple $\widehat{\bm{u}} = (\widehat{u}_1, \widehat{u}_2, \ldots , \widehat{u}_n) \in \mathbb{W}_n$ by setting $\widehat{u}_j = u_j  + j \varsigma / (4n^{4/3})$ for each $j \in \llbracket 1, n \rrbracket$. Define the processes $\bm{\lambda} = (\lambda_1, \lambda_2, \ldots , \lambda_n) \in \llbracket 1, n \rrbracket \times \mathcal{C} (\mathbb{R}_{\ge 0})$, $\widetilde{\bm{\lambda}} = \big( \widetilde{\lambda}_1, \widetilde{\lambda}_2, \ldots , \widetilde{\lambda}_n \big) \in \llbracket 1, n \rrbracket \times \mathcal{C} (\mathbb{R}_{\ge 0})$, and $\widehat{\bm{\lambda}} = \big( \widehat{\lambda}_1, \widehat{\lambda}_2, \ldots , \widehat{\lambda}_n \big)$ by letting $\bm{\lambda} (s)$, $\widetilde{\bm{\lambda}} (s)$, and $\widehat{\bm{\lambda}} (s)$ denote Dyson Brownian motions with initial data $\bm{\lambda} (0) = \bm{u}$, $\widetilde{\bm{\lambda}} (0) = \widetilde{\bm{u}}$, and $\widehat{\bm{\lambda}} (0) = \widehat{\bm{u}}$, respectively.

		Since $\widetilde{u}_j - \widetilde{u}_{j+1} \le u_j - u_{j+1} + \varsigma / (4n^{4/3}) = \widehat{u}_j - \widehat{u}_{j+1}$, \Cref{gapmotion} yields a coupling between $\widetilde{\bm{\lambda}}$ and $\widehat{\bm{\lambda}}$ such that $\widetilde{\lambda}_j (s) - \widetilde{\lambda}_{j+1} (s) \le \widehat{\lambda}_j (s) - \widehat{\lambda}_{j+1} (s)$, for each $(j,s) \in \llbracket 1, n \rrbracket \times \mathbb{R}_{\ge 0}$. Since we also have $|\widehat{u}_j - u_j| \le j \varsigma / (4n^{4/3}) \le \varsigma /(4n^{1/3})$, \Cref{couplemotion} yields a coupling between $\bm{\lambda}$ and $\widehat{\bm{\lambda}}$ such that $\big| \lambda_j (s) - \widehat{\lambda}_j (s) \big| \le \varsigma / (4n^{1/3})$ for each $(j, s) \in \llbracket 1, n \rrbracket \times \mathbb{R}_{\ge 0}$. Combining these couplings, it follows for any $1 \le j < k \le n$ and $s \in \mathbb{R}_{\ge 0}$ that 
		\begin{flalign*} 
			\widetilde{\lambda}_k (s) - \widetilde{\lambda}_j (s) \le \widehat{\lambda}_j (s) - \widehat{\lambda}_k (s) \le \lambda_k (s) - \lambda_j (s) + \displaystyle\frac{\varsigma}{2n^{1/3}} \le \lambda_k (s) - \lambda_j (s) + \varsigma (k^{2/3} - j^{2/3}),
		\end{flalign*} 
		
		\noindent where in the last inequality we used the bound $k^{2/3} - j^{2/3} \ge 1 / (2n^{1/3})$ for $1 \le j < k \le n$. Together with \eqref{tlambdaj}, this implies that $\widetilde{\bm{u}}$ is $\big( [a, b]; \omega + \varsigma; \xi \big)$-packed. 
	\end{proof}

	The following lemma indicates that applying Dyson Brownian motion to a packed $n$-tuple likely yields a packed $n$-tuple.
	
	\begin{lem}
		
		\label{pacxi2} 
		
		Adopt the notation of \Cref{eventpacked}, assuming that $\bm{u}$ is $\big( [a, b]; c; \omega \big)$-packed. For any real number $s_0 \in [0, b]$ the $n$-tuple $\bm{\lambda} (s_0 n^{1/3})$ is $\big( [ \min \{ a-s_0, 0 \}, b-s_0]; \xi/2; \omega \big)$-packed with probability at least $1 - 2\xi^{-1} e^{-\xi (\log n)^2 / 2}$. 
		
	\end{lem}

	\begin{proof}
		
		Let us assume in what follows that $s_0 \le a$, as the proof when $s_0 \in (a, b]$ is entirely analogous. Since $\bm{u}$ is $\big( [a, b]; c; \omega \big)$-packed, we have 
		\begin{flalign*}
			\mathbb{P} \Bigg[ \bigcup_{t \in [a-s_0, b-s_0]} \bigcup_{1 \le j < k \le \lfloor n/2 \rfloor} \Big\{ \lambda_j \big( & (t+s_0) n^{1/3} \big) - \lambda_k \big( (t+s_0) n^{1/3} \big) \\
			& > \omega (k^{2/3} - j^{2/3})  + (\log n)^{25} j^{-1/3} \Big\} \Bigg] \le  \xi^{-1} e^{-\xi (\log n)^2}.
		\end{flalign*}
		
		\noindent  Together with the Markov estimate \Cref{fg0g}, this yields the existence of an event $\mathscr{E}$ measurable with respect to $\bm{\lambda} (s_0 n^{1/3})$, with $\mathbb{P} [\mathscr{E}] \ge 1 - (\xi^{-1} e^{-\xi (\log n)^2})^{1/2} \ge 1 - 2 \xi^{-1} e^{-\xi (\log n)^2 / 2}$, such that the following holds. Conditioning on $\bm{\lambda} (s_0 n^{1/3})$ and restricting to $\mathscr{E}$, we have
		\begin{flalign*}
			\mathbb{P} \Bigg[ \bigcup_{t \in [a-s_0, b-s_0]} & \bigcup_{1 \le j < k \le \lfloor n/2 \rfloor} \Big\{ \lambda_1 \big( (t+s_0) n^{1/3} \big) - \lambda_{\lfloor n/2 \rfloor} \big((t+s_0) n^{1/3} \big) \\
			& \ge \omega (k^{2/3} - j^{2/3}) + (\log n)^{25} j^{-1/3} \Big\} \Bigg] \le (\xi^{-1} e^{-\xi (\log n)^2})^{1/2}  \le 2 \xi^{-1} e^{-\xi (\log n)^2 / 2}. 
		\end{flalign*}
		
		\noindent This, together with the fact that $\bm{\lambda} (s+s_0 n^{1/3})$ has the same law as Dyson Brownian motion run for time $s$ with initial data $\bm{\lambda} (s_0 n^{1/3})$ for $s \ge 0$ (and \Cref{eventpacked}), implies the lemma.
	\end{proof}

		\begin{figure}
	\center
\includegraphics[scale = 1.25, trim=0 0.5cm 0 0.5cm, clip]{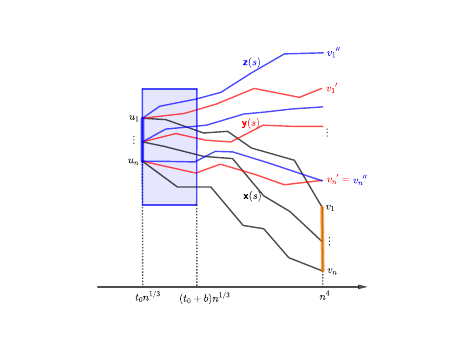}

\caption{Shown above is the setup for the proof of \Cref{kgapc}.}
\label{f:DBM_coupling}
	\end{figure}

	Now we can establish \Cref{kgapc}. 
	
	\begin{proof}[Proof of \Cref{kgapc}]

		We only show \eqref{eventa10}, as the proof of \eqref{eventa20} is entirely analogous. So, assuming that $t \in \mathcal{I}_2$ and denoting the event $\mathscr{E} = \mathscr{A} \cap \textbf{GAP}_n ( n^4 ; n )$, it then suffices to show 
		\begin{flalign}
			\label{xjtb} 
			\begin{aligned}
			& \mathbb{P} \bigg[ \textbf{PAC}_n \Big( t; [0, b-2a]; 2 \omega; \frac{\xi}{2} \Big)^{\complement} \cap \mathscr{E} \Big] \le c^{-1} e^{-c(\log n)^2}; \\ 
			& \mathbb{P} \Big[ \textbf{GAP}_{\lfloor n/2 \rfloor} \big( n^{1/3} \cdot \mathcal{I}_2; 2\omega \big)^{\complement} \cap \mathscr{E} \Big] \le  c^{-1} e^{-c (\log n)^2}. 
			\end{aligned}
		\end{flalign}
		
		To this end, restrict to $\mathscr{E}$ and denote the $n$-tuples $\bm{u} = \bm{\mathsf{x}}_{\llbracket 1, n \rrbracket} (t_0 n^{1/3}) \in \mathbb{W}_n$ and $\bm{v} = \bm{\mathsf{x}}_{\llbracket 1, n \rrbracket} (n^4) \in \mathbb{W}_n$. Define the line ensemble $\bm{\mathsf{y}} = (\mathsf{y}_1, \mathsf{y}_2, \ldots , \mathsf{y}_n) \in \llbracket 1, n \rrbracket \times \mathcal{C} \big( [t_0 n^{1/3}, n^4] \big)$ by, for all $s \in [t_0 n^{1/3}, n^4]$, letting $\bm{\mathsf{y}} (s)$ denote Dyson Brownian motion (recall \Cref{LambdaEquation}) run for time $s - t_0 n^{1/3}$ with initial data $\bm{\mathsf{y}} (t_0 n^{1/3}) = \bm{u}$. Condition on $\bm{u}$ and $\bm{\mathsf{y}} (n^4)$, and define the $n$-tuples $\bm{v}', \bm{v}'' \in \mathbb{W}_n$ by setting $v_j' = \mathsf{y}_j (n^4)$ and $v_j'' = v_j' + (n-j)n$ for each $j \in \llbracket 1, n \rrbracket$. Then sample $n$ non-intersecting Brownian bridges $\bm{\mathsf{z}} = (\mathsf{z}_1, \mathsf{z}_2, \ldots , \mathsf{z}_n) \in \mathcal{C} \big( [t_0 n^{1/3}, n^4] \big)$ from the measure $\mathsf{Q}^{\bm{u}; \bm{v}''}$. See \Cref{f:DBM_coupling}. 
		
		We will first use gap and height monotonicity to compare the gaps of $\bm{\mathsf{x}}$ and $\bm{\mathsf{y}}$, through $\bm{\mathsf{z}}$. To do this, first observe that the law of $\bm{\mathsf{x}}$ is given by $\mathsf{Q}_{\mathsf{x}_{n+1}}^{\bm{u}; \bm{v}}$. For any $j \in \llbracket 1, n-1 \rrbracket$, we have $v_j - v_{j+1} = x_j (n^4) - x_{j+1} (n^4) \le n \le v_j'' - v_{j+1}''$, where the first statement holds by the definition of $\bm{v}$; the second by the fact that we have restricted to the event $\mathscr{E} \subseteq \textbf{GAP}_n (n^4; n)$ (and the fact that $n \big( (j+1)^{2/3} - j^{2/3} \big) + (\log n)^{25} \le n$ for sufficiently large $n$); and the third by the definition of $\bm{v}''$. Hence, it follows from gap monotonicity \Cref{monotonedifference} that we may couple $\bm{\mathsf{x}}$ and $\bm{\mathsf{z}}$ such that 
		\begin{flalign}
			\label{xjtzjtn4} 
			\mathsf{x}_j (t) - \mathsf{x}_{j+1} (t) \le \mathsf{z}_j (t) - \mathsf{z}_{j+1} (t), \qquad \text{for $(j, t) \in \llbracket 1, n-1 \rrbracket \times [t_0 n^{1/3}, n^4]$}. 
		\end{flalign} 
		
		\noindent Moreover observe from the second part of \Cref{lambdat} that the law of $\bm{\mathsf{y}}$ is given by $\mathsf{Q}^{\bm{u}; \bm{v}'}$. Since that of $\bm{\mathsf{z}}$ is given by $\mathsf{Q}^{\bm{u}; \bm{v}''}$ and since $0 \le v_j'' - v_j' \le n^2$ for each $j \in \llbracket 1, n \rrbracket$, it follows from the second part of \Cref{uvv} that we may couple $\bm{\mathsf{y}}$ and $\bm{\mathsf{z}}$ in such a way that 
		\begin{flalign*} 
			\big| \mathsf{y}_j (t) - \mathsf{z}_j (t) \big| \le \displaystyle\frac{b n^{1/3}}{n^4 - t_0 n^{1/3}} \cdot n^2 \le 2bn^{-5/3}, \qquad \text{for $(j, t) \in \llbracket 1, n \rrbracket \times \big[t_0 n^{1/3}, (t_0 + b) n^{1/3} \big]$},
		\end{flalign*} 
		
		\noindent where we used the fact that $n^4 - t_0 n^{1/3} \ge n^4 / 2$ for $n \ge 2$ (as $|t_0| \le n$) for sufficiently large $n$. Combining this bound with \eqref{xjtzjtn4} yields
		\begin{flalign}
			\label{xjtyjtn4} 
			\mathsf{x}_j (t) - \mathsf{x}_k (t) \le \mathsf{y}_j (t) - \mathsf{y}_k (t) + 4b n^{-5/3}, \qquad \text{for $1 \le j < k \le n$ and $t \in \big[t_0 n^{1/3}, (t_0 + b) n^{1/3} \big]$}.
		\end{flalign}	
		
		Next, by \Cref{pacxi2}, for any real number $s \in [an^{1/3}, 2an^{1/3}]$ the $n$-tuple $\bm{\mathsf{y}} (t_0 n^{1/3} + s)$ is $\big( [0, b-sn^{-1/3}]; \omega; \xi/2 \big)$-packed with probability at least $1 - 2 \xi^{-1} e^{-\xi (\log n)^2 / 2}$; in particular, it is $\big( [0, b-2a]; \omega; \xi/2 \big)$-packed with at least the same probability $1 - 2 \xi^{-1} e^{-\xi (\log n)^2 / 2}$. This, \eqref{xjtyjtn4}, \Cref{gap2}, and the fact that $16 b n^{-2/3} \le \omega n^{-1/3}$ for sufficiently large $n$ together imply that the first bound in \eqref{xjtb} holds. Moreover, the fact that $\bm{\mathsf{y}} (t_0 n^{1/3} + s)$ is $\big( [0, b-2a]; \omega; \xi / 2 \big)$-packed; \Cref{gappacked}; \eqref{xjtyjtn4}; and the fact that $16 b n^{-2/3} \le \omega n^{-1/3}$ together imply the second bound in \eqref{xjtb}. This establishes \Cref{kgapc}.
	\end{proof}

	\chapter{Limit Shapes Near the Edge}
	
	\label{EDGESHAPE}

	\section{Limit Shapes for Non-intersecting Brownian Bridges}
	\label{LimitBridges}
	
	In this section we collect some results on limit shapes for families of non-intersecting Brownian bridges, without upper and lower boundaries. We begin by introducing them and their properties in \Cref{Limit0}; we then provide examples in \Cref{ProcessExample}, and continuous variants of monotonicity for them in \Cref{MonotoneContinuous}. In \Cref{EquationHG} and \Cref{EstimatesEquation0} we recall an elliptic partial differential equation and regularity estimates satisfied by these limit shapes. In \Cref{ConcentrationSmooth0} we recall a concentration bound for non-intersecting Brownian bridges, indicating conditions under which they are closely approximated by their limit shapes.
	
	 The results from this section are primarily due to \cite{U,LDASI,FOAMI,LDSI}. Most of the statements we use from \cite{U} are mild modifications of those proven earlier in \cite{LDASI,FOAMI,LDSI}, sometimes together with (conventional) arguments from elliptic partial differential equations. The main exception is the concentration bound \Cref{gh} (used in \Cref{APPROXIMATECURVE} below), constituting \cite[Sections 4-7]{U}. 
	
	\subsection{Limit Shapes} 
	
	\label{Limit0} 
	
	In this section we recall results from \cite{LDASI,FOAMI,LDSI} concerning the limiting macroscopic behavior of non-intersecting Brownian bridges under given starting and ending data, with no upper and lower boundaries. Throughout, we use coordinates $(t, x)$ or sometimes $(t, y)$ for $\mathbb{R}^2$ (instead of $(x, y)$).

	 Fix an interval $I \subseteq \mathbb{R}$. A \emph{measure-valued process} (on the time interval $I$) is a family $\bm{\mu} = (\mu_t)_{t \in I}$ of measures $\mu_t \in \mathscr{P}_{\fin}$ for each $t \in I$. Given a real number $A > 0$, we say that $\bm{\mu}$ has \emph{constant total mass $A$} if $\mu_t (\mathbb{R}) = A$, for each $t \in I$. If $\bm{\mu}$ has constant total mass $1$ (so each $\mu_t \in \mathscr{P}$), we call $\bm{\mu}$ a \emph{probability measure-valued process}. Measure-valued processes can be interpreted as elements of $I \times \mathscr{P}_{\fin}$ and probability measure-valued processes as ones of $I \times \mathscr{P}$. We let $\mathcal{C} ( I; \mathscr{P}_{\fin})$ and $\mathcal{C} ( I; \mathscr{P})$\index{C@$\mathcal{C} (I; \mathscr{P}_{\fin})$, $\mathcal{C} (I; \mathscr{P})$} denote the sets of measure-valued processes and probability measure-valued processes that are continuous in $t \in I$, under the topology of weak convergence on $\mathscr{P}_{\fin}$ and $\mathscr{P}$, respectively. 
	
	Given two measures $\mu, \nu \in \mathscr{P}_{\fin}$ of finite total mass, the L\'{e}vy distance\index{D@$d_{\dL}$; L\'{e}vy distance} between them is 
	\begin{flalign}
		\label{munudistance1}
		d_{\dL} (\mu, \nu) = \inf \Bigg\{ a > 0 : \displaystyle\int_{-\infty}^{y-a} \mu (dx) - a \le \displaystyle\int_{-\infty}^y \nu (dx) \le \displaystyle\int_{-\infty}^{y+a} \mu (dx) + a, \quad \text{for all $y \in \mathbb{R}$} \Bigg\}.
	\end{flalign}
	
	\noindent Given an interval $I \subseteq \mathbb{R}$ and two measure-valued processes $\bm{\mu} = (\mu_t)_{t \in I} \in I \times \mathscr{P}_{\fin}$ and $\bm{\nu} = (\nu_t)_{t \in I} \in I \times \mathscr{P}_{\fin}$ on the time interval $I$, the L\'{e}vy distance between them is defined to be 
	\begin{flalign}
		\label{munut} 
		d_{\dL} (\bm{\mu}, \bm{\nu}) = \sup_{t \in I} d_{\dL} (\mu_t, \nu_t).
	\end{flalign}

	The following lemma from \cite{LDSI} (based on results from \cite{FOAMI,LDASI}) states that, as $n$ tends to $\infty$, the empirical measure (recall \eqref{aemp}) for $n$ non-intersecting Brownian bridges (whose starting and ending data converge in a certain way) has a limit. The following lemma was stated in \cite{LDSI} in the case when $[a, b] = [0, 1]$ and $A = 1$ but, by the scaling invariance (\Cref{scale}) for non-intersecting Brownian bridges, it also holds for any interval $[a, b]$ and real number $A > 0$, as below.  In what follows, we recall the notation $\emp$ from \eqref{aemp}.

	\begin{lem}[{\cite[Claim 2.13]{LDSI}}]
		
		\label{rhot} 
		
		Fix real numbers $a < b$ and compactly supported measures $\mu_a, \mu_b \in \mathscr{P}_{\fin}$, both of total mass $\mu_a (\mathbb{R}) = A = \mu_b (\mathbb{R})$ for some real number $A > 0$. There is a measure-valued process $\bm{\mu}^{\star} = (\mu_t^{\star})_{t \in [a, b]} \in \mathcal{C} \big( [a, b]; \mathscr{P}_{\fin} \big)$ on $[a, b]$ of constant total mass $A$, which is continuous in the pair $(\mu_a, \mu_b) \in \mathscr{P}_{\fin}^2$ under the L\'{e}vy metric, such that the following holds. For each integer $n \ge 1$, let $\bm{u} = \bm{u}^n \in \overline{\mathbb{W}}_n$ and $\bm{v} = \bm{v}^n \in \overline{\mathbb{W}}_n$ denote sequences such that $A \cdot \emp(\bm{u}^n)$ and $A \cdot \emp(\bm{v}^n)$ converge to $\mu_a$ and $\mu_b$ under the L\'{e}vy metric as $n$ tends to $\infty$, respectively. Sample $n$ non-intersecting Brownian bridges $\bm{x}^n = (x_1^n, x_2^n, \ldots , x_n^n) \in \llbracket 1, n \rrbracket \times \mathcal{C} \big( [a, b] \big)$ from $\mathsf{Q}^{\bm{u}; \bm{v}} (An^{-1})$; for any $t \in [a, b]$, denote $\nu_t^n = A \cdot \emp \big(\bm{x}^n (t) \big) \in \mathscr{P}$. For any real number $\varepsilon > 0$, we have $\lim_{n \rightarrow \infty} \mathbb{P} \big[ d_{\dL} (\bm{\nu}^n, \bm{\mu}^{\star}) > \varepsilon \big] = 0$. 
		
	\end{lem}

	Terminology for the limit shape provided by \Cref{rhot} is given by the following definition.

	\begin{figure}
\center
\includegraphics[width=0.9\textwidth]{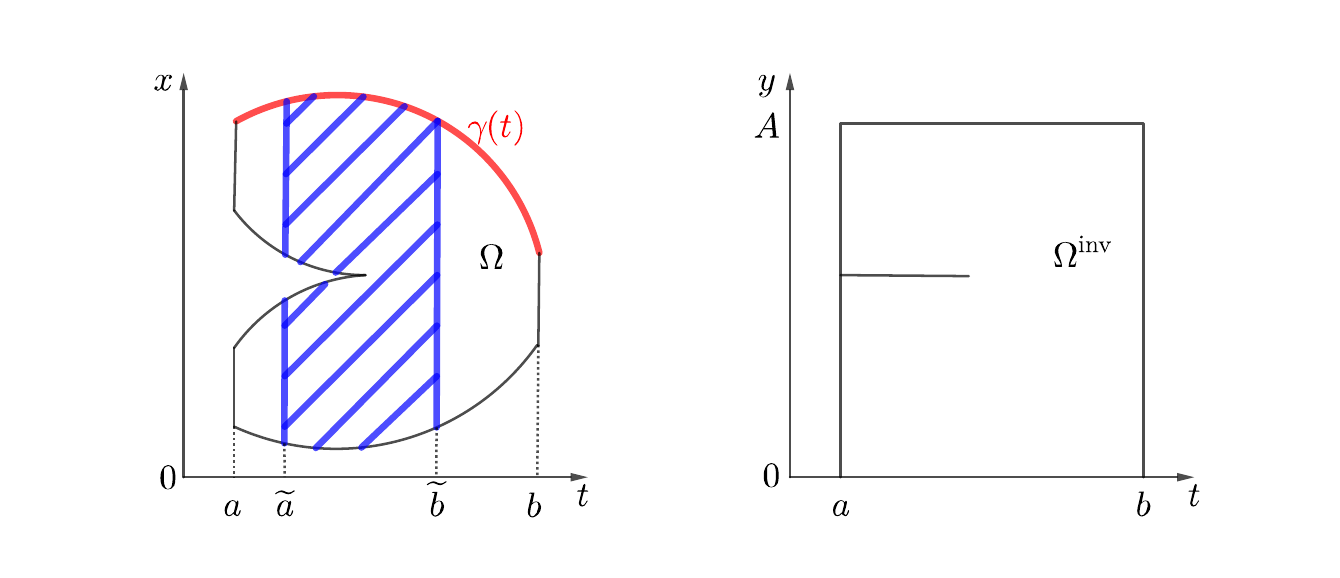}
\caption{Shown to the left is the liquid region $\Omega$ associated with a bridge-limiting measure on $[a,b]$; restricting the latter to the (striped blue) shorter interval $[\widetilde a, \widetilde b]$ again gives a bridge-limiting measure. The (red) curve $\gamma (t)$ traces the north boundary of $\Omega$, called the arctic boundary. Shown to the right is the associated inverted liquid region $\Omega^{\inv}$. As shown there, $\Omega^{\inv}$ does not contain a certain line segment, as $\varrho_t (x) = 0$ at the corresponding $(x, t)$ coordinates. }
\label{f:liquidregion}
\end{figure}

	\begin{definition}
		
		\label{mutmu0mu1}
		
		Adopting the notation of \Cref{rhot}, the measure-valued process $\bm{\mu}^{\star} = (\mu_t^{\star})_{t \in [a, b]}$ is called the \emph{bridge-limiting measure process}\index{B@bridge-limiting measure process} (on the interval $[a, b]$) with boundary data $(\mu_a; \mu_b)$.
		
	\end{definition}

	The following lemma from \cite[Lemma 3.3]{U} indicates how bridge-limiting measure processes restrict to others; see the left side of \Cref{f:liquidregion} for a depiction.
	
	\begin{lem}[{\cite[Lemma 3.3]{U}}]
		
		\label{gtabab}
		
		Adopt the notation and assumptions of \Cref{rhot}, and let $\widetilde{a}, \widetilde{b} \in \mathbb{R}$ be real numbers such that $a \le \widetilde{a} < \widetilde{b} \le b$. Then, the bridge-limiting measure process on the interval $[\widetilde{a}, \widetilde{b}]$ with boundary data $(\mu_{\widetilde{a}}^{\star}; \mu_{\tilde{b}}^{\star})$ is given by $(\mu_t^{\star})_{t \in [\tilde{a}, \tilde{b}]}$.
	\end{lem}

	We will often make use of a height function and inverted height functions associated with a measure-valued process, defined as follows.
	
	\begin{definition}
		
		\label{hrhot} 
		
		Fix an interval $I = [a, b] \subseteq \mathbb{R}$ and a measure-valued process $\bm{\mu} = (\mu_t)_{t \in I}$ of constant total mass $A > 0$. The \emph{height function associated with $\bm{\mu}$}\index{H@height function} is defined to be the function $H = H^{\bm{\mu}} : I \times \mathbb{R} \rightarrow \mathbb{R}$ obtained by setting
		\begin{flalign}
			\label{htxintegral}
			H (t, x) = \displaystyle\int_x^{\infty} \mu_t (dw), \qquad \text{for each $(t, x) \in I \times \mathbb{R}$}.
		\end{flalign}
		
		\noindent The \emph{inverted height function}\index{I@inverted height function} associated with $\bm{\mu}$ is $G = G^{\bm{\mu}}: I \times [0, A] \rightarrow \mathbb{R}$, defined by setting $G(t,0) = \inf \big\{ x : H(t,x) = 0 \big\}$ and 
		\begin{flalign}
			\label{gty}
			G(t, y)  = \sup \Bigg\{ x \in \mathbb{R} : H(t, x) = \displaystyle\int_x^{\infty} \mu_t (dw) \ge y \Bigg\}, \qquad \text{for each $y \in (0, A]$}.
		\end{flalign}
		
		\noindent Thus, in analogy with \eqref{gammaj}, we may view $G (t, y)$ as a classical location of the measure $\mu_t \in \mathscr{P}$.
		
		If $\mu_t = \varrho_t (x) dx$ has a density with respect to Lebesgue measure for each $t \in (a, b)$, then we sometimes associate $H$ (or $G$) with $\varrho = (\varrho_t)$. Moreover, if $\bm{\mu}$ is the bridge-limiting measure process with boundary data $(\mu_a; \mu_b)$ we say that $H$ (or $G$) is associated with boundary data $(\mu_a; \mu_b)$.
	\end{definition} 
	
	\begin{rem} 
		
		Let us briefly translate the continuum notation from \Cref{hrhot} (for limit shapes) to the discrete context (for non-intersecting Brownian bridges). The height function $H(t,x)$ is the continuum analog of the index for a Brownian bridge, around (rescaled) location $x$ at (rescaled) time $t$. The inverted height function $H(t,y)$ is the continuum analog of the location for a Brownian bridge, with (rescaled) index $y$ at (rescaled) time $t$. The (negative inverse of the) density $-\varrho_t (x)^{-1}$ is the continuum analog of the gap between a pair of consecutive Brownian bridges and (rescaled) location $x$ at (rescaled) time $t$.

	\end{rem}

	The following lemma essentially due to \cite{FOAMI} (but appearing as stated below in \cite{U}) indicates that the measures $\mu_t^{\star}$ have a density, and it also discusses properties of this density. In what follows, we recall the free convolution and semicircle distribution $\mu_{\semci}^{(t)}$ from \Cref{TransformConvolution}.

	\begin{lem}[{\cite[Lemma 3.7 and Remark 3.14]{U}}] 
		
		\label{mutrhot}

		Adopting the notation and assumptions of \Cref{rhot}, the following statements hold for each real number $t \in (a, b)$. 
		
		\begin{enumerate} 
			\item There  exists a measurable function $\varrho_t^{\star} : \mathbb{R} \rightarrow \mathbb{R}_{\ge 0}$ such that $\mu_t^{\star} (dx) = \varrho_t^{\star} (x) dx$. 
			\item There exists some compactly supported measure $\nu_t \in \mathscr{P}_{\fin}$ of total mass $\nu_t (\mathbb{R}) = A$, dependent on $\mu_a$ and $\mu_b$ with $\supp \nu_t \subseteq \supp \mu_a + \supp \mu_b$, such that $\varrho_t^{\star}=\nu_t\boxplus \mu_{\semci}^{((t-a)(b-t)/(b-a))}$.
			\item We have $\varrho_t^{\star} (x)^2 \le  A (b-a) \big( (t-a)(b-t) \big)^{-1}$, for any $x \in \mathbb{R}$. 
			\item The function $\varrho_t (x)$ is continuous on $(a, b) \times \mathbb{R}$.
		\end{enumerate}
		
	\end{lem}
	
	The following definition provides notation for the region on which the density $\varrho_t$ is positive (both in terms of the $(t, x)$ coordinates of the height function and the $(t, y)$ coordinates of the inverted height function). See \Cref{f:liquidregion} for a depiction.
	
	\begin{definition} 
		
	\label{omegafg} 
	
	Fix an interval $(a, b) \subseteq \mathbb{R}$ and a family of measures $\bm{\mu} = (\mu_t)_{t \in (a, b)} \in [a, b] \times \mathscr{P}_{\fin}$ of constant total mass $A > 0$. Assume for each $t \in (a, b)$ that each $\mu_t$ has a density $\varrho_t$ with respect to Lebesgue measure,  for some continuous function $\varrho_t : \mathbb{R} \rightarrow \mathbb{R}_{\ge 0}$ that is also continuous in $t$. Recalling the associated height and inverted height functions $H = H^{\bm{\mu}} : [a, b] \times \mathbb{R} \rightarrow \mathbb{R}$ and $G = G^{\bm{\mu}} : [a, b] \times [0, A] \rightarrow \mathbb{R}$ from \Cref{hrhot}, we define the associated \emph{liquid region} $\Omega \subset (a, b) \times \mathbb{R}$ and \emph{inverted liquid region} $ \Omega^{\inv} \subseteq (a, b) \times (0, A)$ by\index{0@$\Omega, \Omega^{\inv}$; (inverted) liquid region}
	\begin{flalign}
		\begin{split}
			\label{omega12}
			&\Omega= \big\{(t, x)\in (a, b) \times \mathbb{R}: \varrho_t (x)>0 \big\}; \\
			& \Omega^{\inv} = \big\{ (t, y) \in (a, b) \times (0, A) : y = H^{\bm{\mu}} (t,x), (t, x) \in \Omega \big\}.
		\end{split}
	\end{flalign} 
	
	\end{definition} 

	\noindent Observe that the map $(t, x) \mapsto \big( t, H(t, x) \big)$ is a bijection to from $\Omega$ to $\Omega^{\inv}$. Moreover, by the continuity of $\varrho$ (which in our context will be verified by \Cref{mutrhot}), the set $\Omega$ is open, which implies that $\Omega^{\inv}$ is also open.

	We next state two lemmas essentially due to \cite{LDSI,LDASI,FOAMI} (but stated as below in \cite{U}). The first reformulates \Cref{rhot} through (inverted) height functions; there, we recall the height function $\mathsf{H}^{\bm{x}}$ associated with a line ensemble $\bm{x}$ from \Cref{htw}. The second indicates that the height and inverted height functions $H^{\star}$ and $G^{\star}$ are smooth on $\Omega$ and $\Omega^{\inv}$, respectively.

	\begin{lem}[{\cite[Corollary 3.6 and Corollary 3.8]{U}}]
		
		\label{convergepathomega} 
		
		Adopt the notation and assumptions of \Cref{rhot}, and fix a real number $\varepsilon > 0$. Let $G^{\star} : [a, b] \times [0, A] \rightarrow \mathbb{R}$ denote the inverted height function associated with $\bm{\mu}^{\star}$, respectively; further denote the associated inverted liquid region by $\Omega^{\inv} \subseteq (a, b) \times (0, A)$. Then the following two statements hold.
		
		\begin{enumerate} 
			\item For any $y \in (0, 1)$, we have
		\begin{flalign*}
			& \displaystyle\lim_{n \rightarrow \infty} \mathbb{P} \Bigg[ \big\{ G^{\star} (t, Ay + \varepsilon) - \varepsilon \le x^n_{\lfloor yn \rfloor} (t) \le G^{\star} (t, Ay - \varepsilon) + \varepsilon \big\} \Bigg] = 1,
		\end{flalign*} 
	
			\item For any $y \in (0, 1)$ such that $(t, Ay) \in \Omega^{\inv}$ holds for each $t \in (a, b)$, we have that $G^{\star} (t, Ay)$ is continuous in $t \in [a, b]$, and
		\begin{flalign*} 
			& \displaystyle\lim_{n \rightarrow \infty} \mathbb{P} \Bigg[ \bigcap_{t \in (a, b)} \big\{ G^{\star} (t, Ay) - \varepsilon \le x_{\lfloor yn \rfloor}^n (t) \le G^{\star} (t, Ay) + \varepsilon \big\} \Bigg] = 1.
		\end{flalign*}
		\end{enumerate}
	\end{lem}

	\begin{lem}[{\cite[Lemma 3.23(1)]{U}}]
		
		\label{urhoderivatives0}
		
		Fix real numbers $a < b$ and $A > 0$, and compactly supported measures $\mu_a, \mu_b \in \mathscr{P}_{\fin}$, satisfying $\mu_a (\mathbb{R}) = A = \mu_b (\mathbb{R})$. Let $\bm{\mu}^{\star} = (\mu_t^{\star})_{t \in [a, b]} \in \mathcal{C} \big( [a, b]; \mathscr{P}_{\fin} \big)$ denote the bridge-limiting measure process on $[a, b]$ with boundary data $(\mu_a; \mu_b)$. Further let $H^{\star} : [a, b] \times \mathbb{R} \rightarrow \mathbb{R}$ and $G^{\star} : [a, b] \times [0, A] \rightarrow \mathbb{R}$ denote the associated height and inverted height functions, respectively. Then, $H^{\star} (t, x)$ is smooth for $(t, x) \in \Omega$ and $G^{\star} (t, y)$ is smooth for $(t, y) \in \Omega^{\inv}$. 
		
	\end{lem} 
	
	We next define a complex slope associated with a limit shape; its imaginary part is given by the associated density $\varrho^{\star}$, and its real part is given by the $t$-derivative of the inverted height function, which we denote by $u^{\star} : \Omega \rightarrow \mathbb{R}$.     
	
	\begin{definition} 
		
		\label{urho}
		
		Adopt the notation and assumptions of \Cref{urhoderivatives0}. Define the function $u^{\star} : \Omega \rightarrow \mathbb{R}$ as follows. For any point $(t, x) \in \Omega$, let $(t, y) \in \Omega^{\inv}$ be the unique point such that $G^{\star} (t, y) = x$ (which is guaranteed to exist since the map $(t, y) \mapsto (t,G^{\star} (t, y))$ is a bijection from $\Omega^{\inv}$ to $\Omega$). Then, define $u^{\star} (t, x) = u_t^{\star} (x)$ by setting
		\begin{flalign}
			\label{u0} 
			u_t^{\star} (x) = \partial_t G^{\star} (t, y),
		\end{flalign}
	
		\noindent and observe that 
		\begin{flalign}
			\label{urho20}
		\varrho_t^{\star} ( x)=	\varrho_t^{\star} \big( G(t, y) \big) = -\partial_x H^{\star} (t, x) = - \big(\partial_y G^{\star} (t, y) \big)^{-1},
		\end{flalign}
	 
	 	\noindent where the last equality holds by \Cref{urhoderivatives0} and the fact from \eqref{gty} that $H^{\star} \big( t, G^{\star} (t, y) \big) = y$ (see also \cite[Remark 3.15]{U}). Further define associated \emph{complex slope} $f = f^{(\mu_a; \mu_b)} : \Omega \rightarrow \overline{\mathbb{H}}$\index{C@complex slope} by setting
	\begin{flalign}
		\label{frhou} 
		f(t, x) = u_t^{\star} (x) + \pi \mathrm{i} \varrho_t^{\star} (x), \qquad \text{for each $(t, x) \in \Omega$}.
	\end{flalign}

	\end{definition}

	The following lemma from \cite{LDSI} (implicitly due to the earlier work \cite{FOAMI}; see also \cite[Remark 3.15 and Lemma 3.23(2)]{U}) indicates that this function $f$ satisfies a complex variant of the Burgers equation.\index{C@complex Burgers equation}
	
	\begin{lem}[{\cite[Lemma 3.23(2)]{LDSI}}]
		
		\label{p:solution}	
		
		Adopting the notation and assumptions of \Cref{urhoderivatives0}, the associated complex slope $f$ satisfies the complex Burgers equation, 
		\begin{equation}
			\label{ftfx}
			\partial_t f(x,t)+f(x,t) \cdot \partial_x f(x,t)=0, \qquad \text{for all $(t, x) \in \Omega$}.
		\end{equation} 
		
	\end{lem}

	\subsection{Examples of Bridge-Limiting Measure Processes} 
	
	\label{ProcessExample} 
	
\begin{figure}
\centering
\begin{subfigure}{.5\textwidth}
  \centering
  \includegraphics[width=1\linewidth]{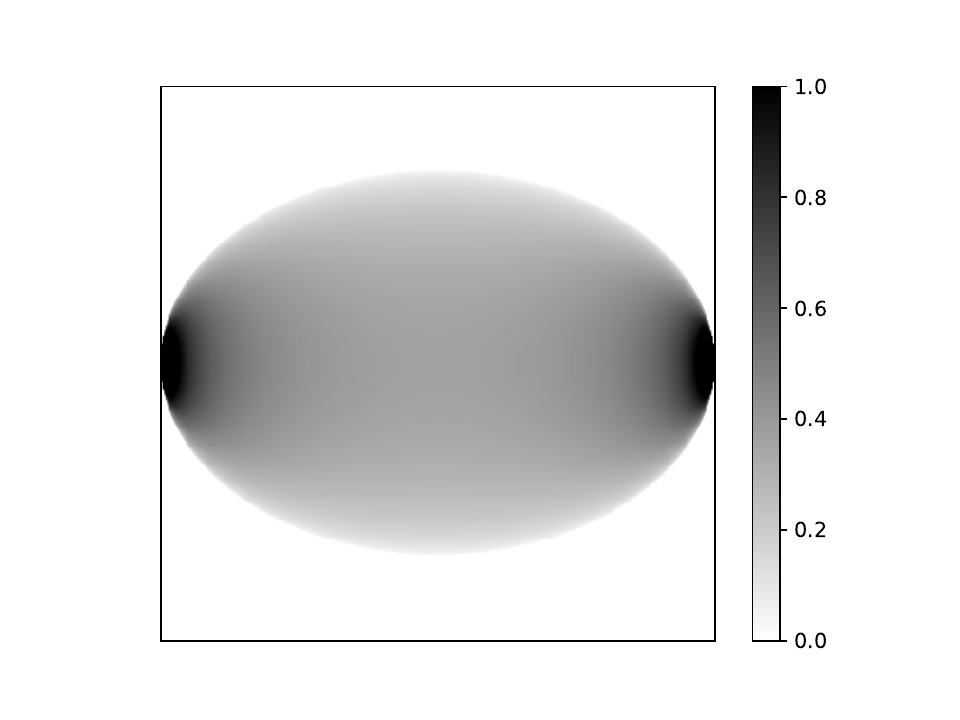}
\end{subfigure}%
\begin{subfigure}{.5\textwidth}
  \centering
  \includegraphics[width=1\linewidth]{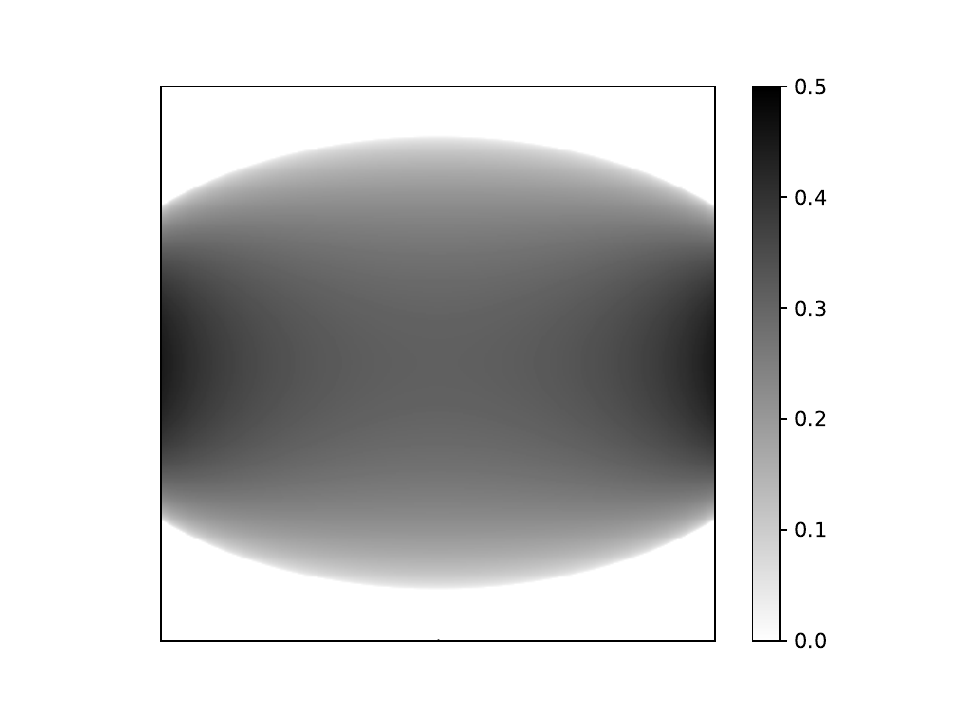}
\end{subfigure}
\caption{Shown on the left is a depiction for \Cref{xuv0} at $(a,b,u,v,A)=(0,6,0,0,2)$.  Shown on the right is a depiction for \Cref{rhoasct} at $(a,b,d, A)=(0,10,2, 2)$. In both, the entire shaded region is the associated liquid region $\Omega$.}
\label{f:density}
\end{figure}

	In this section we describe several examples of the bridge-limiting measure processes from \Cref{Limit0}. The first concerns the case when $\mu_0$ and $\mu_1$ are delta measures, in which the associated non-intersecting Brownian bridges form a Brownian watermelon (recall \Cref{PathsUV0}); see the left side of \Cref{f:density} for a depiction.

	\begin{example}
		
		\label{xuv0} 
		
		Fix real numbers $a < b$; $u, v \in \mathbb{R}$; and $A > 0$. Assume that $(\mu_a, \mu_b) = (A\cdot \delta_u, A \cdot \delta_v)$, where $\delta_x \in \mathscr{P}_0$ denotes the delta measure at $x \in \mathbb{R}$. Then, it follows from \Cref{estimatexj} (multiplying its results by $(An^{-1})^{1/2}$ to account for the fact that the Brownian motions have variance $An^{-1}$ here) and the second statement of \Cref{convergepathomega} (and also the continuity of $\gamma_{\semci} (y)$ below in $y$) that the inverted height function $G^{\star} : [0, 1] \times [0, A] \rightarrow \mathbb{R}$ associated with boundary data $(\mu_0; \mu_1)$ is given by
		\begin{flalign*}
			G^{\star} (t, y) = \bigg( \displaystyle\frac{A (b-t) (t-a)}{b-a} \bigg)^{1/2} \cdot \gamma_{\semci} \Big( \displaystyle\frac{y}{A} \Big) + \displaystyle\frac{b-t}{b-a} \cdot u + \displaystyle\frac{t-a}{b-a} \cdot v, 
		\end{flalign*} 
		
		\noindent  where $\gamma_{\semci} (y)$ is the classical location of the semicircle distribution given by \eqref{gammaj}.  Together with \eqref{htxintegral} and \eqref{gty}, it follows that the associated density process $(\varrho_t^{\star})$ and the height function $H^{\star} : [0, 1] \times \mathbb{R}$ are given by
		\begin{flalign*}
			& \varrho_t^{\star} (x) = A \cdot \varrho_{\semci}^{(\frac{A (b-t)(t-a)}{(b-a)})}  \bigg(x - \displaystyle\frac{b-t}{b-a} \cdot u - \displaystyle\frac{t-a}{b-a} \cdot v \bigg),
		\end{flalign*} 
		
		\noindent and $H^{\star} (t, x) = \int_x^{\infty} \varrho_t^{\star} (y) dy$, where we recall the rescaled semcircle density $\varrho_{\semci}^{(t)}$ from \eqref{rhosct}.
		
	\end{example}

	The second example from \cite{U} concerns the case when $\mu_a$ and $\mu_b$ are rescaled semicircle distributions (recall \eqref{rhosct}), which can be obtained by restricting a watermelon to a smaller time interval; see the right side of \Cref{f:density} for a depiction.
	
	\begin{example}[{\cite[Corollary 3.10]{U}}]
		
		\label{rhoasct}
		
		Fix real numbers $a < b$ and $d, A > 0$; assume that $\mu_a = A \cdot \mu_{\semci}^{(d)} = \mu_b$. Then, the inverted height function $G^{\star} : [a, b] \times [0,A] \rightarrow \mathbb{R}$ and density process $(\varrho_t^{\star})$ associated with boundary data $(\mu_a; \mu_b)$ are given by 
		\begin{flalign}
			\label{rhoa} 
			\varrho_t^{\star} (x) = A \cdot \varrho_{\semci}^{(d + \frac{A(b-t)(t-a)}{b-a+2\kappa})} (x); \qquad G^{\star} (t, y) = \bigg( d + \displaystyle\frac{A(b-t)(t-a)}{b-a + 2\kappa} \bigg)^{1/2} \cdot \gamma_{\semci} \Big( \displaystyle\frac{y}{A} \Big),
		\end{flalign} 
		
		\noindent where $\kappa = \kappa (a, b, d) > 0$ is defined by
		\begin{flalign}
			\label{kappa} 
			\kappa = \displaystyle\frac{d}{A} + \displaystyle\frac{a-b}{2} + \bigg( \Big( \displaystyle\frac{b-a}{2} \Big)^2 + \Big( \displaystyle\frac{d}{A} \Big)^2 \bigg)^{1/2}.
		\end{flalign}
	\end{example}

			\begin{rem}
				
				\label{airy}

				Let us consider the limiting profile associated with affine shifts of the scaled parabolic Airy line ensemble $\bm{\mathcal{S}}$. Fix real numbers $\mathfrak{a}, \mathfrak{b}, \mathfrak{c}$ with $\mathfrak{c} > 0$, and set $\sigma = 2^{1/6} \mathfrak{c}^{1/3}$. For any integer $n \ge 1$, define the affine shift $\bm{\mathcal{S}}^{(\sigma; \mathfrak{a}, \mathfrak{b}; n)} = \big( \mathcal{S}_1^{(\sigma; \mathfrak{a}, \mathfrak{b}; n)}, \mathcal{S}_2^{(\sigma; \mathfrak{a}, \mathfrak{b}; n)}, \ldots \big) \in \mathbb{Z}_{\ge 1} \times \mathcal{C} (\mathbb{R})$ of the rescaled parabolic Airy line ensemble $\bm{\mathcal{S}}^{(\sigma)}$ (recall from \eqref{sigmar}), by for each $(j, t) \in \mathbb{Z}_{\ge 1} \times \mathbb{R}$ setting 
				\begin{flalign*}
					\mathcal{S}_j^{(\sigma; a, b; n)} (t) = \mathcal{S}_1^{(\sigma)} (t) + \mathfrak{a} n^{2/3} + \mathfrak{b} n^{1/3} t.
				\end{flalign*} 
			
				\noindent Observe from \Cref{sigmascale} and \Cref{linear} that $\bm{\mathcal{S}}^{(\sigma; \mathfrak{a}, \mathfrak{b}; n)}$ satisfies the Brownian Gibbs property. Define the \emph{limiting Airy profile} to be the function\index{G@$\mathfrak{G}_{\Ai}$; limiting Airy profile} $\mathfrak{G}_{\Ai} = \mathfrak{G}_{\Ai; \mathfrak{a}, \mathfrak{b}, \mathfrak{c}} : \mathbb{R} \times \mathbb{R}_{\ge 0} \rightarrow \mathbb{R}$ by setting 
				\begin{align}\label{e:Airyprofile}
						\mathfrak{G}_{\Ai} (t, y) =\fa+\fb t- \fc t^2- \bigg( \displaystyle\frac{3 \pi}{4 \mathfrak{c}^{1/2}} \bigg)^{2/3} y^{2/3},
					\end{align}

				\noindent for each $(t, y) \in \mathbb{R} \times \mathbb{R}_{\ge 0}$. By \Cref{kdeltad}, a union bound, and the definition  $\sigma = 2^{1/6} \mathfrak{c}^{1/3}$, we have for any real numbers $a < b$ and $\varepsilon > 0$ that
				\begin{flalign*}
					\displaystyle\lim_{n \rightarrow \infty} \mathbb{P} \Bigg[ \displaystyle\sup_{t \in [a, b]} \displaystyle\sup_{y \in [0, 1]} \Big| n^{-2/3} \cdot \mathcal{S}_{\lfloor yn \rfloor}^{(\sigma; \mathfrak{a}, \mathfrak{b}; n)} (tn^{1/3})  - \mathfrak{G}_{\Ai} (t, y)  \Big| \le \varepsilon \Bigg] = 1.
				\end{flalign*}
			
				\noindent Define the process $\bm{\mu}_{\Ai} = \bm{\mu}_{\Ai; \mathfrak{a}, \mathfrak{b}, \mathfrak{c}} = (\mu_t) = (\mu_t^{\Ai; \mathfrak{a}, \mathfrak{b}, \mathfrak{c}})$ (over $t \in \mathbb{R})$; the density process $(\varrho_t) = (\varrho_t^{\Ai; \mathfrak{a}, \mathfrak{b}, \mathfrak{c}})$; and the function $\mathfrak{H}_{\Ai} = \mathfrak{H}_{\Ai; \mathfrak{a}, \mathfrak{b}, \mathfrak{c}} : \mathbb{R}^2 \rightarrow \mathbb{R}$ by setting
				\begin{flalign}
					\label{grhotx} 
					\varrho_t (x) = \displaystyle\frac{2 \mathfrak{c}^{1/2}}{\pi} (\mathfrak{a} + \mathfrak{b} t - \mathfrak{c} t^2 - x)^{1/2} \cdot \textbf{1}_{x \le \mathfrak{a} + \mathfrak{b} t - \mathfrak{c} t^2}; \quad \mu_t (dx) = \varrho_t (x) dx; \quad \mathfrak{H}_{\Ai} (t, x) = \displaystyle\int_x^{\infty} \varrho_t (w) dw.
				\end{flalign}
			
				\noindent By \eqref{htxintegral} and \eqref{gty}, it is quickly confirmed that $\mathfrak{H}_{\Ai}$ and $\mathfrak{G}_{\Ai}$ are the height and inverted height functions associated with the process $\bm{\mu}_{\Ai}$ (as in \Cref{hrhot}, whose notions are also well-defined if $\bm{\mu}$ has infinite mass). In this way, one can view $\mathfrak{H}_{\Ai}$ and $\mathfrak{G}_{\Ai}$ as the large scale limits of the rescaled parabolic Airy line ensemble. Since each $\mu_t$ has infinite mass, it is not a bridge-limiting measure process in the sense of \Cref{mutmu0mu1}, but we will see that it will satisfies many of the same properties as one. 
				
			\end{rem}

	\subsection{Continuous Variants of Monotonicity}
	
	\label{MonotoneContinuous}
	
	In this section we discuss continuous variants of both height monotonicity (\Cref{monotoneheight}) and gap monotonicity (\Cref{monotonedifference}), which apply to bridge-limiting measure processes (\Cref{mutmu0mu1}). They are given by the first and second lemmas below, respectively. The first statements of these lemmas assume some type of (either inverted height or gap) comparison between two families of boundary data, along their entire west and east boundaries, and deduce that the comparison continues to hold in the interior of the domain. The second statements assume a comparison between these boundary data, but only along the parts of their west and east boundaries that lie above a given ``level line.'' It then shows the comparison continues to hold in the interior of the domain above this level line, if one further assumes a certain comparison between the level lines of the two processes (one lies above the other for height monotonicity, and one is ``more concave'' than the other for gap monotonicity, parallel to \Cref{monotoneheight} and \Cref{monotonedifference}, respectively). The proofs of these two lemmas, which are quick consequences of the discrete variants of monotonicity (\Cref{monotoneheight} and \Cref{monotonedifference}), with the convergence of non-intersecting Brownian bridges to their limit shapes (\Cref{convergepathomega}), are provided in \Cref{ProofContinuousCompare} below. In what follows, we recall the inverted height function and inverted liquid region associated with a measure-valued process from \Cref{hrhot} and \eqref{omega12}, respectively.

	\begin{figure}
\center
\includegraphics[scale = .7, trim =0 1cm 0 0.5cm, clip]{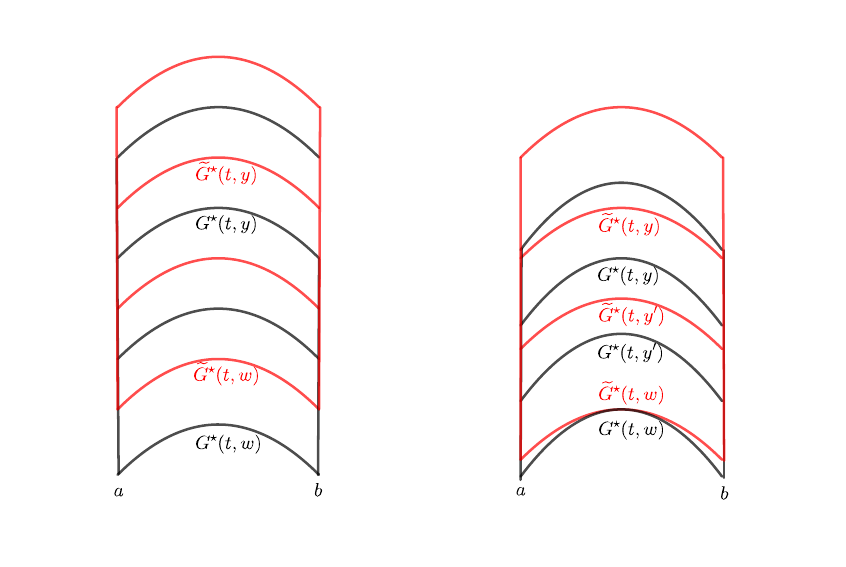}
\caption{Shown to the left is a depiction for the continuous variant of height monotonicity; shown to the right is a depiction for the continuous variant of gap monotonicity. 
 }
\label{f:continuousgap}
\end{figure}

	\begin{lem}
		
		\label{limitheightcompare}
		
		Fix real numbers $a < b$ and $A, \widetilde{A} > 0$; set $A_0 = \min \{ A, \widetilde{A} \}$; and fix compactly supported measures $\mu_a, \widetilde{\mu}_a, \mu_b, \widetilde{\mu}_b \in \mathscr{P}_{\fin}$ such that $\mu_a (\mathbb{R}) = A = \mu_b (\mathbb{R})$ and $\widetilde{\mu}_a (\mathbb{R}) = \widetilde{A} = \widetilde{\mu}_b (\mathbb{R})$. Let $\bm{\mu}^{\star}$ and $\widetilde{\bm{\mu}}^{\star}$ denote the bridge-limiting measure processes on $[a, b]$ with boundary data $(\mu_a; \mu_b)$ and $(\widetilde{\mu}_a; \widetilde{\mu}_b)$, respectively. Also denote the associated inverted height functions by $G^{\star} : [a, b] \times [0, A] \rightarrow \mathbb{R}$ and $\widetilde{G}^{\star} : [a, b] \times [0, \widetilde{A}] \rightarrow \mathbb{R}$ and the associated inverted liquid regions by $\Omega^{\inv}$ and $\widetilde{\Omega}^{\inv}$, respectively.

		\begin{enumerate}
			\item Assume $\widetilde A\geq A$ and for each $(t, y) \in \{ a, b \} \times [0, A_0]$ that $G^{\star} (t, y) \le \widetilde{G}^{\star} (t, y)$. Then, $G^{\star} (t, y) \le \widetilde{G}^{\star} (t, y)$ holds for each $(t, y) \in [a, b] \times [0, A_0]$.

			\item Fix a real number $w \in (0, A_0)$ such that $(t, w) \in \Omega^{\inv} \cap \widetilde{\Omega}^{\inv}$ for each $t \in (a, b)$. Assume for each $(t, y) \in \{ a, b \} \times [0, w]$ that $G^{\star} (t, y) \le  \widetilde{G}^{\star} (t, y)$, and further assume for each $t \in [a,b]$ that $G^{\star} (t, w) \le \widetilde{G}^{\star} (t, w)$. Then, $G^{\star} (t, y) \le \widetilde{G}^{\star} (t, y)$ holds for all $(t, y) \in [a, b] \times [0, w]$. 
			
		\end{enumerate}
	\end{lem}

	\begin{lem}
		
		\label{limitdifferencecompare}
		
		Adopt the notation and assumptions of \Cref{limitheightcompare}.

		\begin{enumerate}
			\item Assume that $A\geq \widetilde A$ and  $G^{\star} (t, y) - G^{\star} (t, y') \le \widetilde{G}^{\star} (t, y) - \widetilde{G}^{\star} (t, y')$ holds for each $t \in \{ a, b \}$ and $y, y' \in [0, A_0]$ with $y < y'$. Then, $G^{\star} (t, y) - G^{\star} (t, y') \le \widetilde{G}^{\star} (t, y) - \widetilde{G}^{\star} (t, y')$ holds for each $t \in [a, b]$ and $y, y' \in [0, A_0]$ with $y < y'$.

			\item Fix a real number $w \in (0, A_0)$ so that $(t, w) \in \Omega^{\inv} \cap \widetilde{\Omega}^{\inv}$ for each $t \in (a, b)$. Assume that 
			\begin{flalign*} 
					& \big| G^{\star} (t, y) - G^{\star} (t, y') \big| \le \big| \widetilde{G}^{\star} (t, y) - \widetilde{G}^{\star} (t, y') \big|,
				\end{flalign*} 
			
				\noindent for all $(t, y), (t, y') \in \{ a, b \} \times [0, w]$. Further assume that $\widetilde{G}^{\star} - G^{\star}$ is convex, that is, for all real numbers $t_1, t_2 \in [a, b]$ and $r \in [0, 1]$ we have
					\begin{flalign} 
						\label{gdifference121}
						\begin{aligned}
					r \cdot G^{\star} (t_1, w & ) - G^{\star} \big( rt_1+ (1-r) t_2, w \big) + (1-r) \cdot G^{\star} (t_2, w) \\
					& \le r \cdot \widetilde{G}^{\star} (t_1, w) - \widetilde{G}^{\star} \big( rt_1 + (1-r) t_2, w \big) + (1-r) \cdot \widetilde{G}^{\star} (t_2, w).
					\end{aligned}
					\end{flalign}
				
					\noindent Then, $\big| G^{\star} (t, y) - G^{\star} (t, y') \big| \le \big| \widetilde{G}^{\star} (t, y) - \widetilde{G}^{\star} (t, y') \big|$  holds for all $(t, y), (t, y') \in [a, b] \times [0, w]$.

		\end{enumerate}

	\end{lem}

	See the left and right sides of \Cref{f:continuousgap} for depictions of \Cref{limitheightcompare} and \Cref{limitdifferencecompare}, respectively.

	While the limiting Airy profiles of \Cref{airy} are not quite bridge-limiting measure processes in the sense of \Cref{mutmu0mu1} (as they have infinite mass), the following analog of \Cref{limitheightcompare} provides a height comparison between bridge-limiting measure processes and limiting Airy profiles. Its proof is similar to that of \Cref{limitheightcompare} (using the concentration bound \Cref{kdeltad} for the rescaled parabolic Airy line ensemble in place of \Cref{convergepathomega}); we also provide it in \Cref{ProofContinuousCompare} below.

	\begin{lem}
		
		\label{airyheightcompare}
		
		Fix real numbers $a < b$ and $A> 0$; and fix measures $\mu_a, \mu_b$ such that $\mu_a (\mathbb{R}) = A = \mu_b (\mathbb{R})$. Let $\bm{\mu}^{\star}$ denote the bridge-limiting measure processes on $[a, b]$ with boundary data $(\mu_a; \mu_b)$. Denote the associated inverted height function by $G^{\star} : [a, b] \times [0, A] \rightarrow \mathbb{R}$ and the associated inverted liquid region by $\Omega^{\inv}$. Let $\widetilde G^{\star}: [a,b]\times [0,\infty]\rightarrow \bR$ be a limiting Airy profile of the form \eqref{e:Airyprofile}.

		\begin{enumerate}
			\item Assume for each $(t, y) \in \{ a, b \} \times [0, A]$ that $G^{\star} (t, y) \leq \widetilde{G}^{\star} (t, y)$. Then, $G^{\star} (t, y) \leq \widetilde{G}^{\star} (t, y)$ holds for each $(t, y) \in [a, b] \times [0, A]$. 
			\item Fix a real number $w \in (0, A)$ such that $(t, w) \in \Omega^{\inv} $ for each $t \in (a, b)$. 
			\begin{enumerate} 
				\item  Assume for each $(t, y) \in \{ a, b \} \times [0, w]$ that $G^{\star} (t, y) \le  \widetilde{G}^{\star} (t, y)$, and for each $t \in [a,b]$ that $G^{\star} (t, w) \le \widetilde{G}^{\star} (t, w)$. Then, $G^{\star} (t, y) \le \widetilde{G}^{\star} (t, y)$ holds for all $(t, y) \in [a, b] \times [0, w]$. 
				\item Assume for each $(t, y) \in \{ a, b \} \times [0, w]$ that $G^{\star} (t, y) \geq \widetilde{G}^{\star} (t, y)$, and for each $t \in [a,b]$ that $G^{\star} (t, w) \geq \widetilde{G}^{\star} (t, w)$. Then, $G^{\star} (t, y) \geq \widetilde{G}^{\star} (t, y)$ holds for all $(t, y) \in [a, b] \times [0, w]$. 
			\end{enumerate} 
		\end{enumerate}
	\end{lem}
	
	It is also possible to state and prove a variant of \Cref{airyheightcompare} that compares the gaps between limiting Airy profiles and those of inverted height function associated with bridge-limiting measure processes. However, we will not pursue this here, since we will not need it.

		\subsection{Elliptic Partial Differential Equations for the Height Function}
	
	\label{EquationHG}

	In this section we state an elliptic partial differential equation satisfied by the inverted height function associated with a bridge-limiting measure process, and related results. The former is provided through the following lemma, shown as stated below in \cite{U} (though implicitly due to the earlier works \cite{LDASI,FOAMI}).

	\begin{lem}[{\cite[Lemma 3.23(3)]{U}}] 
		
	\label{gequation}

	Adopting the notation and assumptions of \Cref{urhoderivatives0}, we have
	\begin{flalign}
		\label{equationxtd}
		\partial_t^2 G^{\star} (t, y) + \pi^2 \big( \partial_y G^{\star} (t, y) \big)^{-4} \cdot \partial_y^2 G^{\star} (t, y) = 0, \qquad \text{for each $(t, y) \in \Omega^{\inv}$}. 
	\end{flalign} 	
	
	\end{lem}

	It will be useful to make use of invariances of the equation \eqref{equationxtd} under the following (linear and multiplicative) transformations. We  first require an additional definition; below, we recall that a locally Lipschitz functions is differentiable almost everywhere (with respect to Lebesgue measure) on its domain, by Rademacher's theorem.
	
	\begin{definition} 
		
		\label{functionsadmissible0} 
		
		For any bounded, open subset $\mathfrak{R} \subset \mathbb{R}^2$, we let $\Adm (\mathfrak{R})$\index{A@$\Adm (\mathfrak{R})$, $\Adm_{\varepsilon} (\mathfrak{R})$} denote the set of locally Lipschitz functions $F \in \mathcal{C}(\overline{\mathfrak{R}})$ such that $\partial_y F (t, y) < 0$, for almost all $(t, y) \in \mathfrak{R}$ (with respect to Lebesgue measure); we call such functions \emph{admissible}. 
	
	\end{definition}

	\begin{lem}[{\cite[Lemma 3.21]{U}}]
		
	\label{invariancesscale} 
	
	Fix a bounded, open subset $\mathfrak{R} \subset \mathbb{R}^2$ and a function $G \in \Adm (\mathfrak{R}) \cap \mathcal{C}^2 (\mathfrak{R})$; assume on $\mathfrak{R}$ that $G$ satisfies \eqref{equationxtd}. Fix nonzero real numbers $\alpha$ and $\beta$, and denote $\widetilde{\mathfrak{R}} = \widetilde{\mathfrak{R}}_{\alpha; \beta} = \big\{ (t, y) \in \mathbb{R}^2  : (\alpha t, \beta y) \in \mathfrak{R} \big\}$. 
	\begin{enumerate}
		
		\item \label{scale21} Assuming $\alpha > 0$ and $\beta > 0$, define $\widetilde{G} \in \mathcal{C}^2 (\widetilde{\mathfrak{R}})$ by $\widetilde{G} (t, y) = (\alpha \beta)^{-1/2} G(\alpha t, \beta y)$. Then $\widetilde{G}$ satisfies \eqref{equationxtd} on $\widetilde{\mathfrak{R}}$.
		
		\item Define $\widehat{G} \in \mathcal{C}^2 (\mathfrak{R})$ by $\widehat{G} (t, y) = G(t, y)+\alpha t$. Then, $\widehat{G}$ satisfies \eqref{equationxtd} on $\mathfrak{R}$.

	\end{enumerate}
	
	\end{lem}
		
	The $\alpha = \beta$ case of \Cref{scale21} in \Cref{invariancesscale} would have held for a solution $G$ to the equation $\sum_{i, j \in \{ t, y \}} \mathfrak{a}_{ij} (\nabla G) \cdot \partial_i \partial_j G = 0$, for any measurable coefficients $\mathfrak{a}_{ij}$. However, that this remains true for all $(\alpha, \beta)$ is special to the specific choice of these coefficients appearing in \eqref{equationxtd}. This more general scaling invariance will be useful in analyzing solutions to \eqref{equationxtd} in \Cref{EDGESHAPE} (see, for example, the proof of \Cref{p:densityrg}).

	\subsection{Regularity Estimates}
	
	\label{EstimatesEquation0}

	In this section we recall from \cite{U} various estimates for  solutions to the partial differential equation \eqref{equationxtd}; throughout, we recall the norms defined in \eqref{e:norms}. We first require the following definition. 
	
	\begin{definition} 
		
		\label{functionadmissible} 
		
		For any real number $\varepsilon \in (0, 1) $ and bounded, open subset $\mathfrak{R} \subset \mathbb{R}^2$, we let $\Adm_{\varepsilon} (\mathfrak{R}) \subset \Adm (\mathfrak{R})$ denote the set of functions $F \in \Adm (\mathfrak{R})$ such that $\varepsilon < -\partial_y F(t, y) < \varepsilon^{-1}$ for almost all $(t, y) \in \mathfrak{R}$ (with respect to Lebesgue measure).
		
	\end{definition}

	Next, we state the maximum principle for solutions of \eqref{equationxtd}.
	
	\begin{lem}[{\cite[Lemma 9.1]{U}}]
		
		\label{maximumboundary}
		
		Fix some open set $\mathfrak{R} \subset \mathbb{R}$, and let $F_1, F_2, F \in \Adm (\mathfrak{R}) \cap \mathcal{C}^2 (\mathfrak{R})$ denote solutions to \eqref{equationxtd} on $\mathfrak{R}$. 
		\begin{enumerate} 
			\item If $F_1 (z) \le F_2 (z)$ for each $z \in \partial \mathfrak{R}$, then $F_1 (z) \le F_2 (z)$ for each $z \in \mathfrak{R}$. 
			\item We have $\sup_{z \in \mathfrak{R}} \big| F_1 (z) - F_2 (z) \big| \le \sup_{z \in \partial \mathfrak{R}} \big| F_1 (z) - F_2 (z) \big|$. In particular, $\sup_{z \in \mathfrak{R}} \big| F (z) \big| = \sup_{z \in \partial \mathfrak{R}} \big| F (z) \big|$.
		\end{enumerate} 
	\end{lem}

		\begin{figure}
	\center
\includegraphics[width=1\textwidth, trim=0 1cm 0 0.5cm, clip]{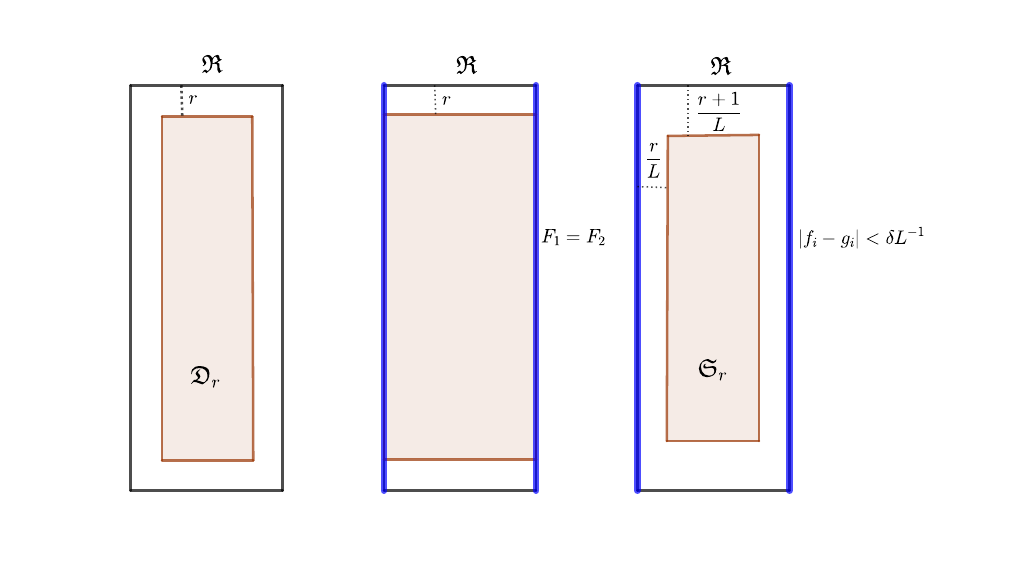}

\caption{Shown to the left is a depiction of $\mathfrak{D}_r$ in \Cref{derivativef}. Shown in the middle is a depiction for \Cref{equationcompareboundary} stating, if $F_1(z)=F_2(z)$ for each $z$ on the blue sides of $\mathfrak{R}$, then $F_1-F_2$ is exponentially small in the shaded region. Shown to the right is a depiction for \Cref{f1f2b}.}
\label{f:setting_1}
	\end{figure}

	We next have the following lemma indicating boundedness of the interior derivatives for a solution to \eqref{equationxtd}; see the left side of \Cref{f:setting_1}. 
	
	\begin{lem}[{\cite[Lemma 9.2]{U}}]
		
		\label{derivativef}
		
		For any integer $m \ge 1$, and real numbers $r > 0$; $\varepsilon \in (0, 1)$; and $B > 1$, there exists a constant $C = C(\varepsilon, r, B, m) > 1$ such that the following holds. Let $\mathfrak{R} \subset \mathbb{R}^2$ be a bounded open set, let $f \in \mathcal{C}(\partial \mathfrak{R})$ be a function satisfying $\| f\|_0 \le B$, and let $F \in \Adm_{\varepsilon} (\mathfrak{R}) \cap \mathcal{C}^2 (\mathfrak{R})$ be a solution to \eqref{equationxtd} on $\mathfrak{R}$ such that $F |_{\partial \mathfrak{R}} = f$. Letting $\mathfrak{D}_r = \big\{ z \in \mathfrak{R} : \dist (z, \partial \mathfrak{R}) > r \big\}$, we have $\| F \|_{\mathcal{C}^m (\overline{\mathfrak{D}}_r)} \le C$.
	\end{lem}

		The following lemma states that the solutions of \eqref{equationxtd} are real analytic.

	\begin{lem}[{\cite[Lemma 9.3]{U}}]
		
		\label{realanalytic}
		
		Fix a real number $\varepsilon\in (0,1)$, some open set $\mathfrak{R} \subset \mathbb{R}^2$, and let $F \in \Adm_{\varepsilon} (\mathfrak{R}) \cap \mathcal{C}^2 (\mathfrak{R})$ denote a solution to \eqref{equationxtd} on $\mathfrak{R}$. Then, $F$ is real analytic on $\mathfrak{R}$. 
	\end{lem}

	The following result states that, given two solutions $F_1, F_2$ to \eqref{equationxtd} on a tall rectangle of aspect ratio $2L$, whose boundary data match on its west and east boundaries, $|F_1 - F_2|$ decays exponentially in $L$ in the middle of the rectangle; see the middle of \Cref{f:setting_1}. It will not be used until \Cref{APPROXIMATECURVE}.
	
	\begin{lem}[{\cite[Proposition 9.5]{U}}]
		
		\label{equationcompareboundary}
		
		For any real numbers $\varepsilon, r \in ( 0, 1 / 4)$ and $B > 1$, there exists a constant $c = c(\varepsilon, r, B) > 0$ such that the following holds. Fix a real number $L > 0$, and define the open rectangle $\mathfrak{R} = (0, L^{-1}) \times (-1, 1)$. Let $F_1, F_2 \in \Adm_{\varepsilon} (\mathfrak{R}) \cap \mathcal{C}^5 (\overline{\mathfrak{R}})$ be two solutions to \eqref{equationxtd} on $\mathfrak{R}$ such that $\| F_i \|_{\mathcal{C}^5 (\mathfrak{R})} \le B$ for each $i \in \{ 1, 2 \}$. Assume that $F_1 (t, x) = F_2 (t, x)$ for any $(t, x) \in \partial \mathfrak{R}$ with $t \in \{ 0, L^{-1} \}$. Then, 
		\begin{align}
			\label{equationcompareboundary2} 
			\big| F_1 (t, x) - F_2 (t, x) \big| \le c^{-1} e^{-cL^{1/8}}, \qquad \text{for any $(t, x) \in [0, L^{-1}] \times [r-1, 1-r]$}.
		\end{align}
	\end{lem}

	\begin{rem} 
		
		The exponential decay in \eqref{equationcompareboundary2} (where the exponent $1/8$ there is not optimal) is typical in the context of uniformly elliptic partial differential equations. It admits the following explanation in the case of the Laplace equation (for this discussion, we will set $B=1$ in \Cref{equationcompareboundary}, as we may by scaling). Let $\mathfrak{R} = (0, L^{-1}) \times (-1, 1)$ be as in \Cref{equationcompareboundary}, but now for each $i \in \{ 1, 2 \}$ let $F_i : \overline{\mathfrak{R}} \rightarrow \mathbb{R}$ denote a solution to $\partial_t^2 F_i (t,x) + \partial_x^2 F_i ( t, x) = 0$. Assume again that $F_1 (t, x) = F_2 (t, x)$ whenever $(t, x) \in \{ 0, L^{-1} \} \times [-1, 1 ]$ is on the west or east edge of $\partial \mathfrak{R}$, and that $\big| F_i (t, x) \big| \le 1$ whenever $(t, x) \in [0, L^{-1}] \times \{ -1, 1 \}$ is on the north or south edge of $\partial \mathfrak{R}$. Then $F_i (t,x) = \mathbb{E} \big[ F_i (W_{\tau}) \big]$, where $W_s : \mathbb{R} \rightarrow \mathbb{R}^2$ denotes a two-dimensional Brownian motion initially at $W_0 = (t, x)$, and $\tau = \inf \{ s : B_s \in \partial \mathfrak{R} \}$ denotes the hitting time where $W_{\tau}$ first meets $\partial \mathfrak{R}$. 
		
		Since $F_1 (W_{\tau}) = F_2 (W_{\tau})$ if $W_{\tau} \in \{ 0, L^{-1} \} \times [-1, 1]$ is on the west or east boundary of $\partial \mathfrak{R}$, and $\big| F_i (W_{\tau}) \big| \le 1$ otherwise, we have that $\big| F_1 (t, x) - F_2 (t, x) \big| = \mathbb{E} \big[ |F_1 (W_{\tau}) - F_2 (W_{\tau})| \big] \le p$, where $p = \mathbb{P} \big[ W_{\tau} \in [0, L^{-1}] \times \{ -1, 1 \} \big]$ is the probability that $W_{\tau}$ is on the north or south edge of $\mathfrak{R}$. If $(t, x) \in [0, L^{-1}] \times [r-1, 1-r]$ for some fixed $r > 0$, then the distance from $W_0 = (t,x)$ to the north and south edges of $\mathfrak{R}$ is bounded away from $0$, while its distance from the west and east edges of $\mathfrak{R}$ is $L^{-1}$. Using this, it can be shown that this probability decays exponentially in $L$, so that for some constant $c = c(r) > 0$ we have $\big| F_1 (t, x) - F_2 (t, x) \big| \le p \le c^{-1} e^{-cL}$. This is analogous to \eqref{equationcompareboundary2} and, indeed, the proof of \Cref{equationcompareboundary} in \cite{U} is based on an adaptation of the above reasoning to the nonlinear equation \eqref{equationxtd}.

	\end{rem}

	We conclude this section with the next lemma (that will also not be used until \Cref{APPROXIMATECURVE}), which states the following. Fix a solution $F$ to \eqref{equationxtd}, bounded in $\mathcal{C}^m$ for some integer $m$, on a rectangle $\mathfrak{R}$, as well some boundary data $g_0$ and $g_1$ on the two vertical sides on the rectangle that are close to $F$. Then, it is possible to find a solution $G$ to \eqref{equationxtd} on a slightly shorter rectangle $\mathfrak{S}$, whose boundary data on the vertical sides of the rectangle are given by $g_0$ and $g_1$ (the first condition in the lemma), and that is close to $F$ (quantified through the second and third conditions of the lemma). The second part of the lemma states that $F$ and $G$ are close in any $\mathcal{C}^k$ norm in the interior of $\mathfrak{S}$, and the third part states that $F$ and $G$ are close in $\mathcal{C}^{m-5}$ (that is, fewer derivatives than the original assumed bound on $F$) up to the boundary of $\mathfrak{S}$. See the right side of \Cref{f:setting_1} for a depiction.
	
	\begin{lem}[{\cite[Lemma 9.6]{U}}]
		
		\label{f1f2b}

	For any integers $m, k \ge 7$, and real numbers $\varepsilon > 0$; $r \in ( 0, 1 / 4)$; and $B > 1$, there exist constants $\delta = \delta(\varepsilon, B) > 0$, $C_1 = C_1 (\varepsilon, r, B, k) > 1$, and $C_2 = C_2 (\varepsilon, B, m) > 1$ such that the following holds. Fix a real number $L > 2$, and define the open rectangles 
	\begin{flalign*} 
		\mathfrak{R} = \Big( 0, \displaystyle\frac{1}{L} \Big) \times (-1, 1); \qquad \mathfrak{S}_r = \Big( \displaystyle\frac{r}{L}, \displaystyle\frac{1-r}{L} \Big) \times \Big( \displaystyle\frac{r+1}{L} - 1, 1 - \displaystyle\frac{r+1}{L} \Big); \qquad \mathfrak{S} = \mathfrak{S}_0. 
	\end{flalign*} 
	
	\noindent Let $F \in \Adm_{\varepsilon} (\mathfrak{R}) \cap \mathcal{C}^m (\overline{\mathfrak{R}})$ denote a solution to \eqref{equationxtd} on $\mathfrak{R}$ such that $\| F \|_{\mathcal{C}^m (\mathfrak{R})} \le B$, and define the functions $f_0, f_1 : [-1, 1] \rightarrow \mathbb{R}$ by setting $f_i (x) = F (iL^{-1}, x)$ for each $(i, x) \in \{ 0, 1 \} \times [-1, 1]$. Further let $g_0, g_1 : [-1, 1] \rightarrow \mathbb{R}$ denote two functions such that $\| g_i \|_{\mathcal{C}^m (-1, 1)} \le B$ and $\big| g_i (x) - f_i (x) \big| \le \delta L^{-1}$ for each $(i, x) \in \{ 0, 1\} \times [-1, 1]$. Then, there exists a solution $G \in \Adm_{\varepsilon/2} (\mathfrak{S}) \cap \mathcal{C}^{m-5} (\overline{\mathfrak{S}})$ to \eqref{equationxtd} on $\mathfrak{S}$ satisfying the following three properties.
	
	\begin{enumerate} 
		\item For each $i \in \{ 0, 1 \}$ and $x \in \big[ L^{-1} - 1, 1 - L^{-1} \big]$, we have $G(iL^{-1}, x) = g_i (x)$. 
		\item We have $\| F - G \|_{\mathcal{C}^k (\mathfrak{S}_r)} \le C_1 L^k \cdot \big( \| f_0 - g_0 \|_{\mathcal{C}^0} + \| f_1 - g_1 \|_{\mathcal{C}^0}  \big)$. 
		\item We have  $\| F - G \|_{\mathcal{C}^{m-5} (\mathfrak{S})} \le C_2 L^{m-5} \cdot \big( \| f_0 - g_0 \|_{\mathcal{C}^0}^{3/m} + \| f_1 - g_1 \|_{\mathcal{C}^0}^{3/m} \big)$.
	\end{enumerate}

	\end{lem}

	\subsection{Concentration Estimates for Non-Intersecting Brownian Bridges} 
	
	\label{ConcentrationSmooth0}
	
	In this section we state a concentration bound from \cite{U} (stronger than \Cref{concentrationbridge} but requiring more stringent hypotheses) for families of non-intersecting Brownian bridges sampled from the measure $\mathsf{Q}_{f;g}^{\bm u; \bm v}$ of \Cref{qxyfg}; it will be used for the proof of \Cref{h0x2} in \Cref{APPROXIMATECURVE}. We begin by specifying the regularity assumption to which we will subject our boundary data (namely, the starting and ending data, $\bm{u}$ and $\bm{v}$, and the lower and upper boundaries, $f$ and $g$, of the paths).

	\begin{assumption}
		
		\label{fgr}
		
		Fix a real number $B > 1$. Let $n \ge 1$ be an integer and $L > 0$ be a real number; define the rectangle $\mathfrak{R} = (0, L^{-1}) \times (0, 1) \subset \mathbb{R}^2$, and let $G \in \mathcal{C}^{50} (\overline{\mathfrak{R}})$ be a function. Assume that 
		\begin{flalign}
			\label{glg} 
			L \in (B^{-1}, n^{1/20000}); \qquad G \in \Adm_{1/10} (\mathfrak{R}); \qquad  \big\| G - G(0,0) \big\|_{\mathcal{C}^{50} (\mathfrak{R})} \le B, 
		\end{flalign} 
		
		\noindent and that $G$ solves \eqref{equationxtd} on $\mathfrak{R}$. Define $f, g: [0, L^{-1}] \rightarrow \mathbb{R}$ by setting $f(s) = G(s, 1)$ and $g(s) = G(s, 0)$, for each $s \in [0, L^{-1}]$. Further let $\bm{u}, \bm{v} \in \overline{\mathbb{W}}_n$ be $n$-tuples, and assume for each $j \in \llbracket 1, n \rrbracket$ that   
		\begin{flalign}
			\label{uvg}
			 u_j = G(0,jn^{-1}); \qquad v_j = G (L^{-1}, jn^{-1}).
		\end{flalign}
		
		\noindent Sample non-intersecting Brownian bridges $\bm{x} = (x_1, x_2, \ldots , x_n) \in \llbracket 1, n \rrbracket \times \mathcal{C} \big( [0,L^{-1}] \big)$ from the measure $\mathsf{Q}_{f;g}^{\bm{u}; \bm{v}} (n^{-1})$. 
	\end{assumption}
	
	Let us briefly explain \Cref{fgr}. The function $G$ will eventually be the limit shape for the family $\bm{x}$ of non-intersecting Brownian bridges, in the sense that we will have $x_j (t) \approx G (t, jn^{-1})$. The conditions that $f(s) = G(s, 0)$ and $g(s) = G(s,1)$ ensure that this holds for the upper and lower boundaries of the model (formally, when $j \in \{ 0, n + 1 \}$), and \eqref{uvg} ensures that this holds when $t \in \{ 0, 1 \}$. The constraint that $G \in \Adm_{1/10} (\mathfrak{R})$ ensures that $G$ has no ``frozen facets'' (macroscopic regions containing no curves), and the constraint that $\| G \|_{\mathcal{C}^{50} (\mathfrak{R})} \le B$ ensures that $G$ has some regularity.

	Under \Cref{fgr}, the below concentration bound holds, stating that $x_j (t) \approx G (t, jn^{-1}) + \mathcal{O} (n^{-23/24})$. It follows from the $(\varkappa, \varepsilon, \delta, m) = (0, 1/10, 1/20000, 49)$ case of \cite[Theorem 1.4]{U}. 
	
	\begin{lem}[{\cite[Theorem 1.4]{U}}] 
		
		\label{gh} 
		
		Adopt \Cref{fgr}. There is a constant $c = c (B) > 0$ with   
		\begin{flalign*} 
			\mathbb{P} \Bigg[ \displaystyle\sup_{s \in [0, L^{-1}]} \bigg( \displaystyle\max_{j \in \llbracket 1, n \rrbracket} \big| x_j (s) - G(s, jn^{-1}) \big| \bigg) > c^{-1} n^{-23/24} \Bigg] \le c^{-1} e^{-c(\log n)^2}.
		\end{flalign*} 
		
	\end{lem}
	
	Let us briefly comment on the constants appearing in \Cref{fgr} and \Cref{gh}. In the assumption $G \in \Adm_{1/10} (\mathfrak{R})$ of \Cref{fgr}, the constant $1/10$ can be replaced by any arbitrarily small, but fixed, positive real number $\varepsilon > 0$. Moreover, the error $n^{-23/24}$ in \Cref{gh} can be replaced by $n^{2/m-1+\delta}$ for any integer $m \ge 4$ and $\delta \in (0, 1/5m^2)$, if we assume that $\big\| G - G(0,0) \big\|_{\mathcal{C}^{m+1} (\mathfrak{R})} \le B$ and $L \le n^{\delta}$; in this way, the error can be made arbitrarily close to $n^{-1}$ (in the exponent), by increasing the regularity of the boundary data. Lemma \ref{gh} comes from taking $(\varepsilon, \delta, m) = (1/10, 20000, 49)$, which will suffice for our purposes.

	\section{Density Estimates for Limit Shapes}
\label{s:density}

\subsection{Free Convolution Estimates} 

In this section we collect some estimates on free convolution measures, which are subject to the following assumption that bounds their integrals. In what follows, we recall notation on free convolutions from \Cref{TransformConvolution}.

\begin{assumption} 
	
\label{blx12} 

Let $B, L \ge 1$ be real numbers, and let $\tau\in[B^{-1}, B]$ be a real number; let $\nu \in \mathscr{P}_{\fin}$ be a measure. Assume that its total mass is $\nu (\mathbb{R}) = L^{3/2}$ and that $\supp \nu \subseteq [-BL,0]$.  We denote the measure $\nu_{\tau} = \nu\boxplus \mu_{\semci}^{(\tau)}$, which is the free convolution of $\nu$ with the rescaled semicircle distribution. As mentioned in \Cref{TransformConvolution}, $\nu_{\tau}$ admits a density with respect to Lebesgue measure, which we denote by $\varrho_{\tau} \in L^1 (\mathbb{R})$. 
\end{assumption}

\begin{assumption}\label{blx122} 
Adopting  \Cref{blx12}, further assume that $\varrho_{\tau}$ satisfies
\begin{align}\label{e:rhotintbound}
	\int_{x}^\infty \varrho_\tau(y) dy \leq B|x|^{3/2},\qquad \text{for each $x \in [-BL, -1]$}.
\end{align}

\end{assumption} 

We then have the following two propositions. The former, established in \Cref{Proofrhox} below, bounds the support of $\varrho_{\tau}$ under \Cref{blx12} and bounds its magnitude under \Cref{blx122}. The latter, established in \Cref{ProofrhoxDerivative} below, bounds the derivatives of $\varrho_{\tau}$ under \Cref{blx122}, assuming a lower bound on $\varrho_{\tau}$ (made precise through the function $\gamma_{\tau}$ in \eqref{e:defgammatau} below).

\begin{prop}\label{p:densityest}

For any real number $B \ge 1$, there exists a constant $C = C(B) > 1$ such that the following holds. 
\begin{enumerate} 
	\item Adopting \Cref{blx12}, we have $\supp \varrho_\tau \subseteq [-CL, CL^{3/4}]$. 
	\item If we further adopt \Cref{blx122} then
\begin{align}\label{e:densitybb0}
\varrho_\tau(x)\leq C\max\{1, -x\}^{3/4},\qquad \text{for each $x\in \bR$}.
\end{align}
\end{enumerate}
\end{prop}

\begin{rem}
	
	Let us briefly (and slightly imprecisely) comment on the intuition behind \eqref{e:densitybb0}. For simplicity, we suppose that $B=1$ in \Cref{blx122}; that $x \le -1$; and that \eqref{e:rhotintbound} also holds for the initial data $\nu$. Then, $\nu \big( [x, \infty) \big) \le |x|^{3/2}$, and so gap monotonicity \Cref{limitdifferencecompare} (or its discrete form, \Cref{gapmotion}, for Dyson Brownian motion) suggests that ``compressing'' the part $\nu |_{[x, \infty)}$ of $\nu$ above $x$ to the delta mass $|x|^{3/2} \cdot \delta_x$ should increase the density $\varrho_{\tau}$. If we were able to neglect the effect the part $\nu |_{(-\infty, x)}$ of $\nu$ below $x$, then this would give rise to the free convolution on $|x|^{3/2} \cdot \delta_x$, run for time $\tau$. By \Cref{mtscale}, the latter is the rescaled semicircle law $|x|^{3/2} \cdot \mu_{\sc}^{(\tau |x|^{3/2})}$ from \eqref{rhosct}, which for $\tau \sim 1$ has density around $|x|^{3/4}$ in its bulk, by \eqref{rho1}; this matches \eqref{e:densitybb0}. However, gap monotonicity implies that neglecting $\nu |_{[-\infty, x)}$ in fact decreases $\varrho_{\tau}$, so this cannot be done directly. We therefore instead prove \Cref{p:densityest} using analytic properties of the free convolution.
	
\end{rem}

	\begin{prop}\label{c:rhoderbound}
		
		For any integer $\ell \ge 1$ and real numbers $A \ge 1$ and $B \ge 3$, there exist constants $\varepsilon=\varepsilon(A, B) > 0$ and $C = C(\ell, A, B) > 1$ such that the following holds. Adopt \Cref{blx12} and \Cref{blx122}. Define the function $\gamma_\tau:[0, L^{3/2}]\mapsto \bR$ by setting
		\begin{align}\label{e:defgammatau}
			\gamma_\tau(y)=\sup \left\{ x \in \mathbb{R} : \int_{x}^\infty \varrho_\tau(u)d u\geq y\right\}, 
		\end{align}
		
		\noindent for each $y \in [0, L^{3/2}]$. We further assume the following two bounds.   
		\begin{enumerate}
			\item We have $\gamma_\tau(B)\geq -A$.
			\item For any $B^{-1}\leq y\leq y'\leq B$ with $y'-y\geq \varepsilon$, we have $\big| \gamma_\tau(y)-\gamma_\tau(y') \big|\leq A|y-y'|$. 
		\end{enumerate}
	
		\noindent Then, we have
		\begin{align}\label{e:lder}
			\big\| \gamma_\tau - \gamma_{\tau} (1) \big\|_{\mathcal{C}^{\ell} ([2/B, B/2])} \leq C.
		\end{align}
		
	\end{prop}

	\begin{rem}
		
		The estimate \eqref{e:lder} reflects the idea that the free convolution has a regularizing effect on measures (in all derivatives). However, for this to be valid, it is necessary to impose some sort of assumption along the lines of the one supposed there, $\big| \gamma_{\tau} (y) - \gamma_{\tau} (y') \big| \le A |y-y'|$. Indeed, otherwise, it would be possible for $\supp \varrho_{\tau}$ to be disconnected; this could cause $\rho_{\tau}$, and thus $\gamma_{\tau}$, to have an infinite second derivative (at the edges of its support). The proof of \eqref{e:lder}, similarly to \eqref{e:densitybb0}, will make use of properties for the free convolution. 
		
	\end{rem}

\subsection{Density Upper and Lower Bound Estimates}

In this section we obtain upper and lower bounds for the density associated with a bridge-limiting measure process. We begin by stating three assumptions, which will be used at various points below (though not necessarily all at once). The first sets notation for the types of boundary measures; bridge-limited measure processes; and associated inverted height functions, inverted liquid regions, and density processes that we will consider in this chapter. In what follows, we recall the inverted height function and density process associated with a bridge-limiting measure process from \Cref{hrhot} (the latter of which exists by the first part of \Cref{mutrhot}), and the associated inverted liquid region from \Cref{omegafg}.

\begin{assumption}

\label{lmu}

Let $L \ge B \ge 10$ be real numbers and $\mu_0, \mu_1 \in \mathscr{P}_{\fin}$ be two measures. Assume that their total masses are $\mu_0 (\mathbb{R}) = L^{3/2} = \mu_1 (\mathbb{R})$, and that they satisfy $\supp \mu_0 \subseteq [-BL, 0]$ and $\supp \mu_1 \subseteq [-BL, 0]$. Let $\bm{\mu} = (\mu_t)$ denote the bridge-limiting measure process on $[0, 1]$ with boundary data $(\mu_0; \mu_1)$. Further denote the associated density process by $(\varrho_t)$; height function by $H : [0, 1] \times \mathbb{R} \rightarrow [0, L^{3/2}]$ as in \eqref{htxintegral}; inverted height function by $G:[0,1]\times[0, L^{3/2}] \rightarrow \bR$ as in \eqref{gty}; liquid region by $\Omega \subseteq (0, 1) \times \mathbb{R}$; inverted liquid region by $\Omega^{\inv} \subseteq (0, 1) \times (0, L^{3/2})$; the function $u : \Omega \rightarrow \mathbb{R}$ as in \eqref{u0}; and the complex slope $f : \Omega \rightarrow \overline{\mathbb{H}}$ as in \eqref{frhou}. Further define the function $\gamma: [0, 1] \rightarrow \mathbb{R}$ by setting $\gamma (t) = G(t, 0)$ for each $t \in [0, 1]$. 

\end{assumption}

Observe by \eqref{gty} that the curve $\gamma (t)$ traces the upper edge for the support of $\bm{\mu}$. We sometimes refer to it as the \emph{arctic boundary}\index{A@arctic boundary} associated with $\bm{\mu}$; see the left side of \Cref{f:liquidregion}. We first show the following result indicating that $G$, $u$, and $\varrho$ are real analytic on $\Omega$.

\begin{lem} 
	
	\label{gderivatives0}
	
	Adopting \Cref{lmu}, the functions $G$, $u$, and $\varrho$ are real anaytic on $\Omega$.
\end{lem} 

\begin{proof} 
	
	Fix a point $(t_0, x_0) \in \Omega$. By \Cref{omegafg}, $\varrho_{t_0} (y_0) > 0$ and so $\partial_y G(t_0, y_0) < 0$. Since $G$ is smooth on $\Omega$ by \Cref{urhoderivatives0}, there exist a real number $\varepsilon = \varepsilon (t_0, x_0) > 0$ and a neighborhood $U = U(t_0, x_0) \subset \Omega$ containing $(t_0, x_0)$ such that $-\varepsilon^{-1} < \partial_y G(t, x) < -\varepsilon$, for each $(t, x) \in U$. Thus, $G \in \Adm_{\varepsilon} (U)$ (recall from \Cref{functionadmissible}). Furthermore, by \Cref{gequation}, $G$ solves \eqref{equationxtd} on $U$, implying by \Cref{realanalytic} that $G$ is real analytic on $U$. By \eqref{u0} and \eqref{urho20}, this implies that $u$ and $\varrho$ are also real analytic on $\Omega$ (where for $\varrho$ we used the fact that $\partial_y G$ is bounded away from $0$ on $U$). Since $(t_0, x_0) \in \Omega$ was arbitrary, this confirms the lemma. 
\end{proof} 

The next two assumptions impose estimates on the boundary measures $\mu_0$ and $\mu_1$; the first states that its integrals are bounded above, and the second states that their densities are bounded below (which we formally express through an upper bound on the gaps of the associated inverted height function). 

\begin{assumption}

\label{integralmu0mu1} 

Adopt \Cref{lmu}; assume that 
	\begin{align}
	\label{e:avdensity}
	\int_{x}^0 \mu_0(dy) \le B |x|^{3/2}, \quad \text{and} \quad  \int_{x}^0  \mu_1(dy)\leq B|x|^{3/2},\qquad \text{for each $x \in [-BL, -1]$}.
\end{align}	

\end{assumption}

\begin{assumption}
	
	\label{gapmu0mu1} 
	
	Adopt \Cref{lmu}; assume that $G(0,0)=G(1,0)=0$. Further assume for all $(t, y) \in \{ 0, 1 \} \times [0, L^{3/2}]$ that $\partial_y G(t,y) \ge - B y^{-1/3}$; equivalently, for any real numbers $0\leq y\leq y'\leq L^{3/2}$, 
	\begin{align}\label{e:rholow}
		G(t,y)-G(t,y')\leq \frac{3B}{2} \big( (y')^{2/3}-y^{2/3} \big),\quad t\in \{0,1\}.
	\end{align}

\end{assumption}

The below result states that, under the integral bound \Cref{integralmu0mu1}, $\varrho_t$ is bounded above at intermediate times $t \in (0, 1)$. Its proof, which appears in \Cref{Proofrho0} below, uses \Cref{p:densityest} and the continuum height comparison \Cref{limitheightcompare}.

	\begin{prop}\label{p:densityub}
		
		Adopting \Cref{integralmu0mu1}, the following two statements hold. 
		\begin{enumerate}
		\item \label{i:Gupb} For each real number $t \in [0, 1]$, we have $\supp \varrho_t \subseteq [-2BL, 4B^2]$ and 
		\begin{align}\label{e:gtbound}
			G(t,r)\leq (2B)^2- \Big( \displaystyle\frac{r}{B} \Big)^{2/3}, \qquad \text{for each $r \in [0, L^{3/2}]$}.
		\end{align}
		\item \label{i:densityup} There exists a constant $C =C(B) > 1$ such that 
		\begin{align}\label{e:densitybb}
			\varrho_{t}(x)\le C  \max \{1,-x\}^{3/4}, \qquad \text{for each $(t, x) \in [B^{-1}, 1 - B^{-1}] \times \mathbb{R}$}.
		\end{align}
		\end{enumerate}
	\end{prop}

	The next result  states that, under the gap bound \Cref{gapmu0mu1}, $\varrho_t$ is bounded below at intermediate times $t \in (0, 1)$. Its proof appears in \Cref{ProofrhoG} and uses the continuum gap comparison \Cref{limitdifferencecompare}.  

	\begin{prop}\label{p:densitylow}
		
		Adopting \Cref{gapmu0mu1}, the following two statements hold. 
		
		\begin{enumerate}
		\item \label{i:Glowb}For any $t \in [0, 1]$, we have $\gamma(t) \geq 0$, and
		\begin{align}\label{e:Glow}
		G(t,r)\geq -3Br^{2/3},\qquad \text{for each $r \in [0,L^{3/2}]$}. 
		\end{align}
		\item \label{i:densitylow} We have $\Omega^{\inv} = (0, 1) \times (0, L^{3/2})$. Moreover,  we have 
		\begin{align}\label{e:rholowb}
			\varrho_t \big(G(t,r) \big)\geq \frac{r^{1/3}}{4B}\geq \bigg( \frac{ \gamma(t)-G(t,r)}{96B^3} \bigg)^{1/2}, \qquad \text{for any $(t, r) \in (0, 1) \times \Big[ 0, \displaystyle\frac{L^{3/2}}{2} \Big]$}. 
		\end{align}
		\end{enumerate}
	\end{prop}
	
	\subsection{Proof of Density Upper Bound} 
	
	\label{Proofrho0} 
	
	In this section we establish \Cref{p:densityub}. We first require the following lemma bounding the inverted height function $G$ at intermediate times by its values on the boundary.
	
	\begin{lem} 
		
		\label{g01y} 
		
		Adopting \Cref{lmu}, we have for each $(t, r) \in [0, 1] \times [0, L^{3/2}]$ that 
		\begin{align}
			\label{e:bound1}
			-L^{3/4}\le G(t,r)- \big((1-t)G(0,r)+t G(1,r) \big)\le L^{3/4}. 
		\end{align}
		
	\end{lem} 
	 
	 \begin{proof} 
	 			We only establish the lower bound in \eqref{e:bound1}, as the proof of the upper bound is entirely analogous. Fixing $r \in [0, L^{3/2}]$, we will compare $G$ with the limiting Brownian watermelon of \Cref{xuv0}, with the $(a, b; A; u, v)$ there equal to $\big( 0, 1; r, G(0, r), G(1, r) \big)$ here, so define (recalling $\gamma_{\semci}$ from \eqref{gammascy})
	 	\begin{flalign}
	 		\label{2gty} 
	 		G^-(t,y)= \big(r(1-t)t \big)^{1/2} \cdot \gamma_{\semci}(r^{-1} y)+(1-t) \cdot G(0,r)+t \cdot G(1,r).
	 	\end{flalign}
	 	
	 	\noindent  Then, for each $(t, y) \in \{ 0, 1 \} \times [0, r]$, we have $G(t,y)\ge G(t, r) = G^-(t,y)$, where the first bound follows from the fact that $G(t, y)$ is non-increasing in $y$ and the second follows from the definition \eqref{2gty} of $G^-$. Thus, by the first statement in \Cref{limitheightcompare} we have for each $(t, y) \in [0, 1] \in [0, r]$ that $G(t, y) \ge G^-  (t, y)$. At $y = r$, this implies
	 	\begin{align}\begin{split}
	 			G(t,r) \geq G^- (t,r) \ge(1-t) \cdot G(0,r)+t \cdot G(1,r) - r^{1/2} \ge (1-t) \cdot G(0,r)+t \cdot G(1,r)-L^{3/4},
	 	\end{split}\end{align}
	 	
	 	\noindent where in the second inequality we used \eqref{2gty} with the bound $\big( r(1-t)t \big)^{1/2} \gamma_{\semci}(1) \geq -r^{1/2}$ (as $\gamma_{\semci} (1) = -2$ by \eqref{gammascy} and $t(1-t) \le 1/4$), and in the third we used the bound $r \leq L^{3/2}$. This confirms \eqref{e:bound1}.
	 \end{proof} 
 
 	Now we can establish \Cref{p:densityub}.

	\begin{proof}[Proof of  \Cref{i:Gupb} in \Cref{p:densityub}]
		Since $\supp \mu_0 \subseteq [-BL, 0]$ and $\supp \mu_1 \subseteq [-BL,0]$, we have by \eqref{gty} that $G(0,y) \le 0$ and $G(1,y)\le 0$ for each $y \in [0, L^{3/2}]$. For $y \in [B, L^{3/2}]$ we have by \eqref{gty} and \eqref{e:avdensity} that $G(0,y) \le -(y/B)^{2/3}$ and $G(1,y)\leq -(y/B)^{2/3}$. Combining these yields 
		\begin{flalign}
			\label{g0yg1yyb}
			G(0,y) \le 1 - \Big( \displaystyle\frac{y}{B} \Big)^{2/3}, \quad \text{and} \quad G(1,y)\leq 1 - \Big( \displaystyle\frac{y}{B} \Big)^{2/3}, \qquad \text{for each $y \in [0, L^{3/2}]$}. 
		\end{flalign}

	 By taking $r=L^{3/2}$ in \eqref{e:bound1} (and using the fact that $G(0, y) \ge -BL$ and $G(1, y) \ge -BL$ for each $y \in [0, L^{3/2}]$, which holds by \eqref{gty} with the facts that $\supp \mu_0 \subseteq [-BL, 0]$ and $\supp \mu_1 \subseteq [-BL, 0]$), the lower bound in \eqref{e:bound1} implies that $G(t,L^{3/2})\ge -BL-L^{3/4}\geq -2BL$ (where in the last bound we used the fact that $L \ge B > 1$). By \eqref{gty}, this implies for each $t \in [0, 1]$ that 
	 \begin{flalign} 
	 	\label{2bl} 
	 	\supp \mu_t \subseteq [-2BL,\infty].
	 \end{flalign}   
	 	
		Next we prove \eqref{e:gtbound} To this end, we will compare $G$ to the limiting Airy profile of \Cref{airy} with the $(a, b; \mathfrak{a}, \mathfrak{b}, \mathfrak{c})$ there equal to $(0, 1; 1, \mathfrak{c}, \mathfrak{c})$ here, where $\fc=(3\pi B/4)^2$. So, for $(t, y) \in [0, 1] \times \mathbb{R}_{\ge 0}$, define
		\begin{align}\begin{split}\label{e:wGtlower1}
		G^+ (t,y)
		:=1+\fc(1-t)t -\left(\frac{3\pi}{4\fc^{1/2}}\right)^{2/3}y^{2/3}=1+\fc(1-t)t -\left(\frac{y}{B}\right)^{2/3}. 
		\end{split}\end{align}
	
	\noindent By \eqref{g0yg1yyb}, we then have the lower bound 
\begin{align}
G^+ (t,y)=1-\left(\frac{y}{B}\right)^{2/3}\geq G(t,y),\qquad \text{for $t\in \{0,1\}$}.
\end{align}
Thus the first statement in continuum height comparison \Cref{airyheightcompare} gives $G(t,y)\leq G^+ (t,y)$ for $(t, y) \in [0, 1] \times [0, L^{3/2}]$. Using the explicit formula \eqref{e:wGtlower1}, we get
\begin{align}\begin{split}\label{e:Gupb1}
G(t,y)
\leq G^+ (t,y)
\leq 1+\frac{\fc}{4}-\left(\frac{y}{B}\right)^{2/3}
\leq 1+2B^2-\left(\frac{y}{B}\right)^{2/3}\leq 4B^2-\left(\frac{y}{B}\right)^{2/3},
\end{split}\end{align}
where in the first inequality we used $\fc t(1-t)\leq \fc/4$; in the second inequality, we used $\fc=(3\pi B/4)^2\leq 8B^2$; and in the last inequality we used $1\leq 2B^2$. This finishes the proof of \eqref{e:gtbound}. By taking $y=0$ in \eqref{e:Gupb1} we get $G(t,0) \le 4B^2$, which with \eqref{gty} implies that $\supp \varrho_t \subseteq (-\infty, 4B^2]$. Together with \eqref{2bl}, this yields $\supp \varrho_t \subseteq [-2BL, 4B^2]$, verifying the first part of the proposition. 
\end{proof}

\begin{proof}[Proof of  \Cref{i:densityup} in \Cref{p:densityub}]

For $y \in \big[ (2B)^4, L^{3/2}]$, \eqref{e:gtbound} implies 
\begin{align}\label{e:Gtyab}
G(t,y)\leq (2B)^2-\left(\frac{y}{B}\right)^{2/3}\leq -\left(\frac{y}{(2B)^4}\right)^{2/3},
\end{align}

\noindent Let $C_1=(2B)^4$. From \Cref{i:Gupb} in \Cref{p:densityub}, we have  $\supp \varrho_t \subseteq [-2BL, 4B^2]$ for each $t \in [0, 1]$. So, by \eqref{gty} (with the fact that $\mu_t (\mathbb{R}) = L^{3/2}$), \eqref{e:Gtyab} implies
		\begin{align}\label{e:rhotd}
			\int_x^\infty \varrho_t(y)d y\leq C_1 |x|^{3/2}, \qquad \text{for $x \in [-C_1^{-2/3} L, -1]$}. 
		\end{align}
	
		 By the second part of \Cref{mutrhot}, for any $t \in [B^{-1}, 1-B^{-1}]$, there exists a measure $\nu_t$ with $\nu_t (\mathbb{R}) = L^{3/2}$ and $\supp \nu_t \subseteq \supp \mu_0 + \supp \mu_1 \subseteq [-2BL, 0]$, such that $\mu_t = \nu_t \boxplus \mu_{\semci}^{(\tau)}$ for $\tau = t(1-t)$. Since $B^{-1}(1-B^{-1})\leq \tau \leq 1/2$, \eqref{e:rhotd} verifies \Cref{blx122} (with the $B$ there equal to $(2B)^4$ here, and using the fact that the left side of \eqref{e:rhotd} is at most $L^{3/2} \le C_1 |x|^{3/2}$ for $x \le -C_1^{-2/3} L$), and so the second part of \Cref{p:densityest} yields \eqref{e:densitybb}.
	\end{proof}

\subsection{Proof of Density Lower Bound}
	
	\label{ProofrhoG} 
	
	In this section we establish \Cref{p:densitylow}.

	\begin{proof}[Proof of  \Cref{i:Glowb} in \Cref{p:densitylow}]
		
	 Since $G(0, 0) = 0 = G(1,0)$, taking $y=0$ and $y'=r$ in the assumption \eqref{e:rholow} gives
\begin{align}\label{e:GGbb}
-G(t,r) = G(t,0)-G(t,r)\leq \frac{3B}{2} \cdot r^{2/3},\qquad \text{for $(t, r) \in  \{0,1\} \times [1, L^{3/2}]$}. 
\end{align}

		Next we prove $\gamma(t)\geq 0$ by comparing $G$ to the limiting Airy profile as in \Cref{airy} with the $(a,b; \mathfrak{a}, \mathfrak{b}, \mathfrak{c})$ there equal to $(0,1; 0; \mathfrak{c}, \mathfrak{c})$ here, where $\fc=\pi^2/(12B^3)$. So, for $(t, y) \in [0, 1] \times \mathbb{R}_{\ge 0}$, define 
		\begin{align}\begin{split}\label{e:wGtlower2}
		G^- (t,y)
		:=\fc(1-t)t -\left(\frac{3\pi}{4\fc^{1/2}}\right)^{2/3} y^{2/3}=\fc(1-t)t -3By^{2/3}.
		\end{split}\end{align}
When $t\in \{0,1\}$, using \eqref{e:GGbb} and \eqref{e:wGtlower2} we deduce the upper bound 
\begin{align}\label{e:Gty}
 G^- (t,y)
=-3By^{2/3}\leq -\frac{3B}{2} \cdot y^{2/3}\leq  G(t,y),\qquad \text{for $(t, y) \in \{0,1\} \times [0, L^{3/2}]$}.
\end{align}

\noindent For $t\in [0, 1]$, we have from \eqref{e:bound1} that 
\begin{align*}
G(t,L^{3/2})
&\geq (1-t) \cdot G(0,L^{3/2})+t \cdot G(1,L^{3/2})- L^{3/4} \\
& \geq -\frac{3BL}{2}-L^{3/4} \geq\frac{\fc}{4} -3BL\geq \fc t(1-t) -3BL= G^- (t, L^{3/2}),
\end{align*} 
where the second statement is from  \eqref{e:GGbb} with $r=L^{3/2}$; the third holds since $B \ge 1$ and $L \ge 1$; the fourth uses $1/4\geq t(1-t)$; and the fifth uses the definition \eqref{e:wGtlower2} of $G^-$. Thus, the second statement in continuum height comparison \Cref{airyheightcompare} gives $G(t,y)\geq G^- (t,y)$ for each $(t, y) \in [0, 1] \times [0, L^{3/2}]$. Using  the explicit formula \eqref{e:wGtlower2} for $G^-$, it follows for each $(t, y) \in [0, 1] \times [0, L^{3/2}]$ that 
\begin{align}\begin{split}\label{e:Gupb2}
G(t,y)
\geq G^- (t,y)= \fc(1-t)t-3 By^{2/3}\geq -3By^{2/3},
\end{split}\end{align}

\noindent verifying \eqref{e:Glow}. Consequently, $\gamma(t)=G(t,0)\geq 0$ by setting $y=0$ in \eqref{e:Gupb2}, verifying the first statement of the first part of the proposition.
\end{proof}

\begin{proof}[Proof of \Cref{i:densitylow} in \Cref{p:densitylow}]
	
To prove \eqref{e:rholowb}, we will compare $G$ to the inverted height function from \Cref{rhoasct} with $(a, b) = (0, 1)$, $A = L^{3/2}$, and $d=8B^2L^2/\pi^2$. So, for $(t, y) \in [0, 1] \times [0, L^{3/2}]$, define 
\begin{flalign}
	\label{3gty}
	\widetilde{G} (t, y) = \bigg( d + \displaystyle\frac{A t(1-t)}{1 + 2\kappa} \bigg)^{1/2} \cdot \gamma_{\semci} (A^{-1} y),
\end{flalign}

\noindent where we recall $\kappa$ from \eqref{kappa} and the classical location $\gamma_{\semci} (y)$ from \eqref{gammascy}. Observe that
\begin{flalign}
	\label{estimatead} 
	d + \displaystyle\frac{At(1-t)}{1 + 2 \kappa } \le d + \displaystyle\frac{A}{8\kappa} \le d + \displaystyle\frac{A^2}{16d} = d + \displaystyle\frac{\pi^2 L}{128 B^2} \le 2d,
\end{flalign}

\noindent where in the first statement we used the bound $t(1-t) \le 1/4$; in the second we used the fact that $\kappa \ge 2A^{-1} d$ (which follows from \eqref{kappa}); in the third we used the definitions of $A$ and $d$; and in the fourth we used the definition of $d$ with the facts that $B \ge 1$ and $L \ge 1$. Thus, for $y \in [0, A/2]$, we obtain 
\begin{align}\label{e:tGderup}
 -\partial_y\widetilde G (t, y) 
 = - A^{-1} \bigg( d + \displaystyle\frac{A(1-t)t}{1 + 2\kappa} \bigg)^{1/2} \gamma'_{\semci} \Big( \displaystyle\frac{y}{A} \Big)
 \leq   
 	\bigg( d + \displaystyle\frac{A(1-t)t}{1 + 2\kappa} \bigg)^{1/2}  \frac{\pi}{A^{2/3}y^{1/3}}\leq \frac{4B}{y^{1/3}},
\end{align}

\noindent where in the first statement we differentiated \eqref{3gty} with respect to $y$; in the second we applied the second part of \Cref{gammaderivative}; and in the third we used \eqref{estimatead} and the definitions of $A$ and $d$. Moreover, for $y \in [0, A]$, we have 
\begin{align}\label{e:tGderlow}
\frac{B}{y^{1/3}}=\frac{d^{1/2}\pi}{2^{3/2} A^{2/3}y^{1/3}}\leq - A^{-1} \bigg( d + \displaystyle\frac{At(1-t)}{1 + 2\kappa} \bigg)^{1/2} \cdot \gamma_{\semci}' \Big( \displaystyle\frac{y}{A} \Big) = -\partial_y\widetilde G (t, y),
\end{align}

\noindent where in the first statement we used the definitions of $A$ and $d$, in the second we used the second part of \Cref{gammaderivative}, and in the third we used the definition \eqref{3gty} of $\widetilde{G}$. Thus, for any $0\leq y\leq y' \le A = L^{3/2}$, we have
\begin{align}
\widetilde G(t,y)-\widetilde G(t,y')\geq \int_{y}^{y'}\frac{Bd r}{r^{1/3}}= \frac{3B}{2} \big((y')^{2/3}-y^{2/3} \big) \geq G(t,y)-G(t,y'),
\end{align}
where the first statement follows from integrating \eqref{e:tGderlow}; the second from performing the integral; and the third from \eqref{e:rholow}.
This verifies the assumption in the first statement of continuum gap comparison \Cref{limitdifferencecompare};  and we conclude for $t \in [0, 1]$ and $0 \le y \le y' \le L^{3/2}$ that 
\begin{align}\label{e:tGG}
\widetilde G(t,y)-\widetilde G(t,y')\geq G(t,y)-G(t,y').
\end{align}

\noindent Since $\widetilde{G} (t, y)$ is differentiable (and has negative derivative) in $y \in (0, L^{3/2})$ for any $t \in (0, 1)$, this implies (by \Cref{hrhot}) that $\varrho_t \big( G(t, y) \big) > 0$ for any $(t, y) \in (0, 1) \times (0, L^{3/2})$, meaning by \eqref{omega12} that $\Omega^{\inv} = (0, 1) \times (0, L^{3/2})$. This establishes the first statement of the second part of the proposition.

Next, by integrating the upper bound \eqref{e:tGderup}, together with \eqref{e:tGG}, we deduce 
\begin{align}\label{e:Gtdiff}
G(t,y)-G(t,y')\leq  \widetilde G(t,y)-\widetilde G(t,y')\leq \int_y^{y'} \frac{4Bd r}{r^{1/3}} = 6B \big( (y')^{2/3}-y^{2/3} \big),
\end{align}
for any $0 < t < 1$ and $0\leq y < y'\leq L^{3/2}/2$. From  \Cref{hrhot} (and recalling from \Cref{mutrhot} that for $t \in (0, 1)$ that the density $\varrho_t$ exists), \eqref{e:Gtdiff} implies for $(t, y) \in (0, 1) \times ( 0, L^{3/2} / 2 ]$ that
	\begin{align}\label{e:drhat2}
			\frac{1}{\varrho_t \big(G(t,y) \big)}=-\partial_y G(t,y)
			=\lim_{y'\rightarrow y^+}\frac{G(t,y)-G(t,y')}{y'-y} \le \frac{4B}{y^{1/3}}.
		\end{align}
By rearranging, this gives the first inequality in \eqref{e:rholowb}.

 To establish the second, recall that $\gamma(t)=G(t,0)$ and take $y=0$ and $y'=r\leq L^{3/2}/2$ in \eqref{e:Gtdiff}; this yields for $t \in [0, 1]$ the bound 
		\begin{align}\label{e:Gloc}
		\gamma(t)-G(t,r)=G(t,0)-G(t,r)\leq 6B r^{2/3}.
		\end{align}
By plugging \eqref{e:Gloc} into \eqref{e:drhat2}, we obtain for $r \in (0, L^{3/2}/2]$ that 	
\begin{align}
\varrho_t \big( G(t,r) \big)\geq \frac{r^{1/3}}{4B}\geq \frac{1}{4B} \bigg( \frac{\gamma(t)-G(t,r)}{6B} \bigg)^{1/2} = \bigg( \frac{\gamma(t)-G(t,r)}{96B^3} \bigg)^{1/2},
\end{align}

\noindent which finishes the proof of the second inequality in \eqref{e:rholowb} for $r \in ( 0, L^{3/2} / 2 ]$. At the endpoint $r = 0$ of this interval, the second bound in \eqref{e:rholowb} continues to hold by the nonnegativity of $\varrho_t$. This verifies the second part of the proposition.
	\end{proof}

	\subsection{Regularity Estimates for Limit Shapes}
	
	\label{GDerivativeConvexgamma} 
	
	In this section we state the following two propositions providing estimates on the inverted height functions $G$ subject to the integral bound \Cref{integralmu0mu1} and gap bound \Cref{gapmu0mu1}. The first provides approximately matching bounds on the $y$-derivatives of $G$, and also shows as \eqref{e:gammaconv} that the arctic  boundary $\gamma$ is uniformly concave (namely, $\gamma''$ is bounded above and below, in the weak sense); we establish it later in this section. The second shows that the functions $u$ and $\varrho$ (recall \Cref{lmu} for their definitions) extend continuously to its arctic boundary; we establish it in \Cref{ProofContinuousu} below.

	\begin{prop}\label{p:checka}
		
		Adopt \Cref{integralmu0mu1} and \Cref{gapmu0mu1}.  There exist constants $c = c(B)>0$ and $C=C(B) > 1$ such that the following two statements hold if $L \ge C$, for any real number $t \in [3B^{-1}, 1 - 3B^{-1}]$. 
		\begin{enumerate}
			\item \label{i:dyG} For any $y \in (0, B^5]$, we have $c y^{-1/3}\leq-\partial_y G(t,y)\leq 4By^{-1/3}$.
			\item \label{i:concave} For any real number $t' \in [3B^{-1}, 1-3B^{-1}]$, we have  $\big| \gamma(t) - \gamma(t') \big| \le C|t-t'|$. Moreover, for any real numbers $t_0, s, \tau_0 \in \mathbb{R}$ with $3B^{-1}\leq t_0-\tau_0 <t_0+\tau_0 \leq 1-3B^{-1}$ and $s \in [-\tau_0, \tau_0]$, we have   
			\begin{align}\label{e:gammaconv}
				C^{-1} (\tau_0^2-s^2) \le \gamma(t_0+s)-\left(\frac{\tau_0-s}{2\tau_0} \cdot \gamma(t_0-\tau_0)+\frac{\tau_0+s}{2\tau_0} \cdot \gamma(t_0+\tau_0) \right)\le C(\tau_0^2-s^2).
			\end{align}
		\end{enumerate}
	\end{prop}
	
	\begin{prop}\label{p:densityrg}
		
		Adopt \Cref{integralmu0mu1} and \Cref{gapmu0mu1}. There exists a constant $C = C(B) > 1$ such that the following two statements hold if $L \ge C$.
	\begin{enumerate}
	\item For any $x_0=G(t_0,y_0)$ with $t_0\in [4B^{-1},1-4B^{-1}]$ and $y_0 \in (0, B]$, we have		
	\begin{align}\label{e:holder}
	\big|\partial_x u(t_0,x_0) \big| + \big|\partial_x \varrho(t_0,x_0) \big| \le C \big(\gamma(t_0)-x_0 \big)^{-1/2}.
	\end{align}
	\item	
	Both $\varrho(t,x)$ and $u(t,x)$ extend continuously to the set $\big\{ (t,\gamma(t)): 4B^{-1}\le t\le 1-4B^{-1} \big\}$, with $\varrho \big(t,\gamma(t) \big)=0$ and $u \big(t,\gamma(t) \big)=\gamma'(t)$. In particular, 
	\begin{align}\label{e:gammadc}
	\gamma'(t) \text{ is continuous in } t \in [4B^{-1}, 1 - 4B^{-1}].
	\end{align}
	\end{enumerate}
	\end{prop}

	\begin{rem}

	Let us comment on \eqref{e:holder}, focusing on the bound $\big| \partial_x \varrho (t_0, x_0) \big| \le C \big( \gamma(t_0) - x_0 \big)^{-1/2}$. Using \eqref{grhotx}, it can be seen that this estimate is optimal (up to the constant $C$) for the explicit limiting Airy profile $\mathfrak{G}_{\Ai}$ given by \eqref{e:Airyprofile}. For general inverted height functions $G$ (satisfying \Cref{gapmu0mu1}), $\varrho (t_0, x_0)$ decays as the square root of the distance from $x_0$ to the arctic boundary, namely, $\varrho (t_0, x_0) \sim \big( \gamma(t_0) - x_0 \big)^{1/2}$ (the upper bound essentially follows from \eqref{e:rholowb}; the lower bound can be deduced as a consequence of \Cref{i:dyG} of \Cref{p:checka}, together with \eqref{e:gtbound}). If we could differentiate this estimate in $x_0$, we would obtain $\big| \partial_x \varrho (t_0, x_0) \big| \le C \big( \gamma(t_0) - x_0 \big)^{-1/2}$. 
	
	\end{rem}

	We now establish \Cref{p:checka}.

\begin{figure}
	\center
\includegraphics[width=1\textwidth]{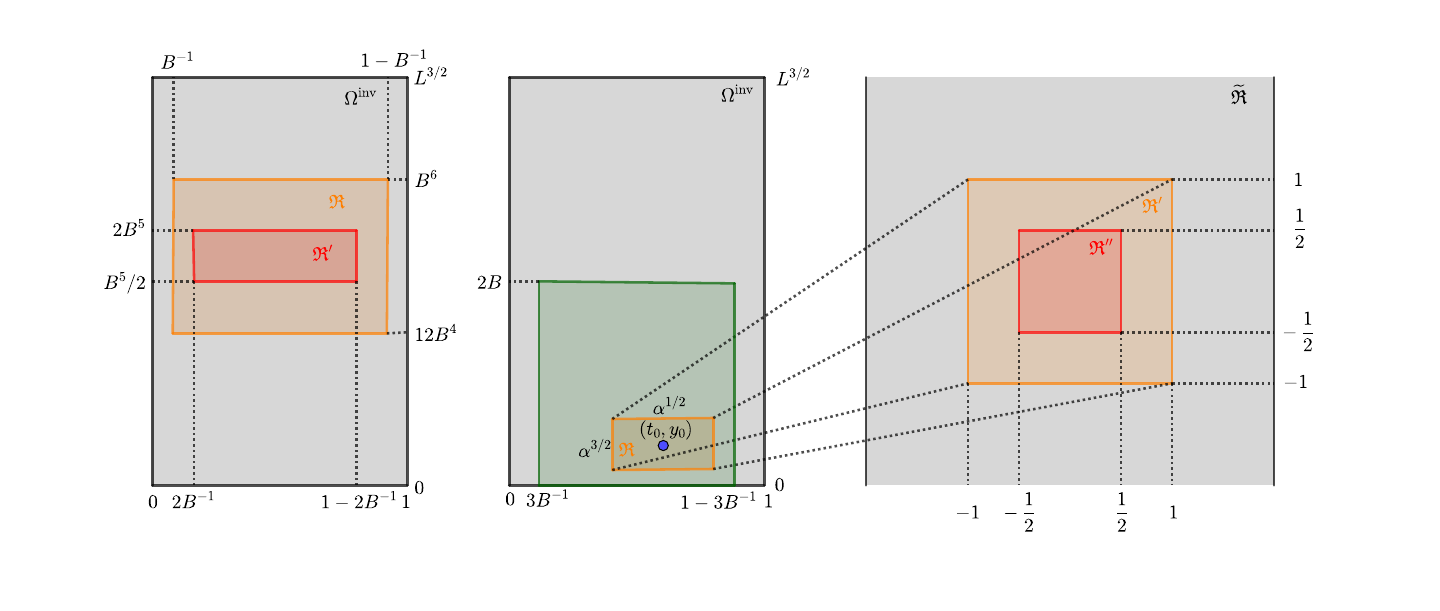}

\caption{Shown on the left is a depiction for the proof of \Cref{i:dyG} in \Cref{p:checka}.  Shown on the right is the ``zooming in'' procedure implemented in the proof of \Cref{p:densityrg}.}

\label{f:Rescale_Omega}
	\end{figure}

	\begin{proof}[Proof of \Cref{i:dyG} in \Cref{p:checka}]

		First observe that for any $(t, y) \in [0, 1] \times [12B^4, B^6]$ we have 
		\begin{align}\label{e:Gbound}
		-3B^5\leq G(t,y)\leq 4B^2- (12B^3)^{2/3}\leq -B^2,
		\end{align}
		
		\noindent where the lower bound is from \eqref{e:Glow} and the upper bound is from \eqref{e:gtbound}. Moreover, for any $(t, y) \in [B^{-1}, 1-B^{-1}] \times [12B^4, B^6]$, we have for some constant $C_1 = C_1 (B) > 1$ that
\begin{align}\label{e:derbound}
		C_1^{-1} \le \frac{1}{\varrho_t \big( G(t, y) \big)} = -\partial_y G(t, y) \leq B^{-1} \cdot (96B^3)^{1/2} = (96B)^{1/2},		
		\end{align}
	
		\noindent where the first statement is from \eqref{e:densitybb} and  \eqref{e:Gbound}; the second is from \eqref{urho20} (and the fact from \Cref{i:densitylow} of \Cref{p:densitylow} that $\Omega^{\inv} = (0, 1) \times (0, L^{3/2})$); and the third is from \eqref{e:rholowb}, as $\gamma(t) - G(t, y) \ge -G(t,y) \ge B^2$ (where in the first inequality we used the fact that $\gamma(t) \ge 0$, from \Cref{i:Glowb} of \Cref{p:densitylow}, and in the second we used \eqref{e:Gbound}).

		 Now define the open rectangle $\fR= (B^{-1}, 1-B^{-1}) \times (12 B^4, B^6)$; see the left side of \Cref{f:Rescale_Omega}. By \Cref{gequation}, $G$ solves \eqref{equationxtd} on $\Omega^{\inv} = (0, 1) \times (0, L^{3/2})$ (where the latter follows from \Cref{i:densitylow} of \Cref{p:densitylow}). Moreover, $G$ and $-\partial_y G(t, y)$ are bounded above and below on $\mathfrak{R}$, by \eqref{e:Gbound} and \eqref{e:derbound}. Hence, recalling \Cref{functionadmissible}, there exists some constant $\varepsilon = \varepsilon (B) > 0$ such that $G \in \Adm_{\varepsilon} (\mathfrak{R})$ (where we observe that $G$ is Lipschitz on $\mathfrak{R}$ since it is real analytic by \Cref{gderivatives0}); this will enable us to apply the regularity results of \Cref{EstimatesEquation0} to $G$. In particular, denoting the open rectangle $\mathfrak{R}' = (2B^{-1}, 1-2B^{-1}) \times (B^5/2, 2 B^5) \subset \mathfrak{R}$, \Cref{derivativef} yields a constant $M(B)=M > 1$ with 
		\begin{align}\label{e:der2Gbound}
		 \displaystyle\sup_{(t, y) \in \overline{\mathfrak{R}}'} \big| \partial_t^2G(t,y) \big|\leq M.
		 \end{align}
	 
	 	\noindent We may assume in what follows that $M \ge 10B^5$. Set $A = M^2/4$.
		
		We next use \eqref{e:der2Gbound} to establish the first statement of the proposition, bounding $-\partial_y G(t, y)$ from above and below. The lower bound is given by the first estimate in \eqref{e:rholowb}, together with the fact (recall \eqref{urho20}) that $-\partial_y G(t, y) = \varrho_t \big( G(t,y) \big)^{-1}$. To prove the upper bound, we use \Cref{limitdifferencecompare} to compare $G$ with the limiting Brownian watermelon of \Cref{xuv0} with the $(a, b; A; u, v)$ there equal to $\big( 2B^{-1}, 1-2B^{-1}; M^2/4; G(2B^{-1}, A), G(1-2B^{-1}, A) \big)$ here. So, define $\widetilde{G} : [a, b] \times [0, A] \rightarrow \mathbb{R}$ by
		\begin{flalign}
			\label{gty02} 
			\widetilde{G} (t, y) = \bigg(  \displaystyle\frac{A(b-t)(t-a)}{b-a} \bigg)^{1/2} \gamma_{\semci} \Big( \displaystyle\frac{y}{A} \Big) + \displaystyle\frac{b-t}{b-a} \cdot G(2B^{-1}, A) + \displaystyle\frac{t-a}{b-a} \cdot G(1-2B^{-1}, A),
		\end{flalign}
	
		\noindent where we recall the classical location $\gamma_{\semci}$ of the semicircle law from \eqref{gammascy}.

		Since we seek to use $\widetilde{G}$ to upper bound $-\partial_y G$, the second part of the continuum gap comparison \Cref{limitdifferencecompare} indicates that we must lower bound $\partial_t^2 G$ by $\partial_t^2 \widetilde{G}$. To that end, observe for $2B^{-1}\le t\le 1-2B^{-1}$ that   
		\begin{align}
			\label{gtb5} 
			\partial^2_t\widetilde G(t,B^5)= -\frac{A^{1/2}(b-a)^{3/2}}{4(b-t)^{3/2}(t-a)^{3/2}} \cdot \gamma_{\semci} \Big( \displaystyle\frac{B^5}{A} \Big)\leq -2A^{1/2} \cdot \gamma_{\semci} \Big( \displaystyle\frac{B^5}{A} \Big)\leq -2A^{1/2} = -M,
		\end{align}
		where the first statement follows from \eqref{gty02}; the second follows from the fact that $(b-t)(t-a) \le (b-a)^2 / 4 \le (b-a) / 4$ (as $b-a = 1 - 4B^{-1} < 1$); the third follows from the facts that $A = M^2 / 4 \geq 25B^5$ (as $M \ge 10 B^5$ and $B > 1$) and $\gamma_{\semci}(1/25)\geq 1$ (by the first part of \Cref{gammaderivative}); and the fourth from the fact that $A = M^2/4$. To upper bound $-\partial_y \widetilde{G}$, observe for $t\in [3B^{-1}, 1-3B^{-1}]$ that 
		\begin{align}\label{e:LtG}
		-\partial_y \widetilde G(t,y)
		=-\left(\frac{ (b-t) (t-a)}{A(b-a)} \right)^{1/2}\gamma'_{\semci} \Big(\frac{y}{A} \Big) 
		\geq \left(\frac{ B^{-1} (1-5B^{-1})}{A(1-4B^{-1})} \right)^{1/2}\frac{\pi A^{1/3}}{2^{3/2} y^{1/3}},
		\end{align}
	where the first statement is from \eqref{gty02}, and the second is from bounding $(b-t)(t-a)\geq B^{-1}(1-5B^{-1})$ (by the fact that $t \in [3B^{-1}, 1-3B^{-1}]$) and $\gamma'_{\semci} (y/A)\geq (\pi A^{1/3})/(2^{3/2} y^{1/3})$ (by the second part of \Cref{gammaderivative}). 
	
		Now, \eqref{e:der2Gbound} and \eqref{gtb5} together yield $\partial^2_t \widetilde G(t,B^5) \le \partial^2_t G(t,B^5)$ for each $t \in [2B^{-1}, 1-2B^{-1}]$.  Since $\widetilde G(2B^{-1},y) = u$ and $\widetilde{G} (1 - 2B^{-1}, y) = v$ are constant in $y$, we also have $\big|\widetilde G(t,y)-\widetilde G(t,y') \big| = 0 \leq \big|G(t,y)-G(t,y') \big|$ for each $t\in\{2B^{-1},1-2B^{-1}\}$ and $0\le y\leq y'\le A$. This verifies the assumptions in the second statement of \Cref{limitdifferencecompare} (with the $(G^{\star}, \widetilde{G}^{\star})$ there equal to $(\widetilde{G}, G)$ here), which gives for each $(t, y) \in [3B^{-1}, 1-3B^{-1}] \times (0, B^5]$ that 
		\begin{align*}
			-\partial_y G(t,y) & = \displaystyle\lim_{y' \rightarrow y^-} \displaystyle\frac{G(t, y') - G(t, y)}{y-y'} \\
			& \ge \displaystyle\lim_{y' \rightarrow y^-} \displaystyle\frac{\widetilde{G} (t, y') - \widetilde{G} (t, y)}{y-y'} = -\partial_y \widetilde G(t,y)\geq \left(\frac{B-5}{AB(B-4)} \right)^{1/2}\frac{\pi A^{1/3}}{2^{3/2} y^{1/3}},
		\end{align*}
		
		\noindent where the last inequality is from \eqref{e:LtG}. This provides the lower bound on $-\partial_y G(t, y)$ and thus finishes the proof of the first statement in \Cref{p:checka}. 
	\end{proof}

	\begin{proof}[Proof of \Cref{i:concave} in \Cref{p:checka}]	
		Fix an interval $[a, b]\subseteq [3B^{-1},1-3B^{-1}]$, and denote $\tau=(b-a)/2$. Define the functions $\breve{G} : [0, 1] \times [0, L^{3/2}] \rightarrow \mathbb{R}$ and $\breve{\gamma} : [0, 1] \rightarrow \mathbb{R}$ by performing an affine shift on $G$ and $\gamma$ respectively, by setting 
		\begin{align}\label{e:breveG}
			\breve G(t,y) =G(t,y)-\left(\frac{b-t}{2\tau} \cdot \gamma(a)+\frac{t-a}{2\tau} \cdot \gamma(b)\right), \qquad \text{and} \qquad \breve{\gamma} (t) = \breve{G}(t, 0),
		\end{align}
	
		\noindent for each $(t, y) \in [0, 1] \times [0, L^{3/2}]$. Then, for each $t \in [0, 1]$, we have 
		\begin{flalign}
			\label{gammaab0} 
			\breve \gamma'(t)=\gamma'(t) + \displaystyle\frac{\gamma(b)-\gamma(a)}{2\tau}, \qquad \text{and} \qquad \breve \gamma(a)=\breve\gamma(b)=0.
		\end{flalign}
	
		\noindent  From the first statements of \Cref{p:densityub} and \Cref{p:densitylow}, we have $0\leq \gamma(t)=G(t,0)\leq (2B)^2$, which upon insertion into \eqref{e:breveG} gives for each $t \in [a, b]$ that 
		 \begin{align}\label{e:breveGup}
		 \breve G(t,0)=\breve \gamma(t)\leq  (2B)^2.
		 \end{align}
		
		Next we show that there exists a constant $C_1 = C_1 (B) > 1$ such that, for each $t \in [a, b]$, 
		\begin{align}\label{e:gammatb}
			C_1^{-1} (b-t)(t-a)  \le\breve\gamma(t)\le C_1 (b-t)(t-a)
		\end{align}
		
		\noindent We only prove the upper bound in \eqref{e:gammatb}, as the proof of the lower bound is entirely analogous. To this end, recall from the first statement in \Cref{p:checka} that there exists some constant $c=c(B)>0$ such that $-\partial_y \breve G(t,y) = -\partial_y G(t, y) \geq c y^{-1/3}$ for each $(t, y) \in [a, b] \times (0, B^5]$. Integrating this estimate from $0$ to $r$ then gives for each $(t, r) \in [a, b] \times [0, B^5]$ that
		\begin{align}\label{e:tGlow0}
			  \breve G(t,r) \le \breve{G} (t, 0) - \displaystyle\frac{3c}{2} \cdot r^{2/3} = \breve \gamma(t)-\frac{3c}{2} \cdot r^{2/3}\leq (2B)^2-\frac{3c}{2} \cdot r^{2/3},
		\end{align}
		where we used \eqref{e:breveGup} in the last inequality. By \eqref{e:breveG} and \eqref{e:gtbound} (with the fact that $\gamma(t) \ge 0$, by \Cref{p:densitylow}), we also have $ \breve G(t,r)\leq G(t,r)\leq (2B)^2-(r/B)^{2/3}$. Together with \eqref{e:tGlow0}, this implies for $(t, r) \in [a, b] \times (0, B^5]$ that  
		\begin{align}\label{e:tGlow1}
		 \breve G(t,r)\leq (2B)^2- r^{2/3} \cdot \max\left\{\frac{3c}{2}, B^{-2/3}\right\}.
		\end{align}
		
		Using this, we compare $\breve{G}$ to the limiting Airy profile $\widetilde G(t,y)$ from \eqref{e:Airyprofile} with $\fc= 243 \pi^2 / (2c^3)$, $\fa=-ab\fc$, and $\fb=(a+b)\fc$, so define the function $\widetilde{G} : [a, b] \times \mathbb{R}_{\ge 0} \rightarrow \mathbb{R}$ by setting
		\begin{align}\begin{split}\label{e:wGtlower}
		\widetilde G(t,y) =\fc(b-t)(t-a) -\left(\frac{3\pi}{4\fc^{1/2}}\right)^{2/3}y^{2/3}=\fc(b-t)(t-a) -\frac{c}{6} \cdot y^{2/3}.
		\end{split}\end{align}
		
		\noindent Then it follows that for each $(t, y) \in \{a,b\} \times [0, B^5]$ we have  
		\begin{align}\begin{split}\label{e:hatGtlow2}
		&\widetilde G(t,y)= -\frac{c}{6} \cdot y^{2/3}
		\geq -\frac{3c}{2} \cdot y^{2/3}\geq  \breve G(t,y),
		\end{split}\end{align}	
		\noindent	where the first equality is from \eqref{e:hatGtlow2}, the second inequality follows from $3c/2\geq c/6$ and the third inequality uses the first bound in \eqref{e:tGlow0} and the equalities $\breve{G} (a, 0) = 0 = \breve{G} (b, 0)$ (by \eqref{gammaab0}), and \eqref{e:wGtlower}. 		Moreover, for each $t\in [a,b]$, we have 
		\begin{align*}
			\widetilde G(t,B^5)
			\geq -\frac{c}{6} \cdot B^{10/3}
			& \ge (2B)^2 - \displaystyle\frac{8}{9} \cdot B^{8/3} - \displaystyle\frac{c}{6} \cdot B^{10/3} \\
			& \geq (2B)^2-B^{10/3} \cdot \max\left\{\frac{3c}{2}, B^{-2/3}\right\}
			 \geq \breve G(t,B^5),
		\end{align*}
	where the first inequality is from \eqref{e:wGtlower} and the fact that $(b-t)(t-a)\geq 0$; the second inequality is from the bound $(2B)^2 \le 8B^{8/3} / 9$, as $B \ge 10$ (recall \Cref{lmu}); the third is from the fact that $c / 6 \leq 3c / 2$; and the fourth is from \eqref{e:tGlow1}.
	This verifies the assumptions in the second statement of \Cref{airyheightcompare} (using \Cref{gtabab}, \eqref{e:breveG}, and \Cref{linear} to confirm that the restriction of $\breve{G}$ to $[a, b] \times [0, L^{3/2}]$ is the inverted height function associated with a bridge-limiting measure process), which yields
		\begin{align*}
			\breve{\gamma} (t, 0) = \breve G(t,0)
			\le \widetilde G(t,0)=\fc(b-t)(t-a)= \frac{243 \pi^2}{2c^3}(b-t)(t-a).
		\end{align*}
	This gives the upper bound in \eqref{e:gammatb}. The proof of the lower bound is very similar, obtained by comparing $\breve{G}$ to a limiting Airy profile from \eqref{e:Airyprofile} with $\mathfrak{c} = \pi^2 / (2^{15} B^3)$, $\mathfrak{a} = -ab\mathfrak{c}$, and $\mathfrak{b} = (a+b) \mathfrak{c}$ (using the bound $\breve{G} (t, B^5) \ge G(t, B^5) - 4B^2 \ge -7 B^{13/3}$, which holds by \eqref{e:breveG}, \eqref{e:Glow} and \eqref{e:gtbound}, and the upper bound in \Cref{i:dyG} of \Cref{p:checka} in place of the lower bound there); further details are therefore omitted.
	
	By \eqref{e:breveG} and the fact that $2 \tau = b - a$, we can rewrite \eqref{e:gammatb} as   
		\begin{align}\label{e:gammabd}
			\begin{aligned} 
				C_1^{-1} (b-t)(t-a) & \le
				\gamma(t)-\left(\frac{b-t}{2\tau} \cdot \gamma(a)+\frac{t-a}{2\tau} \cdot \gamma(b)\right)\\
				&= \gamma(t)-\gamma(a)-(t-a) \cdot \frac{\gamma(b)-\gamma(a)}{2\tau}
				\le C_1 (b-t)(t-a).
				\end{aligned} 
		\end{align}
	
		\noindent The claim \eqref{e:gammaconv} follows from this by taking the $C_1$ here equal to $C$ there; the $(a, b)$ here equal to $(t_0 - \tau_0, t_0 + \tau_0)$ there (so that the $\tau$ here is equal to $\tau_0$ there); and the $t$ here equal to $t_0 + s$ there. This in particular implies that $\gamma(t)$ is concave.
		
		It thus remains to verify the bound $\big| \gamma (t) - \gamma (t') \big| \le C |t - t'|$ for any $t, t' \in [3B^{-1}, 1-3B^{-1}]$. We may suppose that $t > t'$ by symmetry, and similarly that $t' < 1 / 2$; set $(a, b) = (t', 1-3B^{-1})$ in \eqref{e:gammabd}, which guarantees that $2 \tau = 1 - 3B^{-1} - t' \ge 1 / 2 - 3 / 10 \ge 1 / 5$, so that $\tau \ge 1 / 10$. Due to the  bound $0\leq \gamma(t)\leq (2B)^2$ (from the first statements of \Cref{p:densityub} and \Cref{p:densitylow}), it follows from \eqref{e:gammabd} that  
		\begin{align}
		\displaystyle\frac{\big| \gamma (t) - \gamma(t') \big|}{t-t'}=\displaystyle\frac{\big| \gamma (t) - \gamma(a) \big|}{t-a} \le \frac{\big| \gamma(b)-\gamma(a) \big|}{2\tau}+C_1 (b-a)\leq \frac{4B^2}{\tau}+C_1 \le 40B^2 + C_1,
		\end{align}
	
		\noindent which establishes the second statement in \Cref{p:checka}. 		
	\end{proof}

\subsection{Continuous Extensions for $u$ and $\varrho$}

\label{ProofContinuousu}
	
	In this section we establish \Cref{p:densityrg}.
	 
	\begin{proof}[Proof of \Cref{p:densityrg}]

		By \Cref{p:checka}, there exists a constant $D = D(B) > 1$ such that the following two statements hold.  
			\begin{enumerate}
			\item \label{i:dyGcopy}For each $(t, y) \in[3B^{-1}, 1-3B^{-1}] \times (0, B^5]$, we have 
			\begin{flalign}
			\label{yd0} 
			D^{-1}y^{-1/3}\leq-\partial_y G(t,y)\leq Dy^{-1/3}.
			\end{flalign}
				
			\noindent Integrating this bound and using the fact that $\gamma(t) = G(t, 0)$, we obtain for each $(t, y) \in [3B^{-1}, 1 - 3B^{-1}] \times [0, B^5]$ that 
			\begin{align}\label{e:gammaGdiff}
			\displaystyle\frac{3}{2D} \cdot y^{2/3}\leq \gamma(t)-G(t,y)\leq  \displaystyle\frac{3D}{2} \cdot y^{2/3}.
			\end{align}
		
			\noindent Together with \eqref{yd0} and \eqref{urho20}, this gives  
			\begin{align}\label{e:densitybound}
			\bigg( \frac{2(\gamma(t)-G(t,y))}{3D^3} \bigg)^{1/2}  \leq D^{-1} y^{1/3} \leq \varrho \big(t, G(t,y) \big)\leq y^{1/3}D\leq \bigg( \frac{2D^3(\gamma(t)-G(t,y))}{3} \bigg)^{1/2}.
			\end{align} 
			\item \label{i:concavecopy} The function $\gamma(t)$ is concave for $t \in [3B^{-1}, 1-3B^{-1}]$; moreover, for any $t, t' \in [3B^{-1}, 1-3B^{-1}]$, we have $\big| \gamma(t) - \gamma (t') \big| \le D|t-t'|$. Furthermore, for any real numbers $t_0, s, \tau \in \mathbb{R}$ with $3B^{-1}\leq t_0-\tau<t_0<t_0+\tau\leq 1-3B^{-1}$ and $s \in [-\tau, \tau]$, we have  
			\begin{align}\label{e:gammaconv2}
				\left| \gamma(t_0+s)-\Big(\frac{\tau-s}{2\tau} \cdot \gamma(t_0-\tau)+\frac{\tau+s}{2\tau} \cdot \gamma(t_0+\tau) \Big)\right|\leq D\tau^2.
			\end{align}
				\end{enumerate}
			
			Now fix some point $(t_0, y_0) \in [4B^{-1}, 1 - 4B^{-1}] \times (0, B]$; set $x_0=G(t_0,y_0)$; denote $\alpha=\min \big\{(y_0/2)^{2/3}, B^{-2} \big\}$; and define the open rectangle 
			\begin{align}
			\label{2r}
			\fR= (t_0-\alpha^{1/2}, t_0+\alpha^{1/2}) \times (y_0-\alpha^{3/2}, y_0+\alpha^{3/2}),
			\end{align} 
			which is centered at $(t_0,y_0)$. Then we have 
		\begin{align}\begin{split}\label{e:ayrelation}
		 \alpha^{3/2}\leq \frac{y_0}{2},\qquad \alpha^{1/2}\leq \frac{1}{B}, \qquad  1\geq \alpha \Big( \displaystyle\frac{y_0}{2} \Big)^{-2/3}\geq\min \bigg\{1,  \displaystyle\frac{2^{2/3}}{B^{8/3}} \bigg\}\geq B^{-3}, 
		\end{split}\end{align}
		which holds by our choice of $\alpha$ and the bound $y_0\leq B$, and $\fR\subset [3B^{-1},1-3B^{-1}]\times [0, 2B]$. Observe from \Cref{gequation} and \Cref{i:densitylow} of \Cref{p:densitylow} that $G$ solves \eqref{equationxtd} on $\Omega^{\inv} = (0, 1) \times (0, L^{3/2})$. We next rescale $G$ by ``zooming into'' the point $(t_0, y_0)$; see the right side of \Cref{f:Rescale_Omega}.  More specifically, define the rescaled rectangles
		\begin{flalign*}
		\widetilde{\mathfrak{R}} = \big\{ (t, y) \in \mathbb{R}^2 : (\alpha^{1/2} t + t_0, \alpha^{3/2} y + y_0) \in [0, 1] \times [0, L^{3/2}] \big\}; \quad \mathfrak{R}' = [-1, 1]^2; \quad \mathfrak{R}'' = \Big[ -\displaystyle\frac{1}{2}, \displaystyle\frac{1}{2} \Big]^2,
		\end{flalign*} 
	
		\noindent and the function $\widetilde{G} : \widetilde{\mathfrak{R}} \rightarrow \mathbb{R}$, by setting 
		\begin{flalign}
			\label{g20} 
			\widetilde{G} (t, y) = \alpha^{-1} \cdot G (\alpha^{1/2} t + t_0, \alpha^{3/2} y + y_0) - a - bt, 
		\end{flalign} 
	
		\noindent where $a$ and $b$ are defined by  
		\begin{flalign}
			\label{abg}  
			a = (2\alpha)^{-1} \big( \gamma ( t_0 + \alpha^{1/2} ) + \gamma ( t_0 - \alpha^{1/2} ) \big); \qquad b = (2\alpha)^{-1} \big( \gamma ( t_0 + \alpha^{1/2} ) - \gamma ( t_0 - \alpha^{1/2} ) \big).
		\end{flalign}
	
		\noindent Observe that  by our construction the rescaling $(t,y)\mapsto (\alpha^{1/2} t + t_0, \alpha^{3/2} y + y_0)$ maps $\mathfrak{R}'$ to $\mathfrak{R} \subset [3B^{-1},1-3B^{-1}]\times [0, 2B]$ (recall from \eqref{2r}). Thus $\mathfrak{R}'' \subset \mathfrak{R}' \subset \widetilde{\mathfrak{R}}$. 
	
		By \Cref{invariancesscale}, with the $(\alpha,\beta)$ in the first part there equal to $(\alpha^{1/2}, \alpha^{3/2})$ here, $\widetilde{G}$ solves \eqref{equationxtd} on $\widetilde{\mathfrak{R}}$. By \eqref{g20}, we have 
		\begin{align}\begin{split}\label{e:1der}
				&\partial_t  \widetilde G{(0,0)}=\alpha^{-1/2} \cdot \partial_t G(t_0,y_0)-b; \qquad -\partial_y \widetilde G{(0,0)}=-\alpha^{1/2} \cdot \partial_y G(t_0,y_0),
		\end{split} \end{align}
		and
		\begin{align}\begin{split}\label{e:2der}
				&\partial_t \partial_y \widetilde G{(0,0)}=\alpha \cdot \partial_y \partial_t G(t_0,y_0); \qquad \partial_t^2 \widetilde G{(0,0)}=\partial_t^2 G(t_0,y_0); \qquad \partial_y^2 \widetilde G{(0,0)}=\alpha^2 \cdot \partial_y^2 G(t_0,y_0).
		\end{split}\end{align}
		
		To estimate the derivatives of $u$ and $\varrho$ (as in \eqref{e:holder}), we will first estimate $\widetilde{G}$ and its derivatives, and then use \eqref{e:1der} and \eqref{e:2der} to deduce regularity bounds on $u$ and $\varrho$. To this end, let us show that $\widetilde{G}$ and its $y$-derivative are bounded on $\mathfrak{R}'$ and use \Cref{derivativef}. To do this, observe that 
		\begin{align}\begin{split}\label{e:dyGasymp}
			-\partial_y\widetilde G(t,y)&=-\alpha^{1/2} \cdot \partial_y G(\alpha^{1/2}t+t_0, \alpha^{3/2}y+y_0)\\
			&\leq D\alpha^{1/2}(\alpha^{3/2}y+y_0)^{-1/3}\leq 2^{1/3} D\alpha^{1/2}y_0^{-1/3} \le D;\\
			-\partial_y\widetilde G(t,y)& = \alpha^{1/2} \cdot \partial_y G(\alpha^{1/2} t + t_0, \alpha^{3/2} y + y_0) \\
			& \geq \alpha^{1/2} D^{-1} (\alpha^{3/2}y+y_0)^{-1/3}\geq 2^{1/3} 3^{-1/3} D^{-1} \alpha^{1/2} y_0^{-1/3} \geq (3^{1/3} B^2 D)^{-1},
		\end{split}\end{align}
	
		\noindent where we used \eqref{g20} for the first statements of both inequalities; \eqref{yd0} for the second; the fact that $y_0/2\leq \alpha^{3/2}y+y_0\leq 3y_0/2$  (from \eqref{e:ayrelation} and the fact that $(t, y) \in [-1, 1]^2$) for the third; and  $(y_0 / 2)^{2/3} B^{-4}  \le \alpha \le (y_0 / 2)^{2/3}$ from \eqref{e:ayrelation} for the fourth. It follows that there exists a constant $\varepsilon = \varepsilon (B) > 0$ such that $\widetilde{G} \in \Adm_{\varepsilon} (\mathfrak{R}')$ (recall \Cref{functionadmissible}). 
	
		Next we bound $\big| \widetilde{G} (t, y) \big|$ on $\partial \mathfrak{R}'$. For any $(t, y) \in \overline{\mathfrak{R}}'$, we have 
		\begin{align}
			\label{e:gytboundary0}
			\begin{aligned} 
			\Big|  \widetilde G(t,y)-\big(\alpha^{-1} \cdot \gamma(\alpha^{1/2}t+t_0)-a-bt \big) \Big|
				&=\alpha^{-1} \big|G(\alpha^{3/2}y+y_0, \alpha^{1/2}t+t_0)-\gamma(\alpha^{1/2}t+t_0) \big|\\
				&\leq \displaystyle\frac{3D}{2\alpha} \cdot (\alpha^{3/2}y+y_0)^{2/3} \\
				&  \leq \displaystyle\frac{3D}{2\alpha} \cdot \Big( \displaystyle\frac{3y_0}{2} \Big)^{2/3} \leq \displaystyle\frac{3D}{2} \cdot 4B^3 =6B^3D,
				 \end{aligned} 
		\end{align}
	
		\noindent where we used \eqref{g20} for the first statement; \eqref{e:gammaGdiff} for the second; the fact that $\alpha^{3/2}y+y_0\leq 3y_0/2$ from \eqref{e:ayrelation} for the third; and the last statement of \eqref{e:ayrelation} for the fourth. Moreover, the definitions \eqref{abg} of $(a, b)$ and \eqref{e:gammaconv2} (with the $(t_0,\tau,s)$ there given by $(t_0, \alpha^{1/2}, \alpha^{1/2}t)$ here) together imply for $(t, y) \in \overline{\mathfrak{R}}'$ that
		\begin{align*}
				\begin{aligned} 
				\big| & \alpha^{-1} \cdot \gamma(\alpha^{1/2}t+t_0)-a-bt \big|\\
				&=\alpha^{-1} \bigg|\gamma(\alpha^{1/2}t+t_0)- \Big( \displaystyle\frac{1-t}{2}  \Big) \cdot \gamma (t_0- \alpha^{1/2} )- \Big( \displaystyle\frac{t+1}{2} \Big) \cdot \gamma (t_0+ \alpha^{1/2} ) \bigg|\leq D.
		\end{aligned} 
	\end{align*}
		\noindent This, with \eqref{e:gytboundary0}, implies that $\big\| \widetilde G \big\|_{\mathcal{C}^0 (\mathfrak{R}')} \le (6B^3 + 1 )D$. Together with \Cref{derivativef} and the fact that $\widetilde{G}\in \Adm_{\varepsilon} (\mathfrak{R})$, this yields a constant $M = M(B) > 1$ such that $\big\| \widetilde{G} \big\|_{\mathcal{C}^2 (\mathfrak{R}'')} \le M$. 
		
		From \Cref{i:concavecopy} above, we have $\big|\gamma(t) - \gamma (t') \big| \le D|t-t'|$, for any $t, t' \in [4B^{-1}, 1-4B^{-1}]$. Together with \eqref{abg}, this implies that $|b|\leq D/(2 \alpha^{1/2})$, which with \eqref{e:1der}, the bound $\big\| \widetilde{G} \big\|_{\mathcal{C}^1 (\mathfrak{R}'')} \le \big\| \widetilde{G} \big\|_{\mathcal{C}^2 (\mathfrak{R}'')} \le M$, and \eqref{e:dyGasymp} yields 
		\begin{align}\label{e:delg1}
			\big| \partial_t G(t_0,y_0) \big|\leq \alpha^{1/2} \big(M+|b| \big)\leq M \alpha^{1/2} +D; \qquad  (2DB^2)^{-1} \leq -\alpha^{1/2}\partial_y G(t_0,y_0)\leq D.
		\end{align}
			
		\noindent Moreover, \eqref{e:2der} (with the bound $\big\| \widetilde{G} \big\|_{\mathcal{C}^2 (\mathfrak{R}'')} \le M$) implies
		\begin{align}\label{e:delg2}
			\big| \partial^2_t G(t_0,y_0) \big| \leq M,\quad \big| \partial_t\partial_y G(t_0,y_0) \big|\leq M \alpha^{-1},\quad \big| \partial_y^2 G(t_0,y_0) \big|\leq M \alpha^{-2}.
		\end{align}
	
		\noindent Thus, there exists a  constant $C_1 = C_1 (B) > 1$ such that
		\begin{flalign}
			\label{derivativeu} 
			\begin{aligned} 
			\big| \partial_x u( t_0, x_0) \big| & = \big| \partial_y G(t_0, y_0) \big|^{-1} \cdot \big| \partial_y \partial_t G(t_0, y_0) \big| \\
			& \le 2DB^2 \alpha^{1/2} \cdot M \alpha^{-1} = 2B^2 D M \alpha^{-1/2} \le C_1 \big( \gamma(t_0) - x_0 \big)^{-1/2},
			\end{aligned} 
		\end{flalign}
	
		\noindent where in the first statement we used the definition \eqref{u0} of $u$ and the fact that $x_0 = G(t_0, y_0)$; in the second we used \eqref{e:delg1} and \eqref{e:delg2}; in the third we evaluated the product; and in the fourth we used the bound  $\alpha^{1/2}\geq B^{-2} (y_0/2)^{1/3}\geq B^{-2} \big( (\gamma(t_0)-y_0)/3D \big)^{1/2}$, which holds by the last statement in \eqref{e:ayrelation} and \eqref{e:gammaGdiff}. Similarly, there exists a constant $C_2 = C_2 (B) > 1$ such that
		\begin{flalign}
			\label{derivativerho} 
			\begin{aligned}  
			 \big| \partial_x \varrho(t_0, x_0) \big| & = \big| \partial_y G(t_0, y_0) \big|^{-1} \cdot \Big| \partial_y \Big( -\displaystyle\frac{1}{\partial_y G (t_0, y_0)} \Big) \bigg| \\
			 & =  |\partial_y G(t_0,y_0) \big|^{-3} \cdot \big| \partial_y^2 G(t_0,y_0) \big| \\
			& \leq (2DB^2)^3 \alpha^{3/2} \cdot M \alpha^{-2} = 8 B^6 D^3 M \alpha^{-1/2} \leq C_2 \big( \gamma(t_0)-x_0 \big)^{-1/2},
			\end{aligned} 
		\end{flalign} 
	
		 \noindent where in the first statement we used \eqref{urho20} and the fact that $x_0 = G(t_0, y_0)$; in the second we performed the differentiation; in the third we used \eqref{e:delg1} and \eqref{e:delg2}; in the fourth we evaluated the product; and in the fifth we again used the bound $\alpha^{1/2}\geq B^{-2} (y_0/2)^{1/3}\geq B^{-2} \big((\gamma(t_0)-y_0)/3D \big)^{1/2}$. Together, \eqref{derivativeu} and \eqref{derivativerho} verify the first statement \eqref{e:holder} of the proposition.

	Finally we show that $\varrho(t,x)$ and $u(t,x)$ extend continuously to the upper boundary $\big\{ \gamma(t) \big\}$ of $\Omega$. The continuity of $\varrho(t,x)$, and that it converges to $0$ as $(t, x)$ tends to $\big( t, \gamma(t) \big)$, follows from \eqref{e:densitybound}. To show the continuity of $u(t,x)$, it suffices by \eqref{u0} to show that $\partial_t G(t,y)$ extends continuously to the set $(t, y) \in [4B^{-1}, 1-4B^{-1}] \times \{ y=0 \}$. For $y \in (0, B)$, \eqref{e:delg2} gives $\big| \partial_t^2 G(t,y) \big|\leq M$ and $\big|\partial_t \partial_y G(t,y) \big|\leq M B^3 (2/y)^{2/3}$ (where in the latter we used the fact that $\alpha \ge B^{-3} ( y_0 / 2 )^{2/3}$ from the last statement in \eqref{e:ayrelation}). Hence, for any $t, t' \in [4B^{-1}, 4B]$ and $y, y' \in (0, B)$ with $t' < t$ and $y' < y$, we have
	\begin{align}
		\label{gtyy} 
		\begin{aligned}
	&\phantom{{}={}}\big|\partial_t G(t,y)-\partial_tG(t',y') \big|  \le \big| \partial_t G(t, y) - \partial_t G(t, y') \big| + \big| \partial_t G(t, y') - \partial_t G (t', y') \big|  \\
	 &\leq \int_{y'}^y \big| \partial_t \partial_r G(t,r) \big| d r + \displaystyle\int_{t'}^t \big| \partial_s^2 G(s, y') \big| ds \\
	& \leq 2^{2/3} MB^3 \int_{y'}^y r^{-2/3} dr + M (t-t') \leq 6MB^3 \big(y^{1/3} - (y')^{1/3} \big) + M (t-t'),
	\end{aligned} 
	\end{align}

	\noindent and so $\partial_t G$ is uniformly continuous on $[4B^{-1}, 1-4B^{-1}] \times (0, B)$. In particular, for $t \in [4B^{-1}, 1-4B^{-1}]$, the function $u(t,x)$ extends uniformly continuously to the north boundary $\big\{ \gamma(t) \big\}$ of $\Omega$, so  
	\begin{flalign*} 
		\gamma' (t) &= \displaystyle\lim_{t' \rightarrow t} \displaystyle\frac{G(t, 0) - G(t', 0)}{t-t'}  = \displaystyle\lim_{t' \rightarrow t} \Bigg( \displaystyle\lim_{y \rightarrow 0^+} \displaystyle\frac{G(t, y) - G(t', y)}{t-t'} \Bigg) \\
		& = \displaystyle\lim_{y \rightarrow 0^+} \Bigg(\displaystyle\lim_{t'\rightarrow t} \displaystyle\frac{G(t, y) - G(t', y)}{t-t'} \Bigg)  = \displaystyle\lim_{y \rightarrow 0^+} \partial_t G(t, y) = \displaystyle\lim_{x \rightarrow 0^+} u \big( t, \gamma(t) \big),
	\end{flalign*} 

	\noindent where the first statement is by the definition of $\gamma$; the second is by the continuity of $G(t, y)$ around $y=0$ (by \Cref{hrhot}); the third is by the uniformity of the convergence of the right side of \eqref{gtyy} to $0$, if $y = y'$ tends to $0$; the fourth is from the fact that $G$ is smooth on $\Omega^{\inv} = (0, 1) \times (0, L^{3/2})$ (recall \Cref{urhoderivatives0} and \Cref{i:densitylow} in \Cref{p:densitylow}); and the fifth is by \eqref{u0}. Together with the uniform continuity \eqref{gtyy} of $\partial_t G$, this implies that $\gamma'$ is continuous on $[4B^{-1}, 1-4B^{-1}]$.		
	\end{proof}

	\section{Limit Shapes on Tall Rectangles}\label{s:shape1}

	In this section we study the inverted height function associated with a bridge-limiting measure process (as in \Cref{mutmu0mu1}) on a tall $1 \times L^{3/2}$ rectangle. We show that, under \Cref{integralmu0mu1} and \Cref{gapmu0mu1}, its inverted height function (recall \Cref{hrhot}) behaves around its arctic boundary approximately as does the one \eqref{e:Airyprofile} associated with the limiting Airy profile, with coefficients $(\mathfrak{a}, \mathfrak{b}, \mathfrak{c})$ bounded above and below, independently of $L$. Throughout this section, we adopt and recall the notation from \Cref{lmu} and recall from \Cref{gequation} that the inverted height function $G$ satisfies the elliptic equation \eqref{equationxtd}. We also recall the sets $\Adm$ and $\Adm_{\varepsilon}$ of admissible functions from \Cref{functionsadmissible0} and \Cref{functionadmissible}, respectively.

	\subsection{Complex Burgers Equation and Characteristic Maps}
	
	\label{Equation0} 
	
	The following theorem, to be established in \Cref{FunctionEdge} below, indicates that the inverted height function $G$ from \Cref{lmu} (under \Cref{integralmu0mu1} and \Cref{gapmu0mu1}) behaves approximately as a limiting Airy profile (recall \eqref{e:Airyprofile}) near the edge of its support.
	
	\begin{thr}\label{p:limitprofile}
		
		Adopting \Cref{integralmu0mu1} and \Cref{gapmu0mu1}, there exist constants $c = c(B) > 0$ and $C = C(B) > 1$ such that the following holds if $L \ge C$. For any $\ft\in [5B^{-1},1-5B^{-1}]$, there are real numbers $\mathfrak{a}, \mathfrak{b} \in [-C, C]$ and $\mathfrak{c} \in [C^{-1}, C]$ satisfying the below property. For any real numbers $\tau \in [-c, c]$ and $y \in [0, c]$, we have 
		\begin{align}\label{e:limitprofile}
			\Bigg|G(\ft+\tau,y)-\bigg(\fa+\fb\tau-\fc\tau^2-\Big(\frac{3\pi}{4\fc^{1/2}}\Big)^{2/3}y^{2/3} \bigg) \Bigg| \leq C \big( |\tau|^3+ |\tau| y^{2/3}+y \big).
		\end{align}
	\end{thr}

	The proof of this theorem will make considerable use of the bounds from \Cref{GDerivativeConvexgamma}, as well as the complex Burgers equation \Cref{p:solution}. In this section we state some results and properties about the latter; throughout, we adopt \Cref{integralmu0mu1} and \Cref{gapmu0mu1} (and hence we recall the notation from \Cref{lmu}). First, observe by \Cref{p:solution} that complex slope $f$ satisfies \eqref{ftfx}, which can be rewritten as 
	\begin{align}\label{e:burger}
		-\partial_t f(t,x)=f(t,x) \cdot \partial_x f(t,x) = \displaystyle\frac{1}{2} \cdot \partial_x  \big( f(t,x)^2 \big).
	\end{align}

	\noindent We further recall from \Cref{p:densityrg} that, for sufficiently large $L$, the complex slope $f$ extends continuously to the part of the arctic boundary given by $\big\{ (t, \gamma(t) ) \in \mathbb{R}^2 : t \in [4B^{-1}, 1-4B^{-1}] \big\}$. Moreover, that proposition also implies for each $t\in [4B^{-1},1-4B^{-1}]$ that   
	\begin{flalign}
		\label{ufrhogammat}
		\lim_{x\rightarrow \gamma(t)} \varrho_t(x)=0, \quad \text{and} \quad \lim_{x\rightarrow \gamma(t)} u_t(x)=\gamma'(t), \qquad \text{so that} \quad f \big( t, \gamma(t) \big) = \gamma'(t).
	\end{flalign}

	A function that will be useful to analyze the complex Burgers equation will be the following characteristic map, which will later provide a complex coordinate on $\Omega$ (see \Cref{p:analytic} below).  
	
	\begin{definition} 
		
		\label{ztx2}
		
		Adopt \Cref{lmu}; fix a real number $t_0 \in (0, 1)$; and define the set $\Omega (t_0) = \big\{ (t, x) \in \Omega : t \ge t_0 \big\}$. Define the \emph{characteristic map}\index{Z@$z_{t_0}$; characteristic map} $z = z_{t_0} : \Omega (t_0) \rightarrow \mathbb{H}^-$ by for each $(t, x) \in \Omega (t_0)$ setting
		\begin{align}\label{e:zxt}
			z(t,x)=x-(t-t_0) \cdot f(t, x)=x-(t - t_0) \cdot u_t(x) - \pi\mathrm{i} \cdot (t-t_0)  \varrho_t(x),
		\end{align}

	\end{definition} 
	
	\begin{rem} 
		
	\label{zb4}
	
		 By \Cref{p:densityrg}, $z$ extends continuously to $\big( [4B^{-1}, 1 - 4B^{-1}] \times \mathbb{R} \big) \cap \overline{\Omega (t_0)}$ (containing part of the arctic boundary) if we adopt \Cref{integralmu0mu1} and \Cref{gapmu0mu1}. The same proposition implies that $z \big( t, \gamma(t) \big) = \gamma(t) - (t-t_0) \cdot \gamma'(t) \in \mathbb{R}$ for each $t \in [4B^{-1}, 1 - 4B^{-1}]$, if $L$ is sufficiently large.

	\end{rem} 

	The next lemma provides some general properties of the characteristic map. In what follows, for any subset $U \subset \mathbb{R}^2$, a differentiable function $g : U \rightarrow \mathbb{C}$ is called \emph{positively oriented} if the map $(\Real g, \Imaginary g) : U \rightarrow \mathbb{R}^2$ has a nonnegative Jacobian determinant everywhere in the interior of $U$. It is \emph{strictly positively oriented} at a point $u \in U$ if this Jacobian determinant is positive at $u$. Any point $u \in U$ at which this Jacobian determinant is equal to $0$ is called a \emph{critical point} for $g$.

	\begin{lem}\label{p:zprop}
		
		Adopting \Cref{integralmu0mu1} and \Cref{gapmu0mu1}, there exists a constant $C = C(B) > 1$ such that the following holds if $L \ge C$. Fix a real number $t_0 \in [4B^{-1}, 1 - 4B^{-1}]$, and let $z = z_{t_0}$ be the characteristic map as in \Cref{ztx2}. Then $f$ is real analytic on $\Omega$, and $z$ is real analytic, positively oriented on $\Omega(t_0)$, and strictly positively oriented away from its critical points. Moreover, if some point point $(t, x) \in \Omega (t_0)$ is a critical point of $z(t,x)$, then $\partial_t z(t,x)=\partial_x z(t,x)=0$. Additionally, for any $(t, x) \in \Omega (t_0)$, we have
	\begin{align}\label{e:dtzdxz}
		\partial_t z(t,x)=-f(t, x) \cdot \partial_x z(t,x).
	\end{align}
	\end{lem}
\begin{proof} 
	
	 The definition \eqref{frhou} of $f$ and \Cref{gderivatives0} together imply that $f$ is real analytic on $\Omega$; by \eqref{e:zxt}, it follows that $z$ is real analytic on $\Omega (t_0)$. Next, by \Cref{ztx2} and \eqref{e:burger}, we have for any $(t, x) \in \Omega (t_0)$ that 
	\begin{align*}
	&\partial_t z(t,x)=-f(t, x)-(t-t_0)\partial_t f(t, x)=-f(t, x)+(t-t_0)f(t, x) \cdot \partial_x f(t, x); \\
	&\partial_x z(t,x)=1-(t-t_0) \partial_x f(t, x).
	\end{align*}
	Thus, $\partial_t z(t,x)=-f(t, x) \cdot \partial_x z(t,x)$, which verifies \eqref{e:dtzdxz}. Moreover, it implies that the determinant of the Jacobian of  $z(t,x)$ is given by 
	\begin{align*}
	\det & \left[
	\begin{array}{cc}
	\Real \partial_t z(t,x)  &\Imaginary \partial_t z(t,x)\\
	\Real \partial_x z(t,x) &\Imaginary \partial_x z(t,x)
	\end{array}
	\right] \\
	& \qquad = \det \left[ \begin{array}{cc} \Imaginary f \cdot \Imaginary \partial_x z - \Real f \cdot \Real \partial_x z & -\Imaginary f \Real \partial_x z - \Real f \cdot \Imaginary \partial_x z \\ \Real \partial_x z & \Imaginary \partial_x z \end{array} \right] \\
	 & \qquad = \Imaginary f(t, x) \cdot \big| \partial_x z(t,x) \big|^2=\pi \varrho(t,x) \cdot \big|\partial_x z(t,x) \big|^2,
	\end{align*}

	\noindent where in the last equality we applied \eqref{frhou}; this implies that $z$ is positively oriented and strictly positively oriented at $(t, z)$ unless $\partial_x z(t, x) = 0$. Since $\varrho(t, x) > 0$ for $(t, x) \in \Omega (t_0) \subseteq \Omega$, it follows that $(t, x)$ is a critical point of $z$ if and only if $\partial_xz(t,x)=\partial_tz(t,x)=0$, establishing the lemma.
\end{proof}

	\subsection{Heuristic for \Cref{p:limitprofile}}
	
	\label{Profile0}
	
	In this section we provide a heuristic for \Cref{p:limitprofile}, before giving its detailed proof. We will be quite informal here, assuming that all relevant limits exist, functions are smooth, and derivatives are uniformly bounded (the last of which will actually constitute the main effort, and is where we require the estimates proven in \Cref{s:density}).

	To that end, let us fix a real number $t_0$ close to $\mathfrak{t}$; recall the characteristic map $z = z_{t_0}$ defined by \eqref{e:zxt}. For $t \ge t_0$ close to $\mathfrak{t}$, define 
	\begin{flalign}
		\label{xitgamma0} 
		\xi(t) = \gamma(t) - (t-t_0) \cdot \gamma' (t) = z \big( t, \gamma(t) \big),
	\end{flalign}
	
	\noindent where the last equality formally follows from \eqref{e:zxt} and \eqref{ufrhogammat}. Observe that $\xi(t) \in \mathbb{R}$; fix a time $s_0 = \mathfrak{t} + \tau$. Since $G(s_0, 0) = \gamma(s_0)$, for \eqref{e:limitprofile} to hold at $y = 0$, we must have 
	\begin{flalign*} 
		\mathfrak{a} = \gamma(\mathfrak{t}); \qquad \mathfrak{b} = \gamma' (\mathfrak{t}); \qquad  \mathfrak{c} = - \displaystyle\frac{\gamma'' (\mathfrak{t})}{2},
	\end{flalign*} 

	\noindent  where we have implicitly assumed that $\gamma''$ exists. This would guarantee \eqref{e:limitprofile} at $y=0$, by a Taylor expansion of $\gamma$. It therefore suffices, for $y > 0$, to approximate 
	\begin{flalign}
		\label{ygs0y} 
		\partial_y G(s_0, y) \approx -\pi^{2/3} \cdot \big( 3 y\gamma''(\mathfrak{t}) \big)^{-1/3}.
	\end{flalign} 

	\noindent As $\partial_y G(s_0, y) = -\varrho \big( s, G(s_0, y) \big)^{-1}$, it can be seen by direct computation that this is equivalent to 
\begin{flalign}
\label{rhos00} 
\varrho (s_0, x) \approx \pi^{-1} \cdot \Big( -2 \gamma'' (\mathfrak{t}) \cdot \big( \gamma (s_0) - x \big) \Big)^{1/2},
\end{flalign}

\noindent for $x = G(s_0, y) < \gamma (s_0)$ close to $\gamma(s_0)$. Indeed, by integration, \eqref{rhos00} would give (see \Cref{gammag}) 
\begin{flalign*} 
	y = H(s_0, x) \approx \bigg( -\displaystyle\frac{8 \gamma'' (\mathfrak{t})}{9 \pi^2} \bigg)^{1/2} \cdot \big( \gamma(s_0) - x \big)^{3/2},
\end{flalign*} 

\noindent which with \eqref{rhos00} yields \eqref{ygs0y}.

To show \eqref{rhos00} (see \Cref{rhoxtgamma} below), observe by \eqref{e:zxt} that $\rho_{s_0} (x) = -\pi^{-1} \cdot (s_0 - t_0)^{-1} \cdot \Imaginary z(s_0, x)$. To evaluate $z(s_0, x)$, we will approximate it by $z \big( s_0, \gamma(s_0) \big) = \xi(s_0)$. In particular, since $\xi (s_0) \in \mathbb{R}$, denoting $w = z(s_0, x) - \xi (s_0)$, we have that 
\begin{flalign}
\label{rhos0w} 
\rho_{s_0} (x) = -\pi^{-1} \cdot (s_0 - t_0)^{-1} \cdot \Imaginary w.  
\end{flalign} 

\noindent We will show as \Cref{p:analytic} below that $z$ is a homeomorphism from a neighborhood in $\overline{\Omega}$ of the point $\big( \mathfrak{t}, \gamma(\mathfrak{t}) \big)$ to a neighborhood in the upper half-plane $\overline{\mathbb{H}}$ of the point $\xi(t)$. Then, define the function $F$ on its image by setting $F \big( z(t,x) \big) = f(t,x)$. This yields   
\begin{flalign}
\label{xfzs0x}
x = z (s_0, x) + (s_0 - t_0) \cdot f(s_0, x) = z(s_0, x) + (s_0 - t_0) \cdot F \big( z(s_0, x) \big),
\end{flalign} 

\noindent where the first statement follows from the definition \eqref{e:zxt} of $z$, and the second from that of $F$.

Now, a basic fact about $F$ that will also be shown in \Cref{p:analytic} (as a consequence of the complex Burgers equation \Cref{p:zprop}) is that it extends to a smooth (in fact, holomorphic) function on the boundary $\partial \Omega$, even though $f$ and $z$ do not (as their derivatives diverge there). Therefore, Taylor expanding $F$ in \eqref{xfzs0x} (assuming its derivatives are uniformly bounded), we obtain (recalling that $z(s_0, x) = \xi(s_0) + w$) that
\begin{flalign*}
x - \xi(s_0) \approx (s_0 - t_0) \cdot F \big( \xi (s_0) \big) + w + (s_0 - t_0) w \cdot F' \big( \xi(s_0) \big) + (s_0 - t_0) \cdot \displaystyle\frac{w^2}{2} \cdot F'' \big( \xi(s_0) \big).
\end{flalign*}

\noindent Applying \eqref{xfzs0x} at $w=0$, meaning that $x = \gamma(s_0)$, gives $(s_0-t_0) \cdot F \big( \xi(s_0) \big) = \gamma(s_0) - \xi(s_0)$. Inserting this into the above approximation yields 
\begin{flalign}
\label{wfxis0}
w + (s_0-t_0) w \cdot F' \big( \xi(s_0) \big) + (s_0-t_0) \cdot \displaystyle\frac{w^2}{2} \cdot F'' \big( \xi(s_0) \big) = x - \gamma(s_0).
\end{flalign}

\noindent Next let us evaluate $F' \big( \xi(s_0) \big)$ and $F'' \big( \xi(s_0) \big)$ (see \Cref{gamma123} below). For the former, we have
\begin{flalign*} 
F \big( \xi(s_0) \big) = f \big( s_0, \gamma(s_0) \big) = \gamma' (s_0), 
\end{flalign*}

\noindent where the first statement follows from the definition of $F$, and the second formally follows from \eqref{ufrhogammat}. Differentiating in $s_0$, and using the fact that $\xi' (s_0) = (t_0 - s_0) \cdot \gamma''(s_0)$ by \eqref{xitgamma0}, we deduce
\begin{flalign}
\label{f1xis0}
F' \big( \xi(s_0) \big) = \xi'(s_0)^{-1} \cdot \gamma'' (s_0) = (t_0 - s_0)^{-1}.
\end{flalign}

\noindent Differentiating again in $s_0$, and using similar reasoning, we obtain
\begin{flalign}
\label{f2xis0}
F'' \big( \xi(s_0) \big) = -\xi'(s_0)^{-1} \cdot (t_0 - s_0)^{-2} = \gamma''(s_0)^{-1} \cdot (s_0 - t_0)^{-3}.
\end{flalign}

\noindent Inserting \eqref{f1xis0} into \eqref{wfxis0}, we find that the coefficient of the linear term $w$ on the right side is equal to $0$. By \eqref{f2xis0}, this therefore yields
\begin{flalign*}
w \approx  (s_0-t_0) \cdot \Big( - 2\gamma'' (s_0) \cdot \big( \gamma(s_0) - x \big) \Big)^{1/2},
\end{flalign*} 

\noindent which using \eqref{rhos0w} implies \eqref{rhos00}.

	\subsection{Holomorphic Coordinate}		
	
	\label{SetCharacteristic}

\begin{figure}
	\center
\includegraphics[width=1\textwidth, trim=0 1cm 0 1cm, clip]{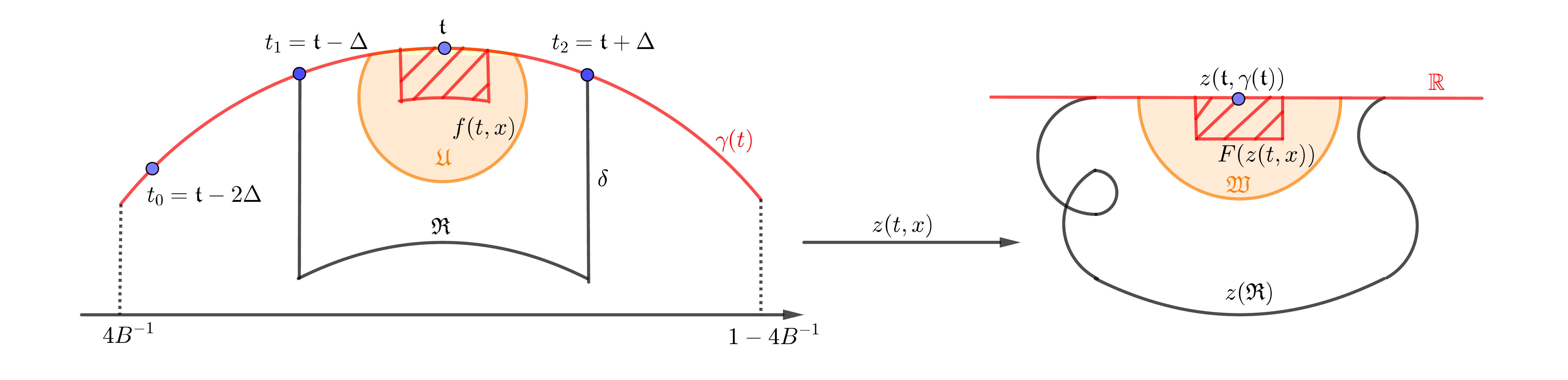}

\caption{ Shown above is a depiction of \Cref{p:analytic}, whose first part indicates that $z$ is a homeomorphism from $\mathfrak{U}$ to $\mathfrak{W}$ (the orange regions) and whose third part indicates that $\mathfrak{U}$ is not too small (contains the red region).}

\label{f:ztx}
	\end{figure}

	The next proposition, to be established in \Cref{SetCharacteristic} below, states that the characteristic map \eqref{e:zxt} is a bijection, at least on an open set in its domain intersecting the arctic boundary; see \Cref{f:ztx} for a depiction. Using this, it defines a holomorphic coordinate $F$ with uniformly bounded derivatives, which (as explained in \Cref{Profile0}) will be used in the proof of \Cref{p:limitprofile}.

	\begin{prop}\label{p:analytic}
		Adopting \Cref{integralmu0mu1} and \Cref{gapmu0mu1}, there exist constants $c_1 = c_1 (B) \in (0, 1)$, $c_2 = c_2 (B) \in (0, 1)$, and $C = C(B) > 1$ such that the following holds if $L \ge C$. Fix real numbers $\mathfrak{t} \in [5B^{-1}, 1 - 5B^{-1}]$ and $\Delta \in ( 0, 1 / 4B]$; set $t_0=\ft-2\Delta$, and let $z = z_{t_0}$ be the characteristic map as in \Cref{ztx2}.
		  \begin{enumerate}
		  \item \label{i:injection}  There exists a neighborhood $\mathfrak{U} \subseteq \Omega$ of $\big( \mathfrak{t}, \gamma (\mathfrak{t}) \big)$ such that the following two statements hold. First, $\mathfrak{U} \subseteq [\mathfrak{t} - \Delta, \mathfrak{t} + \Delta] \times \mathbb{R}$. Second, $z$ is a homeomorphism from $\mathfrak{U}$ to the set $\mathfrak{W} = \big\{ w \in \mathbb{H}^- : \big| w - z(\mathfrak{t}, \gamma(\mathfrak{t})) \big| < 2 c_1 \Delta^2 \big\}$.
		
		\item \label{i:analytic} Define $F : \mathfrak{W} \rightarrow \mathbb{H}^-$ by setting $F \big( z(t, x) \big) = f(t, x)$, for each $(t, x) \in \mathfrak{U}$. Then $F$ extends to a holomorphic function on the set $\big\{ w\in \mathbb{C}: \big| w -z(\ft, \gamma(\ft))|\le c_1 \Delta^2 \big\}$. We have $F (\overline{z}) = \overline{f(z)}$ and, for any integer $k \ge 0$, we have 
		\begin{align}\label{e:fkbound}
		\big| \partial_w^k F (w) \big|\leq \frac{C}{ k! c_1^{k}\Delta^{2k}}.
		\end{align}
		\item \label{i:injectiondomain}The characteristic map $z$ is an injection from 
		\begin{flalign*}
		\big\{(t,x)\in \Omega: |t-\ft|\leq c_2 \Delta, \gamma(t) \ge x\geq \gamma(t)-c_2 \Delta^2 \big\} \quad \text{into} \quad \bigg\{z\in  \mathbb{H}^-: \Big|z- z \big(\ft, \gamma(\ft) \big) \Big|\le c_1 \Delta^2 \bigg\}.
		\end{flalign*}
		\end{enumerate}
	\end{prop}

\begin{figure}
	\center
\includegraphics[scale = .8, trim=0 1cm 0 1cm, clip]{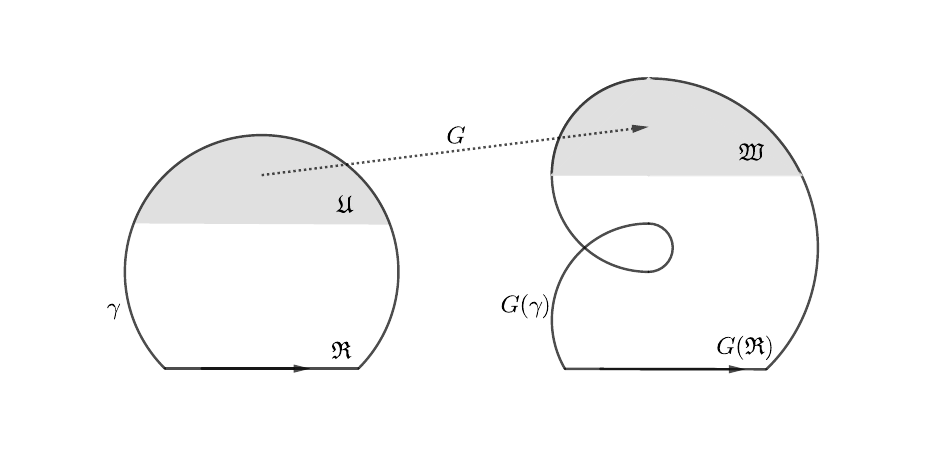}

\caption{Shown above is a depiction for the setup of \Cref{degreed}.}
\label{f:homeo}
	\end{figure}

 To establish \Cref{p:analytic}, we first require the following topological fact providing a sufficient condition for a positively oriented, real analytic map to be a homeomorphism; see \Cref{f:homeo} for a depiction. It is similar to known results (see, for example, \cite[Section 2.5]{DT}), though we have not seen it in the literature stated as written here; so, it is shown in \Cref{ProofFunctionG} below. 

\begin{prop}
	
	\label{degreed}

	Let $\mathfrak{R} \subset \mathbb{R}^2$ denote a bounded, simply-connected, open set, whose boundary $\gamma = \partial \mathfrak{R}$ is a piecewise differentiable Jordan curve. Let $G : \overline{\mathfrak{R}} \rightarrow \mathbb{R}^2$ denote a nonconstant, real analytic function that is strictly positively oriented away from its critical points, and let $\mathfrak{W} \subset \mathbb{R}^2$ denote a connected, bounded, open set, satisfying the following four properties.
	\begin{enumerate}
		\item The set $\mathfrak{W} \cap G(\mathfrak{R})$ is nonempty.
		\item The set $\mathfrak{W}$ is disjoint from the curve $G(\gamma)$. 
		\item The winding number of $G(\gamma)$, with respect to any point $w \in \mathfrak{W}$, is equal to one. 
		\item For any point $w \in \overline{\mathfrak{W}} \cap G(\gamma)$, there is only one point $u \in \overline{\mathfrak{R}}$ such that $G(u) = w$.
	\end{enumerate}
	
	\noindent Let $\mathfrak{U} = G^{-1} (\mathfrak{W}) \subseteq \mathfrak{R}$. Then, $G$ is a homeomorphism from $\mathfrak{U}$ to $\mathfrak{W}$. 
	
\end{prop}

In the remainder of this section, we adopt the notation and assumptions of \Cref{p:analytic}. Observe from \Cref{i:dyG} of \Cref{p:checka} (with the fact that $\gamma(t) = G(t, 0)$) that there exist constants $c_3 = c_3 (B) \in (0, 1)$ and $c_4 = c_4 (B) \in (0, 1)$ such that 
\begin{flalign}
	\label{gammagtb} 
	c_4 \le \displaystyle\frac{3c_3}{2} \cdot B^{2/3} = c_3 \displaystyle\int_0^B y^{-1/3} dy \le \gamma(t) - G(t, B) \le 4B \displaystyle\int_0^B y^{-1/3} dy = 6 B^{5/3} \le c_4^{-1},
\end{flalign}

\noindent for any $t \in [3B^{-1}, 1 - 3B^{-1}]$. Next observe from \Cref{p:checka}, \Cref{p:densityrg}, and \eqref{e:densitybound} that there exists a constant $M=M(B) > 576 B^2$ such that the following three statements hold, for any real numbers $t\in [4B^{-1}, 1-4B^{-1}]$ and $x \in \big[ G(t, B), \gamma(t) \big]$.

 First, we have $\big| \gamma' (t) \big| \le M / 8$ and, for any $\tau \in \mathbb{R}$ such that $3B^{-1} \le t - \tau \le t + \tau \le 1-3B^{-1}$, we have (from the $s=0$ case of \eqref{e:gammaconv}) that
		\begin{align}\label{e:strongconvex}
			 M^{-1} \tau^2 \le \gamma(t)- \bigg( \displaystyle\frac{\gamma(t-\tau) + \gamma(t+\tau)}{2} \bigg) \le M\tau^2.
		\end{align}

\noindent Second, by integrating \eqref{e:holder} (and replacing the constant $C$ there by $M^{1/2} / 2$ here), we have the H\"{o}lder bounds 
		\begin{align}\begin{split}\label{e:rhouholder}
		& \Big|\varrho_{t}(x)-\varrho_{t}\big( \gamma(t) \big) \Big| \leq M^{1/2} \big( \gamma(t)-x \big)^{1/2}; \qquad \Big|u_{t}(x)-u_{t} \big(\gamma(t) \big) \Big| \leq M^{1/2} \big(\gamma(t)-x \big)^{1/2}.
		\end{split}\end{align}
		
		\noindent In particular, since $\varrho_t \big( \gamma(t) \big) = 0$ and $\big| u_t (\gamma(t)) \big| = \big| \gamma'(t) \big| \le M / 8$ (where the first equality is from the second part of \Cref{p:densityrg}), it follows from the bound $M > 576B^2$ and the fact (from \eqref{gammagtb}) that $\gamma(t) - G(t, B) \le 6B^2$ that $\big| \varrho_t (x) \big| \le 3BM^{1/2} \le M / 8$ and $\big| u_t (x) \big| \le M / 8 + 3BM^{1/2} \le M / 2$.  With \eqref{frhou}, this gives
		\begin{flalign}
			\label{festimatetx} 
			\big| f(t, x) \big| \le \big| u_t (x) \big| + \pi \cdot \big| \varrho_t (x) \big| \le M. 
		\end{flalign} 
	\noindent Third, from \eqref{e:densitybound} and the fact that $G$ is a bijection from $\Omega^{\inv} = (0, 1) \times (0, L^{3/2})$ (recall \Cref{i:densityup} of \Cref{p:densitylow}) to $\Omega$, we have  
\begin{align}\label{e:densitylowb}
M^{-1} \big( \gamma(t)-x \big)^{1/2} \leq \varrho_t(x)\leq M \big( \gamma(t) - x \big)^{1/2}.
\end{align}

We now establish \Cref{p:analytic} using these bounds; we refer to \Cref{f:ztx} for a depiction. For the remainder of this section, we recall that $t_0 = \mathfrak{t} - 2 \Delta$, and define 
		\begin{flalign}
			\label{deltatt} 
		\delta= \Big( \displaystyle\frac{c_4 \Delta}{600 \pi M^4} \Big)^2; \qquad c_1= \displaystyle\frac{1}{600 M^5}; \qquad t_1=\ft-\Delta; \qquad t_2=\ft+\Delta.
			\end{flalign}
		
		\noindent Observe that $4B^{-1} \le t_0 \le t_2 \le 1 - 4B^{-1}$ (as $\mathfrak{t} \in [5B^{-1}, 1 - 5B^{-1}]$ and $\Delta \le (2B)^{-1}$), and that $\gamma(t) - \delta \ge \gamma(t) - c_4 \ge G(t, B)$ by \eqref{gammagtb}. Further define  the domain     
		\begin{flalign}\label{e:defO}
			\fR = \big\{(t,x)\in \overline{\Omega}: t_1 <  t < t_2, \gamma(t)-\delta < x < \gamma(t) \big\} \subseteq \big( [4B^{-1}, 1 - 4B^{-1}] \times \mathbb{R} \big) \cap \overline{\Omega (t_0)}.
		\end{flalign}
		
		\noindent We also set notation for its boundary $\partial \mathfrak{R}$, defining 
		\begin{flalign*}
			& \partial_{\no} \mathfrak{R} = \Big\{ \big( t, \gamma(t) \big) : t \in [t_1, t_2] \Big\}; \qquad \qquad \partial_{\ea} \mathfrak{R} = \Big\{ (t_2, x) : x \in \big[ \gamma(t_2) - \delta, \gamma(t_2) \big] \Big\}; \\
			& \partial_{\so} \mathfrak{R} = \Big\{ \big( t, \gamma(t) - \delta \big) : t \in [t_1, t_2] \Big\}; \qquad \partial_{\we} \mathfrak{R} = \Big\{ (t_1, x) : x \in \big[ \gamma(t_1) - \delta, \gamma (t_1) \big] \Big\},
		\end{flalign*} 
		
		\noindent and observing that $\partial \mathfrak{R} = \partial_{\no} \mathfrak{R} \cup \partial_{\ea} \mathfrak{R} \cup \partial_{\so} \mathfrak{R} \cup \partial_{\we} \mathfrak{R}$.

		\begin{proof}[Proof of \Cref{i:injection} in \Cref{p:analytic}]
			
			We will eventually apply \Cref{degreed} to the domain $\mathfrak{R}$, to which end we must study the image of its boundary under the characteristic map $z$. We begin by analyzing the image $z(\partial_{\no} \mathfrak{R})$ of north boundary $\partial_{\no} \mathfrak{R}$ of $\mathfrak{R}$. By \Cref{zb4}, $z$ maps this curve to the subset of the real axis given by 
		\begin{align*}
			\Big\{ z \big( t, \gamma(t) \big): t \in [t_1, t_2] \Big\} = \Big\{  \big( t, \gamma(t) - (t-t_0) \cdot \gamma'(t) \big) : t \in [t_1, t_2] \Big\} \subseteq \mathbb{R}.
		\end{align*}  
		
		\noindent To analyze this subset, let $a, b \in [t_1, t_2]$ be any real numbers with $a < b$. We claim that 
		\begin{flalign} 
			\label{zgammaab} 
			z \big( b, \gamma(b) \big) \ge z \big( a, \gamma (a) \big) + M^{-1} (b-a)^2.
		\end{flalign} 
	
		\noindent To verify \eqref{zgammaab},\footnote{If $\gamma''$ were known to exist (which we will only show later, in \Cref{gamma123} below), then \eqref{zgammaab} can be deduced more quickly as follows. Remark \ref{zb4} states $z \big( t, \gamma(t) \big) = \gamma(t) - (t-t_0) \cdot \gamma' (t)$, which by differentiating in $t$ gives $\partial_t \big( z (t, \gamma(t)) \big) = (t_0 - t) \cdot \gamma'' (t)$. Integrating this over $t \in [a, b]$; using by \eqref{e:strongconvex} that $\gamma'' (t) \le -2M^{-1}$; and observing that $b - t_0 \ge b-a$ (as $t_0 \le t_1 \le a$) gives $z \big( b, \gamma (b) \big) \ge z \big( a, \gamma(a) \big) + 2 M^{-1} \int_a^b (t-t_0) dt = z \big( a, \gamma(a) \big) + M^{-1} (b-a)^2$.} observe since $\gamma$ is concave (by \eqref{e:strongconvex}) that 
		\begin{flalign*}
			\gamma (b) - \gamma \Big( \displaystyle\frac{a+b}{2} \Big) \ge \displaystyle\frac{b-a}{2} \cdot \gamma' (b),
		\end{flalign*} 
	
		\noindent which implies that 
		\begin{flalign}
			\label{gamma1b} 
			\gamma (b) - (b-t_0) \cdot \gamma' (b) \ge \gamma \Big( \displaystyle\frac{a+b}{2} \Big) - \Big( \displaystyle\frac{a+b}{2} - t_0 \Big) \cdot \displaystyle\frac{2}{b-a} \cdot \bigg( \gamma (b) - \gamma \Big( \displaystyle\frac{a+b}{2} \Big) \bigg).
		\end{flalign}
	
		\noindent By similar reasoning, we also have
		\begin{flalign}
			\label{gamma1a} 
			\gamma \Big( \displaystyle\frac{a+b}{2} \Big) - \Big( \displaystyle\frac{a+b}{2} - t_0 \Big) \cdot \displaystyle\frac{2}{b-a} \cdot \bigg( \gamma \Big( \displaystyle\frac{a+b}{2} \Big) - \gamma(a) \bigg) \ge \gamma (a) - (a - t_0) \cdot \gamma' (a).
		\end{flalign}
		
		\noindent Thus,
		\begin{align*}
			\begin{aligned}
			z \big(b,\gamma(b) \big) & = \gamma(b)-(b-t_0) \cdot \gamma'(b) \\
			& \ge \gamma\left(\frac{a+b}{2}\right)-\Big(\frac{a+b}{2}-t_0\Big) \cdot  \displaystyle\frac{2}{b-a} \cdot \bigg(\gamma(b)-\gamma \Big(\frac{a+b}{2} \Big) \bigg)\\
			&\ge \gamma\left(\frac{a+b}{2}\right)-\Big(\frac{a+b}{2}-t_0 \Big) \cdot \displaystyle\frac{2}{b-a} \cdot \bigg( \gamma \Big(\frac{a+b}{2}\Big) - \gamma (a) \bigg)+ M^{-1} (b-a)^2 \\
			&\ge \gamma(a)-(a-t_0) \cdot \gamma'(a) +M^{-1} (b-a)^2= z \big(\gamma(a),a \big)+M^{-1} (b-a)^2,
			\end{aligned}
		\end{align*}
	
		\noindent which confirms \eqref{zgammaab}. Here, the first statement follows from \Cref{zb4}; the second is from \eqref{gamma1b}; the third from  \eqref{e:strongconvex} (with the fact that $t_0 \le t_1 \le a$); the fourth from \eqref{gamma1a}; and the fifth from \Cref{zb4}. In particular, $z$ is increasing in $t$ (and thus injective) on $\partial_{\no} \mathfrak{R}$. Moreover, setting $(a, b)$ equal to $(t_1, \mathfrak{t})$ and $(\mathfrak{t}, t_2)$ in \eqref{zgammaab} and using the fact that $t_2 - \mathfrak{t} = \Delta = \mathfrak{t} - t_1$, we deduce
		\begin{align}\label{e:endpoint}
			z \big(t_2,\gamma(t_2) \big)- z \big(\ft, \gamma(\ft) \big) \ge M^{-1} \Delta^2; \qquad z \big(\ft, \gamma(\ft) \big)-z \big(t_1,\gamma(t_1) \big) \ge M^{-1} \Delta^2.
		\end{align}

		Now let us analyze $z (\partial_{\we} \mathfrak{R})$ and $z (\partial_{\ea} \mathfrak{R})$. To this end, observe for any $t \in \{ t_1, t_2 \}$ and $x \in \big[ \gamma(t) - \delta, \gamma(t) \big]$ that 
		\begin{align}
			\label{e:ztxdiff}
			\begin{aligned}
			\Big|z(t,x)-z \big(t,\gamma(t) \big) \Big| &\le \big|x-\gamma(t) \big| +(t-t_0) \cdot \bigg( \Big| u_t (x) - u_t \big( \gamma(t) \big) \Big|+\pi \Big| \varrho_{t}(x)- \varrho_{t}(\gamma(t) \big) \Big| \bigg)\\
			&\le \delta+ 3 \Delta \cdot (\pi + 1) M^{1/2} \delta^{1/2} \leq \displaystyle\frac{\Delta^2}{4M},
		\end{aligned}
		\end{align}
	
		\noindent where in the first inequality we used \eqref{e:zxt} and \Cref{zb4}; in the second we used the facts that $0 \le \gamma(t) - x \le \delta$ and $0 \le t - t_0 \le 3 \Delta$ for $t \in \{ t_1, t_2 \}$, and \eqref{e:rhouholder}; and in the third we used the definition \eqref{deltatt} of $\delta$. Thus, 
		\begin{flalign}
			\label{zr1}
			\begin{aligned}
			\displaystyle\sup_{(t, x) \in \partial_{\we} \mathfrak{R}} \Real \Big( z(t, x) - z \big( \mathfrak{t}, \gamma (\mathfrak{t}) \big) \Big) \le -\displaystyle\frac{3\Delta^2}{4M} \le -3c_1 \Delta^2; \\
			\displaystyle\sup_{(t, x) \in \partial_{\ea} \mathfrak{R}} \Real \Big( z(t, x) - z \big( \mathfrak{t}, \gamma (\mathfrak{t}) \big) \Big) \ge \displaystyle\frac{3 \Delta^2}{4M} \ge 3c_1 \Delta^2, 
			\end{aligned}
		\end{flalign}
		
		\noindent where in the first inequality we used \eqref{e:endpoint} and \eqref{e:ztxdiff}, and in the second we used \eqref{deltatt}. 
		
		Finally, to analyze $z(\partial_{\so} \mathfrak{R})$, observe using \Cref{ztx2} and \eqref{e:densitylowb} that
		\begin{align}
			\label{zr2} 
			\begin{aligned}
			\Imaginary z \big( t, \gamma(t) - \delta \big) &  =  (t_0-t)\pi \cdot \varrho_t \big(\gamma(t) - \delta \big) \leq -3\pi M^{-1} \Delta \delta^{1/2} \leq -3c_1\Delta^2.
			\end{aligned} 
		\end{align}
		
		\noindent Hence, the image $z(\partial_{\so} \mathfrak{R})$ of the south boundary of $\mathfrak{R}$ is distance $3c_1\Delta^2$ away from the real axis, and thus from $z (\partial_{\no} \mathfrak{R})$. 
		
		We will apply \Cref{degreed} with the $(G; \mathfrak{R}; \mathfrak{W})$ there equal to $(z; \mathfrak{R}; \mathfrak{W})$ here (recall $\mathfrak{W}$ from \Cref{i:analytic} in \Cref{p:analytic}). That $z$ is real analytic and strictly positively oriented away from its critical points follows from \Cref{p:zprop}, so we must verify that $\mathfrak{W}$ satisfies the four properties listed in \Cref{degreed}. By \eqref{zgammaab}, $z \big( t, \gamma(t) \big) \in \mathbb{R}$ is increasing in $t$. Since \eqref{e:endpoint} implies that 
		\begin{flalign} 
			\label{zt1t2t3} 
			z \big( t_1, \gamma(t_1) \big)< z \big( \mathfrak{t}, \gamma(\mathfrak{t}) \big) < z \big( t_2, \gamma(t_2) \big),
		\end{flalign}  the continuity of $z$ yields some $t_3 \in [t_1, t_2]$ such that $z \big( t_3, \gamma(t_3) \big) = \mathfrak{t}$. By the continuity of $z$ (and the fact that $\varrho_t (x) > 0$ for $\gamma(t) - B < x < \gamma(t)$, since $\Omega^{\inv} = (0, 1) \times (0, L^{3/2})$, by \Cref{p:densitylow}), it follows that $\big( t_3, \gamma(t_3) - \varepsilon \big) \in \mathfrak{R}$ and $z \big( t_3 , \gamma(t_3) - \varepsilon \big) \in \mathfrak{W}$ for a sufficiently small real number $\varepsilon > 0$. Hence, $G(\mathfrak{R}) \cap \mathfrak{W}$ is nonempty, verifying the first property in \Cref{degreed}. 
		
		From \eqref{e:endpoint}, \eqref{zr1}, and \eqref{zr2}, we deduce that $\dist \big( z (\partial_{\ea} \mathfrak{R}) \cup z (\partial_{\so} \mathfrak{R}) \cup z (\partial_{\we} \mathfrak{R}); (\mathfrak{t}, \gamma(\mathfrak{t})) \big) \ge 3 c_1 \Delta^2$, and so $\overline{\mathfrak{W}} = \big\{ w \in \mathbb{H}^- : \big| w - z(\mathfrak{t}, \gamma(\mathfrak{t})) \big| \le 2c_1 \Delta^2 \big\}$ is disjoint from $z(\partial_{\ea} \mathfrak{R}) \cup z (\partial_{\so} \mathfrak{R}) \cup z(\partial_{\we} \mathfrak{R})$. Moreover, since $z(\partial_{\no} \mathfrak{R}) \subseteq \mathbb{R}$ (by \eqref{zb4}) and $\mathfrak{W} \subset \mathbb{H}^-$, we also have $z(\partial_{\no} \mathfrak{R})$ is disjoint from $\mathfrak{W}$. Hence, $z(\partial \mathfrak{R})$ is disjoint from $\mathfrak{W}$, verifying the second property in \Cref{degreed}. 
		
		The bounds \eqref{e:endpoint}, \eqref{zr1}, and \eqref{zr2} with the fact that $z(\partial_{\no} \mathfrak{R}) \subset \mathbb{R}$ also quickly imply that $z(\partial \mathfrak{R})$ can be continuously deformed in $\mathbb{C} \setminus \mathfrak{W}$ to the boundary of the rectangle with corners $\big\{ z(t_1, \gamma(t_1)), z(t_1, \gamma(t_1)) - \mathrm{i}, z(t_2, \gamma(t_2)) - \mathrm{i}, z (t_2, \gamma(t_2)) \big\}$. Consequently, the winding number of $G(\partial \mathfrak{R})$ around any point in $\mathfrak{W}$ is equal to one, confirming the third property in \Cref{degreed}.
		
		To establish the fourth, first observe from \eqref{e:endpoint}, \eqref{zr1}, and \eqref{zr2} that $\overline{\mathfrak{W}} \cap G(\partial \mathfrak{R}) = \overline{\mathfrak{W}} \cap G(\partial_{\no} \mathfrak{R})$. By \Cref{zb4}, $G(\partial_{\no} \mathfrak{R}) \subset \mathbb{R}$ and, by \eqref{e:zxt} with the fact that $\varrho_t (x) > 0$ for $\gamma(t) - B \le x < \gamma(t)$, we have $z(t, x) \in \mathbb{R}$ if and only if $(t, x) \in \partial_{\no} \mathfrak{R} = \big\{ (t, \gamma(t)) : t \in [t_1, t_2] \big\}$. Since \eqref{zgammaab} implies that $z \big( t, \gamma(t) \big)$ is increasing in $t$, it follows that $z$ is injective onto its image in $\overline{\mathfrak{W}} \cap z(\partial \mathfrak{R}) = \overline{\mathfrak{W}} \cap z(\partial_{\no} \mathfrak{R})$, verifying the fourth property in \Cref{degreed}. 		
				
		Thus, \Cref{degreed} applies; denoting $\mathfrak{U} = z^{-1} (\mathfrak{W}) \cap \overline{\mathfrak{R}}$, it implies that the map $z : \mathfrak{U} \rightarrow \mathfrak{W}$ is a homeomorphism. Since $\mathfrak{U} \subseteq \mathfrak{R}$, we also have that $(t, x) \in \mathfrak{U}$ implies $t \in [t_1, t_2] = [\mathfrak{t} - \Delta, \mathfrak{t} + \Delta]$, from which the first statement of the first part of the proposition follows.
	\end{proof}
	\begin{proof}[Proof of \Cref{i:analytic} in \Cref{p:analytic}]	
		
		Let us first show that $F : \mathfrak{W} \rightarrow \mathbb{C}$ is holomorphic. To this end, fix some point $w \in \mathfrak{W}$, and let $(t', x') \in \Omega(t_0)$ be such that $w = z(t', x')$. We claim that, if $(t', x')$ is not a critical point of $z$, then $\partial_{\bar z} F (w) =0$. We first establish the holomorphicity of $F$ assuming this claim. 
		
		To this end, observe that, since $z$ is real analytic by \Cref{p:zprop}, the image of its critical points is discrete; thus, $F$ is holomorphic away from a discrete set of points. Moreover, since $z : \mathfrak{U} \rightarrow \mathfrak{W}$ is a homeomorphism (by \Cref{i:injection} of the proposition) and $f$ is continuous on $\Omega(t_0)$ (by \eqref{frhou} and the fact that $\varrho$ and $u$ are smooth on $\Omega$, due to \Cref{gderivatives0}), the function $F$ is continuous on $\mathfrak{W}$. By Riemann's theorem on removable singularities, it follows that $F$ is a holomorphic function on $\mathfrak{W}$.  

		To show that $\partial_{\bar z} F(w) = 0$ unless $(t', x')$ is a critical point of $z$, suppose that $\partial_{\overline{z}} F(w) \ne 0$. Taking derivatives with respect to $t$ and $x$ of the relation $F \big(z(t,x) \big)=f(t,x)$, we get
		\begin{align}\begin{split}\label{e:xt}
				&\partial_z F \big( z(t,x) \big) \cdot \partial_t z (t,x)+\partial_{\bar z} F \big(z(t,x) \big) \cdot \partial_t \bar z(t,x)=\partial_t f(t, x),\\
				&\partial_z F \big(z(t,x) \big) \cdot \partial_x z(t,x)+\partial_{\bar z} F \big(z(t,x) \big) \cdot \partial_x \bar z(t,x)=\partial_x f(t, x).
		\end{split}\end{align}
		
		\noindent By \eqref{e:burger}, \eqref{e:dtzdxz}, multiplying the second relation in \eqref{e:xt} by $f(t, x)$, and summing with the first relation there, we deduce
		\begin{align}\label{e:barz}
			\partial_{\bar z} F \big(z(t,x) \big) \big(\partial_t \bar z(t,x)+f(t, x) \cdot \partial_x \bar z(t,x) \big)=0.
		\end{align}
		
		\noindent Thus, setting $(t, x) = (t', x')$ (and using the definition $z(t', x') = w$), we deduce since $\partial_{\bar z} F(w) \ne 0$ that 
		\begin{flalign}
			\label{tzfxz}
			\partial_t \bar z(t',x')+f(t', x') \cdot \partial_x \bar z(t',x')=0.
		\end{flalign}
	
		\noindent Taking the complex conjugate of \eqref{e:dtzdxz} yields $\partial_t \bar z(t',x')+\bar f(t, x) \cdot \partial_x \bar z(t',x')=0$, which upon subtraction from \eqref{tzfxz} yields $2 \Imaginary f(t', x') \cdot \partial_x \overline{z} (t', x') = 0$. Since $\Imaginary f(t', x') = \varrho_{t'} (x') > 0$ (as $(t', x') \in \mathfrak{U} \subseteq \Omega$), it follows that $\partial_x \overline{z} (t', x') = 0$, meaning by conjugation that $\partial_x z(t', x') = 0$. Together with \Cref{p:zprop}, this implies that $\partial_t z(t', x') = 0$ and that $(t', x')$ is a critical point of $z$, verifying the claim. 
		
		Thus, $F$ is holomorphic on $\mathfrak{W}$. Since \Cref{zb4} implies that $F(w)$ is real for $w\in \big\{z(t,\gamma(t)): t \in [t_1, t_2] \big\}\subset \mathbb{R}$, the reflection principle indicates we can extend $F$ to a holomorphic function on $\big\{z: |z-z(\ft, \gamma(\ft))|\le 2c_1\Delta^2 \big\}$ such that $\overline{F(z)} = F(\overline{z})$. This confirms the first statement of \Cref{i:analytic} of the proposition. To establish the second (given by \eqref{e:fkbound}), define the contour 
		\begin{flalign*} 
			\mathcal{C}=\bigg\{w \in \mathbb{C}: \Big|w-z \big(\ft, \gamma(\ft) \big) \Big|=2c_1\Delta^2 \bigg\}.
		\end{flalign*} 
	
		\noindent Since $z^{-1}(\mathcal{C} \cap \mathbb{H}^-) \subseteq \mathfrak{R}$, we have from \eqref{festimatetx} that $\big| F(w) \big| = \big| f (z(w)) \big|\leq M$ for $w\in \mathcal C$. Moreover, for any $z$ with $\big|z-z(\ft, \gamma(\ft)) \big|\leq c_1\Delta^2$, we have $\dist(z, \mathcal C)\geq c_1\Delta^2$. Together with Cauchy's formula, this gives 
		\begin{align}\label{e:fder}
			\big| \partial_z^k F(z) \big|=\left|\frac{1}{2\pi k!}\oint_{\mathcal{C}} \frac{F(w)d w}{(w-z)^{k+1}}\right|\leq \frac{2M}{ k! c_1^{k}\Delta^{2k}},
		\end{align}
		
		\noindent for any $z \in \mathbb{C}$ satisfying $\big|z-z(\ft,\gamma(\ft)) \big|\leq c_1\Delta^2$, verifying the second part of the proposition.
		\end{proof}
	
		\begin{proof}[Proof of \Cref{i:injectiondomain} in \Cref{p:analytic}]
		
		By \Cref{i:injection} of the proposition and \eqref{deltatt}, it suffices to show for any $(t, x) \in \mathfrak{R}$, with $t \in \big[ \mathfrak{t} - c_1 \Delta / (24 M), \mathfrak{t} + c_1 \Delta / (24 M) \big]$ and $x \in \big[ \gamma(t) - c_1^2 \Delta^2 / 900M, \gamma(t) \big]$, that $\big| z(t, x) - z ( \mathfrak{t}, \gamma(\mathfrak{t})) \big| \le c_1 \Delta^2$. To this end, observe that 
		\begin{flalign}
			\label{ztxz} 
			\Big| z(t,x) - z \big( \mathfrak{t}, \gamma (\mathfrak{t}) \big) \Big| \le \Big| z \big( t, \gamma(t) \big) - z \big( \mathfrak{t}, \gamma(\mathfrak{t}) \big) \Big| +	\Big| z(t, x) - z \big( t, \gamma(t) \big) \Big|.
		\end{flalign}
		
		\noindent We first estimate the second term in \eqref{ztxz}. Denote $a = \min \{ t, \mathfrak{t} \}$ and $b = \max \{ t, \mathfrak{t} \}$, and also let $a'=a-(b-a)$ and $b'=b+(b-a)$. Then, we have
		\begin{align*}
			z \big(b,\gamma(b) \big)
			=\gamma \big(b)-(b-t_0) \cdot \gamma'(b)
			& \leq \gamma(b)-\left(b-t_0\right) \cdot \frac{\gamma(b')-\gamma(b)}{b-a}\\
			&\leq \gamma(b)-(b-t_0) \cdot \frac{\gamma(b)-\gamma(a)}{b-a}
			+6\Delta M(b-a) \\
			& = \gamma(a)-(a-t_0) \cdot \frac{\gamma(b)-\gamma(a)}{b-a}
			+6\Delta M(b-a)
			\\
			&\leq \gamma(a)-(a-t_0) \cdot \frac{\gamma(a)-\gamma(a')}{b-a}
			+12\Delta M(b-a)\\
			&\leq \gamma(a)-(a-t_0) \cdot \gamma'(a)+12\Delta M(b-a) \\
			& = z \big(\gamma(a),a \big)+12\Delta M(b-a),
		\end{align*} 
		where the first and seventh statements follow from \Cref{zb4}; the second and sixth from the fact (by \eqref{e:strongconvex}) that $\gamma(t)$ is concave (observe that $b' \le b + 2 \Delta \le \mathfrak{t} + 3 \Delta \le 1 - 4B^{-1}$, and similarly $a \ge 4B^{-1}$, so these bounds apply); the third and fifth by \eqref{e:strongconvex} and the fact that $a - t_0 \le b - t_0 \le t_2 - t_0 \le 3 \Delta$; and the fourth by performing the addition. Together with the facts that $z \big( t, \gamma(t) \big)$ is increasing in $t \in [t_1, t_2]$ (by \eqref{zgammaab}) and the bound $|t-\ft| \le c_1 \Delta/(24 M)$, this gives
		\begin{align}
			\label{ztxz1} 
			\Big| z \big(t,\gamma(t) \big)-z \big(\ft, \gamma(\ft)\big) \Big| \leq 12 \Delta M|t-\ft|\leq \displaystyle\frac{c_1\Delta^2}{2}.
		\end{align}
	
		To bound the second term on the right side of \eqref{ztxz}, we follow \eqref{e:ztxdiff} to obtain, for any $t \in [t_1, t_2]$ and $x \in \big[ \gamma(t) - c_1^2 \Delta^2 / 900M, \gamma(t) \big]$ that
		\begin{align*}
			\Big|z(t,x)-z \big(t,\gamma(t) \big) \Big| &\le \big|x-\gamma(t) \big|+(t-t_0) \cdot \bigg( \Big| u_{t}(x) - u_t \big( \gamma(t) \big) \Big|+\pi \Big| \varrho_{t}(x)-\varrho_t \big(\gamma(t) \big) \Big| \bigg) \\
			&\le \displaystyle\frac{c_1^2 \Delta^2}{900M} +3 \Delta \cdot (\pi + 1) M^{1/2} \cdot \displaystyle\frac{c_1 \Delta}{30M^{1/2}}  \leq \displaystyle\frac{c_1\Delta^2}{2}.
		\end{align*}
	
		\noindent Together with \eqref{ztxz} and \eqref{ztxz1}, this gives $\big| z(t, x) - z (\mathfrak{t}, \gamma(\mathfrak{t})) \big| \le c_1 \Delta^2$, which as mentioned above yields the proposition. 	
	\end{proof}

\subsection{Proof of Edge Behavior of $G$}
	
	\label{FunctionEdge} 
	
	In this section we establish \Cref{p:limitprofile}. In what follows, we fix some real number $\Delta \in (0, 1 / 4B]$ and set $t_0 = \mathfrak{t} - 2 \Delta$, $t_1 = \mathfrak{t} - \Delta$, and $t_2 = \mathfrak{t} + \Delta$. Define the characteristic map $z = z_{t_0} : \Omega(t_0) \rightarrow \mathbb{H}$ as in \Cref{ztx2}. Then, \Cref{p:analytic} applies, and we adopt the notation of that proposition in what follows, but write $C_1 = C_1 (B) > 1$ for the constant $C = C(B) > 1$ appearing there. Observe by \eqref{e:gammaconv}, and by the fact that on $[4B^{-1}, 1-4B^{-1}]$ we have $\gamma$ is continuously differentiable (from \eqref{e:gammadc}), that \Cref{p:checka} yields a constant $M = M(B) > 1$ such that 
		\begin{align}\label{e:deribb}
		\big| \gamma'(t) \big| \le M, \quad \text{and} \quad M^{-1} \leq \frac{\gamma'(t)-\gamma'(t')}{t' - t}\leq M,\qquad \text{for each $4B^{-1} \le t \le t' \le 1-4B^{-1}$}.
		\end{align}
	
	\noindent In what follows, we define the function $\xi : [4B^{-1}, 1-4B^{-1}]$ by setting
		\begin{flalign}
			\label{xifunctiont} 
			\xi (t) = \gamma (t) - (t-t_0) \cdot \gamma' (t), \qquad \text{for each $t \in [4B^{-1}, 1 - 4B^{-1}]$}.
		\end{flalign}

	Now, fix some $(t, x) \in \overline{\mathfrak{U}}$ (recall from \Cref{i:injection} in \Cref{p:analytic}) such that $z = z(t, x) \in \overline{\mathbb{H}}^-$ satisfies $\big| z-z(\ft, \gamma(\ft)) \big|\le c_1 \Delta^2$. Given such a $z$ we can recover $(t,x)$, as follows. We distinguish two cases, the first being if $z \in \mathbb{H}^-$ (meaning $\varrho_t (x) > 0$ by \eqref{e:zxt}, so $(t, x) \in \Omega$ by \eqref{omega12}) and the second being if $z \in \mathbb{R}$ (meaning $\varrho_t (x) = 0$ by \eqref{e:zxt}, so $x = G(t, 0) = \gamma(t)$ since $\Omega^{\rm inv} = (0, 1) \times (0, L^{3/2})$ by \Cref{i:densitylow} of \Cref{p:densitylow}). 
	
	In the first case (by the fact that $F(z) = f(t, x)$ and by \eqref{e:zxt}) we can solve for $(t,x)$ by
		\begin{align}\label{e:constructtx}
		t = t_0- \displaystyle\frac{\Imaginary z}{\Imaginary F(z)}; \qquad x = \Real z + (t - t_0) \cdot \Real F(z) = z + (t-t_0) \cdot F(z), \qquad \text{if $z \in \mathbb{H}$},
		\end{align}
		
		\noindent and so    
		\begin{align}\label{e:constructurho}
			u_t(x)+\pi \mathrm{i} \varrho_t(x)= F(z) = f(t, x) = f \big(t, z+(t-t_0) F(z) \big).
		\end{align}
		
		\noindent In the second case, we have $x = \gamma(t)$, and \Cref{zb4} and \eqref{ufrhogammat} imply that  
		\begin{align}\label{e:gammat}
			z=\gamma(t)-(t-t_0) \cdot \gamma'(t) = \xi(t); \qquad \gamma'(t)= f \big(t,\gamma(t) \big)=F(z)=F\big( \xi(t) \big).
		\end{align}
	
		We next have the following lemma that evaluates the derivatives of $F$, $\xi$, and $\gamma$. It also Taylor expands $\gamma$, which will be used to show the $y=0$ case of \Cref{p:limitprofile}.
		
	\begin{lem} 
		
		\label{gamma123} 
		
		The following hold for any $t \in [t_1, t_2]$. 
		
		\begin{enumerate} 
			\item The functions $\xi$ and $\gamma$ are both smooth on $[t_1, t_2]$, and $\xi$ is moreover increasing on $[t_1, t_2]$. 
			\item We have $F' \big( \xi(t) \big) = -(t-t_0)^{-1}$ and $\gamma(t) = \xi(t) + (t-t_0) \cdot F \big( \xi(t) \big)$.  
			\item We have  	
			\begin{align}\label{e:xider}
				 \xi'(t)=\frac{1}{(t-t_0)^2 F'' \big( \xi(t) \big)},\quad  \gamma''(t)=-\frac{1}{(t-t_0)^3 F'' \big( \xi(t) \big)},\quad F'' \big( \xi(t) \big)=-\frac{1}{(t-t_0)^3 \gamma''(t)}.
			\end{align}
			\item  We  have 
		\begin{flalign}\label{e:gamma'''}
			 M^{-1} \Delta^{-3}  \le F'' \big( \xi(t) \big) \le 27 M \Delta^{-3}; \qquad M^{-1} \Delta \le \xi' (t) \le 3M \Delta; \qquad \big| \gamma'''(t) \big| \le \displaystyle\frac{2 C_1 M^3}{c_1^3 \Delta^2}.
		\end{flalign}
	
		\item For any real number $\tau \in [-\Delta, \Delta]$, we have   
		\begin{flalign*}
			\bigg| \gamma(\ft+\tau)-\Big(\gamma(\ft)+\gamma'(\ft) \cdot \tau+ \frac{\gamma'' (\mathfrak{t})}{2} \cdot \tau^2\Big) \bigg| \leq \frac{2 C_1 M^3}{c_1^3\Delta^2} \cdot |\tau|^3,
		\end{flalign*}
		\end{enumerate} 
		
	\end{lem} 

	\begin{proof}

		By the second statement in \eqref{e:deribb}, $\gamma(t)$ is (strictly) concave on $[t_1, t_2]$, and so $\gamma'$ is (strictly) decreasing. Hence $\xi$ is (strictly) increasing on $[t_1, t_2]$, as for $t_1 \le t' < t \le t_2$ we have 
		\begin{flalign*} 
			\xi (t) - \xi (t') = \gamma(t) - \gamma(t') - (t-t_0) \cdot \gamma' (t) + (t'-t_0) \cdot \gamma' (t') \ge \gamma(t) - \gamma(t') - (t-t') \cdot \gamma' (t)  > 0,
		\end{flalign*} 
	
		\noindent where in the first statement we used \eqref{xifunctiont}; in the second we used the fact that $\gamma'$ is decreasing (and that $t' \ge t_0$); and in the third we used the fact that $\gamma$ is concave. Next we show that $\gamma(t)$ is  smooth and compute its derivatives in terms of $F$. We first compute the derivative of $F \big( \xi(t) \big)$, obtaining 
		\begin{align}
			\label{e:f'bb}
			\begin{aligned}
		F' \big( \xi(t) \big) &  =\lim_{t'\rightarrow t}\frac{F \big( \xi(t) \big)-F \big(\xi(t') \big)}{\xi(t)-\xi(t')} \\ 
		& =\lim_{t'\rightarrow t}\frac{\gamma'(t)-\gamma'(t')}{\gamma(t)-\gamma(t')-(t-t') \cdot \gamma'(t')-(t-t_0) \cdot \big(\gamma'(t)-\gamma'(t') \big)}\\
		&=\lim_{t'\rightarrow t} \bigg( \frac{\gamma(t)-\gamma(t')-(t-t') \cdot \gamma'(t')}{\gamma'(t)-\gamma'(t')}-(t-t_0) \bigg)^{-1} =-(t-t_0)^{-1},
		\end{aligned}
	\end{align}
		where the first equality holds by the holomorphicity of $F$ (by \Cref{p:analytic}) and the continuity of $\xi$ (by \eqref{xifunctiont} and the fact that $\gamma'$ is continuous); the second by \eqref{e:gammat} and \eqref{xifunctiont}; the third by dividing the numerator and denominator by $\gamma'(t)-\gamma'(t')$; and the last by the continuous differentiability of $\gamma$ (from \eqref{e:gammadc} and the fact that $t, t' \in [t_1, t_2] \subseteq [4B^{-1}, 1 - 4B^{-1}]$), with the fact that $\gamma'(t) - \gamma'(t') \ge M^{-1} (t-t')$ (by \eqref{e:strongconvex}). The equality \eqref{e:f'bb}, together with \eqref{xifunctiont} and the second statement of \eqref{e:gammat}, yields the second statement of the lemma.
		
		Taking a further $t$-derivative on both sides of \eqref{e:f'bb} gives 
		\begin{align*}
	\frac{1}{(t-t_0)^2}
	&=F'' \big( \xi(t) \big) \cdot \lim_{t'\rightarrow t}\frac{\xi(t)-\xi(t')}{t-t'}\\
	&=F'' \big( \xi(t) \big) \cdot \lim_{t'\rightarrow t}\frac{\gamma(t)-\gamma(t')-(t-t') \cdot \gamma'(t')-(t-t_0) \cdot \big(\gamma'(t)-\gamma'(t') \big)}{t-t'}\\
	&= - F'' \big( \xi(t) \big) \cdot \lim_{t'\rightarrow t} \frac{(t-t_0) \cdot \big(\gamma'(t)-\gamma'(t') \big)}{t-t'},
\end{align*}
	where the second equality follows from \eqref{xifunctiont}, and the third equality follows from that $\gamma(t)$ is continuously differentiable (by \eqref{e:gammadc}). Together with the second statement of \eqref{e:deribb} and the holomorphicity of $F$, this implies that $F'' \big( \xi(t) \big) \ne 0$; hence, $\gamma$ is twice-differentiable on $[t_1, t_2]$, and $\gamma'' (t) = - F'' \big( \xi(t) \big)^{-1} \cdot (t - t_0)^{-3}$. This implies the second equality of 
	\eqref{e:xider}, of which the third is a consequence. The first equality there follows from first differentiating \eqref{xifunctiont}, which yields $\xi'(t)=-(t-t_0) \cdot \gamma''(t)$, and then applying the second equality of \eqref{e:xider}.
	
	The fact that $F'' \big( \xi(t) \big) \ne 0$, together with the holomorphicity of $F$, the identity $F' \big( \xi(t) \big) = -(t-t_0)^{-1}$ (from the second part of the lemma), and the Inverse Function Theorem, implies that $\xi$ is smooth on $[t_1, t_2]$. From this, it follows by differentiating \eqref{xifunctiont} (and using the continuous differentiability of $\gamma$ from \eqref{e:gammadc}) that $\gamma$ is continuously twice-differentiable on $[t_1, t_2]$; repeatedly differentiating \eqref{xifunctiont} then yields that $\gamma$ is smooth on $[t_1, t_2]$. Since we confirmed above that $\xi$ is increasing on $[t_1, t_2]$, this yields the first statement of the lemma.
	
	It thus remains to establish the last two statements of the lemma. To show the fourth, observe by the second statement of \eqref{e:deribb} that $M^{-1} \leq -\gamma''(t)\leq M$. Moreover, we have from the third and first statements of \eqref{e:xider} (the latter being equivalent to $\xi'(t) = -(t - t_0) \cdot \gamma''(t)$, by the second equality of \eqref{e:xider}) that  
		\begin{align}\label{e:fttbb}
		M^{-1} \Delta^{-3} \leq F'' \big( \xi(t) \big)\leq 27 M\Delta^{-3}; \qquad M^{-1} \Delta \leq \xi'(t)\leq 3 M \Delta,
		\end{align}
	
		\noindent where we used the bounds $\Delta \le t - t_0 \le 3 \Delta$ for $t \in [t_1, t_2]$. Thus, 
		\begin{align}\label{e:fttbb2}
			\big|\gamma'''(t) \big|=\bigg|\frac{3}{(t-t_0)^4 F'' \big( \xi(t) \big)} \bigg| + \bigg| \frac{F''' \big(\xi(t) \big)}{(t-t_0)^5 F'' \big( \xi(t) \big)^3} \bigg|\leq \displaystyle\frac{3 M}{\Delta} + \displaystyle\frac{C_1 M^3}{6c_1^3 \Delta^2} \le \frac{2 C_1 M^3}{c_1^3\Delta^2}.
		\end{align}
		 
		 \noindent Here, to deduce the first inequality we differentiated the second equality of \eqref{e:xider}, and also used the first equality there; to deduce the second, we used \eqref{e:fttbb} and the facts that $t - t_0 \ge t_1 - t_0 \ge \Delta$ and $\big| F'''(\xi(t)) \big|\leq C_1 /(6c_1^3\Delta^6)$ (the latter by \eqref{e:fkbound}); and to deduce the third we used the fact that $4\Delta \le B^{-1} \le 1$. The fourth part of the lemma then follows from \eqref{e:fttbb} and \eqref{e:fttbb2}. The fifth then follows from the last bound in the fourth, together with a Taylor expansion, thereby establishing the lemma.
	\end{proof}

	Recalling the density process $(\varrho_t)$ associated with $\bm{\mu}$ (and thus with $G$) from \Cref{lmu}, we next approximate $\varrho_t(x)$ around $(t, x) = (\ft, \gamma(\ft))$ through the following lemma. As a corollary, we deduce a bound on $\gamma(t) - G(t, x)$, from which \Cref{p:limitprofile} quickly follows.

	\begin{lem} 
		
		\label{rhoxtgamma}

		There exist constants $c_3 = c_3 (B, \Delta) \in (0, c_2 \Delta^2) \subset (0, 1)$ and $C_2 = C_2 (B, \Delta) > 1$ such that, for any real numbers $\tau \in [-c_3, c_3]$ and $x \in \big[ \gamma(\mathfrak{t} + \tau) - c_3, \gamma(\mathfrak{t} + \tau) \big)$, we have 
		\begin{flalign}
			\label{rhott} 
			 \varrho_{\mathfrak{t} + \tau} (x) = \pi^{-1} \big( 1 + \mathcal{E} (\tau, x) \big) \cdot \Big| 2 \gamma''(\mathfrak{t}) \cdot \big( \gamma (\mathfrak{t} + \tau) - x \big) \Big|^{1/2},
		\end{flalign}
	
		\noindent for some quantity $\mathcal{E} (\tau, x) \in \mathbb{R}$ satisfying 
		\begin{flalign*} 
			\big| \mathcal{E} (\tau, x) \big| \le C_2  \Big( |\tau| + \big| \gamma(\mathfrak{t} + \tau) - x \big|^{1/2} \Big).
		\end{flalign*} 
	\end{lem} 

	\begin{proof}
		
		Throughout this proof, we set $s_0 = \mathfrak{t} + \tau \in [\mathfrak{t} - c_2 \Delta^2, \mathfrak{t} + c_2 \Delta^2] \subset [\mathfrak{t} - \Delta, \mathfrak{t} + \Delta] = [t_0, t_1]$. Since $\gamma(s_0) - c_2 \Delta^2 \le x \le \gamma(s_0)$, \Cref{i:injectiondomain} in \Cref{p:analytic} (with the bound $|\mathfrak{t} - s_0| \le c_2 \Delta^2 \le c_2 \Delta$) gives $\big| z(s_0, x)-z(\ft, \gamma(\ft)) \big|\leq c_1 \Delta^2$. Denote 
		\begin{flalign}
			\label{w} 
			w = z(s_0, x) - \xi (s_0). 
		\end{flalign} 
	
		\noindent Observe that $w \in \mathbb{H}^-$, since $z (s_0) \in \mathbb{H}^-$ (and $\xi (s_0) \in \mathbb{R})$), which follows from \eqref{e:zxt} and the fact that $\varrho_{s_0} (x) > 0$ (the latter since $x < \gamma(s_0)$ and $\Omega^{\rm inv} = (0, 1) \times (0, L^{3/2})$ by \Cref{p:densitylow}). From \eqref{e:constructtx} and \eqref{e:constructurho} (with the fact that $2 \Delta + \tau = s_0 - t_0$, as $t_0 = \mathfrak{t} - 2 \Delta$), $w$ satisfies the relations
		\begin{align}\begin{split}\label{e:xdensity}
			&x=\xi(s_0)+w+(2\Delta+\tau) \cdot F \big(\xi(s_0)+w \big),\\
			& \varrho_{s_0} (x)=\pi^{-1} \cdot \Imaginary F \big(\xi(s_0)+w \big) = -\frac{\Imaginary w }{\pi(2\Delta+\tau)},
		\end{split}\end{align} 
	
		\noindent where in the last equality we used the fact that $\Imaginary w + (2 \Delta + \tau) \cdot \Imaginary F \big( \xi (s_0) + w \big) = 0$ (which follows from the first part of \eqref{e:constructtx}).
		
		Next let us bound $w$, to which end observe that there exists a constant $C_3 = C_3 (B) > 1$ so that
		\begin{flalign}
			\label{westimate0} 
			\begin{aligned}
			|w| & = \Big| z (s_0, x) - z \big( s_0, \gamma (s_0) \big) \Big| \\
			& \le \big| x - \gamma(s_0) \big| +  (s_0 - t_0) \cdot \Big| f(s_0, x) - f \big( s_0, \gamma (s_0) \big) \Big| \\
			& = \big| x - \gamma(s_0) \big| + 3\Delta \bigg( \big| u_{s_0} (x) - u_{s_0} \big( \gamma (s_0) \big) \Big| + \pi \Big| \varrho_{s_0} (x) - \varrho_{s_0} \big( \gamma (s_0) \big) \Big| \bigg)  \le C_3 \big| x - \gamma(s_0) \big|^{1/2}.
			\end{aligned}
		\end{flalign}
	
		\noindent Here, in the first statement we used \eqref{w} and the first equality in \eqref{e:gammat}; in the second we used \eqref{e:zxt}; in the third  we used \eqref{frhou} and the fact that $|s_0 - t_0| = 2 \Delta + \tau \le 3 \Delta$ (as $|\tau| \le c_2 \Delta^2 \le \Delta$); and in the fourth we used (the integral of) \eqref{e:holder}, with the facts that $\big| x - \gamma (s_0) \big| \le c_2 \Delta^2 < 1$ and that $\Delta < 1$. Then, by Taylor expanding $F$ around $\xi(s_0)$ and using \eqref{e:fkbound} (again with the fact that $|2 \Delta + \tau| \le 3 \Delta$, as $|\tau| \le \Delta$), the first relation in \eqref{e:xdensity} gives
			\begin{flalign}
				\label{e:xexp}
			x&=\xi(s_0)+w+(2\Delta+\tau) \Big( F \big(\xi(s_0) \big)+F' \big(\xi(s_0) \big) \cdot w +\displaystyle\frac{F'' \big( \xi(s_0) \big)}{2} \cdot w^2 \Big) + \cE(w), 
			\end{flalign} 
		
			\noindent for some complex number $\mathcal{E} (w) \in \mathbb{C}$ satisfying 
			\begin{flalign}
				\label{eestimatew} 
			 \big| \cE(w) \big| \le \displaystyle\frac{C_1 (2 \Delta + \tau) }{6 c_1^3 \Delta^6} \cdot |w|^3 \le \frac{C_1 |w|^3}{2c_1^3\Delta^5} \le \displaystyle\frac{C_1 C_3^3}{2 c_1^3 \Delta^5} \cdot \big( \gamma(s_0) - x \big)^{3/2},
		\end{flalign}
		 	 
			\noindent where in the last inequality we used \eqref{westimate0}. Moreover, by the second part of \Cref{gamma123} (and again the fact that $2 \Delta + \tau = s_0 - t_0$), we have 
		\begin{align}\label{e:gt3}
			\gamma(s_0)=\xi(s_0)+(2\Delta+\tau) \cdot F \big(\xi(s_0) \big); \qquad F' \big(\xi(s_0) \big) = -(2\Delta+\tau)^{-1}.
		\end{align}		
	
		Taking the difference between \eqref{e:gt3} and \eqref{e:xexp}, we find
		\begin{align*}
			\gamma(s_0)-x =-\Big( \Delta + \frac{\tau}{2} \Big) F'' \big(\xi(s_0) \big) \cdot w^2- \cE(w).
		\end{align*}

		\noindent This, together with the facts that $w \in \mathbb{H}^-$ and that $F'' \big( \xi(t) \big) > 0$ for $t \in [t_1, t_2]$ (by the third statement of \eqref{e:xider} and the concavity from \eqref{e:deribb} of $\gamma$), yields a constant $c_4 = c_4 (B, \Delta) \in (0, 1)$ such that for $x \in \big[ \gamma (s_0) - c_4, \gamma(s_0) \big]$ (implying by \eqref{westimate0} and \eqref{eestimatew} that $|w|$ and $\big| \mathcal{E} (w) \big|$ are sufficiently small) we have
		\begin{flalign*}
			w&=- 2^{1/2} \mathrm{i} \cdot  \Bigg( \frac{\gamma(s_0)-x + \cE(w)}{(2\Delta+\tau) \cdot F'' \big(\xi(s_0) \big)} \Bigg)^{1/2}.
		\end{flalign*}
	
		\noindent Hence, denoting $\xi = \xi (\mathfrak{t})$, we have 
		\begin{flalign*} 
			w =- 2^{1/2} \mathrm{i} \cdot \Bigg(\frac{\gamma(s_0)-x}{(2\Delta+\tau) \cdot F''(\xi)} \Bigg)^{1/2}  \Bigg( \frac{F''(\xi)}{F'' \big(\xi(s_0) \big)} \Bigg)^{1/2} \Bigg(1+\frac{\cE(w)}{\gamma(s_0)-x} \bigg)^{1/2},
		\end{flalign*}
	
		\noindent so the second statement of \eqref{e:xdensity} (together with the fact that $F'' \big( \xi(t) \big) > 0$ for $t \in [t_1, t_2]$) gives 
		\begin{align*}
		\varrho_{s_0}(x)= \pi^{-1} \bigg( \frac{\gamma(s_0)-x}{4\Delta^3 \cdot F''(\xi)} \bigg)^{1/2} \bigg( \frac{(2\Delta)^3 \cdot F''(\xi)}{(2\Delta+\tau)^3 \cdot F'' \big(\xi(s_0) \big)} \bigg)^{1/2} \Imaginary \bigg( 1  +\frac{\cE(w)}{\gamma(s_0) - x} \bigg)^{1/2}.
		\end{align*}
		
		\noindent Observe by the second statement of \eqref{e:xider} that $F'' (\xi)^{-1} = - (2 \Delta)^3 \cdot \gamma'' (\mathfrak{t})$ (as $\mathfrak{t} - t_0 = 2\Delta$), and so it follows that 
		\begin{flalign}
			\label{rhogammatt}	
			\varrho_{s_0}(x)= \pi^{-1} \Big| 2 \gamma''(\mathfrak{t}) \cdot \big( \gamma(s_0)-x \big) \Big|^{1/2} \bigg( \frac{(2\Delta)^3 \cdot F''(\xi)}{(2\Delta+\tau)^3 \cdot F'' \big(\xi(s_0) \big)} \bigg)^{1/2} \Imaginary \bigg(1+\frac{\cE(w)}{\gamma(s_0) - x} \bigg)^{1/2}.
		\end{flalign} 
	
		\noindent The first (of three) terms in the above product is in agreement with \eqref{rhott}; we must therefore approximate last two terms in this product by $1$.
		
		To this end observe, since for each $t \in [t_1, t_2]$ we have $\big|F'''(\xi(t)) \big|\leq C_1 /(6 c_1^3\Delta^6)$ (by \eqref{e:fkbound}) and $M^{-1} \Delta \le \big| \xi' (t) \big| \le 3 M \Delta$ (by \eqref{e:fttbb}), that 
		\begin{flalign*}
			\Big|F'' \big(\xi(s_0) \big)-F''(\xi) \Big| \le \displaystyle\frac{C_1 \big|\xi(s_0)-\xi (\mathfrak{t}) \big|}{6c_1^3\Delta^6} \leq \displaystyle\frac{ C_1 M |\tau|}{2c_1^3\Delta^5},
		\end{flalign*}
	
		\noindent where we have used the fact that $s_0 = \mathfrak{t} + \tau$. This, with the bounds $|\tau| \le c_2^2 \Delta < 1$ and $\big| F''(\xi(s_0)) \big | \ge  M^{-1} \Delta^{-3}$ (the latter of which holds by \eqref{e:fttbb}), yields 
		\begin{flalign*}
			 2 \bigg| \Big( 1 + \displaystyle\frac{\tau}{2\Delta} \Big)^{-3} -1 \bigg| \le \displaystyle\frac{6 |\tau|}{\Delta^3}; \qquad \frac{2 \Big|F''(\xi)-F'' \big( \xi(s_0) \big) \Big|}{\Big| F'' \big(\xi(s_0) \big) \Big|} \le \frac{C_1 M^2 |\tau|}{c_1^3\Delta^2},
		\end{flalign*}
		
		\noindent Together with the bound $|ab-1|\leq 2 \big(|a-1|+|b-1| \big)$ if $|a-1| \le 1$ and $|b-1|\leq 1$, this implies for sufficiently small $|\tau|\leq c_1^3 \Delta^3/(6 C_1 M^2)$ that
		\begin{align}\begin{split}\label{e:diyixiang}
		\left|\frac{(2\Delta)^3 \cdot F''(\xi)}{(2\Delta+\tau)^3 \cdot F'' \big( \xi(s_0) \big)}-1\right|
		&\leq 2 \bigg| \Big( 1 + \displaystyle\frac{\tau}{2\Delta} \Big)^{-3} -1 \bigg|+\frac{2 \Big|F''(\xi)-F'' \big( \xi(s_0) \big) \Big|}{\Big| F'' \big(\xi(s_0) \big) \Big|}\\
		 &\leq \frac{6|\tau|}{\Delta^3}+\frac{C_1 M^2 |\tau|}{c_1^3\Delta^2},
		\end{split}\end{align}
		     
		  \noindent which addresses the second term in \eqref{rhogammatt}. To address the third, observe from \eqref{eestimatew} that
		\begin{flalign*}
			\bigg| \displaystyle\frac{\mathcal{E} (w)}{\gamma(s_0) - x} \bigg| \le \displaystyle\frac{C_1 C_3^3}{2c_1^3 \Delta^5} \cdot \big( \gamma(s_0) - x \big)^{1/2}. 
		\end{flalign*}
	
		\noindent Applying this, with \eqref{e:diyixiang}, we deduce that there exist constants $c_5 = c_5 (B, \Delta) \in (0, 1)$ and $C_4 = C_4 (B, \Delta) > 1$ such that for $|\tau|\leq  c_5$ and  $x \in \big[ \gamma(s_0) - c_5, \gamma(s_0) \big]$ we have
	\begin{flalign*}  
	\Bigg| \bigg(\frac{(2\Delta)^3 \cdot F''(\xi)}{(2\Delta+\tau)^3 \cdot F'' \big(\xi(s_0) \big)} \bigg)^{1/2} \bigg(1+\frac{\cE(w)}{\gamma(s_0)-x} \bigg)^{1/2} -1 \Bigg| \le C_4 \Big( |\tau|+ \big(\gamma(s_0)-x \big)^{1/2} \Big).
		\end{flalign*}
	
	\noindent Together with \eqref{rhogammatt}, this yields the lemma. 
\end{proof}

\begin{cor} 
	
	\label{gammag}
	
	There exist constants $c_3 = c_3 (B, \Delta) \in (0, \Delta)$ and $C_3 = C_3 (B, \Delta) > 1$ such that, for any real numbers $\tau \in [-c_3, c_3]$ and $y \in [0, c_3]$, we have 
	\begin{flalign*}
		\bigg| G(\mathfrak{t} + \tau,y) - \gamma (\mathfrak{t} + \tau) - \Big( -\displaystyle\frac{9 \pi^2}{8 \gamma''(\mathfrak{t})} \Big)^{1/3} y^{2/3} \bigg| \le C_3 \big( |\tau| y^{2/3} + y \big).
	\end{flalign*}
\end{cor} 

\begin{proof} 
	
	Setting $s_0 = \mathfrak{t} + \tau$ and integrating \eqref{rhott} yields constants $c_3 = c_3 (B, \Delta) \in (0, 1)$ and $C_2 = C_2 (B, \Delta) > 1$ such that, for $|\tau| \le c_3$ and $x \in \big[ \gamma(s_0) - c_3, \gamma(s_0) \big]$, we have   
	\begin{flalign}
		\label{rhos0} 
		\Bigg| \displaystyle\int_x^{\gamma(s_0)} \varrho_{s_0} (x) dx - \displaystyle\frac{2}{3 \pi} \big( -2 \gamma'' (\mathfrak{t}) \big)^{1/2} \big( \gamma(s_0) - x \big)^{3/2} \Bigg| \le C_2 \Big( |\tau| \big( \gamma(s_0) - x \big)^{3/2} + \big( \gamma(s_0) - x \big)^2 \Big),
	\end{flalign}

	\noindent where we also used the fact that $-M \le \gamma'' (s_0) \le -M^{-1}$ (by \eqref{e:deribb}). Fix a real number $R = R (B, \Delta) > 1$, to be determined later, and set 
	\begin{flalign}
		\label{x03} 
		x_0 = \gamma(s_0) - \Big( -\displaystyle\frac{9\pi^2}{8 \gamma''(s_0)} \Big)^{1/3} y^{2/3}; \quad x_0^- = x_0 - R \big( |\tau| y^{2/3} + y \big); \quad x_0^+ = x_0 + R \big( |\tau| y^{2/3} + y \big),
	\end{flalign}

	\noindent Then, by a Taylor expansion (again using the fact that $M^{-1} \le -\gamma'' (s_0) \le M$, by \eqref{e:deribb}), there exist constants $c_4 = c_4 (B, \Delta, R) \in (0, 1)$, $c_5 = c_5 (B, \Delta) \in (0, 1)$, and $c_6 (B, \Delta) \in (0, 1)$ such that for $y \in [0, c_4]$ we have 
	\begin{flalign*}
		\displaystyle\frac{2}{3 \pi} \big( -2 \gamma''( & s_0) \big)^{1/2} \cdot \big( \gamma(s_0) - x_0^+ \big)^{3/2} \\
		& = \displaystyle\frac{2}{3 \pi} \big( -2\gamma'' (s_0) \big)^{1/2} \cdot \bigg( \Big( -\displaystyle\frac{9 \pi^2}{8 \gamma'' (s_0)} \Big)^{1/3} y^{2/3} - R \big( |\tau| y^{2/3} + y \big)  \bigg)^{3/2} \\
		& < y - c_5 R \big( |\tau| y + y^{4/3} \big) < y - c_6 R \Big( |\tau| \big( \gamma(s_0) - x_0^+ \big)^{3/2} + \big( \gamma(s_0) - x_0^+ \big)^2 \Big),
	\end{flalign*}

	\noindent where in the last bound we applied \eqref{x03} (which implies for some constant $c_7 = c_7 (B) > 0$ that $c_7 \big( \gamma(s_0) - x_0^+ \big)^{3/2} \le y \le c_7^{-1} \big( \gamma(s_0) - x_0^+ \big)^{3/2}$). For $R > c_6^{-1} C_2$, this bound with \eqref{rhos0} implies that
	\begin{flalign} 
		\label{yrhos0} 
		\displaystyle\int_{x_0^+}^{\gamma(s_0)} \varrho_{s_0} (x) dx < y.
	\end{flalign}

	\noindent By similar reasoning, we have (after increasing $R$ if necessary) that 
	\begin{flalign*}
		\displaystyle\int_{x_0^-}^{\gamma(s_0)} \varrho_{s_0} (x) dx > y. 
	\end{flalign*}

	\noindent Together with \eqref{yrhos0} and \eqref{gty} (with \eqref{htxintegral}), this implies that $x_0^- \le G(s_0, y) \le x_0^+$. By \eqref{x03}, this establishes the lemma.
\end{proof}

\begin{proof}[Proof of \Cref{p:limitprofile}]
	
		Define the real numbers $\mathfrak{a}$, $\mathfrak{b}$, and $\mathfrak{c}$ by setting
		\begin{flalign*} 
			\mathfrak{a} = \gamma (\mathfrak{t}); \qquad \mathfrak{b} = \gamma' (\mathfrak{t}); \qquad \mathfrak{c} = - \displaystyle\frac{\gamma''(\mathfrak{t})}{2}.
		\end{flalign*} 
	
		\noindent By \eqref{e:deribb}, there exists a constant $C_1 = C_1 (B) > 4B^2 > 1$ such that $|\mathfrak{b}| \le C_1$ and $C_1^{-1} \le \mathfrak{c} \le C_1$. We also have by the $r = 0$ case of \eqref{e:gtbound} that $\mathfrak{a} = \gamma(\mathfrak{t}) = G(\mathfrak{t}, 0) \le 4B^2 \le C_1$ and by the $r = 0$ case of \eqref{e:Glow} that $-\mathfrak{a} = -\gamma(\mathfrak{t}) \le 0 < C_1$, meaning that $|\mathfrak{a}| < C_1$. Equation \eqref{e:limitprofile} follows from combining the fifth part of \Cref{gamma123} and \Cref{gammag}. This finishes the proof of \Cref{p:limitprofile}.
	\end{proof}

	\chapter{Couplings on Tall Rectangles}
	
	\label{RectangleLCouple}

	Although the proofs of \Cref{p:globallaw2} and \Cref{p:closerho0}, indicating that a line ensemble $\bm{\mathcal{L}}$ satisfying \Cref{l0} likely satisfies the global law and regular profile events, will appear in \Cref{GlobalRegular} below, let us briefly mention one aspect of them. They will proceed by first restricting $\bm{\mathcal{L}}$ to a tall rectangle; this gives rise to a family of non-intersecting Brownian bridges with lower boundary. However, many of our previous results (such as those appearing in \Cref{Limit0} and \Cref{EDGESHAPE} for limit shapes) analyzed non-intersecting Brownian bridges without lower boundary. Thus, we will require a coupling that compares a family of non-intersecting Brownian bridges on a tall rectangle with lower boundary to one with the same starting and ending data but without a lower boundary; in this way, it ``removes'' the lower boundary condition of the first family, so we sometimes refer to it as a ``boundary removal coupling.'' The purpose of this chapter is to provide such a coupling, which will be stated as \Cref{c:finalcouple} in \Cref{s:couplingandholder} below.

	We begin in \Cref{s:linelong} by establishing several miscellaneous concentration estimates. We then state and prove the  boundary removal coupling in \Cref{RectangleCouple}, assuming the existence of particular ``preliminary couplings'' and certain improvements of the H\"{o}lder regularity bounds (from \Cref{eventtsregular2} and \Cref{sclprobability}). The former will be verified in \Cref{Couple0Proof} and the latter in \Cref{RegularImproved}.

	\section{Concentration Bounds and Extreme Path Estimates}
	
	\label{s:linelong}
	
	In this section we collect several results that will be used to establish the existence of the boundary removal coupling later in this chapter. These include concentration bounds for non-intersecting Brownian bridges in \Cref{SmoothRho} below (which will be proven in \Cref{ProofrhoRandom} and \Cref{ProofRhoDeterministic}), and estimates for the locations of the extreme paths of these bridges in \Cref{Processn} below. 
	
	\subsection{Concentration Around Smooth Profiles}
	
	\label{SmoothRho} 
	
	In this section we state several results indicating that non-intersecting Brownian bridges concentrate around smooth profiles. We begin with the following assumption, indicating that the boundary data for these bridges is ``on-scale'' (analogously to the $\textbf{MED}$ event in \Cref{eventsregular1}).
	
	\begin{assumption}

		\label{a:boundary} 
		
		Fix integers $k, n > 1$, and real numbers $D >1$ and $L \in [1, k^D]$, such that $n=L^{3/2} k$. Further let $A > 0$, $B \ge 2A^{-1}$, and $t \in [B^{-1}, A - B^{-1}]$ be real numbers; set $\mathsf{t} = tk^{1/3}$; and let $\bm{u}, \bm{v} \in \mathbb{W}_n$ be $n$-tuples. Assume that, for each $j \in \llbracket 1, n \rrbracket$, we have 
		\begin{align}\label{e:xidensity}
			-Bk^{2/3} - Bj^{2/3} \leq u_j \le Bk^{2/3} -B^{-1} j^{2/3}; \quad -Bk^{2/3} -Bj^{2/3} \le v_j \leq Bk^{2/3} -B^{-1} j^{2/3}.   
		\end{align}
		
		\noindent  Sample $n$ non-intersecting Brownian bridges $\bm{\mathsf{x}} = (\mathsf{x}_1, \mathsf{x}_2, \ldots , \mathsf{x}_n) \in \llbracket 1, n \rrbracket \times \mathcal{C} \big( [0, Ak^{1/3}] \big)$ from the measure $\mathsf{Q}^{\bm{u}; \bm{v}}$. 
		
	\end{assumption}

	The following proposition indicates the existence of a random (that is, measurable with respect to $\bm{\mathsf{x}}$) measure $\mu_t$ satisfying the following properties. First, recalling the classical locations with respect to a measure from \Cref{gammarho}, $\mathsf{x}_j (\mathsf{t})$ is very close to the $j$-th classical location of $\mu_t$, which is of order $- (j/k)^{2/3}$. Second, $\mu_t$ admits a density $\varrho_t$ with respect to Lebesgue measure, which satisfies bounds similar to those imposed in and implied by \Cref{blx122} (as in \Cref{p:densityest}). Third, assuming an upper bound on the difference between the classical locations of $\mu_t$, the inverse of the cummulative density function for $\mu_t$ is smooth (as in \Cref{c:rhoderbound}). We establish the following proposition in \Cref{ProofrhoRandom} below.

	\begin{prop}\label{p:densityub0}
		
		Adopt \Cref{a:boundary}. There exist constants $c = c(A, B) > 0$, $C_1 = C_1 (A, B) > 1$, and $C_2 = C_2 (A, B, D) > 1$ such that, with probability at least $1 - C_2 e^{-c(\log n)^2}$, there exists a random measure $\mu_t \in \mathscr{P}_{\fin}$ satisfying $\supp \mu_t \subseteq [-C_1 L, C_1 L^{3/4}]$, $\mu_t (\mathbb{R}) = L^{3/2}$, and the following three properties. In what follows, we denote the classical locations (recall \Cref{gammarho})  of $\mu_t$ by $\gamma_j = \gamma_{j;n}^{\mu_t}$, for each $j \in \llbracket 1, n \rrbracket$; we also define the function $\gamma : [0, L^{3/2}] \rightarrow \mathbb{R}$ by for each $y \in [0, L^{3/2}]$ setting
		\begin{align}\label{e:defgammataucopy}
			\gamma(y)=\displaystyle\sup \left\{ x \in \mathbb{R} : \int_{x}^\infty \mu_t (du) \geq y\right\}.
		\end{align}
		
		\begin{enumerate}
			\item  We have
			\begin{align}\label{e:locallaw}
				\gamma_{j+ \lfloor (\log n)^6 \rfloor} - n^{-D} \leq k^{-2/3} \cdot \sfx_j (\sft) \leq \gamma_{j- \lfloor (\log n)^6 \rfloor} + n^{-D}, \qquad \text{for each $j \in \llbracket 1, n \rrbracket$},
			\end{align}
			
			\noindent and 
			\begin{align}\label{e:gammaest0}
				-C_1 \Big( \displaystyle\frac{j}{k} \Big)^{2/3} - C_1 \leq \gamma_j \leq C_1 - C_1^{-1} \Big( \displaystyle\frac{j}{k} \Big)^{2/3}, \qquad \text{for each $j \in \big\llbracket (\log n)^6, n \big\rrbracket$}.
			\end{align}
			\item The measure $\mu_t$ has a density $\varrho_t : \mathbb{R} \rightarrow \mathbb{R}_{\ge 0}$ with respect to Lebesgue measure, satisfying
			\begin{align}\label{e:rhotint}
				\int_{x}^\infty \varrho_t(y) d y \leq C _1|x|^{3/2},\quad \text{for any $x \le -1$}; \qquad \int_{C_1}^\infty \varrho_t(y) d y\leq C_1 k^{-1} (\log n)^6,
			\end{align}
			and
			\begin{align}\label{e:rhotup}
				\varrho_t(x)\leq C_1  \max\{1,-x\}^{3/4}, \qquad \text{for any $x \in \mathbb{R}$}.
			\end{align}
			\item For any integer $\ell \ge 1$ and real number $R > 1$, there exists a constant $C_3 = C_3 (\ell, A, B, R) > 1$ such that the following holds. If for any $y, y' \in [B^{-1}, B]$, with $y'-y\geq 10 k^{-1} (\log n)^{50}$, we have $\big|\gamma(y)-\gamma(y') \big|\leq R|y-y'|$, then $\gamma \in \mathcal{C}^{\ell} \big( [2/B, B/2] \big)$ and
			\begin{align}\label{e:ldercopy}
				\|\gamma\|_{\mathcal{C}^{\ell} ([2/B, B/2])} \leq C_3.  
			\end{align}
			
		\end{enumerate}
	\end{prop}

	The following corollary, to be established in \Cref{ProofRhoDeterministic} below, is a variant of \Cref{p:densityub0} that makes the measure $\mu_t$ deterministic but provides the weaker concentration bound \eqref{e:determine}.
	
	\begin{cor}
		
		\label{p:concentration}
		
		Adopt \Cref{a:boundary}. There exist three constants $c = c(A, B) > 0$, $C_1 = C_1 (A, B) > 1$, and $C_2 = C_2 (A, B, D) > 1$, and a deterministic measure $\mu_t \in \mathscr{P}_{\fin}$, such that $\mu_t (\mathbb{R}) = L^{3/2}$ and the following holds if $n > C_2$. Below, we denote the classical locations (recall \Cref{gammarho})  of $\mu_t$ by $\gamma_j = \gamma_{j;n}^{\mu_t}$ and set $\mathfrak{m}_j = \big\lceil C_1 \log n \cdot \max \{j ^{1/2}, k^{1/2} \} \big\rceil$ for each $j \in \llbracket 1, n \rrbracket$.

		\begin{enumerate}
			\item We have $\supp \mu_t \subseteq [-C_1 L, C_1 L^{3/4}]$, and $\mu_t$ admits a density $\varrho_t: \mathbb{R} \rightarrow \mathbb{R}_{\ge 0}$ with respect to Lebesgue measure that satisfies \eqref{e:rhotint} and \eqref{e:rhotup}.
			\item The bound	\eqref{e:gammaest0} holds for each $j \in \big\llbracket (\log n)^6 + 1, n \big\rrbracket$, and we have 
			\begin{align}\label{e:determine}
				\mathbb{P} \Bigg[ \bigcap_{j = 1}^n \big\{ \gamma_{j + \mathfrak{m}_j}  \leq  k^{-2/3} \cdot \sfx_j (\sft)  \leq \gamma_{j-\mathfrak{m}_j} \big\} \Bigg] \ge 1 - C_2 e^{-c (\log n)^2}. 
			\end{align}
		\end{enumerate}
		
	\end{cor}

	\subsection{Approximation by Random Profiles} 
	
	\label{ProofrhoRandom}

	In this section we establish \Cref{p:densityub0}, whose notation we adopt throughout. We will assume that $A = 1$ and $u_1 = v_1 = 0$ (as we may, the former by the scaling invariance \Cref{scale} and the latter by the affine invariance \Cref{linear}); we will also assume (by replacing $B$ with $B + 10$, if necessary) that $B > 10$. We begin with the following lemma bounding the $\mathsf{x}_j (\mathsf{t})$ with high probability. Set $M = B + 9B^3\pi^2/64 + 1$ and define the event $\mathscr{E} = \mathscr{E}_1 \cap \mathscr{E}_2$, where
	\begin{flalign}
		\label{e1e2gamma} 
		\mathscr{E}_1 = \bigcap_{j=1}^n \big\{ \mathsf{x}_j (\mathsf{t}) \le M k^{2/3} - B^{-1} j^{2/3} \big\}; \qquad \mathscr{E}_2 = \bigcap_{j=1}^n \big\{ \mathsf{x}_j (\mathsf{t}) \ge -M k^{2/3} -2Bj^{2/3}  \big\}.
	\end{flalign}
	
	\begin{lem} 
		
		\label{eprobabilityuv}
		
		There are constants $c = c(B) > 0$ and $C = C(B, D) > 1$ with $\mathbb{P} [\mathscr{E}] \ge 1 - C e^{-c(\log n)^2}$. 
	\end{lem} 
	
	\begin{proof} 
		
		This will follow from \Cref{p:compareAiry}. In particular, apply the first part of that lemma, with the $(f; a, b)$ there equal to $(-\infty; 0, k^{1/3})$ here and the $(d, M, D)$ there equal to $(B^{-1}, Bk^{2/3}, 1)$ here. Its assumptions are verified by the upper bounds in \eqref{e:xidensity}, and so it yields constants $c_1 = c_1 (B) > 0$ and $C_1 = C_1 (B, D) > 0$ such that 
		\begin{flalign}
			\label{e10} 	
			\mathbb{P} [\mathscr{E}_1] \ge \mathbb{P} \Bigg[ \bigcap_{j=1}^n \bigg\{ \mathsf{x}_j (\mathsf{t}) \le \Big( B + \displaystyle\frac{9 B^3 \pi^2}{64} \Big) k^{2/3} - B^{-1} j^{2/3} + 2 (\log n)^2 \bigg\} \Bigg]  \ge  1 - C_1 e^{-c_1 (\log n)^2}, 
		\end{flalign}
		
		\noindent where in the first inequality we used the definition of $M$ and the fact that $2(\log n)^2 \le k^{2/3}$ for sufficiently large $n$ (as $k^{3D/2+1} \ge L^{3/2} k = n$). Next, apply the second part of \Cref{p:compareAiry}, with the $(a, b)$ there equal to $(0, k^{1/3})$ here and the $(A, B, M)$ there equal to $( 1, B, Bk^{2/3})$ here. Its assumptions are verified by the lower bounds in \eqref{e:xidensity}, and so it yields constants $c_2 = c_2 (B) > 0$ and $C_2 = C_2 (B, D) > 1$ such that denoting $A_0 = B + 5 \le 2B$ (as $B \ge 10$) we have
		\begin{flalign*}
			\mathbb{P} [\mathscr{E}_2] \ge	\mathbb{P} \Bigg[ \bigcap_{j=1}^n \bigg\{ \mathsf{x}_j (\mathsf{t}) \ge \Big( \displaystyle\frac{9 \pi^2}{16 A_0^3} t(1-t) - B \Big) k^{2/3} - 2 (\log n)^2 - A_0 j^{2/3} \bigg\} \Bigg] \ge 1 - C_2 e^{-c_2 (\log n)^2},
		\end{flalign*} 
		
		\noindent where in the first inequality we again used the definition of $M$ and the fact that $2 (\log n)^2 \le k^{2/3}$ (and that $t(1-t) \ge 0$). This, together with \eqref{e10} and a union bound, yields the lemma.
	\end{proof}

	Next, we apply \Cref{tuvwx} to equate the law of $\bm{\mathsf{x}}(\mathsf{t})$ with Dyson Brownian motion run under certain (random) initial data. More specifically, recalling the notation from \Cref{MotionCurves}, define the $n\times n$ diagonal matrices $\bm{U} = \diag (\bm{U})$ and $\bm{V} = \diag (\bm{v})$; let $\bm{W}$ denote a random $n \times n$ unitary matrix with law \eqref{wuv}; and define the random Hermitian $n \times n$ matrix
	\begin{align}\label{e:lawycopy}
		\bm A =(1-t) \cdot \bm U+ t \cdot \bm W \bm V \bm W^*
	\end{align}
	
	\noindent Set $\tau = t(1-t)$, and denote the eigenvalues of $\bm{A}$ by $\eig (\bm{A}) = \bm{\mathsf{a}} = (\mathsf{a}_1, \mathsf{a}_2, \ldots , \mathsf{a}_n) \in \overline{\mathbb{W}}_n$. By \Cref{tuvwx} (with the $(t, \mathsf{T})$ there given by $(\mathsf{t}, k^{1/3})$ here) and the fact that $\mathsf{t} = tk^{1/3}$, $\bm{\mathsf{x}} (\mathsf{t})$ has the same law as $\bm{\lambda} (\tau k^{1/3})$, where $\bm{\lambda} (s)$ is Dyson Brownian motion with initial data $\bm{\lambda}(0) = \bm{\mathsf{a}}$, run for time $s$. 
	
	Since \eqref{e:xidensity} (and the fact that $n = L^{3/2} k$) implies that $-2B L k^{2/3}\leq -Bn^{2/3}-Bk^{2/3}\leq u_n \leq u_1\leq 0$ and $-2B L k^{2/3} \le v_n \le v_1 \le 0$, the Weyl interlacing inequality yields $-4B Lk^{2/3}\leq \min \bm{\mathsf{a}} \le \max \bm{\mathsf{a}} \leq 0$. We then set (recalling the notation from \eqref{aemp})
	\begin{flalign}
		\label{snu1} 
		\nu=L^{3/2} \cdot \emp (k^{-2/3} \cdot \bm{\mathsf{a}}) = \frac{1}{k}\sum_{j =1}^n \delta_{\mathsf{a}_j / k^{2/3}}, \quad \text{so} \quad \nu(\mathbb{R}) = L^{3/2} \quad \text{and} \quad \supp \nu \subseteq  [-4B L, 0].
	\end{flalign} 
	
	\noindent Recalling the notation on free convolutions from \Cref{TransformConvolution}, for any real number $s\geq 0$, let $\nu_s = \nu \boxplus \mu_{\semci}^{(s)} \in \mathscr{P}_{\fin}$. Denote the classical locations (recall \Cref{gammarho}) of $\nu_s$ by $\gamma_j (s) = \gamma_{j;n}^{\nu_s}$. 
	
	The following lemma indicates that the $\mathsf{x}_j (\mathsf{t})$ concentrate around these classical locations.		
	
	\begin{lem}
		
		\label{xgamma0}
		
		There exists a constant $C = C(D) > 1$ and a coupling between $\bm{\mathsf{x}}$ and $\nu_{\tau}$ such that
		\begin{flalign*}
			\mathbb{P} \Bigg[ \bigcap_{j = 1}^n \big\{ \gamma_{j+\lfloor (\log n)^6 \rfloor}(\tau)-n^{-50D}\leq k^{-2/3} \cdot \sfx_i(\sft) \leq \gamma_{j- \lfloor (\log n)^6 \rfloor}(\tau)+n^{-50D} \big\} \Bigg] \ge 1 - C e^{ - (\log n)^2}.
		\end{flalign*}
		
	\end{lem}

	\begin{proof}
		
		Recall by \Cref{tuvwx} that $\bm{\mathsf{x}}(\mathsf{t})$ has the same law as $\bm{\lambda} (\tau k^{1/3})$. So, it suffices to show 
		\begin{flalign}
			\label{e:locallawcopy2}
			\begin{aligned} 
			\mathbb{P} \Bigg[ \bigcap_{j = 1}^n \big\{  \gamma_{i+ \lfloor (\log n)^6 \rfloor }(\tau)-n^{-50D}\leq k^{-2/3} \cdot \lambda_i(\tau k^{1/3}) \leq  \gamma_{i- \lfloor (\log n)^6 \rfloor}(& \tau)+n^{-50D} \big\} \Bigg] \\
			&  \ge 1 - C e^{ - (\log n)^2},
			\end{aligned}
		\end{flalign}
		
		\noindent which will follow from \Cref{concentrationequation} and rescaling. More specifically, define for any $s \ge 0$ the probability measures 
		\begin{flalign}
			\label{snu2}
			\widetilde{\nu} = \emp (n^{-1} \cdot \bm{\mathsf{a}}), \qquad \text{and} \qquad  \widetilde{\nu}_s = \widetilde{\nu} \boxplus \mu_{\semci}^{(s)},
		\end{flalign}
		
		\noindent and denote the classical locations of $\widetilde{\nu}_s$ by $\widetilde{\gamma}_j (s) = \gamma_{j;n}^{\tilde{\nu}_s}$. By \Cref{concentrationequation}, there exists a constant $C = C(D) > 1$ such that  
		\begin{flalign}
			\label{e:tildegammarigidity} 
			\begin{aligned}
				\mathbb{P} \Bigg[ \bigcap_{j = 1}^n \big\{&  \widetilde \gamma_{j+ \lfloor (\log n)^6 \rfloor}( \tau k^{1/3} n^{-1})-n^{-50D-1} \\
				& \qquad \le n^{-1} \cdot \lambda_j (\tau k^{1/3}) \le \widetilde\gamma_{j - \lfloor (\log n)^6 \rfloor}(\tau k^{1/3} n^{-1})+n^{-50D-1} \big\} \Bigg] \ge 1 - C e^{ - (\log n)^2}. 
			\end{aligned}
		\end{flalign}
		
		Comparing \eqref{snu1} with \eqref{snu2}, we have for any interval $I \subseteq \mathbb{R}$ that $\nu_0(I)= L^{3/2} \cdot \widetilde\nu_0( n^{-1} k^{2/3} \cdot I)$. By the scaling relations for free convolutions given by \Cref{mtscale} (with the $A$ there equal to $L^{-3/2}$ here) and \Cref{mtscalebeta} (with the $\beta$ there equal to $k^{-1/3} n$ here), we have $\nu_s(I)= L^{3/2} \cdot \widetilde\nu_{k^{1/3} s / n}(n^{-1} k^{2/3} \cdot I)$ for any real number $s\geq 0$ (as for $\beta = k^{-1/3} n$ we would have $\beta^{-1/2} \cdot L^{-3/4} = n^{-1} k^{2/3}$, since $n = L^{3/2} k$). By \Cref{gammarho}, the classical locations therefore satisfy $\gamma_j(s)=n k^{-2/3} \cdot \widetilde\gamma_j (s k^{1/3} n^{-1})$. This, together with \eqref{e:tildegammarigidity}, implies \eqref{e:locallawcopy2} and thus the lemma.						
	\end{proof} 
	
	Now we can establish \Cref{p:densityub0}.

	\begin{proof}[Proof of \Cref{p:densityub0}]
		
		Recalling that $\tau = t(1-t)$, set $\mu_t = \nu_{\tau} = \nu \boxplus \mu_{\semci}^{(\tau)}$, and denote $\gamma_j = \gamma_j (\tau)$; by \Cref{mtscale} and \eqref{snu1}, we have $\mu_t (\mathbb{R}) = \nu(\mathbb{R}) = L^{3/2}$. Moreover, as explained below \Cref{mz}, $\mu_t$ admits a density $\varrho_t$ with respect to Lebesgue measure. Observe by \eqref{snu1} that $\varrho_t$ satisfies \Cref{blx12} (with the $B$ there equal to $4B$ here), so the first statement in  \Cref{p:densityest} implies that there exists a constant $C_3 = C_3 (B) > 1$ such that $\supp \mu_t \subseteq [-C_3 L, C_3 L^{3/4}]$. 
		
		Let us next verify the first statement of the proposition. Observe that \Cref{xgamma0} implies the bound \eqref{e:locallaw}. By \Cref{eprobabilityuv}, the fact that $\mathscr{E} = \mathscr{E}_1 \cap \mathscr{E}_2$, \eqref{e1e2gamma}, \Cref{xgamma0}, and a union bound, there exist constants $c_1 = c_1 (B) > 0$ and $C_0 = C_0 (B, D) > 1$ such that
		\begin{flalign*}
			\mathbb{P} \Bigg[ & \bigcap_{j=\lfloor (\log n)^6 \rfloor}^{n - \lfloor (\log n)^6 \rfloor} \big\{ k^{2/3} \gamma_{j - \lfloor (\log n)^6 \rfloor} - n^{-10} \le M k^{2/3} - B^{-1} j^{2/3} \big\} \\
			& \qquad \qquad \qquad \bigcap \big\{ k^{2/3} \gamma_{j + \lfloor (\log n)^6 \rfloor} + n^{-10} \ge -M k^{2/3} - 2 Bj^{2/3} \big\} \Bigg] \ge 1 - C_0 e^{-c_1 (\log n)^2}.
		\end{flalign*}
		
		\noindent This, with the facts that $Mk^{2/3} + n^{-10} - B^{-1} \big( j +  (\log n)^6  \big)^{2/3} \le 2Mk^{2/3} - (3B)^{-1} j^{2/3}$ and $-Mk^{2/3} - 2 B\big( j +  (\log n)^6  \big)^{2/3} - n^{-10} \ge -2Mk^{2/3} - (3B) j^{2/3}$ for sufficiently large $k$ (as $n = Lk \in [k, k^{D+1}]$), yields (after decreasing $c_1 = c_1 (D) > 0$ and increasing $C_0 = C_0 (B, D) > 1$ if necessary) that $\mathbb{P}[\mathscr{E}_0] \ge 1 - C_0 e^{-c_1 (\log n)^2}$, where
		\begin{flalign}
			\label{e0gamma}
			\mathscr{E}_0 =  \bigcap_{j = \lfloor (\log n)^6 \rfloor}^{n - \lfloor (\log n)^6 \rfloor} \bigg\{ -2M - (3B) \Big( \displaystyle\frac{j}{k} \Big)^{2/3} \le \gamma_j \le 2M - (3B)^{-1} \Big( \displaystyle\frac{j}{k} \Big)^{2/3} \bigg\}.
		\end{flalign}
		
		\noindent This confirms \eqref{e:gammaest0} for $j \in \big\llbracket  (\log n)^6 , n -  (\log n)^6  \big\rrbracket$ with the $C_1$ there equal to $C_4 = \max \{ 2M, 3B \}$ here. The fact that it also holds for $j \in \big\llbracket n -  (\log n)^6 , n \big\rrbracket$ follows from the fact that for such $j$ we have $\gamma_j \ge \gamma_n \ge \inf \supp \mu_t \ge -C_3 L \ge 2 C_3 (k^{-1} j)^{2/3}$, establishing the first statement of the proposition.
		
		We next establish the second, to which end we restrict to the event $\mathscr{E}_0$ for the remainder of this proof. To show the first bound in \eqref{e:rhotint}, fix a real number $x\leq -1$ as stated there; we may assume that $L \ge 4C_4^2$, as otherwise $\mu_t (\mathbb{R}) = L^{3/2} \le 8C_4^3 \le 8C_4^3 |x|^{3/2}$, and that $x \ge -(4C_4)^{-1} L$, as otherwise $\mu_t (\mathbb{R}) = L^{3/2} \le 8C_4^{3/2} |x|^{3/2}$; this verifies the first estimate in \eqref{e:rhotint} in both cases. Then, let $j_0 = j_0 (x) \in \big\llbracket (\log n)^6, n - (\log n)^6 \big\rrbracket$ denote the smallest integer such that $x > C_4 - C_4^{-1} (j_0/k)^{2/3}$, implying by \eqref{e:gammaest0} that $x > \gamma_{j_0}$; observe that such an integer $j_0$ exists, since $x \ge -(4C_4)^{-1} L \ge C_4 - L / (2C_4) \ge C_4 - C_4^{-1} \big( k^{-1} (n- (\log n)^6 )  \big)^{2/3}$ for sufficiently large $n$. This yields
		\begin{align}\label{e:d1}
			\int_x^\infty\varrho_t(x) dx \leq \displaystyle\frac{2j_0 - 1}{2n} \cdot L^{3/2} \le \frac{j_0}{k} \le 2 \big(C_4(C_4-x) \big)^{3/2}\leq 2 \big( C_4 (C_4 +1) \big)^{3/2}|x|^{3/2}\leq 8 C_4^3 |x|^{3/2},
		\end{align}
		
		\noindent where in the first statement we used the fact that $x > \gamma_{j_0}$; in the second we used the fact that $n = L^{3/2} k$; in the third we used the fact that $k^{-2/3} (j_0-1)^{2/3} \le C_4 (C_4 - x) \le k^{-2/3} j_0^{2/3}$ (unless $j_0 \le  (\log n)^6  + 1$, in which case $j_0 \le k$ and so $j_0 k^{-1} \le 1 \le 2 \big( C_4 (C_4 - x) \big)^{3/2}$ again holds); in the fourth used the fact that $|C_4 - x| \le (C_4 + 1) |x|$ (as $x \le -1$); and in the fifth we used the fact that $(C_4 + 1)^{3/2} \le 4 C_4^{3/2}$ (as $C_4 \ge 1$). This confirms the first bound in \eqref{e:rhotint}. Further observe on $\mathscr{E}_0$ that $\gamma_{\lfloor (\log n)^6 \rfloor} \le 2M$, and so very similar reasoning as implemented to deduce \eqref{e:d1} (using $2M$ in place of $x$ there) yields
		\begin{flalign*}
			\displaystyle\int_{2M}^{\infty} \varrho_t (x) dx \le k^{-1} (\log n)^6,
		\end{flalign*}
		
		\noindent verifying the second statement of \eqref{e:rhotint}.
		
		The remaining parts of the lemma will follow from applying \Cref{p:densityest} and \Cref{c:rhoderbound} to the measure $\mu_t =\nu\boxplus \mu_{\semci}^{(\tau)}$. Restricting to $\mathscr{E}_0$, \eqref{e:d1} holds, verifying \Cref{blx122}. Thus, the second part of \Cref{p:densityest} (using the fact that $\nu$ satisfies \Cref{blx12} by \eqref{snu1}) yields \eqref{e:rhotup}, proving the second part of the proposition.    
		
		To show the third, we apply \Cref{c:rhoderbound}. By \eqref{e0gamma}, we have that $\gamma (B)\geq \gamma_{\lceil Bk\rceil} \geq -3(B+M)$, which verifies the first assumption in \Cref{c:rhoderbound}, with the $A$ there equal to $3(B+M)$ here. The second follows from the condition imposed in the third part of \Cref{p:densityub0}, with the $A$ there equal to $R$ here (if $k$ is sufficiently large so that the $\varepsilon$ of \Cref{c:rhoderbound} is less than $10k^{-1} (\log n)^{50}$ here). Thus, \Cref{c:rhoderbound} applies and shows (together with the fact that $\| \gamma \|_{\mathcal{C}^0 ([2/B, B/2])}$ is uniformly bounded, by \eqref{e:gammaest0}) the third part of the proposition.
	\end{proof}

	\subsection{Approximation by Deterministic Profiles} 
	
	\label{ProofRhoDeterministic}
	
	In this section we establish \Cref{p:concentration}, which will follow from \Cref{p:densityub0} together with \Cref{concentrationbridge}.

	\begin{proof}[Proof of \Cref{p:concentration}]
		
		Throughout this proof, we will assume (by replacing $B$ by $B + 10$, if necessary) that $B > 10$. First observe by  \Cref{p:densityub0} that there exist constants $c_1 = c_1 (A, B) \in (0, 1)$, $C_3 = C_3 (A, B) > 1$, and $C_4 = C_4 (A, B, D) > 1$, and an event $\mathscr{E}_0$ with $\mathbb{P} \big[ \mathscr{E}_0^{\complement} \big] \ge 1 - C_4 e^{-c_1 (\log n)^2}$, such that on $\mathscr{E}_0$ there exists a random measure $\widetilde{\mu}_t$ satisfying $\supp \widetilde{\mu}_t \subseteq [-C_3 L, C_3 L^{3/4}]$; $\widetilde{\mu}_t (\mathbb{R}) = L^{3/2}$; and the following two properties. First, denoting the classical locations of $\widetilde{\mu}_t$ by $\widetilde{\gamma}_j = \gamma_{j;n}^{\tilde{\mu}_t}$, \eqref{e:locallaw} and \eqref{e:gammaest0} both hold (with $\gamma_j$ there replaced by $\widetilde{\gamma}_j$ here). Second, $\widetilde{\mu}_t$ has a density $\widetilde{\varrho}_t \in L^1 (\mathbb{R})$ with respect to Lebesgue measure, such that \eqref{e:rhotint} and \eqref{e:rhotup} hold (with $\varrho_t$ there replaced by $\widetilde{\varrho}_t$ here). 
		
		Now define $\mu_t \in \mathscr{P}_{\fin}$ by setting $\mu_t (dx) = \varrho_t (x) dx$, where $\varrho_t : \mathbb{R} \rightarrow \mathbb{R}_{\ge 0}$ is given by the conditional expectation (recalling below that $\mathbb{P} \big[ \mathscr{E}_0^{\complement} \big] \ge 1 - C_4 e^{-c_1 (\log n)^2} > 0$ for sufficiently large $n$)
		\begin{flalign}
			\label{rho1x} 
			\varrho_t(x)=\bE \big[  \widetilde \varrho_t(x) \mid \mathscr{E}_0 \big], \qquad \text{for each $x \in \mathbb{R}$}. 
		\end{flalign}
		
		\noindent Observe that $\mu_t (\mathbb{R}) = L^{3/2}$, since $\widetilde{\mu}_t (\mathbb{R}) = L^{3/2}$. We moreover claim that $\mu_t$ satisfies \eqref{e:rhotint} and \eqref{e:rhotup}. Indeed, for any $x \le -1$ and sufficiently large $n$, we have
		\begin{flalign*}
			\displaystyle\int_x^{\infty} \varrho_t (x) dx \le C_3 |x|^{3/2},
		\end{flalign*}      
		
		\noindent for some constant $C_5 = C_5 (A, B) > 1$, by \eqref{rho1x} and the fact that $\widetilde{\varrho}_t$ satisfies \eqref{e:rhotint}. This verifies the first bound in \eqref{e:rhotint}. The proof of the second is entirely analogous, as is that of \eqref{e:rhotup}, so they are omitted.
		
		We next verify that the $\gamma_j$ satisfy \eqref{e:gammaest0} if $j \ge (\log n)^6 + 1$. Denoting $\gamma_j^- = -C_3 (j / k)^{2/3} - C_3$ and $\gamma_j^+ = C_3 - C_3^{-1} ( j / k )^{2/3}$ for each integer $j \in \llbracket 1, n \rrbracket$, we have, for $j \in \big\llbracket (\log n)^6 + 1, n - (\log n)^6 - 1 \big\rrbracket$,
		\begin{flalign*}
			\displaystyle\int_{\gamma_j^+}^{\infty} \varrho_t (x) dx & \le \displaystyle\int_{\gamma_j^+}^{\infty} \mathbb{E} \big[ \widetilde{\varrho}_t (x) \mid \mathscr{E}_0 \big] dx  \le  \displaystyle\frac{2j-1}{2n} \cdot L^{3/2},
		\end{flalign*}
		
		\noindent where in the first inequality we used \eqref{rho1x}; in the second we used the fact that $\widetilde{\gamma}_j$ satisfies \eqref{e:gammaest0} (with $C_1$ there replaced by $C_3$) and \Cref{gammarho}. This, together with \Cref{gammarho}, implies that $\gamma_j \le \gamma_j^+$. Since $\gamma_j^+ = C_3 - C_3^{-1} ( j / k )^{2/3}$, this shows that $\gamma_j$ satisfies the upper bound in \eqref{e:gammaest0} (with the $C_1$ there equal to $C_3$ here); the proof of the lower bound is entirely analogous and thus omitted.
		
		It therefore remains to verify \eqref{e:determine}; in what follows, we recall from \Cref{htw} the height function $\mathsf{H}^{\bm{\mathsf{x}}}$ associated with the line ensemble $\bm{\mathsf{x}}$. Observe by \eqref{e:locallaw}, \Cref{gammarho} (with the fact that $n^{-1} \cdot \widetilde{\mu}_t (\mathbb{R}) = n^{-1} L^{3/2} = k^{-1}$), and \eqref{e:rhotup} (with the facts that $\supp \widetilde{\varrho}_t \subseteq [-C_3 L, C_3 L^{3/4}]$, that $L \le n$, and that $B > 10$) that, on $\mathscr{E}_0$, for any $x\in \mathbb{R}$ we have 
		\begin{flalign}
			\label{hx4} 
			\begin{aligned} 
			\mathsf{H}^{\bm{\mathsf{x}}} ( \mathsf{t}, k^{2/3} x )  & \le  k\int_{x - n^{-D}}^\infty\widetilde \varrho_t(y) d y +  (\log n)^6 \\
			& \le k \displaystyle\int_x^{\infty} \widetilde{\varrho}_t (y) dy + C_3 (C_3 L)^{3/4} n^{-D} + (\log n)^6 \le k \displaystyle\int_x^{\infty} \widetilde{\varrho}_t (y) dy + 2 (\log n)^6,
			\end{aligned} 
		\end{flalign} 
		
		\noindent and similarly 
		\begin{flalign*}
			& \mathsf{H}^{\bm{\mathsf{x}}} (\mathsf{t}, k^{2/3} x)  \ge k \displaystyle\int_x^{\infty} \widetilde{\varrho}_t (y) dy - 2 (\log n)^6.
		\end{flalign*}
		
		\noindent Taking expectations, we deduce  
		\begin{align}\label{23kx}
			\Bigg|\bE \big[ \mathsf{H}^{\bm{\mathsf{x}}} ( \mathsf{t}, k^{2/3} x) \big] - k \int_x^\infty \mathbb{E} \big[  \widetilde{\varrho}_t (y) \mid \mathscr{E}_0 \big] d y \Bigg|\leq 2 (\log n)^6 + n \cdot \mathbb{P} \big[ \mathscr{E}_0^{\complement} \big] \le 3 (\log n)^6,
		\end{align}
		
		\noindent where in the last bound we used the facts that $\mathsf{H}^{\bm{\mathsf{x}}} (\mathsf{t}, k^{2/3} x) \le n$, that $\mathbb{P} \big[ \mathscr{E}_0^{\complement} \big] \le C_4 e^{-c_1 (\log n)^2}$, and that $n$ is sufficiently large.
		
		We next define the event 
		\begin{flalign*}
			\mathscr{F}= \Big\{ \mathsf{H}^{\bm{\mathsf{x}}} (\mathsf{t}, k^{2/3} x) \le 2C_3 k\big( |x| + 1 \big)^{3/2} \text{ for all } x\Big\}.
		\end{flalign*}
		
		\noindent Since $\widetilde{\varrho}_t$ satisfies \eqref{e:rhotint} (and $3 (\log n)^6 \le C_3 k$, as $n = L^{3/2} k$ and $L \le k^D$), \eqref{hx4} implies that $\mathscr{E}_0 \subseteq \mathscr{F}$ for $n$ sufficiently large; in particular, $\bP \big[ \mathscr{F}^\complement \big] \le \mathbb{P} \big[ \mathscr{E}_0^{\complement} \big] \leq C_4 e^{-c_1 (\log n)^2}$. Thus, \Cref{concentrationbridge} (with $(f,g; w; B; r)$ there equal to $\big(-\infty, \infty; k^{2/3} x; 2C_3 k (|x|+1)^{3/2}; 2\log n \big)$) yields a deterministic number $\mathfrak{Y} = \mathfrak{Y} (\bm{u}; \bm{v}; k; t; x; B) \in \mathbb{R}$ such that 
		\begin{align}\label{e:heightconcentration}
			\bP \bigg[ \big| \mathsf{H}^{\bm{\mathsf{x}}} (\mathsf{t}, k^{2/3} x) -\mathfrak{Y} \big|\geq \big( |x| + 1 \big)^{3/4} (8C_3 k)^{1/2} \log n \bigg] & \leq 4 e^{-(\log n)^2} + C_4 e^{-c_1 (\log n)^2} \le 2C_4 e^{-c_1 (\log n)^2}.
		\end{align}
		
		\noindent Thus, 
		\begin{align}\label{e:Yvalue}
			\begin{aligned}
				\Bigg| \mathfrak{Y} -k\int_x^\infty \varrho_t(y)d y \Bigg| & \le \mathbb{E} \Big[ \big| \mathfrak{Y} - \mathsf{H}^{\bm{\mathsf{x}}} (\mathsf{t}, k^{2/3} x)  \big| \Big] + \Bigg| \mathbb{E}[\mathsf{H}^{\bm{\mathsf{x}}} (\mathsf{t}, k^{2/3} x)] - k \displaystyle\int_x^{\infty} \mathbb{E} \big[  \widetilde{\varrho}_t (y) \mid \mathscr{E}_0 \big] dy \Bigg| \\
				& \leq \big( |x| + 1 \big)^{3/4} (8C_3 k)^{1/2} \log n + 2n C_4 e^{-c_1 (\log n)^2} + 3 (\log n)^6 \\ 
				& \le  3 C_3 k^{1/2} \big( |x| + 1 \big)^{3/4} \log n,
			\end{aligned}
		\end{align}
		
		\noindent where in the first bound we used \eqref{rho1x}; in the second bound we used \eqref{23kx}, \eqref{e:heightconcentration}, and the fact that $\mathsf{H}^{\bm{\mathsf{x}}} (t, k^{2/3} x) \le n$; and in the third we used the facts that $\log n \le (3D/2+1) \log k$ (as $n = L^{3/2} k \le k^{3D/2+1}$) and that $n$ is sufficiently large.

		By inserting \eqref{e:Yvalue} into \eqref{e:heightconcentration}, we get 
		\begin{flalign}
			\label{xfprobability} 
			\mathbb{P} \big[ \mathscr{F} (x)^{\complement} \big] \le 3C_4 e^{-c_1 (\log n)^2}
		\end{flalign}
		
		\noindent  where 
		\begin{flalign}
			\label{xfevent} 
			\mathscr{F} (x) = \Bigg\{ \bigg| \mathsf{H}^{\bm{\mathsf{x}}} (\mathsf{t}, k^{2/3} x)  - k\int_x^\infty \varrho_t(y)d y \bigg| \le 6C_3 k^{1/2} \big( |x| + 1 \big)^{3/4} \log n \Bigg\}.
		\end{flalign}
		
		\noindent  Fix some integer $j \in \big\llbracket (\log n)^6 + 1 , n - (\log n)^6 - 1 \rrbracket$. Then \eqref{e:gammaest0} (which holds for $\gamma_j$) yields a constant $C_6 = C_6 (A, B) > 1$ such that $\big( 1 + |\gamma_j| \big)^{3/4} \le C_6 \max \{ k^{-1/2} j^{1/2}, 1 \}$. Together with \eqref{xfevent} (and \Cref{gammarho}), this implies on $\mathscr{F} (\gamma_j)$ that  
		\begin{align}\begin{split}\label{e:heightcc2}
				\big| \mathsf{H}^{\bm{\mathsf{x}}} ( \mathsf{t}, k^{2/3} \gamma_j) - j \big| \le 9 C_3 k^{1/2} \big( 1+|\gamma_j| \big)^{3/4} \log n + j^{-1} \le 10 C_3 C_6 \log n \cdot \max \{ j^{1/2}, k^{1/2} \}.
		\end{split}\end{align}
		
		Setting $\mathfrak{m}_j = \big\lceil  10 C_3 C_6 \log n \cdot \max \{ j^{1/2}, k^{1/2} \} \big\rceil$ for each integer $j \in \llbracket 1, n \rrbracket$, \eqref{e:heightcc2} implies on $\mathscr{F} (\gamma_j)$ that $\mathsf{x}_{j +  \mathfrak{m}_j} \le k^{2/3} \gamma_j \le \mathsf{x}_{j - \mathfrak{m}_j}$ if $j \in \big\llbracket (\log n)^6, n - (\log n)^6 - 1 \big\rrbracket$. Hence, since $ \mathfrak{m}_j > (\log n)^6 + 1$ for sufficiently large $n$ (again, as $n = L^{3/2} k \le k^{3D/2+1}$), it follows that
		\begin{align*}
			\mathscr{F}_j \subseteq \big\{ \gamma_{j - \mathfrak{m}_j} \le k^{-2/3} \cdot \mathsf{x}_j (\mathsf{t}) \le \gamma_{j + \mathfrak{m}_j} \big\}, \qquad \text{for each $j \in \llbracket 1, n \big\rrbracket$}.
		\end{align*}
		
		\noindent This, together with \eqref{xfprobability} and a union bound over $j \in \llbracket 1 , n \rrbracket$, yields the corollary.	
	\end{proof}

	\subsection{Lowest Path Estimates}
	
	\label{Processn} 
	
	In this section we establish the following proposition, which will be used in the proof of \Cref{p:dentoReg} (through that of \Cref{l:holder}) below. It bounds how far the lowest path in a family of non-intersecting Brownian bridges can go down in a given interval of time. We will require that the boundary data for these non-intersecting Brownian bridges has some regularity, which is prescribed through \eqref{vntrho} and \eqref{e:gammabb}; observe their similarity to \Cref{p:concentration}. In what follows, we recall the classical locations with respect to a measure from \Cref{gammarho}.

	\begin{prop}
		
		\label{uvrho} 
		
		For any real numbers $A > 0$ and $B, D > 1$, there exist constants $c = c(A, B) > 0$, $C_1 = C_1 (A, B) > 1$, and $C_2 =C_2 (A, B, D) > 1$ such that the following holds. Let $k \ge 1$ and $n \in \llbracket k, k^D \rrbracket$ be integers; $\mathsf{T}, \mathfrak{m} \ge 1$ be real numbers; and $\bm{u}, \bm{v} \in \overline{\mathbb{W}}_n$ be $n$-tuples. Sample $n$ non-intersecting Brownian bridges $\bm{\mathsf{x}} = (\mathsf{x}_1, \mathsf{x}_2, \ldots , \mathsf{x}_n) \in \llbracket 1, n \rrbracket \times \mathcal{C} \big( [0, \mathsf{T}] \big)$ from the measure $\mathsf{Q}^{\bm{u}; \bm{v}}$. Let $\mu \in \mathscr{P}_{\fin}$ denote a measure with $\mu (\mathbb{R}) = k^{-1} n$, which admits a density $\varrho \in L^1 (\mathbb{R})$ with respect to Lebesgue measure. Assume that
		\begin{flalign}
			\label{vntrho} 
			 \mathsf{T} \in [ A n^{1/3}, 3A n^{1/3}]; \quad v_n - u_n \ge -B \Big( \displaystyle\frac{n}{k} \Big)^{1/6} n^{2/3}; \quad \sup_{x \in \mathbb{R}} \varrho(x) \leq B \Big( \displaystyle\frac{n}{k} \Big)^{1/2}; \quad \mathfrak{m} \le B n^{1/2} \log n.
		\end{flalign}
		
		\noindent Denoting the classical locations (recall \Cref{gammarho})  of $\mu$ by $\gamma_j = \gamma_{j;n}^{\mu}$ for each $j \in \llbracket 1, n \rrbracket$, further assume for some real number $M \ge 1$ that  
		\begin{align}\label{e:gammabb}
			\gamma_{j+\mathfrak{m}}-M\leq k^{-2/3} \cdot u_j \leq \gamma_{j-\mathfrak{m}}+M, \qquad \text{for each $j \in \llbracket 1, n \rrbracket$}.
		\end{align}
		
		\noindent Then, for any $t \in \big[0, (1 - B^{-1}) A \big]$, we have
		\begin{align}\label{e:xiholder}
			\begin{aligned}
				\mathbb{P} \Bigg[  \displaystyle\frac{\sfx_n(tk^{1/3} ) - \mathsf{x}_n (0)}{k^{2/3}} \le -C_1  \bigg(  & t \Big( \displaystyle\frac{n}{k} \Big)^{1/2}  \big| \log (A t^{-1}) \big|^2 + t^{1/2} M^{1/2} \Big( \displaystyle\frac{n}{k} \Big)^{1/4}  \bigg)  \Bigg]  \le C_2 e^{-c(\log n)^2}.
			\end{aligned}
		\end{align}
	\end{prop}

	\begin{proof}
		
		We will establish this proposition by using \Cref{p:extreme}, together with \Cref{tobridge} (to express Dyson Brownian motion through non-intersecting Brownian bridges ending at the same point). Throughout this proof, we will assume (by translation) that $u_n = 0$ and (by the scaling invariance \Cref{scale}) that $A = B^{-1}$.
		
		 Sample a family of $n$ non-intersecting Brownian bridges $\bm{\mathsf{y}} = (\mathsf{y}_1, \mathsf{y}_2, \ldots , \mathsf{y}_n) \in \llbracket 1, n \rrbracket \times \mathcal{C} \big( [0, \mathsf{T}] \big)$ from the measure $\mathsf{Q}^{\bm{u}; \bm{0}_n}$. By the second statement of \eqref{vntrho}, we have $v_j \ge v_n \geq -B (n/k)^{1/6} n^{2/3}$ for each $j \in \llbracket 1, n \rrbracket$. Hence, setting $\bm{v}' = (v', v', \ldots , v')$, where $v' = -B(n/k)^{1/6} n^{2/3}$ appears with multiplicity $n$, and sampling $n$ non-intersecting Brownian bridges $\bm{\mathsf{z}} = (\mathsf{z}_1, \mathsf{z}_2, \ldots , \mathsf{z}_n) \in \llbracket 1, n \rrbracket \times \mathcal{C} \big( [0, \mathsf{T}] \big)$ from the measure $\mathsf{Q}^{\bm{u}; \bm{v}'}$, \Cref{monotoneheight} and \Cref{linear} (with the fact that $\mathsf{T}^2 \le 9A^2 n^{2/3}$) together yield a coupling between $\bm{\mathsf{x}}$, $\bm{\mathsf{y}}$, and $\bm{\mathsf{z}}$ such that 
		\begin{align}\label{e:coupbdd}
			\sfx_j (\sft)\geq \mathsf{z}_j (\mathsf{t}) = \mathsf{y}_j (\mathsf{t}) - \displaystyle\frac{\mathsf{t}}{\mathsf{T}} \cdot v' \ge \sfy_j (\sft)- 9A^2 B  \Big( \displaystyle\frac{n}{k} \Big)^{1/6}\sft \sfT,\qquad \text{for each $j \in \llbracket 1, n \rrbracket$}.   
		\end{align}
		
		Next, define the process $\widetilde{\bm{\mathsf{y}}} (s) = \big( \widetilde{\mathsf{y}}_1 (s), \widetilde{\mathsf{y}}_2 (s), \ldots , \widetilde{\mathsf{y}}_n (s) \big) \in \llbracket 1, n \rrbracket \times \mathcal{C} (\mathbb{R}_{\ge 0})$ for $s \ge 0$ by setting 
		\begin{align}\label{e:tyDBM}
			\widetilde \sfy_j (s)=\frac{s + \sfT}{\sfT} \cdot \sfy_j \left(\frac{s \mathsf{T}}{s + \mathsf{T}}\right), \qquad \text{for each $(j, s) \in \llbracket 1, n \rrbracket \times \mathbb{R}_{\ge 0}$}.
		\end{align} 
		
		\noindent By \Cref{tobridge}, $\widetilde{\bm{\mathsf{y}}} (s)$ has the same law as Dyson Brownian motion, run for time $s$, with initial data $\widetilde{\bm{\mathsf{y}}} (0) = \bm{u}$. Also, combining \eqref{e:coupbdd} and \eqref{e:tyDBM}, we find for each $\mathsf{t} \in   [0,  k^{1/3}]$ that 
		\begin{flalign}
			\label{xjtyjt1}
			\sfx_n (\sft) - \mathsf{x}_n (0) \geq \sfy_n (\sft)- 9A^2 B \Big( \displaystyle\frac{n}{k} \Big)^{1/6}\sft\sfT 
			=\frac{\sfT-\sft}{\mathsf{T}} \cdot \bigg(\widetilde \sfy_n \Big(\frac{\sft \sfT}{\sfT-\sft} \Big)-\widetilde\sfy_n (0) \bigg) - 9A^2 B \Big( \displaystyle\frac{n}{k} \Big)^{1/6} \mathsf{t} \mathsf{T},
		\end{flalign} 
		
		\noindent where we have used the fact that $\mathsf{x}_n (0) = \mathsf{y}_n (0) = \widetilde{\mathsf{y}}_n (0) = u_n = 0$. 
		
		Let us verify that $\widetilde{\bm{\mathsf{y}}} (0) = \bm{u}$ satisfies an instance of \eqref{e:xi-xj2}, more specifically, that 
		\begin{flalign}\label{e:xi-xj}
			u_i - u_j \geq \bigg( \displaystyle\frac{j-i-2\mathfrak{m}}{B (nk)^{1/2}} - 2M \bigg) k^{2/3}, \qquad \text{for any $1 \le i \le j \le n$}.
		\end{flalign}
		
		\noindent Indeed, if $0 \le j-i \le 2\mathfrak{m}$, then $u_i - u_j \ge 0$, which implies \eqref{e:xi-xj}. If instead $j - i \ge 2\mathfrak{m}$, then 
		\begin{flalign*} 
			u_i - u_j \ge (\gamma_{i+\mathfrak{m}}- \gamma_{j-\mathfrak{m}}-2M)k^{2/3} \geq \left(\frac{j-i-2\mathfrak{m}}{B (nk)^{1/2}}-2M \right)k^{2/3}. 
		\end{flalign*}
		
		\noindent where the first bound follows from \eqref{e:gammabb}, and the second follows from the fact that $\sup_{x \in \mathbb{R}} \varrho(x) \le B (n/k)^{1/2}$ (with \Cref{gammarho} for the classical locations); this again verifies \eqref{e:xi-xj}. Hence, \Cref{p:extreme} (with the $M$ there given by $2M+2\mathfrak{m}/(Bn^{1/2} k^{1/2})$ here) applies to $\widetilde{\bm{\mathsf{y}}}$ and yields constants $c = c(B) > 1$, $C_3 = C_3 (B) > 1$, and $C_4 = C_2 (B, D) > 1$ such that the following holds. For each $s \in [0, 1]$, we have with probability $1 - C_4 e^{- c (\log n)^2}$ that 
		\begin{align*}
			\widetilde{\mathsf{y}}_n (sk^{1/3}) & \ge \widetilde{\mathsf{y}}_n(0)  -C_3 k^{2/3} \Bigg(s (nk^{-1})^{1/2} \big|\log (2s^{-1}) \big|^2 + s^{1/2} \left(2M (nk^{-1})^{1/2} +\frac{2\mathfrak{m}}{Bk}\right)^{1/2} \\ 
			& \qquad \qquad \qquad \qquad \qquad \qquad \qquad \qquad \qquad \qquad \qquad \qquad \qquad + (sk^{-1})^{1/2} (\log n)^3 \Bigg) \\
			&\geq \widetilde{\mathsf{y}}_n (0) -4C_3 k^{2/3} \bigg(s (nk^{-1})^{1/2} \big| \log (2s^{-1}) \big|^2 + s^{1/2} M^{1/2} \Big(  \displaystyle\frac{n}{k} \Big)^{1/4}  \bigg),
		\end{align*}
		
		\noindent where in the last bound we also used the facts that $\mathfrak{m} (Bk)^{-1} \le n^{1/2} k^{-1} \log n \le M (nk^{-1})^{1/2}$ (the first by the last bound in \eqref{vntrho} and the second by the facts that $M \ge 1$ and $k$ is sufficiently large). This at $s = t\mathsf{T} / (\mathsf{T} - tk^{1/3}) \in [t, Bt] \subseteq [0, 1]$ (where in the second statement we used the fact that $\mathsf{T} - tk^{1/3} \ge B^{-1} \mathsf{T} $, as $tk^{1/3} \le (1-B^{-1}) Ak^{1/3}  \le (1 - B^{-1}) \mathsf{T}$ by \Cref{uvrho}, and in the third we used the fact that $t \in \big[ 0, (1-B^{-1}) A \big] \subset [0, B^{-1}]$), together with \eqref{xjtyjt1} at $\mathsf{t} = tk^{1/3}$ and the fact that $27A^2 B (nk^{-1})^{1/6} \mathsf{t} \mathsf{T} \le 12A^3 B t (nk^{-1})^{1/2} k^{2/3}$ (as $\mathsf{T} \le 3A n^{1/3}$ by \eqref{vntrho}), finishes the proof of \eqref{e:xiholder}. 
	\end{proof}

	\section{Boundary Removal Coupling}
	
	\label{RectangleCouple}
	
	In this section we state and establish the existence of the boundary removal coupling. We first state this coupling in \Cref{s:couplingandholder}; it relies on a certain event, called a boundary tall rectangle event $\textbf{BTR}$ (see \Cref{ftrbtr} below). In \Cref{EventC}, we introduce and discuss properties of a stronger variant of this $\textbf{BTR}$ event that will be useful for us. We then state several preliminary couplings in \Cref{s:couple} (which will be proved in \Cref{Couple0Proof} below). We will use these, assuming a certain improved H\"{o}lder estimate (to be shown in \Cref{RegularImproved} below), to prove \Cref{c:finalcouple} in \Cref{ProofCouple1}. Throughout this section, we let $\bm{\mathsf{x}} = (\mathsf{x}_1, \mathsf{x}_2, \ldots )$ denote a $\mathbb{Z}_{\ge 1} \times \mathbb{R}$ indexed line ensemble satisfying the Brownian Gibbs property; we also recall the $\sigma$-algebra $\mathcal{F}_{\ext}$ from \Cref{property}. In this section and the remaining ones in this chapter, we will set $\chi = 2^{-5000}$.

	\subsection{Coupling}
	
	\label{s:couplingandholder}
	
	In this section we state a result indicating the existence of a coupling between a family of non-intersecting Brownian bridges with lower boundary, and one with the same starting and ending data but without a lower boundary. We will assume that these families are subject to certain conditions, to which end we must first introduce several events. We begin with the following location events, which are similar to the medium position ones of \Cref{eventsregular1}. 
	
	\begin{definition}
		
		\label{eventlocation}
		
		For any integer $k \ge 1$; real numbers $b \le B$ and $t \in \mathbb{R}$; and subset $\mathcal{T} \subseteq \mathbb{R}$, define the \emph{location events} $\textbf{LOC}_k (t; b; B) = \textbf{LOC}_k^{\bm{\mathsf{x}}} (t; b; B)$\index{L@$\textbf{LOC}$; location event} and $\textbf{LOC}_k (\mathcal{T}; b; B) = \textbf{LOC}_k^{\bm{\mathsf{x}}} (\mathcal{T}; b; B)$ by 
		\begin{flalign*}
			\textbf{LOC}_k (t; b; B) = \big\{ b \le \mathsf{x}_k (t) \le B \big\}; \qquad \textbf{LOC}_k  (\mathcal{T}; b; B) = \bigcap_{s \in \mathcal{T}} \textbf{LOC}_k (s; b; B). 
		\end{flalign*}  
	\end{definition}

	We next define an event, which constrains the paths $\big( \mathsf{x}_j (t) \big)$ for $(j, t) \in \llbracket 1, n+1 \rrbracket \times [-Ak^{1/3}, Ak^{1/3}]$; we imagine $A$ as bounded and $L^{3/2} = k^{-1} n$ as large, making the rectangle $\llbracket 1, n+1 \rrbracket \times [-Ak^{1/3}, Ak^{1/3}]$ ``tall.'' The following event imposes that the $\mathsf{x}_j (t)$ are ``on-scale'' (so that $-\mathsf{x}_j (t)$ is of order $j^{2/3}$) for $(j, t)$ on the boundary of the rectangle $\llbracket 1, n+1  \rrbracket \times [-Ak^{1/3}, Ak^{1/3}]$. It also imposes some weak bounds on $\mathsf{x}_j (-2Aj^{1/3})$ and $\mathsf{x}_j (2Aj^{1/3})$ for $j \in \llbracket k, n+1 \rrbracket$; this in particular constrains $\bm{\mathsf{x}} (t)$ for $t$ outside of the strip $[-2Ak^{1/3}, 2Ak^{1/3}]$, which will only be needed in the proof of \Cref{l:holder} below (see \Cref{proofd00}).

	\begin{figure}
	\center
\includegraphics[scale = .6, trim=0 1cm 0 1cm, clip]{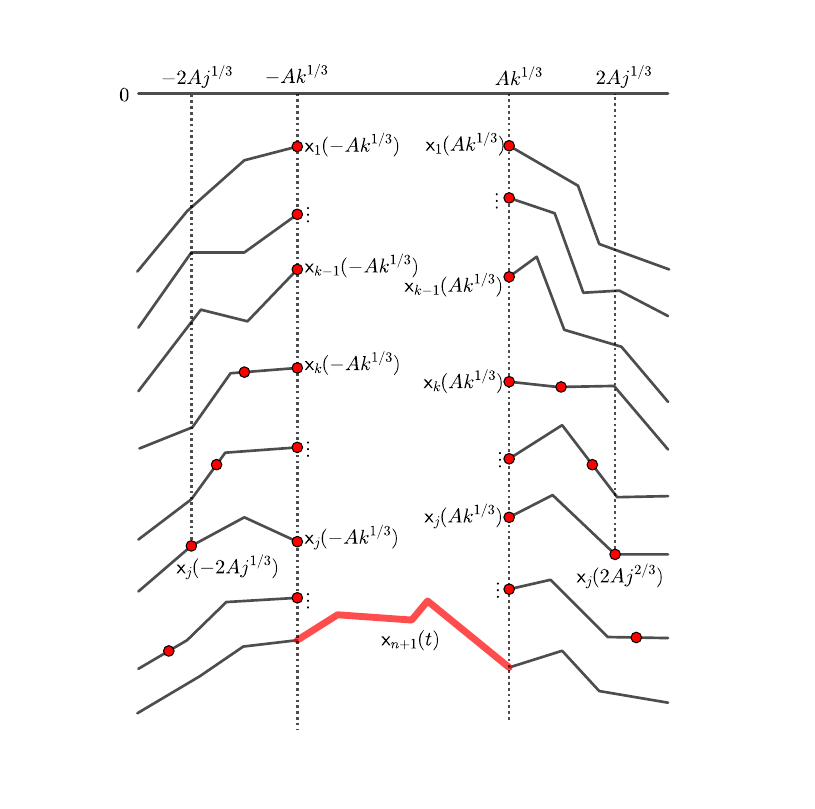}

\caption{Shown above, the red points and curves constitute what are being constrained in the $\textbf{BTR}$ event from  \Cref{ftrbtr}.}
\label{f:BTR}
	\end{figure}

	\begin{definition} 
		
		\label{ftrbtr} 
		
		Recall $\chi = 2^{-5000}$. Fix positive integers $n \ge k \ge 1$, and real numbers $A > 0$ and $B, L \ge 1$, such that $n = L^{3/2} k$. Recalling \Cref{eventlocation}, define the \emph{boundary tall rectangle event} $\textbf{BTR}_n (A; B) = \textbf{BTR}_n^{\bm{\mathsf{x}}} (A, B; k, L)$, measurable with respect to $\mathcal{F}_{\ext}^{\bm{\mathsf{x}}} \big(\llbracket 1, n \rrbracket \times (-Ak^{1/3}, Ak^{1/3}) \big)$, by\index{B@$\textbf{BTR}$; boundary tall rectangle event}
		\begin{flalign} \begin{split}\label{e:cEn}
				\textbf{BTR}_n  (A, B) & =\bigcap_{j = 1}^{n}\textbf{{LOC}}_j \big( \{-Ak^{1/3}, Ak^{1/3}\}; -B j^{2/3}- Bk^{2/3}; B k^{2/3} - B^{-1} j^{2/3} \big)\\
				& \qquad \cap \textbf{{LOC}}_{n+1} \big( [-Ak^{1/3}, Ak^{1/3}]; -B (n+1)^{2/3}-Bk^{2/3}; Bk^{2/3} - B^{-1} (n+1)^{2/3}  \big)\\
				& \qquad \cap \bigcap_{j =k}^{n+1} \big\{\sfx_j(-2A j^{1/3})\geq -L^{\chi/2} j^{2/3}\big\} \cap \big\{ \mathsf{x}_j (2Aj^{1/3}) \ge -L^{\chi/2} j^{2/3} \big\}.
			\end{split}
		\end{flalign}
	
		\noindent See \Cref{f:BTR} for a depiction.

	\end{definition} 
	
	The following assumption imposes that the $\textbf{BTR}$ event has positive probability (in fact, we will only require it to be nonempty, just so that we can restrict to it).

	\begin{assumption}
		
		\label{a:nkrelation}
		
		Fix integers $n \ge k \ge 1$ and real numbers $A, B, D, L \ge 1$, such that $n = L^{3/2} k$ and $L \in [1, k^D]$. Recalling \Cref{ftrbtr}, assume that $\bP \big[ \textbf{BTR}_n^{\bm{\mathsf{x}}} (A;B) \big] > 0$.
		
	\end{assumption}

\begin{figure}
	\center
\includegraphics[scale = .5, trim=0 1cm 0 0.5cm, clip]{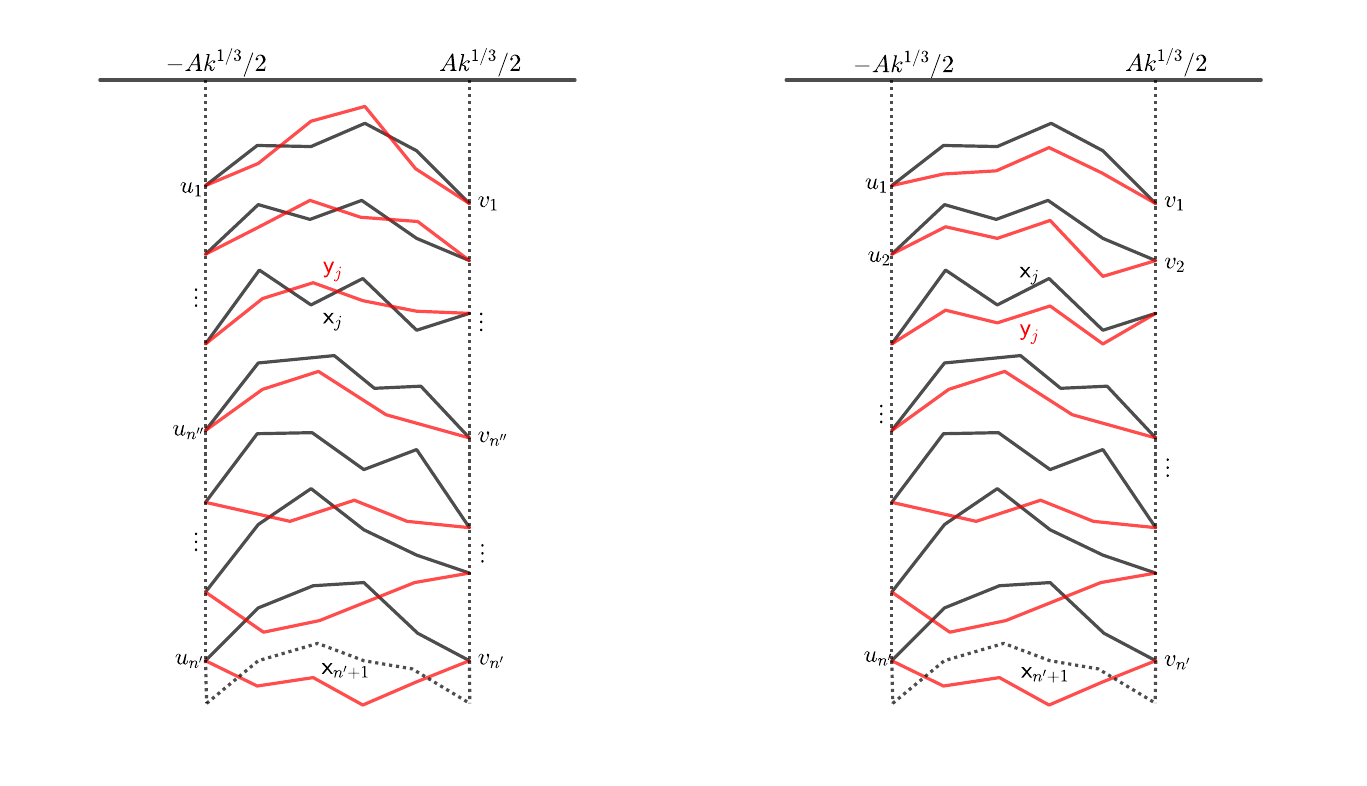}

\caption{Theorem \ref{c:finalcouple} is depicted above. Its first part exhibits a coupling between $\bm{\mathsf{x}}$ and $\bm{\mathsf{y}}$ such that, with high probability, $\mathsf y_j \geq \mathsf x_j -L^{-1/2^{4000}k^{2/3}}$ for each $j \in \llbracket 1, n'' \rrbracket$; this is shown on the left. Its second part exhibits one such that $\mathsf y_j\leq \mathsf x_j$ for $j \in \llbracket 1, n' \rrbracket$; this is shown on the right.}
\label{f:remove_boundary}
	\end{figure}

	The following theorem, to be established in \Cref{ProofCouple1} below, indicates that we may couple $\bm{\mathsf{x}}$ satisfying \Cref{a:nkrelation} with the line ensemble $\bm{\mathsf{y}}$ obtained by restricting it to the time interval $[ -Ak^{1/3} / 2,  Ak^{1/3} / 2]$ and removing its lower curves, so that with high probability the top paths in $\bm{\mathsf{x}}$ and $\bm{\mathsf{y}}$ are ``close to each other'' if $L$ is large (more precisely, we provide two couplings between $\bm{\mathsf{x}}$ and $\bm{\mathsf{y}}$, so that $\bm{\mathsf{x}}$ is almost below $\bm{\mathsf{y}}$ in the former and $\bm{\mathsf{x}}$ is above $\bm{\mathsf{y}}$ in the latter). See \Cref{f:remove_boundary} for a depiction.

	\begin{thr}[Boundary removal coupling]
		
		\label{c:finalcouple}

		Adopt \Cref{a:nkrelation}, and suppose that $A \ge 2$. There exist constants $c = c(A, B) > 0$ and $C = C(A, B, D) > 1$ such that the following holds if $L \ge C$. Set
		\begin{flalign*} 
			n' = \big\lceil L^{1/2^{3000}} k \big\rceil; \qquad n'' = \big\lceil L^{1/2^{4000}} k \big\rceil.
		\end{flalign*} 
		
		\noindent There is an event $\mathscr{A} \subseteq \textbf{\emph{BTR}}_n^{\bm{\mathsf{x}}} (A; B)$, measurable with respect to $\mathcal{F}_{\ext}^{\bm{\mathsf{x}}} \big( \llbracket 1, n' \rrbracket \times (-Ak^{1/3} / 2, Ak^{1/3} / 2) \big)$, satisfying $\mathbb{P} \big[ \textbf{\emph{BTR}}_n^{\bm{\mathsf{x}}} (A; B) \setminus \mathscr{A} \big] \le C e^{-c(\log k)^2}$ and the following. Condition on $\mathcal{F}_{\ext}^{\bm{\mathsf{x}}} \big( \llbracket 1, n' \rrbracket \times ( -Ak^{1/3} / 2, Ak^{1/3} / 2) \big)$; restrict to $\mathscr{A}$; and define the $n'$-tuples $\bm{u} = \bm{\mathsf{x}}_{\llbracket 1, n' \rrbracket} ( -Ak^{1/3} / 2) \in \mathbb{W}_{n'}$ and $\bm{v} = \bm{\mathsf{x}}_{\llbracket 1, n' \rrbracket} ( Ak^{1/3} / 2) \in \mathbb{W}_{n'}$. Sample $n'$ non-intersecting Brownian bridges $\bm{\mathsf{y}} = (y_1, y_2, \ldots ,y_{n'}) \in \llbracket 1, n' \rrbracket \times [ -Ak^{1/3} / 2, Ak^{1/3} / 2]$ from the measure $\mathsf{Q}^{\bm{u}; \bm{v}}$.

		\begin{enumerate}
			\item There exists a coupling between $\bm{\mathsf{x}}$ and $\bm{\mathsf{y}}$ such that 
			\begin{flalign*} 
				\mathbb{P} \Bigg[ \bigcap_{j =1}^{n''} \bigcap_{|t| \le Ak^{1/3}/2} \big\{ \mathsf{y}_j (t) \ge \mathsf{x}_j (t) - L^{-1/2^{4000}} k^{2/3} \big\} \Bigg] \ge 1 - C e^{-c(\log k)^2}.  
			\end{flalign*} 
			
			\item There exists a coupling between $\bm{\mathsf{x}}$ and $\bm{\mathsf{y}}$ such that $\mathsf{y}_j (t) \le \mathsf{x}_j (t)$ for each $(j, t) \in \llbracket 1, n' \rrbracket \times [ -Ak^{1/3} / 2, Ak^{1/3} / 2 ]$.

		\end{enumerate}

	\end{thr}

	\subsection{Completed Tall Rectangle Events}
	
	\label{EventC} 
	
	In this section we introduce and show properties of a variant of the boundary tall rectangle event $\textbf{BTR}$ from \Cref{ftrbtr}. In addition to imposing that $\textbf{BTR}$ holds, it further imposes that the $\mathsf{x}_j (t)$ satisfy the location events (recall \Cref{eventlocation}) on the complete rectangle $\llbracket 1, n+1 \rrbracket \times [-Ak^{1/3}, Ak^{1/3}]$, as opposed to only on its boundary. 
	
	\begin{definition} 
		
		\label{ftrbtr2} 
		
		Adopting the notation of \Cref{ftrbtr} (and recalling \Cref{eventlocation}), define the \emph{complete tall rectangle event} $\textbf{CTR}_n (A; B) = \textbf{CTR}_n^{\bm{\mathsf{x}}} (A, B; k, L)$\index{C@$\textbf{CTR}$; complete tall rectangle event} by
		\begin{flalign} \label{eventcompleter0}
			\textbf{CTR}_n(A; B) & = \textbf{BTR}_n (A; B) \cap \bigcap_{j = 1}^n \textbf{{LOC}}_j \big( [-Ak^{1/3}, Ak^{1/3}] ; -Bj^{2/3}-Bk^{2/3}; -B^{-1} j^{2/3}+B k^{2/3}  \big).
		\end{flalign}
		
	\end{definition}
	
	To prove \Cref{c:finalcouple}, we will frequently make use of the following lemma, indicating that the boundary tall rectangle event of \Cref{ftrbtr} likely implies the complete one (with different constants).

	\begin{lem}
		
		\label{ctr} 
		
		Adopting \Cref{a:nkrelation} and assuming that $A \ge 1$, there exist constants $c = c(A, B) > 0$ and $C = C(A, B, D) > 1$ such that, setting $\widetilde B=  12 A^2 B^3$, we have
		\begin{align}\label{e:loc2}
			\bP \big[  \textbf{\emph{BTR}}_n^{\bm{\mathsf{x}}} (A;B) \cap \textbf{\emph{CTR}}_n^{\bm{\mathsf{x}}} (A;\widetilde B)^{\complement} \big]\leq C e^{-c(\log k)^2}.
		\end{align}
	\end{lem}

	\begin{proof}
		
		Condition on $\mathcal{F}_{\ext}^{\bm{\mathsf{x}}} \big( \llbracket 1, n \rrbracket \times (-Ak^{1/3}, Ak^{1/3}) \big)$ and restrict to the event $\textbf{BTR}_n (A; B)$. Define the $n$-tuples $\bm{u} = \bm{\mathsf{x}}_{\llbracket 1, n \rrbracket} (-Ak^{1/3}) \in \mathbb{W}_n$ and $\bm{v} = \bm{\mathsf{x}}_{\llbracket 1, n \rrbracket} (Ak^{1/3}) \in \mathbb{W}_n$, and the function $f : [-Ak^{1/3}, Ak^{1/3}] \rightarrow \mathbb{R}$ by setting $f(s) = \mathsf{x}_{n+1} (s)$ for each $s \in [-Ak^{1/3}, Ak^{1/3}]$. Then, the law of $\bm{\mathsf{x}}$ is given by $\mathsf{Q}_f^{\bm{u}; \bm{v}}$. 
		
		By \Cref{ftrbtr} (and \Cref{eventlocation} for the $\textbf{LOC}$ events), we have $\max \{ u_j, v_j \} \le Bk^{2/3} - B^{-1} j^{2/3}$ and $f(s) \le Bk^{2/3} - B^{-1} (n+1)^{2/3}$ for each $(j, s) \in \llbracket 1, n \rrbracket \times [-Ak^{1/3}, Ak^{1/3}]$. Hence, the first part of \Cref{p:compareAiry} (applied with the $(b-a, d, M)$ there equal to $(2Ak^{1/3},  B^{-1}, Bk^{2/3}$ here), yields $c_1 = c_1 (A, B) > 0$ and $C_1 = C_1 (A, B, D) > 1$ such that 
		\begin{flalign}
			\label{xjta1} 
			\mathbb{P} \Bigg[ \bigcap_{j=1}^n \bigcap_{|t| \le Ak^{1/3}} \big\{ \mathsf{x}_j (t) \le \widetilde{B} k^{2/3} - \widetilde{B}^{-1} j^{2/3} \big\} \Bigg] \ge 1 - C_1 e^{-c_1 (\log k)^2}.
		\end{flalign}
		
		\noindent where we have used the fact that for $(b-a, d, M) = (2Ak^{1/3}, B^{-1}, Bk^{1/3})$ we have 
		\begin{flalign*} 
			\displaystyle\frac{9 \pi^2 (b-a)^2}{64 d^3} + M + (\log n)^2   \le (6A^2 B^3 + B + 2) k^{2/3} \le \widetilde{B} k^{2/3}, 
		\end{flalign*} 
		
		\noindent for sufficiently large $n$ (as $k^{3D/2+1} \le L^{3/2} k = n$), and that $\widetilde{B}^{-1} < B^{-1}$ (using $A\geq 1$). 
		
		Similarly, \Cref{ftrbtr} (and \Cref{eventlocation}) yields $\min \{ u_j, v_j \} \le -Bj^{2/3} - Bk^{2/3}$ for each $j \in \llbracket 1, n \rrbracket$. Hence, the second part of \Cref{p:compareAiry} (applied with the $(A, B, M)$ there equal to $(2A, B, Bk^{2/3})$ here) yields constants $c_2 = c_2 = (A, B) > 0$ and $C_2 = C_2 (A, B, D) > 1$ such that 
		\begin{flalign*}
			\mathbb{P} \Bigg[ \bigcap_{j=1}^n \bigcap_{|t| \le Ak^{1/3}} \big\{ \mathsf{x}_j (t) \ge -\widetilde{B} k^{2/3} - \widetilde{B} j^{2/3} \big\} \Bigg] \le C_2 e^{-c_2 (\log k)^2},
		\end{flalign*} 
		
		\noindent where we used the facts that for $n$ sufficiently large we have $Bk^{2/3} + 2 (\log n)^2 \le (B+1) k^{2/3} \le \widetilde{B} k^{2/3}$ and $2 (2A)^2 + B + 3 \le 12 A^2 B^3 = \widetilde{B}$. Combining this with \eqref{xjta1}, and using the definition \eqref{eventcompleter0} of the $\textbf{CTR}$ event (with \Cref{eventlocation}), yields the lemma.
	\end{proof}

	\subsection{Preliminary Couplings }
	
	\label{s:couple}
	
	In this section we state several preliminary couplings that will be used to prove \Cref{c:finalcouple}. We first define an improved variant of the H\"{o}lder regularity event from \Cref{eventtsregular2}, which can have a smaller error term than what appears there. The error in the $\textbf{FHR}$ event (which we view as a weaker bound) in \eqref{fhr1} below is analogous to, but slightly different from, the one appearing in \eqref{eventsregular0}; that in the $\textbf{SHR}$ event there (which we view as stronger) is analogous to the one from \eqref{e:xiholder}. The improved H\"{o}lder regularity event will impose the first (weaker) bound on all top $n'$ curves and the second (stronger) bound on some of the top $n'$ ones.

	\begin{definition}
		
		\label{eventregularityimproved} 
		
		Recall $\chi = 2^{-5000}$. Fix integers $n \ge k \ge 1$ and real numbers $A, D, L, S > 0$; $R \ge 2A$; and $\beta \in [0, 1)$, with $n = L^{3/2} k$ and $L \in [1, k^D]$. For any integer $j \in \llbracket 1, n \rrbracket$, define the \emph{first H\"{o}lder event} and \emph{second H\"{o}lder event}, denoted by\index{F@$\textbf{FHR}$; first H\"{o}lder regular event}\index{S@$\textbf{SHR}$; second H\"{o}lder regular event} $\textbf{FHR}_j (A) = \textbf{FHR}_j^{\bm{\mathsf{x}}} (A)$ and $\textbf{SHR}_j (A; \beta; R) = \textbf{SHR}_j^{\bm{\mathsf{x}}} (A; \beta; R)$, respectively, as   
		\begin{flalign}
			\label{fhr1} 
			\begin{aligned}
				& \textbf{FHR}_j (A) =  \bigcap_{\substack{|s| \le Ak^{1/3} \\ |s + tk^{1/3}| \le Ak^{1/3}}} \bigg\{ k^{-2/3} \cdot \big| \mathsf{x}_j (s + tk^{1/3}) - \mathsf{x}_j (s)\big| \\
				& \qquad \qquad \qquad \qquad \qquad \qquad \qquad \qquad \qquad  \le  L^{\chi} \Big(\frac{j \vee k}{k} \Big)^{1/3} |t| +4 \Big(\frac{j \vee k}{k} \Big)^{1/2} |t|^{1/2}  + k^{-D} \bigg\}; \\
				&\textbf{SHR}_j (A; \beta; R) = \bigcap_{\substack{|s| \le Ak^{1/3} \\ |s + tk^{1/3}| \le Ak^{1/3}}} \Bigg\{ k^{-2/3} \cdot \big| \sfx_j(s + tk^{1/3})-\sfx_j(s) \big| \\
				& \qquad \qquad \qquad \qquad \qquad \qquad \qquad  \qquad \quad \leq R \bigg( \Big( \frac{j}{k}\Big)^{1/2} |t| \big( \log  (R|t|^{-1})  \big)^2 + \Big(\frac{j}{k}\Big)^{2\beta/3}  |t|^{1/2} +k^{-D}\bigg) \Bigg\}.
			\end{aligned} 
		\end{flalign}
		
		\noindent For any integer $n' \in \llbracket L^{3S\chi/2} k, n \rrbracket$, define the \emph{improved H\"{o}lder event}\index{I@$\textbf{IHR}$; improved H\"{o}lder event} denoted by $\textbf{IHR}_{n'} (A;\beta; R;S) = \textbf{IHR}_{n'}^{\bm{\mathsf{x}}} (A; \beta; R; S)$ as
		\begin{flalign}
			\label{ihrnn} 
			& \textbf{IHR}_{n'} (A; \beta; R; S)  = \bigcap_{j=1}^{n'} \textbf{FHR}_j (A) \cap \bigcap_{j = \lceil L^{3S \chi/2} k \rceil}^{n'} \textbf{SHR}_j (A; \beta; R) 
		\end{flalign}
	
		\noindent Observe that the events above also depend implicitly on the parameters $L$, $D$, and $k$, but we will omit this from the notation for brevity. 
	\end{definition}

	We next define an event, which is nearly the one on which we will be able to formulate the preliminary couplings.
	
	\begin{definition}
		
		\label{nyy} 
		
		Adopting the notation of \Cref{eventregularityimproved}, further let $B \ge 1$ be a real number, and define the \emph{initial coupling event}\index{I@$\textbf{ICE}$; initial coupling event} $\textbf{ICE}_{n'} = \textbf{ICE}_{n'}^{\bm{\mathsf{x}}} = \textbf{ICE}_{n'}^{\bm{\mathsf{x}}} (A, B; \beta; R, S)$ by
		\begin{flalign*}
			\textbf{ICE}_{n'} & = \textbf{IHR}_{n'} (A; \beta; R; S) \cap  \bigcap_{j=1}^{n'}  \textbf{LOC}_j \big( [-Ak^{1/3}, Ak^{1/3}]; -Bj^{2/3} - Bk^{2/3}; Bk^{2/3} - B^{-1} j^{2/3} \big).
		\end{flalign*} 
	
		\noindent Observe that this event also implicitly depends on the parameters $L$, $D$, and $k$, but we will omit this from the notation for brevity. 
		
	\end{definition}
	
\begin{figure}
	\center
\includegraphics[scale = .6, trim=0 1cm 0 1cm, clip]{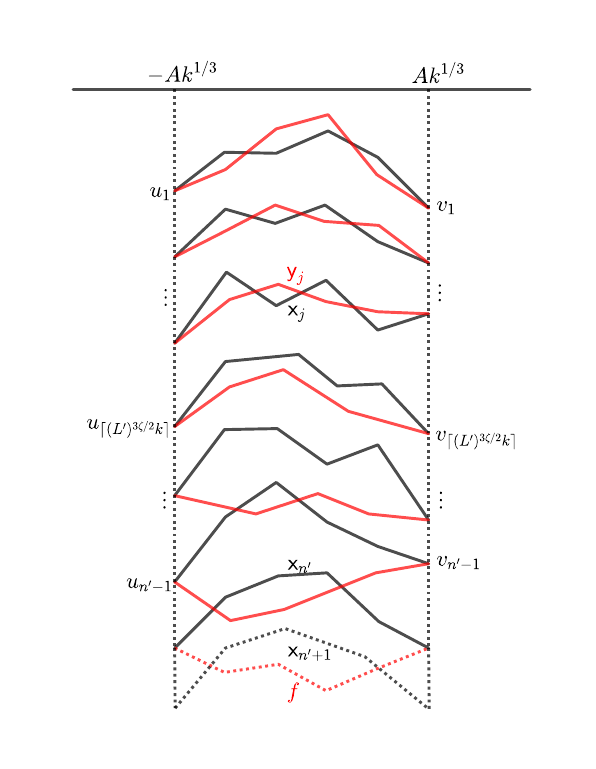}

\caption{Shown above is a depiction of \Cref{p:comparison2}.}
\label{f:preliminary_couple}
	\end{figure}

	The following proposition constitutes the preliminary coupling we will use to establish \Cref{c:finalcouple}; its proof is given in \Cref{ProofCouple0} below (see \Cref{Compare00} for a heuristic). Let us briefly explain this proposition. It considers a family $\bm{\mathsf{y}}$ of non-intersecting Brownian bridges on the interval $[-Ak^{1/3}, Ak^{1/3}]$ with the same starting and ending data as $\bm{\mathsf{x}}$, but with a different lower boundary $f$. Fix an integer $n' = (L')^{3/2} k$ with $L'$ not too small (see \eqref{ldeltal}), and assume that $\textbf{ICE}^{\bm{\mathsf{x}}}$ is likely (see \eqref{probabilityeventc0}) and that $f \ge \mathsf{x}_{n'} - (L^{\alpha} k^{2/3})$ for some $\alpha \le 9/10$ (see \eqref{yxpl}). Then, it provides a coupling between $\bm{\mathsf{x}}$ and $\bm{\mathsf{y}}$ so that the top several paths in the latter nearly bounds those in the former from above (the corresponding lower bound will be a quick consequence of height monotonicity \Cref{monotoneheight}).  See \Cref{f:preliminary_couple} for a depiction.

	\begin{prop}\label{p:comparison2}
		
		Recall $\chi = 2^{-5000}$. For any real numbers 
		\begin{flalign*} 
			A, B, \Xi \ge 1; \qquad \beta \in \Big[ \frac{3}{8}, \frac{7}{8} \Big]; \qquad \xi \in \Big( 0, \displaystyle\frac{1}{2} \Big); \qquad \alpha \in \Big[ 2\beta - \displaystyle\frac{9}{10}, \displaystyle\frac{9}{10} \Big]; \qquad R \ge 2A,
		\end{flalign*} 
		
		\noindent there exist two constants $\zeta = \zeta (\alpha, \beta) \in [2^{-250}, 1]$ and $C = C (A, B, R, \xi, \Xi) > 1$ such that the following holds. Fix real numbers $D, L \ge 1$ and $S \in [1, 2^{500}]$. Further let $n \ge k \ge 1$ be integers, such that $n = L^{3/2} k$ and $L \in [C, k^D]$. Let $n' \in \llbracket k, n \rrbracket$ be an integer, and set $L' = (k^{-1} n')^{2/3}$. Condition on $\mathcal{F}_{\ext}^{\bm{\mathsf{x}}} \big( \llbracket 1, n'-1 \rrbracket \times (-Ak^{1/3}, Ak^{1/3}) \big)$; define the $(n'-1)$-tuples $\bm{u} = \bm{\mathsf{x}}_{\llbracket 1, n' -1 \rrbracket} (-Ak^{1/3}) \in \mathbb{W}_{n'-1}$ and $\bm{v} = \bm{\mathsf{x}}_{\llbracket 1, n' - 1 \rrbracket} (Ak^{1/3}) \in \mathbb{W}_{n'-1}$; and let $f : [-Ak^{1/3}, Ak^{1/3}] \rightarrow \mathbb{R}$ denote a function. Assume that 
		\begin{flalign} 
			\label{ldeltal} 
			L' \ge L^{4S\chi / \zeta},
		\end{flalign}
		
		\noindent that
		\begin{flalign}
			\label{probabilityeventc0} 
			\mathbb{P} \big[ \textbf{\emph{ICE}}_{n'}^{\bm{\mathsf{x}}} (A, B; \beta; R, S) \big] \ge 1 - \Xi e^{-\xi (\log k)^2},
		\end{flalign} 
		
		\noindent and that
		\begin{flalign}
			\label{yxpl}
			f (s) \ge \mathsf{x}_{n'} (s) - (L')^{\alpha} k^{2/3}, \qquad \text{for each $s \in [-Ak^{1/3}, Ak^{1/3}]$}. 
		\end{flalign}
		
		\noindent Sample non-intersecting Brownian bridges $\bm{\mathsf{y}} = (\mathsf{y}_1, \mathsf{y}_2, \ldots , \mathsf{y}_{n'-1}) \in \llbracket 1, n' - 1 \rrbracket \times \mathcal{C} \big( [-Ak^{1/3}, Ak^{1/3}] \big)$ from the measure $\mathsf{Q}_f^{\bm{u}; \bm{v}}$. There exists a coupling between $\bm{\mathsf{x}}$ and $\bm{\mathsf{y}}$ such that 
		\begin{align}\label{e:zilow}
			\mathbb{P} \Bigg[ \bigcap_{j=1}^{\lceil (L')^{3\zeta/2} k \rceil} \bigcap_{|s| \le Ak^{1/3}} \big\{ \mathsf{y}_j (s)\geq \mathsf{x}_j (s)- (L')^{\zeta (2\beta-7/8)}k^{2/3} \big\} \Bigg] \ge 1 - 3^{250} \Xi \cdot  e^{-\xi (\log k)^2/2}.
		\end{align}
		
	\end{prop}

	To obtain a lower bound on the $\bm{\mathsf{y}}$ paths, \Cref{p:comparison2} imposes a lower bound \eqref{yxpl} on the lower boundary $f$. The next corollary, to be established in \Cref{Proofyf} below, removes this constraint, enabling upper and lower bounds on the $\bm{\mathsf{y}}$ paths assuming that $f = -\infty$.

	\begin{cor}
		
		\label{xycouple2}
		
		For any real numbers $\xi \in ( 0, 1 / 2 )$ and $A, B, D, \Xi \ge 1$, and $R \ge 2A$, there exist three constants $\zeta = \zeta (\beta) \in [2^{-250}, 1]$, $c = c(A, B, \xi) > 0$, and $C = C (A, B, D, R, \xi, \Xi) > 1$ such that the following holds. Fix real numbers $\beta \in [ 3 / 8, 7 / 8 ]$,  $L \ge 1$ and $S \in [1, 2^{500}]$; suppose that $L \in [C, k^D]$. Further let $n \ge k \ge 1$ be integers, such that $n = L^{3/2} k$. Let $n' \in \llbracket k, n \rrbracket$ be an integer, and set $L' = (k^{-1} n')^{2/3}$. Condition on $\mathcal{F}_{\ext}^{\bm{\mathsf{x}}} \big( \llbracket 1, n'  \rrbracket \times (-Ak^{1/3}, Ak^{1/3}) \big)$; define the $n'$-tuples $\bm{u} = \bm{\mathsf{x}}_{\llbracket 1, n' \rrbracket} (-Ak^{1/3}) \in \mathbb{W}_{n'}$ and $\bm{v} = \bm{\mathsf{x}}_{\llbracket 1, n' \rrbracket} (Ak^{1/3}) \in \mathbb{W}_{n'}$; and assume that \eqref{ldeltal} and \eqref{probabilityeventc0} both hold. Sample $n'$ non-intersecting Brownian bridges $\bm{\mathsf{y}} = (\mathsf{y}_1, \mathsf{y}_2, \ldots , \mathsf{y}_{n'}) \in \llbracket 1, n' \rrbracket \times \mathcal{C} \big( [-Ak^{1/3}, Ak^{1/3}] \big)$ from the measure $\mathsf{Q}^{\bm{u}; \bm{v}}$. 
		
		\begin{enumerate}
			\item There exists a coupling between $\bm{\mathsf{x}}$ and $\bm{\mathsf{y}}$ such that 
			\begin{flalign}
				\label{yx1}
				\mathbb{P} \Bigg[ \bigcap_{j=1}^{\lceil (L')^{3\zeta/2} k \rceil} \bigcap_{|s| \le Ak^{1/3}} \big\{ \mathsf{y}_j (s) \ge \mathsf{x}_j (s) -  (L')^{\zeta (2\beta-7/8)} k^{2/3} \big\} \Bigg] \ge 1 - C e^{-c (\log k)^2}. 
			\end{flalign}
			
			\item There exists a coupling between $\bm{\mathsf{x}}$ and $\bm{\mathsf{y}}$ such that 
			\begin{flalign} 
				\label{yx2} 
				\mathsf{y}_j (s) \le \mathsf{x}_j (s), \qquad \text{for each $(j, s) \in \llbracket 1, n' \rrbracket \times [-Ak^{1/3}, Ak^{1/3}]$}.
			\end{flalign}
		\end{enumerate}

	\end{cor}

	The coupling from \Cref{xycouple2} that guarantees the lower bound \eqref{yx1} for the $\mathsf{y}_j$ is not necessarily the same as that guaranteeing the upper bound \eqref{yx2}. The next corollary, which will be used in \Cref{RegularImproved}, provides concentration upper and lower bounds for these paths. In the following, we recall the classical locations with respect to a measure from \Cref{gammarho}.

	\begin{cor}
		
		\label{p:closerho}
		
		Adopt the notation and assumptions of \Cref{xycouple2}, and let $b \in (0, 1)$ be a real number. For any $t \in \big[ (b-1) A, (1-b) A \big]$, there exist constants $c_1 = c_1 (b, A, B, \xi) > 0$, $C_1 = C_1 (b, A, B) > 1$, and $C_2 = C_2 (b, A, B, D, R, \xi, \Xi) > 1$, and a (deterministic) measure $\mu_t \in \mathscr{P}_{\fin}$, such that $\mu_t (\mathbb{R}) = (L')^{3/2}$ and the following holds if $L \ge C_2$. Below, we denote the classical locations of $\mu_t$ by $\gamma_j = \gamma_{j; n'}^{\mu_t}$ and set $\mathfrak{m}_j = \big\lceil C_1 \log n \cdot \max \{ j^{1/2}, k^{1/2} \} \big\rceil$ for each $j \in \llbracket 1, n' \rrbracket$. 
		
		\begin{enumerate}
			\item  The measure $\mu_t$ admits a density $\varrho_t: \mathbb{R} \rightarrow \mathbb{R}_{\ge 0}$ with respect to Lebesgue measure that satisfies $\supp \varrho_t \subseteq \big[ -C_1 L', C_1 (L')^{3/4} \big]$ and $\varrho_t(x)\leq C_1 \max \{ 1, - x \}^{3/4}$, for any $x \in \mathbb{R}$.
			\item For each $j \in \big\llbracket (\log n')^6 + 1, n \big\rrbracket$, we have $C_1 (j/k)^{2/3} - C_1 \le \gamma_j \le C_1 - C_1^{-1} (j/k)^{2/3}$.
			\item We have 
			\begin{align}\label{e:sfyi}
				\begin{aligned} 
					\mathbb{P} \Bigg[ \bigcap_{j = 1}^{\lceil (L')^{3\zeta/2} k \rceil} \big\{ \gamma_{j + \mathfrak{m}_j} - (L')^{\zeta(2\beta-7/8)} & \leq  k^{-2/3} \cdot \mathsf{x}_j (tk^{1/3}) \\
					&  \leq \gamma_{j-\mathfrak{m}_j} + (L')^{\zeta (2\beta-7/8)} \big\} \Bigg] \ge 1 - C_2 e^{-c_1 (\log k)^2}. 
				\end{aligned} 
			\end{align}
		\end{enumerate}
		
	\end{cor}
	
	\begin{proof}
		
		Throughout this proof, we abbreviate $\textbf{ICE}_{n'} = \textbf{ICE}_{n'}^{\bm{\mathsf{x}}} (A; B; \beta; R; S)$. Define the $n'$-tuples $\bm{u} = \bm{\mathsf{x}}_{\llbracket 1, n' \rrbracket} (-Ak^{1/3}) \in \mathbb{W}_{n'}$ and $\bm{v} = \bm{\mathsf{x}}_{\llbracket 1, n' \rrbracket} (Ak^{1/3}) \in \mathbb{W}_{n'}$, and sample $n'$ non-intersecting Brownian bridges $\bm{\mathsf{y}} = (\mathsf{y}_1, \mathsf{y}_2, \ldots , \mathsf{y}_{n'}) \in \llbracket 1, n' \rrbracket \times \mathcal{C} \big( [-Ak^{1/3}, Ak^{1/3}] \big)$ from $\mathsf{Q}^{\bm{u}; \bm{v}}$. 
		
		Observe from \Cref{nyy} (and \Cref{eventlocation} for the $\textbf{LOC}$ event) that on $\textbf{ICE}_{n'}$ we have for each $j \in \llbracket 1, n' \rrbracket $ that
		\begin{flalign}
			\label{ujvjn} 
			\begin{aligned} 
				& -Bj^{2/3} - Bk^{2/3} \le u_j \le Bk^{2/3} - B^{-1} j^{2/3}; \qquad  -B j^{2/3} - Bk^{2/3} \le v_j \le Bk^{2/3} - B^{-1} j^{2/3}.
			\end{aligned} 
		\end{flalign}
		
		\noindent Since we have conditioned on $\mathcal{F}_{\ext}^{\bm{\mathsf{x}}} \big( \llbracket 1, n' \rrbracket \times (-Ak^{1/3}, Ak^{1/3}) \big)$, and since \eqref{probabilityeventc0} gives $\mathbb{P} [\textbf{ICE}_{n'}] > 0$ for sufficiently large $n$, it follows (from the $\textbf{LOC}$ event in \Cref{nyy} for $\textbf{ICE}_{n'}$) that \eqref{ujvjn} holds. In particular, $\bm{\mathsf{y}}$ satisfies \Cref{a:boundary} with the $(n; A; B)$ there equal to $(n'; 2A; B + b^{-1})$ here. Hence, by  \Cref{p:concentration}, there exist constants $c_2 = c_2 (b, A, B) > 0$, $C_1 = C_1 (b, A, B) > 1$, and $C_2 = C_2 (b, A, B, D) > 1$, and a measure $\mu_t \in \mathscr{P}_{\fin}$ with $\mu_t (\mathbb{R}) = (L')^{3/2}$, such that the following hold if $L \ge C_2$. First, we have $\supp \mu_t \subseteq \big[-C_1 L', C_1 (L')^{3/4} \big]$, and $\mu_t$ admits a density $\varrho_t \in L^1 (\mathbb{R})$ with respect to Lebesgue measure, satisfying $\varrho_t (x) \le C_1 \max \{ 1, -x \}^{3/4}$ (by \eqref{e:rhotup}). Second, we have $-C_1 (j/k)^{2/3}-C_1 \le \gamma_j \le C_1 - C_1^{-1} (j/k)^{2/3}$ for each $j \in \big\llbracket (\log n')^6 + 1, n' \big\rrbracket$. Third, we have 
		\begin{align}
			\label{y20} 
			\mathbb{P} \Bigg[ \bigcap_{j=1}^n \big\{ \gamma_{j + \mathfrak{m}_j} \le k^{-2/3} \cdot \mathsf{y}_j (tk^{1/3})\leq \gamma_{j - \mathfrak{m}_j} \big\} \Bigg] \ge 1 - C_2 e^{-c_2 (\log k)^2},
		\end{align}
		
		\noindent where we have used the fact that $\log n' \ge \log k$ (as $n' = (L')^{3/2} k \ge k$). The first and second statements confirm the first and second parts of the corollary. The third \eqref{y20}, together with the two couplings of \Cref{xycouple2}, implies the third part.
	\end{proof}

	\subsection{Proof of \Cref{c:finalcouple}} 
	
	\label{ProofCouple1}

	In this section we establish \Cref{c:finalcouple}. This will be done using \Cref{p:comparison2}, to which end we must verify \eqref{probabilityeventc0} there. To this end, we have the following quick lemma, indicating that this estimate holds if the boundary tall rectangle event likely implies the improved H\"{o}lder one (see \eqref{btrr2}). Here, we recall the events $\textbf{BTR}$, $\textbf{IHR}$, and $\textbf{ICE}$ from \Cref{ftrbtr}, \Cref{eventregularityimproved}, and \Cref{nyy}, respectively.	
	
	\begin{lem} 
		
		\label{gevent} 
		
		Adopting \Cref{a:nkrelation}, for any real numbers $\xi > 0$ and $\Xi \ge 1$, there exist constants $c = c(A, B, \xi) > 0$ and $C = C(A, B, D, \Xi)$ such that the following holds. Recall $\chi = 2^{-5000}$. Let $A' \in [0, A]$; $R > 1$; $S \in [1, 2^{500}]$; and $\beta \in [0, 1)$ be real numbers, and let $n', n'', n''' \in \llbracket L^{3S\chi/2} k, n \rrbracket$ be integers with $n''' \le n'' \le n'$. Assume 
		\begin{flalign}
			\label{btrr2} 
			\mathbb{P} \big[ \textbf{\emph{BTR}}_n^{\bm{\mathsf{x}}} (A; B) \cap \textbf{\emph{IHR}}_{n'}^{\bm{\mathsf{x}}} (A'; \beta; R; S)^{\complement} \big] \le \Xi e^{-\xi (\log k)^2}.
		\end{flalign}
		
		\noindent Define the event $\mathscr{G}_0 (c; C)$, measurable with respect to $\mathcal{F}_{\ext}= \mathcal{F}_{\ext}^{\bm{\mathsf{x}}} \big( \llbracket 1, n''' \rrbracket \times (-A'k^{1/3}, A' k^{1/3}) \big)$, by
		\begin{flalign}
			\label{g0event} 
			\mathscr{G}_0 (c; C) = \Big\{ \mathbb{P} \big[ \textbf{\emph{ICE}}_{n''}^{\bm{\mathsf{x}}} (A', 12A^2 B^3; \beta; R, S) \big| \mathcal{F}_{\ext} \big] \ge 1 - C e^{-c(\log k)^2} \Big\},
		\end{flalign}
		
		\noindent where we conditioned on $\mathcal{F}_{\ext}$ in the probability. Then, $\mathbb{P} [ \textbf{\emph{BTR}}_n^{\bm{\mathsf{x}}} (A; B) \cap \mathscr{G}_0 (c, C)^{\complement} \big] \le C e^{-c(\log k)^2}$. 
		
	\end{lem} 
	
	\begin{proof}

		Set $\widetilde{B} = 12 A^3 B^3$; abbreviate $\textbf{ICE}_m (A_0) = \textbf{ICE}_m^{\bm{\mathsf{x}}} (A_0, \widetilde{B}; \beta; R, S)$ for any integer $m \in \llbracket L^{3S\chi/2} k, n \rrbracket$ and real number $A_0 \in [0, A]$; and abbreviate $\textbf{ICE}_m = \textbf{ICE}_m (A')$. We will also assume that $\xi \le 1$ (by replacing $\xi$ with $\min \{ \xi, 1 \}$ if necessary). First observe that there exist constants $c_1 = c_1 (A, B) \in (0, 1)$ and $C_1 = C_1 (A, B, D) > 1$ such that 
		\begin{flalign*}
			\mathbb{P} \Bigg[ \textbf{BTR}_n^{\bm{\mathsf{x}}} (A; B) \cap \bigcup_{j=1}^n \textbf{LOC}_j \big( [-Ak^{1/3}, & Ak^{1/3} ]; \widetilde{B}^{-1} j^{2/3} - \widetilde{B}k^{2/3}; \widetilde{B} j^{2/3} + \widetilde{B}k^{2/3} \big)^{\complement} \Bigg] \\
			&  \le \mathbb{P} \big[ \textbf{BTR}_n^{\bm{\mathsf{x}}} (A; B) \cap \textbf{CTR}_n^{\bm{\mathsf{x}}} (A; \widetilde{B})^{\complement} \big] \le C_1 e^{-c_1 (\log n)^2},
		\end{flalign*} 
		
		\noindent where the first inequality follows from \eqref{eventcompleter0} and the second from \Cref{ctr}. Hence, by \eqref{btrr2}, \Cref{nyy} and a union bound, we have 
		\begin{flalign*} 
			\mathbb{P} \big[ \textbf{BTR}_n (A; B) \cap \textbf{ICE}_{n'} (A)^{\complement} \big] \le C_1 e^{-c_1 (\log n)^2} + \Xi e^{-\xi (\log n)^2} \le (C_1 + \Xi) e^{-c_1 \xi (\log n)^2}.
		\end{flalign*} 
	
		\noindent  Since $\textbf{ICE}_{n'} (A) \subseteq \textbf{ICE}_{n''} (A') = \textbf{ICE}_{n''}$ (as $n'' \le n'$ and $A' \le A$), we deduce $\mathbb{P} \big[ \textbf{BTR}_n^{\bm{\mathsf{x}}} (A; B) \cap \textbf{ICE}_{n''}^{\complement} \big] \le (C_1 + \Xi) e^{-c_1 \xi (\log k)^2}$. Therefore,
		\begin{flalign*}
			\mathbb{E} \Big[ \textbf{1}_{\textbf{BTR}_n^{\bm{\mathsf{x}}} (A; B)} \cdot \mathbb{P} [ \textbf{ICE}_m^{\complement} | \mathcal{F}_{\ext}] \big] = \mathbb{P} \big[ \textbf{BTR}_n^{\bm{\mathsf{x}}} (A; B) \cap \textbf{ICE}_m^{\complement} \big] \le (C_1 + \Xi) e^{-c_1 \xi (\log k)^2},
			\end{flalign*}
			
		\noindent since $\textbf{BTR}_n^{\bm{\mathsf{x}}} (A; B)$ is measurable with respect to $\mathcal{F}_{\ext}$. Thus, for $C = (C_1 + \Xi)^{1/2}$ and $c = c_1 \xi / 2$, it follows from a Markov estimate that 
				\begin{flalign}
					\label{gevent0}
					\mathbb{P} \big[ \textbf{1}_{\textbf{BTR}_n^{\bm{\mathsf{x}}} (A; B)} \cdot \mathbb{P} [ \textbf{ICE}_m^{\complement} | \mathcal{F}_{\ext}] \ge C e^{-c (\log k)^2} \big] \le C e^{-c (\log k)^2}. 
					\end{flalign} 
					
					\noindent By \eqref{g0event} (and again the fact that $\textbf{BTR}_n^{\bm{\mathsf{x}}} (A; B)$ is measurable with respect to $\mathcal{F}_{\ext}$), the left side of \eqref{gevent0} is $\mathbb{P} \big[ \textbf{BTR}_n (A; B) \cap \mathscr{G}_0 (c; C)^{\complement} \big]$, and so the lemma follows by \eqref{gevent0}.
	\end{proof}
	
	Although \Cref{gevent} can be used to verify \eqref{probabilityeventc0} in \Cref{p:comparison2}, we must then confirm that \eqref{btrr2} holds. This will be done through the following proposition, to be established in \Cref{RegularImproved} below. It indicates the boundary tall rectangle event of \Cref{ftrbtr} likely implies the improved H\"{o}lder regularity one of \Cref{eventregularityimproved} (on a smaller time interval); we establish it in \Cref{RegularImprovedProof} below.
	
	\begin{prop}
		
		\label{p:Lip}
		
		Adopting \Cref{a:nkrelation} and recalling \Cref{eventregularityimproved}, there exist constants $c = c(A, B) > 0$, $C_1 = C_1 (A, B) > 1$, and $C_2 = C_2 (A, B, D) > 1$ such that the following holds for $L > C_2$. Letting $n'= \lceil L^{1 / 2^{2000}} k \rceil$, we have 
		\begin{align}\label{e:hdd}
			\bP \bigg[\textbf{\emph{BTR}}_n(A;B) \cap  \textbf{\emph{IHR}}_{n'} \Big(\displaystyle\frac{A}{2}; \displaystyle\frac{3}{8}; C_1; 2^{14} \Big)^{\complement} \bigg] \leq C_2 e^{-c(\log k)^2}.
		\end{align}
	\end{prop} 
	
	Given this result, we can now quickly establish \Cref{c:finalcouple}.

	\begin{proof}[Proof of \Cref{c:finalcouple}]

		Set $\widehat{n} = \lceil L^{1/2^{2000}} k \rceil \ge n'$ and $\chi = 2^{-5000}$. By \Cref{p:Lip}, there exist constants $c_1 = c_1 (A, B) > 0$, $C_1 = C_1 (A, B) > 2A$, and $C_2 = C_2 (A, B, D) > 1$ such that 
		\begin{flalign*}
			\mathbb{P} \bigg[ \textbf{BTR}_n (A; B) \cap \textbf{IHR}_{\hat{n}} \Big( \displaystyle\frac{A}{2}; \displaystyle\frac{3}{8}; C_1; 2^{14} \Big)^{\complement} \bigg] \le C_2 e^{-c_1 (\log k)^2}. 
		\end{flalign*}
		
		\noindent This verifies \eqref{btrr2}, with the integers $(n', n'', n''')$ there equal to $(\widehat{n}, n', n')$ here and the real numbers $(\beta; A', R, S; \xi_0; \Xi_0)$ there equal to $( 3 / 8; A / 2, C_1, 2^{14}; c_1, C_2)$ here (observing that $n' \ge L^{2^{20} \chi} k$, for $n' \ge L^{1/2^{3000}} k$). Hence, by \Cref{gevent}, there exist constants $c_2 = c_2 (A, B) > 0$ and $C_3 (A, B, D) > 1$, and an event $\mathscr{A} \subseteq \textbf{BTR}_n (A; B)$ (obtained by intersecting the event $\mathscr{G}_0$ in \eqref{g0event} with $\textbf{BTR}_n (A; B)$) measurable with respect to $\mathcal{F}_{\ext} = \mathcal{F}_{\ext} \big( \llbracket 1, n' \rrbracket \times (-Ak^{1/3} / 2, Ak^{1/3} / 2) \big)$ such that the following holds. First, we have 
		\begin{flalign*} 
			\mathbb{P} \big[ \textbf{BTR}_n (A; B) \setminus \mathscr{A} \big] \le C_3 e^{-c_2 (\log k)^2}.	
		\end{flalign*} 
		
		\noindent Second, conditioning on $\mathcal{F}_{\ext}$ and restricting to the event $\mathscr{A}$, we have 
		\begin{flalign}
			\label{na2x}
			\mathbb{P} \bigg[ \textbf{ICE}_{n'}^{\bm{\mathsf{x}}} \Big( \displaystyle\frac{A}{2}, 12A^2 B^3; \displaystyle\frac{3}{8}; C_1, 2^{14} \Big) \bigg] \ge 1 - C_3 e^{-c_1 (\log k)^2}.
		\end{flalign}
		
		Now we apply \Cref{xycouple2}, with the $(\beta; R, S)$ there equal to $(3 / 8; C_1, 2^{14})$ here; recall $\zeta \in (0, 2^{-250})$ from that corollary and denote $L' = (k^{-1} n')^{2/3} \ge L^{1/2^{3100}}$, so that $n' = (L')^{3/2} k $. The assumption \eqref{probabilityeventc0} is verified by \eqref{na2x}, and \eqref{ldeltal} is confirmed by the fact that $L' \ge L^{1/2^{3100}} \ge L^{2^{20} \chi / \zeta}$ (as $\chi = 2^{-5000}$ and $\zeta \le 2^{-250}$). Hence, \Cref{xycouple2} applies (with the $(A, B; \beta; R, S)$ there given by $(A / 2, 12A^2 B^3; 3 / 8; C_1, 2^{14})$ here); its second part gives the second part of \Cref{c:finalcouple}. Its first part yields constants $c_2 = c_2 (A, B) > 0$, $C_4 = C_4 (A, B) > 1$, and $C_5 = C_5 (A, B, D) > 1$, and a coupling between $\bm{\mathsf{x}}$ and $\bm{\mathsf{y}}$ such that 
		\begin{flalign*}
			\mathbb{P} \Bigg[ \bigcap_{j=1}^{\lceil (L')^{3\zeta/2} k \rceil} \bigcap_{|s| \le Ak^{1/3}/2} \big\{ \mathsf{y}_j (s) \ge \mathsf{x}_j (s) - C_4 (L')^{- \zeta / 8} k^{2/3} \big\} \Bigg] \ge 1 - C_5 e^{-c_2 (\log k)^2}. 
		\end{flalign*}
		
		\noindent This, together with the facts that $n'' = \lceil L^{1/2^{4000}} k \rceil \le (L')^{3\zeta/2} k$ and $C_4 (L')^{-\zeta/8} \le L^{-1/2^{4000}}$ for sufficiently large $L$ (both of which hold since $L' \ge L^{1/2^{2000}}$ and $\zeta \ge 2^{-250}$), yields \Cref{c:finalcouple}. 
	\end{proof}

	\section{Existence of Preliminary Couplings} 
	
	\label{Couple0Proof}
	
	In this section we establish the preliminary couplings from \Cref{s:couple}. We begin by showing \Cref{xycouple2} using \Cref{p:comparison2} in \Cref{Proofyf}. In \Cref{ProofCouple0} we establish \Cref{p:comparison2} assuming an additional result, which is proven in \Cref{Proofxy1}. Throughout this section, we let $\bm{\mathsf{x}} = (\mathsf{x}_1, \mathsf{x}_2, \ldots )$ denote a $\mathbb{Z}_{\ge 1} \times \mathbb{R}$ indexed line ensemble satisfying the Brownian Gibbs property; we also recall the $\sigma$-algebra $\mathcal{F}_{\ext}$ from \Cref{property}, the location events from \Cref{eventlocation}, and the boundary tall rectangle event $\textbf{BTR}$ from \Cref{ftrbtr}. We further  set $\chi = 2^{-5000}$.

	\subsection{Proof of \Cref{xycouple2}}
	
	\label{Proofyf}
	
	In this section we establish \Cref{xycouple2}. We use the notation of that corollary throughout, and we also abbreviate $\textbf{ICE}_m = \textbf{ICE}_m^{\bm{\mathsf{x}}} (A; B; \beta; R; S)$ for any integer $m \in \llbracket L^{3S \chi / 2} k, n \rrbracket$, recalling that $\chi = 2^{-5000}$.
	
	Define the (random) function $f : [-Ak^{1/3}, Ak^{1/3}]$ by setting $f(s) = \mathsf{y}_{n'} (s)$ for each $s \in [-Ak^{1/3}, Ak^{1/3}]$, and define $\bm{\mathsf{z}} = (\mathsf{z}_1, \mathsf{z}_2, \ldots , \mathsf{z}_{n'-1}) \in \llbracket 1, n'-1 \rrbracket \times \mathcal{C} \big( [-Ak^{1/3}, Ak^{1/3}] \big)$ by setting $\mathsf{z}_j (s) = \mathsf{y}_j (s)$ for each $(j, s) \in \llbracket 1, n'-1 \rrbracket \times [-Ak^{1/3}, Ak^{1/3}]$. We will apply \Cref{p:comparison2}, with the $\bm{\mathsf{y}}$ there equal to $\bm{\mathsf{z}}$ here, to which end we must verify the assumptions of that proposition. To do this, define the event 
	\begin{flalign*}
		\mathscr{E} = \bigcap_{|s| \le Ak^{1/3}} \Big\{ \big| \mathsf{x}_{n'} (s) - f(s) \big| \le (L')^{8/9} k^{2/3} \Big\}. 
	\end{flalign*} 
	
	\begin{lem} 
		
		\label{xfe} 
		
		There exist constants $c = c(\xi) > 0$ and $C = C(A, R) > 1$ such that, if $L' > C$, then $\mathbb{P} [ \mathscr{E} ] \ge 1 - (C + \Xi) e^{-c(\log k)^2}$.
	\end{lem} 
	
	\begin{proof}

		Define the events
		\begin{flalign}
			\label{e1s}
			\mathscr{E}_1 = \bigcap_{|s| \le Ak^{1/3}} \bigg\{ \Big| f(s) - \displaystyle\frac{Ak^{1/3} - s}{2Ak^{1/3}} \cdot u_{n'} - \displaystyle\frac{Ak^{1/3} + s}{2Ak^{1/3}} \cdot v_{n'} \Big| \le  (L')^{7/8} k^{2/3} \bigg\}; \quad \mathscr{E}_2 = \mathscr{E}_1 \cap \textbf{ICE}_{n'}.
		\end{flalign}
		
		\noindent We first claim that $\mathscr{E}_2 \subseteq \mathscr{E}$. To show this observe that, by \Cref{nyy} (with the definitions \eqref{ihrnn} and \eqref{fhr1} of the $\textbf{IHR}$ and $\textbf{SHR}$ events), the fact that $n' = (L')^{3/2} k$, and the fact that $t \big| \log (R|t|^{-1}) \big|^2 \le R$ for $|t| \le R$, we have on $\textbf{ICE}_{n'}$ that 
		\begin{flalign*}
			\displaystyle\sup_{|s| \le Ak^{1/3}} \big| \mathsf{x}_{n'} (s) - u_{n'} \big| \le R k^{2/3} \big( R (L')^{3/4} + (L')^{\beta} (2A)^{1/2} + k^{-D} \big) \le 2 R^2 k^{2/3} \big( (L')^{3/4} + (L')^{\beta} \big),
		\end{flalign*}
		
		\noindent where we have used the facts that $R \ge 2A$, that $u_j = \mathsf{x}_{n'} (-Ak^{1/3})$, and $|s + Ak^{1/3}| \le 2Ak^{1/3}$. Together with the definition \eqref{e1s} of $\mathscr{E}_1$ and the facts that $v_{n'} = \mathsf{x}_{n'} (Ak^{1/3})$ and $\beta \le 7 / 8$, it follows that on $\mathscr{E}_2$ we have for any $s \in [-Ak^{1/3}, Ak^{1/3}]$ that 
		\begin{flalign*}
			\big| \mathsf{x}_{n'} (s) - f(s) \big| & \le \big| \mathsf{x}_{n'} (s) - u_{n'} \big| + \big| u_{n'} - v_{n'} \big| + \Big| f(s) - \displaystyle\frac{Ak^{1/3} - s}{2Ak^{1/3}} \cdot u_{n'} - \displaystyle\frac{Ak^{1/3} + s}{2Ak^{1/3}} \cdot v_{n'} \Big| \\
			& \le 4 R^2 k^{2/3} \big( (L')^{3/4} + (L')^{\beta} \big) + (L')^{7/8} k^{2/3} \le 9 R^2 (L')^{7/8} k^{2/3} \le (L')^{8/9} k^{2/3},
		\end{flalign*}
		
		\noindent meaning that $\mathscr{E}$ holds. 
		
		Thus, it remains to show that $\mathbb{P} [\mathscr{E}_2] \ge 1 - (C + \Xi) e^{-c(\log k)^2}$. To this end, by \eqref{probabilityeventc0} and a union bound, it suffices to show for sufficiently large $L'$ that 
		\begin{flalign}
			\label{e1cn} 
			\mathbb{P} [\mathscr{E}_1] \ge 1 - C e^{-c(\log k)^2}. 
		\end{flalign} 
		
		\noindent We apply \Cref{estimatexj2}, with the $B$ there equal to $ (L')^{1/20} (n')^{1/2}$. Since
		\begin{flalign*}
			\displaystyle\sup_{s \in [0, 2Ak^{1/3}]} s^{1/2} \log (4Ak^{1/3} s^{-1}) \le 2 A^{1/2} k^{1/6},
		\end{flalign*} 
		
		\noindent that lemma yields a constant $C_1 > 0$ such that
		\begin{flalign}
			\label{k13sa} 
			\mathbb{P} \Bigg[ \displaystyle\sup_{|s| \le Ak^{1/3}} \Big| \mathsf{y}_{n'} (s) - \displaystyle\frac{Ak^{1/3} - s}{2Ak^{1/3}} \cdot u_{n'} - \displaystyle\frac{Ak^{1/3} + s}{2Ak^{1/3}} \cdot v_{n'} \Big|  \ge 2 A^{1/2} k^{1/6} (L')^{1/20} (n')^{1/2}  \Bigg] \le C e^{-n'},
		\end{flalign}
		
		\noindent for sufficiently large $L'$ (so that $(L')^{1/10} > 2 c^{-1} C$ for the constants $c$ and $C$ in \Cref{estimatexj2}). Since $n' = (L')^{3/2} k$, we have for sufficiently large $L'$ that $2 A^{1/2} k^{1/6} (L')^{1/20} (n')^{1/2} = 2 A^{1/2} (L')^{4/5} k^{2/3} \le (L')^{7/8} k^{2/3}$, and so \eqref{k13sa} (with the definition \eqref{e1s} of $\mathscr{E}_1$ and the facts that $f = \mathsf{y}_{n'}$ and $e^{-n'} \le e^{-(\log k)^2}$, since $n' \ge k$) verifies \eqref{e1cn} and thus the lemma.
	\end{proof} 
	
	Now we can establish \Cref{xycouple2}.

	\begin{proof}[Proof of \Cref{xycouple2}]
		
		Since the laws of $\bm{\mathsf{x}}$ and $\bm{\mathsf{y}}$ are given by $\mathsf{Q}_{\mathsf{x}_{n+1}}^{\bm{u}; \bm{v}}$ and $\mathsf{Q}^{\bm{u}; \bm{v}}$, respectively, the second statement of the corollary follows from height monotonicity \Cref{monotoneheight}. So, it remains to establish the first. We recall $\bm {\mathsf z}$ from the beginning of \eqref{Proofyf}.
		
		To this end, observe by \eqref{probabilityeventc0}, \Cref{xfe}, and a union bound that there exist constants $c_1 = c_1 (\xi) > 0$ and $C_1 = C_1 (A, R) > 1$ 
		\begin{flalign}
			\label{fsfxns}
			\mathbb{P} \bigg[ \textbf{ICE}_{n'} \cap \bigcap_{|s| \le Ak^{1/3}} \Big\{ \big| \mathsf{x}_{n'} (s) - f(s) \big| \le (L')^{8/9} k^{2/3} \Big\} \bigg] \le (C_1 + 2 \Xi) e^{-c_1 (\log k)^2}.
		\end{flalign} 
	
		\noindent Thus, by the Markov inequality \Cref{fg0g} (with the event $\mathscr{E}$ there equal to that in the probability on the left side of \eqref{fsfxns} here) that there exist constants $c_2 = c_2 (\xi) > 0$ and $C_2 = C_2 (A, R, \Xi) > 1$, and an event $\mathscr{E}_0$, measurable with respect to the $\sigma$-algebra $\mathcal{F}_{n'-1}$ generated by $\mathcal{F}_{\ext}^{\bm{\mathsf{x}}} \big( \llbracket 1, n' -1 \rrbracket \times (-Ak^{1/3}, Ak^{1/3}) \big)$ and $\mathcal{F}_{\ext}^{\bm{\mathsf{y}}} \big( \llbracket 1, n' -1 \rrbracket \times [-Ak^{1/3}, Ak^{1/3}] \big)$, such that the following properties hold. First, we have 
		\begin{flalign}
			\label{e01}
			\mathbb{P} [\mathscr{E}_0] \ge 1 - C_2 e^{-c_2 (\log k)^2}.
		\end{flalign}
		
		\noindent Second, conditioning on $\mathcal{F}_{n'-1}$ and restricting to $\mathscr{E}_0$, we have (since $\mathsf{x}_{n'}$ is measurable with respect to $\mathcal{F}_{n'-1}$) that
		\begin{flalign}
			\label{xnfxn} 
			f(s)=\mathsf{y}_{n'}(s) \ge \mathsf{x}_{n'} (s) - (L')^{8/9} k^{2/3}.
		\end{flalign}
		
		\noindent Third, again restricting to $\mathscr{E}_0$, we have 
		\begin{flalign}
			\label{nn2}
			\mathbb{P} [ \textbf{ICE}_{n'} | \mathcal{F}_{n'-1}] \ge 1 - C_2 e^{-c_2 (\log k)^2}. 
		\end{flalign}
		
		By \eqref{e01} and again the Markov estimate \Cref{fg0g}, it follows that there exists an event $\mathscr{E}_1$ measurable with respect to the $\sigma$-algebra $\mathcal{F}_{n'}$ generated by $\mathcal{F}_{\ext}^{\bm{\mathsf{x}}} \big( \llbracket 1, n' \rrbracket \times (-Ak^{1/3}, Ak^{1/3}) \big)$ and $\mathcal{F}_{\ext}^{\bm{\mathsf{y}}} \big( \llbracket 1, n' \rrbracket \times (-Ak^{1/3}, Ak^{1/3}) \big)$ such that, denoting $c_3 = c_2 / 2$, the following two properties hold. First, we have 
		\begin{flalign} 
			\label{e11} 
			\mathbb{P} [\mathscr{E}_1] \ge 1 - C_3 e^{-c_3 (\log k)^2}. 
		\end{flalign} 
	
		\noindent  Second, restricting to $\mathscr{E}_1$, we have 
		\begin{flalign} 
			\label{e112} 
			\mathbb{P} [\mathscr{E}_0 | \mathcal{F}_{n'} ] \ge 1 - C_3 e^{-c_3 (\log k)^2}.
		\end{flalign} 
		
		Now, condition on $\mathcal{F}_{n'}$, and restrict to the event $\mathscr{E}_1$. If $\mathscr{E}_0$ holds, then apply \Cref{p:comparison2} with the $(\alpha; \xi, \Xi; \bm{\mathsf{x}}; \bm{\mathsf{y}})$ there equal to $(8/9; c_1, C_3; \bm{\mathsf{x}}, \bm{\mathsf{z}})$ here, which yields constants $\zeta$, $C_1$, and $C_2$ satisfying the conditions stated there. Observe that the hypotheses \eqref{probabilityeventc0}, and \eqref{yxpl} of \Cref{p:comparison2} hold by \eqref{nn2} and \eqref{xnfxn}, respectively (and we assumed \eqref{ldeltal} to hold).
		
		Hence, \Cref{p:comparison2} applies and yields on $\mathscr{E}_0 \cap \mathscr{E}_1$ a coupling between $\bm{\mathsf{x}}$ and $\bm{\mathsf{y}}$ such that \eqref{yx1} holds (where here we use the fact that $\mathsf{z}_j = \mathsf{y}_j$ for $j \in \llbracket 1, n - 1 \rrbracket$). Together with the probability bounds \eqref{e11} and \eqref{e112} for $\mathscr{E}_1$ and $\mathscr{E}_0$, this establishes the first statement of the corollary.     
	\end{proof}

	\subsection{Heuristic for \Cref{p:comparison2}}
	
	\label{Compare00}
	
	In this section we provide a heuristic for the preliminary coupling \Cref{p:comparison2}, before giving its careful proof; this will consist of a rescaling, together with an induction. We will be quite informal here, by ignoring various error terms; imagining that high probability events occur deterministically; and assuming that certain recursions behave as expected (none of which will end up requiring much conceptual effort to justify).
	
	Throughout, we set $A=1$ to ease notation. We also assume that we likely have $\mathsf{x}_j (s) \le -j^{2/3}$ for all $(j, s) \in \llbracket 1, n \rrbracket \times [-k^{1/3} k^{1/3}]$. This is almost implied by the fact that $\bm{\mathsf{x}}$ likely satisfies (see \Cref{ctr}) the $\textbf{CTR}$ event \eqref{eventcompleter0}, but is slightly different since in \eqref{eventcompleter0} there is a constant ${B}$ and an additive shift ${B} k^{2/3}$, which we are disregarding (neither of these are too serious, though they simplify some of the computations below).

	\subsubsection{Rescaling}

	As briefly explained in \Cref{CoupleBoundary2}, the proof of \Cref{p:comparison2} will proceed by comparing $\bm{\mathsf{x}}$ to a rescaled and translated version of $\bm{\mathsf{y}}$. In what follows, we let $\vartheta \in (0, 1)$ and $\mathfrak{V} > 0$ be real numbers (determining the rescaling and translation, respectively); both will be fixed later. Define the line ensemble $\bm{\mathsf{z}} = (\mathsf{z}_1, \mathsf{z}_2, \ldots , \mathsf{z}_{n'-1}) \in \llbracket 1, n'-1 \rrbracket \times \mathcal{C} \big( [-(1-\vartheta)^2 k^{1/3}, (1-\vartheta)^2 k^{1/3}] \big)$ by setting 
	\begin{flalign*}
		\mathsf{z}_j (s) = (1 - \vartheta) \cdot \mathsf{y}_j \big( (1 - \vartheta)^{-2} s \big) + \mathfrak{V}; \qquad \widetilde{f} (s) = (1 - \vartheta) \cdot f \big( (1 - \vartheta)^{-2} s \big) + \mathfrak{V},
	\end{flalign*}
	
	\noindent for each $j \in \llbracket 1, n'-1 \rrbracket$ and $s \in \big[ -(1-\vartheta)^2 k^{1/3}, (1-\vartheta)^2 k^{1/3} \big]$. Denote $\bm{u}' = \bm{\mathsf{x}}_{\llbracket 1, n'-1 \rrbracket} \big(-(1-\vartheta)^2 k^{1/3} \big)$; $\bm{v}' = \bm{\mathsf{x}}_{\llbracket 1, n'-1 \rrbracket} \big( (1-\vartheta)^2 k^{1/3} \big)$; $\widetilde{\bm{u}} = \bm{\mathsf{z}} \big(- (1-\vartheta)^2 k^{1/3} \big)$; and $\widetilde{\bm{v}} = \bm{\mathsf{z}} \big( (1-\vartheta)^2 k^{1/3} \big)$. By the invariance of Brownian bridges under diffusive scaling, the law of $\bm{\mathsf{z}}$ is given non-intersecting Brownian bridges on $\big[ -(1-\vartheta)^2 k^{1/3}, (1-\vartheta)^2 k^{1/3}]$, conditioned to start at $\widetilde{\bm{u}}$, end at $\widetilde{\bm{v}}$, and stay above $\widetilde{f}$. 
		
		For an appropriate choice of $(\vartheta, \mathfrak{V}, n')$, let us show that $\bm{\mathsf{x}}$ and $\bm{\mathsf{z}}$ can be coupled so that $\mathsf{z}_j (s) \ge \mathsf{x}_j (s)$ for each $j \in \llbracket 1, n'-1 \rrbracket$ and $s \in \big[ -(1-\vartheta)^2 k^{1/3}, (1-\vartheta)^2 k^{1/3} \big]$, with high probability. By height monotonicity \Cref{monotoneheight}, it suffices to show that we likely have
		\begin{flalign}
			\label{uuvvxf}
			\widetilde{\bm{u}} \ge \bm{u}'; \qquad \widetilde{\bm{v}} \ge \bm{v}'; \qquad \widetilde{f} \ge \mathsf{x}_{n'}.
		\end{flalign}
		
		To that end, observe that 
		\begin{flalign}
			\label{ujuj0} 
			\begin{aligned} 
			\widetilde{u}_j & = (1-\vartheta) \cdot \mathsf{x}_j ( -k^{1/3} ) + \mathfrak{V}  \ge u_j' - \vartheta j^{2/3} + \mathfrak{V} - \Big| \mathsf{x}_j (-k^{1/3}) - \mathsf{x}_j \big( -(1-\vartheta)^2 k^{1/3} \big) \Big|,
			\end{aligned}
		\end{flalign}
		
		\noindent where in the first statement we used the definition of $\widetilde{\bm{u}}$ and $\bm{\mathsf{z}}$, with the fact that $\mathsf{y}_j (-k^{1/3}) = \mathsf{x}_j (-k^{1/3})$; in the second, we used the upper bound $\mathsf{x}_j (-k^{1/3}) \le -j^{2/3}$, with the definition of $u_j' = \mathsf{x}_j \big( -(1-\vartheta)^2 k^{1/3} \big)$. To estimate the right side of \eqref{ujuj0}, we must upper bound $\big| \mathsf{x}_j (-k^{1/3}) - \mathsf{x}_j (-(1-\vartheta)^2 k^{1/3}) \big|$, which explains why we introduced (and imposed) the H\"{o}lder regularity events from \Cref{eventregularityimproved}. If $j \le L^{3S\chi/2} k$, then since $\bm{\mathsf{x}}$ satisfies the $\textbf{FHR}$ event from \eqref{fhr1}, we deduce  
		\begin{flalign*}
			\Big| \mathsf{x}_j (-k^{1/3}) - \mathsf{x}_j \big( -(1-\vartheta)^2 k^{1/3} \big)\Big| = \mathcal{O} ( (L^{(S/2+1)\chi}\vartheta+  L^{3S \chi/4}\vartheta^{1/2})k^{2/3}) = \mathcal{O} \big( (L')^{\beta} k^{2/3} \vartheta^{1/2} \big),
		\end{flalign*}
		
		\noindent where in the last inequality we used the facts that $S\geq 1$, $L' \ge L^{4S\chi}$ (by \eqref{ldeltal}) and that $\beta \ge 3/8$. If $j \ge L^{3S\chi/2} k$, then since $\bm{\mathsf{x}}$ satisfies the $\textbf{SHR}$ event from \eqref{fhr1}, we deduce that 
		\begin{flalign}
			\label{xj0}
			\begin{aligned} 
			\Big| \mathsf{x}_j (-k^{1/3}) - \mathsf{x}_j \big( -(1-\vartheta)^2 k^{1/3} \big) \Big| & = \mathcal{O} \big( k^{2/3} (jk^{-1})^{1/2} \vartheta |\log \vartheta|^2 + k^{2/3} (jk^{-1})^{2\beta/3} \vartheta^{1/2} \big) \\ 
			& = \mathcal{O} \big( (L')^{3/4} \vartheta |\log \vartheta|^2 k^{2/3} + (L')^{\beta} \vartheta^{1/2} k^{2/3} \big),
			\end{aligned}
		\end{flalign}		
		
		\noindent where in the last inequality we used the fact that $(jk^{-1})^{2/3} \le L'$. Thus, \eqref{xj0} holds in either case (see the second part of \Cref{xnf2} below). In order for $\widetilde{u}_j \ge u_j'$, by \eqref{ujuj0}, it therefore suffices to take (see \eqref{e:defz1} below), for some sufficiently large constant $C>1$,
		\begin{flalign}
			\label{r0} 
			\mathfrak{V} = Ck^{2/3} \big( (L')^{3/4} \vartheta |\log \vartheta|^2 + (L')^{\beta} \vartheta^{1/2} \big).
		\end{flalign}  
		
		 Similarly, \eqref{r0} guarantees that $\widetilde{v}_j \ge v_j'$. To show that $\widetilde{f} \ge \mathsf{x}_{n'}$, observe that 
		\begin{flalign}
			\label{fsxns0} 
			\begin{aligned}
			\widetilde{f}(s) & = (1 - \vartheta) \cdot f \big( (1-\vartheta)^{-2} s \big) + \mathfrak{V}   \\
			& \ge (1 - \vartheta) \cdot \mathsf{x}_{n'} (s) -  (L')^{\alpha} k^{2/3}  + \mathfrak{V} -  \Big| \mathsf{x}_{n'} \big( (1-\vartheta)^{-2} s \big) - \mathsf{x}_{n'} (s) \Big|.
			\end{aligned} 
		\end{flalign}
	
		\noindent where in the first statement we used the definition of $\widetilde{f}$ and in the second we used \eqref{yxpl}. By \eqref{r0} and (reasoning similar to that used to derive) \eqref{xj0}, we have 
		\begin{flalign}
			\label{rxnxn} 
			\mathfrak{V} \ge \Big| \mathsf{x}_{n'} \big( (1 - \vartheta)^{-2} s \big) - \mathsf{x}_{n'} (s) \Big|.
		\end{flalign}
	
		\noindent We also have $-\vartheta \cdot \mathsf{x}_{n'} (s) \ge \vartheta (n')^{2/3} = \vartheta L' k^{2/3}$ (by the upper bound $\mathsf{x}_j (s) \le -j^{2/3}$ and the fact that $n' = (L')^{3/2} k$). Inserting this with \eqref{rxnxn} into \eqref{fsxns0} yields that $\widetilde{f}(s) \ge \mathsf{x}_{n'} (s)$, if for some sufficiently large constant $C > 1$ we take (see \eqref{bplalpha})
		\begin{flalign}
			\label{thetaalpha0} 
			\vartheta = C (L')^{\alpha-1}.
		\end{flalign}
		
		 Therefore, under \eqref{thetaalpha0} and \eqref{rxnxn}, we have \eqref{uuvvxf} and thus that there exists a coupling between $\bm{\mathsf{x}}$ and $\bm{\mathsf{z}}$ such that $\mathsf{x}_j (s) \le \mathsf{z}_j (s)$ for each $(j,s) \in \llbracket 1, n' \rrbracket \times \big[ - (1-\vartheta)^2 k^{1/3}, (1-\vartheta)^2 k^{1/3}]$ (see \Cref{zxcouple0}). This and \eqref{rxnxn} together imply that, for any $(j, s) \in \llbracket 1, n' \rrbracket \times [-k^{2/3},  k^{2/3} ]$,
		\begin{flalign*}
			\mathsf{x}_j (s) \le \mathsf{x}_j \big( (1-\vartheta)^2 s \big) + \mathfrak{V}  \le \mathsf{z}_j \big( (1-\vartheta)^2 s \big) + \mathfrak{V} = (1-\vartheta) \cdot \mathsf{y}_j (s) + 2\mathfrak{V}.
		\end{flalign*}	
	
		\noindent Since $\vartheta \ll 1$, we find $\mathsf{y}_j (s) \ge (1 - \vartheta)^{-1} \cdot \mathsf{x}_j (s) - 3 \mathfrak{V}$, and so (ignoring the $|\log \vartheta|^2$ term in \eqref{xj0})
		\begin{flalign}
			\label{yjxj0} 
			\mathsf{y}_j (s) & \ge \mathsf{x}_j (s) - \mathcal{O} \Big(  (L')^{\alpha-1} \cdot j^{2/3} + k^{2/3} \big( (L')^{\alpha-1/4} + (L')^{\beta + \alpha/2-1/2} \big) \Big),
		\end{flalign}
	
		\noindent where we used that $-\mathsf{x}_j (s) = \mathcal{O} (j^{2/3})$ (as $\bm{\mathsf{x}}$ likely satisfies the $\textbf{CTR}$ event \eqref{eventcompleter0}, by \Cref{ctr}). 
		
		Now, recall that \Cref{p:comparison2} stated, if $L'' \le (L')^{\zeta}$ (for some small but fixed $\zeta > 0$) and $j \le (L'')^{3/2} k$, then we likely have $\mathsf{y}_j (s) \ge \mathsf{x}_j (s) - \mathcal{O} \big( (L'')^{2\beta-7/8} k^{2/3} \big)$. The estimate \eqref{yjxj0} does not immediately attain this, as the error on its right side at least $(L')^{\alpha-1/4} k^{2/3}$, which can be much larger than $(L')^{2\beta-7/8} \gg (L'')^{2\beta-7/8}$ (say, if\footnote{This is in fact a potentially relevant choice of parameters, as \Cref{xfe} indicates that $\alpha$ ``starts'' at $7/8$, while we will eventually take $\beta \le 1/2$ (as suggested by \Cref{p:Lip}, where $\beta=3/8$).} $(\alpha, \beta) = (7/8, 1/2)$).
		
		\subsubsection{Induction} 
		
		Still, \eqref{yjxj0} does yield an improvement over a direct application of height monotonicity. Indeed, since the lower boundary $f$ of $\bm{\mathsf{y}}$ satisfies $f(s) \ge \mathsf{x}_{n'} (s) - \mathcal{O} \big( (L')^{\alpha} k^{2/3} \big)$ by \eqref{yxpl}, \Cref{uvv} indicates that $\mathsf{y}_j (s) \ge \mathsf{x}_j (s) - \mathcal{O} \big( (L')^{\alpha} k^{2/3} \big)$. This error is much larger than that in \eqref{yjxj0} if $j \ll (L')^{3/2} k$ (and $\alpha \ge 2\beta-1$). This suggests that one might be able to inductively apply \eqref{yjxj0} to improve that bound further.
		
		To implement this, $\widetilde{L} = (L')^d$, for some parameter $d = d(\alpha)$ to be fixed later, and denote $\widetilde{n} = \big\lceil (\widetilde{L})^{3/2} k \big\rceil$. Then, taking $j = \widetilde{n}$ in \eqref{yjxj0} gives 
		\begin{flalign*}
			\mathsf{y}_{\tilde{n}} (s) \ge \mathsf{x}_{\tilde{n}} (s) - \mathcal{O} \big( ( \widetilde{L}^{(\alpha+d-1)/d} + \widetilde{L}^{(4\alpha-1) / 4d}  + \widetilde{L}^{(2\beta+\alpha-1)/2d}) \cdot k^{2/3} \big).
		\end{flalign*}
	
		\noindent Optimizing the exponent of $\widetilde{L}$ in the above error, over the choice of $d$, leads us to set
		\begin{flalign*}
			d = \displaystyle\max \Big\{ \displaystyle\frac{3}{4}, \beta + \displaystyle\frac{1-\alpha}{2} \Big\}.
		\end{flalign*} 
	
		\noindent This gives the bound $\mathsf{x}_{\tilde{n}} (s) \ge \mathsf{y}_{\tilde{n}} (s) - \mathcal{O} ( \widetilde{L}^{\alpha'} k^{2/3})$, where
		\begin{flalign}
			\label{alpha1alpha} 
			\alpha' = \displaystyle\frac{4}{3} \Big( \alpha - \displaystyle\frac{1}{4} \Big), \quad \text{if $\alpha \ge 2\beta - \displaystyle\frac{1}{2}$}; \qquad \alpha' = \displaystyle\frac{2\beta+\alpha-1}{2\beta-\alpha+1}, \quad \text{if $\alpha < 2\beta - \displaystyle\frac{1}{2}$}.
		\end{flalign} 
		
		 One can verify that $\alpha' \le \alpha - 1/250$ if\footnote{Observe that $\alpha = 2\beta-1$ is a fixed point of the recursion \eqref{alpha1alpha}. This suggests that, for any $\omega > 0$, there exists a constant $c = c(\omega) > 0$ such that $\alpha' \le \alpha - c$ whenever $\alpha \in (2\beta-1+\omega, 1-\omega)$. This is indeed true, and one can use it to show that the exponent $2\beta - 7/8$ appearing in \eqref{e:zilow} can be replaced by any real number greater than $2\beta - 1$ (at the expense of modifying some of the constants in the statement of that proposition). However, the exponent $2\beta - 7/8$ will suffice for our purposes, so we will not pursue this further.} $\alpha \in [2\beta-9/10, 9/10]$. This yields (see \Cref{p:comparison1}) 
		 \begin{flalign}
		 	\label{xnsyns2} 
		 	\mathsf{x}_{\tilde{n}} (s) \ge \mathsf{y}_{\tilde{n}} (s) - \mathcal{O} ( \widetilde{L}^{\alpha - 1/250} k^{2/3} ),
		 \end{flalign}
	 
	 	\noindent So, \eqref{xnsyns2} may be viewed as an improvement of \eqref{yxpl}, upon realizing $\mathsf{y}_{\tilde{n}}$ as the lower boundary for $\bm{\mathsf{y}}_{\llbracket 1, \widetilde{n}-1 \rrbracket}$ (analogous to $f$ being the lower boundary for $\bm{\mathsf{y}}_{\llbracket 1, n' \rrbracket}$ in \Cref{p:comparison2}). Therefore, we may again use \eqref{xnsyns2}, replacing the $(\alpha; L', \widetilde{L})$ there with $( \alpha - 1/250; \widetilde{L}; \widetilde{L}^{d(\alpha - 1/250)})$. Applying this bound inductively, we successively decrease the exponent $\alpha$ to $2\beta - 9/10 < 2\beta - 7/8$, which yields \Cref{p:comparison2}.

	\subsection{Proof of the Preliminary Coupling}
	
	\label{ProofCouple0}
	
	In this section we prove the preliminary coupling from \Cref{s:couple}, given by \Cref{p:comparison2}. To this end, we require the following proposition, which provides a weaker coupling than the one stated in \Cref{p:comparison2}. It will be established in \Cref{Proofxy1} below. Below, we recall that $\chi = 2^{-5000}$.

	\begin{prop}\label{p:comparison1}
		
		Adopt the notation and assumptions of \Cref{p:comparison2}, except for \eqref{ldeltal}, assuming instead that $L' \ge L^{3 S\chi}$. There exists a coupling between $\bm{\mathsf{x}}$ and $\bm{\mathsf{y}}$, such that the following holds. Setting
		\begin{flalign}
			\label{dalpha} 
			d(\alpha) = \displaystyle\max \Big\{ \displaystyle\frac{3}{4}, \beta + \displaystyle\frac{1-\alpha}{2} \Big\}; \qquad \widetilde{L} = (L')^{d(\alpha)}; \qquad \widetilde{n} = \lceil \widetilde{L}^{3/2} k \rceil,
		\end{flalign}
		
		\noindent we have 
			\begin{align}
				\label{e:case1}
				\mathbb{P} \Bigg[ \bigcap_{j = 1}^{\tilde{n}} \bigcap_{|s| \le Ak^{1/3}} \big\{ \mathsf{y}_j (s)\geq \mathsf{x}_j (s)- \widetilde L^{\alpha - 1/250} k^{2/3} \big\} \Bigg] \ge 1 - 2 \Xi e^{-\xi (\log k)^2}.
			\end{align}
			
	\end{prop}

	We establish \Cref{p:comparison2} by repeated application of \Cref{p:comparison1}.

	\begin{figure}
	\center
\includegraphics[scale = .75, trim=0 1cm 0 1cm, clip]{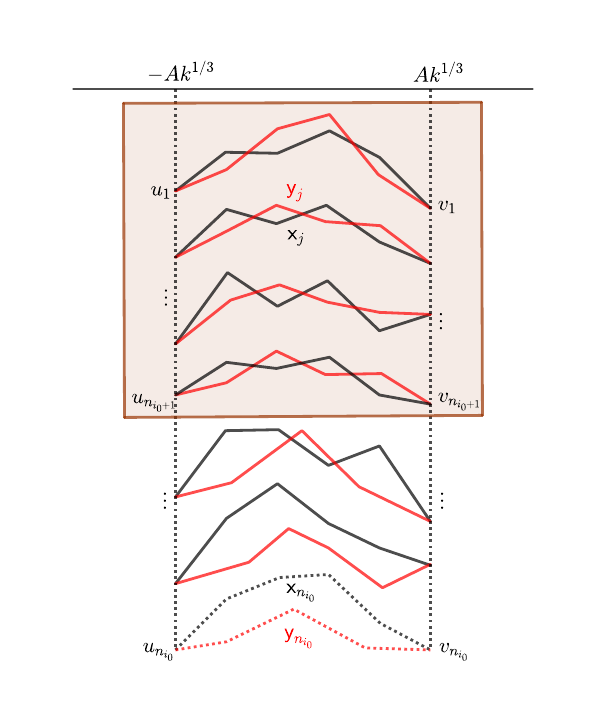}

\caption{
Shown above is a depiction of the inductive argument used in the proof of \Cref{p:comparison2}.
}
\label{f:iteration}
	\end{figure}

	\begin{proof}[Proof of \Cref{p:comparison2}]
		
		Throughout this proof, for each integer $m \in \llbracket L^{3 S \chi/2} k, n \rrbracket$, we abbreviate $\textbf{ICE}_m = \textbf{ICE}_m^{\bm{\mathsf{x}}} (A, B; \beta; R; S)$, and set $\xi'=\xi/2$. Define the sequence $\bm{\alpha} = (\alpha_0, \alpha_1, \ldots )$ by setting $\alpha_j = \alpha - j / 250$ for each $j \ge 0$. Also define the sequence $\bm{d} = (d_0, d_1, \ldots)$ by setting $d_0 = 1$ and 
		\begin{flalign*} 
			d_j = \max \Big\{ \displaystyle\frac{3}{4}, \beta + \displaystyle\frac{1 - \alpha_{j-1}}{2} \Big\}, \qquad \text{for each $j \ge 1$}. 
		\end{flalign*} 
	
		\noindent Observe that $d_i \in [3/4, 1)$ for each $i \ge 0$ satisfying $\alpha_{i-1} > 2\beta-1$. For each integer $i \ge 0$, set
		\begin{flalign}
			\label{ljnjxij}
			L_i = (L')^{ \prod_{j=0}^i d_j}; \qquad n_i = \lceil L_i^{3/2} k \rceil; \qquad \Xi_i = 3^i \Xi. 
		\end{flalign}
		
		\noindent Further let $M \ge 1$ denote the smallest integer such that $\alpha_M < 2 \beta - 7/8$; observe that $M \le 500$, as $M / 250 = \alpha_0 - \alpha_M \le 2$. Fix $\zeta = \prod_{j=0}^M d_j$, so that $\zeta \ge (3/4)^M \ge 2^{-250}$ and $L_M = (L')^{\zeta}$.

		In what follows, we omit the ceilings in the definition \eqref{ljnjxij} of $n_i$, assuming that $n_i = L_i^{3/2} k$, as this will barely affect the proofs below. We claim for each integer $i \in \llbracket 0, M \rrbracket$ that it is possible to couple $\bm{\mathsf{x}}$ and $\bm{\mathsf{y}}$ such that 
		\begin{flalign}
			\label{probabilityjxy} 
			\mathbb{P} \Bigg[ \bigcap_{j=1}^{n_i} \bigcap_{|s| \le Ak^{1/3}} \big\{ \mathsf{y}_j (s) \ge \mathsf{x}_j (s) -  L_i^{\alpha_i} k^{2/3} \big\} \Bigg] \ge 1 - \Xi_i e^{-\xi' (\log k)^2}. 
		\end{flalign}
		
		\noindent To this end, we induct on $i \in \llbracket 0, M \rrbracket$, beginning with the case $i = 0$. Since the laws of $\bm{\mathsf{x}}_{\llbracket 1, n' \rrbracket}$ and $\bm{\mathsf{y}}$ are given by $\mathsf{Q}_{\mathsf{x}_{n'+1}}^{\bm{u}; \bm{v}}$ and $\mathsf{Q}_f^{\bm{u}; \bm{v}}$, respectively, \eqref{yxpl} and \Cref{monotoneheight} together yield a coupling between $\bm{\mathsf{x}}$ and $\bm{\mathsf{y}}$ such that $\mathsf{y}_j (s) \ge \mathsf{x}_j (s) - (L')^{\alpha} k^{2/3}$ for each $(j, s) \in \llbracket 1, n' \rrbracket \times [-Ak^{1/3}, Ak^{1/3}]$. Since $(n_0, L_0, \alpha_0) = (n', L', \alpha)$, this verifies \eqref{probabilityjxy} at $i = 0$. 
		
		Thus, fix some integer $i_0 \in \llbracket 0, M-1 \rrbracket$, and assume that \eqref{probabilityjxy} holds whenever $i \le i_0$; we will show it holds for $i = i_0 + 1$. By \eqref{probabilityeventc0}, the fact that $\textbf{ICE}_{n'} \subseteq \textbf{ICE}_{n_{i_0}}$ (by \Cref{nyy}, \eqref{ihrnn}, and the fact that $n = n_0 \ge n_{i_0}$) and the Markov inequality \Cref{fg0g}, there exists an event $\mathscr{E}_{i_0}$, measurable with respect to  $\mathcal{F}_{n_{i_0}}^{\bm{\mathsf{x}}} = \mathcal{F}_{\ext}^{\bm{\mathsf{x}}} \big(\llbracket 1, n_{i_0} - 1 \rrbracket \times (-Ak^{1/3}, Ak^{1/3}) \big)$, satisfying the following two properties. First, recalling that $\xi' = \xi / 2$, we have $\mathbb{P} [\mathscr{E}_{i_0}] \ge 1 - \Xi^{1/2} e^{-\xi' (\log k)^2}$. Second, conditioning on $\mathcal{F}_{n_{i_0}}^{\bm{\mathsf{x}}}$ and restricting to $\mathscr{E}_{i_0}$, we have
		\begin{flalign}
			\label{ci01} 
			\mathbb{P} \big[ \textbf{ICE}_{n_{i_0}} \big| \mathcal{F}_{n_{i_0}}^{\bm{\mathsf{x}}} \big] \ge (1 - \Xi^{1/2} e^{-\xi' (\log k)^2}) \cdot \textbf{1}_{\mathscr{E}_{i_0}}. 
		\end{flalign}
		
		Intersecting this $\mathscr{E}_{i_0}$ with the event that $\mathsf{y}_{n_{i_0}} (s) \ge \mathsf{x}_{n_{i_0}} (s)$ for each $s \in [-Ak^{1/3}, Ak^{1/3}]$ under the coupling in \eqref{probabilityjxy}, we obtain an event $\mathscr{E}_{i_0}'$, measurable with respect to the $\sigma$-algebra $\mathcal{F}_{n_{i_0}}^{\bm{\mathsf{x}}; \bm{\mathsf{y}}}$ generated by $\mathcal{F}_{\ext}^{\bm{\mathsf{x}}} \big( \llbracket 1, n_{i_0} - 1\rrbracket \times (-Ak^{1/3}, Ak^{1/3}) \big)$ and $\mathcal{F}_{\ext}^{\bm{\mathsf{y}}} \big( \llbracket 1, n_{i_0} - 1 \rrbracket \times (-Ak^{1/3}, Ak^{1/3}) \big)$, such that the following two statements hold. First, 
		\begin{flalign*}
			\mathbb{P} [ \mathscr{E}_{i_0}' ] \ge 1 - (\Xi^{1/2} + \Xi_{i_0}) e^{-\xi' (\log k)^2}.
		\end{flalign*}
		
		\noindent Second, conditioning on $\mathcal{F}_{n_{i_0}}^{\bm{\mathsf{x}}; \bm{\mathsf{y}}}$ and restricting to $\mathscr{E}_{i_0}'$, we have \eqref{ci01} and
		\begin{flalign*}
			\mathsf{y}_{n_{i_0}} (s) \ge \mathsf{x}_{n_{i_0}} (s) - L_{i_0}^{\alpha_{i_0}} k^{2/3}, \qquad \text{for each $s \in [-Ak^{1/3}, Ak^{1/3}]$}.
		\end{flalign*} 
		
		Under the same conditioning and restriction, this verifies on the event $\mathscr{E}_{i_0}'$ the bounds \eqref{probabilityeventc0} and \eqref{yxpl}, with the $(n'; L'; f; \alpha)$ there equal to $(n_{i_0}; L_{i_0}; \mathsf{y}_{n_{i_0}}; \alpha_{i_0})$ here. Thus, since $L_{i_0}  \ge (L')^{\zeta} \ge L^{3S \chi}$ is sufficiently large (using the facts that $L' \ge L^{4 S \chi / \zeta}$ by \eqref{ldeltal} and that $L \ge C$ is sufficiently large), \Cref{p:comparison1} applies, with the $(\widetilde{n}, \widetilde{L})$ there equal to $(n_{i_0+1}, L_{i_0+1})$ here. This proposition yields a coupling between $\bm{\mathsf{x}}$ and $\bm{\mathsf{y}}$ such that 
		\begin{flalign*}
			\mathbb{P} \Bigg[ \bigcap_{j=1}^{n_{i_0+1}} & \bigcap_{|s| \le Ak^{1/3}} \big\{ \mathsf{y}_j (s) \ge \mathsf{x}_j (s) - L_{i_0+1}^{\alpha_{i_0+1}} k^{2/3}  \big\} \Bigg] \\ 
			& \ge 1 - 2 (\Xi^{1/2} + \Xi_{i_0}) e^{-\xi' (\log k)^2} \ge 1 - 3 \Xi_{i_0} e^{-\xi' (\log k)^2} = 1 - \Xi_{i_0 + 1} e^{-\xi' (\log k)^2},
		\end{flalign*}
		
		\noindent where we used the definitions $\alpha_{i_0+1} = \alpha_{i_0} - 1/250$, $L_{i_0+1} = (L_{i_0})^{d_{i_0+1}}$, and $\Xi_{i_0+1} = 3 \Xi_{i_0}$. See \Cref{f:iteration} for a depiction. This verifies \eqref{probabilityjxy}.
		
		Taking $j = M$ in \eqref{probabilityjxy} and using the facts that 
		\begin{flalign*} 
			& \Xi_M = 3^M \Xi \le 3^{250} \Xi; \qquad \zeta \ge 2^{-250}; \quad n_M = \big\lceil (L')^{3 \zeta/2} k \big\rceil; \qquad L_M^{\alpha_M} = (L')^{\zeta \alpha_M} \le (L')^{\zeta (2\beta - 7 / 8)},
		\end{flalign*} 
		
		\noindent this verifies \eqref{e:zilow} and thus the proposition.
	\end{proof}

	\subsection{Proof of \Cref{p:comparison1}}
	
	\label{Proofxy1}  
	
	In this section we establish \Cref{p:comparison1}; we adopt the notation and assumptions of that proposition throughout. First observe that, since $L'\ge C_3$, we may assume that $L'$ is sufficiently large; in particular, $L' \ge (6B)^8 \ge (6B)^{1/(1-\alpha)}$. Throughout this section, we fix the real number 
	\begin{flalign}
		\label{bplalpha}
		\vartheta=2B (L')^{\al-1} \le \displaystyle\frac{1}{3}.
	\end{flalign} 
	
	\noindent To prove \Cref{p:comparison1}, we will first couple two ensembles on the interval $\big[ -(1-\vartheta)^2 Ak^{1/3}, (1-\vartheta)^2 A k^{1/3} \big]$; the first is the restriction of $\bm{\mathsf{x}}$ to this interval, and the second is a rescaled variant of $\bm{\mathsf{y}}$. 
	
\begin{figure}
	\center
\includegraphics[scale = .8]{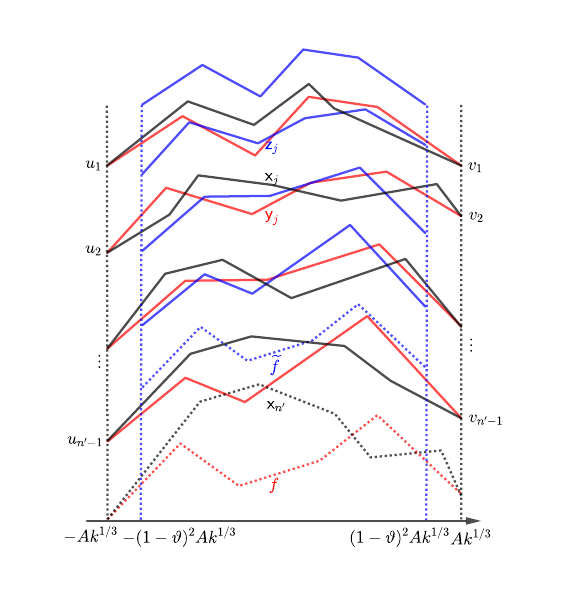}

\caption{Shown above is a depiction of \Cref{Proofxy1}, where $\bm{\mathsf z}$ is a rescaled variant of $\bm{\mathsf y}$, and with high probability we can couple $\bm{\mathsf z}\geq \bm{\mathsf x}$.}
\label{f:rescaling}
	\end{figure}

	To make the latter more explicit, define the line ensemble $\bm{\mathsf{z}} = (\mathsf{z}_1, \mathsf{z}_2, \ldots , \mathsf{z}_{n'-1}) \in \llbracket 1, n'-1 \rrbracket \times \mathcal{C} \big( [-(1-\vartheta)^2 Ak^{1/3}, (1-\vartheta)^2 Ak^{1/3}] \big)$ (see \Cref{f:rescaling}) by for each $(j, s) \in \llbracket 1, n'-1 \rrbracket \times \big[ -(1-\vartheta)^2 Ak^{1/3}, (1-\vartheta)^2 Ak^{1/3} \big]$ setting
	\begin{align}\label{e:defz1}
		&\mathsf{z}_j (s) =(1-\vartheta) \cdot \mathsf{y}_j \big( (1-\vartheta)^{-2} s \big)+ 2 k^{2/3} (R^2 + B) \cdot  \big( (L')^{3/4}\vartheta |\log \vartheta|^2 + (L')^{\beta} \vartheta^{1/2} \big).
	\end{align}
	
	\noindent See \Cref{f:rescaling} for a depiction.
	
	\begin{lem} 
		
		\label{zxcouple0} 
		
		There exists a constant $C = C(A, B, R) > 1$ such that, if $L' > C$, then there exists a coupling between $\bm{\mathsf{x}}$ and $\bm{\mathsf{z}}$ such that
		\begin{flalign*}
			\mathbb{P} \Bigg[ \bigcap_{j=1}^{n'-1} \bigcap_{|s| \le (1-\vartheta)^2 Ak^{1/3}} \big\{ \mathsf{z}_j (s) \ge \mathsf{x}_j (s) \big\} \Bigg] \ge 1 - \Xi e^{-\xi(\log k)^2}. 
		\end{flalign*}
	\end{lem}

	To establish \Cref{zxcouple0}, we require the following lemma. Below, we abbreviate the event $\textbf{ICE}_{n'} = \textbf{ICE}_{n'}^{\bm{\mathsf{x}}} (A; B; \beta; R; S)$.

	\begin{lem} 
		
		\label{xnf2} 
		
		There exists a constant $C = C(A, B, R) > 1$ such that the following hold on the event $\textbf{\emph{ICE}}_{n'}$ if $L' > C$. 
		\begin{enumerate} 
			\item For each $s \in [-Ak^{1/3}, Ak^{1/3}]$, we have $(1 - \vartheta) \cdot f(s) \ge \mathsf{x}_{n'} (s)$. 
			\item For each $(j, s) \in \llbracket 1, n' \rrbracket \times [-Ak^{1/3}, Ak^{1/3}]$, we have 
			\begin{align}\label{e:fbb1}
				\Big| \sfx_j \big ((1-\vartheta)^2 s \big)-\sfx_j (s) \Big| \le 2R^2 k^{2/3} \cdot \big( (L')^{3/4}\vartheta |\log \vartheta|^2+ (L')^{\beta} \vartheta^{1/2} \big),
			\end{align}
		\end{enumerate} 
		
	\end{lem}
	
	\begin{proof} 
		First observe by \Cref{nyy} (and \Cref{eventlocation} for the $\textbf{LOC}$ event) that, on $\textbf{ICE}_{n'}$, we have 
		\begin{flalign} 
			\label{xnsb} 
			\mathsf{x}_{n'} (s) \le Bk^{2/3} - B^{-1} (n')^{2/3} = (B - B^{-1} L') k^{2/3}, \qquad \text{for each $s \in [-Ak^{1/3}, Ak^{1/3}]$}.
		\end{flalign} 
		
		\noindent Hence, on $\textbf{ICE}_{n'}$ we have that
		\begin{align*}
			(1-\vartheta) \cdot f (s)
			&\geq (1-\vartheta) \cdot \big( \sfx_{n'} (s)-  (L')^{\alpha} k^{2/3} \big) \\
			& \geq \sfx_{n'} (s)-2B (L')^{\alpha-1} \cdot \sfx_{n'} (s)-(L')^\alpha k^{2/3}\\
			&\geq \sfx_{n'} (s)+2B (L')^{\alpha-1} k^{2/3} \cdot (B^{-1}L'-B)-(L')^\alpha k^{2/3}\geq \sfx_{n'} (s),
		\end{align*}
		
		\noindent where the first bound follows from \eqref{yxpl}, the second from \eqref{bplalpha}, the third from \eqref{xnsb}, and the fourth from the fact that $L' \ge 2B^2$ is sufficiently large. This verifies the first statement of the lemma. 
		
		By \Cref{nyy} (and \eqref{ihrnn} and \eqref{fhr1} for the $\textbf{IHR}$ and $\textbf{SHR}$ events), we have on $\textbf{ICE}_{n'}$ that, for each $(j, s) \in \llbracket L^{3S\chi/2} k, n' \rrbracket \times [-Ak^{1/3}, Ak^{1/3}]$ (recall that $\chi = 2^{-5000}$), 
		\begin{align*}
			\Big| \sfx_j \big ((1-\vartheta)^2 s \big)-\sfx_j (s) \Big| & \leq R k^{2/3} \cdot \bigg( \Big( \displaystyle\frac{j}{k} \Big)^{1/2} t_0 \big| \log (Rt_0^{-1}) \big|^2 + \Big(\displaystyle\frac{j}{k} \Big)^{2\beta/3} t_0^{1/2} + k^{-D} \bigg) \\
			& \le  2 R^2 k^{2/3} \cdot \big( (L')^{3/4}\vartheta |\log \vartheta|^2+ (L')^{\beta} \vartheta^{1/2} \big),
		\end{align*}
		
		\noindent where we have denoted $t_0 = k^{-1/3} \big|s - (1-\vartheta)^2 s  \big|$. Here, we used the facts that $jk^{-1} \le (L')^{3/2}$, that $t_0 k^{1/3} = \big|s - (1 - \vartheta)^2 s \big|\le 2 \vartheta |s| \le 2 \vartheta A k^{1/3} \le \vartheta R k^{1/3}$, that $t_0 \big| \log (R t_0^{-1}) \big|^2 \le \vartheta R |\log \vartheta|^2$ for $t_0 \in [0, \vartheta R]$ and $L'$ sufficiently large (so that $\vartheta$ is sufficiently small), and that $(L')^{\beta} \vartheta^{1/2} \ge \vartheta^{1/2} \ge (L')^{-1/2} \ge k^{-D}$ (as $L' \le L \le k^D$). This verifies \eqref{e:fbb1} for $j \in \llbracket L^{3S\chi/2} k, n' \rrbracket \times [-Ak^{1/3}, Ak^{1/3}]$. Similarly, by \Cref{nyy} (and \eqref{ihrnn} and \eqref{fhr1} for the $\textbf{IHR}$ and $\textbf{FHR}$ events), we have on $\textbf{ICE}_{n'}$ that, for each $(j, s) \in \llbracket 1, L^{3S\chi/2} k \rrbracket \times [-Ak^{1/3}, Ak^{1/3}]$, 
		\begin{align*}
			& \Big| \sfx_j \big( (1-\vartheta)^2 s \big)-\sfx_j (s) \Big| \leq 2A k^{2/3} \cdot \big( L^{3S\chi/2}\vartheta + 5 L^{S\chi} \vartheta^{1/2} \big)	\leq R k^{2/3} \cdot \big( (L')^{3/4}\vartheta + (L')^{\beta} \vartheta^{1/2}  \big).
		\end{align*}
		
		\noindent Here, in the first estimate we used the facts that we have $\big| (1-\vartheta)^2 s - s \big| \le 2 A \vartheta k^{1/3}$; that $L^{\chi} \cdot \max \{ j^{1/3} k^{-1/3}, 1 \} \le L^{(S/2+1)\chi} \le L^{3S\chi/2}$ (since $j \le L^{3S\chi/2} k$ and $S \ge 1$); that $\max \{ j^{1/2} k^{-1/2} , 1 \} \le (L')^{S \chi}$; and that $L^{S \chi} \vartheta^{1/2} \ge \vartheta^{1/2} \ge (L')^{-1/2} \ge k^{-D}$. In the second, we used the facts that $R \ge 2A$, that $\beta\geq 3/8$, and that $L'\geq L^{3S \chi} \ge L^{S\chi/ \beta}$ (with the fact that $L'$ is sufficiently large). It follows that \eqref{e:fbb1} holds on $\textbf{ICE}_{n'}$ for each $(j, s) \in \llbracket 1, n' \rrbracket \times [-Ak^{1/3}, Ak^{1/3}]$, verifying the second part of the lemma.
	\end{proof} 
	
	\begin{proof}[Proof of \Cref{zxcouple0}]
		
		Set $\mathsf{T} = (1-\vartheta)^2 Ak^{1/3}$, and condition on $\mathcal{F}_{\ext}^{\bm{\mathsf{x}}} \big( \llbracket 1, n' -1 \rrbracket \times (-\mathsf{T}, \mathsf{T}) \big)$. Define the $(n'-1)$-tuples $\bm{u}', \bm{v}', \widetilde{\bm{u}}, \widetilde{\bm{v}} \in \mathbb{W}_{n'-1}$ and the function $\widetilde{f} : [-\mathsf{T}, \mathsf{T}] \rightarrow \mathbb{R}$ by for each $s \in [-\mathsf{T}, \mathsf{T}]$ setting 
		\begin{flalign}
			\label{u2v2f}
			\begin{aligned}
				& \bm{u}' = \bm{\mathsf{x}}_{\llbracket 1, n'-1 \rrbracket} (-\mathsf{T}); \qquad \bm{v}' = \bm{\mathsf{x}}_{\llbracket 1, n'-1 \rrbracket} (\mathsf{T}); \qquad \widetilde{\bm{u}} = \bm{\mathsf{z}} ( -\mathsf{T}); \qquad \widetilde{\bm{v}} = \bm{\mathsf{z}} (\mathsf{T}); \\ 
				& \widetilde{f}(s) = (1-\vartheta) \cdot f \big( (1-\vartheta)^{-2} s \big) + 2k^{2/3} (R^2+B) \cdot \big( (L')^{3/4} \vartheta |\log \vartheta|^2 + (L')^{\beta} \vartheta^{1/2} \big).
			\end{aligned} 
		\end{flalign}
		
		\noindent Then, by \eqref{e:defz1} and \Cref{scale}, the law of $\bm{\mathsf{z}}$ is given by $\mathsf{Q}_{\tilde{f}}^{\tilde{\bm{u}}; \tilde{\bm{v}}}$. Hence, by \eqref{probabilityeventc0} and height monotonicity \Cref{monotoneheight} (with the fact that the law of $\bm{\mathsf{x}}_{\llbracket 1, n'-1 \rrbracket}$ is given by $\mathsf{Q}_{\mathsf{x}_{n'}}^{\bm{u}'; \bm{v}'}$), it suffices to show for sufficiently large $L'$ that 
		\begin{flalign}
			\label{uuvvff0} 
			\bm{u}' \le \widetilde{\bm{u}}, \quad \bm{v}' \le \widetilde{\bm{v}}, \quad \text{and} \quad \mathsf{x}_n \le \widetilde{f}, \qquad \text{all hold on the event $\textbf{ICE}_{n'}$}.
		\end{flalign}
	
		The proofs of the first and second statements in \eqref{uuvvff0} are entirely analogous, so we only provide that of the former. To this end, observe on $\textbf{ICE}_{n'}$ that, for each $j \in \llbracket 1, n' \rrbracket$,
		\begin{align*}
			\widetilde{u}_j & = (1-\vartheta) \cdot \sfy_j (-A k^{1/3})+ 2k^{2/3} (R^2 + B) \cdot \big( (L')^{3/4}\vartheta | \log \vartheta|^2 + (L')^{\beta} \vartheta^{1/2} \big) \\
			&= (1-\vartheta) \cdot \mathsf{x}_j (-A k^{1/3})+2k^{2/3} (R^2 + B) \cdot  \big( (L')^{3/4}\vartheta | \log \vartheta|^2 + (L')^{\beta} \vartheta^{1/2} \big) \\
			&\geq \sfx_j (- A k^{1/3})+2R^2 k^{2/3} \cdot \big((L')^{3/4}\vartheta |\log \vartheta|^2+ (L')^{\beta} \vartheta^{1/2} \big) \geq \sfx_j (-\mathsf{T}),
		\end{align*}
		
		\noindent where the first statement holds by \eqref{u2v2f} and \eqref{e:defz1}; the second by the fact $\sfy_j (-A k^{1/3})= u_j = \sfx_j (-A k^{1/3})$; the third by the bound $\vartheta \cdot \sfx_j (- A k^{1/3}) \le \vartheta B k^{2/3} \le B (L')^{\beta} \vartheta^{1/2} k^{2/3}$ on $\textbf{ICE}_{n'}$ (where the former is due to the $\textbf{LOC}$ events in \Cref{nyy}, and the latter is due to the facts that $\vartheta \le 1$ and $B, L' \ge 1$); and the fourth by (the $s = -Ak^{1/3}$ case of) \eqref{e:fbb1}. This verifies the first bound in \eqref{uuvvff0}; the proof of the second is entirely analogous and is thus omitted.
		
		To confirm the third, observe on $\textbf{ICE}_{n'}$ that, for $s \in [-Ak^{1/3}, Ak^{1/3}]$, 
		\begin{align*}
			\widetilde{f} \big( (1-\vartheta)^2 s \big) & = (1-\vartheta) \cdot f(s)+ 2 k^{2/3} (R^2 + B) \cdot \big( (L')^{3/4}\vartheta |\log \vartheta|^2 + (L')^{\beta} \vartheta^{1/2} \big)\\
			&\geq \sfx_{n'}(s)+ 2 k^{2/3} (R^2 + B) \cdot \big( (L')^{3/4}\vartheta |\log \vartheta|^2 + (L')^{\beta} \vartheta^{1/2} \big) \geq \sfx_{n'} \big((1-\vartheta)^2 s \big) ,
		\end{align*}
		where the first statement holds by \eqref{u2v2f}; the second by the first statement of \Cref{xnf2}; and the third by \eqref{e:fbb1}. This establishes \eqref{uuvvff0} and thus the lemma. 
	\end{proof}

	Now we can establish \Cref{p:comparison1}.

	\begin{proof}[Proof of \Cref{p:comparison1}]
		
		First observe that we may couple $\bm{\mathsf{x}}$ and $\bm{\mathsf{y}}$ such that the following holds with probability at least $1 - \Xi e^{-\xi (\log k)^2}$. For each $(j, s) \in \llbracket 1, n' - 1 \rrbracket \times [-Ak^{1/3}, Ak^{1/3}]$, we have 
		\begin{align*}
			(1-\vartheta) \cdot \sfy_j (& s)+ 2 k^{2/3} (R^2+B) \cdot \big( (L')^{3/4}\vartheta |\log \vartheta|^2 + (L')^{\beta} \vartheta^{1/2} \big) \\
			& =\sfz_j \big( (1-\vartheta)^2 s \big) \ge  \sfx_j \big( (1-\vartheta)^2 s \big)
			\geq \sfx_j (s)-2R^2 k^{2/3} \cdot \big ((L')^{3/4}\vartheta |\log \vartheta|^2 + (L')^{\beta} \vartheta^{1/2} \big),
		\end{align*}
		
		\noindent where the first statement holds by \eqref{e:defz1}; the second by \Cref{zxcouple0}; and the third by  \eqref{e:fbb1}. Hence, with probability at least $1 - 2 \Xi e^{-\xi (\log k)^2}$, for each $(j, s) \in \llbracket 1, n' - 1 \rrbracket \times [-Ak^{1/3}, Ak^{1/3}]$ we have
		\begin{align}\begin{split}\label{e:xiyi}
				\sfy_j (s)&\geq (1 - \vartheta)^{-1} \cdot \sfx_j (s) - 2 (1 - \vartheta)^{-1} R^2  k^{2/3} \cdot \big( (L')^{3/4}\vartheta |\log \vartheta|^2 + (L')^{\beta} \vartheta^{1/2} \big) \\
				&\geq
				\sfx_j (s)-6(R^2+B) \Big( k^{2/3} \big((L')^{3/4}\vartheta |\log \vartheta|^2  + (L')^{\beta} \vartheta^{1/2} \big) + \vartheta (k^{2/3}+j^{2/3}) \Big),
		\end{split}\end{align}
		
		\noindent where for the second inequality, we used the facts that $\vartheta \le 1 / 3$ (by \eqref{bplalpha}) and that $(1 - \vartheta)^{-1} \cdot \sfx_j (s) \ge \mathsf{x}_j (s) - 2 \vartheta B (k^{2/3}+j^{2/3})$ holds with probability at least $1 - \Xi e^{-\xi (\log k)^2}$ (and a union bound). Indeed, the latter follows from the bound $(1 - \vartheta)^{-1} \le 1 + 2\vartheta$ for $\vartheta \le 1 / 3$, the fact that $\mathsf{x}_j (s) \ge -B(k^{2/3} + j^{2/3})$ on $\textbf{ICE}_{n'}$ (by the $\textbf{LOC}$ events in \Cref{nyy}), and \eqref{probabilityeventc0}.

		Now, set $\mathfrak{d} = 10^{-10}$. First assume that $\alpha \ge 2\beta-1/2$, in which case $d(\alpha) = 3/4$, by \eqref{dalpha}. Then, for sufficiently large $L'$ we have by \eqref{bplalpha} that
		\begin{align}\begin{split}\label{e:someup}
				&(L')^\beta \vartheta^{1/2}= (2B)^{1/2} (L')^{\beta+(\alpha-1)/2}\leq (2B)^{1/2}  (L')^{\alpha-1/4}\leq (L')^{\alpha-1/4 + \mathfrak{d}}; \\
				& (L')^{3/4}\vartheta |\log \vartheta|^2 =  2B (L')^{\alpha-1/4} |\log \vartheta|^2 \leq (L')^{\alpha-1/4 + \mathfrak{d}}.
		\end{split}\end{align}
		
		\noindent Taking $\widetilde{L} = (L')^{3/4}$ and $ \widetilde n=\lceil \widetilde{L}^{3/2}k\rceil$, we also have for sufficiently large $L'$ that
		\begin{flalign*} 
			\vartheta(k^{2/3}+ j^{2/3})\leq 2B (L')^{\alpha-1} \cdot \big( (L')^{3/4} + 2 \big) k^{2/3} \le (L')^{\alpha - 1/4 + \mathfrak{d}} k^{2/3}.
		\end{flalign*} 
		
		\noindent for sufficiently large $L'$ and any integer $j \in \llbracket 1, \widetilde{n} \rrbracket$. This, together with \eqref{e:xiyi} and \eqref{e:someup}, gives for each $(j, s) \in \llbracket 1, \widetilde{n} \rrbracket \times [-Ak^{1/3}, Ak^{1/3}]$ that
		\begin{align*}
			\sfy_j (s)\geq \sfx_j (s)-  18(R^2+B) (L')^{\alpha-1/4+\fd} k^{2/3} = \mathsf{x}_j (s) - 18(R^2+B) \widetilde{L}^{(4\alpha-1+4 \mathfrak{d})/3} k^{2/3}.
		\end{align*}
		
		\noindent Since $(4\alpha - 1 + 4 \mathfrak{d}) / 3 \le \alpha - 1/200$ (as $\alpha < 9/10$ and $\mathfrak{d} = 10^{-10}$), this establishes \eqref{e:case1} when $\alpha \ge 2\beta - 1/2$.
		
		Next assume that $\alpha < 2\beta-1/2$, in which case $d(\alpha) = \beta + (1-\alpha)/2$, by \eqref{dalpha}. Then, for sufficiently large $L'$ we have by \eqref{bplalpha} that 
		\begin{align}\label{e:someup2}
			(L')^{\beta + \mathfrak{d}} \vartheta^{1/2} = (2B)^{1/2}  (L')^{\beta+(\alpha-1)/2 + \mathfrak{d}} \geq  (2B)^{1/2} (L')^{\alpha-1/4 + \mathfrak{d}} \geq (L')^{3/4}\vartheta |\log \vartheta|^2. 
		\end{align}
		
		\noindent Taking $\widetilde{L} = (L')^{\beta + (1-\alpha)/2}$ and $ \widetilde n=\lceil \widetilde{L}^{3/2} k\rceil$, we also have 
		\begin{flalign*} 
			\vartheta(k^{2/3}+j^{2/3})\leq 2B (L')^{\alpha-1} \cdot \big( (L')^{\beta + (1 - \alpha)/2} + 2 \big)k^{2/3} \le 3B (L')^{\beta+(\alpha-1)/2}k^{2/3},
		\end{flalign*} 
		
		\noindent for sufficiently large $L'$ and any integer $j \in \llbracket 1, \widetilde{n} \rrbracket$.  This, together with \eqref{e:xiyi}, \eqref{e:someup2}, and \eqref{bplalpha} (and the fact that $d(\alpha) = \beta + (1-\alpha)/2 \ge 1/4$ for $\alpha \le 9/10$ and $\beta \ge 3/8$), gives for each $(j, s) \in \llbracket 1, \widetilde{n} \rrbracket \times [-Ak^{1/3}, Ak^{1/3}]$ that 
		\begin{align}
			\sfy_j (s) \geq 
			\sfx_j (s)- 36B (R^2+B) (L')^{\beta+(\alpha-1) / 2 + \mathfrak{d}} k^{2/3} \ge \mathsf{x}_j (s) -  \widetilde{L}^{(2\beta + \alpha-1) / (2\beta -\alpha + 1) + 5 \mathfrak{d}}.
		\end{align}
	
		\noindent Since $(2 \beta + \alpha - 1) / (2\beta - \alpha + 1) + 5 \mathfrak{d} \le \alpha - 1/225$ for $\alpha \in [2\beta-9/10, 9/10]$ and $\beta \in [3/8, 7/8]$ (as this is equivalent to $(5\mathfrak{d} + 1/225)(2\beta - \alpha + 1) \le (1 - \alpha)(\alpha - 2\beta + 1)$, which holds as $\mathfrak{d} = 10^{-10}$, $(1-\alpha)(\alpha-2\beta+1) \ge 1/100$, and $2\beta - \alpha + 1 \le 2$), this proves \eqref{e:case1} when $\alpha < 2\beta - 1/2$, thereby establishing the proposition.
	\end{proof}

	\section{Improved H\"{o}lder Estimates}
	
	\label{RegularImproved} 
	
	In this section we establish the improved H\"{o}lder estimate \Cref{p:Lip}, which will be based on three results. The first indicates, under \Cref{a:nkrelation}, that the improved H\"{o}lder event $\textbf{IHR}$ (from \Cref{eventregularityimproved}) likely holds at $\beta = 3 / 4$; to establish \Cref{p:Lip}, we must improve this value of $\beta$ to $3 / 8$. To this end, we define a ``density regularity event'' $\textbf{DEN}$, on which the paths in the line ensemble $\bm{\mathsf{x}}$ are well approximated by a measure with regular density; this event will also involve a parameter $\beta$, prescribing the error in the approximation. The second result we will show indicates that $\textbf{DEN}$ likely implies $\textbf{IHR}$ with an improved value of $\beta$; the third indicates that $\textbf{IHR}$ likely implies $\textbf{DEN}$ with an improved value of $\beta$. By inductively applying the latter two statements, we will improve the $\beta$ in $\textbf{IHR}$ from $3 / 4$ to $3 / 8$, establishing \Cref{p:Lip}. 
	
	We begin in \Cref{RegularImproved} by defining the regular density event $\textbf{DEN}$, formulating these three statements, and establishing \Cref{p:Lip} assuming them. We then establish the first, second, and third results mentioned above in \Cref{ProofB}, \Cref{ProofDensity1}, and \Cref{ProofDensity2}, respectively. Throughout this section, we let $\bm{\mathsf{x}} = (\mathsf{x}_1, \mathsf{x}_2, \ldots )$ denote a $\mathbb{Z}_{\ge 1} \times \mathbb{R}$ indexed line ensemble satisfying the Brownian Gibbs property; we also recall the $\sigma$-algebra $\mathcal{F}_{\ext}$ from \Cref{property}, and the location event $\textbf{LOC}$ from \Cref{eventlocation} and the boundary tall rectangle event $\textbf{BTR}$ from \Cref{ftrbtr}.

	\subsection{Proof of the Improved H\"{o}lder Estimate}
	
	\label{RegularImprovedProof}
	
	In this section we establish \Cref{p:Lip}. We begin with the following lemma, to be established in \Cref{ProofB} below, indicating that the boundary tall rectangle event $\textbf{BTR}$ of \Cref{ftrbtr} likely implies the first H\"{o}lder events $\textbf{FHR}$ of \eqref{fhr1}. We further set $\chi = 2^{-5000}$.

	\begin{lem}

		\label{p:initial}
		
		Adopting \Cref{a:nkrelation}, there exist constants $c = c(A, B) > 0$ and $C = C(A, B, D) > 1$ such that the following holds if $L \ge (2B)^{2/\chi}$. For any real number $A' \in [0, A - k^{-1/3}]$, we have (recalling \Cref{eventregularityimproved}) that
		\begin{flalign*}
			\mathbb{P} \Bigg[ \textbf{\emph{BTR}}_n^{\bm{\mathsf{x}}} (A; B) \cap \bigcup_{j=1}^n \textbf{\emph{FHR}}_j^{\bm{\mathsf{x}}} (A')^{\complement}  \Bigg] \le Ce^{-c(\log k)^2}.
		\end{flalign*}
		
	\end{lem}

	The following lemma indicates that intersections of the $\textbf{FHR}$ events are equal to the $\beta = 3/4$ cases of the improved H\"{o}lder events $\textbf{IHR}$ from \eqref{ihrnn}. So, \Cref{p:initial} can be viewed as the inital case of the induction outlined at the beginning of this section.
	
	\begin{lem} 
		
		\label{fhrihr}

		Recall $\chi = 2^{-5000}$. Fix integers $n \ge k \ge 1$, and real numbers $A, B, D, L \ge 1$; $S \in [4, 2^{1000}]$; and $R \ge \max \{ 2A, 5 \}$, such that $n = L^{3/2} k$ and $L \in [1, k^D]$. Recalling \Cref{eventregularityimproved}, we have for any integer $n' \in \llbracket L^{3S\chi/2} k, n \rrbracket$ that 
		\begin{flalign*}
			\textbf{\emph{IHR}}_{n'} \Big( A; \displaystyle\frac{3}{4}; R; S \Big) = \bigcap_{j=1}^{n'} \textbf{\emph{FHR}}_j (A). 
		\end{flalign*}
		
	\end{lem} 
	
	\begin{proof} 
		
		By \eqref{ihrnn}, it suffices to show that $ \textbf{FHR}_j (A)\subseteq \textbf{SHR}_j (A; 3 / 4; R)$ for each integer $j \in \llbracket L^{3S \chi / 2} k, n \rrbracket$. Setting $\beta = 3 / 4$, by \eqref{fhr1}, this follows from the facts that $R \ge 5$ and that, for any $(j, t) \in \llbracket L^{3S \chi / 2} k, n \rrbracket \times [-A, A]$, we have 
		\begin{flalign*}
			L^{\chi} \Big( \displaystyle\frac{j}{k} \Big)^{1/3} t + 4 \Big(\displaystyle\frac{j}{k} \Big)^{1/2} t^{1/2} + k^{-D} & \le 5 \Big( \displaystyle\frac{j}{k} \Big)^{1/2} t^{1/2} + k^{-D} \\
			& \le 5 \Bigg( \bigg( \displaystyle\frac{j}{k} \Big)^{1/2} t \Big (\log ( 5|t|^{-1} \big) \Big)^2 + \Big( \displaystyle\frac{j}{k} \Big)^{2\beta/3} t^{1/2} + k^{-D} \Bigg),
		\end{flalign*} 
		
		\noindent where the first bound holds as $L^{\chi} (k^{-1} j)^{1/3} \le (k^{-1} j)^{1/2}$ for $j \ge L^{3S\chi/2} \ge L^{6 \chi} k$ (recall $S \ge 4$).
	\end{proof}

	We next introduce the following event on which the $\bm{\mathsf{x}}_{\llbracket 1, i \rrbracket} (tk^{1/3})$ are, for each $(t, i)$, well-approximated by the classical locations with respect to a measure with a regular density (in a form similar to what is guaranteed by \Cref{p:closerho}). In what follows, we recall the classical locations with respect to a measure from \Cref{gammarho}. 
	
	\begin{definition}
		
		\label{eventrho} 
		
		Recall $\chi = 2^{-5000}$. Fix integers $n \ge k \ge 1$ and real numbers $A, D, L, R, S  \ge 1$ and $\beta\in [-1,3/4]$, with $n = L^{3/2} k$ and $L \in [1, k^D]$. For any integer $i \in \llbracket L^{3S\chi/2} k, n \rrbracket$ and real number $t \in [-A, A]$, define the \emph{regular density event} $\textbf{DEN}_{i} (t;\beta;R) = \textbf{DEN}_{i}^{\bm{\mathsf{x}}} (t; \beta; R)$ to be that on which the following holds. There exists a measure $\mu = \mu_t^{(i)}$ with $\mu (\mathbb{R}) = k^{-1} i$, satisfying the following properties. Below, we denote the classical locations of $\mu$ by $\gamma_j = \gamma_{j; i}^{\mu}$ and set $\mathfrak{m}_j = \mathfrak{m}_j (R) = \big\lceil R \log n \cdot \max \{ j^{1/2}, k^{1/2} \} \big\rceil$, for each $j \in \llbracket 1, i \rrbracket$.\index{D@$\textbf{DEN}$; regular density event}
		
		\begin{enumerate} 
			
			\item The measure $\mu$ admits a density $\varrho \in L^1 (\mathbb{R})$ with respect to Lebesgue measure satisfying $\supp \varrho \subseteq [-RL, RL^{3/4}]$ and $\varrho (x) \le R  (k^{-1} i)^{1/2}$ for each $x \in \mathbb{R}$.
			\item  For each integer $j \in \llbracket 1, i \rrbracket$, we have  
			\begin{align}\label{e:defDen}
				\gamma_{j + \mathfrak{m}_j}-R\left(\frac{i}{k}\right)^{2\beta/3} \le k^{-2/3} \cdot \sfx_j (tk^{1/3}) \le  \gamma_{j - \mathfrak{m}_j}+R\left(\frac{i}{k}\right)^{2\beta/3}.
			\end{align}
			
		\end{enumerate} 
		
		\noindent For any integer $n' \in \llbracket L^{3S\chi/2} k, n \rrbracket$, also define $\textbf{DEN}_{n'} (A; \beta; R; S) = \textbf{DEN}_{n'}^{\bm{\mathsf{x}}} (A; \beta; R; S)$ by 
		\begin{flalign*}
			\textbf{DEN}_{n'} (A; \beta; R; S) = \bigcap_{i = \lceil L^{3S\chi/2} k \rceil}^{n'} \bigcap_{|t| \le A} \textbf{DEN}_i (t; \beta; R).
		\end{flalign*}
		
		\noindent Observe that these events also implicitly depend on $k$ and $L$, but we will omit this from the notation for brevity.
		
	\end{definition}

	\begin{rem} 
		
		\label{event0density} 
		
		As in \Cref{pfl0event}, observe that those $\bm{\mathsf{x}}$ satisfying $\textbf{DEN}_i^{\bm{\mathsf{x}}} (t; \beta; R)$ define a closed (and thus measurable) subset of $\llbracket 1, n \rrbracket \times \mathcal{C} \big( [-Ak^{1/3}, Ak^{1/3}] \big)$. Indeed, the set $\mathcal{Y}$ of measures $\mu \in \mathscr{P}_{\fin}$ satisfying all properties except for \eqref{e:defDen} in \Cref{eventrho} forms a compact set (where $\mathscr{P}_{\fin}$ is given the topology of weak convergence). Therefore, the function mapping $\bm{\mathsf{x}} (tk^{1/3})$ to $\inf_{\mu \in \mathcal{Y}} \max_{j \in \llbracket 1, i \rrbracket} \big\{ k^{-2/3} \cdot x_j (tk^{1/3}) - \gamma_{j + \mathfrak{m}_j}^{\mu}, \gamma_{j - \mathfrak{m}_j}^{\mu} - k^{-2/3} \cdot \mathsf{x}_j (tk^{1/3}) \big\}$ is continuous on $\llbracket 1, n \rrbracket \times \mathcal{C} \big( [-Ak^{1/3}, Ak^{1/3}] \big)$. Since $\textbf{DEN}_i^{\bm{\mathsf{x}}} (t; \beta; R)$ is the event on which this function is at most $R(ik^{-1})^{2\beta/3}$, it is closed. By the measurable selection theorem, we can also choose $\mu \in \mathcal{Y}$ satisfying \eqref{e:defDen} as an argument infimum of the above continuous function (see \cite[Theorem 18.13]{IDA}) to be Borel measurable in $\bm{\mathsf{x}} (tk^{1/3})$. 
	
	\end{rem}

	The following two propositions provide implications between the regular density event $\textbf{DEN}$ and improved H\"{o}lder one $\textbf{IHR}$. The first, to be established in \Cref{ProofDensity1} below (and eventually amounting from \Cref{uvrho}), indicates (upon restricting to $\textbf{BTR}$) that $\textbf{DEN}$ likely implies $\textbf{IHR}$ with a different value of $\beta$ (see \Cref{proofd00} for an explanation of this value), on a slightly thinner rectangle.

	\begin{prop}\label{p:dentoReg}
		Adopting \Cref{a:nkrelation} and letting $R \ge 1$ be a real number, there exist constants $c_1 = c_1 (A, B, R) > 0$, $C_1 = C_1 (A, B, R) > 1$, and $C_2 = C_2 (A, B, R, D) > 1$ such that the following holds if $L \ge C_2$. Fix real numbers $\beta\in [0, 3/4]$, $A' \in [0, A-k^{-1/3}]$, and $S \in [4, 2^{1000}]$. For any integer $n' \in \llbracket L^{3S\chi/2} k, n \rrbracket$, we have (recalling \Cref{eventregularityimproved} and \Cref{eventrho}) that
		\begin{align}\begin{split}
				\bP \bigg[\textbf{\emph{BTR}}_n(A;B)\cap  \textbf{\emph{DEN}}_{n'}(A';\beta; R;S) \cap  \textbf{\emph{IHR}}_{n'} \Big(A'; \displaystyle\frac{\beta}{2} + \displaystyle\frac{3}{8};  C_1; S \Big)^{\complement} \bigg] \le C_2 e^{-c(\log k)^2}.
		\end{split} \end{align} 
	\end{prop}

	The second, to be established in \Cref{ProofDensity2} below, indicates that $\textbf{IHR}$ likely implies $\textbf{DEN}$ also with a different value of $\beta$ (but does so by stating that, if $\textbf{BTR}$ likely implies $\textbf{IHR}$, then it also likely implies $\textbf{DEN}$), on a slightly thinner rectangle. It will essentially be a quick consequence of \Cref{p:closerho} (which exhibits the same exponent $2\beta-7/8$).

	\begin{prop}\label{p:RegtoDen}

		Adopting \Cref{a:nkrelation} and letting $A' \in [1, A - k^{-1/3}]$; $b, \xi \in ( 0, 1 / 4)$; and $R, \Xi \ge 1$ be real numbers, there exist constants $\zeta \in [2^{-250}, 1]$, $c = c (b, A, A', B, \xi) > 0$, $C_1 = C_1 (b, A, A', B, R) > 1$, and $C_2 = C_2 (b, A, A', B, D, R, \xi, \Xi) > 0$ such that the following holds if $L \ge C_2$. Fix real numbers $\beta \in [ 3 / 8, 3 / 4]$ and $S \in [1, 2^{500}]$. Assume for some integer $n' \in \llbracket L^{3S\chi/2} k, n \rrbracket$ and real number $L' \in [C_2, L]$, such that $n' = (L')^{3/2} k$ and $(L')^{3\zeta/2} \ge L^{6S\chi}$, that  we have    
		\begin{align}
			\label{btrr} 
			&\bP \big[ \textbf{\emph{BTR}}_{n}(A;B) \cap \textbf{\emph{IHR}}_{n'}(A'; \beta; R; S)^{\complement} \big] \leq \Xi e^{-\xi (\log k)^2}.
		\end{align} 
		
		\noindent Then, denoting $\widetilde{n} = \big\lceil (L')^{3\zeta / 2} k \big\rceil$, we have (recalling \Cref{eventregularityimproved} and \Cref{eventrho}) that 
		\begin{align*} 
			&\bP \bigg[ \textbf{\emph{BTR}}_n(A;B) \cap \textbf{\emph{DEN}}_{\widetilde n} \Big((1-b)A'; 2\beta- \displaystyle\frac{7}{8}; C_1; 4S \Big)^{\complement} \bigg] \le C_2 e^{-c (\log k)^2}.
		\end{align*} 
	\end{prop}

	Given the above results, we can establish \Cref{p:Lip}. Underlying its proof is the fact that composing the functions $\beta \mapsto 2\beta - 7/8$ (from \Cref{p:RegtoDen}) and $\beta \mapsto \beta/2 + 3/8$ (from \Cref{p:dentoReg}) yields the funciton $\beta \mapsto \beta - 1/16$. In this way, applying \Cref{p:RegtoDen} and then \Cref{p:dentoReg} decreases the H\"{o}lder regularity parameter $\beta$ to $\beta - 1/16$. Implementing this six times reduces $\beta$ from $\beta = 3/4$ (recall \Cref{fhrihr} and \Cref{p:initial}) to $\beta = 3/8$, which is what is stated in \Cref{p:Lip}. Throughout the process, some of the other parameters in the $\textbf{DEN}$ and $\textbf{IHR}$ events will change, and this must be tracked (though it will not cause any issues).

	\begin{proof}[Proof of \Cref{p:Lip}] 
		
		Let $\zeta \in [2^{-250}, 1]$ be as in \Cref{p:RegtoDen}, and set $b_0 = 1 - 2^{-1/7}$. For each integer $i \in \llbracket 0, 7 \rrbracket$, set 
		\begin{flalign*} 
			A_i = (1 - b_0)^i A; \qquad \beta_i = \displaystyle\frac{3}{4} - \displaystyle\frac{i-1}{16}; \qquad S_i = 4^i; \qquad L_i = L^{\zeta^{i-1}}; \qquad n_i = \lceil L_i^{3/2} k \rceil.
		\end{flalign*} 
		
		\noindent We will omit the ceilings in what follows, assuming that $n_i = L_i^{3/2} k$, as this will barely affect the proofs; we may also suppose that $k$ is sufficiently large so that $A_i - k^{-1/3} \ge A_{i+1}$ for each $i \in \llbracket 0, 6 \rrbracket$. Observe that $n_i \ge n_7 \ge L^{3S_7 \chi/2} k$, since $2S_7 \chi = 2^{15} \chi < 2^{-1500} \le \zeta^6$ (as $\chi = 2^{-5000}$) for each $i \in \llbracket 1, 7 \rrbracket$. 
		
		We claim for each integer $i \in \llbracket 1, 7 \rrbracket$ that there exist constants $\xi_i = \xi_i (A, B) > 0$, $R_i = R_i (A, B) > 1$, and $\Xi_i = \Xi_i (A, B, D) > 1$ such that for $L > \Xi_i$ we have
		\begin{flalign}
			\label{ibtr17}
			\mathbb{P} \big[ \textbf{BTR}_n (A; B) \cap \textbf{IHR}_{n_i} (A_i; \beta_i; R_i; S_i )^{\complement} \big] \le \Xi e^{-\xi_i (\log k)^2}.
		\end{flalign}
		
		\noindent The proposition would then follow from taking $i = 7$ in \eqref{ibtr17} and using the inclusion of events $\textbf{IHR}_{n_7} (A_7; \beta_7; R_7; S_7) \subseteq \textbf{IHR}_{n'} ( A / 2; 3 / 8; R_7; 2^{14})$, which holds since $(A_7, \beta_7, S_7) = ( A / 2, 3 / 8, 2^{14})$,  since $\textbf{IHR}_{n_7} \subseteq \textbf{IHR}_{n'}$ by \Cref{eventregularityimproved} and the fact that $n_7 \ge n'$ (as $\zeta^7 \ge 2^{-1750}$). 
		
		It therefore remains to verify \eqref{ibtr17}, which we do by induction on $i$. We begin with the case $i = 1$. To this end, first observe by \Cref{p:initial} that there exist constants $\xi_1 = c_1 (A, B) > 0$ and $\Xi_1 = \Xi_1 (A, B, D) > 1$ such that  for $L > \Xi_1$ we have
		\begin{align*}
			\bP \Bigg[ \textbf{BTR}_n(A;B) \cap \bigcup_{j=1}^n \textbf{FHR}_j (A_1)^{\complement} \Bigg] \le \Xi_1 e^{-\xi_1 (\log k)^2},
		\end{align*}
		
		\noindent which together with \Cref{fhrihr} yields \eqref{ibtr17} at $i = 1$ (with $R_1 = \max \{ 2A, 5 \}$).

		Now, assume \eqref{ibtr17} holds for some integer $i \in \llbracket 1, 6 \rrbracket$, and we will show it continues to hold upon replacing $i$ with $i+1$. By \Cref{p:RegtoDen} (with the parameters $(n', L'; \widetilde{n})$ there equal to $(n_i, L_i; n_{i+1})$ here; the $(b; \beta)$ there equal to $(b_0; \beta_i)$ here; and the $\big(A', (1 - b) A'; R, S, 4S; \xi, \Xi \big)$ there equal to $(A_i, A_{i+1}; R_i, S_i, S_{i+1}; \xi_i, \Xi_i)$ here, observing that $L_i^{3\zeta/2} \ge L^{\zeta^6} \ge L^{-2^{1/1500}} \ge L^{2^{20} \chi} \ge L^{6 S_i \chi}$ as $\chi = 2^{-5000}$), there exist constants $c_1 = c_1 (A, B, \xi_i) > 0$, $C_1 = C_1 (A, B, R_i) > 1$, and $C_2 = C_2 (A, B, D, R_i, \xi_i, \Xi_i) > 1$ such that, for $L > C_2$, we have
		\begin{flalign}
			\label{btr0} 
			\mathbb{P} \bigg[ \textbf{BTR}_n (A; B) \cap \textbf{DEN}_{n_{i+1}} \Big( A_{i+1}; 2\beta_i - \displaystyle\frac{7}{8}; C_1; S_{i+1} \Big)^{\complement} \bigg] \le C_2 e^{-c_1 (\log k)^2}. 
		\end{flalign}
		
		\noindent Moreover, \Cref{p:dentoReg} (with the $(n'; \beta; A', R, S)$ there equal to $( n_{i+1}; 2\beta - 7 / 8; A_{i+1}, C_1, S_{i+1})$ here, where we observe that $2\beta_i - 7 / 8 \geq 2 \beta_6 - 7 / 8 = 0$) yields constants $c_2 = c_2 (A, B, R_i) > 0$, $R_{i+1} = R_{i+1} (A, B, C_1) > 1$, and $C_3 = C_3 (A, B, C_1, D) > 1$ such that, for $L > C_3$, we have
		\begin{flalign*} 
			\mathbb{P} \bigg[ \textbf{BTR}_n (A; B) & \cap \textbf{DEN}_{n_{i+1}} (A_{i+1}; 2\beta_i - \displaystyle\frac{7}{8}; C_1; S_{i+1}) \\
			& \cap \textbf{IHR}_{n_{i+1}} \Big( A_{i+1}; \beta_i - \displaystyle\frac{1}{16}; R_{i+1}; S_{i+1} \Big)^{\complement} \bigg] \le C_3 e^{-c_2 (\log k)^2}.
		\end{flalign*}  
		
		\noindent This together with \eqref{btr0} and union bound (with the fact that $\beta_{i+1} = \beta_i - 1 / 16$) yields \eqref{ibtr17} with the $i$ there given by $i+1$. This establishes \eqref{ibtr17} and thus the propostion.
	\end{proof}

	\subsection{Likelihood of $\textbf{FHR}$ Restricted to $\textbf{BTR}$}
	
	\label{ProofB}

	In this section we establish \Cref{p:initial}, which is a consequence of the below ``pointwise'' variant of it. Below, we recall that $\chi = 2^{-5000}$.
	
	\begin{lem} 
		\label{l:holder0}
		Adopting \Cref{a:nkrelation}, there exist constants $c = c(A, B) >0$ and $C_1 = C(A, B, D) > 1$ such that the following holds if $L \ge (2B)^{2/\chi}$. For any integer $j \in \llbracket 1, n \rrbracket$ and real numbers $\mathsf{s}, \mathsf{s} + tk^{1/3} \in [-Ak^{1/3}, Ak^{1/3}]$, we have
		\begin{align*}
			\begin{aligned} 
				\mathbb{P} \Bigg[ \textbf{\emph{BTR}}_n (A; B) \cap \bigg\{ \frac{\sfx_j(\mathsf{s}+tk^{1/3})-\sfx_j(\mathsf{s})}{k^{2/3}} \le -L^\chi |t| \Big(\frac{j\vee k}{k} \Big)^{1/3}-4 |t|^{1/2} \Big(\frac{j\vee k}{k} & \Big)^{1/2}  -k^{-D} \bigg\}  \Bigg] \\
				& \le C e^{-c(\log k)^2}.
			\end{aligned}
		\end{align*}
	\end{lem} 
	
	\begin{proof}
		
		  The proof of this lemma will be similar to that of \Cref{probabilityeventr}. In what follows, we will assume that $t \ge 0$, as we may by symmetry under reflection through the line $\{ t = 0 \}$. 
		
		Let $\sfT=2A(j\vee k)^{1/3}$ and $\widetilde{B} = 12 A^2 B^3$, and define the event 
		\begin{flalign*}
			\mathscr{E} = \big\{ \mathsf{x}_j (\mathsf{s}) \le \widetilde{B} k^{2/3} - \widetilde{B}^{-1} j^{2/3} \big\} \cap \big\{ \mathsf{x}_j (\mathsf{T}) \ge -L^{\chi/2} (j\vee k)^{2/3} \big\}, \qquad \text{so that} \qquad \textbf{CTR}_n (A; \widetilde{B}) \subseteq \mathscr{E},
		\end{flalign*} 
		
		\noindent where the last inclusion follows from \Cref{ftrbtr} and \Cref{ftrbtr2} (with the fact that $\mathsf{x}_j (\mathsf{T}) \ge \mathsf{x}_k (\mathsf{T})$ if $j \le k$). Further recall by \Cref{ctr} that there exist constants $c_1 = c_1 (A, B) > 0$ and $C_1 = C_1 (A, B) > 1$ such that 
		\begin{flalign*} 
			\mathbb{P} \big[ \textbf{BTR}_n (A; B) \cap \textbf{CTR}_n (A; \widetilde{B})^{\complement} \big] \le C_1 e^{-c_1 (\log k)^2}, 
		\end{flalign*} 
		
		\noindent Thus, by a union bound, it suffices to show 
		\begin{align}
			\label{xje2}
			\mathbb{P} \Bigg[ \mathscr{E}  \cap \bigg\{ \frac{\sfx_j(\mathsf{s}+tk^{1/3})-\sfx_j(\mathsf{s})}{k^{2/3}} \le -L^\chi t \Big(\frac{j\vee k}{k} \Big)^{1/3}-4 t^{1/2} \Big(\frac{j\vee k}{k} & \Big)^{1/2}  -k^{-D} \bigg\}  \Bigg] \le C e^{-c(\log k)^2}.
		\end{align}
		
		To this end, condition on $\mathcal{F}_{\ext}^{\bm{\mathsf{x}}} \big( \llbracket 1, j \rrbracket \times (\mathsf{s}, \mathsf{T}) \big)$ and restrict to the event $\mathscr{E}$. Let $u = \mathsf{x}_j (\mathsf{s})$ and $v = \mathsf{x}_j (\mathsf{T})$, and set $j_0 = j \vee k$. Sample $j_0$ non-intersecting Brownian bridges $\bm{\mathsf{y}} = (\mathsf{y}_1, \mathsf{y}_2, \ldots , \mathsf{y}_{j_0}) \in \llbracket 1, j_0 \rrbracket \times \mathcal{C} \big( [\mathsf{s}, \mathsf{T}] \big)$ from the measure $\mathsf{Q}^{\bm{u}; \bm{v}}$, where $\bm{u} = (u, u, \ldots, u)$ and $\bm{v} = (v, v, \ldots , v)$ (with each entry appearing with multiplicity $j_0$). Then, $\mathsf{x}_{i+j-j_0} (\mathsf{s}) \ge \mathsf{x}_j (\mathsf{s}) = u = \mathsf{y}_i (\mathsf{s})$ and $\mathsf{x}_{i+j-j_0} (\mathsf{T}) \ge \mathsf{x}_j (\mathsf{T}) = v = \mathsf{y}_i (\mathsf{T})$ for each $i \in \llbracket 1, j \rrbracket$, where we have set $\mathsf{x}_m = \infty$ if $m < 1$. Hence, by \Cref{monotoneheight}, we may couple $\bm{\mathsf{x}}$ and $\bm{\mathsf{y}}$ in such a way that 
		\begin{flalign}
			\label{xjsyjs} 
			\mathsf{x}_j (\mathsf{s} + tk^{1/3}) \ge \mathsf{y}_{j_0} (\mathsf{s} + tk^{1/3}).
		\end{flalign}
		
		Next, by the second part of \Cref{estimatexj} (and using the facts that $\log j_0 \ge \log k$ and that $\big( (tk^{1/3}) (\mathsf{T} - \mathsf{s} - tk^{1/3}) (\mathsf{T} - \mathsf{s})^{-1} \big)^{1/2} \le t^{1/2} k^{1/6}$), there exists constants $c_2 = c_2 (A) > 0$ and $C_2 = C_2 (A, D) > 1$ such that
		\begin{flalign*} 
			\mathbb{P} \bigg[ \mathsf{y}_{j_0} (\mathsf{s} + tk^{1/3}) - u \le \displaystyle\frac{tk^{1/3}}{\mathsf{T} - \mathsf{s}} \cdot (v-u) - (8j_0 t )^{1/2} k^{1/6} - k^{-D} \bigg] \le C_2 e^{-c_2 (\log k)^5},
		\end{flalign*}
		
		\noindent Since $\mathscr{E}$ holds, we have $v-u = \mathsf{x}_j (\mathsf{T}) - \mathsf{x}_j (\mathsf{s}) \ge \widetilde{B}^{-1} j^{2/3} - \widetilde{B} k^{2/3} -L^{\chi/2} j_0^{2/3}  \ge - j_0^{2/3} (L^{\chi/2} + B)$. Thus, since $\mathsf{T} - \mathsf{s} \ge 2Aj_0^{1/3} - Ak^{1/3} \ge Aj_0^{1/3}$ and $u = \mathsf{x}_j (\mathsf{s})$, we have
		\begin{flalign*}
			\mathbb{P} \bigg[ \mathsf{y}_{j_0} (\mathsf{s} + tk^{1/3}) - \mathsf{x}_j (\mathsf{s}) \leq - A^{-1} j_0^{1/3} k^{1/3} (L^{\chi/2} + B) \cdot t - 4 j_0^{1/2} k^{1/6} t^{1/2} - k^{-D} \bigg] \le C_2 e^{-c_2 (\log k)^5}.
		\end{flalign*} 
		
		\noindent Together with \eqref{xjsyjs} and the facts that $j_0 = j \vee k$ and $L^{\chi/2} + B \le L^{\chi} A$ (as $L^{\chi/2} \ge 2B \ge B (A^{-1} + 1)$), this yields \eqref{xje2} and thus the lemma.
	\end{proof}

	\begin{proof}[Proof of \Cref{p:initial}]
		
		The proof of this lemma given \Cref{l:holder0} is similar to that of \Cref{probabilityeventr0} given \Cref{probabilityeventr}. In particular, by \Cref{eventregularityimproved}, it suffices to show that 
		\begin{flalign}
			\label{btrka} 
			\begin{aligned}
				\mathbb{P} \Bigg[ \bigcup_{j=1}^{n'} \bigcup_{\substack{|s| \le A' k^{1/3} \\ |s + tk^{1/3}| \le A' k^{1/3}}} \bigg\{ \displaystyle\frac{\mathsf{x}_j (s+tk^{1/3}) - \mathsf{x}_j (s)}{k^{2/3}} < - L^{\chi} \Big( \displaystyle\frac{j \vee k}{k} & \Big)^{1/3} |t| - 4 \Big( \displaystyle\frac{j \vee k}{k}  \Big)^{1/2} |t|^{1/2} - k^{-D} \bigg\} \\ 
				& \cap \textbf{BTR}_n (A; B) \Bigg] \le C_2 e^{-c_1 (\log k)^2},
			\end{aligned} 
		\end{flalign}
		
		\noindent observing that $t$ can be either positive or negative above (and for any $M \ge 0$ that $\big| \mathsf{x}_j (s+tk^{1/3}) - \mathsf{x}_j (s) \big| \le M$ holds if and only if we have both $\mathsf{x}_j (s + tk^{1/3}) - \mathsf{x}_j (s) \ge -M$ and $\mathsf{x}_j (s) - \mathsf{x}_j (s + tk^{1/3}) \ge -M$). 
		
		Now, denote the $n^{-50(D+1)}$-mesh $\mathcal{T} = [-A' k^{1/3}, A'k^{1/3}] \cap (n^{-50(D+1)} \cdot \mathbb{Z})$, and take a union bound in \Cref{l:holder0} (with the $D$ there given by $2D$ here) over all $i \in \llbracket 1, n \rrbracket$ and $s, t \in \mathcal{T}$; this consists of at most $9(A')^2 k^{2/3} n^{100(D+1)+1} \le 9A^2 n^{300D} \le 9A^2 k^{750D^2}$ (as $n = L^{3/2} k \le k^{3D/2+1} \le k^{5D/2}$) triples $(i, s, t)$. Hence, it yields constants $c_1 = c_1 (A, B) \in (0, 1)$ and $C_1 = C_1 (A, B, D) > 1$ such that $\mathbb{P} \big[ \textbf{BTR}_n (A; B) \cap \mathscr{E}_1^{\complement} \big] \le  C_1 e^{-c_1 (\log k)^2}$, where we have defined the event 
		\begin{flalign*}
			\mathscr{E}_1 = \bigcap_{j=1}^n \bigcap_{s, s+tk^{1/3} \in \mathcal{T}} \bigg\{ \displaystyle\frac{\mathsf{x}_j (s+tk^{1/3}) - \mathsf{x}_j (s)}{k^{2/3}} \ge - L^{\chi} \Big( \displaystyle\frac{j \vee k}{k} & \Big)^{1/3} |t| - 4 \Big( \displaystyle\frac{j \vee k}{k} \Big)^{1/2} |t|^{1/2} - k^{-2D} \bigg\}.
		\end{flalign*}
		
		\noindent Also define the  event 
		\begin{flalign}
			\label{e2e0} 
			\mathscr{E}_2 = \bigcap_{j = 1}^n \bigcap_{|s|, |s'| \le A'k^{1/3}} \Big\{ \big| \mathsf{x}_j (s) - \mathsf{x}_j (s') \big| \le n^5 |s-s'|^{1/3} \Big\}.
		\end{flalign}
		
		We claim that there exist constants $c_2 = c_2 (A, B) \in (0, 1)$ and $C_2 = C_2 (A, B, D) > 1$ such that 
		\begin{flalign}
			\label{e2e1e2} 
			\mathbb{P} \big[ \mathscr{E}_2^{\complement} \cap \textbf{BTR}_n (A; B) \big] \le  C_2 e^{-c_2 (\log k)^2}.
		\end{flalign}

		\noindent Let us establish \eqref{btrka} assuming \eqref{e2e1e2}. First observe by \eqref{e2e1e2}, the estimate $\mathbb{P} \big[ \textbf{BTR}_n (A; B) \cap \mathscr{E}_1^{\complement} \big] \le C_1 e^{-c_1 (\log k)^2}$, and a union bound that we have
		\begin{flalign}
			\label{btrn2}
			\mathbb{P} \big[ \textbf{BTR}_n (A; B) \cap (\mathscr{E}_1^{\complement} \cup \mathscr{E}_2^{\complement}) \big] \le 2C_1 C_2 e^{-c_1 c_2 (\log k)^2}.
		\end{flalign}
		
		\noindent Next, restrict to the event $\mathscr{E}_1 \cap \mathscr{E}_2$; by \eqref{btrn2}, it suffices to show that the event on the left side of \eqref{btrka} does not hold for sufficiently large $k$. To do this, fix an integer $j \in \llbracket 1, n \rrbracket$ and real numbers $s, s + tk^{1/3} \in [-A'k^{1/3}, A'k^{1/3}]$. Set $r = s + tk^{1/3}$, and let $s_0, r_0 \in \mathcal{T}$ be such that $|s-s_0| \le n^{-50(D+1)}$ and $|r - r_0| \le n^{-50(D+1)}$. Then,
		\begin{flalign*}
			\displaystyle\frac{\mathsf{x}_j (s+tk^{1/3}) - \mathsf{x}_j (s)}{k^{2/3}} & \ge \displaystyle\frac{\mathsf{x}_j (r_0) -\mathsf{x}_j (s_0)}{k^{2/3}} - k^{-2/3} \Big( \big| \mathsf{x}_j (r) - \mathsf{x}_j (r_0) \big| - \big| \mathsf{x}_j (s) - \mathsf{x}_j (s_0) \big| \Big) \\
			& \ge -L^{\chi} \Big( \displaystyle\frac{j \vee k}{k} \Big)^{1/3} |t| - 4 \Big( \displaystyle\frac{j \vee k}{k} \Big)^{1/2} |t|^{1/2} - k^{-2D} - 2 n^5 k^{-15D - 15} \\
			& \ge - L^{\chi} \Big( \displaystyle\frac{j \vee k}{k} \Big)^{1/3} |t| - 4 \Big( \displaystyle\frac{j \vee k}{k} \Big)^{1/2} |t|^{1/2} - k^{-D}.
		\end{flalign*}
		
		\noindent where the second bound follows from the facts that we have restricted to $\mathscr{E}_1 \cap \mathscr{E}_2$ and that $|s - s_0|^{1/3} + |r-r_0|^{1/3} \le 2k^{-15D-15}$, and the third holds since $k^{-2D} + 2n^5 k^{-15D-15} \le k^{-2D} + 2 n^{-5D-5} < k^{-D}$ for $k \ge 2$ (as $k^{15D+15} \ge L^{15} k^{15} = n^{10} k^5$ and $n \ge k$). This confirms that the event on the left side of \eqref{btrka} cannot hold on $\mathscr{E}_1 \cap \mathscr{E}_2$, which (as mentioned above) implies the lemma. 
		
		It therefore suffices to verify \eqref{e2e1e2}, which will follow from \Cref{estimatexj3}. In particular, condition on $\mathcal{F}_{\ext}^{\bm{\mathsf{x}}} \big( \llbracket 1, n \rrbracket \times (-Ak^{1/3}, Ak^{1/3}) \big)$ and restrict to $\textbf{BTR}_n (A; B)$. Then, for any $t_0 \in \{ -Ak^{1/3}, Ak^{1/3} \}$ and $s \in [-Ak^{1/3}, Ak^{1/3}]$, we have $\mathsf{x}_{n+1} (s) - \mathsf{x}_1 (t_0) \le 2Bk^{2/3} + B$, due the $\textbf{LOC}$ events (recall \Cref{eventlocation}) in the definition \eqref{e:cEn} of $\textbf{BTR}$. Moreover, since $|Ak^{1/3} - A'k^{1/3}| \ge 1$ (as $A' \in [0, A-k^{1/3}]$), \Cref{estimatexj3} (with the $(a, b)$ there equal to $(-Ak^{1/3}, Ak^{1/3})$ here, and the $ (A, B; \kappa)$ there equal to $(3Bk^{1/2}, n^2; 1 / 4)$ here) applies and yields constants $c_3 = c_3 (A, B) > 0$ and $C_3 = C_3 (A, B) > 1$ such that
		\begin{flalign*} 
			\mathbb{P} \Bigg[ \bigcap_{|s|, |s'| \le A'k^{1/3}} \bigg\{ \big| \mathsf{x}_j (s) - \mathsf{x}_j (s') \big| & \le |s-s'|^{1/2} \Big( n^2 \log \big( 4A k^{1/3} |s-s'|^{-1}  \big)  + 8 A k^{1/3} (n^2 + 3B k^{1/2}) \Big)^2 \\
			& \qquad \qquad \qquad \qquad + 4B n^{2/3} |s-s'| \bigg\} \Bigg] \ge 1 - C_3 e^{C_3 n -c_3 n^4},
		\end{flalign*} 
		
		\noindent where we also used the fact that $\big| \mathsf{x}_j (Ak^{1/3}) - \mathsf{x}_j (-Ak^{1/3}) \big| \le 2B (j^{2/3} + k^{2/3}) \le 4B n^{2/3}$ (again due to the $\textbf{LOC}$ events in $\textbf{BTR}$). This, together with the definition \eqref{e2e0} of $\mathscr{E}_2$ and the fact that for sufficiently large $n \ge k$ we have
		\begin{flalign*} 
			|s-s' &|^{1/2} \Big( n^2 \log \big( 4Ak^{1/3} |s-s'|^{-1} \big) + 8 Ak^{1/3} (n^2 + 3 B k^{1/2}) \Big)^2 + 4B n^{2/3} |s-s'| \\
			& \le |s-s'|^{1/2} \big(  20 Ak^{1/3} n^2 |s-s'|^{-1/12} + 64 A^2 B n^{29/12}  |s-s'|^{-1/12} \big)^2 + 4ABn |s-s'|^{1/3} \\
			& \le |s-s'|^{1/3} ( 400 A^2 n^{14/3} + 4096 A^4 B^2 n^{29/6} + 4ABn) \le n^5 |s-s'|^{1/3}
		\end{flalign*} 
		
		\noindent yields \eqref{e2e1e2} and thus the lemma. Here, in the first inequality, we used the facts that $|s-s'| \le 2An^{1/3} |s-s'|^{1/3}$ (as $|s-s'| \le 2Ak^{1/3} \le 2An^{1/3}$), that $w^{1/12} \log (Mw^{-1}) \le 12 e^{-1} M^{1/12} \le 5 M$ for any real numbers $M \ge 1$ and $w \in (0, M]$ (in particular, at $(w, M) = \big( |s-s'|, 4Ak^{1/3} \big)$), and that $k^{1/3} (n^2 + 3B k^{1/2}) \le 4B n^{7/3} \le 8 AB n^{29/12} |s-s'|^{-1/12}$ (the last since $|s-s'| \le 2Ak^{1/3}$); and in the third, we used the facts that $n \ge k$ and that $(x+y)^2 \le 2(x^2 + y^2)$ for any $x, y \in \mathbb{R}$; and in the fifth we used that $n$ is sufficiently large.
	\end{proof}

	\subsection{Likelihood of $\textbf{IHR}$ Restricted to $\textbf{DEN}$ and $\textbf{BTR}$}
	
	\label{ProofDensity1} 
	
	In this section we establish \Cref{p:dentoReg}, which will again be a consequence of its below pointwise variant.

	\begin{lem} 
		\label{l:holder}

		Adopt the notation and assumptions of \Cref{p:dentoReg}. Setting ${\widetilde \beta}=\beta/2+3/8$, we have for any integer $i \in \llbracket L^{3S\chi/2} k, n \rrbracket$ and real numbers $\mathsf{s}, \mathsf{s} + tk^{1/3} \in [-A'k^{1/3}, A'k^{1/3}]$ that 
		\begin{align}\label{e:pointhold}
			\begin{aligned} 
				\mathbb{P} \Bigg[ \bigg\{ \frac{\sfx_i(\sfs+tk^{1/3})-\sfx_i(\sfs)}{k^{2/3}} \le  -C_1 \bigg( \Big( & \displaystyle\frac{i}{k} \Big)^{1/2} |t| \big| \log(C_1|t|^{-1}) \big|^2 + \Big( \displaystyle\frac{i}{k} \Big)^{2{\widetilde \beta}/3} |t|^{1/2} \bigg) \bigg\} \\
				& \cap \textbf{\emph{BTR}}_n (A; B) \cap \textbf{\emph{DEN}}_{n'} (A'; \beta; R; S) \Bigg] \le  C_2 e^{-c(\log k)^2}.
			\end{aligned}
		\end{align}
	\end{lem} 
	
	\begin{proof}[Proof of \Cref{p:dentoReg} (Outline)]
		
		By \Cref{p:initial}, \eqref{ihrnn}, \eqref{fhr1}, and a union bound, denoting $\widetilde{\beta} = \beta / 2 + 3\beta / 8$, it suffices to show that 
		\begin{flalign*}
			\mathbb{P} \Bigg[  \bigcup_{i=\lceil L^{3S\chi/2} k \rceil}^{n'} \bigcup_{\substack{|\mathsf{s}| \le A' k^{1/3} \\ |\mathsf{s}+tk^{1/3}| \le A'k^{1/3}}} & \bigg\{ \frac{\sfx_i(\sfs+tk^{1/3})-\sfx_i(\sfs)}{k^{2/3}} <  -C_1 \bigg( \Big(  \displaystyle\frac{i}{k} \Big)^{1/2} |t| \big| \log( C_1 |t|^{-1}) \big|^2 \\
			& \qquad \qquad \qquad \qquad \qquad \qquad \qquad + \Big( \displaystyle\frac{i}{k} \Big)^{2{\widetilde \beta}/3} |t|^{1/2} +k^{-D} \bigg) \bigg\} \\
			& \qquad \qquad \cap \textbf{BTR}_n (A; B) \cap \textbf{DEN}_{n'} (A'; \beta; R; S) \Bigg] \le  C_2 e^{-c(\log k)^2}.
		\end{flalign*}
		
		\noindent Given \Cref{l:holder}, the proof of this bound is very similar to that of \Cref{p:initial} given \Cref{l:holder0}, by taking a union bound of \Cref{l:holder} over all $i \in \llbracket L^{3S\chi/2} k, n' \rrbracket$ and $s, t$ in an $n^{-50(D+1)}$ mesh to $[-A'k^{1/3}, A'k^{1/3}]$, and then using the high probability H\"{o}lder regularity of $\bm{\mathsf{x}}$ on $[-A' k^{1/3}, A'k^{1/3}]$ guaranteed by \Cref{estimatexj3} to conclude. We omit further details.		
	\end{proof}

	\begin{figure}
	\center
\includegraphics[scale=.6, trim = 0 1cm 0 1cm]{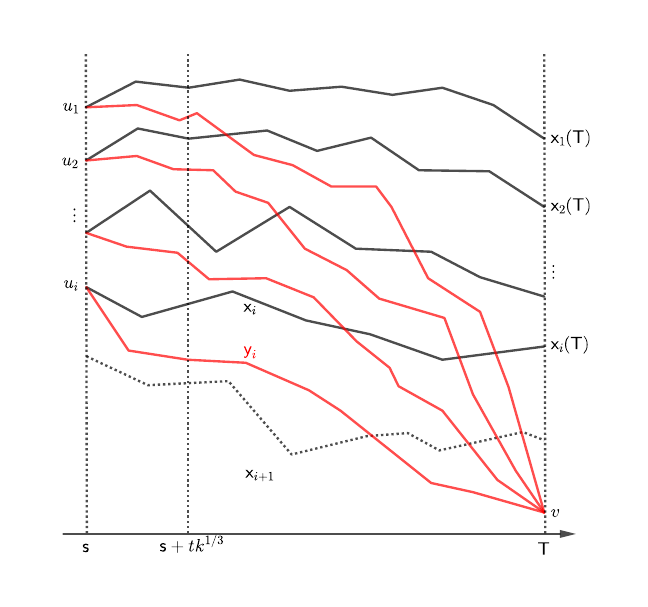}

\caption{Shown above is a depiction of the proof of \Cref{l:holder}.}
\label{f:IHR}
	\end{figure}

	\begin{rem} 
		
		\label{proofd00} 
		
		Before giving the detailed proof of \Cref{l:holder}, let us provide a brief explanation of how it will proceed. By symmetry, we may assume that $t > 0$, so we must bound how low $\mathsf{x}_i$ can reach in a given time interval. Let $\mathsf{T} = 2Ai^{1/3}$; while $\mathsf{T} \notin [-Ak^{1/3}, Ak^{1/3}]$, under the $\textbf{BTR}$ event (recall \Cref{ftrbtr}), the location $\mathsf{x}_i (\mathsf{T}) \ge -L^{\chi/2} i^{2/3}$ is not too low. Denoting $v = -L^{\chi/2} i^{2/3}$, we may therefore by height monotonicity \Cref{monotoneheight} bound $\bm{\mathsf{x}}_{\llbracket 1, i \rrbracket}$ below by a family $\bm{\mathsf{y}} = (\mathsf{y}_1, \mathsf{y}_2, \ldots , \mathsf{y}_i)$ of non-intersecting Brownian bridges on $[\mathsf{s}, \mathsf{T}]$, with the same starting data $\bm{\mathsf{y}} (\mathsf{s}) = \bm{\mathsf{x}}_{\llbracket 1, i \rrbracket} (\mathsf{s})$ as $\bm{\mathsf{x}}_{\llbracket 1, i \rrbracket}$, but with ending data $\bm{\mathsf{y}} (\mathsf{T}) = (v, v, \ldots , v)$. We can estimate $\bm{\mathsf{y}}$ from below using \Cref{uvrho}, where the main condition \eqref{e:gammabb} imposed there is verified at $M = R(k^{-1} i)^{2\beta/3}$ (and $n=i$), due to our restriction to the $\textbf{DEN}$ event (recall \eqref{e:defDen}). The error $M^{1/2} (k^{-1} n)^{1/4}$ in \eqref{e:xiholder} then becomes of order $(k^{-1} i)^{\beta/3 + 1/4}$, which recovers the second error term on the right side of the definition \eqref{fhr1} for the $\textbf{SHR}$ event, at $\beta \mapsto \beta/2 + 3/8$; this is the reason for the choice of $\widetilde{\beta}$ in \Cref{l:holder}. 
		
		Let us mention that our reason for taking $\mathsf{T}$ outside of the interval $[-Ak^{1/3}, Ak^{1/3}]$ (which is the only place in the proof of \Cref{c:finalcouple} where our arguments require doing this) is to control the drifts of the paths in $\bm{\mathsf{y}}$ (this arises in a more precise sense through the lower bound on $v_n - u_n$ imposed through \eqref{vntrho} in \Cref{uvrho}). Indeed, the $\textbf{BTR}$ event allows the $i$-th path in $\bm{\mathsf{x}}$ to travel down by of order $i^{2/3}$ in time $k^{1/3}$, which would yield a drift of around $-k^{-1/3} i^{2/3}$. By increasing the length of the time interval to $\mathsf{T} \sim i^{1/3}$, this reduces the drift to around $-L^{\chi/2} i^{1/3}$, which is now within the first error term on the right side of the definition \eqref{fhr1} for the $\textbf{SHR}$ event (while the original factor of $i^{2/3}$ was not).
		
	\end{rem}

	\begin{proof}[Proof of \Cref{l:holder}]

		The proof of this lemma will follow that of \Cref{l:holder0}, replacing the use of \Cref{estimatexj} (when comparing to a Brownian watermelon) by one of \Cref{uvrho}. In what follows, we will assume that $t \ge 0$, as we may by symmetry under reflection through the line $\{ t = 0 \}$; we also assume (by the scaling invariance \Cref{scale}) that $A = 1$.
		
		Set $L'' = (ik^{-1})^{2/3}$ and $\mathsf{T} = 2Ai^{1/3} = 2i^{1/3}$.  Then define the event
		\begin{flalign*}
			\mathscr{E} = \textbf{DEN}_{i} (\mathsf{s} k^{-1/3}; \beta; R) \cap \big\{ \mathsf{x}_i (\mathsf{s}) \le 0 \big\} \cap \big\{ \mathsf{x}_i (\mathsf{T}) \ge -L^{\chi/2} i^{2/3} \big\}, 
		\end{flalign*} 
	
		\noindent so if $\widetilde{B} = 12A^2 B^3$ (as in \Cref{ctr}), the event $\mathscr{E}$ holds on $\textbf{DEN}_{n'} (A; \beta; R; S) \cap \textbf{CTR}_n (A; \widetilde{B})$, where the last inclusion follows from \Cref{eventrho}, \Cref{ftrbtr2}, and \Cref{ftrbtr} (using the fact that $\widetilde{B} k^{2/3} - \widetilde{B}^{-1} i^{2/3} \le 0$ for $L$ sufficiently large, as $i \ge L^{3S\chi/2} k$). By \Cref{ctr}, it follows that upon restricting to $\textbf{BTR}_n (A; B) \cap \textbf{DEN}_{n'} (A; \beta; R; S)$, the event $\mathscr{E}$ holds off of an event of probability at most $1 - c^{-1} e^{-c(\log n)^2}$. So, conditioning on $\mathcal{F}_{\ext}^{\bm{\mathsf{x}}} \big( \llbracket 1, i \rrbracket \times (\mathsf{s}, \mathsf{T}) \big)$, and restricting to $\mathscr{E}$, to verify \eqref{e:pointhold} it suffices to show that
		\begin{flalign}
			\label{xistk13x}
			\mathbb{P} \Bigg[  	\frac{\sfx_i(\sfs+tk^{1/3})-\sfx_i(\sfs)}{k^{2/3}} \le  -C_1 \bigg( \Big( \displaystyle\frac{i}{k} \Big)^{1/2} \big| \log(2|t|^{-1}) \big|^2  t + \Big( \displaystyle\frac{i}{k} \Big)^{2{\widetilde \beta}/3} t^{1/2} \bigg) \Bigg] \le  C_2 e^{-c(\log k)^2}.
		\end{flalign}

		To this end, let $v = -L^{\chi/2} i^{2/3}$, so that
		\begin{flalign} 
			\label{vxi} 
			v = -L^{\chi/2} i^{2/3} \le \mathsf{x}_i (\mathsf{T}),
		\end{flalign} 
		
		\noindent where the last inequality holds by our restriction to $\mathscr{E}$. Then, since $k^{-1} i \ge L^{3S\chi/2} \ge L^{3\chi}$ (as $S \ge 2$),
		\begin{flalign}
			\label{vxi2} 
			v = -L^{\chi/2} i^{2/3} \ge -(k^{-1} i)^{1/6} i^{2/3}, 
		\end{flalign}
		
		\noindent Also define the $i$-tuples $\bm{u} = \mathsf{x}_{\llbracket 1, i \rrbracket} (\mathsf{s}) \in \mathbb{W}_i$ and $\bm{v} = (v, v, \ldots , v) \in \overline{\mathbb{W}}_i$ (where $v$ appears with multiplicty $i$). Sample $i$ non-intersecting Brownian bridges $\bm{\mathsf{y}} = (\mathsf{y}_1, \mathsf{y}_2, \ldots , \mathsf{y}_i) \in \llbracket 1, i \rrbracket \times \mathcal{C} \big( [\mathsf{s}, \mathsf{T}] \big)$ from the measure $\mathsf{Q}^{\bm{u}; \bm{v}}$; see \Cref{f:IHR}. Since for any $j \in \llbracket 1, i \rrbracket$ we have $\mathsf{x}_j (\mathsf{s}) = u_j$ and $\mathsf{x}_j (\mathsf{T}) \ge \mathsf{x}_i (\mathsf{T}) \ge v = v_i$ by \eqref{vxi}, \Cref{monotoneheight} yields a coupling between $\bm{\mathsf{x}}$ and $\bm{\mathsf{y}}$ such that $\mathsf{x}_i (\mathsf{s} + tk^{1/3}) \ge \mathsf{y}_i (\mathsf{s} + tk^{1/3})$. Since $\mathsf{x}_i (\mathsf{s}) = \mathsf{y}_i (\mathsf{s})$, to prove \eqref{xistk13x} it therefore suffices to show
		\begin{flalign}
			\label{yistk13y}
			\mathbb{P} \Bigg[  	\frac{\mathsf{y}_i(\sfs+tk^{1/3}) - \mathsf{y}_i(\sfs)}{k^{2/3}} \le  -C_1 \bigg( \Big( \displaystyle\frac{i}{k} \Big)^{1/2} \big| \log(2t^{-1}) \big|^2 t + \Big( \displaystyle\frac{i}{k} \Big)^{2{\widetilde \beta}/3} t^{1/2}  \bigg) \Bigg] \le  C_2 e^{-c(\log k)^2}.
		\end{flalign}
		
		This will follow from an application of \Cref{uvrho}, to which end we must verify that $\bm{\mathsf{x}}$ satisfies the hypotheses stated there. This will be a consequence of the fact that we have restricted to the event $\mathscr{E}\subseteq \textbf{DEN}_i (\mathsf{s} k^{-1/3}; \beta; R) $. Indeed, observe from \Cref{eventrho} that there exists a measure $\mu \in \mathscr{P}_{\fin}$ with $\mu(\mathbb{R}) =k^{-1} i$, satisfying the following two properties. First, $\mu$ admits a density $\varrho \in L^1 (\mathbb{R})$ with respect to Lebesgue measure satisfying $\varrho (x) \le R (k^{-1} i)^{1/2}$. Second, denoting its classical locations by $\gamma_j = \gamma_{j; i}^{\mu}$ (recall \Cref{gammarho}) and setting $\mathfrak{m} =  \mathfrak{m}_i = \lceil R \log n \cdot i ^{1/2} \rceil$, we have
		\begin{align}\label{e:dendef}
			\gamma_{i + \mathfrak{m}_i} - R \Big(\frac{i}{k} \Big)^{2\beta/3} \leq k^{-2/3} \cdot \sfx_i(\mathsf{s}) \leq \gamma_{i + \mathfrak{m}_i} +R \Big( \frac{i}{k} \Big)^{2\beta/3}.
		\end{align}
		
		Let us confirm the hypotheses stated in \Cref{uvrho}, with the $(n, K)$ there equal to $\big(i, (i/k)^{2\beta/3} \big)$ here. We have that $\mathsf{T} - \mathsf{s} \in [i^{1/3}, 3i^{1/3}]$ (as $\mathsf{T} = 2i^{1/3}$ and $|\mathsf{s}| \le k^{1/3} \le i^{1/3}$); this verifies the first statement in \eqref{vntrho}. We also have that $\mathsf{y}_n (\mathsf{T}) - \mathsf{y}_n (\mathsf{s}) \ge v \ge - (k^{-1} i)^{1/6} i^{2/3}$ (by \eqref{vxi2} and the fact that $\mathsf{y}_n (\mathsf{s}) \le 0$, as we restricted to $\mathscr{E}$); this verifies the second statement in \eqref{vntrho}. The fact that $\varrho (x) \le R (k^{-1} i)^{1/2}$ verifies the third statement in \eqref{vntrho} (with the $B$ there equal to $R$ here). The definition of $\mathfrak{m}$ verifies the fourth statement of \eqref{vntrho} (with the $B$ there equal to $R$ here). Moreover, the bound \eqref{e:dendef} verifies \eqref{e:gammabb}, with the $M$ there equal to $R(i/k)^{2\beta/3}$ here.
		
		This confirms the hypotheses stated in \Cref{uvrho}, with the $(n, K; \mathsf{T})$ there equal to $\big(i; (i/k)^{2\beta/3}; \mathsf{T} - \mathsf{s} \big)$ here (and the $(A, B)$ there equal to $(1, R)$ here). Hence, \Cref{uvrho} yields constants $c = c(A, B, R) > 0$, $C_3 = C_3 (A, B, R) > 1$, and $C_4 = C_4 (A, B, R, D) > 1$ such that
		\begin{flalign*}
			\mathbb{P} \Bigg[ &  \frac{\sfy_i(\sfs+tk^{1/3})-\sfy_i(\sfs)}{k^{2/3}}  \le -C_3 \bigg( t \Big( \displaystyle\frac{i}{k} \Big)^{1/2} \big| \log (2t^{-1}) \big|^2    + t^{1/2}  \Big( \displaystyle\frac{i}{k} \Big)^{\beta/3 + 1/4}  \bigg) \Bigg] \le C_4 e^{-c (\log k)^2},
		\end{flalign*}
		
		\noindent which implies \eqref{yistk13y} upon recalling the definition of $\widetilde{\beta} = \beta/2+3/8$. 	 
	\end{proof}

	\subsection{Likelihood of $\textbf{DEN}$ Restricted to $\textbf{IHR}$ and $\textbf{BTR}$}
	
	\label{ProofDensity2}
	
	In this section we prove \Cref{p:RegtoDen}, which will be a consequence of its below pointwise variant.

	\begin{lem} 
		
		\label{inrd}
		
		Adopting the notation and assumptions of \Cref{p:RegtoDen}, we have for any integer $i \in \llbracket L^{6S \chi} k, \widetilde{n} \rrbracket$ and real number $t \in \big[ (b - 1) A', (1 - b) A'\big]$ that
		\begin{flalign}
			\label{xbtr1}
			\mathbb{P} \bigg[ \textbf{\emph{BTR}}_n (A; B) \cap \textbf{\emph{DEN}}_i \Big( t; 2 \beta - \displaystyle\frac{7}{8}; C_1 \Big)^{\complement} \bigg] \le C_2 e^{-c (\log k)^2}. 
		\end{flalign}
		
	\end{lem}

	\begin{proof}[Proof of \Cref{p:RegtoDen} (Outline)]
		
		The proof of this proposition given \Cref{inrd} is (as that of \Cref{p:dentoReg} given \Cref{l:holder}) very similar to that of \Cref{p:initial} given \Cref{l:holder0}. In particular, we first take a union bound in \eqref{xbtr1} over all $i \in \llbracket L^{6S \chi} k, \widetilde{n} \rrbracket$ and $t \in \mathcal{T}$, for some $n^{-50}$-mesh $\mathcal{T}$ of $\big[ (b-1) A', (1-b) A' \big]$. For any integer $i \in \llbracket L^{6S\chi} k, n' \rrbracket$ and real number $t \in \mathcal{T}$, this yields a measure $\mu_t^{(i)} \in \mathscr{P}_{\fin}$ satisfying the properties in \Cref{eventrho}, with the $(\beta, R)$ there equal to $( 2\beta - 7 / 8, C_1)$ here. For $t \in \big[ (b-1) A', (1-b) A' \big] \setminus \mathcal{T}$, set $\mu_t^{(i)} = \mu_{t'}^{(i)}$, where $t'$ is an arbitrary element of $\mathcal{T}$ such that $|t-t'| \le n^{-50}$; this $\mu_t^{(i)}$ satisfies the first property in \Cref{eventrho}, since $\mu_{t'}^{(i)}$ does. Using the high probability H\"{o}lder bound for $\bm{\mathsf{x}}$ guaranteed by \Cref{estimatexj3}, it is then quickly verified that $\mu_t^{(i)}$ likely satisfies the second property \eqref{e:defDen} in \Cref{eventrho}, with the $(\beta, R)$ there equal to $( 2\beta - 7 / 8, 2C_1)$ here, for any $t \in \big[ (b-1) A', (1-b) A' \big]$; this confirms that $\textbf{DEN}_{\tilde{n}} \big( (1-b) A'; 2\beta - 7 / 8; 2C_1; 4S \big)$ holds with high probability. We omit further details.	
	\end{proof} 
	
	\begin{figure}
	\center
\includegraphics[scale=.6, trim = 0 1cm 0 1cm]{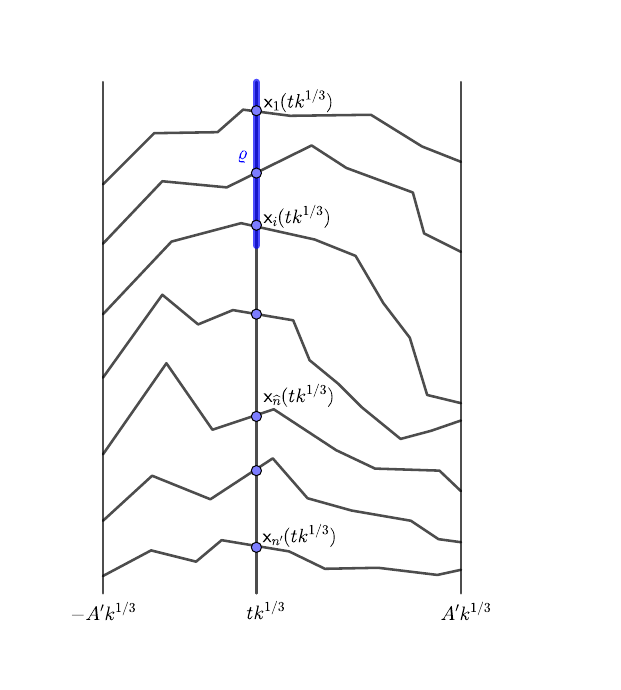}

\caption{Shown above is a depiction of the proof of \Cref{inrd}.}
\label{f:DEN}
	\end{figure}

	\begin{proof}[Proof of \Cref{inrd}] 
		
		 This lemma will follow from an application of \Cref{p:closerho}. Let $\zeta \in [2^{-250}, 1]$ denote the constant defined there (see also \Cref{xycouple2}), and define $\widehat{L} = (k^{-1} i)^{2/3 \zeta}$ and $\widehat{n} = \lceil \widehat{L}^{3/2} k \rceil$; we will omit the ceilings in what follows, assuming that $\widetilde{n} = (L')^{3\zeta/2} k$ and $\widehat{n}= \widehat{L}^{3/2} k$, as this will barely affect the proofs. Observe that $\widehat{L} = (k^{-1} i)^{2/3\zeta} \le (k^{-1} \widetilde{n})^{2/3\zeta} = L'$ (as $i \le \widetilde{n})$; in particular, $\widehat{n} = \widehat{L}^{3/2} k \le (L')^{3/2} k = n'$. Throughout, we abbreviate the event $\textbf{ICE}_{\widehat{n}} = \textbf{ICE}_{\widehat{n}}^{\bm{\mathsf{x}}} (A', 12A^2 B^3; \beta; R, S)$. 
		
		We will apply \Cref{p:closerho} with the $(n', L', k)$ there equal to $(\widehat{n}, \widehat{L}, k)$ here, to which end we must verify the assumptions imposed there. First observe since $i \ge L^{6 S \chi} k$ that $\widehat{L} = (k^{-1} i)^{2/3\zeta} \ge L^{4S\chi / \zeta}$, confirming \eqref{ldeltal}. To show \eqref{probabilityeventc0}, we apply \Cref{gevent} (with the $(n', n'', n''')$ there given by $(n', \widehat{n}, \widehat{n})$ here), whose hypothesis \eqref{btrr2} is verified by \eqref{btrr}. That lemma yields constants $c_1 = c_1 (A, B, \xi) > 0$ and $C_3 = C_3 (A, B, D, \Xi) > 1$, and an event $\mathscr{G} \subseteq \textbf{BTR}_n (A; B)$ (obtained by intersecting the $\mathscr{G}_0$ of \eqref{g0event} with $\textbf{BTR}_n (A; B)$) measurable with respect to $\mathcal{F}_{\ext} = \mathcal{F}_{\ext}^{\bm{\mathsf{x}}} \big( \llbracket 1, \widehat{n} \rrbracket \times (-A' k^{1/3}, A' k^{1/3}) \big)$ satisfying the following two properties. First, we have
		\begin{flalign}
			\label{btrg}
			\mathbb{P} \big[ \textbf{BTR}_n (A; B) \setminus \mathscr{G} \big] \le C_3 e^{-c_1 (\log k)^2}.
		\end{flalign} 
		
		\noindent Second, conditioning on $\mathcal{F}_{\ext}$ and restricting to $\mathscr{G}$, we have $\mathbb{P} [ \textbf{ICE}_{\widehat{n}} ] \ge 1 - C_3 e^{-c_1 (\log k)^2}$. The latter verifies \eqref{probabilityeventc0}. 
		
		Thus, conditioning on $\mathcal{F}_{\ext}$ and restricting to $\mathscr{G}$, \Cref{p:closerho} applies and (since $t \in \big[ (b-1) A', (1-b) A' \big]$) yields constants $c_2 = c_2 (b, A', B, \xi) > 0$ and $C_1 = C_1 (b, A', B, R) > 1$, and $C_4 = C_4 (b, A', B, D, R, \xi, \Xi) > 1$, and a measure $\widehat{\mu} \in \mathscr{P}_{\fin}$ with $\widehat{\mu} (\mathbb{R}) = \widehat{L}^{3/2}$, satisfying the following two properties. First, $\widehat{\mu}$ admits a density $\widehat{\varrho} \in L^1 (\mathbb{R})$ with respect to Lebesgue measure satisfying $\supp \widehat{\varrho} \subseteq [-C_1 L, C_1 L^{3/4}]$ and $\widehat{\varrho} (x) \le C_1 \max \{ 1, -x \}^{3/4}$. Second, denoting the classical locations of $\widehat{\mu}$ by $\gamma_j = \gamma_{j; \widehat{n}}^{\widehat{\mu}}$, we have $-C_1 (j/k)^{2/3} - C_1 \le \gamma_j \le C_1 - C_1^{-1} (j/k)^{2/3}$ for each $j \in \big\llbracket (\log \widehat{n})^6 + 1, \widehat{n} \big\rrbracket$. Third, denoting $\mathfrak{m}_j = \big\lceil C_1 \log n \cdot \max \{ j^{1/2}, k^{1/2} \} \big\rceil$ for each integer $j \in \llbracket 1, \widehat{n} \rrbracket$, we have 
		\begin{flalign}
			\label{k23zetagammax} 
			\begin{aligned} 
				\mathbb{P} \Bigg[ \bigcap_{j=1}^{\lceil \widehat{L}^{3\zeta/2} k \rceil} \big\{  \gamma_{j+\mathfrak{m}_j} - C_1 \widehat{L}^{\zeta (2\beta - 7/8)} \le k^{-2/3} \cdot \mathsf{x}_j (tk^{1/3}) & \le \gamma_{j - \mathfrak{m}_j} + C_1 \widehat{L}^{\zeta (2\beta-7/8)} \big\} \Bigg] \\
				& \qquad \qquad \quad \ge 1 - C_4 e^{-c_2 (\log k)^2}.
			\end{aligned} 
		\end{flalign}
		
		Now, define $\gamma \in \mathbb{R}$ to be the minimal real number such that $\widehat{\mu} \big( [\gamma, \infty) \big) = k^{-1} i$. The existence of at least one $\gamma$ follows from the fact that $\widehat{\varrho}$ is bounded and that $\widehat{\mu} (\mathbb{R}) = \widehat{L}^{3/2} = (k^{-1} i)^{1/\zeta} > k^{-1} i$; it moreover satisfies $-\gamma \le -\gamma_i \le C_1 (i/k)^{2/3} + C_1 \le 2C_1 (i/k)^{2/3}$. Define $\varrho \in L^1 (\mathbb{R})$ by setting $\varrho (x) = \widehat{\varrho} (x) \cdot \textbf{1}_{x \ge \gamma}$, and define the measure $\mu = \varrho (x) dx \in \mathscr{P}_{\fin}$. Then, $\mu(\mathbb{R}) = k^{-1} i$, and the classical locations of $\mu$ are given by $\gamma_j = \gamma_{j; i}^{\mu} = \gamma_{j; \widehat{n}}^{\widehat{\mu}}$, for each $j \in \llbracket 1, i \rrbracket$; moreover, $\supp \mu \subseteq \supp \widehat{\varrho} \subseteq [-C_1 L, C_1 L^{3/4}]$. See \Cref{f:DEN}. Hence, the fact that $\varrho (x) \le \widehat{\varrho} (x) \le C_1 \max \{ 1, -x \}^{3/4} \le C_1 \max \{ \gamma, 1 \}^{3/4} \le 2C_1^2 (i/k)^{1/2}$; \eqref{k23zetagammax}; and the fact that $\widehat{L}^{\zeta} = (k^{-1} i)^{2/3}$ (so $\widehat{L}^{3\zeta/2} k = i$), together imply by \Cref{eventrho} that 
		\begin{flalign*} 
			\mathbb{P} \bigg[ \textbf{DEN}_i^{\bm{\mathsf{x}}} \Big( t; 2\beta - \displaystyle\frac{7}{8}; 2C_1^2 \Big) \bigg] \ge 1 - C_4 e^{-c_2 (\log k)^2},
		\end{flalign*} 
		
		\noindent after conditioning on $\mathcal{F}_{\ext}$ and restricting to the event $\mathscr{G} \subseteq \textbf{BTR}_n (A; B)$. This, together with \eqref{btrg} and a union bound, establishes the lemma.
	\end{proof}

		\chapter{Global Law and Regularity}
		
		\label{GlobalRegular} 
		
		Throughout this chapter, we let $\bm{\mathsf{x}} = (\mathsf{x}_1, \mathsf{x}_2, \ldots )$ denote a $\mathbb{Z}_{\ge 1} \times \mathbb{R}$ indexed line ensemble satisfying the Brownian Gibbs property, and we recall the $\sigma$-algebra $\mathcal{F}_{\ext}$ from \Cref{property}. We also recall the events $\textbf{TOP}$, $\textbf{GAP}$, and $\textbf{BTR}$ from \Cref{eventsregular1}, \Cref{gap}, and \Cref{ftrbtr}, respectively.

	\section{Likelihood of Regular Profile Events} 
	
	\label{ProofRegular}

	In this section we prove \Cref{p:closerho0}, which indicates that regular profile events are likely, upon restricting to the intersection of several $\textbf{TOP}$ events (from \Cref{eventsregular1}). Recall from \Cref{c:finalcouple} that the boundary removal coupling held under restricting to the boundary tall rectangle event $\textbf{BTR}$ (from \Cref{ftrbtr}) is likely; we verify that this holds upon restricting to several $\textbf{TOP}$ events in \Cref{EventProofB}. We then establish \Cref{p:closerho0} in \Cref{ProofRegular0} and \Cref{Proofx2}; it will eventually amount to being a consequence of the boundary removal coupling, together with \Cref{p:densityub0}.

		\subsection{Likelihood of $\textbf{BTR}$ Restricted to $\textbf{TOP}$ Events}

	\label{EventProofB}

	In this section we lower bound the $\textbf{BTR}$ event, through the following proposition. In what follows, we recall the events $\textbf{GAP}$ from \Cref{gap}, and $\textbf{BTR}$ from \Cref{ftrbtr}.

	\begin{prop} 
		
		\label{ndeltal} 
		
		Adopt \Cref{l0}. For any real numbers $A \ge 3$ and $\varepsilon \in ( 0,  1 / 2)$, there exist constants $B = B(A) > 1$, $R = R(A) > 1$, and $C = C(A, \varepsilon) > 1$, such that the following holds for any integer $k \ge C$.  Setting $L = k^{2^{5000}}$ and $n = L^{3/2} k$, we have
		\begin{flalign*}
			\mathbb{P} \Big[ \textbf{\emph{BTR}}_n^{\bm{\mathcal{L}}} (A, B; k, L) \cap \textbf{\emph{GAP}}_n^{\bm{\mathcal{L}}} \big( [-Ak^{1/3}, Ak^{1/3}]; R \big) \Big] \ge 1 - \varepsilon.
		\end{flalign*}
	\end{prop} 

	 We will quickly deduce \Cref{ndeltal} as a consequence of the following proposition (together with \Cref{x1lsmall}), which states that the boundary tall rectangle event $\textbf{BTR}$ and gap event $\textbf{GAP}$ are likely upon restricting to the intersection of several $\textbf{TOP}$ events; it applies to any $\mathbb{Z}_{\ge 1} \times \mathbb{R}$ indexed line ensemble $\bm{\mathsf{x}}$ satisfying the Brownian Gibbs property (as fixed at the beginning of this chapter). In what follows, we recall the $\textbf{TOP}$ event from \Cref{eventsregular1}.
	
		\begin{prop}\label{p:cEpbound}
			
			For any real number $A \ge 3$, there exist constants $c = c(A) > 0$, $\vartheta = \vartheta (A) > 0$, $B = B (A) > 1$, and $R = R(A) > 1$, such that the following holds. Set $D = 2^{5000}$ and $\omega = 2^{-6000}$; let $k \ge 1$ be an integer; and set $L = k^D$ and $n = L^{3/2} k$. For each integer $j \ge 0$, set $K_j = \lceil L^{j\omega} k \rceil$. Then,  
			\begin{align}\begin{split}\label{e:large}
					\bP\Bigg[
				\bigcap_{j=0}^{ 3/\omega }\textbf{\emph{TOP}}^{\bm{\mathsf{x}}} & \big( [-\vartheta^{-1} K_j^{1/3}, \vartheta^{-1} K_j^{1/3}]; \vartheta K_j^{2/3} \big)\cap \textbf{\emph{TOP}}^{\bm{\mathsf{x}}} \big( [-\vartheta^{-1} K_j^{30D}, \vartheta^{-1} K_j^{30D}]; \vartheta K_j^{60D}\big) \\
					& \cap \left( \textbf{\emph{BTR}}_n^{\bm{\mathsf{x}}} (A, B;  k, L)\cap \textbf{\emph{GAP}}_n^{\bm{\mathsf{x}}} \big( [-Ak^{1/3}, Ak^{1/3}]; R \big) \right)^{\complement}  \Bigg]	
					\le c^{-1} e^{-c (\log k)^2}. 
			\end{split}\end{align}

		\end{prop}

	\begin{proof}[Proof of \Cref{ndeltal}]

		Fix $D = 2^{5000}$ and $\omega = 2^{-6000}$, denote $K_j = \lceil L^{j\omega} k \rceil$ for each integer $j \ge 0$, and let $\vartheta = \vartheta (A)$ be as in \Cref{p:cEpbound}. By \Cref{x1lsmall}, we have for any sufficiently large real number $m > 1$ that $\mathbb{P} \big[ \textbf{{TOP}}^{\bm{\mathcal{L}}} \big( [-\vartheta^{-1} m^{1/3}, \vartheta^{-1} m^{1/3}]; \vartheta m^{2/3} \big) \big] \ge 1 - \omega \varepsilon / 24$. Taking  a union bound over $m \in \{ K_0, K_1, \ldots , K_{ 3 / \omega } \} \cup \{ K_0^{90D}, K_1^{90D}, \ldots , K_{ 3 / \omega }^{90D} \}$ (observing that the size of this set is at most $12 \omega^{-1}$) and taking $k$ to be sufficiently large, we deduce that 
		\begin{flalign*} 
			\mathbb{P} \Bigg[ \bigcap_{j = 0}^{ 3 / \omega }  \textbf{TOP}^{\bm{\mathcal{L}}} \big( [-\vartheta^{-1} K_j^{1/3}, \vartheta^{-1} K_j^{1/3} ]; \vartheta K_j^{2/3} \big) \cap \textbf{TOP}^{\bm{\mathcal{L}}} \big( [-\vartheta^{-1} K_j^{30D}, \vartheta^{-1} K_j^{30D} & ]; \vartheta K_j^{60D} \big)  \Bigg] \\ 
			& \ge 1 - \displaystyle\frac{\varepsilon}{2}.
		\end{flalign*} 
	
		\noindent This, together with \Cref{p:cEpbound} and a union bound, implies the proposition. 		
	\end{proof} 

	To prove \Cref{p:cEpbound}, we use the following lemma, stating that the intersection of certain $\textbf{GAP}$ and $\textbf{IMP}$ events implies the $\textbf{BTR}$ event; its proof follows quickly from the definitions of these events. Below, we recall the improved medium position event $\textbf{{IMP}}$ from \Cref{eventimproved}; observe that they can be expressed through the $\textbf{LOC}$ events of \Cref{eventlocation} by 
	\begin{flalign}
		\label{locationimproved}
		\textbf{IMP}_n (A; B; C; R) = \bigcap_{j = \lceil n/B \rceil}^{\lfloor Rn \rfloor} \textbf{LOC}_j \big( [-An^{1/3}, An^{1/3}]; C^{-1} n^{2/3} - Cj^{2/3}; Cn^{2/3} - C^{-1} j^{2/3} \big).
	\end{flalign}

		\begin{lem}\label{l:cEandIMP}
			
			Set $D = 2^{5000}$ and $\chi = 2^{-5000}$; let $k \ge 2^{500}$ be an integer; set $L = k^D$; and denote $n = L^{3/2}k $. Fix real numbers $A,B, C, R \geq 2$, and assume that $L \ge (2C)^{4/\chi}$. Defining $n_j= \lceil L^{3j\chi/8}k \rceil$ for each integer $j\geq 0$, we have
			\begin{align*} 
				\textbf{\emph{GAP}}_n^{\bm{\mathsf{x}}} \big( [-Ak^{1/3}, Ak^{1/3}]; R \big) \cap \bigcap_{j=0}^{4/\chi} \textbf{\emph{IMP}}_{n_j}^{\bm{\mathsf{x}}} (2A; B; C; n_j^{3D}) \subseteq \textbf{\emph{BTR}}_n^{\bm{\mathsf{x}}} (A, C+R+2; k, L).
			\end{align*}
		\end{lem}

		\begin{proof}
			
			Set $B_0 = C+R+2$. By \Cref{ftrbtr} (and the fact that $\llbracket k, n \rrbracket \subseteq \bigcup_{j=1}^{\lceil 4 / \chi\rceil} \llbracket n_{j-1}, n_j \rrbracket$), to establish the lemma, it suffices to verify the following two statements. First (recalling that $k = n_0$), for any integer $j \in \llbracket 1, n+1 \rrbracket$, we have 
			\begin{align}\label{e:IMPbelong}
				\begin{aligned}
					\textbf{{IMP}}_{k} (2A; B; C; k^{3D} & ) \cap \textbf{{GAP}}_n \big( [-Ak^{1/3}, Ak^{1/3}]; R \big)\\
					&\subseteq\textbf{{LOC}}_j \big( [-Ak^{1/3}, Ak^{1/3}] ; -B_0j^{2/3} - B_0k^{2/3}; -B_0^{-1} j^{2/3}+B_0 k^{2/3}   \big).
				\end{aligned} 
			\end{align}
			
			\noindent Second, for any integers $i \in \big\llbracket 1, \lceil 4 \chi^{-1} \rceil \big\rrbracket$ and $j \in \llbracket n_{i-1}, n_i \rrbracket$, we have
			\begin{align}\label{e:IMPni}
				\textbf{{IMP}}_{n_i} (2A; B; C; n_i^{3D})   
				\subseteq \big\{ \sfx_j(-2A j^{1/3})\geq -L^{\chi/2} j^{2/3} \big\} \cap \big\{ \mathsf{x}_j (2Aj^{1/3}) \ge -L^{\chi/2} j^{2/3} \big\}. 
			\end{align}
			
			For $j \in \llbracket k, n + 1 \rrbracket$, the inclusion \eqref{e:IMPbelong} follows from \eqref{locationimproved}, the facts that $\llbracket k, n+1 \rrbracket \subseteq \llbracket B^{-1} k, k^{3D+1} \rrbracket$, that $B_0 \ge C$, and that $\textbf{LOC}_j (\mathcal{T}; b; B) \subseteq \textbf{LOC}_j (\mathcal{T}'; b'; B')$ whenever $\mathcal{T}' \subseteq \mathcal{T}$ and $b' \le b \le B \le B'$. To confirm this inclusion when $j \in \llbracket 1, k-1 \rrbracket$, restrict to the event $\textbf{{IMP}}_{k} (2A; B; C; k^{3D}) \cap \textbf{{GAP}}_n \big( [-Ak^{1/3}, Ak^{1/3}]; R \big)$. We must verify that $-B_0 j^{2/3} - B_0 k^{2/3} \le \mathsf{x}_j (t) \le B_0 k^{2/3} - B_0^{-1} j^{2/3}$ for any $t \in [-Ak^{1/3}, Ak^{1/3}]$. To establish the lower bound on $\mathsf{x}_j$, observe for any such $t$ that 
			\begin{flalign}
				\label{xjtb0} 
				\sfx_j(t)
				&\geq \sfx_{k}(t)
				\geq  -C k^{2/3}\geq -B_0 j^{2/3}- B_0 k^{2/3},
			\end{flalign}
			
			\noindent where in the first bound we used the fact that $\mathsf{x}_j \ge \mathsf{x}_{j'}$ for $j \le j'$; in the second we used the fact that we restricted to the $\textbf{LOC}_k$ event contained in the $\textbf{IMP}_k$ one by \eqref{locationimproved}; and in the third we used the fact that $B_0 \ge C$. To establish the upper bound on $\mathsf{x}_j$, observe that 
			\begin{flalign*} 
				\mathsf{x}_j (t) &\leq \big| \sfx_{j}(t)-\sfx_k (t) \big| + \mathsf{x}_k (t) 
				\le Rk^{2/3} + (\log k)^{25} -C^{-1} k^{2/3}  \leq (R+1) k^{2/3} \leq  B_0 k^{2/3} - B_0^{-1}j^{2/3},
			\end{flalign*}
			
			\noindent where in the second bound we used our restriction to $\textbf{GAP}$ event (recall \Cref{gap}) and the $\textbf{LOC}_k$ event contained in the $\textbf{IMP}_k$ one (by \eqref{locationimproved}); in the third we used the fact that $(\log k)^{25} \le k^{2/3}$ for $k \ge 2^{500}$; and in the fourth we used the facts that $B_0 \ge R+2$ and $j \le k$. By this and \eqref{xjtb0}, the event $\textbf{LOC}_j \big( [-Ak^{1/3}, Ak^{1/3}]; -B_0 j^{2/3} - B_0 k^{2/3}; B_0 k^{2/3} - B_0^{-1} j^{2/3} \big)$ holds, verifying \eqref{e:IMPbelong}. 
			
			The inclusion \eqref{e:IMPni} follows from the fact that, fixing $j \in \llbracket n_{i-1}, n_i \rrbracket$ and restricting to the event $\textbf{{IMP}}_{n_i} (2A; B; C; n_i^{3D})$, we have for any $t \in \{ -2Aj^{1/3}, 2Aj^{1/3} \}$ that
			\begin{align*}
				\sfx_j(t)\geq \sfx_{n_i} (t)\geq  -Cn_i^{2/3} \ge -2 CL^{\chi/4} n_{i-1}^{2/3}\geq -L^{\chi / 2}j^{2/3}.
			\end{align*}
			
			\noindent Here, in the first bound we used the fact that $\mathsf{x}_j \ge \mathsf{x}_{j'}$ for $j \le j'$; in the second we used our restriction to the $\textbf{LOC}_{n_i}$ event contained in the $\textbf{IMP}_{n_i}$ one, by \eqref{locationimproved}; in the third we used the fact that $n_i = \lceil L^{3i\chi/8} k \rceil \le 2 L^{3\chi/8} \cdot \lceil L^{3(i-1)\chi/8} k \rceil = 2 L^{3\chi/8} n_{i-1}$; and in the fourth we used the facts that $j \ge n_{i-1}$ and $L\geq (2C)^{4/\chi}$.
		\end{proof}

		\begin{proof}[Proof of \Cref{p:cEpbound}]	
			
			Throughout, we recall the medium position event $\textbf{MED}$ from \Cref{eventsregular1} and the on-scale event $\textbf{SCL}$ from \Cref{eventscl}. Set $\mathfrak{C}_1$ to be the constant $C_1 (450)$ from \Cref{probabilityimproved}; set $\mathfrak{C}_2$ to be $C_2 (\mathfrak{C}_1 + A, 1800)$, where $C_2$ denotes the constant from \Cref{sclprobability}; and set $\vartheta$ to be such that $\vartheta^{-1}$ is $\mathfrak{C}_1 + \mathfrak{C}_2 + 2A + C_1 (1800) + 200$, where $C_1$ denotes the constant from \Cref{sclprobability}. Then, define the event 
			\begin{flalign}
				\label{te} 
				\mathscr{E} = \bigcap_{j =0}^{ 3 /\omega } \textbf{TOP} \big( [-\vartheta^{-1} K_j^{1/3}, \vartheta^{-1} K_j^{1/3}]; \vartheta K_j^{2/3} \big) \cap \textbf{TOP} \big( [-\vartheta^{-1} K_j^{30D}, \vartheta^{-1} K_j^{30D}]; \vartheta K_j^{60D} \big).
			\end{flalign}
			
			For the remainder of this proof, we restrict to $\mathscr{E}$. By \Cref{sclprobability}, there exist constants $c_1 = c_1 (A)$ and $R = R(A)$ such that the event 
			\begin{flalign}
				\label{e2} 
				\mathscr{E}' = \bigcap_{j=0}^{3/\omega} \textbf{SCL}_{K_j} (\mathfrak{C}_1 + A; 1800; R) \cap \textbf{SCL}_{K_j^{90D}} (\mathfrak{C}_1 + A; 1800; R),
			\end{flalign} 
		
			\noindent  holds outside an event of probability at most $c_1^{-1} e^{-c_1 (\log k)^2}$. So, we further restrict to $\mathscr{E}'$ below. Since for $K_{3/2\omega} = n $, we have (from \Cref{eventscl}) that $\textbf{GAP}_n \big( [-n^{1/3} (\mathfrak{C}_1 + A), n^{1/3} (\mathfrak{C}_1 + A)]; R \big) \subseteq \textbf{GAP}_n \big( [-An^{1/3}, An^{1/3}];  R \big)$ holds on $\mathscr{E}'$. Thus, it remains to show for some constant $B = B(A) > 1$ that $\textbf{BTR}_n^{\bm{\mathsf{x}}} (A, B; k, L)$ holds with probability at least $1 - c^{-1} e^{-c(\log k)^2}$. To that end, setting $\chi = 2^{-5000}$ and $n_j = \lceil L^{3 j \chi / 8} \rceil$ for each integer $j \ge 0$, it suffices by \Cref{l:cEandIMP} to show for some sufficiently large constant $M = M(A) > 1$ that $\bigcap_{j=0}^{4/\chi} \textbf{IMP}_{n_j}^{\bm{\mathsf{x}}} (2A; 450; M; n_j^{3D})$ holds with probability at least $1 - c^{-1} e^{-c(\log k)^2}$. 
			
			This will follow from \Cref{probabilityimproved}. Indeed, fixing $j \in \llbracket 0, 4/\chi \rrbracket$ and denoting $j' = 2^{1000} \cdot 3 j / 8$, we have $n_j = K_{j'}$ (and $j' \in \llbracket 0, 3 / \omega \rrbracket$). Therefore, the event $\textbf{TOP} \big( [-\vartheta^{-1} n^{1/3}, \vartheta^{-1} n^{1/3}]; \vartheta n_j^{2/3} \big) \subseteq \textbf{TOP} \big( [-2An_j^{1/3}, 2An_j^{1/3}]; n_j^{2/3}/200 \big)$ holds on $\mathscr{E}$, as does $\textbf{TOP} \big( [-\vartheta^{-1} n_j^{30D}, \vartheta^{-1} n_j^{30D}]; \vartheta n_j^{60D} \big) \subseteq \textbf{TOP} \big( [-\mathfrak{C}_1 n_j^{30D}, \mathfrak{C}_1 n_j^{30D}]; 450 n_j^{60D} \big)$ (due to the $j'$-th $\textbf{TOP}$ events on the right side of \eqref{te}). Similarly, on $\mathscr{E}'$, the event $\textbf{MED}_{\lfloor n_j / 1800 \rfloor} \big( [ -  2An_j^{1/3}, 2A n^{1/3}]; n_j^{2/3} /100; 450 n_j^{2/3} \big)$ holds, as does $\textbf{MED}_{n_j^{90D}} \big( [-\mathfrak{C}_1 n_j^{30D}, \mathfrak{C}_1 n_j^{30D}]; -450 n_j^{60D}; 450 n_j^{60D} \big)$ (due to the $j'$-th $\textbf{SCL}$ events on the right side of \eqref{e2}). Therefore, \Cref{probabilityimproved} (with the $(A, B, b, D)$ there equal to $(2A, 450, 1/200, 3D)$ here) yields constants $c = c(A) > 1$ and $M = M(A) > 1$ such that $\textbf{IMP}_{n_j} (2A; 450; M; n_j^{3D})$ holds with probability at least $1 - c_2^{-1} e^{-c_2 (\log k)^2}$, upon restricting to $\mathscr{E} \cap \mathscr{E}'$. Taking a union bound over $j \in \llbracket 0, 4/\chi \rrbracket$, we (as mentioned previously) deduce the proposition. 	
		\end{proof}

		\subsection{Likelihood of Regular Profile Events}
		
		\label{ProofRegular0}
			
		In this section we establish \Cref{p:closerho0}, which will be a consequence of the below proposition (together with \Cref{ndeltal} and \Cref{p:cEpbound}). The latter is a general result stating that, if $\bm{\mathsf{x}}$ is a $\mathbb{Z}_{\ge 1} \times \mathbb{R}$ indexed line ensemble satisfying the Brownian Gibbs property (as fixed at the beginning of this chapter) for which both a boundary tall rectangle event $\textbf{BTR}$ and gap event $\textbf{GAP}$ are likely, then $\bm{\mathsf{x}}$ also satisfies a regular profile event $\textbf{PFL}$ with high probability. In what follows, we recall the events $\textbf{GAP}$, $\textbf{PFL}$, and $\textbf{BTR}$ from \Cref{gap}, \Cref{eventpfl}, and \Cref{ftrbtr}, respectively.
		
		\begin{prop}
		
			\label{x2p} 
			
			For any real numbers $A, B, R \ge 4$, there exists a constant $C = C (A, B, R) > 1$ such that the following holds. Let $k \ge 1$ be an integer; setting $L = k^{2^{5000}}$, denote $n = L^{3/2} k$. Define $\bm{x} = (x_1, x_2, \ldots , x_k) \in \llbracket 1, k \rrbracket \times \mathcal{C} \big( [-A / 2, A / 2] \big)$ from $\bm{\mathsf{x}} \in \mathbb{Z}_{\ge 1} \times \mathcal{C}(\mathbb{R})$ by setting $x_j (s) = k^{-2/3} \cdot \mathsf{x}_{j+k} (sk^{1/3})$ for each $(j, s) \in \llbracket 1, k \rrbracket \times \big[ -A / 2, A / 2 \big]$. Then,
			\begin{flalign*}
				\mathbb{P} \Bigg[ \bigcup_{|t| \le A / 4} \textbf{\emph{PFL}}^{\bm{x}} \big( t; k^{-1} (\log k)^7; C \big)^{\complement} & \cap \textbf{\emph{BTR}}_n^{\bm{\mathsf{x}}} (A, B; k, L) \\
				&  \cap \textbf{\emph{GAP}}_n^{\bm{\mathsf{x}}} \big( [-Ak^{1/3}, Ak^{1/3}]; R \big)  \Bigg] \le C k^{-100}. 
			\end{flalign*} 
		\end{prop}

		\begin{proof}[Proof of \Cref{p:closerho0}]
			
			Throughout this proof, we recall the $\textbf{GAP}$ and $\textbf{BTR}$ events from \Cref{gap} and \Cref{ftrbtr}, respectively. We further fix 
			\begin{flalign*} 
				D = 2^{5000}; \qquad \chi = 2^{-5000}; \qquad \omega = 2^{-6000}; \qquad L = n^D; \qquad N = L^{3/2} n.
			\end{flalign*} 
		
			\noindent Additionally, let $\vartheta$ denote the constant $\vartheta (4A, 2^{5000})$ from \Cref{p:cEpbound}; set $\vartheta_0 = \vartheta/2$; and restrict to the event 
			\begin{flalign}
				\label{eomega2}
				\mathscr{E} = \bigcap_{j=1}^{1/\omega^2} \textbf{TOP}^{\bm{\mathcal{L}}} \big( [-\vartheta_0^{-1} n^{j\omega/3}, \vartheta_0^{-1} n^{j \omega/3}]; \vartheta_0 n^{2j\omega/3} \big).
			\end{flalign}
		
			\noindent For each integer $j \ge 0$, also set $K_j = \lceil L^{j \omega} n \rceil$, and denote $m_j = jD + \omega^{-1}$ and $m_j' = 90 m_j$. Observe for all $j \in \llbracket 0, 3/\omega \rrbracket$ that $\omega^{-2} \ge m_j' \ge m_j$ (as $\omega = 2^{-6000}$ and $D = 2^{5000}$). Therefore, on $\mathscr{E}$, the event 
			\begin{flalign}
				\label{3omega2} 
				\bigcap_{j=0}^{3/\omega} \textbf{TOP}^{\bm{\mathcal{L}}} \big( [-\vartheta^{-1} K_j^{1/3}, \vartheta^{-1} K_j^{1/3}]; \vartheta K_j^{2/3} \big) \cap \textbf{TOP}^{\bm{\mathcal{L}}} \big( [-\vartheta^{-1} K_j^{30D}, \vartheta^{-1} K_j^{30D}]; K_j^{60D} \big),
			\end{flalign} 
		
			\noindent holds, since the $j$-th event in the intersection in \eqref{3omega2} contains the intersection between the $m_j$-th and $m_j'$-th events in \eqref{eomega2} (as $L = n^D$).
			
			Thus, by \Cref{p:cEpbound} (with the $(n, k, A)$ there equal to $(N, n, 4A)$ here), there exist constants $B = B(A) > 1$, $R = R(A) > 1$, and $c = c(A) > 0$ such that the event $\mathscr{E}' = \textbf{BTR}_n^{\bm{\mathcal{L}}} (4A, B; n, L) \cap \textbf{GAP}_n^{\bm{\mathcal{L}}} \big( [-4An^{1/3}, 4An^{1/3}]; R \big)$ holds with probability at least $1 - c^{-1} e^{-c(\log n)^2}$. Hence, we may further restrict to $\mathscr{E}'$. Then, \Cref{x2p} (with the $(A, k)$ there equal to $(4A, n)$ here) yields a constant $C = C (A) > 1$ such that $\bigcap_{|t| \le A} \textbf{PFL}^{\bm{l}} \big( t; n^{-1} (\log n)^7; C)$ holds with probability at least $1 - C n^{-100}$. This establishes the theorem.
		\end{proof}

		The proof of \Cref{x2} will be a quick consequence of the following proposition (together with a high-probability H\"{o}lder bound on the paths in $\bm{\mathsf{x}}$, guaranteed by \Cref{p:initial}), to be established in \Cref{Proofx2} below. Instead of showing that $\textbf{PFL}^{\bm{x}}$ holds for all $t \in [ -A / 4, A / 4 ]$ simultaneously, it shows this statement for a fixed time $t \in [ -A / 4, A / 4 ]$.

		\begin{prop}
			
			\label{x2}
			
			Adopt the notation and assumptions in \Cref{x2p}. For any real number $t \in [ -A / 4, A / 4]$, we have   
			\begin{flalign*} 
				\mathbb{P} \bigg[ \textbf{\emph{PFL}}^{\bm{x}} \Big( t; \displaystyle\frac{(\log k)^7}{2k}; C_1 \Big)^{\complement} \cap \textbf{\emph{BTR}}^{\bm{\mathsf x}}_n (A, B; k, L) \cap \textbf{\emph{GAP}}^{\bm{\mathsf x}}_n \big( [-Ak^{1/3}, Ak^{1/3}]; R \big) \bigg] \le C_2 k^{-200}. 
			\end{flalign*} 			
		
		\end{prop}

		\begin{proof}[Proof of \Cref{x2p}] 
			
			Throughout this proof, we set $\chi = 2^{-5000}$ and $D = 2^{5000}$, and we abbreviate the events $\textbf{BTR}_n = \textbf{BTR}_n^{\bm{\mathsf{x}}} (A, B; k, L)$ and $\textbf{GAP}_n = \textbf{GAP}^{\bm{\mathsf x}}_n \big( [-Ak^{1/3}, Ak^{1/3}]; R \big)$. The proposition will follow from applying \Cref{x2} over $t$ in a $k^{-10}$-mesh of the interval $[ -Ak^{1/3} / 4, Ak^{1/3} / 4 ]$, and using a H\"{o}lder type bound for the paths in $\bm{\mathsf{x}}$ guaranteed by \Cref{p:initial}. More specifically, define the $k^{-10}$-mesh $\mathcal{T} = [ -A / 4, A / 4] \cap (k^{-10} \cdot \mathbb{Z})$, which satisfies $|\mathcal{T}| \le Ak^{10}$, and let $C_1 = C_1 (A, B, R) > 1$ denote the constant $C_1$ from \Cref{x2}.  Then, define event $\mathscr{E} = \mathscr{E}_1 \cap \mathscr{E}_2$, where
			 \begin{flalign*}
			  \mathscr{E}_1 = \bigcap_{s \in \mathcal{T}} \textbf{PFL}^{\bm{x}} \Big( t; \displaystyle\frac{(\log k)^7}{2k}; C_1 \Big); \quad \mathscr{E}_2 =  \bigcap_{j=1}^{2k} \bigcap_{\substack{|s| \le Ak^{1/3}/4 \\ |s+t| \le Ak^{1/3} / 4}} \Big\{ \big| \mathsf{x}_j (s+t) - \mathsf{x}_j (t) \big| \le 10 A k^2 t^{1/2} + k^{-D} \Big\}.   
			 \end{flalign*} 
		 	
		 	Observe that there exist constants $c = c (A, B) > 0$ and $C_2 = C_2 (A, B, D) > 0$ such that 
		 	\begin{flalign}
		 		\label{e1btr} 
		 		\mathbb{P} \big[ \mathscr{E}_1^{\complement} \cap \textbf{BTR}_n \cap \textbf{GAP}_n \big] \le C_2 k^{-150}; \qquad \mathbb{P} \big[ \mathscr{E}_2^{\complement} \cap \textbf{BTR}_n \big] \le C_2 e^{-c(\log k)^2}.
		 	\end{flalign}
	 	 
		 	\noindent Indeed, the first follows from taking a union bound in \Cref{x2} over $t \in \mathcal{T}$. The second follows from the $A' = A / 4$ case of \Cref{p:initial}, using the definition \eqref{fhr1} of the event $\textbf{FHR}$; and the facts that $L^{\chi} k^{1/3} (j \vee k)^{1/3} t \le AL^{\chi} (j \vee k)  t^{1/2} \le 2 Ak^2 t^{1/2}$ (as $L^{\chi} = k$ for $\chi = 2^{-5000}$ and $L = k^{2^{5000}}$) and $4 k^{1/6} (j \vee k)^{1/2} t^{1/2} \le 4 (j \vee k) t^{1/2} \le 8 Ak^2 t^{1/2}$ for $|t| \le Ak^{1/3} / 4$ and $j \in \llbracket 1, 2k \rrbracket$. 
		 	
		 	Now, restrict to the event $\mathscr{E}$. By \eqref{e1btr} and a union bound, it suffices to show for sufficiently large $k$ that $\textbf{PFL}^{\bm{x}} \big( t; k^{-1} (\log k)^7; C_1 \big)$ holds, for any $t \in [-A/4, A/4]$. So, fix $t_0 \in [-A/4, A/4]$ and an arbitrary element $s \in \mathcal{T}$ such that $|s - t_0| \le k^{-10}$. Since we have restricted to $\mathscr{E} \subseteq \mathscr{E}_1$, \Cref{eventpfl} for the $\textbf{PFL}$ event yields a function $\gamma_s : [0, 1] \rightarrow \mathbb{R}$ such that $\big| x_j (s) - \gamma_s (jk^{-1}) \big| \le (2k)^{-1} (\log k)^7$ for each $j \in \llbracket 1, k \rrbracket$ and $\big\| \gamma_s - \gamma_s (0) \big\|_{\mathcal{C}^{50}} \le C_1$. Set $\gamma_{t_0} = \gamma_s$, which satisfies $\big| \gamma_{t_0} - \gamma_{t_0} (0) \big\| \le C_1$ since $\gamma_s$ does. Moreover, for any integer $j \in \llbracket 1, k \rrbracket$, for sufficiently large $k$ we have 
			\begin{align*}
				\big| x_j (t_0) - \gamma_{t_0} (jk^{-1}) \big| & \le \big| x_j (s) - \gamma_s (jk^{-1}) \big| + k^{-2/3} \cdot \big| \mathsf{x}_{j+k} (t_0 k^{1/3}) - \mathsf{x}_{j+k} (sk^{1/3}) \big| \\ 
				& \le (2k)^{-1} (\log k)^7  + 11 A k^{-3} \le k^{-1} (\log k)^7,
			\end{align*}
		
			\noindent where in the first statement we used the facts that $x_j (t_0) = k^{-2/3} \cdot \mathsf{x}_{j+k} (t_0 k^{1/3})$ and $\gamma_s = \gamma_{t_0}$; in the second we used the fact that $\big| x_j (s) - \gamma_s (jk^{-1}) \big| \le (2k)^{-1} (\log k)^7$, our restriction to the event $\mathscr{E} \subseteq \mathscr{E}_2$, and the fact that $|s-t_0| \le k^{-10}$ (and that $D = 2^{5000}$); and in the third we used the fact that $k$ is sufficiently large. 
		
			This verifies that $\gamma_{t_0}$ satisfies the first bound in \eqref{xjtgammat}. Thus, $\textbf{PFL}^{\bm{x}} \big( t_0; k^{-1} (\log k)^7; C_1 \big)$ holds. Since $t_0 \in [ -A / 4, A / 4  ]$ was arbitrary, this establishes the proposition. 
		\end{proof}

		\subsection{Proof of \Cref{x2}} 
		
		\label{Proofx2} 
		
		In this section we establish \Cref{x2}, guaranteeing that $\bm{x}$ likely satisfies a regular profile event. The third part of \Cref{p:densityub0} provides a way of ensuring that families of non-intersecting Brownian bridges without lower boundary satisfy these events. We therefore require a way of coupling $\bm{\mathsf{x}}$ with such a family $\bm{\mathsf{y}}$, in such a way that their upper paths are close to each other; \Cref{c:finalcouple} does not quite do this, since the two couplings it provides are not necessarily the same. The following lemma indicates that if $L$ is sufficiently large with respect to $k$ (namely, $L \ge k^{2^{6000}}$), then there exists a coupling between $\bm{\mathsf{x}}$ and $\bm{\mathsf{y}}$ guaranteeing that their top $k^2$ paths are likely close.

		\begin{lem}
			
			\label{lklarge} 
			
			Adopt the notation and assumptions of \Cref{c:finalcouple}; assume that $L \ge k^{2^{5000}}$. For any real number $\mathsf{t} \in [ -Ak^{1/3} / 2, Ak^{1/3} / 2]$, there exists a coupling between $\bm{\mathsf{x}}$ and $\bm{\mathsf{y}}$ such that 
			\begin{flalign}
				\label{xyk2}
				\mathbb{P} \Bigg[ \bigcap_{j=1}^{k^2}  \Big\{ \big| \mathsf{x}_j (\mathsf{t}) - \mathsf{y}_j (\mathsf{t}) \big| \le k^{-2} \Big\} \Bigg] \ge 1 - C k^{-200}.
			\end{flalign}
		\end{lem} 
		
		To prove \Cref{lklarge}, we will apply a Markov estimate to the quantity $\mathsf{x}_j (\mathsf{t}) - \mathsf{y}_j (\mathsf{t})$. This will require a (weak) tail bound on the latter random variable, which is provided by the following lemma, to be established in \Cref{ProofxyB} below.
		
		\begin{lem} 
			
			\label{estimateb}
			
			Adopting \Cref{a:nkrelation}, there exist constants $c = c(A, B) > 0$ and $C = C(A, B, D) > 1$ such that the following holds. Set $n' = \big\lceil L^{1/2^{3000}} k \big\rceil$. For any real number $\mathsf{t} \in [ -Ak^{1/3} / 2, Ak^{1/3} / 2 ]$, there is an event $\mathscr{A} = \mathscr{A}_{\mathsf{t}} \subseteq \textbf{\emph{BTR}}_n^{\bm{\mathsf{x}}} (A; B)$, that is measurable with respect to $\mathcal{F}_{\ext}^{\bm{\mathsf{x}}} \big( \llbracket 1, n' \rrbracket \times (-Ak^{1/3} / 2, Ak^{1/3} / 2) \big)$, satisfying $\mathbb{P} \big[ \textbf{\emph{BTR}}_n^{\bm{\mathsf{x}}} (A; B) \setminus \mathscr{A} \big] \le c^{-1} e^{-k}$ and the following. Condition on $\mathcal{F}_{\ext}^{\bm{\mathsf{x}}} \big( \llbracket 1, n' \rrbracket \times (Ak^{1/3} / 2, Ak^{1/3} / 2) \big)$; restrict to $\mathscr{A}$; and define the $n'$-tuples $\bm{u} = \bm{\mathsf{x}}_{\llbracket 1, n' \rrbracket} ( -Ak^{1/3} / 2) \in \mathbb{W}_{n'}$ and $\bm{v} = \bm{\mathsf{x}}_{\llbracket 1, n' \rrbracket} ( Ak^{1/3} / 2) \in \mathbb{W}_{n'}$. Sample $n'$ non-intersecting Brownian bridges $\bm{\mathsf{y}} = (y_1, y_2, \ldots ,y_{n'}) \in \llbracket 1, n' \rrbracket \times [ -Ak^{1/3} / 2, Ak^{1/3} / 2]$ from the measure $\mathsf{Q}^{\bm{u}; \bm{v}}$. Then, under any coupling between $\bm{\mathsf{x}}$ and $\bm{\mathsf{y}}$, we have
			\begin{flalign}
				\label{ckixy} 
				\mathbb{P} \bigg[ \displaystyle\max_{j \in \llbracket 1, n' \rrbracket} \big| \mathsf{x}_j (\mathsf{t}) - \mathsf{y}_j (\mathsf{t}) \big| \ge c^{-1} i k^{4D} \bigg] \le C e^{-ik}, \qquad \text{for every integer $i \ge 0$}.
			\end{flalign} 
			
		\end{lem}

		\begin{proof}[Proof of \Cref{lklarge}] 
			
			Throughout this proof, we abbreviate $\textbf{BTR}_n = \textbf{BTR}_n^{\bm{\mathsf{x}}} (A; B)$. By \Cref{estimateb} and \Cref{c:finalcouple}, we deduce the existence of constants $c = c (A, B) \in (0, 1)$ and $C_1 = C_1 (A, B, D) > 1$, and events $\mathscr{A}', \mathscr{A}'' \subseteq \textbf{BTR}_n$, both measurable with respect to $\mathcal{F}_{\ext}^{\bm{\mathsf{x}}} \big( \llbracket 1, n' \rrbracket \times (-Ak^{1/3} / 2, Ak^{1/3} / 2) \big)$ (recalling $n' = \lceil L^{1/2^{3000}} k \rceil$, as we have adopted the notation of \Cref{c:finalcouple}), satisfying the following three properties. First, we have 
			\begin{flalign} 
				\label{btrnaa} 
				\mathbb{P} [\textbf{BTR}_n \setminus \mathscr{A}'] \le C e^{-k}; \qquad \mathbb{P} \big[ \textbf{BTR}_n \setminus \mathscr{A}''] \le C e^{-c (\log k)^2}.
			\end{flalign}
			
			\noindent Second, conditioning on $\mathcal{F}_{\ext}^{\bm{\mathsf{x}}} \big( \llbracket 1, n' \rrbracket \times (-Ak^{1/3} / 2, Ak^{1/3} / 2) \big)$ and restricting to $\mathscr{A}'$, we have for any coupling between $\bm{\mathsf{x}}$ and $\bm{\mathsf{y}}$ that
			\begin{flalign}
				\label{xykd4} 
				\mathbb{P} \bigg[ \displaystyle\max_{j \in \llbracket 1, n \rrbracket} \big| \mathsf{x}_j (\mathsf{t}) - \mathsf{y}_j (\mathsf{t}) \big| \ge c^{-1} ik^{4D} \bigg] \le C_1 e^{-ik}, \qquad \text{for each integer $i \ge 0$}.
			\end{flalign}
			
			\noindent Third, again conditioning on $\mathcal{F}_{\ext}^{\bm{\mathsf{x}}} \big( \llbracket 1, n' \rrbracket \times (-Ak^{1/3} / 2, Ak^{1/3} / 2) \big)$ and now restricting to $\mathscr{A}''$, there exist two coupling between $\bm{\mathsf{x}}$ and $\bm{\mathsf{y}}$. Under the first, we have $\mathsf{x}_j (s) \le \mathsf{y}_j (s)$ for each $(j, s) \in \llbracket 1, n' \rrbracket \times [ -Ak^{1/3} / 2, Ak^{1/3} / 2 ]$, almost surely. Under the second, we have (recalling $n'' = \lceil L^{1/2^{4000}} k \rceil$)
			\begin{flalign}
				\label{exjyj} 
				\mathbb{P} \big[ \mathscr{E}^{\complement} \big] \le C_1 e^{-c (\log k)^2}, \qquad \text{where} \qquad \mathscr{E} = \bigcap_{j=1}^{n''}  \big\{ \mathsf{y}_j (\mathsf{t}) \ge \mathsf{x}_j (\mathsf{t}) - L^{-1/2^{4000}} k^{2/3}\big\}.
			\end{flalign}
			
			Set $\mathscr{A} = \mathscr{A}' \cap \mathscr{A}''$, which by \eqref{btrnaa} and a union bound satisfies $\mathbb{P} [\textbf{BTR} \setminus \mathscr{A}] \ge 1 - C_1 e^{-k} - C_1 e^{-c (\log k)^2} \ge 1 - 2C_1 e^{-c (\log k)^2}$. As in the statement of the lemma (see also that of \Cref{c:finalcouple}), condition on $\mathcal{F}_{\ext}^{\bm{\mathsf{x}}} \big( \llbracket 1, n' \rrbracket \times (-Ak^{1/3} / 2, Ak^{1/3} / 2) \big)$ and restrict to $\mathscr{A}$. We will show that the coupling between $\bm{\mathsf{x}}$ and $\bm{\mathsf{y}}$ provided by the first part of \Cref{c:finalcouple} satisfies \eqref{xyk2}. This will proceed by using a Markov estimate. 
			
			In particular, we claim for sufficiently large $k$ that
			\begin{flalign}
				\label{xjtyjtk} 
				\mathbb{E} \big[ \mathsf{x}_j (\mathsf{t}) \big] - \mathbb{E} \big[ \mathsf{y}_j (\mathsf{t}) \big] \le k^{-250}, \qquad \text{for each integer $j \in \llbracket 1, k^2 \rrbracket$}.
			\end{flalign}
			
			\noindent Let us establish the lemma assuming \eqref{xjtyjtk}. Since we have restricted to $\mathscr{A}$, there exists a coupling between $\bm{\mathsf{x}}$ and $\bm{\mathsf{y}}$ such that $\mathsf{x}_j (s) \le \mathsf{y}_j (s)$ for each $(j, s) \in \llbracket 1, n' \rrbracket \times [ -Ak^{1/3} / 2, Ak^{1/3} / 2 ]$, almost surely. Hence, under this coupling, we have for sufficiently large $k$ and any integer $j \in \llbracket 1, k^2 \rrbracket$ that
			\begin{flalign*} 
				\mathbb{P} \Big[ \big| \mathsf{x}_j (\mathsf{t}) - \mathsf{y}_j (\mathsf{t}) \big| \ge k^{-2} \Big] = \mathbb{P} \big[ \mathsf{x}_j (\mathsf{t}) - \mathsf{y}_j (\mathsf{t}) \ge k^{-2} \big] \le k^2 \cdot \mathbb{E} \big[ \mathsf{x}_j (\mathsf{t}) - \mathsf{y}_j (\mathsf{t}) \big] \le k^{-210},
			\end{flalign*} 
			
			\noindent where in the first statement we used the fact that $\mathsf{x}_j \ge \mathsf{y}_j$; in the second we used a Markov bound; and in the third we used \eqref{xjtyjtk}. Taking a union bound over all $j \in \llbracket 1, k^2 \rrbracket$ then yields the lemma.
			
			It therefore remains to establish \eqref{xjtyjtk}; in what follows, we fix an integer $j \in \llbracket 1, k^2 \rrbracket$. Since $L = k^{2^{5000}}$, we have $L^{-1/2^{4000}} = k^{-2^{1000}} \le k^{-300}$; in particular, $n'' \ge L^{1/2^{4000}} k \ge k^{300} \ge k^2$, so $j \in \llbracket 1, n'' \rrbracket$. Hence, for sufficiently large $k$, \eqref{exjyj} yields  
			\begin{flalign}
				\label{xjeyje}
				\mathbb{E} \Big[ \textbf{1}_{\mathscr{E}} \cdot \big( \mathsf{x}_j (\mathsf{t}) - \mathsf{y}_j (\mathsf{t}) \big) \Big] \le L^{-1/2^{4000}} k \le C_1 k^{-299} \le k^{-298}. 
			\end{flalign}
			
			It thus remains to bound the expectation of $\mathsf{x}_j (\mathsf{t}) - \mathsf{y}_j (\mathsf{t})$ off of $\mathscr{E}$, which will make use of the tail bound \eqref{xykd4}. In particular, observe for $c' = c / 2$ and $k$ sufficiently large that
			\begin{flalign*}
				\mathbb{E} & \Big[ \textbf{1}_{\mathscr{E}^{\complement}} \cdot \big( \mathsf{x}_j (\mathsf{t}) - \mathsf{y}_j (\mathsf{t}) \big) \Big] \\
				& \le \mathbb{E} \big[ \textbf{1}_{\mathscr{E}^{\complement}} \big]^{1/2} \cdot \mathbb{E} \Big[ \big| \mathsf{x}_j (\mathsf{t}) - \mathsf{y}_j (\mathsf{t}) \big|^2 \Big]^{1/2} \\ 
				& \le C_1^{1/2} e^{-c' (\log k)^2} \cdot \Bigg( \displaystyle\sum_{i=0}^{\infty} c^{-2} (i+1)^2 k^{8D} \cdot \mathbb{P} \Big[ \big| \mathsf{x}_j (\mathsf{t}) - \mathsf{y}_j (\mathsf{t}) \big| \in \big[ c^{-1} ik^{4D}, c^{-1} (i+1) k^{4D} \big] \Big] \Bigg)^{1/2} \\
				& \le c^{-1} C_1 k^{4D} e^{-c' (\log k)^2} \cdot \bigg( \displaystyle\sum_{i=0}^{\infty} e^{-ik} (i+1)^2 \bigg)^{1/2} \le 3 c^{-1} C_1 k^{4D} e^{-c' (\log k)^2} \le k^{-298},
			\end{flalign*}
			
			\noindent where in the second inequality we used \eqref{exjyj}; in the third we used \eqref{xykd4}; in the fourth we used the bound $\sum_{i=0}^{\infty} e^{-ik} (i+1)^2 \le \sum_{i=0}^{\infty} e^{-i} (i+1)^2 < 9$; and in the fifth we used the fact that $k$ is sufficiently large. This, together with \eqref{xjeyje} and the fact that $2k^{-298} \le k^{-250}$ confirms \eqref{xjtyjtk}. 
		\end{proof}
		
		Now we can establish \Cref{x2}.

		\begin{proof}[Proof of \Cref{x2}]
			
			Recall the location event $\textbf{LOC}$ from \Cref{eventlocation} and the complete rectangle event $\textbf{CTR}$ from \Cref{ftrbtr2}. Throughout this proof, we define $n' = \lceil L^{1/3000} k \rceil$ and $\widetilde{B} = 12 A^2 B^3$; abbreviate the events $\textbf{BTR}_n = \textbf{BTR}_n^{\bm{\mathsf{x}}} (A, B; k, L)$ and $\textbf{CTR}_n = \textbf{CTR}_n^{\bm{\mathsf{x}}} (A, \widetilde{B}; k, L)$; abbreviate the $\sigma$-algebra $\mathcal{F}_{\ext} =  \mathcal{F}_{\ext}^{\bm{\mathsf{x}}} \big( \llbracket 1, n' \rrbracket \times (-Ak^{1/3} / 2, Ak^{1/3} / 2) \big)$; and define the $n'$-tuples  $\bm{u} = \bm{\mathsf{x}}_{\llbracket 1, n' \rrbracket} ( -Ak^{1/3} / 2) \in \mathbb{W}_n$ and $\bm{v} = \bm{\mathsf{x}}_{\llbracket 1, n' \rrbracket} (Ak^{1/3} / 2) \in \mathbb{W}_n$. Further sample $n'$ non-intersecting Brownian bridges $\bm{\mathsf{y}} = (\mathsf{y}_1, \mathsf{y}_2, \ldots , \mathsf{y}_{n'}) \in \llbracket 1, n' \rrbracket \times \mathcal{C} \big( [-Ak^{1/3}/ 2, Ak^{1/3} / 2] \big)$ from the measure $\mathsf{Q}^{\bm{u}; \bm{v}}$. We may assume throughout that $\mathbb{P} [ \textbf{BTR}_n^{\bm{\mathsf{x}}}] > 0$, as otherwise the proposition holds.  
			
			Define the $\mathcal{F}_{\ext}$-measurable event  
			\begin{flalign}
				\label{a1a} 
				\mathscr{A}_1 = \textbf{BTR}_n \cap \bigcap_{j=1}^{n'} \textbf{LOC}_j^{\bm{\mathsf{x}}} \Bigg( \Big\{ -\frac{Ak^{1/3}}{2}, \displaystyle\frac{Ak^{1/3}}{2} \Big\}; -\widetilde{B} j^{2/3} - \widetilde{B} k^{2/3}; \widetilde{B} k^{2/3} - \widetilde{B}^{-1} j^{2/3} \bigg).
			\end{flalign} 
			
			\noindent Since $\textbf{BTR}_n \cap \textbf{CTR}_n \subseteq \mathscr{A}_1$ by \Cref{ftrbtr2},  \Cref{ctr} yields constants $c_1 = c_1 (A, B) > 0$ and $C_3 = C_3 (A, B, D) > 1$ such that $\mathbb{P} [\textbf{BTR}_n \setminus \mathscr{A}_1] \le C_3 e^{-c_1 (\log k)^2}$. Next, by the $\mathsf{t} = tk^{1/3}$ case of \Cref{lklarge} (and altering the constants $c_1 > 0$ and $C_3 > 1$ if necessary), there exists an $\mathcal{F}_{\ext}$-measurable event $\mathscr{A}_2 \subseteq \textbf{BTR}_n$ satisfying $\mathbb{P} [\textbf{BTR}_n \setminus \mathscr{A}_2] \le C_3 e^{-c_1 (\log k)^2}$ and the following. Conditioning on $\mathcal{F}_{\ext}$ and restricting to $\mathscr{A}_2$, there exists a coupling between $\bm{\mathsf{x}}$ and $\bm{\mathsf{y}}$ such that
			\begin{flalign}
				\label{xjk2yj}
				\mathbb{P} \Bigg[ \bigcap_{j=1}^{k^2} \Big\{ \big| \mathsf{x}_j (tk^{1/3}) - \mathsf{y}_j (tk^{1/3}) \big| \ge k^{-2} \Big\} \Bigg] \le C_3 k^{-200}.
			\end{flalign}
			
			\noindent Define the $\mathcal{F}_{\ext}$-measurable event $\mathscr{A} = \mathscr{A}_1 \cap \mathscr{A}_2 \subseteq \textbf{BTR}_n$, which by a union bound satisfies 
			\begin{flalign}
				\label{btrna}  
				\mathbb{P} [\textbf{BTR}_n \setminus \mathscr{A}] \le 2C_3 e^{-c_1 (\log k)^2}.
			\end{flalign} 
			
			Condition on $\mathcal{F}_{\ext}$ and restrict to the event $\mathscr{A}$. By \eqref{btrna}, it suffices to show for some constants $C_1 = C_1 (A, B, D) > 1$ and $C_2 = C_2 (A, B, D, R) > 1$ that
			\begin{flalign}
				\label{xtk2} 
				\mathbb{P} \bigg[ \textbf{PFL}^{\bm{x}} \Big( t; \displaystyle\frac{(\log k)^7}{2k}; C_1 \Big) \cap \textbf{GAP}_n \big( [-Ak^{1/3}, Ak^{1/3}]; R \big) \bigg] \le C_2 k^{-200}. 
			\end{flalign}
			
			This will follow from \Cref{p:densityub0}; we must first verify \Cref{a:boundary} for that proposition. Denote the $n$-tuples $\bm{u} = \bm{\mathsf{x}}_{\llbracket 1, n \rrbracket} ( -Ak^{1/3} / 2 ) \in \mathbb{W}_n$ and $\bm{v} = \bm{\mathsf{x}}_{\llbracket 1, n \rrbracket} ( Ak^{1/3} / 2 ) \in \mathbb{W}_n$. Observe (by \eqref{a1a} and \Cref{eventlocation}) since we have restricted to $\mathscr{A} \subseteq \mathscr{A}_1$ that, for each integer $j \in \llbracket 1, n' \rrbracket$,
			\begin{flalign} 
				\label{uvb} 
				-\widetilde{B}j^{2/3} - \widetilde{B}k^{2/3} \le u_j \le \widetilde{B}k^{2/3} - \widetilde{B}^{-1} j^{2/3}; \quad -\widetilde{B}j^{2/3} - \widetilde{B}k^{2/3} \le v_j \le \widetilde{B}k^{2/3} - \widetilde{B}^{-1} j^{2/3}. 
			\end{flalign} 
			
			\noindent This verifies \eqref{e:xidensity} of \Cref{a:boundary}. Since we moreover have $t \in [-A / 4, A / 4]$, \Cref{a:boundary} holds with the $(\bm{\mathsf{x}}; L, n, B; t)$ there equal to $\big( \bm{\mathsf{y}}; (k^{-1} n')^{2/3}, n'; \widetilde{B}; t + A / 2 \big)$ here (with the arguments of the paths in $\bm{\mathsf{y}}$ shifted by $Ak^{1/3} / 2$). Thus, \Cref{p:densityub0} applies; its first and third parts will be the ones of relevance for us. 
			
			Its first part yields constants $c_2 = c_2 (A, B) > 0$, $C_4 = C_4 (A, B) > 1$, and $C_5 = C_5 (A, B, D) > 1$ and an event $\mathscr{E}$, with 
			\begin{flalign} 
				\label{ce} 
				\mathbb{P} [\mathscr{E}] \ge 1 - C_5 e^{-c_2 (\log k)^2},
			\end{flalign} 
			
			\noindent on which there exists a (random) measure $\mu$ with $\mu (\mathbb{R}) = k^{-1} n'$, satisfying the following property. Denoting the classical locations of $\mu$ (recall \Cref{gammarho}) by $\gamma_j = \gamma_{j; n'}^{\mu}$ for each $j \in \llbracket 1, n' \rrbracket$, we have 
			\begin{flalign}
				\label{ygamma}
				\gamma_{j + \lfloor (3D \log k)^6 \rfloor} - k^{-2} \le k^{-2/3} \cdot \mathsf{y}_j (tk^{1/3}) \le \gamma_{j - \lfloor (3D \log k)^6 \rfloor} + k^{-2}, \qquad \text{for each $j \in \llbracket 1, n' \rrbracket$},
			\end{flalign}
			
			\noindent where we have used the facts that $(n')^{-D} \le k^{-D} \le k^{-2}$ and that $3D \log k = \log k^{3D} \ge \log k^{3D/2+1} \ge \log (L^{3/2} k) = \log n \ge \log n'$. 
			
			We would also like to use the third statement in Proposition \ref{p:densityub0}, to which end we must verify its hypothesis, which we will do upon further restricting to the event 
			\begin{flalign}
				\label{eak}
				\mathscr{E}' = \bigcap_{j=1}^{k^2} \Big\{ \big| \mathsf{x}_j (tk^{1/3}) - \mathsf{y}_j (tk^{1/3}) \big| \le k^{-2} \Big\} \cap \textbf{GAP}_n \big( [-Ak^{1/3}, Ak^{1/3}]; R \big).
			\end{flalign}
			
			\noindent So, as in \eqref{e:defgammataucopy}, define the function $\gamma : [0, k^{-1} n'] \rightarrow \mathbb{R}$ by 
			\begin{flalign}
				\label{gammaygamma} 
				\gamma (y) = \sup \Bigg\{ x \in \mathbb{R} : \displaystyle\int_x^{\infty} \mu (dx) \ge y \Bigg\}, \qquad \text{so that} \qquad \gamma_j = \gamma \Big( \displaystyle\frac{2j-1}{2k} \Big), \quad \text{for each $j \in \llbracket 1, n' \rrbracket$}.
			\end{flalign}
			
			\noindent On the event $\mathscr{E}'$, \eqref{ygamma}, the first event in \eqref{eak}, and the second statement in \eqref{gammaygamma} give for sufficiently large $k$ (using the facts that $\gamma$ is non-increasing and that $(6D \log k)^6 \ge 2 \big\lfloor (3D \log k)^6 \big\rfloor + 1$) that
			\begin{align}\label{xgamma2}
				\gamma \Big( \displaystyle\frac{2j + (6D \log k)^6}{2k} \Big) - 2k^{-2} \leq k^{-2/3} \cdot \mathsf{x}_j(tk^{1/3}) \leq \gamma \Big( \displaystyle\frac{j - (6D \log k)^6}{k} \Big) + 2k^{-2}, \quad \text{for $j \in \llbracket 1, k^2 \rrbracket$}. 
			\end{align}
			
			\noindent On the event $\mathscr{E}' \subseteq \textbf{GAP}_n \big( [-A k^{1 /3}, A k^{1/3}];   R \big)$, we further have by \Cref{gap} that (as $tk^{1/3} \in [-Ak^{1/3}, Ak^{1/3}]$), for any integers $1 \le i \le j \le n$,
			\begin{align}\label{xgamma3}
				\big| \sfx_i(tk^{1/3})-\sfx_j(tk^{1/3}) \big|\leq R (j^{2/3}-i^{2/3})+ i^{-1/3} (\log k)^{25}.
			\end{align}
			
			\noindent  Hence, for $k$ sufficiently large and any real numbers $\widetilde{B}^{-1} \le y \le y' \le \widetilde{B}$ with $y' - y \ge 10 k^{-1} (\log n')^{50}$, it follows that
			\begin{align}\label{e:gapyy}
				\begin{aligned} 
					\big| \gamma (y')-\gamma (y) \big| = \gamma (y) - \gamma(y') & \leq k^{-2/3} \cdot \big( \sfx_{\lfloor yk\rfloor- (6D \log k)^6}(tk^{1/3})-\sfx_{\lceil y'k \rceil+(6D \log k)^6} \big) (tk^{1/3}) + 4k^{-2} \\
					& \le R \Big( \big(y' + k^{-1} (12 D \log k)^6  \big)^{2/3} - \big( y - k^{-1} (12 D \log k)^6 \big)^{2/3} \Big) \\
					& \qquad  + \displaystyle\frac{(\log k)^{25}}{k^{2/3} \big( yk - (6D \log k)^6 - 1 \big)^{1/3}} + 4k^{-2} \\
					& \le R \big( (y')^{2/3} - y^{2/3} \big) + 10 \widetilde{B} k^{-1} (12 D \log k)^{25} \le 2R \big( (y')^{2/3}-y^{2/3} \big).
				\end{aligned} 
			\end{align}
			
			\noindent Here, in the first statement, we used the fact that $\gamma$ is non-increasing; in the second we used \eqref{xgamma2} (which applies since  $\lceil yk \rceil + (6D \log k)^6 \le \lceil y' k \rceil + (6D \log k)^6 \le \lceil \widetilde{B} k \rceil + (6D \log k)^6 < k^2$ for sufficiently large $k$); and, in the third, we used \eqref{xgamma3}. In the fourth we used the facts that $a^{2/3} - b^{2/3} \le b^{-1/3} (a-b)$ for any real numbers $a \ge b > 0$ (applied for $(a, b) = \big( y' + k^{-1} (12 D \log k)^6, y' \big)$ and $(a, b) = \big( y, y - k^{-1} (12 D \log k)^6 \big)$), that $(2\widetilde{B})^{-1} \le y - k^{-1} (12 D \log k)^6 \le \widetilde{B}$ for sufficiently large $k$, and that $4k^{-2} \le 4 \widetilde{B} k^{-1} (12 D \log k)^{25}$. In the fifth, we used the bound $(y')^{2/3} - y^{2/3} \ge 2|y'-y| / 3 \widetilde{B}^{1/3} \ge 20 (3 \widetilde{B} k)^{-1} (\log n')^{50} \ge 4 \widetilde{B} k^{-1} (12 D \log k)^{25}$ for sufficiently large $k$ (as $n' \ge k$). 
			
			The estimate \eqref{e:gapyy} verifies the assumptions in the third statement in Proposition \ref{p:densityub0} (with the $R$ there equal to $2R$ here). Since $\widetilde{B} = 12 A^2 B^3 \ge 12$, we have $[ 1 / 6, 6 ] \subseteq [ 2 / \widetilde{B}, \widetilde{B} / 2]$, and so it follows from the third part of \Cref{p:densityub0} that there exists a constant $C_6 = C_6 (A, B, D, R) > 1$ such that
			\begin{align}\label{gammaxy}
				\big\|\gamma |_{[1/6, 6]} - \gamma (0) \big\|_{\mathcal{C}^{50}} \leq C_6
			\end{align}
			
			We now claim for sufficiently large $k$ that $\textbf{PFL}^{\bm{x}} \big( t; (2k)^{-1} (\log k)^7; C_6 \big)$ holds on $\mathscr{E} \cap \mathscr{E}'$, with the associated function $\gamma_t : [0, 1] \rightarrow \mathbb{R}$ of \Cref{eventpfl} given by $\gamma_t (x) = \gamma(x+1)$ for each $x \in [0, 1]$. That this choice of $\gamma_t$ satisfies the second bound in \eqref{xjtgammat} follows from \eqref{gammaxy}. To verify that it satisfies the first, observe for sufficiently large $k$ and any $j \in \llbracket 1, k \rrbracket$ that
			\begin{align*}
				\big| x_j (t) - \gamma_t (k^{-1} j) \big| = \bigg| k^{-2/3} \cdot \mathsf{x}_{j+k} (tk^{1/3}) - \gamma \Big( \displaystyle\frac{j+k}{k} \Big) \bigg| \le C_6 \cdot \displaystyle\frac{(6D \log k)^6}{2k} + 2k^{-2} \le \displaystyle\frac{(\log k)^7}{2k}, 
			\end{align*}
			
			\noindent where in the first statement we used the facts that $x_j (t) = k^{-2/3} \cdot \mathsf{x}_{j+k} (tk^{1/3})$ and $\gamma_t (x) = \gamma (x + 1)$; in the second we used \eqref{xgamma2} and \eqref{gammaxy}; and in the third we used that $k$ is sufficiently large. Thus, $\gamma_t$ satisfies the first bound in \eqref{xjtgammat}, and so $\textbf{PFL}^{\bm{x}} \big( t; (2k)^{-1} (\log k)^7; C_6 \big)$ holds.
			
			Hence, $\mathscr{E} \cap \mathscr{E}' \subseteq \textbf{PFL}^{\bm{x}} \big( t; (2k)^{-1} (\log k)^7; C_6 \big)$. With the bound $\mathbb{P} [\mathscr{E}] \ge 1 - C_5 e^{-c_2 (\log k)^2}$ from \eqref{ce}, the definition \eqref{eak} of $\mathscr{E}'$, \eqref{xjk2yj}, and a union bound, this gives \eqref{xtk2} and thus the proposition.
		\end{proof}

		\subsection{Proof of \Cref{estimateb}}
		
		\label{ProofxyB}
		
		In this section we establish \Cref{estimateb}, which will follow as an application of \Cref{estimatexj2} and \Cref{monotoneheight}. 
		
		\begin{proof}[Proof of \Cref{estimateb}]
			
			Throughout this proof, we condition on the $\sigma$-algebra given by $\mathcal{F}_{\ext} = \mathcal{F}_{\ext}^{\bm{\mathsf{x}}} \big( \llbracket 1, n \rrbracket \times (-Ak^{1/3}, Ak^{1/3}) \big)$, and we abbreviate the event $\textbf{BTR}_n = \textbf{BTR}_n^{\bm{\mathsf{x}}} (A; B)$. Observe that it suffices to show that there exists constants $c = c(A, B) > 0$ and $C = C(A, B, D) > 1$ such that, for any integers $i \ge 1$ and $j \in \llbracket 1, n' \rrbracket$, we have
			\begin{flalign*}
				\mathbb{P} \bigg[ \Big\{  \big| \mathsf{x}_j (\mathsf{t}) \big| \ge c^{-1} ik^{4D} \Big\} \cap \textbf{BTR}_n \bigg] \le C e^{-4ik}; \qquad \mathbb{P} \bigg[ \Big\{  \big| \mathsf{y}_j (\mathsf{t}) \big| \ge c^{-1} ik^{4D} \Big\} \cap \textbf{BTR}_n \bigg] \le C e^{-4ik}.
			\end{flalign*}
		
			\noindent Indeed, by a union bound this would yield 
			\begin{flalign}
				\label{jxy} 
				\mathbb{P} \Bigg[ \bigg\{ \displaystyle\max_{j \in \llbracket 1, n' \rrbracket} \big| \mathsf{x}_j (\mathsf{t}) - \mathsf{y}_j (\mathsf{t}) \big| \ge 2 c^{-1} ik^{4D} \bigg\} \cap \textbf{BTR}_n \Bigg] \le 2 n' C e^{-4ik}.
			\end{flalign}
			
			\noindent The bound \eqref{jxy}, together with the Markov estimate \Cref{fg0g} then yields for each integer $i \ge 1$ an event $\mathscr{A} (i) = \mathscr{A}_t (i) \subseteq \textbf{BTR}_n$, measurable with respect to $\mathcal{F}_{\ext}^{\bm{\mathsf{x}}} \big( \llbracket 1, n' \rrbracket \times (-Ak^{1/3} / 2, Ak^{1/3} / 2) \big)$, satisfying the following two properties. First, $\mathbb{P} \big[ \textbf{BTR}_n \setminus \mathscr{A} (i) \big] \le 2C n' e^{-2ik}$. Second, conditioning on $\mathcal{F}_{\ext}^{\bm{\mathsf{x}}} \big( \llbracket 1, n' \rrbracket \times (-Ak^{1/3} / 2, Ak^{1/3} / 2) \big)$ and restricting to $\mathscr{A}(i)$, we have
			\begin{flalign}
				\label{xy2a} 
				\mathbb{P} \Bigg[ \bigg\{ \displaystyle\max_{j \in \llbracket 1, n' \rrbracket} \big| \mathsf{x}_j (\mathsf{t}) - \mathsf{y}_j (\mathsf{t}) \big| \ge 2c^{-1} ik^{4D} \bigg\} \cap \textbf{BTR}_n \Bigg] \le 2Cn' e^{-2ik}.
			\end{flalign}
			
			\noindent The lemma then follows from taking $\mathscr{A} = \bigcap_{i = 1}^{\infty} \mathscr{A} (i)$, which by a union bound satisfies 
			\begin{flalign*} 
				\mathbb{P} [\textbf{BTR}_n \setminus \mathscr{A}] \le \sum_{i=1}^{\infty} \mathbb{P} \big[ \textbf{BTR}_n \setminus \mathscr{A}(i) \big] \le 2 C n \sum_{i=1}^{\infty} e^{-2ik} \le 4n' C e^{-2 k} \le C' e^{-k},
			\end{flalign*} 
			
			\noindent for some constant $C' = C' (A, B, D) > 1$, where in the last bound we used the fact that $n' \le n = L^{3/2} k \le k^{3D/2+1}$. By \eqref{xy2a} satisfies \eqref{ckixy} for any integer $i \ge 1$, with the $(c, C)$ there equal to $( c / 2, C')$ here; observe then that \eqref{ckixy} also holds for $i = 0$ (as $C' > 1$), establishing the lemma.
			
			It therefore remains to establish \eqref{jxy}; in what follows, we fix integers $i \ge 1$ and $j \in \llbracket 1, n' \rrbracket$. Let us only verify the first bound in \eqref{jxy}, as the proof of the second is entirely analogous. Since $\mathsf{x}_1 \ge \mathsf{x}_2 \ge \cdots \ge \mathsf{x}_n$, we must then show for sufficiently small $c = c(A, B) > 0$ and large $C = C(A, B, D) > 1$ that 
			\begin{flalign}
				\label{2x2}
				\mathbb{P} \Big[ \big\{  \mathsf{x}_1 (\mathsf{t}) \ge c^{-1} ik^{4D} \big\} \cap \textbf{BTR}_n \Big] \le C e^{-4ik}; \quad \mathbb{P} \Big[ \{ \mathsf{x}_n (\mathsf{t}) \le - c^{-1} ik^{4D} \big\} \cap \textbf{BTR}_n \Big] \le C e^{-4ik}.
			\end{flalign} 
			
			\noindent We only confirm the first bound in \eqref{2x2}, as the proof of the second is again very similar. Recall that we have conditioned on $\mathcal{F}_{\ext} = \mathcal{F}_{\ext}^{\bm{\mathsf{x}}} \big( \llbracket 1, n \rrbracket \times (-Ak^{1/3}, Ak^{1/3}) \big)$; we further restrict to the event $\textbf{BTR}_n$. Set $u = C_0 ik^{4D}$, for a constant $C_0 = C_0 (A, B) > B \ge 1$ to be fixed later. Denote the $n$-tuple $\bm{u} = (u, u, \ldots , u) \in \overline{\mathbb{W}}_{n}$ (where $u$ appears with multiplicity $n$), and define the function $f : [-Ak^{1/3}, Ak^{1/3}] \rightarrow \mathbb{R}$ by setting $f(s) = \mathsf{x}_{n + 1} (s)$ for each $s \in [-Ak^{1/3}, Ak^{1/3}]$. Then, sample two families of non-intersecting Brownian bridges $\bm{\mathsf{z}} = (\mathsf{z}_1, \mathsf{z}_2, \ldots , \mathsf{z}_{n}) \in \llbracket 1, n \rrbracket \times \mathcal{C} \big( [-Ak^{1/3}, Ak^{1/3}] \big)$ and $\widetilde{\bm{\mathsf{z}}} = (\widetilde{\mathsf{z}}_1, \widetilde{\mathsf{z}}_2, \ldots , \widetilde{\mathsf{z}}_{n}) \in \llbracket 1, n \rrbracket \times \mathcal{C} \big([-Ak^{1/3}, Ak^{1/3}] \big)$, from the measures $\mathsf{Q}^{\bm{u}; \bm{u}}$ and $\mathsf{Q}_f^{\bm{u}; \bm{u}}$, respectively. 
			
			We will compare $\bm{\mathsf{x}}$ to $\bm{\mathsf{z}}$ through $\widetilde{\bm{\mathsf{z}}}$, and then use \Cref{estimatexj2} to analyze $\bm{\mathsf{z}}$. To implement this, first observe since we have restricted to the event $\textbf{BTR}_n$ that \Cref{ftrbtr} implies for each $r \in \{ -Ak^{1/3}, Ak^{1/3} \}$ that $\mathsf{x}_j (r) \le \mathsf{x}_1 (r) \le Bk^{2/3} \le C_0 i k^{3D} = u$. Hence, since the law of $\bm{\mathsf{x}}$ is given by $\mathsf{Q}_f^{\bm{\mathsf{x}} (-Ak^{1/3}); \bm{\mathsf{x}} (Ak^{1/3})}$, it follows from \Cref{monotoneheight} that we may couple $\bm{\mathsf{x}}$ and $\widetilde{\bm{\mathsf{z}}}$ such that $\mathsf{x}_1 (\mathsf{t}) \le \widetilde{\mathsf{z}}_1 (\mathsf{t})$. 
			
			Next, we compare $\bm{\mathsf{z}}$ and $\widetilde{\bm{\mathsf{z}}}$ by showing that the paths in $\bm{\mathsf{z}}$ are with high probability already above the lower boundary $f$. To do this, observe by \Cref{estimatexj2} (applied with the $(t, s; a, b; n; B)$ there given by $(\mathsf{s}, Ak^{1/3}; -Ak^{1/3}, Ak^{1/3}; n; R)$ here) that there exist constants $c_1 \in (0, 1)$ and $C_1 > 2$ such that, for any real number $R \ge 1$, we have
			\begin{flalign}
				\label{zjr} 
				&\mathbb{P} \Bigg[ \displaystyle\sup_{|s| \le Ak^{1/3}} \mathsf{z}_{n} (s) \le u - 2An R \Bigg] \le  C_1 e^{C_1 n - c_1 R^2}; \qquad \mathbb{P} \Bigg[ \mathsf{z}_1 (\mathsf{t}) \ge u + 2AnR \Bigg] \le C_1 e^{C_1 n - c_1 R^2}, 
			\end{flalign} 
			
			\noindent where we used the fact that $(Ak^{1/3} - s)^{1/2} \log \big(4 Ak^{1/3} (Ak^{1/3} - s)^{-1} \big) \le 2A^{1/2} k^{1/6} \le 2An$ for each $s \in [-Ak^{1/3}, Ak^{1/3}]$. Now set $R = 3c_1^{-1} C_1 ik^{3D/2}$ and fix $C_0 = 7c_1^{-1} C_1 AB$. Observe that 
			\begin{flalign}
				\label{uanr} 
				\begin{aligned} 
					u - 2AnR & = 7c_1^{-1} C_1 AB ik^{4D} - 6 c_1^{-1} C_1 A i k^{3D/2} n  \ge AB ik^{4D} \ge \displaystyle\sup_{|s| \le Ak^{1/3}} f(s), 
				\end{aligned} 
			\end{flalign}
			
			\noindent where in the first statement we used the definitions of $u = C_0 ik^{4D}$, of $C_0$, and of $R$; in the second we used the fact that $n = L^{3/2} k \le L^{3D/2} k \le k^{5D/2}$ (and that $c_1 \in (0, 1)$ and $B, C_1 \ge 1$); and in the third we used the fact that $f(s) = \mathsf{x}_{n+1} (s) \le Bk^{2/3} \le ABik^{4D}$, which holds (by \Cref{ftrbtr}) since we have restricted to the event $\textbf{BTR}_n$. Inserting this into the first bound in \eqref{zjr} (and using the bound $C_1 n - c_1 R^2 \le -4ik$, since $C_1 n = C_1 L^{3/2} k \le C_1 k^{3D/2+1} \le C_1 k^{5D/2}$ and $c_1 R^2 = 9c_1^{-1} C_1^2 i^2 k^{3D} \ge 9 C_1 (i+1) k^{5D/2}$) yields
			\begin{flalign*}
				\mathbb{P} \Bigg[ \bigcap_{|s| \le Ak^{1/3}} \big\{ \mathsf{z}_{n} (s) \ge f(s) \big\} \Bigg] \le C_1 e^{-4i k}.
			\end{flalign*}
			
			Thus, the paths in $\bm{\mathsf{z}}$ are above the lower boundary $f$ with probability at least $1 - C_1 e^{-4ik}$, and so $\mathbb{P} \big[ \widetilde{\mathsf{z}}_1 (\mathsf{t}) \ge a \big] \le (1 - C_1 e^{-4ik})^{-1} \cdot \mathbb{P} \big[ \widetilde{\mathsf{z}}_1 (\mathsf{t}) \ge a \big]$ for any real number $a \in \mathbb{R}$. Together with the second bound in \eqref{zjr}, the fact that $u + 2AnR = 7 c_1^{-1} C_1 AB ik^{4D} + 6 c_1^{-1} C_1 Aikn \le 13 c_1^{-1} C_1 AB ik^{4D}$ (by following the same reasoning as used to deduce \eqref{uanr}), and the bound $C_1 n - c_1 R^2 \le -4ik$, this gives
			\begin{flalign*}
				\mathbb{P} \big[ \widetilde{\mathsf{z}}_1 (\mathsf{t}) \ge 13 c_1^{-1} C_1 AB ik^{4D} \big] & \le (1 - C_1 e^{-4ik})^{-1} \cdot \mathbb{P} \big[ \mathsf{z}_1 (\mathsf{t}) \le 13c_1^{-1} C_1 AB ik^{4D} \big] \\
				& \le 2 \cdot \mathbb{P} \big[ \mathsf{z}_1 (\mathsf{t}) \ge u + 2AnR \big] \le  2C_1 e^{-4ik}.
			\end{flalign*}
			
			\noindent Combining this with the existence of a coupling between $\bm{\mathsf{x}}$ and $\bm{\mathsf{z}}$ such that $\mathsf{x}_1 (\mathsf{t}) \le \widetilde{\mathsf{z}}_1 (\mathsf{t})$, this yields \eqref{2x2} and thus the lemma.
		\end{proof}

		\section{Proof of the Global Law}
		
		\label{ProofGlobal}

		In this section we establish \Cref{p:globallaw2}, indicating that a line ensemble satisfying \Cref{l0} satisfies a global law. As in \Cref{ProofRegular0}, we first in \Cref{ProofxGlobal} reduce it to a general statement, given by \Cref{p:globallaw} below, about line ensembles $\bm{\mathsf{x}}$ that likely satisfy a boundary tall rectangle and gap event. To show the latter, we will restrict $\bm{\mathsf{x}}$ to a tall rectangle, and use \Cref{c:finalcouple} to couple it to a family of non-intersecting Brownian bridges without lower boundary; we will then analyze the latter using \Cref{rhot} (stating it converges to a limit shape) and \Cref{p:limitprofile} (to analyze the edge behavior of this limit shape). To implement this, we require a variant of \Cref{rhot} that has uniform rates of convergence; we show such a statement in \Cref{MeasuremuUniform} using compactness arguments. We then establish \Cref{p:globallaw} in \Cref{Proofxabc}.
		
		Throughout this section, we recall from the beginning of this chapter that $\bm{\mathsf{x}} = (\mathsf{x}_1, \mathsf{x}_2, \ldots )$ is a $\mathbb{Z}_{\ge 1} \times \mathbb{R}$ indexed line ensemble satisfying the Brownian Gibbs property. 
		
		\subsection{Likelihood of the Global Law Event}
		
		\label{ProofxGlobal} 
		
		In this section we prove \Cref{p:globallaw2} as a consequence of the following general result, to be established in \Cref{Proofxabc} below. It states that, if $\bm{\mathsf{x}}$ satisfies a boundary tall rectangle event $\textbf{BTR}$ and a gap event $\textbf{GAP}$ then, for small $\theta$ and sufficiently large $k$, its top $\lceil \theta^3 k \rceil$ curves on the time interval $[-\theta k^{1/3}, \theta k^{1/3}]$ approximate a limiting Airy profile $\mathfrak{G}_{\Ai; \mathfrak{a}, \mathfrak{b}, \mathfrak{c}}$ from \eqref{e:Airyprofile} (with random coefficients $\mathfrak{a}, \mathfrak{b}, \mathfrak{c}$), up to error $\mathcal{O} (\theta^3 k^{2/3})$. Observe that, unlike in \Cref{x2p} where $L$ was growing faster than $k$, in the below proposition $L$ is bounded independently of $k$ (though the constant $C$ prescribing the error\footnote{However, the lower bound $K_0$ on $k$ depends on all parameters $(\theta, \varpi, B, R, L)$ involved} below does not depend on $L$). In what follows, we recall the events $\textbf{GAP}$ and $\textbf{BTR}$ from \Cref{gap} and \Cref{ftrbtr}, respectively.
		
		\begin{prop}\label{p:globallaw}
			
			For any real numbers $\theta, \varpi \in ( 0, 1 / 2)$ and $B, R, L \ge 1$, there exist constants $C = C (B, R) > 1$ and $K_0 = K_0 (\theta, \varpi, B, R, L ) > 1$ such that the following holds for any integer $k \ge K_0$. Fix integers $n \ge k$; assume that $n = L^{3/2} k$, that $L \ge C + \theta^{-2^{6000}}$, and that $\theta < C^{-1}$. If 
			\begin{flalign}
				\label{bg2}
				\mathbb{P} \Big[ \textbf{\emph{BTR}}_n^{\bm{\mathsf{x}}} (4, B; k, L) \cap \textbf{\emph{GAP}}_n^{\bm{\mathsf{x}}} \big( [-4k^{1/3}, 4k^{1/3}]; R \big) \Big] \ge 1 - \varpi,
			\end{flalign} 
		
			\noindent then there exist random variables $\mathfrak{a}, \mathfrak{b} \in \mathbb{R}$ and $\mathfrak{c} \in [C^{-1}, C]$ such that 
			\begin{align}\label{e:globallawuse}
				\mathbb{P} \left[ \bigcap_{j=1}^{\lfloor \theta^3 k \rfloor} \bigcap_{|t| \le \theta} \Bigg\{ \bigg| \mathsf{x}_j (tk^{1/3}) - k^{2/3} \cdot (\mathfrak{a} + \mathfrak{b} t - \mathfrak{c} t^2) +  \Big( \displaystyle\frac{3 \pi}{4 \mathfrak{c}^{1/2}} \Big)^{2/3} j^{2/3}  \bigg| \le C \theta^3 k^{2/3}  \Bigg\} \right] \ge 1 - 3 \varpi.
			\end{align}
		\end{prop}

		\begin{proof}[Proof of \Cref{p:globallaw2}]
				
				Throughout this proof, we recall the events $\textbf{TOP}$, $\textbf{GAP}$, and $\textbf{BTR}$ from \Cref{eventsregular1}, \Cref{gap}, and \Cref{ftrbtr}, respectively. Let $B_0, R_0 > 1$ denote the constants $B (4)$ and $R(4)$ from \Cref{ndeltal} (with the parameters $(A, D)$ there given by $(4, 4)$ here). Further let $\mathfrak{C} > 1$ denote the constant $ C (B_0, R_0)$ from \Cref{p:globallaw}, and set 
			\begin{flalign}
				\label{thetan} 
				\varpi = \displaystyle\frac{\delta}{4}; \quad \theta = \displaystyle\frac{\delta}{90 \mathfrak{C} B^3}; \quad k = \Big( \displaystyle\frac{B}{\theta} \Big)^3 n; \quad L = B + \mathfrak{C} + \theta^{-2^{6000}}; \quad N =  L^{3/2} k,   
			\end{flalign}
		
			\noindent assuming in what follows that $k$ and $N$ are integers (for otherwise we may increase $B$ and $\mathfrak{C}$, and decrease $\theta$, to ensure this to hold).
			
			By (the $A = 4$ case of) \Cref{ndeltal}, $\mathbb{P} \big[ \textbf{BTR}_N^{\bm{\mathcal{L}}} (4, B_0; k, L) \cap \textbf{GAP}_N^{\bm{\mathcal{L}}} ( [-4k^{1/3}, 4k^{1/3}]; R_0) \big] \ge 1 - \varpi$ holds for sufficiently large $k \ge n$, verifying \eqref{bg2} (with the $(\bm{\mathsf{x}}, n)$ there given by $(\bm{\mathcal{L}}, N)$ here). Hence, \Cref{p:globallaw} applies and yields random variables $\mathfrak{a}, \mathfrak{b} \in \mathbb{R}$ and $\mathfrak{c} \in [\mathfrak{C}^{-1}, \mathfrak{C}]$ such that, for sufficiently large $k$, we have $\mathbb{P}[\mathscr{E}_1] \ge 1 - 3 \varpi$, where we have defined the event
			\begin{flalign}
				\label{te1} 
				 \mathscr{E}_1 = \bigcap_{j=1}^{\lfloor \theta^3 k \rfloor} \bigcap_{|t| \le \theta} \Bigg\{ \bigg| \mathsf{x}_j (tk^{1/3}) - k^{2/3} \cdot (\mathfrak{a} + \mathfrak{b} t - \mathfrak{c} t^2) + \Big( \displaystyle\frac{3\pi}{4 \mathfrak{c}^{1/2}} \Big)^{2/3} j^{2/3}  \bigg| \le \mathfrak{C} \theta^3 k^{2/3} \Bigg\}.	
			\end{flalign}
			
			\noindent Moreover, from \Cref{x1lsmall} (with the $(B, \vartheta)$ there given by $(B\theta^{-1}, \theta)$ here), we have for sufficiently large $k \ge n$ that $\mathbb{P} [\mathscr{E}_2] \ge 1 - \varpi$, where         
			\begin{flalign}
				\label{te2}
				\mathscr{E}_2 = \textbf{TOP}^{\bm{\mathcal{L}}} \big( [- B \theta^{-1} n^{1/3}, B \theta^{-1} n^{1/3}]; \theta n^{2/3} \big).
			\end{flalign}
		
			\noindent Hence, denoting $\mathscr{E} = \mathscr{E}_1 \cap \mathscr{E}_2$, we have by a union bound that $\mathbb{P} [\mathscr{E}] \ge 1 - 4 \varpi = 1 - \delta$. It therefore suffices to show that $\textbf{GBL}_n^{\bm{\mathcal{L}}} (\delta; B)$ holds on $\mathscr{E}$, for sufficiently large $n$.
			
			To this end, restrict to the event $\mathscr{E}$; by \Cref{functiong} (and \eqref{gfunction0}), we must show that 
			\begin{flalign}
				\label{xjt212} 
				\big| \mathsf{x}_j (tn^{1/3}) + 2^{-1/2} t^2 n^{2/3} + 2^{-7/6} (3\pi)^{2/3} j^{2/3} \big| \le \delta n^{2/3}, \quad \text{for all $(j, t) \in \llbracket 1, Bn \rrbracket \times [-Bn^{1/3}, Bn^{1/3}]$}.
			\end{flalign}
		
			\noindent We will show that $\mathfrak{a} + \mathfrak{b} t - \mathfrak{c}^2 t \approx -2^{-1/2} t^2$ and $( 3\pi / 4 \mathfrak{c}^{1/2} )^{2/3} \approx 2^{-7/6} (3\pi)^{2/3}$ in \eqref{te1}, by comparing the $j=1$ case of \eqref{te1} with \eqref{te2}. First observe for sufficiently large $n$ that 
			\begin{flalign*}
				\mathfrak{C} \theta^3 k^{2/3} + \Big( \displaystyle\frac{3 \pi}{4 \mathfrak{c}^{1/2}} \Big)^{2/3} \le \mathfrak{C} \theta^3 k^{2/3} + 3 \mathfrak{C} \le 2 \mathfrak{C} \theta^3 k^{2/3} \le \displaystyle\frac{\delta n^{2/3}}{45},
			\end{flalign*}
		
			\noindent where in the first bound we used the fact that $\mathfrak{c} \in [\mathfrak{C}^{-1}, \mathfrak{C}]$; in the second we used the fact that $k \ge n$ is sufficiently large; and in the third we used \eqref{thetan}. Thus, since we have restricted to $\mathscr{E} \subseteq \mathscr{E}_1$, applying \eqref{te1} at $j = 1$ gives for $t \in [-\theta, \theta]$ that  
			\begin{flalign}
				\label{x1abc} 
				 \big| \mathsf{x}_1 (t k^{1/3}) - k^{2/3} (\mathfrak{a} + \mathfrak{b} t - \mathfrak{c} t^2) \big| \le \displaystyle\frac{\delta n^{2/3}}{45}. 
			\end{flalign}
			
			Since we have also restricted to $\mathscr{E} \subseteq \mathscr{E}_2$, \eqref{te2} (with \Cref{eventsregular1}), \eqref{thetan}, and the fact that $B \theta^{-1} n^{1/3} = k^{1/3} \ge  \theta k^{1/3}$ together imply for each $t \in [-\theta, \theta]$ that
			\begin{flalign*}
				\big| \mathsf{x}_1 (tk^{1/3}) + 2^{-1/2} t^2 k^{2/3} \big| \le \theta n^{2/3} < \displaystyle\frac{\delta n^{2/3}}{90}.
			\end{flalign*}
		
			\noindent Together with \eqref{x1abc} (and \eqref{thetan}), this gives
			\begin{flalign}
				\label{tabc}
				\displaystyle\sup_{|t| \le \theta} \big| \mathfrak{a} + \mathfrak{b} t - t^2 (\mathfrak{c} - 2^{-1/2}) \big| \le \displaystyle\frac{\delta}{30} \cdot (k^{-1} n)^{2/3} = \displaystyle\frac{\delta \theta^2}{30B^2} .
			\end{flalign} 
		
			\noindent Adding \eqref{tabc} at $t \in \{ -\theta, \theta \}$ and subtracting twice of it at $t = 0$ yields $2 \theta^2 |\mathfrak{c} + 2^{1/2}| \le 2 \delta \theta^2 / 15B^2$, so that
			\begin{flalign*}
				\displaystyle\frac{3}{5} \le 2^{-1/2} - \displaystyle\frac{\delta}{15B^2} \le \mathfrak{c} \le 2^{-1/2} + \displaystyle\frac{\delta}{15B^2} \le \displaystyle\frac{4}{5}.
			\end{flalign*}
		
			\noindent In particular, since $|a^{-1/3} - b^{-1/3}| \le \big( |a-b| / 3 \big) \cdot \max \{ a^{-4/3}, b^{-4/3} \}$ for $a, b > 0$, it follows that 
			\begin{flalign}
				\label{c2}
				\bigg|  2^{-7/6} (3 \pi)^{2/3} - \Big( \displaystyle\frac{3 \pi}{4\mathfrak{c}^{1/2}} \Big)^{2/3} \bigg| \le \Big( \displaystyle\frac{3\pi}{4} \Big)^{2/3} \cdot \displaystyle\frac{1}{3} \Big( \displaystyle\frac{3}{5} \Big)^{-4/3}\cdot \displaystyle\frac{\delta}{15B^2} \le \displaystyle\frac{\delta}{5B^2}.
			\end{flalign}
		
			 Since by \eqref{thetan} we have $\lfloor \theta^3 k \rfloor = \lfloor B^3 n \rfloor$ it follows that, for any $(j, t) \in \llbracket 1, Bn \rrbracket \times [-\theta, \theta]$, 
			\begin{flalign*}
				\big| \mathsf{x}_j & (tk^{1/3}) + 2^{-1/2} t^2 k^{2/3} + 2^{-7/6} (3\pi)^{2/3} j^{2/3} \big| \\
				& \le \bigg| \mathsf{x}_j (tk^{1/3}) - k^{2/3} \cdot (\mathfrak{a} + \mathfrak{b} t - \mathfrak{c} t^2) + \Big( \displaystyle\frac{3\pi}{4 \mathfrak{c}^{1/2}} \Big)^{2/3} j^{2/3} \bigg| + j^{2/3} \cdot \bigg| 2^{-7/6} (3 \pi)^{2/3} - \Big( \displaystyle\frac{3 \pi}{4 \mathfrak{c}} \Big)^{2/3} \bigg| \\ 
				& \qquad  + k^{2/3} \cdot \big| \mathfrak{a} + \mathfrak{b} t - t^2 (\mathfrak{c} - 2^{-1/2}) \big| \\
				& \le \mathfrak{C} \theta^3 k^{2/3} + \displaystyle\frac{\delta j^{2/3}}{5B^2} + \displaystyle\frac{\delta \theta^2 k^{2/3}}{30B^2} \le \displaystyle\frac{\delta n^{2/3}}{90B} + \displaystyle\frac{\delta n^{2/3}}{5B} + \displaystyle\frac{\delta n^{2/3}}{30} \le \delta n^{2/3},
			\end{flalign*}
		
			\noindent In the second bound we used the fact that we are restricting to $\mathscr{E} \subseteq \mathscr{E}_1$ (with \eqref{te1}), \eqref{c2}, and \eqref{tabc}; in the third we used \eqref{thetan} and the fact that $j \le Bn$; and in the fourth we used the fact that $B > 1$. Since $[-\theta k^{1/3}, \theta k^{1/3}] = [-Bn^{1/3}, Bn^{1/3}]$, this (upon replacing $t$ by $B^{-1} \theta t = (k^{-1} n)^{1/3} t$) establishes \eqref{xjt212} and thus the theorem.
		\end{proof}

		\subsection{Uniform Convergence to Bridge-Limiting Measure Processes} 
		
		\label{MeasuremuUniform} 
		
		 To show \Cref{p:globallaw}, we will use some form of \Cref{rhot}, stating convergence of non-intersecting Brownian bridges without upper and lower boundaries to a limit shape. Observe that \Cref{rhot} assumes that the limiting starting and ending data for this family (denoted there by $\mu_a$ and $\mu_b$) have been fixed in advance; this will not be the case in our context. So, in this section we provide a variant of that result applying to all boundary data, subject to certain conditions, uniformly.
		 		
		 To state this result, we first require the following set of measures. The second condition below may be viewed as a continuum analog for the first intersection of $\textbf{LOC}$ events (recall \Cref{eventlocation}) appearing in the $\textbf{BTR}$ event of \Cref{ftrbtr}), and the third as one for the $\textbf{GAP}$ event from \Cref{gap}. Observe that these conditions also serve as reformulations of those in \Cref{integralmu0mu1} and \Cref{gapmu0mu1}, which will enable us to use results from \Cref{EDGESHAPE} (such as \Cref{p:limitprofile}). 
		 		
		\begin{definition} 
			
			\label{plu} 
			
			For any real numbers $L \ge 1$, $U > 0$, and $\varkappa \ge 0$, let  $\mathscr P(L;U; \varkappa)\subset \mathscr P_{\fin}$\index{P@$\mathscr{P} (L; U;\varkappa)$} denote the set of measures $\mu$ satisfying the following three properties.
		\begin{enumerate}
			\item \label{e:it1} We have $\mu(\mathbb{R}) = L^{3/2}$ and $\supp \mu \subseteq [-UL, U]$.
			\item \label{e:it2}For any real number $x \le -1$, we have $\mu \big( [x, \infty) \big) \leq U|x|^{3/2}$.
			\item \label{e:it3} Define (analogously to \eqref{gty}) the function $G = G^{\mu}: [0, L^{3/2}] \rightarrow \mathbb{R}$ by setting 
			\begin{flalign}
				\label{gy}
				G(y) = \displaystyle\sup \Big\{ x \in \mathbb{R} : \mu  \big( [x, U] \big) \ge y \Big\}, \qquad \text{for each $y \in [0, L^{3/2}]$}.
				\end{flalign} 
			
				\noindent  Then, for any real numbers $0\leq x\leq y\leq L^{3/2}$, we have $G(x)-G(y)\leq U(y^{2/3}-x^{2/3}) + \varkappa$.
		\end{enumerate}
	
		\noindent Further set $\mathscr{P}(L;U) = \mathscr{P}(L;U;0)$\index{P@$\mathscr{P}(L;U)$} in the case $\varkappa = 0$.
		
		\end{definition} 
	
		Now we can state the following proposition, to be established at the end of this section, which indicates the following. Given some $n$-tuples $\bm{u}$ and $\bm{v}$ satisfying variants of the gap event (from \Cref{gap}) and of the boundary tall rectangle event (from \Cref{ftrbtr}), one can find a limiting bridge-limiting measure process $\bm{\mu}$, with boundary data in some $\mathscr{P} (L; U)$, that approximates non-intersecting Brownian bridges with starting data $\bm{u}$ and ending data $\bm{v}$. Below, we recall the notation $\emp$ from \eqref{aemp}, the L\'{e}vy metric $d_{\dL}$ from \eqref{munudistance1}, and notation on measure-valued processes and bridge-limiting measure processes from \Cref{Limit0} (and \Cref{mutmu0mu1}). 
		
		\begin{prop} 
			
			\label{pluconverge} 
			
			For any real numbers $\theta > 0$ and $B, L \ge 1$, there exists a constant $C = C(B, L, \theta) > 1$ such that the following holds. Let $n \ge k \ge C$ be integers with $n = L^{3/2} k$. Also let $\bm{u}, \bm{v} \in \overline{\mathbb{W}}_n$ be $n$-tuples such that, for any index $a \in \{ u, v \}$ and integers $1 \le i \le j \le n$, we have
			\begin{flalign}
				\label{bbaa} 
				-Bj^{2/3} - Bk^{2/3} \le a_j \le Bk^{2/3} - B^{-1} j^{2/3}, \quad \text{and} \quad a_i - a_j \le B (j^{2/3}-i^{2/3}) + (\log n)^{30} i^{-1/3}.
			\end{flalign}
			
			\noindent Sample $n$ non-intersecting Brownian bridges $\bm{\mathsf{x}} = (\mathsf{x}_1, \mathsf{x}_2, \ldots , \mathsf{x}_n) \in \llbracket 1, n \rrbracket \times \mathcal{C} \big( [-2k^{1/3}, 2k^{1/3}] \big)$ from the measure $\mathsf{Q}^{\bm{u}; \bm{v}}$, and define the measure-valued process $\bm{\mu} = (\mu_t)_{t \in [0, 1]} \in \mathcal{C} \big( [0, 1]; \mathscr{P}_{\fin} \big)$ by setting 
			\begin{flalign}
				\label{lmut}
				\mu_t = L^{3/2} \cdot \emp \Big( (2k^{2/3})^{-1}  \cdot \bm{\mathsf{x}} \big( 2(2t-1) \cdot k^{1/3} \big) \Big), \qquad \text{for each $t \in [0, 1]$}.
			\end{flalign}
		
			\noindent Then, there exist measures $\nu_0, \nu_1 \in \mathscr{P} (L; 10B^3)$ such that $\mathbb{P} \big[ d_{\dL} (\bm{\mu}, \bm{\nu}) \le \theta \big] \ge 1 - \theta$, where $\bm{\nu} \in \mathcal{C} \big( [0, 1]; \mathscr{P}_{\fin} \big)$ denotes the bridge-limiting measure process on $[0, 1]$ with boundary data $(\nu_0; \nu_1)$.

		\end{prop} 
		
		The proof of \Cref{pluconverge} proceeds by combining \Cref{rhot} with the following three results. The first indicates that the set of measures $\mathscr{P} (L; U; \varkappa)$ from \Cref{plu} is compact; this quickly implies the second, which states that any measure in $\mathscr{P} (L; U; \varkappa)$ is uniformly approximated by one in $\mathscr{P}(L; U)$, if $\varkappa$ is sufficiently small. The third indicates that, given any sequence $\bm{a}$ satisfying \eqref{bbaa}, there exists a measure in $\mathscr{P} (L; 4B^3)$ approximating the (shifted and rescaled) empirical measure associated with $\bm{a}$.

		\begin{lem}\label{l:compact}
			The set $\mathscr P(L;U; \varkappa)$ from \Cref{plu} is compact under the L{\'e}vy metric.
		\end{lem}
		\begin{proof}
			
			Fix some sequence of measures $\mu_1, \mu_2, \ldots \in \mathscr{P} (L; U; \varkappa)$. Since each $\supp \mu_j \subseteq [-UL, U]$, this sequence is tight and therefore admits a weak limit $\mu \in \mathscr{P}_{\fin}$ with $\lim_{j \rightarrow \infty} d_{\dL} (\mu_j, \mu) = 0$. To show $\mathscr{P} (L; U; \varkappa)$ is compact, we must verify that $\mu \in \mathscr{P} (L; U; \varkappa)$. That $\mu$ satisfies \Cref{e:it1} of \Cref{plu} follows from the fact that each $\mu_j$ does. Moreover, by weak convergence, we have 
			\begin{flalign*} 
				\mu \big( [x,\infty) \big) =\lim_{x'\rightarrow x^-}\mu \big((x',U] \big)\leq \lim_{x'\rightarrow x^-} \bigg( \lim_{m\rightarrow \infty}\mu_m \big((x',U] \big) \bigg) \le U \cdot \displaystyle\lim_{x' \rightarrow x^-} |x'|^{3/2} = U|x|^{3/2},
			\end{flalign*} 
		
			\noindent and so $\mu$ also satisfies \Cref{e:it2} of \Cref{plu}. Defining $G_m = G^{\mu_m} : [0, L^{3/2}] \rightarrow \mathbb{R}$ and $G = G^{\mu} : [0, L^{3/2}] \rightarrow \mathbb{R}$ as in \eqref{gy}, we have by weak convergence that, for any $0 \le x \le y \le L^{3/2}$,
			\begin{flalign*} 
				G(x)-G(y) & \leq \lim_{\varepsilon\rightarrow 0} \bigg( \lim_{m\rightarrow \infty} \big( G_m(x-\varepsilon)-G_m(y+\varepsilon) \big) \bigg) \\
				& \le U \cdot \lim_{\varepsilon \rightarrow 0} \big( (y+\varepsilon)^{2/3} - (x - \varepsilon)^{2/3} \big) + \varkappa \leq U(y^{2/3}-x^{2/3}) + \varkappa,
			\end{flalign*} 
		
			\noindent Thus, $\mu$ also satisfies \Cref{e:it3} of \Cref{plu}, and we conclude that $\mu \in \mathscr P(L;U; \varkappa)$.
		\end{proof}

		\begin{cor} 
			
			\label{plukappa0}
			
			Fix real numbers $L \ge 1$ and $U, \theta > 0$. There exists a constant $\varkappa = \varkappa (L, U, \theta)$ such that, for any measure $\mu \in \mathscr{P} (L; U; \varkappa)$, there exists a measure $\nu \in \mathscr{P}(L; U)$ with $d_{\mathrm{L}} (\mu, \nu) < \theta$. 
		\end{cor} 
		
		\begin{proof} 
			
			Assume to the contrary that this is false. Then, for each integer $j \ge 1$, there exists a measure $\mu_j \in \mathscr{P} (L; U; j^{-1})$ such that $d_{\mathrm{L}} (\mu_j, \nu) \ge \theta$ for all $\nu \in \mathscr{P}(L; U)$. Since $\mathscr{P}(L; U; 1)$ is compact under the L\'{e}vy metric by \Cref{l:compact}, the $(\mu_j)$ admit a limit $\nu \in \mathscr{P} (L; U; 1)$, which satisfies $d_{\mathrm{L}} (\mu_j, \nu) < \theta$ for sufficiently large $j$. Since $\nu \in \bigcap_{j=1}^{\infty} \mathscr{P} (L; U; j^{-1}) = \mathscr{P} (L; U)$ by the compactness \Cref{l:compact} of each $\mathscr{P} (L; U; j^{-1})$, this is a contradiction, which confirms the corollary.
		\end{proof}
	
		\begin{lem}\label{l:constructmeasure}
			
			For any real numbers $B, L \ge 1$ and $\theta > 0$, the following holds for any sufficiently large integer $k \ge 1$. Let $n \ge k$ be an integer with $n = L^{3/2} k$, and let $\bm{a} = (a_1, a_2, \ldots , a_n) \in \overline{\mathbb{W}}_n$ be an $n$-tuple such that \eqref{bbaa} holds for any integers $1 \le i \le j \le n$. Defining the measure $\mu = L^{3/2} \cdot \emp \big( (2k^{2/3})^{-1} \cdot \bm{a} \big) \in \mathscr{P}_{\fin}$, there exists a measure $\nu \in \mathscr{P} (L; 4B^3)$ satisfying $d_{\rL} (\mu, \nu) \le \theta$. 
			
		\end{lem}
	
		\begin{proof}

			Abbreviating $U = 4B^3$, let us verify that 
			\begin{flalign} 
				\label{mu0} 
				\mu \in \mathscr{P} (L; U; 3Bk^{-2/3} (\log n)^{30} \big).
			\end{flalign}

			\noindent To verify the first statement in \Cref{plu}, observe that $\supp \mu \subseteq [-BL-B, B] \subseteq [ -UL, U )$, where the first statement holds by the first statement in \eqref{bbaa} and the last statement holds since $B, L \ge 1$. To verify the second, recalling \eqref{gy}, it suffices (as $\supp \mu \subseteq [-UL, U)$) to show for any $y \ge U \ge 1$ that we have $G(y) < -(y/U)^{2/3}$. To that end, observe that 
			\begin{flalign*} 
				G (y) = k^{-2/3} \cdot a_{\lceil ky \rceil} \le B - B^{-1} y^{2/3} < -(yU^{-1})^{2/3}.
			\end{flalign*} 
		
			\noindent Here, the first statement follows from the definitions of $G$ and $\mu$, the second from \eqref{bbaa}, and the third from the fact that $(B^{-1} - U^{-2/3}) y^{2/3} \ge (2B)^{-1} y^{2/3} > B$ (where the first holds since $U = 4B^3 \ge (2B)^{3/2}$, as $B \ge 1$, and the second holds since $y^{2/3} \ge U^{2/3} > 2B^2$). To verify the third statement in \Cref{plu}, observe that 
			\begin{flalign*}
				G(x) - G(y) = k^{-2/3} \cdot (a_{\lceil kx \rceil} - a_{\lceil ky \rceil}) & \le B k^{-2/3} \big(\lceil ky \rceil^{2/3} - \lceil kx \rceil^{2/3} \big) + k^{-2/3} (\log n)^{30} \\
				& \le B (y^{2/3} - x^{2/3}) + 3Bk^{-2/3} (\log n)^{30},
			\end{flalign*}
		
			\noindent where the first statement holds by the definitions of $G$ and $\mu$, the second by the second statement of \eqref{bbaa}, and the third by the bound $Bk^{-2/3} \cdot \big( \lceil r \rceil^{2/3} - r^{2/3} \big) \le Bk^{-2/3}$ for any real number $r \ge 0$. 
			
			This confirms \eqref{mu0}. Now, the lemma follows from \Cref{plukappa0}, upon taking $k$ sufficiently large so that the $3Bk^{-2/3} (\log n)^{30}$ here is smaller than the $\varkappa (L, U, \theta)$ there.				
		\end{proof}
	
		Now we can establish \Cref{pluconverge}.

		\begin{proof}[Proof of \Cref{pluconverge}]
			
			We first set some notation. Abbreviate the constant $U = 4B^3$. Since $\mathscr{P} (L; U)$ is compact by \Cref{l:compact}, there exists an integer $K = K(B, L, \theta) \ge 1$ and measures $\nu^{(1)}, \nu^{(2)}, \ldots , \nu^{(K)} \in \mathscr{P} (L; U)$ so that, for any measure $\nu \in \mathscr{P} (L; U)$, there exists some $i = i(\nu) \in \llbracket 1, K \rrbracket$ with $d_{\dL} \big( \nu, \nu^{(i)} \big) < \theta / 4$. For any integers $i_1, i_2 \in \llbracket 1, K \rrbracket$, let $\bm{\nu}^{(i_1, i_2)} \in \mathcal{C} \big( [0, 1]; \mathscr{P}_{\fin} \big)$ denote the bridge-limiting measure process on $[0, 1]$ with boundary data $\big( \nu^{(i_1)}; \nu^{(i_2)} \big)$. 
			
			Given any integers $j \ge 1$ and $i \in  \llbracket 1, K \rrbracket$, set $J = \lceil L^{3/2} j \rceil$ and fix a $J$-tuple $\bm{a}^{(i; j)} \in \overline{\mathbb{W}}_J$ such that the following holds. Defining the measure $\mu^{(i; j)} = L^{3/2} \cdot \emp \big( (2j^{2/3})^{-1} \cdot \bm{a}^{(i; j)} \big)$, we have 
			\begin{flalign}
				\label{muijnui}
				\lim_{j \rightarrow \infty} d_{\dL} \big( \mu^{(i; j)}, \nu^{(i)} \big) = 0.
			\end{flalign}

			\noindent We will omit the ceilings in what follows, assuming that $J = L^{3/2} j$, as this will barely affect the proof. Next, for any integers $i_1, i_2 \in \llbracket 1, K \rrbracket$, sample a family of $J$ non-intersecting Brownian bridges $\bm{\mathsf{y}}^{(i_1, i_2; j)} = \big( \mathsf{y}_1^{(i_1, i_2; j)}, \mathsf{y}_2^{(i_1, i_2; j)}, \ldots , \mathsf{y}_J^{(i_1, i_2; j)} \big) \in \llbracket 1, J \rrbracket \times \mathcal{C} \big( [-2j^{1/3}, 2j^{1/3}] \big)$ from the measure $\mathsf{Q}^{\bm{a}^{(i_1; j)}; \bm{a}^{(i_2; j)}}$. Then the rescaled paths $\big( y_1^{(i_1, i_2; j)}, y_2^{(i_1, i_2; j)}, \ldots , y_J^{(i_1, i_2; j)} \big) \in \llbracket 1, J \rrbracket \times \mathcal{C} \big( [0, 1] \big)$, defined by setting $y_h^{(i_1, i_2; j)} (t)  = (2j^{2/3})^{-1} \cdot \mathsf{y}_h^{(i_1, i_2; j)} \big(  2(2t-1) j^{1/3}  \big)$ for each $(h, t) \in \llbracket 1, J \rrbracket \times [0, 1]$, are non-intersecting Brownian bridges with variances $j^{-1} = L^{3/2} J^{-1}$. So, by \eqref{muijnui}, \Cref{rhot} yields  
			\begin{flalign*} 
				\lim_{j \rightarrow \infty} \mathbb{P} \Big[ d_{\dL} \big( \bm{\mu}^{(i_1, i_2; j)}, \bm{\nu}^{(i_1, i_2)} \big) < \displaystyle\frac{\theta}{4} \Big] = 1,
			\end{flalign*} 
		
			\noindent where we have defined the measure-valued process $\bm{\mu}^{(i_1, i_2; j)} = \big( \mu_t^{(i_1, i_2; j)} \big) \in \mathcal{C} \big( [0, 1]; \mathscr{P}_{\fin} \big)$ by setting 
			\begin{flalign}
				\label{lmut2} 
				\mu_t^{(i_1, i_2; j)} = L^{3/2} \cdot \emp \Big( (2j^{2/3})^{-1} \cdot \bm{\mathsf{y}} \big(2(2t-1) j^{1/3} \big) \Big), \qquad \text{for each $t \in [0, 1]$}.
			\end{flalign}
		
			\noindent  This yields a constant $C_1 = C_1 (B, L, \theta) > 1$ such that, for any integers $i_1, i_2 \in \llbracket 1, K \rrbracket$,  
			\begin{flalign}
				\label{mui1i2j} 
				 \mathbb{P} \Big[ d_{\dL} \big( \bm{\mu}^{(i_1, i_2; j)}, \bm{\nu}^{(i_1, i_2)} \big) < \displaystyle\frac{\theta}{4} \Big] \ge 1 - \displaystyle\frac{\theta}{2}, \qquad \text{whenever $j \ge C_1$}. 
			\end{flalign}
		
			Now, recall from \eqref{lmut} that $\mu_0$ and $\mu_1$ are given by $\mu_0 = L^{3/2} \cdot \emp \big( (2k^{2/3})^{-1} \cdot \bm{u} \big)$ and $\mu_1 = L^{3/2} \cdot \emp \big( (2k^{2/3})^{-1} \cdot \bm{v} \big)$. By \Cref{l:constructmeasure}, there exists a constant $C_2 = C_2 (B, L, \theta) > 1$ and measures $\nu_0', \nu_1' \in \mathscr{P} (L; U)$ such that 
			\begin{flalign} 
				\label{muanua} 
				d_{\dL} (\mu_0, \nu_0') < \displaystyle\frac{\theta}{4}, \quad \text{and} \quad  d_{\dL} (\mu_1, \nu_1') < \displaystyle\frac{\theta}{4},
			\end{flalign}
		
			\noindent whenever $k \ge C_2$. Setting $C = \max \{ C_1, C_2 \}$, we assume for the remainder of this proof that $k \ge C$. 
			
			Fix integers $i_1, i_2 \in \llbracket 1, K \rrbracket$ satisfying 
			\begin{flalign}
				\label{nu01i1i2}  
				d_{\dL} \big( \nu_0', \nu^{(i_1)} \big) < \displaystyle\frac{\theta}{4}; \qquad d_{\dL} \big( \nu_1', \nu^{(i_2)} \big) < \displaystyle\frac{\theta}{4}. 
			\end{flalign}
		
			\noindent Also observe that, by \eqref{mui1i2j} and \eqref{munut} (taken at $t \in \{ 0, 1 \}$), we have $d_{\dL} \big( \mu^{(i_1; k)}, \nu^{(i_1)} \big) < \theta / 4$ and $d_{\dL} \big( \mu^{(i_2; k)}, \nu^{(i_2)} \big) < \theta / 4$ (where these events hold deterministically, as $\big( \mu^{(i_1; k)}; \mu^{(i_2; k)} \big)$ constitutes the deterministic boundary data for $\bm{\mathsf{y}}^{(i_1, i_2; k)}$). Together with \eqref{muanua} and \eqref{nu01i1i2}, this gives  
			\begin{flalign}
				\label{mumua}
				d_{\dL} \big( \mu_0, \mu^{(i_1; k)} \big) < \displaystyle\frac{3 \theta}{4}, \qquad \text{and} \qquad d_{\dL} \big( \mu_1, \mu^{(i_2; k)} \big) < \displaystyle\frac{3 \theta}{4}.
			\end{flalign}
		
			It then suffices to show that it is possible to couple $\bm{\mathsf{x}}$ and $\bm{\mathsf{y}}^{(i_1, i_2; k)}$ in two ways, so that under the first coupling we have for each $(x, t) \in \mathbb{R} \times [0, 1]$ that
			\begin{flalign}
				\label{mumui1i2}
				\displaystyle\int_x^{\infty} \mu_t (dr) \le \displaystyle\int_{x - 3 \theta/4}^{\infty} \mu_t^{(i_1, i_2; k)} (dr) + \displaystyle\frac{3 \theta}{4}.
			\end{flalign}
		
			\noindent and under the second we have for each $(x, t) \in \mathbb{R} \times [0, 1]$ that 
			\begin{flalign} 
				\label{0mumui1i2} 
					\displaystyle\int_x^{\infty} \mu_t (dr) \ge \displaystyle\int_{x+3 \theta/4}^{\infty} \mu_t^{(i_1, i_2; k)} (dr) - \displaystyle\frac{3 \theta}{4}. 
			\end{flalign} 
		
			\noindent Indeed, assuming \eqref{mumui1i2} and \eqref{0mumui1i2}, set $\nu_0 = \nu^{(i_1)} \in \mathscr{P} (L; U)$ and $\nu_1 = \nu^{(i_2)} \in \mathscr{P} (L; U)$, so that $\bm{\nu} = \bm{\nu}^{(i_1, i_2)} \in \mathcal{C} \big( [0, 1]; \mathscr{P}_{\fin} \big)$ is the bridge-limiting measure process on $[0, 1]$ with boundary data $(\nu_0; \nu_1)$. Then, \eqref{mui1i2j} (with the definition \eqref{munudistance1} of $d_{\dL}$), \eqref{mumui1i2}, and \eqref{0mumui1i2} together imply that 
			\begin{flalign*} 
				& \mathbb{P} \Bigg[ \bigcap_{x \in \mathbb{R}} \bigcap_{t \in [0, 1]} \bigg\{ \displaystyle\int_x^{\infty} \mu_t (dr)  \le \displaystyle\int_{x - \theta}^{\infty} \nu_t^{(i_1, i_2)} (dr) + \theta \bigg\}\Bigg] \ge 1 - \displaystyle\frac{\theta}{2}; \\ 
				& \mathbb{P} \Bigg[ \bigcap_{x \in \mathbb{R}} \bigcap_{t \in [0, 1]} \bigg\{ \displaystyle\int_x^{\infty} \mu_t (dr) \ge \displaystyle\int_{x + \theta}^{\infty} \nu_t^{(i_1, i_2)} (dr) - \theta \bigg\} \Bigg] \ge 1 - \displaystyle\frac{\theta}{2}.
			\end{flalign*} 
		
			\noindent Together with a union bound, the fact that $\bm{\nu} = \bm{\nu}^{(i_1, i_2)}$, and the definition \eqref{munudistance1} of $d_{\dL}$, this implies the proposition.
			
			It therefore remains to find a coupling between $\bm{\mathsf{x}}$ and $\bm{\mathsf{y}}^{(i_1, i_2; k)}$ such that \eqref{mumui1i2} holds, and one such that \eqref{0mumui1i2} does. Both follow in a very similar way from \Cref{monotoneheight}, so let us only implement the former. To this end, observe since $\mu_0 = L^{3/2} \cdot \emp \big( (2k^{2/3})^{-1} \cdot \bm{u} \big)$ and $\mu_0^{(i_1, i_2; k)} = \mu^{(i_1; k)} = L^{3/2} \cdot \emp \big( (2k^{2/3})^{-1} \cdot \bm{a}^{(i_1; k)} \big)$ that \eqref{mumua} (with \eqref{munudistance1}) yields 
			\begin{flalign*} 
				u_h \le a_{h - \lfloor 3 \theta L^{3/2} / 4n \rfloor}^{(i_1; k)} +  \displaystyle\frac{3\theta}{4} \cdot 2 k^{2/3} = a_{h - \lfloor 3 \theta k / 4 \rfloor}^{(i_1; k)} + \displaystyle\frac{3\theta}{4} \cdot  (2A)^{1/2} k^{2/3}, 
			\end{flalign*} 
		
			\noindent for any integer $h \in \llbracket 1, n \rrbracket$ (where we also used the fact that $L^{-3/2} n = k$). Similarly, $v_h \le a_{h - \lfloor 3 \theta k / 4 \rfloor}^{(i_2; k)} + 3 \theta k^{2/3} / 2$ for any $h \in \llbracket 1, k \rrbracket$. Thus, since the laws of $\bm{\mathsf{x}}$ and $\bm{\mathsf{y}}$ are given by $\mathsf{Q}^{\bm{u}; \bm{v}}$ and $\mathsf{Q}^{\bm{a}^{(i_1; k)}; \bm{a}^{(i_2; k)}}$, respectively, \Cref{monotoneheight} yields a coupling between these two ensembles such that $\mathsf{x}_h (tk^{1/3}) \le \mathsf{y}_{h - \lfloor 3 \theta k / 4 \rfloor}^{(i_1, i_2; k)} (tk^{1/3}) + 3 \theta k^{2/3} / 2$, or equivalently
			\begin{flalign*} 
			(2k^{2/3})^{-1} \cdot \mathsf{x}_h (tk^{1/3}) \le  (2k^{2/3})^{-1} \cdot \mathsf{y}_{h - \lfloor 3 \theta k / 4 \rfloor}^{(i_1, i_2; k)} (tk^{1/3}) + \displaystyle\frac{3\theta}{4}, \quad \text{for all $(h, t) \in \llbracket 1, n \rrbracket \times [0, 1]$}.
			\end{flalign*} 
		
			\noindent This, together with \eqref{lmut} and \eqref{lmut2},  establishes \eqref{mumui1i2} and thus the proposition.
		\end{proof}

		\subsection{Proof of \Cref{p:globallaw}}
		
		\label{Proofxabc}

		In this section we establish \Cref{p:globallaw}.

		\begin{proof}[Proof of \Cref{p:globallaw}]

			Throughout this proof, we recall the notation $\emp$ from \eqref{aemp}; that on measure-valued processes, bridge-limiting measures, and inverted height functions from \Cref{Limit0}; the set of measures $\mathscr{P} (L; U)$ from \Cref{plu}; and the completed rectangle event $\textbf{CTR}$ from \Cref{ftrbtr2}. Let us briefly outline how we will proceed. First, using \Cref{c:finalcouple}, we will exhibit a coupling between $\bm{\mathsf{x}}$ and an ensemble $\bm{\mathsf{y}}$ of non-intersecting Brownian bridges without lower boundary, so that their upper paths are close with high probability. Using \Cref{pluconverge}, we will show that the empirical measure of the latter converges under the L\'{e}vy metric to a bridge limiting measure $\bm{\nu}$, with boundary data in some $\mathscr{P} (L; U)$. By restricting to the gap event $\textbf{GAP}$, we will show that this implies the paths in $\bm{\mathsf{x}}$ are approximated by the inverted height function $G = G^{\bm{\nu}}$ associated with $\bm{\nu}$; we will then use \Cref{p:limitprofile} to show that the edge behavior of the latter behaves as the function appearing in \eqref{e:globallawuse}. 
			
			To implement this, set $\widetilde{B} = 192 B^3$; denote the integers $n' = \lceil L^{1/2^{4000}} k \rceil$ and $n'' = \lceil L^{2^{1/5000}} \rceil$. Also abbreviate the events $\textbf{BTR}_n = \textbf{BTR}_n (4, B; k, L)$, $\textbf{CTR}_n = \textbf{CTR}_n (4, \widetilde{B}; k, L)$, and $\textbf{GAP}_n = \textbf{GAP}_n \big( [-4k^{1/3}, 4k^{1/3}]; R \big)$; further abbreviate the $\sigma$-algebra $\mathcal{F}_{\ext} = \mathcal{F}_{\ext}^{\bm{\mathsf{x}}} \big( \llbracket 1, n' \rrbracket \times (-2k^{1/3}, 2k^{1/3}) \big)$. Analogously to \eqref{a1a}, define the $\mathcal{F}_{\ext}$-measurable event 
			\begin{flalign}
				\label{a10}
				\mathscr{A}_1 = \textbf{BTR}_n \cap \bigcap_{j=1}^{n'} \textbf{LOC}_j \big( \{ -2k^{1/3}, 2k^{1/3} \}; -\widetilde{B} j^{2/3} - \widetilde{B} k^{2/3}; \widetilde{B} k^{2/3} - \widetilde{B} j^{2/3} \big).
			\end{flalign} 
		
			\noindent By \Cref{ftrbtr2}, we have $\textbf{CTR}_n \subseteq \mathscr{A}_1$, and so \Cref{ctr} (applied at $D = 1$) yields constants $c_1 = c_1 (B) > 0$ and $C_1 = C_1 (B) > 1$ such that $\mathbb{P} [\textbf{BTR}_n \setminus \mathscr{A}_1] \le C_1 e^{-c_1 (\log k)^2}$. Next, by the $A = 4$ case of \Cref{c:finalcouple} (after altering $c_1$ and $C_1$ if necessary), there exists an $\mathcal{F}_{\ext}$-measurable event $\mathscr{A}_2 \subseteq \textbf{BTR}_n$ satisfying $\mathbb{P} [\textbf{BTR}_n \setminus \mathscr{A}_2] \le C_1 e^{-c_1 (\log k)^2}$ and the following. Condition on $\mathcal{F}_{\ext}$ and restrict to $\mathscr{A}_2$. Denote the $n'$-tuples $\bm{u} = \bm{\mathsf{x}}_{\llbracket 1, n' \rrbracket} (-2k^{1/3}) \in \mathbb{W}_{n'}$ and $\bm{v} = \bm{\mathsf{x}}_{\llbracket 1, n' \rrbracket} (2k^{1/3}) \in \mathbb{W}_{n'}$, and sample a family of $n'$ non-intersecting Brownian bridges $\bm{\mathsf{y}} = (\mathsf{y}_1, \mathsf{y}_2, \ldots , \mathsf{y}_{n'}) \in \llbracket 1, n' \rrbracket \times \mathcal{C} \big( [-2k^{1/3}, 2k^{1/3}] \big)$. There exist two couplings between $\bm{\mathsf{x}}$ and $\bm{\mathsf{y}}$ such that, under the first we have 
			\begin{flalign}
				\label{xy1} 
				\mathbb{P} \Bigg[ \bigcap_{j = 1}^{n''} \bigcap_{|t| \le 2k^{1/3}} \big\{ \mathsf{x}_j (t) \le \mathsf{y}_j (t) - L^{-1/2^{5000}} k^{2/3}\big\}  \Bigg] \ge 1 - C_1 e^{-c_1 (\log k)^2}, 
			\end{flalign}
		
			\noindent and under the second we almost surely have 
			\begin{flalign}
				\label{xy2} 
				\mathsf{x}_j (t) \ge \mathsf{y}_j (t), \qquad \text{for each $(j, t) \in \llbracket 1, n' \rrbracket \times [-2k^{1/3}, 2k^{1/3}]$}.
			\end{flalign}
			
			Denote the $\mathcal{F}_{\ext}$-measurable event $\mathscr{A} = \mathscr{A}_1 \cap \mathscr{A}_2 \cap \textbf{GAP}_n \big( \{ -2k^{1/3}, 2k^{1/3} \}; R \big)$. In view of the inclusion $\textbf{GAP}_n  \subseteq \textbf{GAP}_n \big( \{ -2k^{1/3}, 2k^{1/3} \}; R \big)$, \eqref{bg2} and a union bound together indicate for sufficiently large $k$ that 
			\begin{flalign}
				\label{a32} 
				\begin{aligned}
				\mathbb{P} [\mathscr{A}] & \ge \mathbb{P} [ \textbf{BTR}_n \cap \textbf{GAP}_n] - \mathbb{P} [\textbf{BTR}_n \setminus \mathscr{A}_1] - \mathbb{P} [\textbf{BTR}_n \setminus \mathscr{A}_2] \\
				& \ge 1 - \varpi - 2C_1 e^{-c_1 (\log k)^2} \ge 1 - \displaystyle\frac{3 \varpi}{2}.
				\end{aligned} 
			\end{flalign}
			
			For the remainder of this proof, we condition on $\mathcal{F}_{\ext}$ and restrict to the event $\mathscr{A}$. By \eqref{a32}, the fact that $\mathbb{P} [\textbf{GAP}_n] \ge 1 - \varpi$ (by \eqref{bg2}), and a union bound, it then suffices to show for some constant $C = C(B, R) > 1$ and sufficiently large $k$ that there exist real numbers\footnote{Observe that this would allow us to take $(\mathfrak{a}, \mathfrak{b}, \mathfrak{c})$, satisfying the hypotheses of the proposition, to be measurable in $\bm{\mathsf{x}}$. Indeed, the function $\inf_{(\mathfrak{a}, \mathfrak{b}) \in [-C, C]} \inf_{\mathfrak{c} \in [1/C, C]} \big| \mathsf{x}_j (tk^{1/3}) - k^{2/3} (\mathfrak{a} + \mathfrak{b} t - \mathfrak{c} t^2) + (3 \pi / 4 \mathfrak{c}^{1/2})^{2/3} j^{2/3} \big|$ is continuous in $\bm{\mathsf{x}}$, and $\mathscr{E}$ defines the event on which this function is at most equal to $C \theta^3 k^{2/3}$. By the measurable selection theorem, we can take $(\mathfrak{a}, \mathfrak{b}, \mathfrak{c})$ as an argument infimum of the above continuous function (see \cite[Theorem 18.13]{IDA}) to be measurable in $\bm{\mathsf{x}}$.}  $\mathfrak{a}, \mathfrak{b} \in [-C, C]$ and $\mathfrak{c} \in [C^{-1}, C]$ such that 
			\begin{flalign}
				\label{ge} 
				\mathbb{P} [\textbf{GAP}_n \cap \mathscr{E}^{\complement}] \le \displaystyle\frac{\varpi}{2}, 
			\end{flalign} 
		
			\noindent where we have defined the event  
			\begin{flalign}
				\label{xeabc} 
				\mathscr{E} = \bigcap_{j =1}^{\lfloor \theta^3 k \rfloor} \bigcap_{|t| \le \theta} \Bigg\{ \bigg| \mathsf{x}_j (tk^{1/3}) - k^{2/3} \cdot (\mathfrak{a} + \mathfrak{b} t - \mathfrak{c} t^2) + \Big( \displaystyle\frac{3 \pi}{4 \mathfrak{c}^{1/2}} \Big)^{2/3} j^{2/3} \bigg|  \le C \theta^3  k^{2/3} \Bigg\}.
			\end{flalign}
			
			To do this, set $U = 16 \widetilde{B}^3 + 12R$; let $L' = (k^{-1} n')^{2/3} \ge L^{1/2^{4500}}$, so that $n' = (L')^{3/2} k$; and define the measure-valued process $\bm{\mu} = (\mu_t)_{t \in [0, 1]} \in \mathcal{C} \big( [0, 1]; \mathscr{P}_{\fin} \big)$ by setting 
			\begin{flalign}
				\label{muty} 
				\mu_t = (L')^{3/2} \cdot \emp \Big( (2k^{2/3})^{-1} \cdot \bm{\mathsf{y}} \big( 2(2t-1) k^{1/3}) \big) \Big), \qquad \text{for each $t \in [0, 1]$}.
			\end{flalign}
				
			\noindent Since we have restricted to $\mathscr{A} \subseteq \textbf{GAP}_n \big( \{-2k^{1/3}, 2k^{1/3} \}; R \big) \cap \mathscr{A}_1$, we have by \Cref{gap} (with the fact that $(\log n)^{25} \le (\log n')^{30}$ for sufficiently large $k$, as $n = L^{3/2} k \le L^{3/2} n'$) and \eqref{a10} that $\bm{u}$ and $\bm{v}$ satisfy \eqref{bbaa}, with the $B$ there equal to $\widetilde{B} + R$ here.  
			
			Thus, \Cref{pluconverge} applies, with the $(\bm{\mathsf{x}}; n; B)$ there equal to $(\bm{\mathsf{y}}; n'; \widetilde{B} + R)$ here. Setting $U = 4(\widetilde{B} + R)^3$, it yields measures $\nu_0, \nu_1 \in \mathscr{P} (L'; U)$ such that $\mathbb{P} \big[ d_{\dL} (\bm{\mu}, \bm{\nu}) \le \theta^6/2 \big] \ge 1 - \varpi / 8$, where $\bm{\nu} \in \mathcal{C} \big( [0, 1]; \mathscr{P}_{\fin} \big)$ is the bridge-limiting measure process on $[0, 1]$ with boundary data $(\mu_0; \mu_1)$. Denote the inverted height function associated with $\bm{\nu}$ (recall \Cref{hrhot}) by $G = G^{\bm{\nu}} : \big[ 0, (L')^{3/2} \big] \rightarrow \mathbb{R}$. By \eqref{muty}, the definitions \eqref{munudistance1} of $d_{\dL}$ and \eqref{gty} of $G$, the bound $\mathbb{P} \big[ d_{\dL} (\bm{\mu}, \bm{\nu}) \le \theta^6 / 2 \big] \ge 1 - \varpi / 8$ is equivalent to 
			\begin{flalign*}
				\mathbb{P} \Bigg[ \bigcap_{j=1}^{n'} \bigcap_{|t| \le 2} \bigg\{ k^{-2/3} \cdot \mathsf{y}_{j + \lfloor \theta^6 k \rfloor} (tk^{1/3})  - \theta^6 \le 2G \Big( \displaystyle\frac{t+2}{4}, \displaystyle\frac{j}{k} \Big) \le k^{-2/3} \cdot \mathsf{y}_{j - \lfloor \theta^6 k \rfloor} (tk^{1/3} & ) + \theta^6 \bigg\} \Bigg] \\
				& \ge 1 - \displaystyle\frac{\varpi}{8},
			\end{flalign*}
		
			\noindent where we have denoted $\mathsf{y}_j = \infty$ if $j < 1$ and $\mathsf{y}_j = -\infty$ if $j > n'$. Together with the couplings \eqref{xy1} and \eqref{xy2}, and a union bound, it follows for $k$ sufficiently large (so that $C_1 e^{-c_1 (\log k)^2} \le \varpi / 4$) that $\mathbb{P} [\mathscr{E}_0] \ge 1 - \varpi / 2$, where we have defined the event 
			\begin{flalign*}
				\mathscr{E}_0 =\bigcap_{j=1}^{n'' - \lfloor \theta^6 k \rfloor} \bigcap_{|t| \le 2} \bigg\{ k^{-2/3} \cdot \mathsf{x}_{j + \lfloor \theta^6 k \rfloor} (&  tk^{1/3}) - \theta^6 - L^{-1/2^{5000}} \\
				& \le 2 G \Big( \displaystyle\frac{t+2}{4}, \displaystyle\frac{j}{k} \Big) \le k^{-2/3} \cdot \mathsf{x}_{j - \lfloor \theta^6 k \rfloor} (tk^{1/3}) + \theta^6 \bigg\}.
			\end{flalign*}
			
			\noindent Thus, to show \eqref{ge}, it suffices to show that there exists a constant $C = C(B, R) > 1$ and real numbers $\mathfrak{a}, \mathfrak{b} \in \mathbb{R}$ and $\mathfrak{c} \in [C^{-1}, C]$ such that 
			\begin{flalign}
				\label{eeg}
				 \mathscr{E}_0 \cap \textbf{GAP}_n \subseteq \mathscr{E}.
			\end{flalign} 
		
			To this end, restrict to the event $\mathscr{E}_0 \cap \textbf{GAP}_n$. Then, for any $(j, t) \in \llbracket 1, n'' - \theta^6 k \rrbracket \times [-2, 2]$, 
			\begin{flalign}
				\label{xjg1} 
				k^{-2/3}\cdot \mathsf{x}_j (tk^{1/3}) \geq k^{-2/3}\cdot   \mathsf{x}_{j + \lfloor \theta^6 k \rfloor} (tk^{1/3}) \ge 2 G \Big( \displaystyle\frac{t+2}{4}, \displaystyle\frac{j}{k} + 2 \theta^6 \Big) - \theta^6,
			\end{flalign}
		
			\noindent where in the first bound we used the fact that $\mathsf{x}_j \ge \mathsf{x}_{j'}$ whenever $j \le j'$, and in the second we used the fact that we are restricting to $\mathscr{E}_0$. Similarly, for any $(j, t) \in \llbracket 1, n'' - \theta^6 k \rrbracket \times [-2, 2]$, we have for sufficiently large $k$ that
			\begin{flalign}
				\label{xjg2} 
				\begin{aligned}
				k^{-2/3} \cdot \mathsf{x}_j (tk^{1/3}) & \le k^{-2/3} \cdot \mathsf{x}_{j + \lfloor \theta^6 k \rfloor} (tk^{1/3}) + k^{-2/3} \big( R( j + \theta^6 k)^{2/3} - R j^{2/3}  + (\log n)^{25} \big) \\ 
				& \le 2 G \Big( \displaystyle\frac{t+2}{4}, \displaystyle\frac{j}{k} \Big)  + \theta^6  + L^{-1/2^{5000}} + (R+1) \theta^4 \le 2 G \Big( \displaystyle\frac{t+2}{4}, \displaystyle\frac{j}{k} \Big) + (R+3) \theta^4.
				\end{aligned} 
			\end{flalign}
		
			\noindent Here, in the first bound we used the fact that we have restricted to $\textbf{GAP}_n$; in the second, we used the fact that we have restricted to $\mathscr{E}_0$, as well as the bounds $(j + \theta^6 k)^{2/3} - j^{2/3} \le (\theta^6 k)^{2/3} = \theta^4 k^{2/3}$ and $(\log n)^{25} \le \theta^4 k^{2/3}$ for sufficiently large $k$ (as $n' \le n \le L^{3/2} k$); and in the third we used the fact that $L^{-1/2^{5000}} \le \theta^4$ (as $L \ge \theta^{-2^{6000}}$). Together with the fact that $n'' - \theta^6 k \ge L^{1/2^{5000}} k - \theta^6 k \ge \theta^3 k$ (as $L^{1/2^{5000}} \ge \theta^{-2^{1000}}$), we find from \eqref{xjg1} and \eqref{xjg2} that, for sufficiently large $k$,
			\begin{flalign}
				\label{xg3}
				 2 G \Big( \displaystyle\frac{t+2}{4}, \displaystyle\frac{j}{k} +2\theta^6 \Big) - \theta^6 \le k^{-2/3} \cdot \mathsf{x}_j (tk^{1/3}) \le  2 G \Big( \displaystyle\frac{t+2}{4}, \displaystyle\frac{j}{k} \Big) + (R+3) \theta^4,
			\end{flalign} 
		
			\noindent for any $(j, t) \in \llbracket 1, \theta^3 k \rrbracket \times [-2, 2]$.
			
			Now, observe that $\bm{\nu}$ is the bridge-limiting measure associated with boundary data $\nu_0, \nu_1 \in \mathscr{P} (L; U)$. We will assume in what follows that $G(0, 0) = 0 = G(1, 0)$ (as we may otherwise apply an affine transformation to $G (t, y)$, using the second part of \Cref{invariancesscale} to replace it by $G(t, y) - (1-t) \cdot G(0, 0) - t \cdot G(1, 0)$; such an affine transformation will only affect the constants $\mathfrak{a}$ and $\mathfrak{b}$ that appear below). Then, by \Cref{plu}, $\nu_0$ and $\nu_1$ satisfy \Cref{integralmu0mu1} and \Cref{gapmu0mu1} (with the $B$ there equal to $U$ here). Thus, \Cref{p:limitprofile} applies and yields a constant $C_2 = C_2 (B, R) > 1$ and real numbers $\mathfrak{a}_0, \mathfrak{b}_0 \in [-C_2, C_2]$ and $\mathfrak{c} \in [C_2^{-1}, C_2]$ such that for $L > C_2$ and $\theta < C_2^{-1}$ we have
			\begin{flalign*} 
				\displaystyle\sup_{|y| \le 2 \theta^3} \displaystyle\sup_{|s| \le \theta} \Bigg| G \Big( \displaystyle\frac{1}{2} + s, y \Big) - (\mathfrak{a}_0 + \mathfrak{b}_0 s - \mathfrak{c}_0 s^2) + \bigg( \displaystyle\frac{3\pi}{4 \mathfrak{c}_0^{1/2} } \bigg)^{1/2} y^{2/3} \Bigg| \le C_2 \theta^3.
			\end{flalign*} 
		
			\noindent Setting $s = t / 4$, it follows for $(\mathfrak{a}, \mathfrak{b}) = (2\mathfrak{a}_0, \mathfrak{b}_0 / 2)$ and $\mathfrak{c} =  \mathfrak{c}_0 / 8 \in \big[ (8C_2)^{-1}, C_2 \big]$ that 
			\begin{flalign}
				\label{gt} 
				\displaystyle\sup_{|y| \le 2 \theta^3} \displaystyle\sup_{|t| \le \theta} \bigg| 2 G \Big( \displaystyle\frac{1}{2} + \displaystyle\frac{t}{4}, y \Big) - (\mathfrak{a} + \mathfrak{b} t - \mathfrak{c} t^2) + \Big( \displaystyle\frac{3\pi}{4 \mathfrak{c}^{1/2}} \Big)^{2/3} y^{2/3} \bigg| \le 2 C_2 \theta^3.
			\end{flalign}
		
			\noindent In particular,   
			\begin{flalign*}
				\displaystyle\sup_{|y| \le 2 \theta^3 - 2\theta^6} \displaystyle\sup_{|t| \le \theta} & \bigg| 2G \Big( \displaystyle\frac{1}{2} + \displaystyle\frac{t}{4}, y \Big) - 2 G \Big( \displaystyle\frac{1}{2} + \displaystyle\frac{t}{4}, y + 2 \theta^6 \Big)  \bigg| \\ 
				& \le 4 C_2 \theta^3 + \Big( \displaystyle\frac{3 \pi}{4 \mathfrak{c}^{1/2}} \Big)^{2/3} \big( (y+\theta^6)^{2/3} - y^{2/3} \big) \le 4 C_2 \theta^3 + 4 C_2^{1/3} \theta^4 \le 8 C_2 \theta^3,
			\end{flalign*}
		
			\noindent where in the first inequality we applied \eqref{gt} twice, and in the second we used the facts that $\mathfrak{c} = \mathfrak{c}_0 / 8 \ge (8C_2)^{-1}$ and that $(y + \theta^6)^{2/3} - y^{2/3} \le \theta^4$. Inserting this with \eqref{gt} into \eqref{xg3} yields 
			\begin{flalign*} 
				\displaystyle\max_{j \in \llbracket 1, \theta^3 k \rrbracket} \displaystyle\sup_{|t| \le \theta} \bigg| k^{-2/3} \cdot \mathsf{x}_j (tk^{1/3}) - (\mathfrak{a} + \mathfrak{b} t - \mathfrak{c} t^2) + \Big( \displaystyle\frac{3\pi}{4 \mathfrak{c}^{1/2}} \Big)^{2/3} \Big( \displaystyle\frac{j}{k} \Big)^{2/3} \bigg| \le (8 C_2 + R + 3) \theta^3.
			\end{flalign*} 
		
			\noindent This verifies that $\mathscr{E}$ holds (by its definition \eqref{xeabc}) at $C =   8C_2 + R + 3$, thereby confirming \eqref{eeg} and establishing the proposition.
		\end{proof}

\chapter{Curvature Approximation} 

\label{APPROXIMATECURVE}

\section{Second Derivative Approximations for Paths} 

\label{DerivativePath2}

\subsection{Proof of \Cref{h0x2}} 

\label{Proofhj}

In this section we establish \Cref{h0x2}. We will frequently make use of the following quick lemma on the function $G$ from \eqref{g2121} (recalling \Cref{functionadmissible}). 

\begin{lem}
	
	\label{g00}
	
	Adopting \Cref{derivativegxpfl}, we have
	\begin{flalign}
		\label{g212} 
		G \in \Adm_{1/2} (\mathfrak{R}), \qquad \text{and} \qquad \text{$G$ satisfies \eqref{equationxtd} on $\mathfrak{R}$}.
	\end{flalign}
	
\end{lem}

\begin{proof}
	
	For any $(t, x) \in \mathfrak{R} = [-\xi, \xi] \times [0, 1]$, we have $-\partial_x G(t,x) = 2^{-1/6} 3^{-1/3} \pi^{2/3} (x+1)^{-1/3} \in [1/2,2]$; this confirms the first statement in \eqref{g212}. That $G$ satisfies \eqref{equationxtd} on $\mathfrak{R}$ follows either by direct verification, or by combining \Cref{airy} and \Cref{gequation}.
\end{proof} 

\begin{figure}
	\center
	\includegraphics[width=1\textwidth]{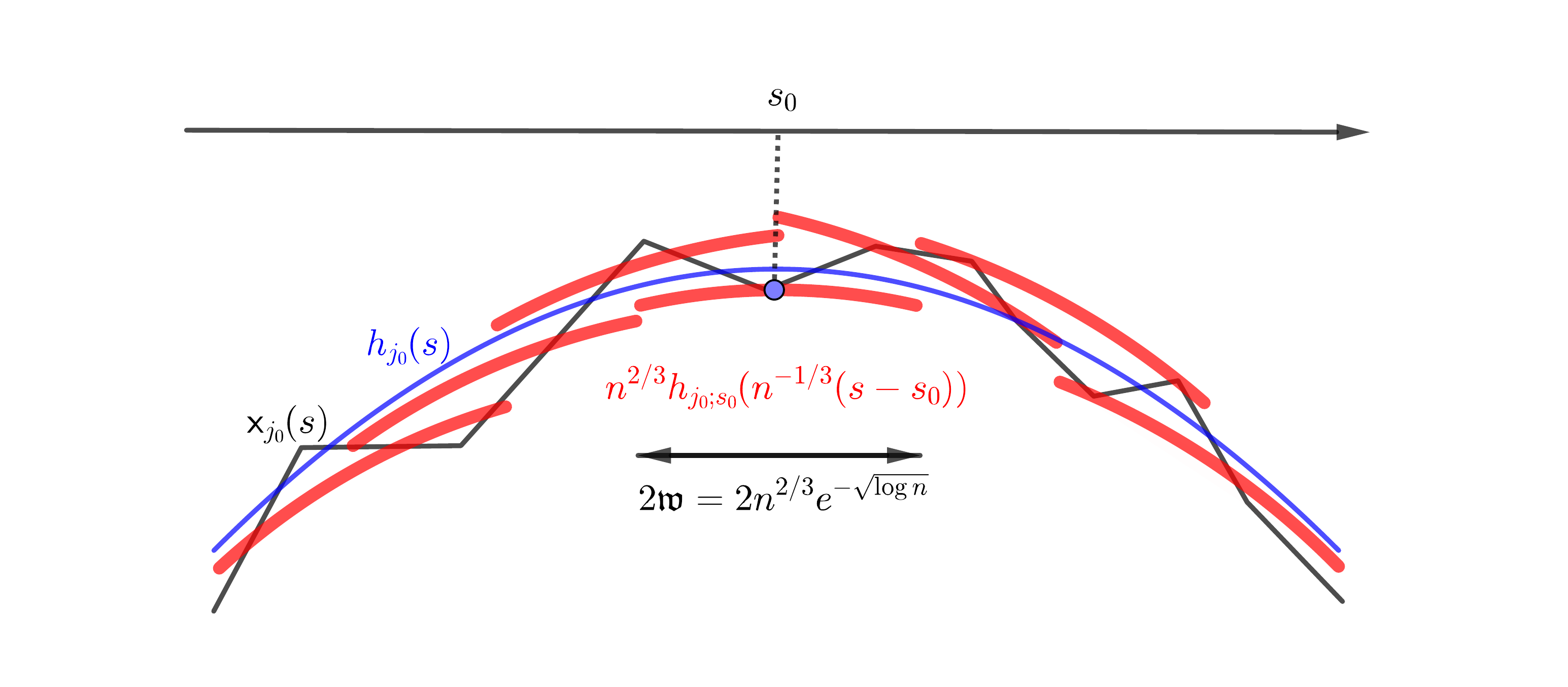}
	
	\caption{Shown above in red are the local approximations $\mathfrak{h}_{j; k}$ to the black curve $x_j$, which ``glue'' together to form the blue global approximation $h_j$.}
	
	\label{f:local_glue}
\end{figure}

To establish \Cref{h0x2} we will first ``locally'' produce the functions $h_j$, on time scales of length $2 e^{-\sqrt{\log n}}$, that satisfy the required properties (up to the rescaling \eqref{xjsn} in \eqref{xj0s}, transforming $\mathsf{x}_j$ into $x_j$); we will then ``glue'' these local functions together to form a global one; see \Cref{f:local_glue} for a depiction. The following proposition implements the first task; its proof is given in \Cref{LOCALHJ} below.

\begin{prop} 
	
	\label{h0xlocal} 
	
	Adopting the notation and assumptions of \Cref{h0x2}, there exist constants $c = c(B) > 0$ and $C = C(B) > 1$ such that the following holds whenever $\delta \in (0, c)$ and $n > C$. Denote $\eta = e^{-\sqrt{\log n}}$; fix a real number $s_0 \in [-\xi/2,  \xi/2]$; and fix an integer $j_0 \in \llbracket n/3, 2n/3 \rrbracket$. Then, with probability at least $1 - n^{-15}$, there exists a (random) twice-differentiable function $h_{j_0; s_0} : [- \eta, \eta] \rightarrow \mathbb{R}$ with 
	\begin{flalign}
		\label{bhs2}
		\displaystyle\sup_{|s| \le \eta} \big| \partial_s^2 h_{j_0; s_0} (s) + 2^{1/2} \big| \le \delta^{1/6} + (\log n)^{-1/3}, \qquad \text{and} \qquad \| h_{j_0; s_0} \|_{\mathcal{C}^1} \le 10 B,
	\end{flalign} 
	
	\noindent such that
	\begin{flalign}
		\label{xj0s} 
		\displaystyle\sup_{|s-s_0| \le \eta} \big| x_{j_0} (s) - h_{j_0; s_0} (s-s_0) \big| < n^{-13/15}. 
	\end{flalign}
	
\end{prop}

\begin{proof}[Proof of \Cref{h0x2}]
	
	Set $K = \lfloor \xi \eta^{-1} \rfloor$ and $s_k = k \eta - \xi/2$, for each $k \in \llbracket 0, K+1 \rrbracket$. Observe that $K \le B \eta^{-1} \le n^{1/10}$ for sufficiently large $n$, as $\xi \le B$ and $\eta = e^{-\sqrt{\log n}}$. Further denote 
	\begin{flalign}
		\label{thetadeltab} 
		\varpi = n^{-13/15}; \qquad \theta = \delta^{1/6} + (\log n)^{-1/3}; \qquad \mathfrak{B} = 10B.
	\end{flalign}
	
	By \Cref{h0xlocal} and a union bound, there exist constants $c = c(B) > 0$ and $C = C(B) > 1$ such that, if $\delta < c$, the following holds with probability at least $1 - Cn^{-10}$. For each integer $j \in \llbracket n/3, 2n/3 \rrbracket$ and $k \in \llbracket 0, K-1 \rrbracket$, there exists a twice-differentiable (random) function $\mathfrak{h}_{j; k} : [s_k, s_{k+2}] \rightarrow \mathbb{R}$ with 
	\begin{flalign}
		\label{hks} 
		\displaystyle\sup_{s \in [s_k, s_{k+2}]} \big| \mathfrak{h}_{j;k}'' (s) + 2^{1/2} \big| < \theta; \qquad \text{and} \qquad \| \mathfrak{h}_{j;k} \|_{\mathcal{C}^1} < \mathfrak{B},
	\end{flalign}
	
	\noindent such that 
	\begin{flalign}
		\label{hks2}
		\displaystyle\sup_{|s - s_k| \le \mathfrak{a}/K} \big| x_j (s) - \mathfrak{h}_{j;k} (s) \big| < \varpi.
	\end{flalign}
	
	Now restrict to the event the above holds, and fix an index $j \in \llbracket n/3, 2n/3 \rrbracket$; we will define a function $h_j : [-\xi/2, \xi/2] \rightarrow \mathbb{R}$ satisfying \eqref{hjxj2} and \eqref{hjxj1}. To that end, for each integer $k \in \llbracket 0, K \rrbracket$, define the intervals $\mathcal{I}_k^-$, $\mathcal{I}_k$, and $\mathcal{I}_k^+$ by setting
	\begin{flalign*}
		& \mathcal{I}_k^- = ( s_k, s_k + \eta / 3 ]; \qquad \mathcal{I}_k = ( s_k + \eta / 3, s_k + 2 \eta / 3]; \qquad \mathcal{I}_k^+ = ( s_k + 2 \eta / 3, s_{k+1}].
	\end{flalign*} 
	
	\noindent Denote $\mathcal{J}_k = \mathcal{I}_k^-  \cup \mathcal{I}_k \cup \mathcal{I}_k^+ = [s_k, s_{k+1}]$. Fix a twice-differentiable function $\psi : [0, 1] \rightarrow [0, 1]$ with  
	\begin{flalign*}
		\psi (s) = 1, \quad \text{for $s \in [ 0, 1/3]$}; \qquad \psi (s) = 0, \quad \text{for $s \in [ 2 / 3, 1]$}; \qquad \| \psi \|_{\mathcal{C}^1} \le 20; \qquad \| \psi \|_{\mathcal{C}^2} \le 200.	
	\end{flalign*} 
	
	\noindent For each integer $k \in \llbracket 0, K \rrbracket$, define $\psi_k : \mathcal{J}_k \rightarrow \mathbb{R}$ by for each $s \in \mathcal{J}_k$ setting 
	\begin{flalign}
		\label{psik} 
		\psi_k (s) = \psi \big( \eta^{-1} (s-s_k) \big), \qquad \text{so that} \quad [\psi_k]_1 \le 20\eta^{-1}, \quad \text{and} \quad [\psi_k]_2 \le 200 \eta^{-2}.
	\end{flalign}
	
	\noindent Then, define $h_j: [-\xi/2, \xi/2] \rightarrow \mathbb{R}$ by for each $s \in [-\xi/2, \xi/2]$ setting 
	\begin{flalign}
		\label{hs0}
		h_j(s) = \displaystyle\sum_{k=0}^{K} \textbf{1}_{s \in \mathcal{J}_k} \cdot \Big( \psi_k (s) \cdot \mathfrak{h}_{j;k-1} (s) + \big( 1 - \psi_k (s) \big) \cdot \mathfrak{h}_{j;k} (s) \Big),
	\end{flalign}	
	
	\noindent where we have set $\mathfrak{h}_{j;-1} = \mathfrak{h}_{j;0}$ and $\mathfrak{h}_K = \mathfrak{h}_{j;K-1}$. Observe in this way that $h_j(s) = \mathfrak{h}_{j;k} (s)$ for each $s$ in a neighborhood of $s_{k+1}$, since $\psi (0^+) = 1$ and $\psi (1^-) = 0$. By the facts that $0 \le \psi_k \le 1$ (as $0 \le \psi \le 1$) and that the intervals $\{ \mathcal{J}_k \}$ are disjoint, we then have by \eqref{hks2} that
	\begin{flalign}
		\label{hjxjs2}
		\displaystyle\sup_{|s| \le \xi/2} \big| h_j(s) - x_j (s) \big| \le \displaystyle\max_{k \in \llbracket 0, K-1 \rrbracket}  \displaystyle\sup_{s \in \mathcal{J}_k}  \max \Big\{ \big| \mathfrak{h}_{j;k-1} (s) - x_j (s) \big|, \big| \mathfrak{h}_{j;k} (s) - x_j (s) \big| \Big\}  < \varpi.
	\end{flalign}
	
	\noindent Thus, 
	\begin{flalign*}
		\displaystyle\sup_{|s| \le \mathsf{T}/2} \big| \mathsf{x}_j (s) - n^{2/3} \cdot h_j (sn^{-1/3}) \big| = n^{2/3} \cdot \displaystyle\sup_{|s| \le \xi/2} \big| x_j (s) - h_j (s) \big| < \varpi n^{2/3} = n^{-1/5},
	\end{flalign*} 

	\noindent where in the first statement we used \eqref{xjsn} and the fact that $\mathsf{T} = \xi n^{1/3}$; in the second we used \eqref{hjxjs2}; and in the third we used the definition \eqref{thetadeltab} of $\varpi$. This confirms \eqref{hjxj1}; it remains to establish \eqref{hjxj2}.
	
	Moreover, for any $s \in [-\xi/2, \xi/2]$, letting $k = k (s) \in \llbracket 0, K \rrbracket$ be the unique integer such that $s \in \mathcal{J}_k$, we have  by \eqref{hs0} that
	\begin{flalign}
		\label{hs2}
		\begin{aligned}
			&  h_j' (s)  =  \big( \mathfrak{h}_{j;k-1} (s) - \mathfrak{h}_{j;k} (s) \big) \cdot \psi_k' (s) + \big( \mathfrak{h}_{j;k-1}' (s) - \mathfrak{h}_{j;k}' (s) \big)\cdot \psi_k (s) + \mathfrak{h}_{j;k}' (s); \\
			&   h_j'' (s)  =  \big( \mathfrak{h}_{j;k-1} (s) - \mathfrak{h}_{j;k} (s) \big) \cdot \psi_k'' (s) + 2 \big( \mathfrak{h}_{j;k-1}' (s) - \mathfrak{h}_{j;k}' (s) \big)\cdot \mathfrak{\psi}_k' (s) \\
			& \qquad \qquad \qquad + \big( \mathfrak{h}_{j;k-1}'' (s) - \mathfrak{h}_{j;k}'' (s) \big) \cdot \psi_k (s) + \mathfrak{h}_{j;k}'' (s).
		\end{aligned} 
	\end{flalign}
	
	To bound the right sides of \eqref{hs2}, first observe by \eqref{hks} that
	\begin{flalign}
		\label{hkhk1s2}
		\displaystyle\sup_{s \in \mathcal{J}_k} \big| \mathfrak{h}_{j;k-1} (s) - \mathfrak{h}_{j;k} (s) \big| \le \displaystyle\sup_{s \in \mathcal{J}_k} \Big( \big| \mathfrak{h}_{j;k-1} (s) - x_j (s) \big| + \big| \mathfrak{h}_{j;k} (s) - x_j (s) \big| \Big) < 2 \varpi.
	\end{flalign}
	
	\noindent Applying \eqref{hkhk1s2} at $s = s_k$ and $s = s_{k+1}$, and using the continuity of $\mathfrak{h}_{j;k-1}' - \mathfrak{h}_{j;k}'$, we find that there exists an $S \in \mathcal{J}_k$ such that 
	\begin{flalign} 
		\label{hk1hk} 
		\big| \mathfrak{h}_{j;k-1}' (S) - \mathfrak{h}_{j;k}' (S) \big| < 4 \eta^{-1} \varpi.
	\end{flalign}
	
	\noindent Again by \eqref{hks}, we have 
	\begin{flalign}
		\label{hk1hk22} 
		\displaystyle\sup_{s \in \mathcal{J}_k} \big| \mathfrak{h}_{j;k-1}'' (s) - \mathfrak{h}_{j;k}'' (s) \big| \le 2 \theta,
	\end{flalign}
	
	\noindent which with \eqref{hk1hk} yields    
	\begin{flalign}
		\label{hk1hk1}
		\displaystyle\sup_{s \in \mathcal{J}_k} \big| \mathfrak{h}_{j;k-1}' (s) - \mathfrak{h}_{j;k}' (s) \big| < 4 \varpi + 2 \eta \theta.
	\end{flalign}
	
	\noindent Inserting \eqref{hkhk1s2}, \eqref{hk1hk22}, \eqref{hk1hk1}, and \eqref{psik} (with the facts that $\| \mathfrak{h}_{j;k} \|_{\mathcal{C}^1} \le \mathfrak{B}$ and $\big| h_j'' (s) + 2^{1/2} \big| \le \theta$, by \eqref{hks}) into \eqref{hs2} yields 
	\begin{flalign*}
		& \| h_j\|_{\mathcal{C}^1}  < 40 \varpi \eta^{-1} + 4  \varpi + 2 \eta \theta + \mathfrak{B} \le 20B,
	\end{flalign*} 

	\noindent and 
	\begin{flalign*} 
		 \big| h_j'' (s) + 2^{1/2} \big| < 400 \varpi \eta^{-2} + 40 \eta^{-1} (4 \varpi + 2 \eta \theta) + 3 \theta & = 400 \varpi \eta^{-2} + 160 \varpi \eta^{-1} + 83 \theta \\
		& \le \delta^{1/8} + (\log n)^{-1/4},
	\end{flalign*}
	
	\noindent where in both of the last bounds we used the facts that $\eta = e^{-\sqrt{\log n}}$; the definition \eqref{thetadeltab} of $(\varpi, \theta, \mathfrak{B})$; and the fact that $n$ and $\delta$ are sufficiently large and small, respectively (and that $B>1$). This verifies \eqref{hjxj2} and thus the theorem. 
\end{proof}

\subsection{Perturbations of Boundary Data for Limit Shapes} 

\label{Proof0Boundary}

To establish \Cref{h0xlocal}, one might seek to apply \Cref{gh} to the restriction of $\bm{x}$ to a $2\eta \times 1$ rectangle centered at $(s_0, j_0 n^{-1})$. To do this, one must verify the assumptions of \Cref{fgr} indicating that the boundary data of $\bm{x}$ along this rectangle are sufficiently regular. That this holds for the starting and ending data would be a consequence of \eqref{pflx2} (indicating that the regular profile event $\textbf{PFL}^{\bm{x}}$ likely holds), but no such guarantee holds for the upper and lower boundaries. 

\begin{figure}
	\center
	\includegraphics[width=1\textwidth]{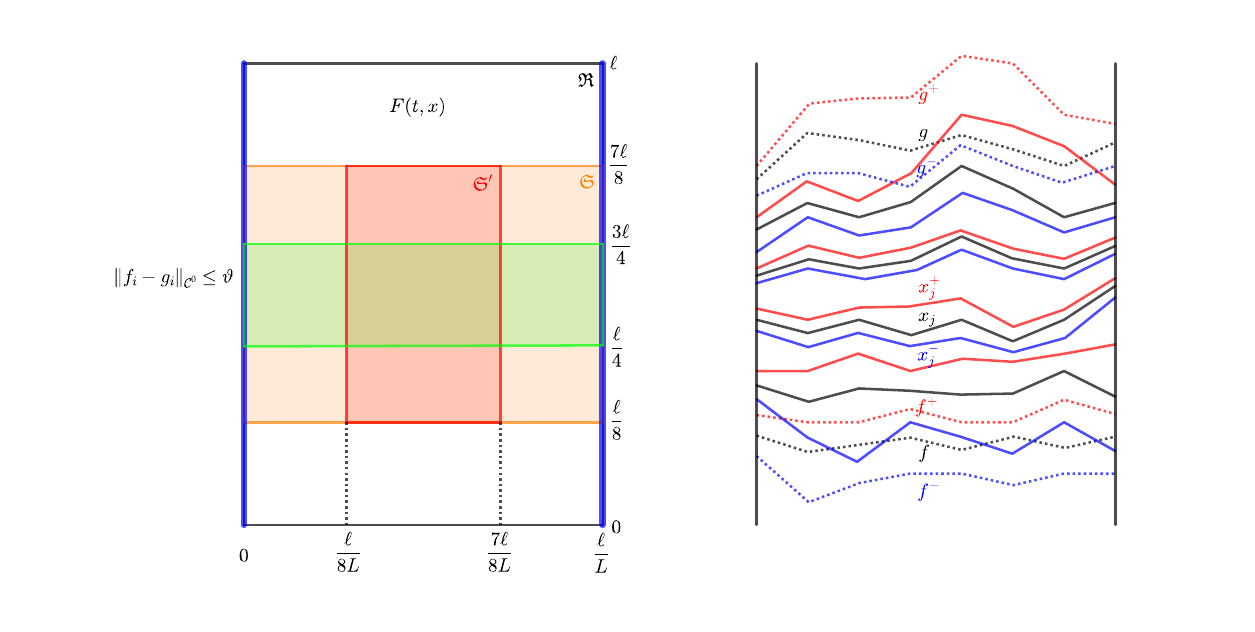}
	
	\caption{Shown to the left is a depiction for \Cref{derivativef}, indicating that there exist two inverted height functions $G^\pm$ on $\mathfrak{S}$ that are close to $F$, such that their difference is exponentially small on the green region. Shown to the right is a depiction for the sandwiching argument.
	}
	
	\label{f:sandwich}
\end{figure}

To circumvent this issue, we will instead introduce two families $\bm{x}^-$ and $\bm{x}^+$ of non-intersecting Brownian bridges and sandwich $\bm{x}$ between $\bm{x}^-$ and $\bm{x}^+$. These families will be defined so that their starting and ending data almost coincide with that of $\bm{x}$ near the middle of the rectangle. However, around its top and bottom, the starting and ending data of $\bm{x}$ will be higher than those of $\bm{x}^-$, and lower than those of $\bm{x}^+$; we will also make the upper and lower boundaries for $\bm{x}^-$ and $\bm{x}^+$ regular. Thus, \Cref{fgr} (and hence the concentration bound \Cref{gh}) will apply to $\bm{x}^-$ and $\bm{x}^+$, giving a bound on $\bm{x}$ due to the sandwiching; see the right side of \Cref{f:sandwich}. For this sandwiching to be effective, we must verify that it is possible to introduce these boundary perturbations in such a way that they do not substantially affect the model in the middle of the rectangle. 

In this section we state the below lemma, showing that this holds for the associated limit shape. Its proof largely follows from \Cref{equationcompareboundary} and \Cref{f1f2b}, and is provided in \Cref{ProofGFG} below. 

\begin{lem} 
	
	\label{fg1g2} 
	
	For any integer $m \ge 7$ and real numbers $\varepsilon > 0$ and $B > 1$, there exist constants $c = c(\varepsilon, B, m) \in (0, 1)$ and $C = C (\varepsilon, B, m) > 1$ such that the following holds. Let $L > 4$ and $\vartheta \in (0, c)$ be real numbers with $|\log \vartheta|^{20} < L < \vartheta^{-1/2m^2}$; also let $\ell \in (B^{-1}, B)$ be a real number. Define the open rectangles    
	\begin{flalign}
		\label{ssr}
		\mathfrak{R} = \Big( 0, \displaystyle\frac{\ell}{L} \Big) \times (0, \ell) ; \qquad \mathfrak{S} = \Big( 0, \displaystyle\frac{\ell}{L} \Big) \times \Big( \displaystyle\frac{\ell}{8}, \displaystyle\frac{7 \ell}{8} \Big); \qquad \mathfrak{S}' = \Big( \displaystyle\frac{\ell}{8L}, \displaystyle\frac{7 \ell}{8L} \Big) \times \Big( \displaystyle\frac{\ell}{8}, \displaystyle\frac{7 \ell}{8} \Big).
	\end{flalign}
	
	\noindent Let $F \in \Adm_{\varepsilon} (\mathfrak{R}) \cap \mathcal{C}^m (\overline{\mathfrak{R}})$ denote a solution to \eqref{equationxtd} such that $\| F \|_{\mathcal{C}^m (\mathfrak{R})} \le B$, and define the functions $f_0, f_1 : [0, \ell] \rightarrow \mathbb{R}$ by setting $f_i (x) = F (i\ell L^{-1}, x)$ for each $i \in \{ 0, 1 \}$ and $x \in [0, \ell]$. Further fix functions $g_0, g_1 : [0, \ell] \rightarrow \mathbb{R}$ such that $\| f_i - g_i \|_{\mathcal{C}^0} \le \vartheta$ and $\| g_i \|_{\mathcal{C}^m} \le B$ for each index $i \in \{ 0, 1 \}$. Then there exist solutions $G^-, G^+ \in \Adm_{\varepsilon/2} (\mathfrak{S}) \cap \mathcal{C}^m (\overline{\mathfrak{S}})$ to \eqref{equationxtd} on $\mathfrak{S}$ satisfying the following properties.
	
	\begin{enumerate}
		\item For each $i \in \{ 0, 1 \}$ and $x \in [ \ell / 5, 4 \ell / 5 ]$, we have $G^- (i \ell L^{-1}, x) = g_i (x) = G^+ (i \ell L^{-1}, x)$. 
		\item For each $i \in \{ 0, 1 \}$ and $x \in [ \ell / 8, 7\ell / 8 ]$, we have $G^- (i \ell L^{-1}, x) \le g_i (x) \le G^+ (i \ell L^{-1}, x)$. 
		\item We have $\| G^- \|_{\mathcal{C}^{m-5} (\mathfrak{S})} + \| G^+  \|_{\mathcal{C}^{m-5} (\mathfrak{S})} \le C$.
		\item We have $\| G^- - F \|_{\mathcal{C}^m (\mathfrak{S}')} + \| G^+ - F \|_{\mathcal{C}^m (\mathfrak{S}')} \le C \vartheta^{3/4}$.
		\item For each $(t, x) \in [0, \ell L^{-1}] \times [\ell / 4, 3\ell / 4]$, we have $\big| G^+ (t, x) - G^- (t, x) \big| \le C e^{-c L^{1/8}}$.
		\item For each $(t, x) \in [0, \ell L^{-1}] \times \{ \ell / 8, 7\ell / 8 \}$, we have $G^- (t, x) \le F(t,x) - \vartheta < F(t, x) + \vartheta \le G^+ (t, x)$. 
	\end{enumerate}
	
\end{lem} 

Let us briefly explain \Cref{fg1g2}; see the left side of \Cref{f:sandwich}. One may view $F$ as the ``original'' function and $G^-$ and $G^+$ as two perturbations of it that have different boundary data along the two vertical sides of $\mathfrak{R}$. The first part of the lemma indicates that the boundary data of $G^-$ and $G^+$ are both given by $g_i$ around the middles of these two sides; the second indicates that $G^-$ and $G^+$ are larger than smaller than $g_i$ around the endpoints of these sides, respectively. The third indicates that $G^-$ and $G^+$ are regular up to the boundary of $\mathfrak{S}$; the fourth indicates that $G^-$ and $G^+$ (and their derivatives) are close to the original function $F$ in the interior $\mathfrak{S}' \subset \mathfrak{S}$ of the rectangle. The fifth indicates that $G^+$ and $G^-$ are quite close in the middle of $\mathfrak{R}$, which will eventually make sandwiching between them effective. The sixth indicates that the boundary data of $G^-$ and $G^+$ along the two horizontal sides of $\mathfrak{S}$ are lower and higher than those of $F$, by at least $\vartheta$, respectively. 

\subsection{Proof of \Cref{h0xlocal}} 

\label{LOCALHJ}

In this section we establish \Cref{h0xlocal}; we adopt the notation of that proposition throughout. The content in \Cref{Proof0Boundary} presented several elements of its proof, but it simplified the discussion on the regularity for the starting and ending data of $\bm{x}$ (along the $2 \eta \times 1$ rectangle centered at $(s_0, j_0 n^{-1})$ described there). Although the likelihood \eqref{pflx2} of the regular profile event indeed indicates that the starting and ending data are each individually regular in the vertical direction, it does not directly forbid the possibility that these data are far from each other; this would make it impossible to find a regular profile (with uniformly bounded $t$-derivative) that interpolates between them in the sense of \eqref{uvg} in \Cref{fgr}. To preclude this possibility, we induct on scales, applying the discussion of \Cref{Proof0Boundary} on thinner rectangles, until eventually reaching width around $2 \eta$.

This requires some additional notation. Let $c_0 = c_0 (B) > 0$ and $C_0 = C_0 (B) > 1$ denote the constants $c ( 1 / 4, 2B, 50 ) > 0$ and $C (1/4, 2B, 50 ) > 1$ from \Cref{fg1g2}, respectively. For any integer $k \ge 0$, define the real numbers $\delta_0, \omega_k, \varsigma_k, \vartheta_k, \Theta_k > 0$ and $L_k > 1$, and integer $n_k \ge 1$, by setting 
\begin{flalign}
	\label{deltaomegal} 
	\begin{aligned}
		& \omega_k = 4^{-k-1}; \qquad n_k = \lfloor \omega_k n \rfloor; \qquad \delta_0 = \delta^{1/2} + ( \log n)^{-1}; \qquad L_k = \delta_0^{-\sqrt{k+1}};   \\
		& \varsigma_k = 5C_0 \exp \Big( -\displaystyle\frac{c_0}{2} \cdot L_k^{1/80} \Big); \qquad \vartheta_k = \varsigma_k + n^{-13/15}; \qquad \Theta_k = \delta_0^{3/4} + \displaystyle\sum_{j=0}^k \omega_j^{-1} \vartheta_j^{3/4}.
	\end{aligned}
\end{flalign}
\begin{figure}
	\begin{center}
    \includegraphics[width=1\textwidth, trim=2cm 6cm 2cm 6cm, clip]{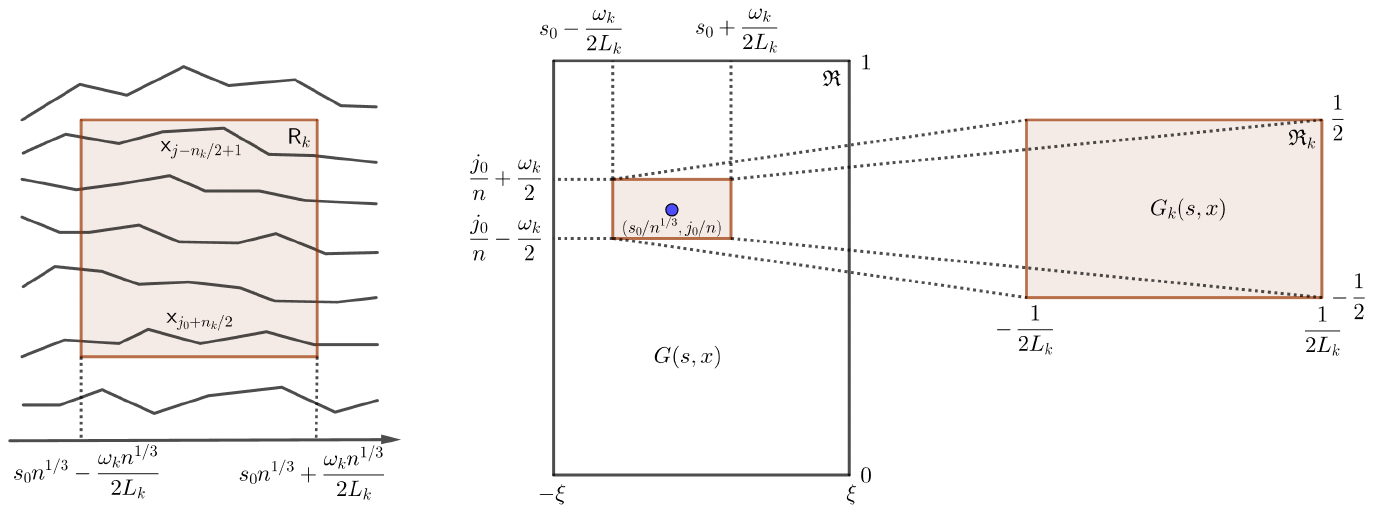}
\end{center}

	\caption{Shown to the left is the rectangle $\mathsf{R}_k$. Shown to the right is the rectangle $\mathfrak{R}_k$ obtained by ``zooming in'' around the point $(s_0, j_0/n)$.}
	\label{f:setting_2}
\end{figure}

\noindent Also let $K_0 \ge 1$ denote the maximal integer such that $\omega_{K_0} L_{K_0}^{-1} \ge 3 \eta$. To ease notation, we will omit the floors in what follows, assuming that each $\omega_k n$ is an integer; this will barely affect the proofs. For each integer $k \in \llbracket 0, K_0 \rrbracket$, define the set $\mathsf{R}_k \subset \mathbb{Z} \times \mathbb{R}$ and open rectangle $\mathfrak{R}_k \subset \mathbb{R}^2$ by
\begin{flalign}
	\label{rkrk} 
	\begin{aligned}
		& \mathsf{R}_k = \Big\llbracket j_0 - \displaystyle\frac{n_k}{2} + 1, j_0 + \displaystyle\frac{n_k}{2} \Big\rrbracket \times  \Big( s_0 n^{1/3} - \displaystyle\frac{\omega_k n^{1/3}}{2 L_k}, s_0 n^{1/3} + \displaystyle\frac{\omega_k n^{1/3}}{2 L_k} \Big); \\ 
		& \mathfrak{R}_k = \Big( -\displaystyle\frac{1}{2L_k}, \displaystyle\frac{1}{2L_k} \Big) \times \Big( -\displaystyle\frac{1}{2}, \displaystyle\frac{1}{2} \Big).
	\end{aligned} 
\end{flalign}		

\noindent See the left side of \Cref{f:setting_2} for a depiction. For each $k \in \llbracket 0, K_0 \rrbracket$, also define the function $G_k : \overline{\mathfrak{R}}_k \rightarrow \mathbb{R}$ by rescaling $G$, namely, by setting 
\begin{flalign} 
	\label{gksx2} 
	\begin{aligned}
		G_k (s, x) = \omega_k^{-1} \cdot G ( s_0 + \omega_k s, j_0 n^{-1} + \omega_k x), \qquad \text{for each $(s, x) \in \overline{\mathfrak{R}}_k$}. 
	\end{aligned}
\end{flalign}

\noindent See the right side of \Cref{f:setting_2} for a depiction. Observe that the $G_k$ satisfy \eqref{equationxtd} on $\mathfrak{R}_k$, by the first part (at $\alpha = \beta = \omega_k$) of \Cref{invariancesscale}. Then, for each integer $k \in \llbracket 0, K_0 \rrbracket$, define the sequence of functions $\bm{x}^{(k)} = \big( x_1^{(k)}, x_2^{(k)}, \ldots , x_{n_k}^{(k)} \big) \in \llbracket 1, n_k \rrbracket \times \mathcal{C} \big( [-1 / 2L_k, 1 / 2L_k] \big)$ by rescaling and reindexing $\bm{\mathsf{x}}$, namely, by setting 
\begin{flalign}
	\label{xks} 
	x_i^{(k)} (s) = \omega_k^{-1} n^{-2/3} \cdot \mathsf{x}_{i+j_0 - n_k/2} ( s \omega_k  n^{1/3} + s_0 n^{1/3}) = \omega_k^{-1} \cdot x_{i+j_0-n_k/2} (s \omega_k + s_0),
\end{flalign} 

\noindent for each $(i, s) \in \llbracket 0, n_k+1 \rrbracket \times [ -1/2L_k, 1/2L_k]$.

\begin{rem} 
	
	\label{delta0omegaklk}
	
	Let us explain the notation above. As mentioned at the beginning of this \Cref{LOCALHJ}, we will analyze the paths in $\bm{\mathsf{x}}$ on a sequence of rectangles of decreasing width, eventually reaching width $2 \eta n^{1/3}$. Transferring estimates on $\bm{\mathsf{x}}$ from one rectangle to the next will necessitate that the heights and widths of these rectangles do not decrease too quickly. In particular, we cannot directly pass from our initial rectangle (whose width is of order $n^{1/3}$) to one of width $2 \eta n^{1/3}$; this is why we study a sequence of rectangles, through the inductive procedure described below. 
	
	The ``height'' of (that is, the number of paths in) each rectangle is reduced by a constant factor $4$ when passing to the subsequent one. Thus, we consider $n_k$ paths of $\bm{\mathsf{x}}$ in the $k$-th rectangle, denoted by $\mathsf{R}_k$. Since the variance of each Brownian bridge in $\bm{\mathsf{x}}$ is equal to one, the scaling \eqref{xks} implies that the variance of each Brownian bridge in $\bm{x}^{(k)}$ is equal to $(\omega_k n)^{-1} = n_k^{-1}$; this is consistent with \Cref{fgr}, which will eventually enable us to apply the concentration bound \Cref{gh} to $\bm{x}^{(k)}$. Implementing the corresponding scaling for $\mathsf{R}_k$ yields its macroscopic variant  $\mathsf{R}_k$ from \eqref{rkrk}, which has aspect ratio $L_k$ from \eqref{deltaomegal}; applying it to $G$ yields the function $G_k$ from \eqref{gksx2}. 
	
	As discussed in the beginning of \Cref{Proof0Boundary} (see also \Cref{Derivative2Approximate}), to apply \Cref{gh} we must modify the upper and lower boundaries for $\bm{x}^{(k)}$ to make them smooth; this will produce two functions $G_k^+$ and $G_k^-$ (see \Cref{gk} below), which will approximately bound $\bm{x}^{(k)}$ from above and below, respectively (see the proof of \Cref{omega42} of \Cref{omega4}). Since the aspect ratio of $\mathfrak{R}_k$ is $L_k$, the fifth property in \Cref{fg1g2} will indicate that $G_k^+$ and $G_k^-$ will differ by at most the parameter $\varsigma_k$ in \eqref{deltaomegal} (see the fifth property in \Cref{gk}). The parameter $\vartheta_k$ incorporates an additional error arising from the concentration estimate \Cref{gh}; it will bound the distance from $\bm{x}^{(k)}$ to $G_k^{\pm}$ (see \eqref{omegak43} and \Cref{omegakestimate} below). It also (using the fourth property in \Cref{fg1g2}) estimates the difference between $G_k^{\pm}$ and (the correspondingly normalized) $G_{k-1}^{\pm}$ (see the fourth property in \Cref{gk}). Summing over $k$ yields the error $\Theta_k$ from \eqref{deltaomegal}, which estimates the error between $G_k$ (the initial choice of $G_0 = G$, suitably normalized) and $G_k^{\pm}$. We must then choose the $L_k$ so that $\Theta_k$ is sufficiently small and $\varsigma_{K_0} \ll n^{-1}$, and we will see that the choice from \eqref{deltaomegal} suffices. 
	
\end{rem} 

Having explained the quantities from \eqref{deltaomegal} in \Cref{delta0omegaklk}, let us bound them through the following lemma. Its proof, which is by direct computation, is provided in \Cref{ProofOmega}.

\begin{lem} 
	
	\label{thetak} 
	
	There exist constants $c = c(B) > 0$ and $C = C (B) > 1$ such that the following hold if $n > C$ and $\delta < c$. First, we have $\vartheta_k \le \Theta_k \le \delta_0^{1/2} / 2$ for each integer $k \in \big\llbracket 0, (\log n)^{3/4} \big\rrbracket$. Second, we have
	\begin{flalign}
		\label{k0} 
		\begin{aligned}
			& (\log n)^{1/3} \le K_0 \le (\log n)^{1/2}; \qquad e^{(\log n)^{1/6}} \le L_{K_0} \le e^{(\log n)^{1/2}}; \\
			& \displaystyle\frac{n_{K_0}}{n} \ge n^{-1/25000}; \qquad \qquad \quad \qquad \varsigma_{K_0-1} \le\vartheta_{K_0-1} \le 2n^{-13/15}.
		\end{aligned} 
	\end{flalign}
	
	\noindent Third, for any integer $k \in \llbracket 1, K_0 \rrbracket$, we have $|\log \vartheta_{k-1}|^{20} \le 4L_k / 3 \le \vartheta_{k-1}^{-1/5000}$. 
	
\end{lem}

For each integer $k \in \llbracket 0, K_0 \rrbracket$, we next inductively define a sequence of events $\Omega_k$, measurable with respect to $\mathcal{F}_{\ext} (\mathsf{R}_{k+1})$ (recall \Cref{property}), and sequences of functions $G_k^-, G_k^+ : \mathfrak{R}_k \rightarrow \mathbb{R}$ satisfying \eqref{equationxtd} on $\mathfrak{R}_k$. At $k = 0$, define the functions $G_0^-, G_0^+ : \overline{\mathfrak{R}}_0 \rightarrow \mathbb{R}$ and event $\Omega_0$ by setting
\begin{flalign}
	\label{gksx} 
	\begin{aligned}
		& G_0^- (s, x) = G_0 (s, x) = G_0^+ (s, x); \qquad \text{for each $(s, x) \in \overline{\mathfrak{R}}_0$},
	\end{aligned} 
\end{flalign} 

\noindent which (by the first part of \Cref{invariancesscale}) solve \eqref{equationxtd} on $\mathfrak{R}_0$, and 
\begin{flalign}
	\label{omega0}
	\Omega_0 = \Bigg\{ \displaystyle\sup_{(s, j/n_0) \in \overline{\mathfrak{R}}_0 \setminus \frac{1}{4} \cdot \overline{\mathfrak{R}}_1} \big| x_{j+n_0/2}^{(0)} (s) - G_0^+ (s, jn_0^{-1}) \big| \le \delta_0^{3/4} \Bigg\}. 	
\end{flalign}

\noindent Observe that $\Omega_0$ is measurable with respect to $\mathcal{F}_{\ext} (\mathsf{R}_1)$ since (by \eqref{rkrk}, \eqref{xks}, and the facts that $\omega_0 = 4 \omega_1$ and $n_0 = 4n_1$) it amounts to constraining $\mathsf{x}_j (t)$  for $(j, t) \notin \mathsf{R}_1$. 

We now let $k \in \llbracket 1, K_0 \rrbracket$ denote an integer. Assume that we have defined the functions $G_{k-1}^-$ and $G_{k-1}^+$, satisfying \eqref{equationxtd} on $\mathfrak{R}_{k-1}$, and the event $\Omega_{k-1}$, measurable with respect to $\mathcal{F}_{\ext} (\mathsf{R}_k)$. We will then define the functions $(G_k^-, G_k^+)$, as well as the event $\Omega_k$. 

The former will eventually be obtained by applying \Cref{fg1g2}, with the $F$ there equal to a rescaling of the $G_{k-1}^+$ here (and the $(g_1, g_2)$ there given by regular profiles approximating rescalings of $\bm{x}^{(k-1)}$, which will be guaranteed to exist by the $\textbf{PFL}$ event). However, to apply \Cref{fg1g2} we must first set some further notation and restrict to an additional event (that will verify its hypotheses). So, the precise prescription of $(G_k^-, G_k^+)$ is deferred to \Cref{gk} below. 

The event $\Omega_k$ will be the intersection of two events $\Omega_k^{(1)}$ and $\Omega_k^{(2)}$. The former is defined by
\begin{flalign}
	\label{1omega1} 
	& \Omega_k^{(1)} = \textbf{PFL}^{\bm{x}^{(k-1)}} \Big( -\displaystyle\frac{1}{8L_k}; \displaystyle\frac{1}{4n^{9/10}}; 2B \Big) \cap \textbf{PFL}^{\bm{x}^{(k-1)}} \Big( \displaystyle\frac{1}{8L_k}; \displaystyle\frac{1}{4n^{9/10}}; 2B \Big) \cap \Omega_{k-1},
\end{flalign}

\noindent where we recall the event $\textbf{PFL}$ from \Cref{eventpfl}. Since $\Omega_{k-1}$ is measurable with respect to $\mathcal{F}_{\ext} (\mathsf{R}_k)$ and since by \eqref{xks} the two regular profile events in \eqref{1omega1} amount to constraining $\bm{\mathsf{x}} (t)$ at $t \in \big\{ s_0 n^{1/3} - \omega_{k-1} n^{1/3} / 8L_k, s_0 n^{1/3}  \omega_{k-1} n^{1/3} / 8L_k \big\} = \{ s_0 n^{1/3} - \omega_k n^{1/3} / 2L_k, s_0 n^{1/3} + \omega_k n^{1/3} / 2L_k \}$, it follows from \eqref{rkrk} that $\Omega_k^{(1)}$ is also measurable with respect to $\mathcal{F}_{\ext} (\mathsf{R}_k)$ (and thus to $\mathcal{F}_{\ext} (\mathsf{R}_{k+1})$). The latter event $\Omega_k^{(2)}$ will be that on which $\bm{x}^{(k)}$ approximates $G_k^+$. Its precise definition must therefore be postponed until we define $G_k^+$ in \Cref{gk}; see \eqref{omegak42} below.

In what follows, we condition on $\mathcal{F}_{\ext} (\mathsf{R}_k)$ and restrict to the event $\Omega_k^{(1)}$. Then, there exists for each real number $t \in \big\{ -1 / 8L_k, 1 / 8L_k \big\}$ a function\footnote{This function may be taken to be measurable with respect to $\bm{x}^{(k-1)} (t)$ as an argument infimum of the continuous function described in \Cref{pfl0event} (see \cite[Theorem 18.13]{IDA}), and thus to $\mathsf{R}_k$, by the measurable selection theorem.} $\gamma_t^{(k-1)}: [-1/2, 1/2] \rightarrow \mathbb{R}$ such that for each $j \in \llbracket 1 - n_{k-1} / 2, n_{k-1} / 2 \rrbracket$ we have 
\begin{flalign}
	\label{xjkgamma2} 
	\big| x_{j+n_{k-1}/2}^{(k-1)} (t) - \gamma_t^{(k-1)} (jn_{k-1}^{-1}) \big| \le \displaystyle\frac{1}{4n^{9/10}}; \qquad \big\| \gamma_t^{(k-1)} - \gamma_t^{(k-1)} (0) \big\|_{\mathcal{C}^{50}} \le 2B;
\end{flalign}

\noindent observe here that we have shifted the argument of $\gamma_t^{(k-1)}$ by $1 / 2$ in comparison to \Cref{eventpfl}. 

\begin{rem}
	
	\label{x4} 
	
	By \eqref{xks}, we have $x_{i-j_0+n_k/2}^{(k)} (s) = 4 x_{i-j_0+n_{k-1}/2}^{(k-1)} (s/4)$, since $\omega_k = 4^{-k-1}$. Thus, we will often use the rescaling $f(z) \mapsto 4 f(z/4)$ to transfer information from scale $k-1$ to scale $k$. We will frequently use two properties of this rescaling. The first is that it preserves the equation \eqref{equationxtd}, by the $\alpha=\beta=4$ case of \Cref{scale21} of \Cref{invariancesscale}. The second is that it preserves gradients and reduces higher derivatives (more specifically, it multiplies the $k$-th derivative by $4^{1-k}$). 
\end{rem} 

Define for each $t \in \{ -1 / 2L_k, 1 / 2L_k \}$ the rescaled function $\widetilde{\gamma}_t^{(k-1)} : [-2, 2] \rightarrow \mathbb{R}$ by setting
\begin{flalign*}
	\widetilde{\gamma}_t^{(k-1)} (x) = 4 \gamma_{t/4}^{(k-1)} \Big( \displaystyle\frac{x}{4} \Big), \qquad \text{for each $x \in [-2, 2]$}.
\end{flalign*} 

\noindent By \eqref{xjkgamma2} and the fact that $x_{j+n_k/2}^{(k)} (s) = 4 x_{j+n_{k-1}/2}^{(k-1)} (s / 4)$ (due to \eqref{xks} and the equality $\omega_{k-1} = 4 \omega_k$), we have for each $t \in \{ -1 / 2L_k, 1 / 2L_k \}$ and $j \in \llbracket 1 - n_{k-1} / 2, n_{k-1} / 2 \rrbracket = \llbracket 1- 2n_k, 2n_k \rrbracket$ that, on $\Omega_k^{(1)}$,
\begin{flalign}
	\label{xjkgamma} 
	\big| x_{j+n_k/2}^{(k)} (t) - \widetilde{\gamma}_t^{(k-1)} (jn_k^{-1}) \big| \le n^{-9/10}; \qquad \big\| \widetilde{\gamma}_t^{(k-1)} - \widetilde{\gamma}_t^{(k-1)} (0) \big\|_{\mathcal{C}^{50}} \le 4B. 
\end{flalign}

Define the open rectangle $\widetilde{\mathfrak{R}}_{k-1} = 4 \cdot \mathfrak{R}_{k-1} = (-2L_{k-1}^{-1}, 2L_{k-1}^{-1}) \times (-2, 2)$ and the function $\widetilde{G}_{k-1}^+ : \widetilde{\mathfrak{R}}_{k-1} \rightarrow \mathbb{R}$ by rescaling $G_{k-1}^+$, namely, by setting
\begin{flalign}
	\label{gk1sx} 
	\widetilde{G}_{k-1}^+ (s, x) = 4 G_{k-1}^+ \Big( \displaystyle\frac{s}{4}, \displaystyle\frac{x}{4} \Big), \qquad \text{for each $(s, x) \in \widetilde{\mathfrak{R}}_{k-1}$},
\end{flalign}

\noindent which satisfies \eqref{equationxtd} on $\widetilde{\mathfrak{R}}_{k-1}$ by the first part (at $\alpha = \beta = 1 / 4$) of \Cref{invariancesscale}. Further define the functions $\mathfrak{y}_0^{(k)}, \mathfrak{y}_1^{(k)}, \mathfrak{z}_0^{(k)}, \mathfrak{z}_1^{(k)}, \widetilde{\mathfrak{z}}_0^{(k)}, \widetilde{\mathfrak{z}}_1^{(k)} : [ -2 / 3, 2 / 3] \rightarrow \mathbb{R}$ by setting 
\begin{flalign}
	\label{figik} 
	\mathfrak{y}_i^{(k)} (x) = \widetilde{\gamma}_{(2i-1)/2L_k}^{(k-1)} (x); \qquad \mathfrak{z}_i^{(k)} (x) = G_k \Big( \displaystyle\frac{2i-1}{2L_k}, x \Big); \qquad \widetilde{\mathfrak{z}}_i^{(k-1)} (x) = \widetilde{G}_{k-1}^+ \Big( \displaystyle\frac{2i-1}{2L_k}, x \Big),
\end{flalign}

\noindent for each index $i \in \{ 0, 1 \}$ and real number $x \in [ -2 / 3, 2 / 3]$.

We then have the following lemma, which states that $\widetilde{G}_{k-1}^+$ (with its derivatives) is close to $G_k$. It will be established in \Cref{ProofOmega123}. If $k = 1$, it follows directly from the fact that $\widetilde{G}_0^+ = G_1$ (by \eqref{gksx}, \eqref{gk1sx}, and \eqref{gksx2}); if $k \ge 2$, its proof will mainly use the fourth property in \Cref{gk} below, together with the rescaling in \Cref{x4}. In what follows, we recall the constant $C_0$ from above \eqref{deltaomegal} (as well as the norms $[f]_{k; \mathfrak{R}}$ from \Cref{FNotation}).

\begin{lem} 
	
	\label{omega3} 
	
	There exist constants $c = c(B) > 0$ and $C = C(B) > 1$ such that the following holds for any integer $k \in \llbracket 1, K_0 - 1 \rrbracket$.  If $n > C$ and $\delta < c$, then on the event $\Omega_k^{(1)}$ we have 
	\begin{flalign}
		\label{omega30} 
		\big[ \widetilde{G}_{k-1}^+ - G_k \big]_{1; \mathfrak{R}_k} + \omega_{k-1}^{-1} \cdot \displaystyle\sum_{d=2}^{50} \big[ \widetilde{G}_{k-1}^+ - G_k \big]_{d; \mathfrak{R}_k} \le 2C_0 \Theta_{k-1}.
	\end{flalign}
	
\end{lem} 

Observe from \Cref{omega3}; the bound $\Theta_{k-1} \le \delta_0^{1/2} / 2 \le 1 / 8C_0$ for sufficiently large $n$ and small $\delta$ (by \Cref{thetak}); the fact that $G_k \in \Adm_{1/2} (\mathfrak{R}_k)$ (which holds since $G \in \Adm_{1/2} (\mathfrak{R})$, by \eqref{g212}, together with the fact that the scaling \eqref{gksx2} producing $G_k$ from $G$ preserves gradients); and the fact that $\big\| G_k - G_k (0, 0) \big\|_{\mathcal{C}^{50} (\mathfrak{R})} \le 2B$ (as $\big\| G - G(0, 0) \big\|_{\mathcal{C}^{50} (\mathfrak{R})} \le B$ by \eqref{h0x2}) that 
\begin{flalign}
	\label{gk1b} 
	\widetilde{G}_{k-1}^+ \in \Adm_{1/4} (\mathfrak{R}_k), \quad \text{and} \quad \big\| \widetilde{G}_{k-1}^+ - \widetilde{G}_{k-1}^+ (0, 0) \big\|_{\mathcal{C}^{50} (\mathfrak{R}_{k-1})} \le 4B + 4 C_0 \Theta_{k-1} \le 5B. 
\end{flalign}

The next lemma states that $\mathfrak{y}_i^{(k)}$ is very close to $\widetilde{\mathfrak{z}}_i^{(k-1)}$. It will also be established in \Cref{ProofOmega123}, using \eqref{xjkgamma}; the definition of the event $\Omega_{k-1}^{(2)}$ in \eqref{omegak42} below (and that \eqref{omega0} of $\Omega_0$ if $k=1$); and the rescaling from \Cref{x4}.

\begin{lem}  
	
	\label{omega2} 
	
	There exist constants $c = c(B) > 0$ and $C = C(B) > 1$ such that the following holds for any integer $k \in \llbracket 1, K_0 - 1 \rrbracket$.  If $n > C$ and $\delta < c$, then on the event $\Omega_k^{(1)}$ we have 
	\begin{flalign}
		\label{omega20} 
		\displaystyle\max_{i \in \{ 0, 1 \}} \big\| \mathfrak{y}_i^{(k)} - \widetilde{\mathfrak{z}}_i^{(k-1)} \big\|_{\mathcal{C}^0 ([-2/3,2/3])} \le \vartheta_{k-1}.
	\end{flalign}
	
\end{lem}

The next definition introduces the functions $G_k^-$ and $G_k^+$, using \Cref{fg1g2}. Here, we define the open rectangles
\begin{flalign}
	\label{sksk} 
	\begin{aligned} 
		& \widehat{\mathfrak{R}}_k = \Big( -\displaystyle\frac{1}{2L_k}, \displaystyle\frac{1}{2L_k} \Big) \times \Big( -\displaystyle\frac{2}{3}, \displaystyle\frac{2}{3} \Big); \qquad \mathfrak{S}_k' = \Big( -\displaystyle\frac{3}{8L_k}, \displaystyle\frac{3}{8L_k} \Big) \times \Big( -\displaystyle\frac{1}{2}, \displaystyle\frac{1}{2} \Big). 
	\end{aligned}
\end{flalign}

\begin{definition} 
	
	\label{gk} 
	
	Apply \Cref{fg1g2} (translated by $( -\ell / 2L_k, -\ell / 2 )$) with the $(\ell, L, \vartheta)$ there equal to $( 4 / 3, 4L_k / 3, \vartheta_{k-1})$ here; the $(\mathfrak{R}, \mathfrak{S}, \mathfrak{S}')$ there equal to $\big(\widehat{\mathfrak{R}}_k, \mathfrak{R}_k, \mathfrak{S}_k' \big)$ here; the $(\varepsilon, B, m)$ there equal to $( 1 / 8, 5B, 50)$ here; and the $(F; g_0, g_1)$ there equal to $\big( \widetilde{G}_{k-1}^+ |_{\widehat{\mathfrak{R}}_k}; \mathfrak{y}_0^{(k)}; \mathfrak{y}_1^{(k)} \big)$ here (implicitly shifting all of these functions by the constant $\widetilde{G}_{k-1}^+ (0, 0)$). The assumptions of this lemma are verified by \Cref{omega2}, \eqref{gk1b}, \eqref{xjkgamma}, and \Cref{thetak}. This yields solutions\footnote{It is again quickly verified that we may choose $G_k^-$ and $G_k^+$, satisfying the constraints in \Cref{gk}, to be measurable with respect to $\big( \mathfrak{y}_i^{(k)}, \widetilde{G}_{k-1}^+$) and thus to $\mathcal{F}_{\ext} (\mathsf{R}_k)$, by the measurable selection theorem \cite[Theorem 18.13]{IDA}.} $G_k^-, G_k^+ \in \Adm_{1/8} (\mathfrak{R}) \cap \mathcal{C}^{50} (\overline{\mathfrak{R}}_k)$ to \eqref{equationxtd}, such that the following six properties hold. 
	
	\begin{enumerate}
		\item For each $i \in \{ 0, 1 \}$ and $x \in [ -2/5, 2/5]$, we have $G_k^- \big( \frac{2i-1}{2L_k}, x \big) = \mathfrak{y}_i^{(k)} (x) = G_k^+ \big( \frac{2i-1}{2L_k}, x \big)$. 
		\item For each $i \in \{ 0, 1 \}$ and $x \in [-1/2, 1/2]$, we have $G_k^- \big( \frac{2i-1}{2L_k}, x \big) \le \mathfrak{y}_i^{(k)} (x) \le G_k^+ \big( \frac{2i-1}{2L_k}, x \big)$. 
		\item We have $\big\| G_k^- - G_k^- (0, 0) \big\|_{\mathcal{C}^{45} (\mathfrak{S}_k)} + \big\| G_k^+ - G_k^+ (0, 0) \big\|_{\mathcal{C}^{45} (\mathfrak{S}_k)} \le 2C_0$.
		\item We have $\| G_k^- - \widetilde{G}_{k-1}^+ \|_{\mathcal{C}^{50} (\mathfrak{S}_k')} + \| G_k^+ - \widetilde{G}_{k-1}^+ \|_{\mathcal{C}^{50} (\mathfrak{S}_k')} \le C_0 \vartheta_{k-1}^{3/4}$.
		\item For each $(t, x) \in [ -1 / 2L_k, 1 / 2L_k ] \times [-1/3, 1/3]$, we have $\big| G_k^+ (t, x) - G_k^- (t, x) \big| \le C_0 e^{-c_0 L_k^{1/8}}$.
		\item For each $(t, x) \in [ -1 / 2L_k, 1 / 2L_k ] \times \{ -1/2, 1/2 \}$, we have $G_k^- (t, x) \le \widetilde{G}_{k-1}^+(t,x) - \vartheta_{k-1} < \widetilde{G}_{k-1}^+ (t, x) + \vartheta_{k-1} \le G_k^+ (t, x)$. 
	\end{enumerate}
	
\end{definition}

Then, define the event
\begin{flalign}
	\label{omegak42}
	\Omega_k^{(2)} & = \Bigg\{ \displaystyle\sup_{\substack{(s,j/n_k) \in \mathfrak{R}_k \setminus \frac{1}{4} \cdot \mathfrak{R}_{k+1} \\ |j| \le n_k/4}} \big| x_{j+n_k/2}^{(k)} (s) - G_k^+ (s, jn_k^{-1}) \big| \le \displaystyle\frac{\vartheta_k}{5} \Bigg\},
\end{flalign}

\noindent which is measurable with respect to $\mathcal{F}_{\ext} (\mathsf{R}_{k+1})$ since (by \eqref{rkrk}, \eqref{xks}, and the facts that $\omega_{k-1} = 4 \omega_k$ and $n_{k-1} = 4 n_k$) it amounts to constraining the paths $\mathsf{x}_j (t)$ for $(j, t) \notin \mathsf{R}_{k+1}$. Also define the event (not measurable with respect to $\mathcal{F}_{\ext} (\mathsf{R}_{k+1})$)
\begin{flalign}
	\label{omegak43} 
	\Omega_k^{(3)} = \Bigg\{ \displaystyle\sup_{\substack{(s, j/n_k) \in \mathfrak{R}_k \\ |j| \le n_k/4}} \Big| x_{j+n_k/2}^{(k)} (s) - G_k^+ \big( s, jn_k^{-1} \big) \Big| \le \displaystyle\frac{\vartheta_k}{5} \Bigg\}.
\end{flalign}

\noindent These events indicate that the random paths in $\bm{x}^{(k)}$ closely approximate the limit shape $G_k^+$. That these two events likely hold (stated as \Cref{omega4}, and shown in \Cref{ProofOmega4}) will follow the sandwiching scheme outlined in the beginning of this section and \Cref{Proof0Boundary}.

Observing that $\Omega_k^{(3)} \subseteq \Omega_k^{(2)}$, further define the events
\begin{flalign}
	\label{omegak0}
	\Omega_k = \Omega_k^{(1)} \cap \Omega_k^{(2)}; \qquad \Omega_k' = \Omega_k^{(1)} \cap \Omega_k^{(3)} \subseteq \Omega_k.
\end{flalign}	

The next lemma, to be shown in \Cref{ProofOmega}, indicates that the final event $\Omega_{K_0}'$ is likely.

\begin{lem}
	
	\label{omegakestimate}
	
	There exist constants $c = c(B) > 0$ and $C = C(B) > 1$ such that $\mathbb{P} [\Omega_{K_0}' ] \ge 1 - n^{-19}$, whenever $n > C$ and $\delta < c$.
	
\end{lem}

Given this result, we can establish \Cref{h0xlocal}.

\begin{proof}[Proof of \Cref{h0xlocal}]

	Throughout this proof, we abbreviate $K = K_0$ and assume in what follows that the event $\Omega_K'$ holds (which we may by \Cref{omegakestimate}). Then, define $h_{j_0; s_0} : [-\eta, \eta] \rightarrow \mathbb{R}$ by setting
	\begin{flalign}
		\label{hj0s0}
		h_{j_0; s_0} (s) = \omega_{K} \cdot G_{K}^+ (\omega_{K}^{-1} s, 0).
	\end{flalign}	
	
	\noindent To verify \eqref{xj0s}, observe that
	\begin{flalign*}
		\displaystyle\sup_{|s| \le \omega_K / 2L_K} \big| x_{j_0} (s+s_0) - h_{j_0; s_0} (s) \big| & = \displaystyle\sup_{|s| \le \omega_K /2L_{K}} \big| x_{j_0} (s+s_0 ) - \omega_{K}  \cdot G_{K}^+ (\omega_{K}^{-1} s, 0) \big| \\ 
		& = \omega_K \cdot \displaystyle\sup_{|s| \le 1 / 2L_K} \big| x_{n_K/2}^{(K)} (s) - G_K^+ (s, 0) \big| \\
		& \le \omega_{K} \vartheta_{K} \le 2 \omega_K n^{-13/15} < n^{-13/15},
	\end{flalign*}
	
	\noindent where the first statement follows from \eqref{hj0s0}; the second from \eqref{xks}; the third from the fact that $\Omega_K^{(3)} \subseteq \Omega_K'$ holds; the fourth from \eqref{k0}; and the fifth from the fact that $\omega_K \le 1/4$. Since $\omega_{K} \ge 3 \eta L_K$ by the definition of $K = K_0$, this verifies \eqref{xj0s}.
	
	Next let us confirm \eqref{bhs2}, starting with the first bound there. Observe from \eqref{g2121} and \eqref{gksx2} that 
	\begin{flalign}
		\label{g2122}
		\displaystyle\max_{|s| < 1/2L_K}  \big| \omega_K^{-1} \cdot \partial_s^2 G_K (s, 0) + 2^{1/2} \big|  = \displaystyle\max_{|s| < \omega_K/2L_K} \big| \partial_t^2 G(s n^{-1/3}_0 + \omega_K s, 0) + 2^{1/2} \big| = 0.
	\end{flalign}
	
	\noindent Therefore,  
	\begin{flalign*}
		\displaystyle\sup_{|s| \le \eta} \big| \partial_s^2 h_{j_0; s_0} (s) + 2^{1/2} \big| & \le \displaystyle\sup_{|s| < 1/3L_K} \big| \omega_K^{-1} \cdot \partial_s^2 G_K^+ (s, 0) + 2^{1/2} \big| \\
		& \le \omega_K^{-1} \cdot \displaystyle\sup_{|s| < 1/3L_K} \big| \partial_s^2 G_K^+ (s, 0) - \partial_s^2 G_K (s, 0) \big| + \delta \\
		& \le \omega_K^{-1} \cdot \displaystyle\sup_{|s| < 1/3L_K} \big| \partial_s^2 \widetilde{G}_{K-1}^+ (s, 0) - \partial_s^2 G_K (s,0) \big| + C_0 \omega_K^{-1} \vartheta_{K-1}^{3/4} + \delta \\
		& \le 8 C_0 \Theta_{K-1} + 4C_0 \omega_{K-1}^{-1} \vartheta_{K-1}^{3/4} + \delta \le 8 C_0 \delta_0^{1/2} + n^{-1/2} + \delta < \delta_0^{1/3}.
	\end{flalign*} 
	
	\noindent Here, in the first bound we applied \eqref{hj0s0}, replaced $s$ by $\omega_K^{-1} s$, and used the fact that $\omega_K^{-1} \eta \le (3L_K)^{-1}$; in the second we applied \eqref{g2122}; in the third we applied the fourth statement of \Cref{gk}; in the fourth we applied 
	\eqref{omega30} (and that $\omega_{K-1} = 4 \omega_K$); in the fifth we applied the facts that $\vartheta_{K-1} \le 2 n^{-13/15}$ and that $\omega_{K-1} = 4 \omega_K \le 4 n^{-1/25000}$ (both as consequences of \eqref{deltaomegal} and \Cref{thetak}); and in the sixth we used the fact that $n$ is sufficiently large and $\delta$ is sufficiently small. Since $\delta_0 = \delta^{1/2} + (\log n)^{-1}$ by \eqref{deltaomegal}, this establishes the first statement of \eqref{bhs2}. 
	
	To establish the second statement of \eqref{bhs2}, observe that 
	\begin{flalign*}
		\| h_{j_0; s_0} \|_{\mathcal{C}^1} \le \| G_K^+ \|_{\mathcal{C}^1 (\mathfrak{S}_K')} \le \| \widetilde{G}_{K-1}^+ \|_{\mathcal{C}^1 (\mathfrak{R}_K)} + C_0 \vartheta_{K-1}^{3/4} \le 5B + C_0 \vartheta_{K-1}^{3/4} < 10B,
	\end{flalign*}
	
	\noindent where the first statement holds by \eqref{hj0s0}, the second by the fourth part of \Cref{gk} (with the fact that $[-\eta, \eta] \times \{ 0 \} \subset \mathfrak{S}_K'$), the third by \eqref{gk1b} (with the fact that $\Omega_K^{(1)}$ holds) and \Cref{thetak}, and the fourth by \eqref{k0} and the fact that $n$ is sufficiently large. This confirms the second part of \eqref{bhs2}, verifying the proposition. 		
\end{proof}

\section{Proofs of Results From Section \ref{LOCALHJ}}

\label{Proof0Omega}

Throughout, we recall the notation from \Cref{LOCALHJ}.   

\subsection{Proofs of \Cref{thetak} and \Cref{omegakestimate}} 

\label{ProofOmega}

We begin by proving \Cref{omegakestimate}, which will proceed inductively, showing that the events $\Omega_k^{(i)}$ are likely if $\Omega_{k-1}$ holds. This is summarized through the following two lemmas; the first will be shown in \Cref{ProofOmega123} and the second in \Cref{ProofOmega4}.

\begin{lem} 
	
	\label{omega1} 
	
	There exist constants $c = c(B) > 0$ and $C = C(B) > 1$ such that the following holds for $n > C$ and $\delta < c$. We have $\mathbb{P} \big[ \Omega_{k-1} \cap (\Omega_k^{(1)})^{\complement} \big] \le 2n^{-20}$ for any integer $k \in \llbracket 1, K_0 \rrbracket$.
	
\end{lem} 

\begin{lem} 
	
	\label{omega4} 
	
	There exist constants $c = c(B) > 0$ and $C = C(B) > 1$ such that the following statements hold for $n > C$ and $\delta < c$. 
	
	\begin{enumerate} 
		\item \label{omega41} We have $\mathbb{P} \big[ \Omega_0^{\complement} \big] \le n^{-20}$. 
		\item \label{omega42} For any integer $k \in \llbracket 1, K_0 \rrbracket$, we have 
	\begin{flalign*} 
		\mathbb{P} \big[ \Omega_k^{(1)} \cap ( \Omega_k^{(2)} )^{\complement} \big] \le \mathbb{P} \big[ \Omega_k^{(1)} \cap (\Omega_k^{(3)})^{\complement} \big] \le n^{-20}. 
	\end{flalign*}  
	\end{enumerate} 
\end{lem} 

Given these results, we can quickly establish \Cref{omegakestimate}. 

\begin{proof}[Proof of \Cref{omegakestimate}]
	
	It suffices to show the bound 
	\begin{flalign}
		\label{omegak1} 
		\mathbb{P} \big[ \Omega_k \big] \ge \mathbb{P} [ \Omega_k'] \ge 1 - 3(k+1) n^{-20}, \qquad \text{for each integer $k \in \llbracket 0, K_0 \rrbracket$},
	\end{flalign}
	
	\noindent from which the lemma follows from taking $k = K_0 \le \log n$ (where the last bound holds by \eqref{k0}). To this end, we induct on $k \in \llbracket 0, K_0 \rrbracket$; for $k = 0$, \eqref{omegak1} holds by the first statement of \Cref{omega4}. We then assume that \eqref{omegak1} holds for some $k \in \llbracket 0, K_0 - 1 \rrbracket$ and show that it continues to hold for $k$ replaced by $k+1$. Since $\Omega_k' \subseteq \Omega_k$, it suffices to show only the second inequality in \eqref{omegak1}. This follows from the estimates 
	\begin{flalign*}
		\mathbb{P} [\Omega_{k+1}'] & \ge \mathbb{P} [\Omega_k] - \mathbb{P} \big[ \Omega_k \cap (\Omega_{k+1}^{(1)})^{\complement} \big] - \mathbb{P} \big[ \Omega_k \cap \Omega_{k+1}^{(1)} \cap ( \Omega_{k+1}^{(3)} )^{\complement} \big] \\
		& \ge  1 - 3(k+1) n^{-20} - \mathbb{P} \big[ \Omega_k \cap \big( \Omega_{k+1}^{(1)} \big)^{\complement} \big] - \mathbb{P} \big[  \Omega_{k+1}^{(1)} \cap ( \Omega_{k+1}^{(3)} )^{\complement} \big] \ge 1 - 3(k+4) n^{-20},
	\end{flalign*}
	
	\noindent where in the first inequality we applied a union bound and \eqref{omegak0}; in the second we applied the inductive hypothesis; and in the third we applied \Cref{omega1} and the second statement of \Cref{omega4}. This yields the lemma.
\end{proof} 

We next establish \Cref{thetak}.

\begin{proof}[Proof of \Cref{thetak}]

	First observe from \eqref{deltaomegal} that for sufficiently large $n$ and small $\delta$ we have $\varsigma_k \le 1/4$ and $n^{-13/15} \le 1/4$. Hence, $\vartheta_k \le 1/2$ and so the definition \eqref{deltaomegal} of $\Theta_k$ yields $\Theta_k \ge \omega_k^{-1} \vartheta_k^{3/4} \ge \vartheta_k$. Thus, for any integer $k \ge 1$, we have 
	\begin{flalign*}
		\vartheta_k \le \Theta_k = \delta_0^{3/4} + \displaystyle\sum_{j=0}^k \omega_j^{-1} \vartheta_j^{3/4} & \le \delta_0^{3/4} + 4 \displaystyle\sum_{j=0}^k 4^j (\varsigma_j^{3/4} + n^{-1/2}) \\
		& \le \delta_0^{3/4} + 4^{k+2} n^{-1/2} + 20 C_0 \displaystyle\sum_{j=0}^k 4^j \exp \Big( -\displaystyle\frac{3c_0}{8} L_j^{1/80} \Big),
	\end{flalign*}
	
	\noindent where we have used the definitions \eqref{deltaomegal} of $\Theta_k$, $\omega_j$, $\vartheta_j$, and $\varsigma_j$, with the bound $\vartheta_j^{3/4} = (\varsigma_j + n^{-13/15})^{3/4} \le \varsigma_j^{3/4} + n^{-1/2}$. Since, for sufficiently small $\delta_0 > 0$, we have $L_j = \delta_0^{-\sqrt{j+1}} \ge (2c_0^{-1} (j+1))^{160} \delta_0^{-1/2}$ for each integer $j \ge 0$, it follows that 
	\begin{flalign*}
		\vartheta_k \le \Theta_k \le \delta_0^{3/4} + 4^{k+2} n^{-1/2} + 20C_0 \displaystyle\sum_{j=0}^k \exp \big( 2j - (j+1)^2 \delta_0^{-1/200} \big) \le 2 \delta_0^{3/4} + 4^{k+2} n^{-1/2},
	\end{flalign*}
	
	\noindent for sufficiently large $n$ and small $\delta$. Hence, since for $k \le (\log n)^{3/4}$ and sufficiently large $n$ we have $4^{k+2} n^{-1/2} \le n^{-1/4} \le (\log n)^{-1} \le \delta_0 \le \delta_0^{3/4}$, it follows that $\vartheta_k \le \Theta_k \le 3\delta_0^{3/4} \le \delta_0^{1/2} / 2$, confirming the first statement of the lemma. 
	
	We next verify \eqref{k0}. To establish the first bound there, on $K_0$, observe for $k \le (\log n)^{1/3} + 1$ that for sufficiently large $n$ we have
	\begin{flalign*}
		\omega_k L_k^{-1} = 4^{-k-1} \delta_0^{\sqrt{k+1}} \ge 4^{-2(\log n)^{1/3}} \cdot (\log n)^{-2 (\log n)^{1/3}} > e^{-(\log n)^{2/5}} > 3 \eta,
	\end{flalign*}
	
	\noindent which indicates that $K_0 \ge (\log n)^{1/3}$. Here in the first statement we used the definition \eqref{deltaomegal} of $\omega_k$ and $L_k$; in the second we used the facts that $\delta_0 \ge (\log n)^{-1}$ and that $\sqrt{k+1} \le k+1 \le 2(\log n)^{1/3}$; in the third we used the fact that $2 (\log n)^{1/3} \cdot \log (4 \log n) < (\log n)^{2/5}$ for sufficiently large $n$; and in the fourth we recalled that $\eta = e^{-\sqrt{\log n}}$. For $k \ge (\log n)^{1/2}$, we have that $\omega_k L_k^{-1} \le 4^{-k-1} < e^{\sqrt{\log n}} < 3 \eta$, indicating that $K_0 \le (\log n)^{1/2}$, verifying the first bound in \eqref{k0}. 
	
	The second bound in \eqref{k0}, given by $e^{(\log n)^{1/6}} \le L_{K_0} \le e^{(\log n)^{1/2}}$, follows from the first bound in \eqref{k0}, together with the facts that $L_k = \delta_0^{-\sqrt{k+1}}$ and that $e^{-(\log n)^{1/5}} \le \delta_0 \le e^{-1}$ for sufficiently large $n$ and small $\delta$. The third follows from the fact that $n_{K_0} n^{-1} = \omega_{K_0} = 4^{-K_0-1} \ge e^{-2\sqrt{\log n}-1} \ge n^{-1/25000}$. The fourth follows from the fact that 
	\begin{flalign*} 
		\vartheta_{K_0-1} = \varsigma_{K_0-1} + n^{-13/15} & = n^{-13/15} + 5C_0 \exp \Big( -\displaystyle\frac{c_0}{2} L_{K_0-1}^{1/80} \Big) \\ 
		& \le n^{-13/15} + 5C_0 \exp \Big( -\displaystyle\frac{c_0}{2} e^{(\log n)^{1/480}} \Big) \\
		& \le n^{-13/15} + 5C_0 e^{-(\log n)^2} \le 2n^{-13/15},
	\end{flalign*} 		
	
	\noindent where in the first and second statements we used the definitions \eqref{deltaomegal} of $\vartheta_{K_0-1}$ and $\varsigma_{K_0-1}$; in the third we used the fact that $L_{K_0-1} = \delta_0^{-\sqrt{K_0}} \ge \delta_0^{(\log n)^{1/6}} \ge e^{(\log n)^{1/6}}$ (by \eqref{deltaomegal} and the first bound in \eqref{k0}); and the fourth and fifth follow since $n$ is sufficiently large. This establishes \eqref{k0} and thus the second statement of the lemma. 	
	
	To establish the third statement of the lemma, observe for any real numbers $a, b \in (0, 1 / 4)$ that $\big| \log (a+b) \big| \le |\log a|$. Applying this with $(a, b) = (\varsigma_k, n^{-13/15})$ yields for sufficiently large $n$ and small $\delta$ that
	\begin{flalign*}
		|\log \vartheta_k|^{20} \le |\log \varsigma_k|^{20} = \Big( \displaystyle\frac{c_0}{2} L_k^{1/80} - \log (5C_0) \Big)^{20} \le 2^{-20} c_0^{20} L_k^{1/4} \le \displaystyle\frac{4L_k}{3},
	\end{flalign*}
	
	\noindent establishing the first bound of the third statement. To establish the second bound there, observe for any real numbers $a, b \in ( 0, 1 / 4)$ and $r \in (0, 1)$ that $(a+b)^{-r} \ge 2^{-r} \min \big\{ (2a)^{-r}, (2b)^{-r} \big\}$. Setting $(a, b; r) = (\varsigma_k, n^{-13/15}; 1 / 5000)$ yields 
	\begin{flalign}
		\label{abr2} 
		\vartheta_k^{-1/5000} \ge \min \big\{ (2\varsigma_k)^{-1/5000}, n^{1/10000} \big\}.
	\end{flalign}
	
	\noindent By the definition \eqref{deltaomegal}, it is quickly verified that for sufficiently large $L_k$ (and hence sufficiently large $n$ and small $\delta)$ we have $(2 \varsigma_k)^{-1/5000} > 4L_k / 3$. Moreover, we have $4L_k / 3 \le 2L_k \le 2L_{K_0} \le 2 e^{\sqrt{\log n}} \le n^{1/10000}$, where in the first and second bounds we used the facts that $L_k$ is positive and increasing in $k$ (by its definition \eqref{deltaomegal}); in the third we used the second statement of \eqref{k0}; and in the fourth we used the fact that $n$ is sufficiently large. Together with \eqref{abr2}, these two bounds verify the third statement of the lemma.
\end{proof}

\subsection{Proofs of \Cref{omega3}, and \Cref{omega2}, and \Cref{omega1}}

\label{ProofOmega123}

In this section we establish first \Cref{omega1}, then \Cref{omega3}, and next \Cref{omega2}. 

\begin{proof}[Proof of \Cref{omega1}]
	
	By \eqref{1omega1}, we must verify that 
	\begin{flalign*} 
		 \mathbb{P} \bigg[ \textbf{PFL}^{\bm{x}^{(k-1)}} \Big( -\frac{1}{8L_k}; \displaystyle\frac{1}{4 n^{9/10}}; 2B \Big)^{\complement} \cup \textbf{PFL}^{\bm{x}^{(k-1)}} \Big( \frac{1}{8L_k}; \displaystyle\frac{1}{4 n^{9/10}}; 2B \Big)^{\complement} \bigg] \le 2n^{-20}.
	\end{flalign*} 	
	
	\noindent To that end, by \eqref{pflx2} and a union bound, it suffices to show that 
	\begin{flalign}
		\label{pflxkpflx}
		\begin{aligned} 
			&  \textbf{PFL}^{\bm{x}} \Big( \displaystyle\frac{s_0}{n^{1/3}} - \displaystyle\frac{\omega_{k-1}}{8 L_k}; n^{-19/20}; B \Big) \subseteq  \textbf{PFL}^{\bm{x}^{(k-1)}} \Big( - \displaystyle\frac{1}{8L_k}; \displaystyle\frac{1}{4n^{9/10}}; 2B \Big); \\
			&  \textbf{PFL}^{\bm{x}} \Big( \displaystyle\frac{s_0}{n^{1/3}} + \displaystyle\frac{\omega_{k-1}}{8 L_k}; n^{-19/20}; B \Big) \subseteq \textbf{PFL}^{\bm{x}^{(k-1)}} \Big( \displaystyle\frac{1}{8L_k}; \displaystyle\frac{1}{4n^{9/10}}; 2B \Big).
		\end{aligned} 
	\end{flalign}
	
	\noindent This will follow quickly from rescaling. We only show the first bound in \eqref{pflxkpflx}, as the proof of the second is entirely analogous. To that end, set $t_1 = -1 / 8L_k$ and $s_1 = s n^{-1/3}_0 + t_1 \omega_k$, and observe on the event $\textbf{PFL}^{\bm{x}} (s_1; n^{-19/20}; B)$ that there exists a function $\gamma_{s_1} : [0, 1] \rightarrow \mathbb{R}$ such that 
	\begin{flalign}
		\label{xjs11}
		\displaystyle\max_{j \in \llbracket 1, n \rrbracket} \big| x_j (s_1) - \gamma_{s_1} (jn^{-1}) \big| < n^{-19/20}; \qquad  \big\| \gamma_{s_1} - \gamma_{s_1} (0) \big\|_{\mathcal{C}^{50}} \le B.
	\end{flalign}
	
	 We then define $\gamma : [0, 1] \rightarrow \mathbb{R}$ by rescaling $\gamma_{s_1}$, namely, by setting
	\begin{flalign}
		\label{gamma0k} 
		\gamma (x) = \omega_{k-1}^{-1} \cdot \gamma_{s_1} \Big( \omega_{k-1} x + \displaystyle\frac{j_0}{n} - \displaystyle\frac{\omega_{k-1}}{2} \Big), 
	\end{flalign} 
	
	\noindent so that 
	\begin{flalign} 
		\label{gamma0k1}
		\big\| \gamma - \gamma (0) \big\|_{\mathcal{C}^{50}} \le \big\| \gamma_{s_1} - \gamma_{s_1}(0) \big\|_{\mathcal{C}^{50}} + \big\| \gamma_{s_1} - \gamma_{s_1} (0) \big\|_{\mathcal{C}^1}  \le 2B,
	\end{flalign}
	
	\noindent where in the last inequalities we used the facts that $[\gamma]_m = \omega_{k-1}^{m-1} \cdot [\gamma_{s_1}]$, for each integer $m \ge 0$, and that $\omega_{k-1} \in [0, 1]$. We further have for sufficiently large $n$ that
	\begin{flalign*}
		\displaystyle\max_{i \in \llbracket 1, n_{k-1} \rrbracket}  \bigg| x_i^{(k-1)} (t_1) -  \gamma \Big(\displaystyle\frac{i}{n_{k-1}} \Big) \bigg| & = \omega_{k-1}^{-1} \cdot \displaystyle\max_{i \in \llbracket 1, n_{k-1} \rrbracket} \bigg| x_{i+j_0 - n_{k-1} / 2} (s_1) - \gamma_{s_1} \Big( \displaystyle\frac{i + j_0}{n} - \displaystyle\frac{\omega_{k-1}}{2} \Big) \bigg| \\
		& \le n^{-11/12} < \displaystyle\frac{1}{4n^{9/10}}.
	\end{flalign*} 
	
	\noindent Here, in the first statement we used \eqref{gamma0k}, with the facts that $n_{k-1} = \omega_{k-1} n$ and that $x_i^{(k-1)} (t_1) = \omega_{k-1}^{-1} \cdot x_{i+j_0-n_{k-1}/2} (s_1)$ (which holds by \eqref{xks}); in the second, we used the first statement of \eqref{xjs11} and the fact that $\omega_{k-1}^{-1} \le \omega_{K_0}^{-1} \le 4^{\sqrt{\log n} + 1} < n^{1/30}$ (by \eqref{k0}) for sufficiently large $n$; and in the third we used the fact that $n$ is sufficiently large. This, together with \eqref{gamma0k1} (and \Cref{eventpfl}) yields \eqref{pflxkpflx} and thus the lemma. 
\end{proof}

\begin{proof}[Proof of \Cref{omega3}] 
	
	We induct on $k \in \llbracket 1, K_0 - 1 \rrbracket$. To verify the lemma in the case $k = 1$, observe by \eqref{gk1sx} and \eqref{gksx} that $\widetilde{G}^+_{k-1} (s, x) = 4G_0^+ ( s / 4, x / 4 ) = 4 G_0 ( s / 4, x / 4)$. Moreover, by \eqref{gksx2} (with the fact that $\omega_0 = 4 \omega_1$), we have $G_1 (s, x) = 4 G_0 ( s / 4, x / 4)$. Thus, $G_1 = \widetilde{G}_0$, which gives \eqref{omega30} at $k = 1$.
	
	So, fix some integer $k \ge 2$ and assume that \eqref{omega30} holds for $k$ there equal to $k-1$ here. In what follows, we restrict to the event $\Omega_k^{(1)}$ and then must show that \eqref{omega30} holds. To this end, first observe since we have restricted to the event $\Omega_{k-1}^{(1)} \subseteq \Omega_{k-1} \subseteq \Omega_k^{(1)}$, the inductive hypothesis (with the definition \eqref{deltaomegal} of $\Theta_{k-2}$) yields
	\begin{flalign}
		\label{gk2gk1} 
		\big[ \widetilde{G}_{k-2}^+ - G_{k-1} \big]_{1; \mathfrak{R}_{k-1}} + \omega_{k-2}^{-1} \cdot \displaystyle\sum_{d=2}^{50} \big[ \widetilde{G}_{k-2}^+ - G_{k-1} \big]_{d; \mathfrak{R}_{k-1}} \le 2 C_0 \delta_0^{3/4} + 2 C_0 \displaystyle\sum_{k=1}^{k-2} \omega_j^{-1} \vartheta_j^{3/4}.
	\end{flalign}
	
	\noindent Now, define $\widehat{G}_{k-2}^+ : 4 \cdot \mathfrak{R}_{k-1} \rightarrow \mathbb{R}$ by rescaling $\widetilde{G}_{k-2}^+$, namely, by setting
	\begin{flalign}
		\label{gk2gk2}
		\widehat{G}_{k-2}^+ (s, x) = 4 \widetilde{G}_{k-2}^+ \Big(\displaystyle\frac{s}{4}, \displaystyle\frac{x}{4} \Big),
	\end{flalign}
	
	\noindent for each $(s, x) \in 4 \cdot \mathfrak{R}_{k-1}$. Since we have from \eqref{gksx2} (with the fact that $\omega_{k-1} = 4 \omega_k$) that $G_k (s, x) = 4 G_{k-1} ( s / 4, x / 4)$, this yields (recall \Cref{x4})
	\begin{flalign*}
		& \nabla \widehat{G}_{k-2}^+ (s, x) - \nabla G_k (s, x) = \nabla \widetilde{G}_{k-2}^+ \Big( \displaystyle\frac{s}{4}, \displaystyle\frac{x}{4} \Big) - \nabla G_{k-1} \Big( \displaystyle\frac{s}{4}, \displaystyle\frac{x}{4} \Big), \qquad \quad \text{for each $(s, x) \in \mathfrak{R}_k$}; \\ 
		& \big[ \widehat{G}_{k-2}^+ - G_k \big]_{d; \mathfrak{R}_k} = 4^{1-d} \cdot \big[ \widetilde{G}_{k-2}^+ - G_{k-1} \big]_{d; \frac{1}{4} \cdot \mathfrak{R}_k},	\qquad \qquad \qquad \qquad \text{for each $d \ge 2$}.
	\end{flalign*}
	
	\noindent Together with \eqref{gk2gk1} and the facts that $\frac{1}{4} \cdot \mathfrak{R}_{k-1} \subseteq \mathfrak{R}_k$ and $\omega_{k-2} = 4 \omega_{k-1}$, this gives
	\begin{flalign*}
		\big[ \widehat{G}_{k-2}^+ - G_k \big]_{1; \mathfrak{R}_k} + \omega_{k-1}^{-1} \cdot \displaystyle\sum_{d=2}^{50} 4^{d-2} \cdot \big[ \widehat{G}_{k-2}^+ - G_k \big]_{d; \mathfrak{R}_k} \le 2 C_0 \delta_0^{3/4} + 2 C_0 \displaystyle\sum_{k=1}^{k-2} \omega_j^{-1} \vartheta_j^{3/4}.
	\end{flalign*}
	
	\noindent Hence, to verify that \eqref{omega30} holds, it suffices to show that 
	\begin{flalign}
		\label{gk1gk1} 
		\big[ \widetilde{G}_{k-1}^+ - \widehat{G}_{k-2}^+ \big]_{1; \mathfrak{R}_k} \le C_0 \vartheta_{k-1}^{3/4}; \qquad \displaystyle\sum_{d=2}^{50} 4^{d-2} \cdot \big[ \widetilde{G}_{k-1}^+ - \widehat{G}_{k-2}^+ \big]_{d; \mathfrak{R}_k}  \le C_0 \vartheta_{k-1}^{3/4}.
	\end{flalign}
	
	\noindent By \eqref{gk2gk2} and the fact from \eqref{gk1sx} that $\widetilde{G}_{k-1}^+ (s, x) = 4 G_{k-1}^+ ( s / 4, x / 4 )$, we have (recall \Cref{x4})
	\begin{flalign*}
		& \nabla \widehat{G}_{k-2}^+ (s, x) - \nabla \widetilde{G}_{k-1}^+ (s, x) = \nabla \widetilde{G}_{k-2}^+ \Big( \displaystyle\frac{s}{4}, \displaystyle\frac{x}{4} \Big) - \nabla G_{k-1}^+ \Big( \displaystyle\frac{s}{4}, \displaystyle\frac{x}{4} \Big), \qquad \text{for each $(s, x) \in 4 \cdot \mathfrak{R}_{k-1}$}; \\ 
		& \big[ \widehat{G}_{k-2}^+ - \widetilde{G}_{k-1}^+ \big]_{d; \mathfrak{R}_k} = 4^{1-d} \cdot \big[ \widetilde{G}_{k-2}^+ - G_{k-1}^+ \big]_{d; \frac{1}{4} \cdot \mathfrak{R}_k},\qquad \qquad \qquad \quad \text{for each $d \ge 2$}.
	\end{flalign*}
	
	\noindent and so (again since $\mathfrak{R}_k \subseteq \mathfrak{R}_{k-1} \subseteq 4 \cdot \mathfrak{R}_{k-1}$) to confirm \eqref{gk1gk1} we may show that 	
	\begin{flalign*}
		\big[ G_{k-1}^+ - \widetilde{G}_{k-2}^+ \big]_{1; \frac{1}{4} \cdot \mathfrak{R}_k} \le C_0 \vartheta_{k-1}^{3/4}; \qquad \displaystyle\sum_{d=2}^{50} \big[ G_{k-1}^+ - \widetilde{G}_{k-2}^+ \big]_{d; \frac{1}{4} \cdot \mathfrak{R}_k}  \le C_0 \vartheta_{k-1}^{3/4}.
	\end{flalign*}
	
	\noindent Both follow from the fourth property in \Cref{gk} (with the $k$ there equal to $k-1$ here), and the fact that $\mathfrak{R}_k \subseteq 4 \cdot \mathfrak{S}_{k-1}'$ (by \eqref{rkrk} and \eqref{sksk}). This verifies \eqref{gk1gk1} and thus the lemma. 
\end{proof}

\begin{proof}[Proof of \Cref{omega2}] 
	
	Throughout this proof, we restrict to the event $\Omega_k^{(1)}$; we then must show that \eqref{omega20} holds. We only verify the bound $\big\| \mathfrak{y}_i^{(k)} - \widetilde{\mathfrak{z}}_i^{(k-1)} \big\|_{\mathcal{C}^0} \le \vartheta_{k-1}$ at $i = 0$, as the proof that it holds at $i = 1$ is entirely analogous. To this end, observe that
	\begin{flalign}
		\label{yz1} 
		\begin{aligned}
			\big\| \mathfrak{y}_0^{(k)} - \widetilde{\mathfrak{z}}_0^{(k-1)} \big\|_{\mathcal{C}^0} & = \displaystyle\sup_{|x| \le 2/3} \bigg| \widetilde{\gamma}_{-1/2L_k}^{(k-1)} (x) - \widetilde{G}_{k-1}^+ \Big( -\displaystyle\frac{1}{2L_k}, x \Big) \bigg| \\
			& \le \displaystyle\sup_{j \in \llbracket - 2n_k/3, 2n_k/3 \rrbracket} \bigg| \widetilde{\gamma}_{-1/2L_k}^{(k-1)} \Big( \displaystyle\frac{j}{n_k} \Big) - \widetilde{G}_{k-1}^+ \Big( -\displaystyle\frac{1}{2L_k}, \displaystyle\frac{j}{n_k} \Big) \bigg| + \displaystyle\frac{10B}{n_k},
		\end{aligned} 
	\end{flalign} 
	
	\noindent where in the statement we used \eqref{figik}, and in the second we used the facts that $\big[ \widetilde{\gamma}_{-1/2L_k}^{(k-1)} \big]_1 \le 4B$ and $\big[ \widetilde{G}_{k-1}^+ \big]_1 \le 5B$ (where the former holds by \eqref{xjkgamma} and the latter by \eqref{gk1b}). We also have 
	\begin{flalign}
		\label{yz2} 
		\begin{aligned}
			& \displaystyle\sup_{j \in \llbracket - 2n_k/3, 2n_k/3 \rrbracket} \bigg| \widetilde{\gamma}_{-1/2L_k}^{(k-1)} \Big( \displaystyle\frac{j}{n_k} \Big) - \widetilde{G}_{k-1}^+ \Big( -\displaystyle\frac{1}{2L_k}, \displaystyle\frac{j}{n_k} \Big) \bigg| \\
			& \qquad \le \displaystyle\sup_{j \in \llbracket -2n_k/3, 2n_k/3 \rrbracket} \bigg| x_{j+n_k/2}^{(k)} \Big( -\displaystyle\frac{1}{2L_k} \Big) - 4 G_{k-1}^+ \Big( -\displaystyle\frac{1}{8L_k}, \displaystyle\frac{j}{n_{k-1}} \Big) \bigg| + n^{-9/10} \\
			& \qquad = 4 \cdot \displaystyle\sup_{j \in \llbracket - n_{k-1}/6, n_{k-1}/6 \rrbracket} \bigg| x_{j+n_{k-1}/2}^{(k-1)} \Big( -\displaystyle\frac{1}{8L_k} \Big) - G_{k-1}^+ \Big( -\displaystyle\frac{1}{8L_k}, \displaystyle\frac{j}{n_{k-1}} \Big) \bigg| + n^{-9/10}.
		\end{aligned} 
	\end{flalign}
	
	\noindent where in the first statement we used \eqref{gk1sx} and \eqref{xjkgamma} (with the facts that $n_{k-1} = 4 n_k$ and that we are restricting to $\Omega_k^{(1)}$), and in the second we used the facts that $x_{j+n_k/2}^{(k)} (s) = 4 x_{j+n_{k-1}/2}^{(k-1)} ( s / 4 )$ (which holds by \eqref{xks} and the facts that $\omega_{k-1} = 4 \omega_k$ and $n_{k-1} = 4n_k$). Next, since we are restricting to the event $\Omega_{k-1}^{(2)} \subseteq \Omega_{k-1} \subseteq \Omega_k^{(1)}$, we have by \eqref{omegak42} (with the fact that $(-1 / 8L_k, jn_{k-1}^{-1}) \notin \frac{1}{4} \cdot \mathfrak{R}_k$, by \eqref{rkrk}), we have 
	\begin{flalign}
		\label{yz3}
		\displaystyle\sup_{|j| \le n_{k-1}/6} \bigg| x_{j+n_{k-1}/2}^{(k-1)} \Big( -\displaystyle\frac{1}{8L_k} \Big) - G_{k-1}^+ \Big( -\displaystyle\frac{1}{8L_k}, \displaystyle\frac{j}{n_{k-1}} \Big) \bigg| \le \displaystyle\frac{\vartheta_{k-1}}{5}.
	\end{flalign}
	
	\noindent Combining \eqref{yz1}, \eqref{yz2}, \eqref{yz3}, and the fact that $4 \vartheta_{k-1} / 5 + n^{-9/10} + 10Bn_k^{-1} \le 4 \vartheta_{k-1} / 5  + 2n^{-9/10} \le \vartheta_{k-1}$ for sufficiently large $n$ (due to \eqref{k0} and the fact that $\vartheta_{k-1} \ge n^{-13/15}$ by \eqref{deltaomegal}) we deduce that \eqref{omega20} holds; this establishes the lemma. 
\end{proof}

\subsection{Proof of \Cref{omega4}} 

\label{ProofOmega4}

In this section we establish \Cref{omega4}, following the ideas outlined in the beginning of \Cref{Proof0Boundary} and \Cref{LOCALHJ}. In particular, both statements of the lemma will follow from a suitable application of \Cref{gh}.

\begin{proof}[Proof of of \Cref{omega41} of \Cref{omega4}]
	
	For each index $\pm \in \{ +, - \}$, define the $n$-tuples $\bm{u}^{\pm} = (u_1^{\pm}, u_2^{\pm}, \ldots , u_n^{\pm}) \in \mathbb{W}_n$ and $\bm{v}^{\pm} = (v_1^{\pm}, v_2^{\pm}, \ldots , v_n^{\pm}) \in \mathbb{W}_n$, and functions $f^{\pm}, g^{\pm} : [-\xi, \xi] \rightarrow \mathbb{R}$ by setting 
	\begin{flalign*}
	u_j^{\pm}  = G ( -\xi; jn^{-1}) \pm \delta; \quad v_j^{\pm} = G (\xi; jn^{-1}) \pm \delta; \quad f^{\pm} (s) = G (s; 1) \pm \delta; \quad g^{\pm} (s) = G(s; 0) \pm \delta,
	\end{flalign*} 

	\noindent for each $j \in \llbracket 1, n \rrbracket$ and $s \in [-\xi, \xi]$. Then sample the line ensemble $\bm{x}^{\pm} = (x_1^{\pm}, x_2^{\pm}, \ldots , x_n^{\pm}) \in \llbracket 1, n \rrbracket \times \mathcal{C} \big( [-\xi, \xi] \big)$ from the measure $\mathsf{Q}_{f^{\pm}; g^{\pm}}^{\bm{u}^{\pm}; \bm{v}^{\pm}} (n^{-1})$. 
	
	Apply \Cref{gh} (translated by $(-\xi, 0)$) to $\bm{x}^{\pm}$, with $(L, B; G)$ there given by $( 1 / 2\xi, 2B; G \pm \delta)$ here; \Cref{fgr} is then verified by \Cref{derivativegxpfl} and \Cref{g00}. This yields a constant $c_1 = c_1 (B) > 0$ such that for each $\pm \in \{ +, - \}$ we have 
	\begin{flalign}
		\label{x00} 
		\mathbb{P} \Bigg[ \displaystyle\sup_{s \in [-\xi, \xi]} \bigg( \displaystyle\max_{j \in \llbracket 1, n \rrbracket} \big| x_j^{\pm} (s) - G^{\pm} (s, jn^{-1}) \big| \bigg) > c_1^{-1} n^{-23/24} \Bigg] \le c_1^{-1} e^{-c_1 (\log n)^2}.
	\end{flalign}	

	\noindent Next, by \eqref{gnuv}, we have that $u_j^- \le u_j \le u_j^+$ and $v_j^- \le v_j \le v_j^+$ for each $j \in \llbracket 1, n \rrbracket$, and that $f^- \le f \le f^+$ and $g^- \le g \le g^+$. Therefore, by height monotonicity \Cref{monotoneheight}, there exists a coupling between $(\bm{x}, \bm{x}^+)$ and between $(\bm{x}^-, \bm{x})$ such that the following holds. For each $(s, j) \in [-\xi, \xi] \times \llbracket 1, n \rrbracket$, we have $x_j (s) \le x_j^+ (s)$ under the former and $x_j^- (s) \le x_j (s)$ under the latter. Together with \eqref{x00}, a union bound, and the fact that $G^{\pm} (t, x) = G (t, x) \pm \delta$, it follows that 
	\begin{flalign}
		\label{e0} 
		\mathbb{P} \big[ \mathscr{E}_0^{\complement} \big] \le 2c_1^{-1} e^{-c_1 (\log n)^2}, \quad \text{where} \quad \mathscr{E}_0 = \Bigg\{ \displaystyle\sup_{|s| \le \xi} \bigg( \displaystyle\max_{j \in \llbracket 1, n \rrbracket} \big| x_j (s) - G(s, jn^{-1}) \big| \bigg) \le \delta + n^{-23/24} \Bigg\}. 
	\end{flalign}
	
	\noindent Then, observe that 
	\begin{flalign*} 
		\mathscr{E}_0 & \subseteq \Bigg\{ \displaystyle\sup_{|s-s_0 n^{-1/3}| \le \omega_0 / 2L_0} \bigg( \displaystyle\max_{-n_0/2 < |j| \le n_0/2} \big| x_{j_0+j} (s) - G(s, jn^{-1} + j_0 n^{-1}) \big| \bigg) \le \delta + n^{-23/24} \Bigg\} \\
		& = \Bigg\{ \displaystyle\sup_{(j/n_0, s) \in \overline{\mathfrak{R}}_0} \big| x_{j+n_0/2}^{(0)} (s) - G_0 (s, jn_0^{-1}) \big| \le \displaystyle\frac{1}{4} \cdot (\delta + n^{-23/24}) \Bigg\} \subseteq \Omega_0,	
	\end{flalign*} 
	
	\noindent where in the first statement we restricted the range of $(j, s)$ in the definition \eqref{e0} of $\mathscr{E}_0$ (using the facts that $n / 3 \le j_0 \le 2n / 3$; that $n_0 = n /4$ by \eqref{deltaomegal}; that $| s n^{-1/3}_0| \le \xi / 2$; and that $\omega_0 / 2L_0 < \delta_0^{1/2} < (2B)^{-1} \le \xi / 2$ for sufficiently small $\delta_0$, again by \eqref{deltaomegal}); in the second we used the fact that $G_0 (s, jn_0^{-1}) = 4 G (s n^{-1/3} + \omega_0 x, j n^{-1} + j_0 n^{-1})$ (by \eqref{gksx2} and the fact that $\omega_0 = 1 / 4$), the fact that $x_{j+n_0/2}^{(0)} (s) = 4 x_{j+j_0} (s n^{-1/3}_0 + \omega_0 s)$ (by \eqref{xks} and the fact that $\omega_0 = 1 / 4$) and the definition \eqref{rkrk} of $\mathfrak{R}_k$; and in the third we used \eqref{gksx}, the definition \eqref{omega0} of $\Omega_0$, and the fact that $\delta_0^{3/4} > \delta + n^{-23/24}$ for sufficiently large $n$ and small $\delta$ (by \eqref{deltaomegal}). Together with \eqref{e0}, this yields the first statement of the lemma.
	\end{proof} 

	\begin{proof}[Proof of \Cref{omega42} of \Cref{omega4}]
	
	 Throughout this proof, we condition on $\mathcal{F}_{\ext} (\mathsf{R}_k)$ and restrict to the event $\Omega_k^{(1)}$ (which, as stated below \eqref{1omega1}, is measurable with respect to $\mathcal{F}_{\ext} (\mathsf{R}_k)$); all of the probabilities in this proof will be with respect to this conditioning and restriction. Since  $\Omega_k^{(3)} \subseteq \Omega_k^{(2)}$, it then suffices to show that $\Omega_k^{(3)}$ holds with probability at least $1 - n^{-20}$. To this end, for each index $\pm \in \{ +, - \}$, define the sequences $\bm{u}^{(k; \pm)}, \bm{v}^{(k; \pm)}, \bm{u}^{(k)}, \bm{v}^{(k)} \in \mathbb{W}_{n_k}$ and functions $f_k^{\pm}, g_k^{\pm}, f_k, g_k  : [ -1 / 2L_k, 1 / 2L_k ] \rightarrow \mathbb{R}$ by for each $j \in \llbracket 1 - n_k / 2, n_k / 2 \rrbracket$ and $s \in [ -1 / 2L_k, 1 / 2L_k ]$ setting 
	\begin{flalign}
		\label{uvfg} 
		\begin{aligned}
			& u_{j+n_k/2}^{(k; \pm)} = G_k^{\pm} \Big( -\displaystyle\frac{1}{2L_k}; \displaystyle\frac{j}{n_k} \Big) \pm n^{-9/10}; \qquad	 \bm{u}^{(k)} = \bm{x}^{(k)} \Big( - \displaystyle\frac{1}{2L_k} \Big); \\ 			
			& v_{j+n_k/2}^{(k; \pm)} = G_k^{\pm} \Big( -\displaystyle\frac{1}{2L_k}; \displaystyle\frac{j}{n_k} \Big) \pm n^{-9/10}; \qquad \bm{v}^{(k)} = \bm{x}^{(k)} \Big( \displaystyle\frac{1}{2L_k} \Big); \\
			&f_k^{\pm} (s) = G_k^{\pm} \Big( s, \displaystyle\frac{1}{2} \Big) \pm n^{-9/10}; \qquad \qquad \qquad g_k^{\pm} (s) = G_k^{\pm} \Big( s, -\displaystyle\frac{1}{2} \Big)\pm n^{-9/10};\\
			&f_k (s) = x_{n_k+1}^{(k)} (s); \qquad \qquad \qquad \qquad \qquad \quad g_k (s) = x_0^{(k)} (s).
		\end{aligned} 
	\end{flalign} 
	For each index $\pm \in \{ +, - \}$, denote the line ensemble $\bm{x}^{(k; \pm)} = \big( x_1^{(k; \pm)}, x_2^{(k; \pm)}, \ldots , x_{n_k}^{(k; \pm)} \big) \in \llbracket 1, n_k \rrbracket \times \mathcal{C} \big( [-1 / 2L_k, 1 / 2L_k] \big)$ sampled from the measure $\mathsf{Q}_{f_k^{\pm}; g_k^{\pm}}^{\bm{u}^{(k;\pm)}; \bm{v}^{(k;\pm)}} (n_k^{-1})$. Further denote the events
	\begin{flalign*}
		\mathscr{E}_k^{\pm} = \Bigg\{ \displaystyle\max_{j \in \llbracket 1 -n_k/2, n_k/2 \rrbracket} \Big| x_{j+n_k/2}^{(k; \pm)} (s) - \big( G_k^{\pm} ( s, jn_k^{-1}) \pm n^{-9/10} \big) \Big| \le n^{-9/10} \Bigg\}; \qquad \mathscr{E}_k = \mathscr{E}_k^- \cap \mathscr{E}_k^+.
	\end{flalign*}
	
	\noindent We next show $\mathscr{E}_k$ is likely, by applying \Cref{gh} (translated by $(-1 / 2L_k, -1 / 2 )$) with the  $(n; G; \bm{u}; \bm{v})$ there equal to $(n_k; G_k^{\pm} \pm n^{-9/10}; \bm{u}^{(k;\pm)}; \bm{v}^{(k;\pm)})$ here. In this way, the $(B_0; L)$ there are equal to $(2C_0; L_k)$ here. To verify \Cref{fgr}, the first statement in \eqref{glg} there follows from \eqref{k0}, the fact that $L_k \ge 1$ is increasing in $k$, and the estimate $e^{\sqrt{\log n}} \le n^{1/20000}$; the second and third statements in \eqref{glg} are verified by \Cref{gk}; and the bound \eqref{uvg} there is verified by \eqref{uvfg}. This yields a constant $c_2 = c_2 (B) > 0$ such that 
	\begin{flalign}
		\label{probabilityek} 
		\mathbb{P} \big[ \mathscr{E}_k^{\complement} \big] \le c_2^{-1} e^{-c_2 (\log n_k)^2} \le n^{-20},		
	\end{flalign} 
	
	\noindent where we have implicitly used the fact that $\log n_k \ge (\log n) / 2$ (as $n_k \ge n_{K_0} \ge n^{1/2}$ by \eqref{k0} and the fact that $n_k = \omega_k n$ is decreasing in $k$). 
	
	Now, observe from \eqref{xks}, \Cref{scale}, and the Brownian Gibbs property that the family $\bm{x}^{(k)}$ of non-intersecting Brownian bridges has law $\mathsf{Q}_{f_k; g_k}^{\bm{u}^{(k)}; \bm{v}^{(k)}} (n_k^{-1})$. We claim that it is possible to couple the three famlies of non-intersecting Brownian bridges $\big( \bm{x}^{(k; -)}, \bm{x}^{(k)}, \bm{x}^{(k; +)} \big)$ so that 
	\begin{flalign}
		\label{xjxj3}
		x_j^{(k;-)} (s) \le x_j^{(k)} (s) \le x_j^{(k;+)} (s), \qquad \text{almost surely, for each $(j, s) \in \llbracket 1, n_k \rrbracket \times \Big[ -\displaystyle\frac{1}{2L_k}, \displaystyle\frac{1}{2L_k} \Big]$}.
	\end{flalign}	
	
	\noindent See the right side of \Cref{f:sandwich} for a depiction. To this end, it suffices by height monotonicity (\Cref{monotoneheight}) to show that
	\begin{flalign}
		\label{uuvvfg} 
		\begin{aligned} 
			\bm{u}^{(k;-)} \le \bm{x}^{(k)} & \Big( - \displaystyle\frac{1}{2L_k} \Big) \le \bm{u}^{(k;+)}; \qquad \bm{v}^{(k;-)} \le \bm{x}^{(k)} \Big( \displaystyle\frac{1}{2L_k} \Big) \le \bm{v}^{(k;+)}; \\
			& f_k^- \le x_{n_k + 1}^{(k)} \le f_k^+; \qquad g_k^- \le x_0^{(k)} \le g_k^+.
		\end{aligned}
	\end{flalign}
	
	\noindent To do this, observe for any $j \in \llbracket 1 - n_k / 2, n_k / 2 \rrbracket$ that 
	\begin{flalign*}
		u_{j+n_k/2}^{(k; -)} (0) \le \mathfrak{y}_0^{(k)} (jn_k^{-1}) - n^{-9/10} & = \widetilde{\gamma}_{-1/2L_k}^{(k-1)} (j n_k^{-1}) - n^{-9/10} \le x_{j+n_k/2}^{(k)} \Big( -\displaystyle\frac{1}{2L_k} \Big),
	\end{flalign*} 
	
	\noindent where in the first statement we applied \eqref{uvfg} and the second statement of \Cref{gk}; in the second we applied \eqref{figik}; and in the third we applied \eqref{xjkgamma}. This shows that $\bm{u}^{(k;-)} \le \bm{x}^{(k)} ( -1 / 2L_k )$. By similar reasoning we also have $\bm{x}^{(k)} (-1 / 2L_k) \le \bm{u}^{(k;+)}$, establishing the first statement of \eqref{uuvvfg}; the proof of the second is entirely analogous and is thus omitted. To establish the third, observe that 
	\begin{flalign*}
		x_{n_k+1}^{(k)} (s) =4x_{5n_{k-1}/8+1}^{(k-1)} \Big( \displaystyle\frac{s}{4} \Big) & \le 4 G_{k-1}^+ \Big( \displaystyle\frac{s}{4}, \displaystyle\frac{1}{8} + \displaystyle\frac{1}{n_{k-1}} \Big) + \displaystyle\frac{4 \vartheta_{k-1}}{5} \\
		& \le 4 G_{k-1}^+ \Big( \displaystyle\frac{s}{4}, \displaystyle\frac{1}{8} \Big) + \displaystyle\frac{4 \vartheta_{k-1}}{5} + \displaystyle\frac{2 C_0}{n_k} \\
		& \le 4 G_{k-1}^+ \Big( \displaystyle\frac{s}{4}, \displaystyle\frac{1}{8} \Big) + \vartheta_{k-1}  \le \widetilde{G}_{k-1}^+ \Big( s, \displaystyle\frac{1}{2} \Big) + \vartheta_{k-1} \le G_k^+ \Big(s, \displaystyle\frac{1}{2} \Big) \leq f_k^+ (s).
	\end{flalign*}
	
	\noindent Here, in the first statement we used that $x_{n_k+1}^{(k)} (s) = 4 x_{n_{k-1}/2+n_k/2+1}^{(k-1)} ( s / 4) = 4 x_{5n_{k-1}/8+1}^{(k-1)} ( s / 4 )$ (which follows from \eqref{xks}, with the equalities $\omega_{k-1}= 4 \omega_k$ and $n_{k-1} = 4n_k$); in the second we used \eqref{omegak42}, the fact that $( s / 4, 1 / 8 + n_{k-1}^{-1} ) \in \mathfrak{R}_{k-1} \setminus \frac{1}{4} \cdot \mathfrak{R}_k$ for $s \in [ -1 / 2L_k, 1 / 2L_k ]$ (by \eqref{rkrk}), and our restriction to the event $ \Omega_k^{(1)}\subseteq \Omega_{k-1}\subseteq \Omega_{k-1}^{(2)}$; in the third we used the third statement of \Cref{gk}; in the fourth we used the bound $\vartheta_{k-1} \ge n^{-13/15} \ge 5 C_0 n_k^{-1}$ (which holds for sufficiently large $n$ by \eqref{deltaomegal} and \eqref{k0}); in the fifth we used \eqref{gk1sx}; in the sixth we used the sixth statement of \Cref{gk}; and in the seventh we used \eqref{uvfg}. Similar reasoning indicates that $x_{n_k + 1}^{(k)} \ge f_k^-$, which yields the third statement of \eqref{uuvvfg}; the proof of the fourth is entirely analogous. This verifies \eqref{uuvvfg} and thus \eqref{xjxj3}. 
	
	In view of \eqref{xjxj3} and \eqref{probabilityek}, we have that 
	\begin{flalign*} 
		\mathbb{P} \Bigg[ \bigcap_{j = 1-n_k/2}^{n_k/2} \bigcap_{|s| \le 1/2L_k} \big\{ G_k^- ( s, jn_k^{-1}) - 2n^{-9/10} \le x_{j+n_k/2}^{(k)} (s) \le G_k^+ (s, jn_k^{-1}& ) + 2 n^{-9/10} \big\} \Bigg] \\
		& \ge \mathbb{P} [ \mathscr{E}_k ] \ge 1 - n^{-20}.
	\end{flalign*} 
	
	\noindent Since by \eqref{deltaomegal} and the fifth statement of \Cref{gk} we have $\big| G_k^- (s, x) - G_k^+ (s, x) \big| \le C_0 e^{-c_0 L_k^{1/8}} \le \varsigma_k / 5$ whenever $|x| \le 1/4$, it follows that 
	\begin{flalign}
		\label{kxjnk}
		\mathbb{P} \Bigg[ \bigcap_{j = -n_k/4}^{n_k/4} \bigcap_{|s| \le 1/2L_k} \Big\{ \big| x_{j+n_k/2}^{(k)} (s) - G_k^+ ( s, jn_k^{-1}) \big| \le \displaystyle\frac{\varsigma_k}{5} + 4n^{-9/10} \Big\} \Bigg] \ge 1 - n^{-20}.
	\end{flalign}
	
	\noindent Due to the bound $\vartheta_k / 5 = (n^{-13/15} + \varsigma_k) / 5 \ge 4n^{-9/10} + \varsigma_k / 5$, the event on the left side of \eqref{kxjnk} is contained in $\Omega_k^{(3)}$; this verifies  \Cref{omega4}.
\end{proof}

\chapter{Airy Statistics} 

\label{STATISTICSBRIDGES}

\section{Gap Convergence to the Airy Point Process}

\label{ProofDifference}

In this section we establish the Airy gaps \Cref{xdifferenceconverge}. We first stochastically bound the gaps of $\bm{\mathsf{x}}$ below by that of an Airy point process in \Cref{DifferenceLower}; then, after recalling a result for edge statistics of Dyson Brownian motion in \Cref{MotionEdge}, we provide the complementary stochastic upper bound for the gaps in \Cref{DifferenceUpper}.

\subsection{Gap Lower Bound} 

\label{DifferenceLower}

\begin{figure}
	\center
	\includegraphics[width=1\textwidth, trim= 1cm 0.5cm 1cm 0.8cm, clip]{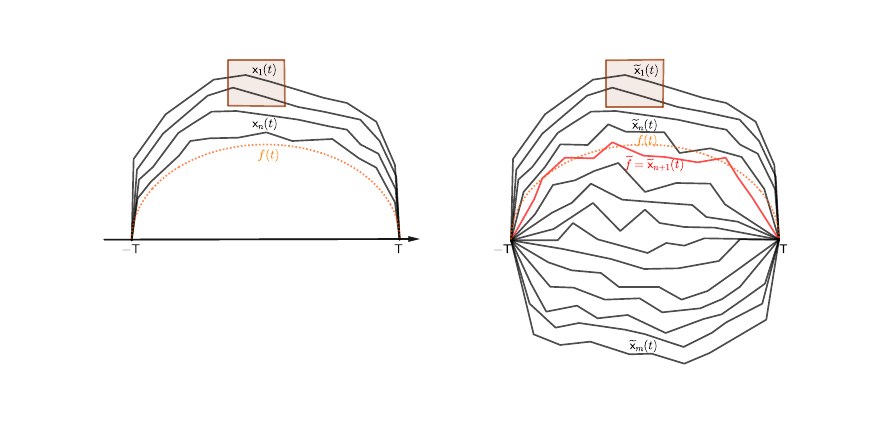}
	
	\caption{Shown on the left is a depiction of the setup for \Cref{f1x2converge}, whose proof follows from coupling $\bm{\mathsf x}$ with the Brownian watermelon $\widetilde{\bm{\mathsf x}}$ as shown to the right. \label{f:Gap_Lower_Bound}}
	
\end{figure}

In this section we establish a lower bound on the gaps considered in \Cref{xdifferenceconverge}. This will follow from a suitable application of gap monotonicity \Cref{monotonedifference}, together with the following result indicating convergence to the Airy line ensemble for non-intersecting Brownian bridges whose lower boundary is given by a rescaled semicircle; see the left side of \Cref{f:Gap_Lower_Bound} for a depiction. In what follows, we recall the classical locations $\gamma_{\semci; n}$ of the semicircle distribution from \eqref{gammaj}, which in the next lemma will be related to the stretching factor for the semicircle (when we later apply the lemma, it will be close to $2$).

\begin{lem}
	
	\label{f1x2converge} 
	
	Let $\bm{\sigma} = (\sigma_1, \sigma_2, \ldots)$ and $\bm{\mathsf{T}} = (\mathsf{T}_1, \mathsf{T}_2, \ldots)$ be two sequences of positive real numbers, so that $\lim_{n \rightarrow \infty} \sigma_n = 1$ and $ \mathsf{T}_n \in [2n^{1/3}, n^{1/2}]$, for each integer $n \ge 64$. For any integer $n \ge 64$, define $f = f_n : [-\mathsf{T}_n, \mathsf{T}_n] \rightarrow \mathbb{R}$ by setting 
	\begin{flalign}
		\label{fn}
		f_n (t) = \sigma_n \mathsf{T}_n \Big( \displaystyle\frac{\mathsf{T}_n^2 - t^2}{2} \Big)^{1/2} \cdot \gamma_{\semci; \lfloor \sigma_n^2 \mathsf{T}_n^3 \rfloor} (n + 1), \qquad \text{for each $t \in [-\mathsf{T}_n, \mathsf{T}_n]$}.
	\end{flalign} 
	
	\noindent Sample non-intersecting Brownian bridges $\bm{\mathsf{x}} = (\mathsf{x}_1, \mathsf{x}_2, \ldots , \mathsf{x}_n) \in \llbracket 1, n \rrbracket \times \mathcal{C} \big( [-\mathsf{T}_n, \mathsf{T}_n] \big)$ under the measure $\mathsf{Q}_f^{\bm{0}_n; \bm{0}_n}$, and define
	\begin{flalign*} 
		\bm{\mathsf{X}}^n = (\mathsf{X}_1^n, \mathsf{X}_2^n, \ldots , \mathsf{X}_n^n) \in \llbracket 1, n \rrbracket \times \mathcal{C} \big( [-\mathsf{T}, \mathsf{T}] \big), \quad \text{where} \quad \mathsf{X}_j^n (t) =  \mathsf{x}_i (t) - 2^{1/2} \sigma_n\mathsf{T}_n^2. 
	\end{flalign*} 
	
	\noindent Then, $\bm{\mathsf{X}}^n$ converges to $\bm{\mathcal{S}}$ on compact subsets of $\mathbb{Z}_{\ge 1} \times \mathbb{R}$, as $n$ tends to $\infty$.  
\end{lem}

\begin{proof}
	
	Throughout this proof, we abbreviate $\mathsf{T} = \mathsf{T}_n$ and $\sigma = \sigma_n$, and we set $m = m_n = \lfloor \sigma^2 \mathsf{T}^3 \rfloor\geq 2n$. To establish the lemma, we will first use \Cref{estimatexj} to approximate the ensemble $\bm{\mathsf{x}}$ by the first $n$ curves of a watermelon $\widetilde{\bm{\mathsf{x}}}$ with $m$ bridges (where its $(n+1)$-th curve closely mimics the shape of $f(t)$); then, we will apply the first part of \Cref{convergea} (with $(n, a,b,j)$ there equal to $(m, -\mathsf T, \mathsf T, n+1)$ here) to show the latter converges to $\bm{\mathcal{S}}$, implying the same for $\bm{\mathsf{x}}$. 
	
	To this end, sample $m$ non-intersecting Brownian bridges $\widetilde{\bm{\mathsf{x}}} = (\widetilde{\mathsf{x}}_1, \widetilde{\mathsf{x}}_2, \ldots , \widetilde{\mathsf{x}}_{m}) \in \llbracket 1, m \rrbracket \times \mathcal{C} \big( [-\mathsf{T}, \mathsf{T}] \big)$ from the measure $\mathsf{Q}^{\bm{0}_m; \bm{0}_m}$. Then, \Cref{estimatexj} gives for sufficiently large $n$ that
	\begin{flalign*}
		\mathbb{P} \Bigg[ \displaystyle\sup_{t \in [-\mathsf{T}, \mathsf{T}]} \big| \widetilde{\mathsf{x}}_{n+1} (t) - \sigma \mathsf{T} \Big( \displaystyle\frac{\mathsf{T}^2 - t^2}{2} \Big)^{1/2} \cdot \gamma_{\semci; m} (n+1) \big| \ge n^{-1/6} \Bigg] \le n^{-5}.
	\end{flalign*}
	
	\noindent So, defining the (random) function $\widetilde{f} = \widetilde{f}_n : [-\mathsf{T}, \mathsf{T}] \rightarrow \mathbb{R}$ by $\widetilde{f}(t) = \widetilde{\mathsf{x}}_{n+1} (t)$, we have
	\begin{flalign}
		\label{gn} 
		\mathbb{P} [\mathscr{E}_n] \le n^{-5}, \qquad \text{where} \qquad \mathscr{E}_n = \Bigg\{ \displaystyle\sup_{t \in [-\mathsf{T}, \mathsf{T}]} \big| \widetilde{f}(t) - f(t) \big| \ge n^{-1/6} \Bigg\}.
	\end{flalign}
	
	Now condition on $\widetilde{f} = \widetilde{\mathsf{x}}_{n+1}$. Then, the law of the first $n$ curves $\big( \widetilde{\mathsf{x}}_j (t) \big) \in \llbracket 1, n \rrbracket \times \mathcal{C} \big( [-\mathsf{T}, \mathsf{T}] \big)$ of $\widetilde{\bm{x}}$ is given by $\mathsf{Q}_{\tilde{f}}^{\bm{0}_n; \bm{0}_n}$. By the first part of \Cref{uvv}, it follows that we can couple $\bm{\mathsf{x}}$ and $\widetilde{\bm{\mathsf{x}}}$ such that 
	\begin{flalign}
		\label{xjxj} 
		\displaystyle\max_{j \in \llbracket 1, n \rrbracket} \displaystyle\sup_{t \in [-\mathsf{T}, \mathsf{T}]} \big| \mathsf{x}_j (t) - \widetilde{\mathsf{x}}_j (t) \big| \le n^{-1/6}, \qquad \text{on the event $\mathscr{E}_n^{\complement}$}.
	\end{flalign}
	
	\noindent Moreover, by \Cref{convergea} (applied with the $(n, T, \sigma)$ there equal to $(m, \sigma_n^{-2/3}, \sigma_n^{-1/3}$ here), the ensemble $\widetilde{\bm{\mathsf{X}}}^n = (\widetilde{\mathsf{X}}_1^n, \widetilde{\mathsf{X}}_2^n, \ldots , \widetilde{\mathsf{X}}_n^n) \in \llbracket 1, n \rrbracket \times \mathcal{C} \big( [-\mathsf{T}, \mathsf{T}] \big)$ defined by $\widetilde{\mathsf{X}}_j^n (t) = \sigma_n^{1/3} \big( \widetilde{\mathsf{x}}_j (\sigma_n^{-2/3} t) - 2^{1/2} \sigma_n \mathsf{T}_n^2 \big)$ converges to $\bm{\mathcal{S}}$ on compact subsets of $\mathbb{Z}_{\ge 1} \times \mathbb{R}$, as $n$ tends to $\infty$. This, together with \eqref{gn}, \eqref{xjxj}, and the fact that $\lim_{n \rightarrow \infty} \sigma_n = 1$, implies the same for $\bm{\mathsf{X}}^n$, thus implying the lemma. 						
\end{proof}

Using \Cref{f1x2converge}, we can lower bound the gaps of the bridges from \Cref{xfdelta2}.

\begin{prop}
	
	\label{xjxj11}
	
	Adopt \Cref{xfdelta2} (recalling the ensemble $\bm{\mathsf{x}} = \bm{\mathsf{x}}^n$ introduced there). Fix an integer $k \ge 1$; a real number $t \in \mathbb{R}$; and nonnegative real numbers $r_1, r_2, \ldots , r_k \ge 0$. Then, recalling the Airy point process $\bm{\mathfrak{a}} = (\mathfrak{a}_1, \mathfrak{a}_2, \ldots )$ from \Cref{a0}, we have
	\begin{flalign*} 
		\displaystyle\liminf_{n \rightarrow \infty} \mathbb{P} \Bigg[ \bigcap_{j=1}^k \Big\{ 2^{1/2} \big( \mathsf{x}_j (t) - \mathsf{x}_{j+1} (t) \big) \ge r_j \Big\} \Bigg] \ge \mathbb{P} \Bigg[ \bigcap_{j=1}^k \{ \mathfrak{a}_j - \mathfrak{a}_{j+1} \ge r_j \} \Bigg].
	\end{flalign*}
	
\end{prop}

\begin{proof}
	
	Set $\sigma_n = 1 + \delta_n^{1/2}$, and denote $m_n = \lfloor \sigma_n^2 \mathsf{T}_n^3 \rfloor$. By the first part of \Cref{gammaderivative} (and the fact that $n \mathsf{T}^{-3} \le \delta^3$ by \Cref{xfdelta2}), there exists a constant $C > 1$ such that 
	\begin{flalign*} 
		0 \le 2 - \gamma_{\semci; m_n} (n+1) \le C \delta_n^2, \qquad \text{so} \qquad \sigma_n \cdot \gamma_{\semci; m_n} (n+1) \ge (1 + \delta_n^{1/2}) (2 - C \delta_n^2) \ge 2 + 2^{1/2}\delta_n,
	\end{flalign*} 
	
	\noindent for sufficiently large $n$. Thus, letting $\widetilde{f}_n$ denote the function $f_n$ from \eqref{fn}, for each $s \in [-\mathsf{T}_n, \mathsf{T}_n]$ we have
	\begin{flalign*}
		\partial_s^2 \widetilde{f}_n (s) = -\displaystyle\frac{\sigma_n \mathsf{T_n}}{(2\mathsf{T}_n^2 - 2s^2)^{1/2}} \bigg( 1 + \displaystyle\frac{s^2}{\mathsf{T}_n^2 - s^2} \bigg) \cdot \gamma_{\semci; m_n} (n+1) \le - \displaystyle\frac{\sigma_n}{2^{1/2}} \cdot \gamma_{\semci; m_n} (n+1) \le -2^{1/2} - \delta_n.
	\end{flalign*}
	
	\noindent In particular, by \eqref{fsfs} it follows for sufficiently large $n$ that
	\begin{flalign}
		\label{fnt22} 
		\partial_s^2 \widetilde{f}_n (s) \le \partial_s^2 h_n (s), \qquad \text{for each $s \in [-\mathsf{T}_n, \mathsf{T}_n]$},
	\end{flalign}
	
	\noindent which will enable us to apply the gap monotonicity \Cref{monotonedifference}. 
	
	To implement this, condition on $\mathcal{F}_{\ext} \big( \llbracket 1, n \rrbracket \times (-\mathsf{T}_n, \mathsf{T}_n) \big)$ and let $\bm{w} = \bm{\mathsf{x}} (\mathsf{T}) \in \mathbb{W}_n$. Identifying $f_n$ with $f_n |_{[-\mathsf{T}_n, \mathsf{T}_n]}$, the law of $\bm{\mathsf{ x}} |_{[-\mathsf{T}_n, \mathsf{T}_n]}$ is then given by $\mathsf{Q}_{f_n}^{\bm{u}; \bm{w}}$. Sample non-intersecting Brownian bridges $\breve{\bm{\mathsf{x}}}=\breve{\bm{\mathsf{x}}}^n = (\breve{\mathsf{x}}_1, \breve{\mathsf{x}}_2, \ldots , \breve{\mathsf{x}}_n) \in \llbracket 1, n \rrbracket \times \mathcal{C} \big( [-\mathsf{T}_n, \mathsf{T}_n] \big)$ and $\widetilde{\bm{\mathsf{x}}}=\widetilde{\bm{\mathsf{x}}}^n = (\widetilde{\mathsf{x}}_1, \widetilde{\mathsf{x}}_2, \ldots , \widetilde{\mathsf{x}}_n) \in \llbracket 1, n \rrbracket \times \mathcal{C} \big( [-\mathsf{T}_n, \mathsf{T}_n] \big)$ from the measures $\mathsf{Q}_{h_n - \delta_n}^{\bm{u}; \bm{w}}$ and $\mathsf{Q}_{\tilde{f}_n}^{\bm{0}_n; \bm{0}_n}$, respectively. Due to the second bound in \eqref{fsfs}, the first part of \Cref{uvv} yields a coupling between $\bm{\mathsf{x}}$ and $\breve{\bm{\mathsf{x}}}$ such that 
	\begin{flalign}
		\label{xjxj23}
		\big| \mathsf{x}_j (t) - \breve{\mathsf{x}}_j (t) \big| \le 2 \delta_n, \qquad \text{for each $(j, t) \in \llbracket 1, n \rrbracket \times [-\mathsf{T}_n, \mathsf{T}_n]$}.
	\end{flalign}
	
	\noindent Moreover, \eqref{fnt22} and gap monotonicity \Cref{monotonedifference} together yield a coupling between $\breve{\bm{\mathsf{x}}}$ and $\widetilde{\bm{\mathsf{x}}}$ so that $\breve{\mathsf{x}}_j (t) - \breve{\mathsf{x}}_{j+1} (t) \ge \widetilde{\mathsf{x}}_j (t) - \widetilde{\mathsf{x}}_{j+1} (t)$, for each $t \in [-\mathsf{T}_n, \mathsf{T}_n]$ and $j \in \llbracket 1, n-1 \rrbracket$. Together with \eqref{xjxj23}, this implies
	\begin{flalign}
		\label{xjxj1r} 
		\mathbb{P} \Bigg[ \bigcap_{j=1}^k \big\{ \mathsf{x}_j (t) - \mathsf{x}_{j+1} (t) \ge 2^{-1/2}r_j\}\Bigg] \ge \mathbb{P} \Bigg[ \bigcap_{j=1}^k \big\{ \widetilde{\mathsf{x}}_j (t) - \widetilde{\mathsf{x}}_{j+1} (t) \ge 2^{-1/2}r_j - 4 \delta_n \big\} \Bigg].
	\end{flalign}
	
	\noindent Moreover, \Cref{f1x2converge}, \Cref{a0}, and the fact that $\lim_{n \rightarrow \infty} \delta_n = 0$ together imply that 
	\begin{flalign*} 
		\displaystyle\lim_{n \rightarrow \infty} \mathbb{P} \Bigg[ \bigcap_{j=1}^k \Big\{ 2^{1/2} \big( \widetilde{\mathsf{x}}_j (t) - \widetilde{\mathsf{x}}_{j+1} (t) \big) \ge r_j - 2^{5/2} \delta_n \Big\} \Bigg] = \mathbb{P} \Bigg[ \bigcap_{j=1}^k \{ \mathfrak{a}_j - \mathfrak{a}_{j+1} \ge r_j \} \Bigg].
	\end{flalign*} 
	
	\noindent This and \eqref{xjxj1r} together yield the lemma. 			
\end{proof}

\subsection{Gap Upper Bound} 

\label{DifferenceUpper}

In this section we establish the following upper bound for the gaps between the bridges in $\bm{\mathsf{x}}$ satisfying \Cref{xfdelta2}, which (together with \Cref{xjxj11}) quickly implies the Airy gaps \Cref{xdifferenceconverge}.

	\begin{figure}[h]
    \centering
    \includegraphics[width=0.7\textwidth, trim=6cm 10.5cm 6cm 11cm, clip]{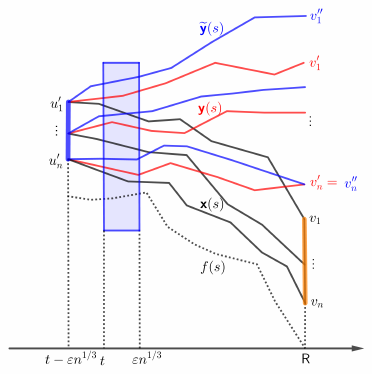}
    \caption{The proof of \Cref{xjxj12} is illustrated above.}
    \label{f:DBM_coupling2}
\end{figure}

\begin{prop}
	
	\label{xjxj12} 
	
	Adopting the notation and assumptions of \Cref{xjxj11}, and also adopting \Cref{xfdelta2}, we have 
	\begin{flalign*}
		\displaystyle\limsup_{n \rightarrow \infty} \mathbb{P} \Bigg[ \bigcap_{j=1}^k \Big\{ 2^{1/2} \big( \mathsf{x}_j (t) - \mathsf{x}_{j+1} (t) \big) \ge r_j \Big\} \Bigg] \le \mathbb{P} \Bigg[ \bigcap_{j=1}^k \{ \mathfrak{a}_j - \mathfrak{a}_{j+1} \ge r_j \} \Bigg]. 
	\end{flalign*}
	
\end{prop}

\begin{proof}[Proof of \Cref{xdifferenceconverge}]
	
	This follows from \Cref{xjxj11} and \Cref{xjxj12}.
\end{proof}

To establish \Cref{xjxj12}, we will use \Cref{yconvergedelta}, to which end we require the following estimate on the constant $\sigma_{\nu;t}$ defined in \Cref{yconvergedelta0}. 

\begin{lem}
	
	\label{sigmanux12}
	
	There exists a constant $C > 1$ such that the following holds. Adopt the notation of \Cref{yconvergedelta0}, and assume that 
	\begin{flalign}
		\label{nuy}
		\nu(dy) = \textbf{\emph{1}}_{y \in [-2^{-7/6} (3\pi)^{2/3}, 0]} \cdot 2^{3/4} \pi^{-1} |y|^{1/2} dy. 
	\end{flalign} 
	
	\noindent Then, for any real number $t \in (0, 1)$, we have $|\sigma_{\nu; t} - 2^{1/2}| \le C t$.
	
\end{lem}

\begin{proof}
	
	Recall the real number $z_0 = z_0 (\nu; t) > 0$ from \eqref{z0tsigma}, and denote $\mathfrak{c} = 2^{-7/6} (3 \pi)^{2/3}$. Let us begin by verifying the approximation $z_0 \approx 2^{-1/2} t^2$. Changing variables from $y$ to $-z_0 w^2$ in the first integral in \eqref{z0tsigma}, we deduce 
	\begin{flalign}
		\label{z0t1} 
		2^{-3/4} \pi t^{-1} z_0^{1/2} = z_0^{1/2} \displaystyle\int_{-\mathfrak{c}}^0 \displaystyle\frac{|y|^{1/2} dy}{(y-z_0)^2} = 2 \displaystyle\int_0^{(\mathfrak{c} / z_0)^{1/2}} \displaystyle\frac{w^2 dw}{(w^2+1)^2} \le 2 \displaystyle\int_0^{\infty} \displaystyle\frac{w^2 dw}{(w^2 + 1)^2} = \displaystyle\frac{\pi}{2},
	\end{flalign}
	
	\noindent from which it follows that $z_0 \le 2^{-1/2} t^2$. Inserting this into \eqref{z0t1}, it follows that 
	\begin{flalign*}
		2^{-3/4} \pi t^{-1} z_0^{1/2} = \displaystyle\frac{\pi}{2} - 2 \displaystyle\int_{(\mathfrak{c} / z_0)^{1/2}}^{\infty} \displaystyle\frac{w^2 dw}{(w^2+1)^2} \ge \displaystyle\frac{\pi}{2} - \displaystyle\int_{(\mathfrak{c} / z_0)^{1/2}}^{\infty} \frac{w^{-2} dw}{2} \ge \displaystyle\frac{\pi}{2} - \frac{z_0^{1/2}}{2} \ge \displaystyle\frac{\pi}{2} - \displaystyle\frac{t}{2},
	\end{flalign*}
	
	\noindent from which we deduce $z_0^{1/2}\geq 2^{-1/4}t -t^2/(2^{1/4}\pi)$ and thus
	\begin{flalign*}
		2^{-1/2} t^2 - z_0 = (2^{-1/4} t - z_0^{1/2}) (2^{-1/4} t + z_0^{1/2}) \le \displaystyle\frac{t^2}{2^{1/4} \pi} \cdot 2^{3/4} t \le t^3,
	\end{flalign*}
	
	\noindent and so
	\begin{flalign}
		\label{z0t}
		0 \le 2^{-1/2} t^2 - z_0 \le t^3.
	\end{flalign}
	
	Next, changing variables from $y$ to $-z_0 w^2$ in the second integral in \eqref{z0tsigma} yields 
	\begin{flalign*}
		2^{-3/4} \pi \sigma_{\nu; t}^{-3} = t^3 \displaystyle\int_{-\mathfrak{c}}^{0} \displaystyle\frac{|y|^{1/2} dy}{(y-z_0)^3} & = 2 t^3 z_0^{-3/2} \displaystyle\int_0^{(\mathfrak{c} / z_0)^{1/2}} \displaystyle\frac{w^2 dw}{(w^2+1)^3} \\
		& = t^3 z_0^{-3/2} \bigg( \displaystyle\frac{\pi}{8} - 2 \displaystyle\int_{(\mathfrak{c} / z_0)^{1/2}}^{\infty} \displaystyle\frac{w^2 dw}{(w^2+1)^3} \bigg).
	\end{flalign*}
	
	\noindent Together with \eqref{z0t} and the fact 
	\begin{flalign*}
		\displaystyle\int_{(\mathfrak{c} / z_0)^{1/2}}^{\infty} \displaystyle\frac{w^2 dw}{(w^2 + 1)^3} \le \displaystyle\int_{(\mathfrak{c} / z_0)^{1/2}}^{\infty} w^{-4} dw \le \displaystyle\frac{z_0^{3/2}}{3} \le \frac{t^3}{3},
	\end{flalign*}
	
	\noindent it follows that 
	\begin{flalign*}
		\sigma_{\nu; t}^3 = \displaystyle\frac{2^{9/4} z_0^{3/2}}{t^3 + \mathcal{O} (t^6)} = 2^{3/2} + \mathcal{O} (t) = 2^{3/2} + \mathcal{O} (t),
	\end{flalign*}
	
	\noindent from which we deduce the lemma.
\end{proof}

We now establish \Cref{xjxj12} using \Cref{yconvergedelta} and gap monotonicity \Cref{monotonedifference} (in a way broadly analogous to the proof of \Cref{kgapc}).

\begin{proof}[Proof of \Cref{xjxj12}]
	We will prove \Cref{xjxj12} if $t\geq 0$; the case when $t\leq 0$ can be proven in the same way by symmetry. Fix some $\varepsilon \in ( 0, 1 / 4 )$; assume that $n$ is sufficiently large so that $t \le \varepsilon n^{1/3} \le \mathsf{R}$; condition on $\mathcal{F}_{\ext} \big( \llbracket 1, n \rrbracket \times (t-\varepsilon n^{1/3}, \mathsf{R}) \big)$; and restrict to the event $\mathscr{F}_n (t-\varepsilon n^{1/3})$ of \eqref{tfn}. Denote the $n$-tuple $\bm{u}' = \bm{\mathsf{x}}^n (t-\varepsilon n^{1/3}) \in \mathbb{W}_n$; the conditional law of $\bm{\mathsf{x}}^n$ is then given by $\mathsf{Q}_f^{\bm{u}'; \bm{v}}$. Further define the process $\bm{\mathsf{z}}^n = (\mathsf{z}_1^n, \mathsf{z}_2^n, \ldots , \mathsf{z}_n^n) \in \llbracket 1, n \rrbracket \times \mathcal{C} \big( [0, n^{2/3} \mathsf{R} + \varepsilon n - t n^{2/3}] \big)$ so that $\bm{\mathsf{z}} (s)$ is obtained by running Dyson Brownian motion for time $s$ with initial data $n^{1/3} \cdot \bm{u}'$. Then define $\bm{\mathsf{y}}^n = (\mathsf{y}_1^n, \mathsf{y}_2^n, \ldots , \mathsf{y}_n^n) \in \llbracket 1, n \rrbracket \times \mathcal{C} \big( [t - \varepsilon n^{1/3}, \mathsf{R}] \big)$ by scaling, so that $\mathsf{y}_j^n (s) = n^{-1/3} \cdot \mathsf{z}_j^n (n^{2/3} s + \varepsilon n - t n^{2/3})$ for each $(j, s) \in \llbracket 1, n \rrbracket \times [t -\varepsilon n^{1/3}, \mathsf{R}]$. Denoting $\bm{v}' = \bm{\mathsf{y}}^n (\mathsf{R}) \in \mathbb{W}_n$, observe from the second part of \Cref{lambdat} with \Cref{scale} that $\bm{y}^n$ has law $\mathsf{Q}^{\bm{u}'; \bm{v}'}$. See \Cref{f:DBM_coupling2}.
	
	Define the $n$-tuple $\bm{v}'' \in \mathbb{W}_n$ by setting $v_j'' = v_j' + (n-j) n$, for each $j \in \llbracket 1, n \rrbracket$. Then, sample non-intersecting Brownian bridges $\widetilde{\bm{\mathsf{y}}}^n = (\widetilde{\mathsf{y}}_1^n, \widetilde{\mathsf{y}}_2^n, \ldots , \widetilde{\mathsf{y}}_n^n) \in \llbracket 1, n \rrbracket \times \mathcal{C} \big( [t - \varepsilon n^{1/3}, \mathsf{R}] \big)$ from the measure $\mathsf{Q}^{\bm{u}'; \bm{v}''}$. Applying the second part of \Cref{uvv} (with the $(\bm{u}, \bm{v}, \widetilde{\bm{v}}, B)$ there equal to the $(\bm{u}', \bm{v}', \bm{v}'', n^2)$ here), yields a coupling between $\mathsf{y}^n$ and $\widetilde{\mathsf{y}}^n$ such that for each $j \in \llbracket 1, n \rrbracket$ we have
	\begin{flalign}
		\label{xnyn2}
		\displaystyle\max_{s \in [t - \varepsilon n^{1/3}, n^{1/3}]}\big| \mathsf{y}_j^n (s) - \widetilde{\mathsf{y}}_j^n (s) \big| \le \frac{(1 + \varepsilon) n^{1/3}}{\mathsf{R}} \cdot n^2 \le 2n^{-2},
	\end{flalign}
	
	\noindent where in the last bound we used the facts that $\varepsilon \le 1/2$ and $\mathsf{R} = n^{20}$. Moreover, since $\bm{\mathsf{x}}^n$ has conditional law $\mathsf{Q}_f^{\bm{u}'; \bm{v}}$, applying gap monotonicity \Cref{monotonedifference} to the measures $\mathsf{Q}_f^{\bm{u}'; \bm{v}}$ and $\mathsf{Q}^{\bm{u}; \bm{v}''}$ (using the fact that $v_j'' - v_{j+1}'' \ge n \ge v_1 - v_n \ge v_j - v_{j+1}$ for each $j \in \llbracket 1, n-1 \rrbracket$, where in the third bound we used \eqref{fvvf}, there is a coupling between $\bm{\mathsf{x}}$ and $\widetilde{\bm{\mathsf{y}}}$ satisfying $\mathsf{x}_j^n (s) - \mathsf{x}_{j+1}^n (s) \le \widetilde{\mathsf{y}}_j^n (s) - \widetilde{\mathsf{y}}_{j+1}^n (s)$, for each $(j, s) \in \llbracket 1, n - 1 \rrbracket \times [t - \varepsilon n^{1/3}, \mathsf{R}]$. Together with \eqref{xnyn2}, this yields a coupling between $\bm{\mathsf{x}}$ and $\bm{\mathsf{y}}$ such that
	\begin{flalign}
		\label{xjxjyjyjn}
		\mathsf{x}_j^n (s) - \mathsf{x}_{j+1}^n (s) \le \mathsf{y}_j^n (s) - \mathsf{y}_{j+1}^n (s) + 4n^{-2}, \qquad \text{for each $(j, s) \in \llbracket 1, n \rrbracket \times [t-\varepsilon n^{1/3}, \mathsf{R}]$}.
	\end{flalign}
	
	Next, define $\nu$ as in \eqref{nuy}, which satisfies \eqref{yndelta1}; also set $\sigma_{\nu;\varepsilon}$ as in \eqref{z0tsigma}. Then, \eqref{tfn} implies on the event $\mathscr{F}_n (t-\varepsilon n^{1/3})$ that the $(u_j' - u_n')$ satisfy \eqref{yndelta} (possibly with a different sequence $\bm{\delta}$ tending to $0$). Hence \Cref{yconvergedelta} (with the $t$ there equal to $\varepsilon$ here), and \Cref{scale}, yield for $n \ge N_0$ sufficiently large (where $N_0 > 1$ only depends on $\varepsilon$, $\bm{\delta}$, $k$, and the $(r_i)$) that 
	\begin{flalign*}
		 & \mathbb{P} \bigg[ \bigcap_{j=1}^k \Big\{ \sigma_{\nu; \varepsilon} \cdot \big( \mathsf{y}_j^n (t) - \mathsf{y}_{j+1}^n (t) \big) \ge r_j \Big\} \bigg] -  \mathbb{P} \bigg[ \bigcap_{j=1}^k \{ \mathfrak{a}_j - \mathfrak{a}_{j+1} \ge r_j \} \bigg]  \\
		& \quad =  \mathbb{P} \bigg[ \bigcap_{j=1}^k \Big\{ \sigma_{\nu; \varepsilon} n^{-1/3} \cdot \big( \mathsf{z}_t^n (\varepsilon n) - \mathsf{z}_{j+1}^n (\varepsilon n) \big) \Big\} \bigg] - \mathbb{P} \Bigg[ \bigcap_{j=1}^k \{ \mathfrak{a}_j - \mathfrak{a}_{j+1} \ge r_j \} \bigg] \le \varepsilon.
	\end{flalign*}
	
	\noindent Together with \eqref{xjxjyjyjn}, this yields for $n \ge N_0$ that
	\begin{flalign*}
		 \mathbb{P} \bigg[ \mathscr{F}(t - \varepsilon n^{1/3}) \cap \bigcap_{j=1}^k \Big\{ 2^{1/2}  \cdot \big( \mathsf{x}_{j+1}^n (t) - \mathsf{x}_j^n (t) & \big) \ge 2^{1/2} \sigma_{\nu; \varepsilon}^{-1} r_j + 4n^{-2} \Big\} \bigg] \\
		 & \le \mathbb{P} \bigg[ \bigcap_{j=1}^k \{ \mathfrak{a}_j - \mathfrak{a}_{j+1} \ge r_j \} \bigg] + \varepsilon.
	\end{flalign*}
	
	\noindent Letting $n$ tend to $\infty$ and using \eqref{probabilityfnt}, it follows that 
	\begin{flalign*}
		\displaystyle\limsup_{n \rightarrow \infty} \mathbb{P} \bigg[ \bigcap_{j=1}^k \Big\{ 2^{1/2}  \cdot \big( \mathsf{x}_{j+1}^n (t) - \mathsf{x}_j^n (t) \big) \ge 2^{1/2} \sigma_{\nu; \varepsilon}^{-1} r_j \Big\} \bigg] \le \mathbb{P} \bigg[ \bigcap_{j=1}^k \{ \mathfrak{a}_j - \mathfrak{a}_{j+1} \ge r_j \} \bigg] + \varepsilon.
	\end{flalign*}
	
	 \noindent Since \Cref{sigmanux12} yields a constant $C > 1$ such that $2^{1/2} \sigma_{\nu; \varepsilon}^{-1} \in [1 - C\varepsilon, 1 + C \varepsilon]$, this implies the proposition upon letting $\varepsilon$ tend to $0$.
\end{proof}

\section{Airy Line Ensembles From Airy Point Processes}

\label{ProofL0}

\subsection{Proof of \Cref{qsf}} 

\label{ProofL00}

\begin{figure}
	\center
	\includegraphics[width=0.7\textwidth, trim=0 0.5cm 0 0.5cm, clip]{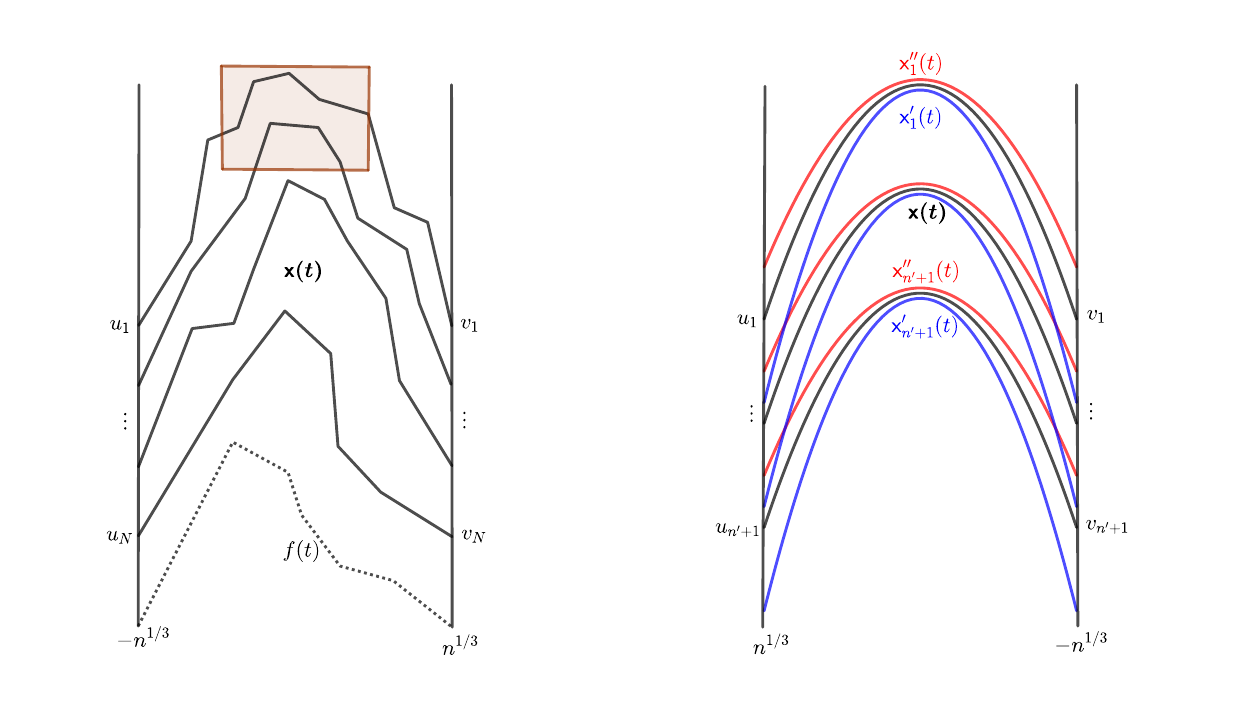}
	
	\caption{Shown to the left is a depiction of \Cref{uvxrconverge}; shown to the right is a depiction of its proof.}
	\label{f:Airy_Line}
\end{figure}

In this section we establish \Cref{qsf}, which will follow as a consequence of the following proposition, to be established in \Cref{ProofxR} below. It states that the edge statistics of a family of $N$ non-intersecting Brownian bridges on a shorter interval $[-n^{1/3}, n^{1/3}]$ (where $n$ is much smaller than $N$), whose boundary data is close to the expected values of the Airy point process (and whose lower boundary is not too irregular), converges to the Airy line ensemble. See the left side of \Cref{f:Airy_Line} for a depiction.

\begin{prop}
	
	\label{uvxrconverge} 
	
	Let $n \ge 1$ denote an integer, and set $N = n^{15}$. Let $\bm{u} = \bm{u}^n \in \mathbb{W}_N$ and $\bm{v} = \bm{v}^n \in \mathbb{W}_N$ denote two $n$-tuples such that
	\begin{flalign}
		\label{xjuv0}
		\big| u_j + 2^{-1/2} n^{2/3} + 2^{-7/6} (3 \pi)^{2/3} j^{2/3} \big| + \big|v_j + 2^{-1/2} n^{2/3} + 2^{-7/6} (3 \pi)^{2/3} j^{2/3} \big| \le (\log n)^3 j^{-1/3},
	\end{flalign}
	
	\noindent for each integer $j \in \llbracket 1, N \rrbracket$. Also let $f = f_n : [-n^{1/3}, n^{1/3}] \rightarrow \mathbb{R}$ denote a function such that 
	\begin{flalign}
		\label{unfs}
		\displaystyle\sup_{|s| \le n^{1/3}} \big| f (s) - u_N \big| \le n^8. 
	\end{flalign}
	
	\noindent Sample $N$ non-intersecting Brownian bridges $\bm{\mathsf{x}}^n = (\mathsf{x}_1^n, \mathsf{x}_2^n, \ldots , \mathsf{x}_N^n) \in \llbracket 1, N \rrbracket \times \mathcal{C} \big( [-n^{1/3}, n^{1/3}] \big)$ from the measure $\mathsf{Q}_f^{\bm{u}; \bm{v}}$. Then, $\bm{\mathsf{x}}^n$ converges to $\bm{\mathcal{S}}$, uniformly on compact subsets of $\mathbb{Z}_{\ge 1} \times \mathbb{R}$, as $n$ tends to $\infty$. 
\end{prop}

\begin{proof}[Proof of \Cref{qsf}]
	
	Let $n \ge 1$ be an integer; set $N = N_n = n^{15}$ and $\mathsf{T} = \mathsf{T}_n = n^{1/3}$; and abbreviate $\mathcal{F}_{\ext}^n = \mathcal{F}_{\ext} \big( \llbracket 1, N \rrbracket \times (-\mathsf{T}, \mathsf{T}) \big)$ (recall \Cref{property}). Define the $\mathcal{F}_{\ext}^n$-measurable random variables 
	\begin{flalign}
		\label{xizeta} 
		\begin{aligned}
			\xi_n & = (2 \mathsf{T})^{-1} \cdot \big( \mathcal{L}_N (\mathsf{T}) - \mathcal{L}_N (-\mathsf{T}) \big); \\
			\zeta_n & = \displaystyle\frac{1}{2} \cdot \big( \mathcal{L}_N (-\mathsf{T}) + \mathcal{L}_N (\mathsf{T}) \big) + 2^{-1/2} n^{2/3} + 2^{-7/6} (3 \pi)^{2/3} N^{2/3}.
		\end{aligned} 
	\end{flalign}
	
	\noindent Also define the family of non-intersecting curves $\bm{\mathsf{L}}^n = (\mathsf{L}_1^n, \mathsf{L}_2^n, \ldots , \mathsf{L}_N^n) \in \llbracket 1, N \rrbracket \times \mathcal{C} \big( [-\mathsf{T}, \mathsf{T}] \big)$; the $N$-tuples $\bm{u} = \bm{u}^n \in \mathbb{W}_N$ and $\bm{v} = \bm{v}^n \in \mathbb{W}_N$; and the function $f = f_n : [-\mathsf{T}, \mathsf{T}] \rightarrow \mathbb{R}$, by setting 
	\begin{flalign}
		\label{ljn} 
		\mathsf{L}_j^n (s) = \mathcal{L}_j (s) - \xi_n s - \zeta_n; \qquad u_j = \mathsf{L}_j^n (-\mathsf{T}); \qquad v_j = \mathsf{L}_j^n (\mathsf{T}); \qquad f(s) = \mathsf{L}_{N+1}^n (s),
	\end{flalign} 
	
	\noindent for each integer $j \ge 1$ and real number $s \in \mathbb{R}$ (so these quantities are still defined for $(j, s) \notin \llbracket 1, N \rrbracket \times [-\mathsf{T}, \mathsf{T}]$). This in particular guarantees that 
	\begin{flalign}
		\label{ttln} 
		\mathsf{L}_N^n (-\mathsf{T}) = -2^{-1/2} n^{2/3} - 2^{-7/6} (3 \pi)^{2/3} N^{2/3} = \mathsf{L}_N^n (\mathsf{T}).
	\end{flalign} 
	
	\noindent Moreover, conditional on $\mathcal{F}_{\ext}^n$, the ensemble $\bm{\mathsf{L}}^n$ is a family of $N$ non-intersecting Brownian bridges sampled from the measure $\mathsf{Q}_f^{\bm{u}; \bm{v}}$ (by \Cref{linear}). To apply \Cref{uvxrconverge} to this ensemble, we must restrict to an event on which its conditions \eqref{xjuv0} and \eqref{unfs} hold.
	
	So, define the $\mathcal{F}_{\ext}^n$-measurable event $\mathscr{E} (n) = \mathscr{E}_1 (n) \cap \mathscr{E}_2 (n) \cap \mathscr{E}_3 (n)$, where
	\begin{flalign*}
		\mathscr{E}_1 (n) & =  \bigcap_{j \in \llbracket 1, N \rrbracket}\Big\{ \big| u_j  + 2^{-1/2} n^{2/3} + 2^{-7/6} (3 \pi)^{2/3} j^{2/3} \big| \le  (\log n)^{5/2} j^{-1/3}\Big\}; \\
		\mathscr{E}_2 (n) & = \bigcap_{j \in \llbracket 1, N \rrbracket} \Big\{ \big| v_j + 2^{-1/2} n^{2/3} + 2^{-7/6} (3 \pi)^{2/3} \big| \le  (\log n)^{5/2} j^{-1/3} \Big\}; \\
		\mathscr{E}_3 (n) & = \Bigg\{ \displaystyle\sup_{|t| \le \mathsf{T}} \big| f (t) - u_N \big| \le n^8 \Bigg\}.
	\end{flalign*}
	
	\noindent Let us show that $\lim_{n \rightarrow \infty} \mathbb{P} \big[ \mathscr{E} (n) \big] = 1$; it suffices by a union bound to show that 
	\begin{flalign}
		\label{e123} 
		\displaystyle\lim_{n \rightarrow \infty} \mathbb{P} \big[ \mathscr{E}_1 (n)^{\complement} \big] = 0; \qquad \displaystyle\lim_{n \rightarrow \infty} \mathbb{P} \big[ \mathscr{E}_2 (n)^{\complement} \big] = 0; \qquad \displaystyle\lim_{n \rightarrow \infty} \mathbb{P} \big[ \mathscr{E}_3 (n)^{\complement} \big] = 0. 
	\end{flalign}  
	
	We begin by confirming the first bound in \eqref{e123}. Observe since the gaps $\big( \mathcal{L}_j (-\mathsf{T}) - \mathcal{L}_{j+1} (-\mathsf{T}) \big)_{j \ge 1}$ of $\bm{\mathcal{L}}$ have the same law as those $2^{-1/2} \cdot (\mathfrak{a}_j - \mathfrak{a}_{j+1})_{j \ge 1}$ of a rescaled Airy point process $2^{-1/2} \cdot \bm{\mathfrak{a}}$, so do those of the gaps $\big( \mathsf{L}_j^n (-\mathsf{T}) \big)_{j \ge 1}$ of $\bm{\mathsf{L}}^n$ (by \eqref{ljn}). Next, by \Cref{kdeltad} (with the $(t, \sigma, k, u)$ there equal to $\big( 0, 1, N, (\log n)^2 \big)$ here), we have 
	\begin{flalign*}
		\mathbb{P} \Big[ \big| 2^{-1/2} \cdot \mathfrak{a}_N + 2^{-7/6} (3 \pi)^{2/3} N^{2/3} \big| \le (\log n)^2 N^{-1/3} \Big] \ge 1 - c_1^{-1} e^{-c_1 (\log n)^2},
	\end{flalign*}  
	
	\noindent for some constant $c_1 > 0$. Together with \eqref{ttln} and the fact that the law of the gaps of $\bm{\mathsf{L}}^n$ coincides with that of $2{-1/2} \cdot \bm{\mathfrak{a}}$, this gives a coupling between $\bm{\mathsf{L}}^n$ and $\bm{\mathfrak{a}}$ such that 
	\begin{flalign*}
		\mathbb{P} \Bigg[ \displaystyle\max_{j \in \llbracket 1, N \rrbracket} \big| 2^{-1/2} \cdot \mathfrak{a}_j - \mathsf{L}_j^n (-\mathsf{T}) - 2^{-1/2} n^{2/3} \big| \le (\log n)^2 N^{-1/3} \Bigg] \ge 1 - c_1^{-1} e^{-c_1 (\log n)^2}.
	\end{flalign*}
	
	\noindent Together with the fact (from \Cref{kdeltad} at $\sigma=1$, $\mathfrak{a}_j=  2^{1/2} \cdot \mathcal{S}_j (0)$, and a union bound) that 
	\begin{flalign*}
		\mathbb{P} \Bigg[ \displaystyle\max_{j \in \llbracket 1, N \rrbracket} \big| 2^{-1/2} \cdot \mathfrak{a}_j + 2^{-7/6} (3 \pi)^{2/3} j^{2/3} \big| \le (\log n)^2 j^{-1/3} \Bigg] \ge 1 - c_2^{-1} e^{-c_2 (\log n)^2},
	\end{flalign*}
	
	\noindent for some constant $c_2 \in (0, c_1]$, this yields
	\begin{flalign*}
		\mathbb{P} \Bigg[ \bigcap_{j \in \llbracket 1, N \rrbracket} \Big\{ \big| \mathsf{L}_j^n (-\mathsf{T}) + 2^{-1/2} n^{2/3} + 2^{-7/6} (3 \pi)^{2/3} j^{2/3} \big| \le (\log n)^2 (& j^{-1/3} + N^{-1/3}) \Big\} \Bigg] \\
		& \ge 1 - 2c_2^{-1} e^{-c_2 (\log n)^2},
	\end{flalign*}
	
	\noindent which implies that $\mathbb{P} \big[ \mathscr{E}_1 (n) \big] \ge 1 - 2c_2^{-1} e^{-c_2 (\log n)^2}$ (as $u_j = \mathsf{L}_j^n (-\mathsf{T})$ and $(\log n)^2 (j^{-1/3} + N^{-1/3}) \le 2 (\log n)^2 j^{-1/3} \le (\log n)^{5/2} j^{-1/3}$, for sufficiently large $n$). This verifies the first statement in \eqref{e123}; the proof of the second is entirely analogous and is thus omitted.
	
	We next confirm the third statement in \eqref{e123}. Let $\mathfrak{c}$, $\mathfrak{C}_1$, and $\mathfrak{C}_2$ denote the constants $c$, $C_1$, and $C_2$ from \Cref{sclprobability} at $(A, B, D, R) = (2, 2, 10, \mathfrak{C}_2)$. By \Cref{x1lsmall} (with the $(n, B, \vartheta)$ there equal to $(N, 2\mathfrak{C}_2, \mathfrak{C}_1^{-1})$ here), there exists a non-increasing sequence $\bm{\delta} = (\delta_1, \delta_2, \ldots )$ of real numbers with $\lim_{j \rightarrow \infty} \delta_j = 0$ such that $\delta_j \ge (\log j + 1)^{-1}$ and
	\begin{flalign}
		\label{eventc1}
		\mathbb{P} \Big[ \textbf{TOP} \big( [-2\mathfrak{C}_2 N^{1/3}, 2\mathfrak{C}_2 N^{1/3}]; \mathfrak{C}_1^{-1} N^{2/3} \big) \Big] \ge 1 - \delta_n. 
	\end{flalign}
	
	\noindent Hence, for sufficiently large $n$, we have     
	\begin{flalign}
		\label{scl2e} 
		\begin{aligned} 
			\mathbb{P} & \Big[ \textbf{REG}_{N, N+1} \big( [-2N^{1/3}, 2N^{1/3}]; 4000; 2N+2 \big) \cap \textbf{GAP}_{N+1} \big( [-2N^{1/3}, 2N^{1/3}]; \mathfrak{C}_2 \big) \Big] \\
			& \ge \mathbb{P} \Big[ \textbf{SCL}_N (2; 2; \mathfrak{C}_2) \Big] \ge \mathbb{P} \Big[ \textbf{TOP} \big( [-\mathfrak{C}_2 N^{1/3}, \mathfrak{C}_2 N^{1/3}]; \mathfrak{C}_1^{-1} N^{2/3} \big) \Big] - \mathfrak{c}^{-1} e^{-\mathfrak{c} (\log n)^2} \ge 1 - 2 \delta_n,
		\end{aligned} 
	\end{flalign}
	
	\noindent where here we set $\textbf{REG}_{N, N+1} = \textbf{REG}_N \cap \textbf{REG}_{N+1}$. Here, to obtain the first bound we used \Cref{eventscl}); to obtain the second we used \Cref{sclprobability}; and the third follows we used \eqref{eventc1} and the fact that $\delta_n \ge (\log n + 1)^{-1}$. Observe on the event in the left side of \eqref{scl2e} that 
	\begin{flalign*}
		\displaystyle\sup_{|t| \le \mathsf{T}} \big| f(t) - u_N \big| & \le \displaystyle\sup_{|t| \le \mathsf{T}} \big| \mathsf{L}_{N+1}^n (t) - \mathsf{L}_{N+1}^n (-\mathsf{T}) \big| + \mathsf{L}_N (-\mathsf{T}) - \mathsf{L}_{N+1} (-\mathsf{T}) \\
		& \le \displaystyle\sup_{|t| \le \mathsf{T}} \big| \mathcal{L}_{N+1} (t) - \mathcal{L}_{N+1} (-\mathsf{T}) \big| + (t + \mathsf{T}) |\xi_n| + 2 \mathfrak{C}_2 N^{-1/3} (\log n)^{30} \\ 
		& \le 4 \big( 4(N+1) \mathsf{T} \big)^{1/2} + 24000 N^{1/3} \mathsf{T} + n^{-10} + 2 \mathfrak{C}_2 N^{-1/3} (\log n)^{30} + 2 \mathsf{T} |\xi_n| \\
		& \le 10 n^{1/6} N^{1/2} + \big| \mathcal{L}_N (\mathsf{T}) - \mathcal{L}_N (-\mathsf{T}) \big|  \le 20 n^{1/6} N^{1/2} \le n^8.
	\end{flalign*}
	
	\noindent Here, in the first statement we used the definitions \eqref{ljn} of $\bm{u}$ and $f$; in the second we used  \Cref{gap} for the $\textbf{GAP}_{N+1}$ event, the fact that $\mathfrak{C}_2 \big( (N+1)^{2/3} - N^{2/3} \big) + N^{-1/3} \big( \log (N+1) \big)^{25} \le \mathfrak{C}_2 N^{-1/3} (\log n)^{30}$ for sufficiently large $n$, and the definition \eqref{ljn} of $\bm{\mathsf{L}}^n$; in the third we used \Cref{eventtsregular2} for the $\textbf{REG}_N$ event (with the fact that $\mathsf{T} + t \le \mathsf{T} + |t| \leq 2\mathsf{T}$ for $|t| \le \mathsf{T}$); in the fourth we used the definition \eqref{xizeta} of $\xi_n$ and the facts that $\mathsf{T} = n^{1/3}$, that $N = n^{15}$, and that $n$ is sufficiently large; and in the fifth and sixth we again used \Cref{eventtsregular2} for the $\textbf{REG}_{N+1}$ event with the fact that $n$ is sufficiently large (and that $N = n^{15}$). Hence, $\sup_{|t| \le \mathsf{T}} \big| f(t) - u_N \big| \le n^8$ holds in the event on the left side of \eqref{scl2e}. Together with \eqref{scl2e} and the fact that $\lim_{j \rightarrow \infty} \delta_j = 0$, this yields the third statement of \eqref{e123}. Hence, \eqref{e123} holds, so $\lim_{n \rightarrow \infty} \mathbb{P} \big[ \mathscr{E} (n) \big] = 1$. 
	
	Now let us condition on $\mathcal{F}_{\ext}^n$ and restrict to the event $\mathscr{E}$. Then apply \Cref{uvxrconverge}, with the $\bm{\mathsf{x}}^n$ there equal to $\bm{\mathsf{L}}^n$ here; observe that its condition \eqref{xjuv0} is verified by $\mathscr{E}_1 \cap \mathscr{E}_2$ and \eqref{unfs} by $\mathscr{E}_3$. Thus, this proposition implies that the conditional law of $\bm{\mathsf{L}}^n$ on $\mathscr{E}$ converges as $n$ tends to $\infty$ to $\bm{\mathcal{S}}$, uniformly on compact subsets of $\mathbb{Z}_{\ge 1} \times \mathbb{R}$. More precisely, for finite intervals $\Sigma \subset \mathbb{Z}_{\ge 1}$ and $I \subset \mathbb{R}$, and any bounded continuous function $F: \mathcal{C} (\Sigma \times I) \rightarrow \mathbb{R}$, we have upon denoting $\bm{\mathcal{L}}^n |_{\Sigma \times I} = \big( \mathsf{L}_j (s) \big)_{(j, s) \in \Sigma \times I}$ and $\bm{\mathcal{S}} |_{\Sigma \times I} = \big( \mathcal{S}_j (s) \big)_{(j,s) \in \Sigma \times I}$ that
	\begin{flalign} 
		\label{fns} 
		\displaystyle\lim_{n \rightarrow \infty} \Big( \mathbb{E} \big[ F ( \bm{\mathsf{L}}^n |_{\Sigma \times I}) \big| \mathcal{F}_{\ext}^n \big] - \mathbb{E} \big[ F (\bm{\mathcal{S}} |_{\Sigma \times I}) \big] \Big) \cdot \textbf{1}_{\mathscr{E}(n)} = 0,
	\end{flalign} 
	
	\noindent where the first expectation above is conditional with respect to $\mathcal{F}_{\ext}^n$. 
	
	Since $\lim_{n \rightarrow \infty} \mathbb{P} \big[ \mathscr{E} (n) \big] = 1$ (and $F$ is bounded), we find upon taking expectation over $\mathcal{F}_{\ext}^n$ that 
	\begin{flalign*}
		\displaystyle\lim_{n \rightarrow \infty} \mathbb{E} \big[ F ( \bm{\mathsf{L}}^n |_{\Sigma \times I}) \big] = \mathbb{E} \big[ F( \bm{\mathcal{S}} |_{\Sigma \times I}) \big].
	\end{flalign*}

	\noindent Since $\mathsf{L}_j^n (t) = \mathcal{L}_j (t) - \xi_n t - \zeta_n$ by \eqref{ljn}, this shows that $\big( \mathsf{L}_1^n (-1), \mathsf{L}_1^n (1) \big) = \big( \mathcal{L}_1 (-1) + \xi_n - \zeta_n, \mathcal{L}_1 (1) - \xi_n - \zeta_n \big)$ converges in law to $\big( \mathcal{S}_1 (-1), \mathcal{S}_1 (1) \big)$, as $n$ tends to $\infty$. Hence, as the variables $\big\{ \mathcal{L}_1 (-1), \mathcal{L}_1 (1), \mathcal{S}_1 (-1), \mathcal{S}_1 (1) \big\}$ are tight, the pair of random variables $(\xi_n, \zeta_n)$ is also tight. 
	
	Now, let us show that the limits $\lim_{n \rightarrow \infty} \xi_n$ and $\lim_{n \rightarrow \infty} \zeta_n$ exist almost surely. To do so, fix any $\delta > 0$. It suffices to verify that 
	\begin{flalign}
		\label{mkdelta}
		 \displaystyle\lim_{k \rightarrow \infty}	\displaystyle\lim_{m \rightarrow \infty} \displaystyle\limsup_{M\rightarrow \infty} \mathbb{P} \bigg[ \displaystyle\max_{n \in \llbracket k, m \rrbracket} \big( | \xi_n - \xi_M | + |\zeta_n - \zeta_M | \big) > \delta \bigg] = 0,
	\end{flalign} 
	
	\noindent as this would imply that
	\begin{flalign*}
		\displaystyle\lim_{k \rightarrow \infty} \mathbb{P} \bigg[ \displaystyle\max_{n, n' \ge k} \big( |\xi_n - \xi_{n'}| + |\zeta_n - \zeta_{n'}| \big) > \delta \bigg] = 0,
	\end{flalign*}
	
	\noindent which would indicate that $(\xi_n, \zeta_n)$ almost surely forms a Cauchy sequence and thus has a limit.
	
	To verify \eqref{mkdelta}, fix a real number $\delta \in (0, 1)$ and an integer $n \ge 1$; as above, set $N = N_n = n^{15}$ and $\mathsf{T} = \mathsf{T}_n = n^{1/3}$. Define the event 
	\begin{flalign}
		\label{gjm} 
		\mathscr{G}_n (M) & = \Big\{ \big| \mathsf{L}_N^M (-\mathsf{T}) + 2^{-1/2} n^{2/3} + 2^{-7/6} (3 \pi)^{2/3} N^{2/3} \big| > \displaystyle\frac{\delta}{2} \Big\} \\
		& \qquad \cap \Big\{ \big| \mathsf{L}_N^M (\mathsf{T}) + 2^{-1/2} n^{2/3} + 2^{-7/6} (3 \pi)^{2/3} N^{2/3} \big| > \displaystyle\frac{\delta}{2} \Big\}.
	\end{flalign} 
	
	\noindent Then there exists a constant $c_3 > 0$ such that
	\begin{flalign}
		\label{probabilitygnm} 
		\begin{aligned} 
		\displaystyle\lim_{M \rightarrow \infty} \mathbb{P} \big[ \mathscr{G}_n (M) \big] & = \mathbb{P} \bigg[ \Big\{ \big| \mathcal{S}_N (-\mathsf{T}) + 2^{-1/2} n^{2/3} + 2^{-7/6} (3 \pi)^{2/3} n^{2/3} \big| > \displaystyle\frac{\delta}{2} \Big\} \\
		& \qquad \qquad \cup \Big\{ \mathcal{S}_N (\mathsf{T}) + 2^{-1/2} n^{2/3} + 2^{-7/6} (3 \pi)^{2/3} N^{2/3} \big| > \displaystyle\frac{\delta}{2} \Big\} \bigg] \le c_3^{-1} e^{-c_3\delta n},
		\end{aligned} 
	\end{flalign}	
		
	\noindent where in the first statement we applied \eqref{fns} and in the second we applied \eqref{kdeltad} (with the $(k, u, \sigma, t)$ there equal to $(N, \delta N^{1/3} / 2, 1, \mathsf{T})$ here, observing that $N^{1/3} \ge n$). Recalling the definition of $\bm{\mathsf{L}}^M$ from \eqref{ljn} and of the $(\xi_i, \zeta_i)$ from \eqref{xizeta}, we find that 
	\begin{flalign*}
		& \xi_j = \xi_M + (2\mathsf{T})^{-1} \cdot \big( \mathsf{L}_N^M (\mathsf{T}) - \mathsf{L}_N^M (-\mathsf{T}) \big); \\
		& \zeta_j = \zeta_M + \displaystyle\frac{1}{2} \cdot \big( \mathsf{L}_N^M (-\mathsf{T}) + \mathsf{L}_N^M (\mathsf{T}) \big) + 2^{-1/2} n^{2/3} + 2^{-7/6} (3 \pi)^{2/3} N^{2/3}.
	\end{flalign*}
	
	\noindent Therefore, by \eqref{gjm}, we have on the event $\mathscr{G}_n (M)^{\complement}$ that $|\xi_n - \xi_M| \le \delta / 2$ and $|\xi_n - \xi_M| \le \delta / 2$. Hence, for any integers $m \ge k \ge 1$, we have by \eqref{probabilitygnm} that
	\begin{flalign*}
		\displaystyle\limsup_{M\rightarrow \infty} \mathbb{P} \bigg[ \displaystyle\max_{n \in \llbracket k, m \rrbracket} \big( | \xi_n - \xi_M | + |\zeta_n - \zeta_M | \big) > \delta \bigg] \le \displaystyle\sum_{n=k}^m \mathbb{P} \big[ \mathscr{G}_n (M)^{\complement} \big] \le c_3^{-1} \displaystyle\sum_{n=k}^{\infty} e^{-c_3 \delta n}.
	\end{flalign*}
	
	\noindent Letting first $m$ tend to $\infty$, then $k$ tend to $\infty$, and next $\delta$ tend to $0$, this yields \eqref{mkdelta}. As explained above, this implies that $(\xi_n, \zeta_n)$ converges almost surely to some pair $(\mathfrak{l}_{\infty}, \mathfrak{c}_{\infty})$, as $n$ tends to $\infty$. 

	Define $\bm{\mathsf{L}}^{\infty} = (\mathsf{L}_1^{\infty}, \mathsf{L}_2^{\infty}, \ldots ) \in \mathbb{Z}_{\ge 1} \times \mathcal{C}(\mathbb{R})$ by setting $\mathsf{L}_j^{\infty} (s) = \mathcal{L}_j (s) - \mathfrak{l}_{\infty} s - \mathfrak{c}_{\infty}$ for each $(j, s) \in \mathbb{Z}_{\ge 1} \times \mathbb{R}$. Sample a pair of random variables $(\mathfrak{l}, \mathfrak{c})$ according to the law of $(\mathfrak{l}_{\infty}, \mathfrak{c}_{\infty})$, and place $(\bm{\mathcal{S}, \mathfrak{l}, \mathfrak{c}})$ on the same probability space, so that $\bm{\mathcal{S}}$ is independent from $(\mathfrak{l}, \mathfrak{c})$. Then, for any bounded intervals $\Sigma \subset \mathbb{Z}_{\ge 1}$ and $I \subset \mathbb{R}$, and bounded continuous functions $F: \mathcal{C} (\Sigma \times I) \rightarrow \mathbb{R}$ and $G : \mathbb{R}^2 \rightarrow \mathbb{R}$, we have upon denoting $\bm{\mathsf{L}}^{\infty} |_{\Sigma \times I} = \big( \mathsf{L}_j^{\infty} (s) \big)_{(j, s) \in \Sigma \times I}$ that
	\begin{flalign*}
		\mathbb{E} \big[ F( \bm{\mathsf{L}}^{\infty} |_{\Sigma \times I}) \cdot G(\mathfrak{l}_{\infty}, \mathfrak{c}_{\infty}) \big] & = \displaystyle\lim_{n \rightarrow \infty} \mathbb{E} \big[ F (\bm{\mathsf{L}}^n |_{\Sigma \times I}) \cdot G(\mathfrak{l}_n, \mathfrak{c}_n) \big] \\
		& = \displaystyle\lim_{n \rightarrow \infty} \mathbb{E} \Big[ \mathbb{E} \big[ F (\bm{\mathsf{L}}^n |_{\Sigma \times I}) | \mathcal{F}_{\ext} \big] \cdot \textbf{1}_{\mathscr{E} (n)} \cdot G(\mathfrak{l}_n, \mathfrak{c}_n) \Big] \\
		&  = \mathbb{E} \big[ F( \bm{\mathcal{S}} |_{\Sigma \times I}) \big] \cdot \mathbb{E} \big[ G(\mathfrak{l}, \mathfrak{c}) \big].
	\end{flalign*} 
	 
	 \noindent Here, the first equality follows from the convergence of $(\mathfrak{l}_n, \mathfrak{c}_n)$ to $(\mathfrak{l}_{\infty}, \mathfrak{c}_{\infty})$ (which implies that of $\bm{\mathsf{L}}^n$ to $\bm{\mathsf{L}}^{\infty}$, uniformly on $\Sigma \times I$); the second follows from the fact that $(\mathfrak{l}_n, \mathfrak{c}_n)$ is measurable with respect to $\mathcal{F}_{\ext}^n$ and that $\lim_{n \rightarrow \infty} \mathbb{P} \big[ \mathscr{E}(n) \big] = 1$; and the third follows from \eqref{fns} (and again the fact that $\lim_{n \rightarrow \infty} \mathbb{P} \big[ \mathscr{E}(n) \big] = 1$). Since $F$ and $G$ were arbitrary, we deduce that $(\bm{\mathsf{L}}^{\infty}, \mathfrak{l}_{\infty}, \mathfrak{c}_{\infty})$ coincides in law with $(\bm{\mathcal{S}}, \mathfrak{l}, \mathfrak{c})$. This yields the theorem, upon recalling that $\bm{\mathsf{L}}^{\infty} (s) = \bm{\mathcal{L}} (s) - \mathfrak{l}_{\infty} s - \mathfrak{c}_{\infty}$.  
\end{proof}

\subsection{Approximate Parabolicity of Paths} 

\label{Pathsx2t} 

\begin{figure}
	\center
	\includegraphics[width=0.7\textwidth, trim=0 0.5cm 0 0.5cm, clip]{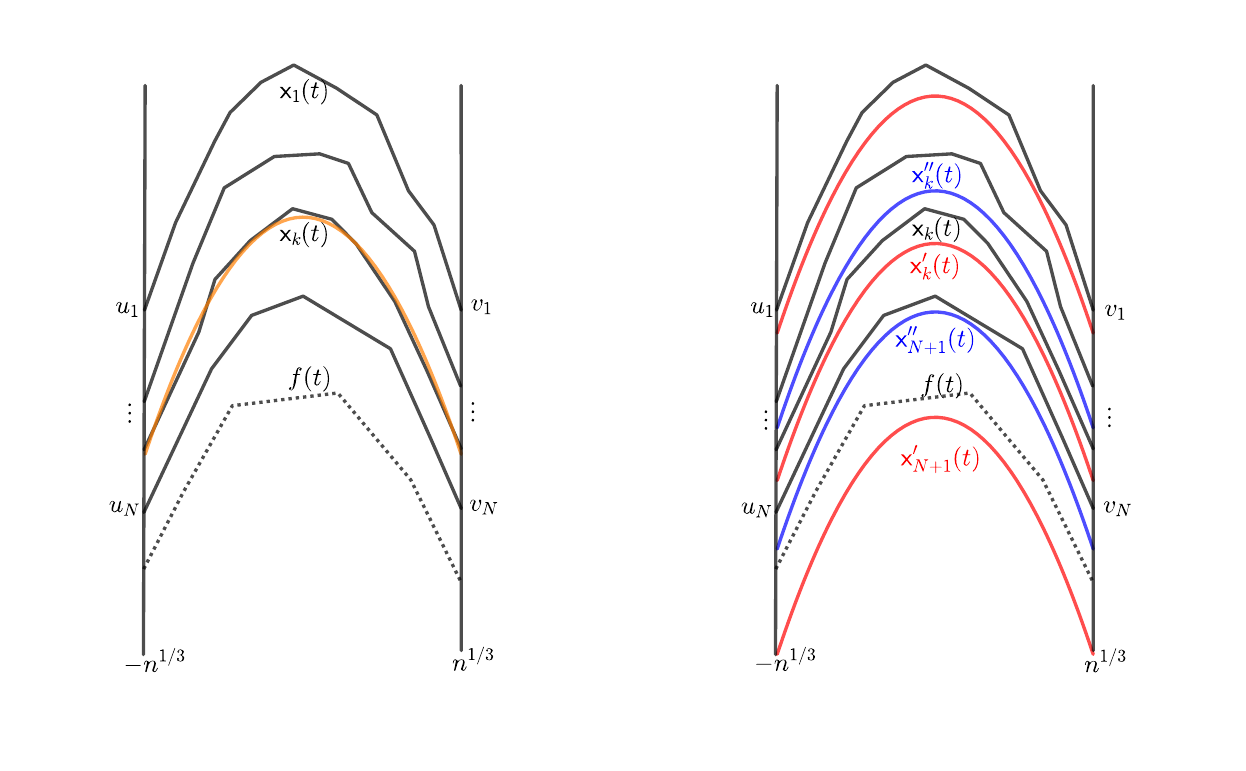}
	
	\caption{Shown to the left is a depiction of \Cref{uvxrn10}, indicating that $\mathsf{x}_k$ is close to the orange parabola. Shown to the right is a depiction of its proof.} 
	\label{f:Airy_Line2}
\end{figure}

To establish \Cref{uvxrconverge}, we will use the following lemma indicating that the paths in $\bm{\mathsf{x}}^n$ closely approximate parabolas; see the left side of \Cref{f:Airy_Line2} for a depiction.

\begin{lem}
	
	\label{uvxrn10}
	
	Adopt the notation and assumptions of \Cref{uvxrconverge}. There exists a constant $c > 1$ such that, denoting $k = \lceil n^{1/6} \rceil$, we have
	\begin{flalign*}
		\mathbb{P} \Bigg[ \displaystyle\sup_{|t| \le n^{1/3}} \big| \mathsf{x}_k^n (t) + 2^{-1/2} t^2 + 2^{-7/6} (3 \pi)^{2/3} k^{2/3} \big| \ge k^{-1/30} \Bigg] \le c^{-1} e^{-c(\log n)^2}.
	\end{flalign*} 
	
\end{lem} 

\begin{proof}
	
	Throughout this proof, we abbreviate $\bm{\mathsf{x}} = \bm{\mathsf{x}}^n$ and $\mathsf{x}_j = \mathsf{x}_j^n$ for each integer $j \in \llbracket 1, N \rrbracket$; we also set $\mathsf{T} = n^{1/3}$. We establish the lemma by bounding $\bm{\mathsf{x}}$ between two rescaled parabolic Airy line ensembles with approximately equal parameters. To this end, define the line ensembles $\bm{\mathsf{y}}' = (\mathsf{y}_1', \mathsf{y}_2', \ldots ) \in \mathbb{Z}_{\ge 1} \times \mathcal{C} (\mathbb{R})$ and $\bm{\mathsf{y}}'' = (\mathsf{y}_1'', \mathsf{y}_2'', \ldots ) \in \mathbb{Z}_{\ge 1} \times \mathcal{C} (\mathbb{R})$ by 
	\begin{flalign} 
		\label{yjr1} 
		\bm{\mathsf{y}}' = \bm{\mathcal{S}}^{(\sigma')}, \quad \text{and} \quad \bm{\mathsf{y}}'' = \bm{\mathcal{S}}^{(\sigma'')}, \qquad \text{where} \quad \sigma' = 1 - n^{-1}, \quad \text{and} \quad \sigma'' = 1 + n^{-1}.
	\end{flalign} 
	
	\noindent where we recall the rescaled parabolic Airy line ensemble $\bm{\mathcal{S}}^{(\sigma)}$ from \eqref{sigmar}; we also set $\bm{\mathsf{y}}_j' = \bm{\mathsf{y}}_j'' = \infty$ if $j \le 0$. We further define $\bm{\mathsf{x}}' = (\mathsf{x}_1', \mathsf{x}_2', \ldots ) \in \mathbb{Z}_{> -n_0} \times \mathcal{C} (\mathbb{R})$ and $\bm{\mathsf{x}}'' = (\mathsf{x}_{n_0+1}'', \mathsf{x}_{n_0+2}'', \ldots ) \in \mathbb{Z}_{> n_0} \times \mathcal{C} (\mathbb{R})$ by shifting the indices of $\bm{\mathsf{y}}'$ and $\bm{\mathsf{y}}''$ respectively, namely, by setting
	\begin{flalign}
		\label{xjr1}
		\mathsf{x}_{j-n_0}' = \mathsf{y}_j' - n^{-1/4}, \quad \text{and} \quad \mathsf{x}_{j+n_0}'' = \mathsf{y}_j'' + n^{-1/4},  \qquad \text{where} \quad n_0 = \lfloor n^{1/50} \rfloor,
	\end{flalign} 
	
	\noindent for any $j \in \mathbb{Z}_{\ge 1}$; we also set $\mathsf{x}_j' = \infty$ if $j \le -n_0$ and $\mathsf{x}_j'' = \infty$ if $j \le n_0$. We will show that it is with high probability possible to couple $\bm{\mathsf{x}}$ to lie between $\bm{\mathsf{x}}'$ and $\bm{\mathsf{x}}''$. See the right side of \Cref{f:Airy_Line2}.

	To this end, we define the event $\mathscr{E} = \mathscr{E}_1 \cap \mathscr{E}_2$, where $\mathscr{E}_1 = \mathscr{E}_1' \cap \mathscr{E}_1'' $ and $\mathscr{E}_2 = \mathscr{E}_2' \cap \mathscr{E}_2''$. Here, 
	\begin{flalign}
		\label{11e22e} 
		\begin{aligned} 
			\mathscr{E}_1' & = \bigcap_{j=1}^N \Bigg\{ \displaystyle\sup_{t \in \{ -\mathsf{T}, \mathsf{T} \}} \big| \mathsf{x}_j' (t) + \mathfrak{p} (j + n_0; t; \sigma') + n^{-1/4} \big| \le (\log n)^3 (j+n_0)^{-1/3} \Bigg\}; \\ 
			\mathscr{E}_1''&  = \bigcap_{j=n_0+1}^N \Bigg\{ \displaystyle\sup_{t \in \{ -\mathsf{T}, \mathsf{T} \}} \big| \mathsf{x}_j'' (t) +  \mathfrak{p} (j-n_0; t; \sigma'') - n^{-1/4} \big| \le (\log n)^3 (j-n_0)^{-1/3} \Bigg\}; \\ 
			\mathscr{E}_2' & = \Bigg\{ \displaystyle\sup_{t \in [ -\mathsf{T}, \mathsf{T}]} \big| \mathsf{x}_{N+1}' (t) + \mathfrak{p} (N+n_0+1; t; \sigma') + n^{-1/4} \big| \le (\log n)^3 N^{-1/3} \Bigg\}; \\ 
			\mathscr{E}_2'' & = \Bigg\{ \displaystyle\sup_{t \in [-\mathsf{T}, \mathsf{T}]} \big| \mathsf{x}_{N+1}'' (t) + \mathfrak{p} (N-n_0+1; t; \sigma'') - n^{-1/4} \big| \le (\log n)^3 N^{-1/3} \Bigg\}. 
		\end{aligned} 
	\end{flalign} 
	
	\noindent where we have denoted
	\begin{flalign}
		\label{jtsigmap} 
		\mathfrak{p} (j; t; \sigma) = 2^{-1/2} \sigma^3 t^2 + 2^{-7/6} (3 \pi)^{2/3} \sigma^{-1} j^{2/3}.
	\end{flalign}
	
	\noindent Applying the definitions \eqref{xjr1} and \eqref{yjr1} of $\bm{\mathsf{x}}$ and $\bm{\mathsf{y}}$ in terms of rescaled parabolic Airy line ensembles; the concentration estimate \Cref{kdeltad} for the latter (with the fact that $N-n_0 \ge N/2$ for sufficiently large $n$); and a union bound yields a constant $c_1 > 0$ such that 
	\begin{flalign}
		\label{probabilitye1122} 
		\max \Big\{ \mathbb{P} \big[ \mathscr{E}_1'^{\complement} \big], \mathbb{P} \big[ \mathscr{E}_1''^{\complement} \big], \mathbb{P} \big[ \mathscr{E}_2'^{\complement} \big], \mathbb{P} \big[ \mathscr{E}_2''^{\complement} \big] \Big\} \le (4c_1)^{-1} e^{-c_1 (\log n)^2}, \qquad \text{so} \qquad \mathbb{P} [\mathscr{E}^{\complement}] \le c_1^{-1} e^{-c_1 (\log n)^2}. 
	\end{flalign}
	
	\noindent Now condition on the curves $\big( \mathsf{x}_j' (t) \big)$ and $\big( \mathsf{x}_j'' (t) \big)$ for $(j, t) \notin \llbracket 1, N \rrbracket \times (-\mathsf{T}, \mathsf{T})$, and restrict to the event $\mathscr{E}$. We claim that 
	\begin{flalign}
		\label{xj3xnfxn} 
		\begin{aligned} 
			& \mathsf{x}_j' (t) \le \mathsf{x}_j (t) \le \mathsf{x}_j'' (t), \qquad \qquad \text{for each $(j, t) \in \llbracket 1, N \rrbracket \times \{ -\mathsf{T}, \mathsf{T} \}$}; \\
			& \mathsf{x}_{N+1}' (t) \le f(t) \le \mathsf{x}_{N+1}'' (t), \qquad   \text{for each $t \in [-\mathsf{T}, \mathsf{T}]$}.
		\end{aligned} 
	\end{flalign}
	
	To this end, observe for any $(j, t) \in \llbracket 1, N \rrbracket \times \{ -\mathsf{T}, \mathsf{T} \}$ and sufficiently large $n$ that 
	\begin{flalign*}
		\mathsf{x}_j (t) - \mathsf{x}_j' (t) & \ge \big( -(\log n)^3 j^{-1/3} - 2^{-1/2} n^{2/3} - 2^{-7/6} (3 \pi)^{2/3} j^{2/3} \big) \\
		& \qquad - \big( (\log n)^3 (j+n_0)^{-1/3} -2^{-1/2} \sigma'^3 n^{2/3} - 2^{-7/6} (3 \pi)^{2/3} \sigma'^{-1} (j+n_0)^{2/3} - n^{-1/4} \big) \\
		& = 2^{-7/6} (3 \pi)^{2/3} \big( (j+n_0)^{2/3} - j^{2/3} \big) - 2 (\log n)^3 j^{-1/3} \\
		& \qquad + 2^{-1/2} (\sigma'^3 - 1) n^{2/3} + 2^{-7/6} (3 \pi)^{2/3} (\sigma'^{-1} - 1) (j+n_0)^{2/3} + n^{-1/4} \\				
		& \ge n_0^{2/3} j^{-1/3} - 2 (\log n)^3 j^{-1/3} -2^{3/2} n^{-1/3} + n^{-1/4} \ge 0,
	\end{flalign*}
	
	\noindent where in the first statement we used \eqref{xjuv0}, the fact that we are restricting to the event $\mathscr{E}_1'$ from \eqref{11e22e}, and the definition \eqref{jtsigmap} of $\mathfrak{p}$; in the second we performed the subtraction; in the third we used the facts that $\sigma'^3 - 1 \ge -4n^{-1}$ and $\sigma'^{-1} - 1 \ge n^{-1} \ge 0$ (by the definition \eqref{yjr1} of $\sigma' = 1 - n^{-1}$), that $(j+n_0)^{2/3} - j^{2/3} \ge n_0^{2/3} / (2 j^{1/3})$, and that $2^{-7/6} (3 \pi)^{2/3} \ge 2$; and in the fourth we used the definition \eqref{xjr1} of $n_0$ and the fact that $n$ is sufficiently large. This verifies the first bound in the first statement of \eqref{xj3xnfxn}; the proof of the second part is entirely analogous (upon taking into account the fact that $\mathsf{x}_j'' = \infty$ for $j \in \llbracket 1, n_0 \rrbracket$) and is therefore omitted.
	
	To verify the second statement in \eqref{xj3xnfxn}, observe for any $t \in [-\mathsf{T}, \mathsf{T}]$ that
	\begin{flalign*}
		f(t &) - \mathsf{x}_{N+1}' (t) \\
		& \ge \big( -2^{-1/2} n^{2/3} - 2^{-7/6} (3 \pi)^{2/3} N^{2/3} - (\log n)^3 N^{-1/3} - n^8 \big) \\
		& \quad - \big( (\log n)^3 N^{-1/3} - 2^{-1/2} \sigma'^3 t^2  - 2^{-7/6} (3 \pi)^{2/3} \sigma'^{-1} (N+n_0+1)^{2/3} - n^{-1/4} \big) \\
		& \ge 2^{-7/6} (3 \pi)^{2/3} \big( (1 + n^{-1}) (N+n_0+1)^{2/3} - N^{2/3}\big) - n^8 -2^{-1/2} n^{2/3} - 2 (\log n)^3 N^{-1/3} \\
		& \ge n^{-1} N^{2/3} - 2n^8 \ge n^9 - 2n^8 \ge 0,
	\end{flalign*}
	
	\noindent where in the first inequality we used \eqref{unfs}, \eqref{xjuv0} (to bound $u_N$), the fact that we are restricting to the event $\mathscr{E}_2'$ from \eqref{11e22e}, and the definition \eqref{jtsigmap} of $\mathfrak{p}$; in the second we used the fact that $\sigma'^{-1} \ge 1+n^{-1}$ (by the definition \eqref{xjr1} of $\sigma' = 1 - n^{-1}$); in the third we used the facts that $2^{-7/6} (3 \pi)^{2/3} \ge 1$, that $(1+n^{-1}) (N+n_0+1)^{2/3} \ge N^{2/3} + n^{-1} N^{2/3}$, and that $n$ is sufficiently large; in the fourth we used the fact that $N = n^{15}$; and in the fifth we used the fact that $n \ge 2$. This verifies the first part of the second bound in \eqref{xj3xnfxn}; the proof of the second part is entirely analogous and is therefore omitted. 
	
	Thus, \eqref{xj3xnfxn} holds. Denote the four $(N-n_0)$-tuples $\bm{u}' = \bm{\mathsf{x}}'_{\llbracket n_0 + 1, N \rrbracket} (-\mathsf{T}) \in \mathbb{W}_{N-n_0}$, $\bm{v}' = \bm{\mathsf{x}}_{\llbracket n_0 + 1, N \rrbracket} '(\mathsf{T}) \in \mathbb{W}_{N-n_0}$, $\bm{u}'' = \bm{\mathsf{x}}''_{\llbracket n_0 + 1, N \rrbracket} (-\mathsf{T}) \in \mathbb{W}_{N-n_0}$, and $\bm{v}'' = \bm{\mathsf{x}}''_{\llbracket n_0 + 1, N \rrbracket} (\mathsf{T}) \in \mathbb{W}_{N-n_0}$. By \eqref{xjr1}, and the fact from \Cref{propertya} that $\bm{\mathsf{y}}'$ and $\bm{\mathsf{y}}''$ satisfy the Brownian Gibbs property, the laws of $\big( \mathsf{x}_j' (s) \big)$ and $\big( \mathsf{x}_j'' (s) \big)$ for $(j, s) \in \llbracket n_0 + 1, N \rrbracket \times [-\mathsf{T}, \mathsf{T}]$ are given by $\mathsf{Q}_{\mathsf{x}_{N+1}'; \mathsf{x}_{n_0}'}^{\bm{u}'; \bm{v}'}$ and $\mathsf{Q}_{\mathsf{x}_{N+1}'', \mathsf{x}_{n_0}''}^{\bm{u}''; \bm{v}''}$, respectively. Thus, by \eqref{xj3xnfxn} and \Cref{monotoneheight}, we may couple $\bm{\mathsf{x}}$, $\bm{\mathsf{x}}'$, and $\bm{\mathsf{x}}''$ so that 
	\begin{flalign}
		\label{xj3} 
		\mathsf{x}_j' (s) \le \mathsf{x}_j (s) \le \mathsf{x}_j'' (s), \qquad \text{for each $(j, s) \in \llbracket n_0+1, N \rrbracket \times [-\mathsf{T}, \mathsf{T}]$}.
	\end{flalign}

	Now recall that $k = \lceil n^{1/6} \rceil > n_0$, and define the event $\mathscr{E}_3 = \mathscr{E}_3' \cap \mathscr{E}_3''$, where
	\begin{flalign*} 
		\mathscr{E}_3' & = \Bigg\{ \displaystyle\sup_{t \in [ -\mathsf{T}, \mathsf{T}]} \big| \mathsf{x}_k' (t) + 2^{-1/2} t^2  + 2^{-7/6} (3 \pi)^{2/3} k^{2/3} \big| \le k^{-1/30} \Bigg\}; \\ 
		\mathscr{E}_3'' & = \Bigg\{ \displaystyle\sup_{t \in [ -\mathsf{T}, \mathsf{T}]} \big| \mathsf{x}_k'' (t) + 2^{-1/2} t^2  + 2^{-7/6} (3 \pi)^{2/3} k^{2/3} \big| \le k^{-1/30} \Bigg\}.  
	\end{flalign*} 
	
	\noindent We claim that    
	\begin{flalign} 
		\label{e30} 
		\max \Big\{	\mathbb{P} \big[ \mathscr{E}_3'^{\complement} \big], \mathbb{P} \big[ \mathscr{E}_3''^{\complement} \big] \Big\} \le c_1^{-1} e^{-c_1 (\log n)^2}, \quad \text{so} \quad \mathbb{P} \big[ \mathscr{E}_3^{\complement} \big] \le 2 c_1^{-1} e^{-c_1 (\log n)^2}.
	\end{flalign} 
	
	\noindent Together with \eqref{probabilitye1122}, \eqref{e30} would imply that $\mathbb{P} \big[ (\mathscr{E} \cap \mathscr{E}_3)^{\complement} \big] \le 3c_1^{-1} e^{-c_1 (\log n)^2}$. Since $\mathsf{x}_k (t) + 2^{-1/2} t^2 + 2^{-7/6} (3 \pi)^{2/3} k^{2/3} \big| \le k^{-1/30}$ holds on $\mathscr{E} \cap \mathscr{E}_3$ by the definitions of $\mathscr{E}_3'$ and $\mathscr{E}_3''$ and \eqref{xj3}, this would imply the lemma. Hence, it suffices to verify \eqref{e30}. 
	
	We only establish the first bound there, as the proof of the second is entirely analogous. To this end, observe from \Cref{kdeltad} (and the definitions \eqref{xjr1} and \eqref{yjr1} of $\mathsf{x}_j'$ and $\mathsf{y}_j'$) that 
	\begin{flalign*}
		\mathbb{P} \Bigg[ \displaystyle\sup_{t \in [-\mathsf{T}, \mathsf{T}]} \big| \mathsf{x}_k' (t) + 2^{-1/2} \sigma'^3 t^2 + 2^{-7/6} (3 \pi)^{2/3} \sigma'^{-1} (k+n_0)^{2/3} \big| \le (\log k)^2 k^{-1/3} \Bigg] \le c_1^{-1} e^{-c_1 (\log n)^2}. 
	\end{flalign*}
	
	\noindent This, together with the fact that, for any $t \in [-\mathsf{T}, \mathsf{T}]$, 
	\begin{flalign*}
		\Big| & \big(2^{-1/2} \sigma'^3 t^2 + 2^{-7/6} (3 \pi)^{2/3} \sigma'^{-1} (k+n_0)^{2/3} \big) - \big(2^{-1/2} t^2 + 2^{-7/6} (3 \pi)^{2/3} k^{2/3} \big) \Big| \\
		& \le |1-\sigma'^3| n^{2/3} + 5 \big( (k+n_0)^{2/3} - k^{2/3} \big) + 5| \sigma'^{-1} - 1| (k+n_0)^{2/3} \\
		& \le 5 \big( (k+n_0)^{2/3} - k^{2/3} \big) + 4n^{-1/3} + 10n^{-1} (k + n_0)^{2/3} \le  5n_0 k^{-1/3} + 25 n^{-1/3} \le n^{-1/30},
	\end{flalign*}
	
	\noindent implies the first bound in \eqref{e30} and thus the lemma. To establish the first statement above, we used the facts that $t^2 \le n^{2/3}$ and that $2^{-7/6} (3\pi)^{2/3} \le 5$; to establish the second we used the facts that $|1-\sigma'^3| \le 4n^{-1}$ and $|\sigma'^{-1} - 1| \le 2n^{-1}$ for sufficiently large $n$ (by the definition \eqref{yjr1} of $\sigma' = 1 - n^{-1}$); to establish the third we used the facts that $(k+n_0)^{2/3} - k^{2/3} \le n_0 k^{-1/3}$ and $k+n_0  \le n$; and to establish in the fourth we used the facts that $k\ge n^{1/6}$ and that $n_0 \le n^{1/50}$ (and that $n$ is sufficiently large).
\end{proof}

\subsection{Proof of \Cref{uvxrconverge}}

\label{ProofxR}

In this section we establish \Cref{uvxrconverge}. Given \Cref{uvxrn10}, its proof is similar to that of \cite[Proposition 3.18]{ESTP}, obtained by locally sandwiching $\bm{\mathsf{x}}^n$ between two rescaled parabolic Airy line ensembles with slightly different curvatures.

\begin{proof}[Proof of \Cref{uvxrconverge}] 
	
	Throughout this proof, we abbreviate $\bm{\mathsf{x}} = \bm{\mathsf{x}}^n$ and $\mathsf{x}_j = \mathsf{x}_j^n$ for each integer $j \in \llbracket 1, N \rrbracket$. We also set $\mathsf{T} = n^{1/3}$, abbreviate $n' = \lceil n^{1/6} \rceil - 1$, and define the real numbers 
	\begin{flalign}
		\label{sigma12}
		\sigma' = 1 + n^{-1/4}; \qquad \sigma'' = 1 - n^{-1/4}.
	\end{flalign}
	
	\noindent Further define the line ensembles $\bm{\mathsf{x}}' = (\mathsf{x}_1', \mathsf{x}_2', \ldots ) \in \mathbb{Z}_{\ge 1} \times \mathcal{C} (\mathbb{R})$ and $\bm{\mathsf{x}}'' = (\mathsf{x}_1'', \mathsf{x}_2'', \ldots ) \in \mathbb{Z}_{\ge 1} \times \mathcal{C} (\mathbb{R})$ by for each integer $j \ge 1$ setting
	\begin{flalign} 
		\label{yjr2} 
		\mathsf{x}_j' = \mathcal{S}_j^{(\sigma')} - n^{-1/75}, \quad \text{and} \quad \mathsf{x}'' = \mathcal{S}_j^{(\sigma'')} + n^{-1/75}.
	\end{flalign} 
	
	\noindent where we recall the rescaled parabolic Airy line ensemble $\bm{\mathcal{S}}^{(\sigma)} = \big( \mathcal{S}_1^{(\sigma)}, \mathcal{S}_2^{(\sigma)} , \ldots  \big)$ from \eqref{sigmar}. Then, $\bm{\mathsf{x}}'$ and $\bm{\mathsf{x}}'$ both satisfy the Brownian Gibbs property by \Cref{propertya}. We will show that it is with high probability possible to couple $\bm{\mathsf{x}}$ to lie between $\bm{\mathsf{x}}'$ and $\bm{\mathsf{x}}''$; see the right side of \Cref{f:Airy_Line}.
	
	To this end, we define the event $\mathscr{E} = \mathscr{E}_1 \cap \mathscr{E}_2$, where $\mathscr{E}_1 = \mathscr{E}_1' \cap \mathscr{E}_1'' $ and $\mathscr{E}_2 = \breve{\mathscr{E}}_2 \cap \mathscr{E}_2' \cap \mathscr{E}_2''$. Here, 
	\begin{flalign}
		\label{211e22e} 
		\begin{aligned} 
			\mathscr{E}_1' & = \bigcap_{j=1}^{n'} \Bigg\{ \displaystyle\sup_{t \in \{ -\mathsf{T}, \mathsf{T} \}} \big| \mathsf{x}_j' (t) + \mathfrak{p} (j; t; \sigma') + n^{-1/75} \big| \le (\log n)^3 j^{-1/3} \Bigg\}; \\ 
			\mathscr{E}_1''&  = \bigcap_{j=1}^{n'} \Bigg\{ \displaystyle\sup_{t \in \{ -\mathsf{T}, \mathsf{T} \}} \big| \mathsf{x}_j'' (t) +  \mathfrak{p} (j; t; \sigma'') - n^{-1/75} \big| \le (\log n)^3 j^{-1/3} \Bigg\}; \\ 
			\breve{\mathscr{E}}_2 & = \Bigg\{ \displaystyle\sup_{t \in [-\mathsf{T}, \mathsf{T}]} \big| \mathsf{x}_{n'+1} (t) + 2^{-1/2} t^2 + 2^{-7/6} (3 \pi)^{2/3} (n'+1)^{2/3} \big| \le n^{-1/30} \Bigg\}; \\
			\mathscr{E}_2' & = \Bigg\{ \displaystyle\sup_{t \in [ -\mathsf{T}, \mathsf{T}]} \big| \mathsf{x}_{n'+1}' (t) + \mathfrak{p} (n'+1; t; \sigma') + n^{-1/75} \big| \le (\log n)^3 n'^{-1/3} \Bigg\}; \\ 
			\mathscr{E}_2'' & = \Bigg\{ \displaystyle\sup_{t \in [-\mathsf{T}, \mathsf{T}]} \big| \mathsf{x}_{n'+1}'' (t) + \mathfrak{p} (n'+1; t; \sigma'') - n^{-1/75} \big| \le (\log n)^3 n'^{-1/3} \Bigg\}. 
		\end{aligned} 
	\end{flalign} 
	
	\noindent where we recall the function $\mathfrak{p}$ from \eqref{jtsigmap}. Applying the definitions \eqref{yjr2} of $\bm{\mathsf{x}}$ in terms of rescaled parabolic Airy line ensembles; the concentration estimate \Cref{kdeltad} for the latter; and a union bound yields a constant $c_1 > 0$ such that 
	\begin{flalign*}
		\max \Big\{ \mathbb{P} \big[ \mathscr{E}_1'^{\complement} \big], \mathbb{P} \big[ \mathscr{E}_1''^{\complement} \big], \mathbb{P} \big[ \mathscr{E}_2'^{\complement} \big], \mathbb{P} \big[ \mathscr{E}_2''^{\complement} \big] \Big\} \le (5c_1)^{-1} e^{-c_1 (\log n)^2}, 
	\end{flalign*} 
	
	\noindent Together with the bound $\mathbb{P} \big[ \breve{\mathscr{E}}_2^{\complement} \big] \le (5c_1)^{-1} e^{-c_1 (\log n)^2}$ (by \Cref{uvxrn10}) and a union bound, this yields 
	\begin{flalign}
		\label{2probabilitye}
		\mathbb{P} \big[\mathscr{E}^{\complement} \big] \le c_1^{-1} e^{-c_1 (\log n)^2}. 
	\end{flalign}
	
	\noindent Now condition on the curves $\big( \mathsf{x}_j (t) \big)$, $\big( \mathsf{x}_j' (t) \big)$, and $\big( \mathsf{x}_j'' (t) \big)$ for $(j, t) \notin \llbracket 1, n' \rrbracket \times (-\mathsf{T}, \mathsf{T})$, and restrict to the event $\mathscr{E}$. We claim that 
	\begin{flalign}
		\label{1xj3xnfxn} 
		\begin{aligned} 
			& \mathsf{x}_j' (t) \le \mathsf{x}_j (t) \le \mathsf{x}_j'' (t), \qquad \qquad \qquad \text{for each $(j, t) \in \llbracket 1, n' \rrbracket \times \{ -\mathsf{T}, \mathsf{T} \}$}; \\
			& \mathsf{x}_{n'+1}' (t) \le \mathsf{x}_{n'+1} (t) \le \mathsf{x}_{n'+1}'' (t), \qquad   \text{for each $t \in [-\mathsf{T}, \mathsf{T}]$}.
		\end{aligned} 
	\end{flalign}
	
	To this end, observe for any $(j, t) \in \llbracket 1, N \rrbracket \times \{ -\mathsf{T}, \mathsf{T} \}$ and sufficiently large $n$ that 
	\begin{flalign*}
		\mathsf{x}_j (t) - \mathsf{x}_j' (t) & \ge \big( -(\log n)^3 j^{-1/3} - 2^{-1/2} n^{2/3} - 2^{-7/6} (3 \pi)^{2/3} j^{2/3} \big) \\
		& \qquad - \big( (\log n)^3 j^{-1/3} -2^{-1/2} \sigma'^3 n^{2/3} - 2^{-7/6} (3 \pi)^{2/3} \sigma'^{-1} j^{2/3} - n^{-1/75} \big) \\
		& = 2^{-1/2} (\sigma'^3 - 1) n^{2/3} - 2 (\log n)^3 j^{-1/3} + 2^{-7/6} (3 \pi)^{2/3} (\sigma'^{-1} - 1) j^{2/3} + n^{-1/75} \\				
		& \ge - 2 (\log n)^3 j^{-1/3} - 5 n^{-1/4} j^{2/3} + n^{-1/75} \ge 0,
	\end{flalign*}
	
	\noindent where in the first statement we used \eqref{xjuv0}, the fact that we are restricting to the event $\mathscr{E}_1'$ from \eqref{211e22e}, and the definition \eqref{jtsigmap} of $\mathfrak{p}$; in the second we performed the subtraction; in the third we used the facts that $\sigma'^3 - 1 \ge 0$, that ${\sigma'}^{-1} - 1 \ge -n^{-1/4} $ (by the definition \eqref{sigma12} of $\sigma' = 1 + n^{-1/4}$), and that $2^{-7/6} (3 \pi)^{2/3} \le 5$; and in the fourth we used the fact that $j \le n' \le n^{1/6}$ and that $n$ is sufficiently large. This verifies the first bound in the first statement of \eqref{1xj3xnfxn}; the proof of the second part is entirely analogous and is therefore omitted.
	
	To verify the second statement in \eqref{1xj3xnfxn}, observe for any $t \in [-\mathsf{T}, \mathsf{T}]$ that
	\begin{flalign*}
		& \mathsf{x}_{n'+1}  (t) - \mathsf{x}_{n'+1}' (t) \\
		& \quad \ge \big( -2^{-1/2} t^2 - 2^{-7/6} (3 \pi)^{2/3} (n'+1)^{2/3}  - n^{-1/50} \big) \\
		& \qquad - \big( (\log n)^3 n'^{-1/3} - 2^{-1/2} \sigma'^3 t^2  - 2^{-7/6} (3 \pi)^{2/3} \sigma'^{-1} (n'+1)^{2/3} - n^{-1/75} \big) \\
		& \quad = n^{1/75} + 2^{-7/6} (3 \pi)^{2/3} (\sigma'^{-1} - 1) (n'+1)^{2/3} + 2^{-1/2} (\sigma'^3 - 1) t^2  - n^{-1/50} -  (\log n)^3 n'^{-1/3} \\
		& \quad \ge n^{-1/75} - 5 n^{-1/4} (n'+1)^{2/3} - n^{-1/50} - (\log n)^3 n'^{-1/3}  \ge 0,
	\end{flalign*}
	
	\noindent where in the first statement we used the fact that we are restricting to the events $\breve{\mathscr{E}}_2$ and $\mathscr{E}_2'$ from \eqref{211e22e}, as well as the definition \eqref{jtsigmap} of $\mathfrak{p}$; in the second we performed the subtraction; in the third we used the facts that $\sigma' \ge 1$, that $\sigma'^{-1}-1 \le -n^{-1/4}$ for $n$ sufficiently large (by the definition \eqref{sigma12} of $\sigma' = 1 + n^{-1/4}$), and that $2^{-7/6} (3 \pi)^{2/3} \le 5$; and in the fourth we used the facts that $n' + 1 \le n^{1/6} + 1 \le 2n^{1/6}$ and that $n$ is sufficiently large. This verifies the first part of the second bound in \eqref{xj3xnfxn}; the proof of the second part is entirely analogous and is therefore omitted. 
	
	Thus, \eqref{1xj3xnfxn} holds. As in the proof of \Cref{uvxrn10}, it follows from \eqref{1xj3xnfxn} and \Cref{monotoneheight} that, on $\mathscr{E}$, we may couple $\bm{\mathsf{x}}$, $\bm{\mathsf{x}}'$, and $\bm{\mathsf{x}}''$ so that 
	\begin{flalign*}
		\mathsf{x}_j' (s) \le \mathsf{x}_j (s) \le \mathsf{x}_j'' (s), \qquad \text{for each $(j, s) \in \llbracket 1, n' \rrbracket \times [-\mathsf{T}, \mathsf{T}]$}. 
	\end{flalign*}			
	
	\noindent Since \eqref{yjr2} and \eqref{sigma12} imply that $\mathsf{x}_j' (s)$ and $\mathsf{x}_j'' (s)$ both converge to $\bm{\mathcal{S}}$, uniformly on compact subsets on $\mathbb{Z}_{\ge 1} \times \mathbb{R}$, as $n$ tends to $\infty$, this and \eqref{2probabilitye} together imply that $\bm{\mathsf{x}}^n$ converges to $\bm{\mathcal{S}}$, establishing the proposition.
\end{proof}

\chapter{Appendices}

\addtocontents{toc}{\protect\setcounter{tocdepth}{0}}

	\section{Proofs of Results From Chapter \ref{BRIDGES0}}
	
	\label{ProofBridgeSum}

	\subsection{Proofs of \Cref{fg0g}, \Cref{monotoneheight}, and \Cref{yconvergedelta}}
	
	\label{Proofydelta}

	In this section we first show \Cref{fg0g}, then \Cref{boundaryconvergeensemble} below, next \Cref{monotoneheight}, and subsequently \Cref{yconvergedelta}.
	
	\begin{proof}[Proof of \Cref{fg0g}]
		
		Condition on $\mathcal{F}_{\ext}$, and let $\mathscr{G} = \big\{ \mathbb{P} [\mathscr{E} | \mathcal{F}_{\ext}] \ge 1-a \big\}$, which is measurable with respect to $\mathcal{F}_{\ext}$ and satisfies the second statement of the lemma. Then, 
		\begin{flalign*}
			\mathbb{P} \big[ \mathscr{G}^{\complement} \big] = \mathbb{P} \big[ \mathbb{P} [ \mathscr{E}^{\complement} | \mathcal{F}_{\ext}] \ge a \big] \le a^{-1} \cdot \mathbb{E} \big[ \mathbb{P} [\mathscr{E}^{\complement}| \mathcal{F}_{\ext}] \big] = a^{-1} \cdot \mathbb{P} \big[\mathscr{E}^{\complement} \big] \le b,
		\end{flalign*}
		
		\noindent where the first statement follows from the definition of $\mathscr{G}$; the second from a Markov estimate, and the third and fourth from our assumption that $\mathbb{P}[\mathscr{E}] \ge 1-ab$. This establishes the lemma.
	\end{proof}

The following lemma indicates that the limiting procedure used to define the measure $\mathsf{Q}_{f; g}^{\bm{u}; \bm{v}}$ (from \Cref{qxyfg}) when either $\bm{u} \in \overline{\mathbb{W}}_n \setminus \mathbb{W}_n$ or $\bm{v} \in \overline{\mathbb{W}}_n \setminus \mathbb{W}_n$ is well-defined.  

\begin{lem} 
	
	\label{boundaryconvergeensemble} 
	
	Fix an integer $n \ge 1$; real numbers $a < b$; two $n$-tuples $\bm{u}, \bm{v} \in \overline{\mathbb{W}}_n$; and two continuous functions $f, g : [a, b] \rightarrow \overline{\mathbb{R}}$ with $f \le g$, such that $f(a) \le u_n \le u_1 \le g(a)$ and $f(b) \le v_n \le v_1 \le g(b)$. Let $\varepsilon > 0$ be a real number, and define $\bm{u}^{\varepsilon}, \bm{v}^{\varepsilon} \in \mathbb{W}_n$ by setting $u_j^{\varepsilon} = u_j - j \varepsilon$ and $v_j^{\varepsilon} = v_j - j \varepsilon$ for each $j \in \llbracket 1, n \rrbracket$; also define $f^{\varepsilon}, g^{\varepsilon} : [a, b] \rightarrow \mathbb{R}$ by setting $f^{\varepsilon} (s) = f(s) - (n+1) \varepsilon$ and $g^{\varepsilon} (s) = g(s)$ for each $s \in [a, b]$. Sample $\bm{\mathsf{x}}^{\varepsilon} = (\mathsf{x}_1^{\varepsilon}, \mathsf{x}_2^{\varepsilon}, \ldots , \mathsf{x}_n^{\varepsilon}) \in \llbracket 1, n \rrbracket \times \mathcal{C} \big( [a, b] \big)$ from the measure $\mathsf{Q}_{f^{\varepsilon}; g^{\varepsilon}}^{\bm{u}^{\varepsilon}; \bm{v}^{\varepsilon}}$. As $\varepsilon$ tends to $0$, $\bm{\mathsf{x}}^{\varepsilon}$ converges in law to a $\llbracket 1, n \rrbracket \times [a, b]$-indexed line ensemble $\bm{\mathsf{x}} \in \llbracket 1, n \rrbracket \times \mathcal{C} \big( [a, b] \big)$.
	
\end{lem}

\begin{proof}
	
	It quickly follows from monotonicity (\Cref{uvv}) that the finite-dimensional distributions of $\bm{\mathsf{x}}^{\varepsilon}$ converge, as $\varepsilon$ tends to $0$. Thus, it suffices to verify the tightness of $(\bm{\mathsf{x}}^{\varepsilon})_{\varepsilon> 0}$, equivalently, we have for any fixed $\delta > 0$ that	 
	\begin{flalign}
		\label{xisxit} 
		\displaystyle\lim_{\omega \rightarrow 0} \displaystyle\sup_{\varepsilon > 0} \mathbb{P} \Bigg[ \displaystyle\max_{i \in \llbracket 1, n \rrbracket} \displaystyle\sup_{\substack{s, t \in [a, b] \\ |s-t| < \omega}} \big| \mathsf{x}_i (s) - \mathsf{x}_i (t) \big| > \delta \Bigg] = 0.
	\end{flalign} 
	
	To do so, define $\bm{u}^{\varepsilon; \delta}, \bm{v}^{\varepsilon; \delta} \in \mathbb{W}_n$ by setting $u_j^{\varepsilon; \delta} = u_j^{\varepsilon} + (n-j) \delta / 4n$ and $v_j^{\varepsilon; \delta} = v_j^{\varepsilon} + (n-j) \delta / 4n$ for each $j \in \llbracket 1, n \rrbracket$, and define $f^{\varepsilon; \delta}, g^{\varepsilon; \delta} : [a, b] \rightarrow \mathbb{R}$ by setting $f^{\varepsilon; \delta} (s) = f^{\varepsilon} (s) - \delta/4n$ and $g^{\varepsilon; \delta} (s) = g^{\varepsilon} (s) + \delta/4$ for each $s \in [a, b]$. Then, we have $f^{\varepsilon; \delta} (a) \le  u_n^{\varepsilon; \delta} - \delta/4n \le u_1^{\varepsilon; \delta} + \delta/4n \le g^{\varepsilon; \delta} (a)$ and $f^{\varepsilon; \delta} (b) \le v_n^{\varepsilon; \delta} - \delta/4n \le v_1^{\varepsilon; \delta} + \delta/4n \le g^{\varepsilon; \delta} (b)$.
	
	Sample non-intersecting Brownian bridges $\bm{\mathsf{x}}^{\varepsilon; \delta} \in \llbracket 1, n \rrbracket \times \mathcal{C} \big( [a,b]\big)$ under $\mathsf{Q}_{f^{\varepsilon; \delta}; g^{\varepsilon; \delta}}^{\bm{u}^{\varepsilon; \delta}; \bm{v}^{\varepsilon; \delta}}$. By \Cref{uvv}, there is a coupling between $\bm{\mathsf{x}}^{\varepsilon}$ and $\bm{\mathsf{x}}^{\varepsilon; \delta}$ so that we almost surely have 
	\begin{flalign}
		\label{xydelta} 
		\big| \mathsf{x}_j^{\varepsilon} (s) - \mathsf{x}_j^{\varepsilon;\delta} (s) \big| \le \displaystyle\frac{\delta}{4}, \qquad \text{for each $s \in [a, b]$}. 
	\end{flalign} 
	
	\noindent Next, let $\bm{\mathsf{B}} = (\mathsf{B}_1, \mathsf{B}_2, \ldots , \mathsf{B}_n) \in \llbracket 1, n \rrbracket \times \mathcal{C} \big( [a, b] \big)$ denote $n$ mutually independent Brownian bridges on $[a, b]$, starting at $\bm{\mathsf{B}}(a) = \bm{u}^{\varepsilon; \delta}$ and ending at $\bm{\mathsf{B}} (b) = \bm{v}^{\varepsilon; \delta}$. Since the points in $\bm{u}^{\varepsilon; \delta} \cup \big\{ f^{\varepsilon; \delta} (a), g^{\varepsilon; \delta} (a) \big\}$ are bounded away from each other, and as are points in $\bm{v}^{\varepsilon; \delta} \cup \big\{ f^{\varepsilon; \delta} (b), g^{\varepsilon; \delta} (b) \big\}$ (all independently of $\varepsilon$), there exists a constant $c = c(\delta, \bm{u}, \bm{v}, f, g) > 0$ (that is in particular independent of $\varepsilon$) such that 
	\begin{flalign*} 
		\mathbb{P} \bigg[ \bigcap_{s \in [a, b]} \big\{ f^{\varepsilon; \delta} (s) < \mathsf{B}_n (s) < \mathsf{B}_{n-1} (s) < \cdots < \mathsf{B}_1 (s) < g^{\varepsilon; \delta} (s) \big\} \bigg] > c.
	\end{flalign*} 
	
	\noindent It follows that  
	\begin{flalign*} 
		\displaystyle\lim_{\omega \rightarrow 0} & \displaystyle\sup_{\varepsilon > 0} \mathbb{P} \Bigg[ \displaystyle\max_{i \in \llbracket 1, n \rrbracket} \displaystyle\sup_{\substack{s, t \in [a, b] \\ |s-t| < \omega}} \big| \mathsf{x}_i^{\varepsilon; \delta} (s) - \mathsf{x}_i^{\varepsilon; \delta} (t) \big| > \displaystyle\frac{\delta}{2} \Bigg] \\
		& \qquad \le c^{-1} \cdot  \displaystyle\lim_{\omega \rightarrow 0} \displaystyle\sup_{\varepsilon > 0}  \mathbb{P} \Bigg[ \displaystyle\max_{i \in \llbracket 1, n \rrbracket} \displaystyle\sup_{\substack{s, t \in [a, b] \\ |s-t| < \omega}} \big| \mathsf{B}_i (s) - \mathsf{B}_i (t) \big| > \displaystyle\frac{\delta}{2} \Bigg] = 0,
	\end{flalign*} 
	
	\noindent where in the last equality we used the fact that a Brownian bridge is almost surely uniformly continuous. Together with \eqref{xydelta}, this yields \eqref{xisxit} upon letting $\delta$ tend to $0$, showing the lemma. 
\end{proof}

We next deduce \Cref{monotoneheight} as a direct consequence of the same monotonicity result (due to \cite{PLE}) when $\bm{u}, \widetilde{\bm{u}}, \bm{v}, \widetilde{\bm{v}}$ are strictly ordered.

\begin{lem}[{\cite[Lemmas 2.6 and 2.7]{PLE}}] 
		
		\label{monotoneheightwn}
		
		If $\bm{u}, \widetilde{\bm{u}}, \bm{v}, \widetilde{\bm{v}} \in \mathbb{W}_n$, then \Cref{monotoneheight} holds.
		
	\end{lem}

	\begin{proof}[Proof of \Cref{monotoneheight}]
		
		For any $\varepsilon > 0$, define the $n$-tuples $\bm{u}^{\varepsilon}, \widetilde{\bm{u}}^{\varepsilon}, \bm{v}^{\varepsilon}, \widetilde{\bm{v}}^{\varepsilon} \in \mathbb{W}_n$ by setting 
		\begin{flalign*} 
			u_j^{\varepsilon} = u_j - j \varepsilon; \qquad \widetilde{u}_j^{\varepsilon} = \widetilde{u}_j - j \varepsilon; \qquad v_j^{\varepsilon} = v_j - j \varepsilon; \qquad \widetilde{v}_j^{\varepsilon} = \widetilde{v}_j - j \varepsilon,
		\end{flalign*}
		
		\noindent for each $j \in \llbracket 1, n \rrbracket$. Similarly, define the functions $f^{\varepsilon}, \widetilde{f}^{\varepsilon}, g^{\varepsilon}, \widetilde{g}^{\varepsilon} : [a, b] \rightarrow \overline{\mathbb{R}}$ by setting 
		\begin{flalign*} 
			f^{\varepsilon} (s) = f(s); \qquad \widetilde{f}^{\varepsilon} (s) = \widetilde{f} (s) - (n+1) \varepsilon; \qquad g^{\varepsilon} = g(s)  ; \qquad \widetilde{g}^{\varepsilon} (s) = \widetilde{g} (s),
		\end{flalign*}
		
		\noindent for each $s \in [a, b]$. Sample non-intersecting Brownian bridges $\bm{\mathsf{x}} = (\mathsf{x}_1^{\varepsilon}, \mathsf{x}_2^{\varepsilon}, \ldots , \mathsf{x}_n^{\varepsilon}) \in \llbracket 1, n \rrbracket \times \mathcal{C} \big( [a, b] \big)$ and $\widetilde{\bm{\mathsf{x}}} = (\widetilde{\mathsf{x}}_1^{\varepsilon}, \widetilde{\mathsf{x}}_2^{\varepsilon}, \ldots , \widetilde{\mathsf{x}}_n^{\varepsilon}) \in \llbracket 1, n \rrbracket \times \mathcal{C} \big( [a, b] \big)$ according to the measures $\mathsf{Q}_{f^{\varepsilon}; g^{\varepsilon}}^{\bm{u}^{\varepsilon}; \bm{v}^{\varepsilon}}$ and  $\mathsf{Q}_{\tilde{f}^{\varepsilon}; \tilde{g}^{\varepsilon}}^{\tilde{\bm{u}}^{\varepsilon}; \tilde{\bm{v}}^{\varepsilon}}$, respectively. 
		
		Since $\bm{u}^{\varepsilon}, \bm{v}^{\varepsilon}, \widetilde{\bm{u}}^{\varepsilon}, \widetilde{\bm{v}}^{\varepsilon} \in \mathbb{W}_n$ (and $\bm{u}^{\varepsilon} \le \widetilde{\bm{u}}^{\varepsilon}$, $\bm{v}^{\varepsilon} \le \widetilde{\bm{v}}^{\varepsilon}$, $f^{\varepsilon} \le \widetilde{f}^{\varepsilon}$, and $g^{\varepsilon} = \widetilde{g}^{\varepsilon}$), \Cref{monotoneheightwn} applies and yields for each $\varepsilon > 0$ a coupling between $\bm{\mathsf{x}}$ and $\widetilde{\bm{\mathsf{x}}}^{\varepsilon}$ such that $\mathsf{x}_j^{\varepsilon} (t) \le \widetilde{\mathsf{x}}_j^{\varepsilon} (t)$ holds for all $(j, t) \in \llbracket 1, n \rrbracket \times [a, b]$, almost surely. Since the measures  $\mathsf{Q}_{f^{\varepsilon}; g^{\varepsilon}}^{\bm{u}^{\varepsilon}; \bm{v}^{\varepsilon}}$ and  $\mathsf{Q}_{\tilde{f}^{\varepsilon}; \tilde{g}^{\varepsilon}}^{\tilde{\bm{u}}^{\varepsilon}; \tilde{\bm{v}}^{\varepsilon}}$ converge to $\mathsf{Q}_{f; g}^{\bm{u}; \bm{v}}$ and  $\mathsf{Q}_{\tilde{f}; \tilde{g}}^{\tilde{\bm{u}}; \tilde{\bm{v}}}$, respectively, the couplings between $( \bm{\mathsf{x}}^{\varepsilon}, \widetilde{\bm{\mathsf{x}}}^{\varepsilon})$ are tight in $\varepsilon > 0$. Taking any limit point of these couplings along a subsequence of $\varepsilon$ tending to $0$ (and using \Cref{boundaryconvergeensemble}), we obtain a coupling between $(\bm{\mathsf{x}}, \widetilde{\bm{\mathsf{x}}})$ satisfying the conditions of the lemma. 
		\end{proof}

	Next we establish \Cref{yconvergedelta} as a quick consequence of a result from \cite{FESD}. To state the latter, we require the following assumption, which is a more qualitative variant of \Cref{yconvergedelta0}.
	
	\begin{assumption} 
		
		\label{nup0} 
		
		Fix a real number $t > 0$ and a measure $\nu \in \mathscr{P}_0$ such that
		\begin{flalign}
			\label{e:b1}
			\displaystyle\inf_{s \in \supp \nu} \displaystyle\lim_{\varepsilon \rightarrow 0} \displaystyle\int_{-\infty}^{\infty} \displaystyle\frac{\nu(dx)}{(s-x)^2 + \varepsilon^2} > t^{-1}.
		\end{flalign} 
		
		\noindent For each integer $n \ge 1$ let $\bm{y} = \bm{y}^n=(y_1, y_2,\cdots, y_n) \in \overline{\mathbb{W}}_n$ be a sequence. Assume that it satisfies the following two conditions.
		\begin{enumerate}
			\item The measures $\nu_n = n^{-1} \sum_{j=1}^n \delta_{y_j/n}$ converge weakly to $\nu$, as $n$ tends to $\infty$.
			\item We have $\lim_{n \rightarrow \infty} \max_{1 \le j \le n} \dist (n^{-1} y_j, \supp \nu) = 0$.
		\end{enumerate}
		
		\noindent For each integer $n \ge 1$, let $\bm{\lambda} = \bm{\lambda}^n \in \overline{\mathbb{W}}_n$ denote Dyson Brownian motion run for time $tn$, with initial data $\bm{y}^n$, and define $\sigma = \sigma_{\nu;t}$ as in \eqref{z0tsigma}.
		
	\end{assumption} 
	
	The following result from \cite{FESD} indicates that the largest particles of $\bm{\lambda}$ converge to the Airy point process. The convergence to the Airy point process follows from  \cite[Theorem 1.1]{FESD} (after scaling the measure and its argument $\nu$ by $t^{1/2}$), and the explicit form of the scaling factor $\sigma$ follows from \cite[Section 4.2.1]{FESD}, together with \cite[Lemmas 3.1 and 3.4]{FESD}.\footnote{In \cite[Theorem 1.1]{FESD}, \eqref{convergelambda} is stated with $\alpha$ replaced by $\alpha^{-1}$. This is a misprint, stemming from a corresponding one when changing of variables to pass from \cite[Equation (49)]{FESD} to the following ones. Numerous other works have also proved edge statistics results in various different regimes, and they showed that the scaling appears as we have written in \eqref{convergelambda}; see, for example, \cite[Equation (17) and Theorem 2(iii)]{ULER} and \cite[Theorem 2.2 and Equation (2.12)]{landon2017edge}.}
	
	\begin{lem}[{\cite{FESD}}]
		
		\label{sumyw}
		
		Adopting \Cref{nup0}, for any integer $n \ge 1$, there exists a real number $E_n$ such that the following holds for any fixed integer $k \ge 1$. As $n$ tends to $\infty$, the sequence 
		\begin{flalign}
			\label{convergelambda0}
			\big(\sigma n^{-1/3} (\lambda_1-E_n), \sigma n^{-1/3} (\lambda_2-E_n), \ldots , \sigma n^{-1/3} (\lambda_k - E_n) \big), \quad \text{converges to} \quad (\mathfrak{a}_1, \mathfrak{a}_2, \ldots , \mathfrak{a}_k),
		\end{flalign}
		
		\noindent in law, where the latter is given by the first $k$ points of the Airy point process (recall \Cref{a0}).
		
	\end{lem}

	\begin{proof}[Proof of \Cref{yconvergedelta}]
		
		Assume to the contrary that \eqref{convergelambda} is false. Then, there exists a real number $\varepsilon > 0$ and a sequence $n_1 < n_2 < \cdots $ of positive integers such that 
		\begin{flalign}
			\label{ynk0}
			\Bigg| \mathbb{P} \bigg[ \bigcap_{j=1}^k \big\{ \sigma n_i^{-1/3} (\lambda_j^{n_i} - \lambda_{j+1}^{n_i}) \ge r_j \big\} \bigg] - \mathbb{P} \bigg[ \bigcap_{j=1}^k \big\{ \mathfrak{a}_j - \mathfrak{a}_{j+1} \ge r_j \}\bigg] \Bigg| > \varepsilon, 
		\end{flalign}
		
		\noindent for each integer $i \ge 1$. By \eqref{yndelta}, the sequence $(\bm{y}^{n_k})$ satisfies the two conditions in \Cref{nup0}, and \eqref{yndelta1} verifies \eqref{e:b1}. Therefore, \Cref{sumyw} applies\footnote{While this sequence only involves elements of $\overline{\mathbb{W}}_{n_k}$, one can define $\widetilde{\bm{y}}^n \in \overline{\mathbb{W}}_n$ for each integer $n \ge 1$ by setting $\widetilde{\bm{y}}^n = \bm{y}_{n_k}$ for each integer $n \in [n_{k-1}, n_k)$ (where $n_{k-1} = 1$) and apply \Cref{sumyw} to the $(\widetilde{\bm{y}}^n)$.} and yields real numbers $E_{n_i}$ such that $\big( \sigma n_i^{-1/3} (\lambda_j^{n_i} - E_{n_i}) \big)_{j \in \llbracket 1, k + 1 \rrbracket}$ converges in law to $(\mathfrak{a}_1, \mathfrak{a}_2, \ldots , \mathfrak{a}_{k+1})$, as $i$ tends to $\infty$. Thus, 
		\begin{flalign*}
			\displaystyle\lim_{i \rightarrow \infty} \mathbb{P} \bigg[ \bigcap_{j=1}^k \big\{ \sigma n_i^{-1/3} (\lambda_j^{n_i} - \lambda_{j+1}^{n_i}) \ge r_j \big\} \bigg] = \mathbb{P} \bigg[ \bigcap_{j=1}^k \{ \mathfrak{a}_j - \mathfrak{a}_{j+1} \ge r_j \} \bigg],
		\end{flalign*} 
		
		\noindent which contradicts \eqref{ynk0}, thereby establishing the lemma.		
	\end{proof}

	\subsection{Proofs of \Cref{t:lawt} and \Cref{gammaderivative}} 
	
	\label{ProofBridgeSum1} 
	
	In this section we establish first \Cref{t:lawt} and then \Cref{gammaderivative}.

	\begin{proof}[Proof of \Cref{t:lawt}]
		
		By the second part of \Cref{lambdat}, the law of $ \bm{\mathsf{x}} (t)$ is given by Dyson Brownian motion run for time $t$, with initial data $\bm{u}$, conditioned to end at $\bm{v}$ at time $\mathsf{T}$. Denoting by $\bm{H} (s) = \bm{H}_n (s)$ an $n \times n$ Hermitian Brownian motion, the first part of \Cref{lambdat} implies that the latter process is given by $\eig \big( \bm{U} + \bm{H} (s) \big)$, where $\bm{U} + \bm{H} (s)$ is conditioned to be of the form $\bm{W} \bm{V} \bm{W}^*$ at time $s = \mathsf{T}$, for some unitary matrix $\bm{W} \in \mathrm{U} (n)$. 
		
		Since the entries of $\bm{H} (s)$ are complex Gaussian random variables of variance $s$, the density of $\bm{U} + \bm{H} (\mathsf T)$ is proportional to 
		\begin{flalign*} 
			\exp \Big( -\displaystyle\frac{1}{2\mathsf T} \Tr \bm{H}(\mathsf{T})^2 \Big) d\bm{H}(\mathsf{T}) = \exp \Big( -\displaystyle\frac{1}{2\mathsf T} \Tr (\bm{W} \bm{V} \bm{W}^* - \bm{U})^2 \Big) d(\bm{W} \bm{V} \bm{W}^*).
		\end{flalign*}
		
		\noindent Upon conditioning on the eigenvalues of $\bm{W} \bm{V} \bm{W}^* $ (and dividing by the constant $e^{\mathsf{T}^{-1} \Tr (\bm{U}^2 + \bm{V}^2)}$), the above density is proportional to \eqref{wuv}, which therefore prescribes the law of the unitary matrix $\bm{W}$.
		
		Hence, denoting the $(i, j)$ entry of any matrix $\bm{M}$ by $(\bm{M})_{ij}$, the law of the upper triangular entries $\big(\bm{U}+ \bm{H} (s) \big)_{ij}$ (for $1 \le i \le j \le n$) are given by Brownian bridges conditioned to start at $(\bm{U})_{ij}$ (at time $t = 0$) and end at $(\bm{W} \bm{V} \bm{W}^*)_{ij}$ (at time $s = \mathsf{T})$. Since any Brownian bridge $B: [0, \mathsf{T}] \rightarrow \mathbb{R}$ with $B(0) = a$ and $B(\mathsf{T}) = b$ can be represented as 
		\begin{flalign*} 
			B(s) = \displaystyle\frac{\mathsf{T} - s}{\mathsf{T}} \cdot a + \displaystyle\frac{s}{\mathsf{T}} \cdot b + \displaystyle\frac{\mathsf{T} - s}{\mathsf{T}^{1/2}} \cdot Y \Big( \displaystyle\frac{s}{\mathsf{T} - s} \Big), 
		\end{flalign*} 
	
		\noindent for some Brownian motion $Y : \mathbb{R}_{\ge 0} \rightarrow \mathbb{R}$, it follows that $\bm{\mathsf{x}} (t)$ has the same law as
		\begin{flalign*}
			\eig \bigg( \displaystyle\frac{\mathsf{T}-t}{\mathsf{T}} \cdot \bm{U} + \displaystyle\frac{s}{\mathsf{T}} \cdot \bm{W} \bm{V} \bm{W}^* + \displaystyle\frac{\mathsf{T}-t}{\mathsf{T}^{1/2}} \cdot \bm{G} \Big( \displaystyle\frac{t}{\mathsf{T}-t} \Big) \bigg),
		\end{flalign*}
		
		\noindent where $\bm{W}$ is sampled under \eqref{wuv} and $\bm{G} \big( t / (\mathsf{T}-t) \big)$ is an independent Hermitian Brownian motion run for time $t (\mathsf{T}-t)^{-1}$. The lemma then follows from the fact that $\bm{G} \big( t / (\mathsf{T}-t) \big)$ has the same law as $t^{1/2} (\mathsf{T} - t)^{-1/2} \cdot \bm{G}$. 
	\end{proof}

			\begin{proof}[Proof of \Cref{gammaderivative}] 
		
		Observe from \eqref{gammascy} that, for $y \in [0, 1]$, we have   
		\begin{flalign}
			\label{gamma1} 
			(2 \pi)^{-1} \displaystyle\int_{\gamma_{\semci} (y)}^2 (4-w^2)^{1/2} dw = y.
		\end{flalign}
		
		\noindent By \eqref{gamma1} and the symmetry of the integrand $(4-w^2)^{1/2}$ there in $w$, we have $0 \le \gamma_{\semci} (y) \le 2$ for $y \in [0, 1/2]$ and $-2 \le \gamma_{\semci} (y) \le 0$ for $y \in [1/2, 1]$. The latter verifies the first statement of the lemma when $y \in [1/2, 1]$, so let assume that $y \in [0, 1/2)$ so $0 \le \gamma_{\semci} (y) \le 2$. Then the first part of the lemma follows from the fact that, for any real number $\theta \in [0, 2]$ we have
		\begin{flalign*}
			\Big( \displaystyle\frac{\theta}{8} \Big)^{3/2} \le \displaystyle\frac{2^{1/2} \theta^{3/2}}{3 \pi} & = 2^{-1/2} \pi^{-1} \displaystyle\int_{2-\theta}^2 (2-w)^{1/2} dw \\ 
			& \le  (2 \pi)^{-1} \displaystyle\int_{2-\theta}^2 (4-w^2)^{1/2} dw \le \pi^{-1} \displaystyle\int_{2-\theta}^2 (2-w)^{1/2} dw = \displaystyle\frac{2\theta^{3/2}}{3\pi} \le \Big( \displaystyle\frac{\theta}{2} \Big)^{3/2},
		\end{flalign*}
		
		\noindent where we used the bound $2(2-w) \le 4-w^2 \le 4(2-w)$ for $w \in [0, 2]$. 
		
		To establish the second part of the lemma, we differentiate  \eqref{gamma1} with respect to $y$ to obtain
		\begin{flalign}
			\label{derivative2}
			-\gamma'_{\semci} (y)=2\pi \big( 4-\gamma_{\semci} (y)^2 \big)^{-1/2}.
		\end{flalign}
		
		\noindent Since $0\leq \gamma_{\semci} (y)\leq 2$ for $y \in [0, 1/2]$, by the first part of the lemma we have $4y^{2/3}\leq 4-\gamma^2_{\semci} (y)\leq 32y^{2/3}$ for $0\leq y\leq 1/2$ and $4-\gamma^2_{\semci} (y)\leq 32y^{2/3}$ for $0\leq y\leq 1$. Together with \eqref{derivative2}, these estimates yield the second part of the lemma.
	\end{proof}

	\subsection{Proofs of Corollary \ref{p:smallinitial} and Lemma \ref{p:compareAiry}}
	
	\label{ProofBridgeSum2} 
	
	In this section we establish first \Cref{p:smallinitial} and then \Cref{p:compareAiry}.

	\begin{proof}[Proof of \Cref{p:smallinitial}]
		
		We claim that there exist constants $c = c(A, B) >  0$ and $C = C (A, B) > 1$ such that, for any fixed real number $t \in [T / 4, 3T / 4]$, we have 
		\begin{flalign}
			\label{xjkn13} 
			\mathbb{P} \Bigg[ \bigcup_{1 \le j < k \le \lfloor n/2 \rfloor} \Big\{ \mathsf{x}_j (tn^{1/3}) - \mathsf{x}_k (tn^{1/3}) \le C (k^{2/3} - j^{2/3}) + (\log n)^{24} j^{-1/3} \Big\} \Bigg] \le c^{-1} e^{-c(\log n)^2}.
		\end{flalign}
		
		\noindent We first establish the corollary assuming \eqref{xjkn13}. To this end, define the set $\mathcal{T} = (n^{-9} \cdot \mathbb{Z}) \cap [0, Tn^{1/3}]$ and the events
		\begin{flalign*}
			& \mathscr{E}_1 = \bigcap_{j=1}^k \bigcap_{0 \le r < r + s \le Tn^{1/3}} \Big\{ \big| \mathsf{x}_j (r+s) - \mathsf{x}_j (r) \big| \le 2n^2 s^{1/3} \Big\}; \\ 
			& \mathscr{E}_2 = \bigcap_{t \in \mathcal{T}} \bigcap_{1 \le j < k \le \lfloor n/2 \rfloor} \big\{ \mathsf{x}_j (tn^{1/3}) - \mathsf{x}_k (tn^{1/3}) \le C (k^{2/3} - j^{2/3}) + (\log n)^{24} j^{-1/3} \big\}.
		\end{flalign*} 
		
		\noindent We then claim that there exists a constant $c_0 = c_0 (A, B) > 0$ such that
		\begin{flalign} 
			\label{probabilitysrxn13} 
			\mathbb{P} [\mathscr{E}_1^{\complement}] \le c_0^{-1} e^{-c_0 (\log n)^2}; \qquad \mathbb{P} \big[ \mathscr{E}_2^{\complement} \big] \le c_0^{-1} e^{-c_0 (\log n)^2}. 
		\end{flalign}
		
		\noindent Indeed, the first bound in \eqref{probabilitysrxn13} follows from the $B = n$ case of \Cref{estimatexj2}, together with the facts that for sufficiently large $n$ we have that $|v_j - u_j| \le 2Bn^{2/3} \le n$ (by \eqref{buvb}), that $s (Tn^{1/3})^{-1} \le s^{1/3} n$ (for $s \in [0, Tn^{1/3}]$), and that $ns^{1/2} \log (2s^{-1} T n^{1/3}) \le n^2 s^{1/3}$. The second bound in \eqref{probabilitysrxn13} follows from taking a union bound in \eqref{xjkn13} over $t \in \mathcal{T}$ (and using the fact that $|\mathcal{T}| \le 3An^{10}$).
		
		Now restrict to the event $\mathscr{E}_1 \cap \mathscr{E}_2$. Fix $s \in [Tn^{1/3}/ 4, 3Tn^{1/3} / 4]$ and let $s' \in \mathcal{T}$ be the closest number in $\mathcal{T}$ to $s$ (if more than one exists, we select one arbitrarily). Then, for any integers $1 \le j < k \le n$, we have
		\begin{flalign*}
			\mathsf{x}_j (s) - \mathsf{x}_k (s) & \le \big| \mathsf{x}_j (s) - \mathsf{x}_j (s') \big| + \big| \mathsf{x}_j (s') - \mathsf{x}_k (s') \big| + \big| \mathsf{x}_k (s') - \mathsf{x}_k (s) \big| \\
			& \le  C (k^{2/3} - j^{2/3}) + (\log n)^{24} j^{-1/3} + 4n^2 |s-s'|^{1/3} \\
			& \le C (k^{2/3} - j^{2/3}) + (\log n)^{24} j^{-1/3} + 4n^{-1} \le C (k^{2/3} - j^{2/3}) + (\log n)^{25} j^{-1/3},
		\end{flalign*} 
		
		\noindent where in the second bound we used the fact that we are restricting to $\mathscr{E}_1 \cap \mathscr{E}_2$; in the third we used the fact that $|s-s'| \le n^{-9}$ (since $s \in [0, Tn^{1/3}]$ and $\mathcal{T} = (n^9 \cdot \mathbb{Z}) \cap [0, Tn^{1/3}]$); and in the fourth we used the fact that $(\log n)^{24} j^{-1/3} + 4n^{-1} \le (\log n)^{25} j^{-1/3}$ for $j \in \llbracket 1, n \rrbracket$ and sufficiently large $n$. This, together with the fact that $\mathbb{P} [\mathscr{E}_1 \cap \mathscr{E}_2] \ge 1 - 2c_0 ^{-1} e^{-c_0 (\log n)^2}$ (by \eqref{probabilitysrxn13} and a union bound), implies \eqref{xjk2n13c} and thus the corollary.
		
		It therefore remains to establish \eqref{xjkn13}. To this end, define the real numbers $t_0, T_0 > 0$; the $n$-tuples $\bm{u}', \bm{v}' \in \overline{\mathbb{W}}_n$; and the ensemble $\bm{\mathsf{y}} = (\mathsf{y}_1, \mathsf{y}_2, \ldots , \mathsf{y}_n) \in \llbracket 1, n \rrbracket \times \mathcal{C} \big( [0, T_0 n^{1/3}] \big)$ by for any $s \in [0, T_0 n^{1/3}]$ setting 
		\begin{flalign}
			\label{ttuvx}
			\begin{aligned}
				t_0 = \bigg( t \Big(1 - \frac{t}{T} \Big) \bigg)^{1/2}; \quad & T_0 = t_0^{-2} T; \quad \bm{u}' = t_0^{-1} \cdot \bm{u}; \quad \bm{v}' = t_0^{-1} \cdot \bm{v}; \quad \mathsf{y}_j (s) = t_0^{-1} \cdot \mathsf{x}_j (t_0^2 s).
			\end{aligned} 
		\end{flalign}
		
		\noindent By \Cref{scale}, the law of $\bm{\mathsf{y}}$ is given by $\mathsf{Q}^{\bm{u}'; \bm{v}'}$. Next, denote the $n \times n$ diagonal matrices $\bm{U} = \diag (\bm{u}')$ and $\bm{V} = \diag (\bm{v}')$; letting the unitary random matrix $\bm{W} \in \mathrm{U}(n)$ have law \eqref{wuv}, set $\bm{A} = (T_0 - t_0) T_0^{-1} \cdot \bm{U} + t_0 T_0^{-1} \cdot \bm{W} \bm{V} \bm{W}^*$ and $\bm{a} = \eig (\bm{A})$. Then, by \Cref{tuvwx}, the law of $\bm{\mathsf{y}} (t_0^{-2} t n^{1/3}) = t_0^{-1} \cdot \bm{\mathsf{x}} (t n^{1/3})$ is given by $\bm{\lambda}(n^{1/3})$, where $\bm{\lambda}(s) = \big( \lambda_1 (s), \lambda_2 (s), \ldots , \lambda_n (s) \big) \in \overline{\mathbb{W}}_n$ denotes Dyson Brownian motion run for time $s$ with initial data $\bm{\lambda} (0) = \bm{a}$.
		
		We analyze $\bm{\lambda} (n^{1/3})$ using \Cref{initialsmall2}. By the Weyl interlacing inequality, we have
		\begin{flalign*}
			\max \bm{a} \le \max \eig (\bm{U}) + \max \eig(\bm{V}) = t_0^{-1} (\max \bm{u} + \max \bm{v}) \le 2 B t_0^{-1} n^{2/3},
		\end{flalign*} 
		
		\noindent and similarly $\min \bm{a} \ge -2B t_0^{-1} n^{2/3}$. Let $c_1$ denote the constant $c(2) > 0$ from \Cref{initialsmall2}. Observe since $t \in [ T / 4, 3T / 4]$ that $t_0 \ge T^{1/2} / 4$ and $T \ge C_1$ that we can make $4Bt_0^{-1} < c_1$ by taking $C_1 = C_1 (B) > 1$ sufficiently large. Then, $\max \bm{a} - \min \bm{a} \le 4 B t_0^{-1} n^{2/3} < c_1 n^{2/3}$, and so \Cref{initialsmall2} applies and yields a constant $c_2 > 0$ such that 
		\begin{flalign*}
			\mathbb{P} \Bigg[ \bigcup_{1 \le j < k \le \lfloor n/2 \rfloor} \Big\{ \big| \lambda_j (n^{1/3}) - \lambda_k (n^{1/3}) \ge 25 (k^{2/3} - j^{2/3}) + (\log n)^{20} j^{-1/3} \Big\} \Bigg] \le c_2^{-1} e^{-c_2 (\log n)^2}. 
		\end{flalign*} 
		
		\noindent Since $\bm{\lambda} (n^{1/3})$ has the same law as $\bm{\mathsf{y}} (t_0^{-2} t n^{1/3}) = t_0^{-1} \bm{\mathsf{x}} (t n^{1/3})$, and since $t_0 < T^{1/2} \le (AC_1)^{1/2}$, it follows that 
		\begin{flalign*}
			\mathbb{P} \Bigg[ \bigcap_{1 \le j < k \le \lfloor n / 2 \rfloor} \big\{ \mathsf{x}_j (tn^{1/3}) - \mathsf{x}_k (tn^{1/3}) \ge 25 (AC_1)^{1/2} (k^{2/3} - j^{2/3}) + (AC_1)^{1/2} (& \log n)^{20} j^{-1/3} \big\} \Bigg] \\
			& \le c_2^{-1} e^{-c_2 (\log n)^2}, 
		\end{flalign*}
		
		\noindent from which \eqref{xjkn13} follows, as $(AC_1)^{1/2} (\log n)^{20} \le (\log n)^{24}$ for sufficiently large $n$.			
	\end{proof}

			\begin{proof}[Proof of \Cref{p:compareAiry}]

				We will establish the lemma by comparing the non-intersecting Brownian bridges $\bm{\mathsf{x}}$ with certain rescaled parabolic Airy line ensembles and Brownian watermelons; we will prove the first part of the lemma in detail and only outline the proof for the second part, as it is fairly similar. In what follows, for any real number $\sigma > 0$, we recall the rescaled parabolic Airy line ensemble $\bm{\mathcal{S}}^{(\sigma)} = \big( \mathcal{S}_1^{(\sigma)}, \mathcal{S}_2^{(\sigma)}, \ldots \big)$ from \eqref{sigmar}. For any integer $n \ge 1$, \Cref{kdeltad} and a union bound (with the $u$ there equal to $(\log n)^2$ here) together yield a constant $c_3 = c_3 (\sigma, D) > 0$ such that 
				\begin{align}
					\label{anrsigmaprobability} 
					\mathbb{P} \Bigg[ \bigcup_{j=1}^{2n} \bigcup_{t \in [-n^D, n^D]} \bigg|  \mathcal{S}_j^{(\sigma)} (s) + 2^{-1/2} \sigma^3 s^2 + \displaystyle\frac{(3 \pi)^{2/3} j^{2/3}}{2^{7/6} \sigma} \bigg| \ge (\log n)^2j^{-1/3} \Bigg] \le c_3^{-1} e^{-c_3 (\log n)^2}.
				\end{align}				
				
				We begin by verifying the first part of the lemma. To this end, take $\sigma_1=2^{-7/6} (3\pi)^{2/3} d^{-1} $, and denote the line ensemble $\widetilde{\bm{\mathcal{S}}} = \big( \widetilde{\mathcal{S}}_1, \widetilde{\mathcal{S}}_2, \ldots \big) \in \mathbb{Z}_{\ge 1} \times \mathcal{C} (\mathbb{R})$ by for each $(j, t) \in \mathbb{Z}_{\ge 1} \times \mathbb{R}$ setting  
				\begin{align}
					\label{rsigma1} 
					\widetilde{\mathcal{S}}_j(t) = \mathcal{S}^{(\sigma_1)}_j \Big( t- \displaystyle\frac{a+b}{2} \Big) +2^{-5/2} \sigma_1^3 (b-a)^2 +M+(\log n)^2,
				\end{align}
				
				\noindent which satisfies the Brownian Gibbs property by \Cref{sigmascale}. Now observe, due to the upper bounds assumed on $\bm{u}$, $\bm{v}$, and $f$, that for each $(j, t) \in \llbracket 1, n \rrbracket \times \{ a, b \}$ we have 
				\begin{flalign*}
					2^{-5/2} \sigma_1^3 (b-a)^2 + M - 2^{-1/2} \sigma_1^3 \Big( t - \displaystyle\frac{b+a}{2} \Big)^2 - \displaystyle\frac{(3\pi)^{2/3} j^{2/3}}{2^{7/6} \sigma_1} = M - d j^{2/3} \ge \max \{ u_j, v_j \},
				\end{flalign*}
			
				\noindent and for each $t \in [a, b]$ (using the bound $\big( t - (a+b) / 2  \big)^2 \le (b-a)^2 / 4$ for $t \in [a, b]$) that
				\begin{flalign*}
					2^{-5/2} \sigma_1^3 (b-a)^2 + M - 2^{-1/2} \sigma_1^3 \Big( t - \displaystyle\frac{b+a}{2} \Big)^2 - \displaystyle\frac{(3\pi)^{2/3} (n+1)^{2/3}}{2^{7/6} \sigma_1} \ge M - d (n+1)^{2/3} \ge f(t).
				\end{flalign*}
			
				\noindent Thus, from \eqref{rsigma1} and \eqref{anrsigmaprobability} (with the translation-invariance of $\bm{\mathcal{S}}^{(\sigma_1)}$, which holds by \Cref{translationa}, to shift the interval $[-n^D, n^D]$ in \eqref{anrsigmaprobability} to one containing $[a, b]$ here), there exists a constant $c_4 = c_4 (d, D) > 0$ such that 
				\begin{flalign} 
					\label{e1probabilityr} 
					\mathbb{P} [\mathscr{E}_1] \ge 1 - c_4^{-1} e^{-c_4 (\log n)^2},
				\end{flalign} 
			
				\noindent where we have defined the event 
				\begin{flalign*}
					\mathscr{E}_1 = \bigcap_{j = 1}^n \big\{ \widetilde{\mathcal{S}}_j (a) \ge u_j \big\} \cap \big\{ \widetilde{\mathcal{S}}_j (b) \ge v_j \big\} \cap \bigcap_{t \in [a, b]} \big\{ \widetilde{\mathcal{S}}_{n+1} (t) \ge f(t) \big\}. 
				\end{flalign*}
				
				\noindent Condition on $\mathcal{S}_j (t)$ for $(j, t) \notin \llbracket 1, n \rrbracket \times [a, b]$ (that is, on $\mathcal{F}_{\ext}^{\bm{\mathcal{S}}} \big( \llbracket 1, n \rrbracket \times [a, b] \big)$ from \eqref{property}), and restrict to $\mathscr{E}_1$. Then, \Cref{monotoneheight} implies on $\mathscr{E}_1$ that we may couple $\bm{\mathsf{x}}$ and $\widetilde{\bm{\mathcal{S}}}$ such that $\mathsf{x}_j (t) \le \widetilde{\mathcal{S}}_j (t)$ for each $(j, t) \in \llbracket 1, n \rrbracket \times [a, b]$. This, \eqref{rsigma1}, \eqref{anrsigmaprobability} (again with the translation-invariance of $\bm{\mathcal{S}}^{(\sigma_1)}$), \eqref{e1probabilityr}, and the bound
				\begin{flalign*} 
					2^{-5/2} \sigma_1^3 (b-a)^2 + M - 2^{-1/2} \sigma_1^3 \Big( t - \displaystyle\frac{a+b}{2} \Big)^2 - \displaystyle\frac{(3\pi)^{2/3} j^{2/3}}{2^{7/6} \sigma_1} \le 2^{-5/2} \sigma_1^3 (b-a)^2 + M - dj^{2/3},
				\end{flalign*} 
			
				\noindent  yields  a constant $c_1 = c_1 (d, D) > 0$ such that 
				\begin{align*} 
					\mathbb{P} \Bigg[ \bigcap_{j = 1}^n \bigcap_{t \in [a, b]} \bigg\{ \mathsf{x}_j (t) \le M + 2^{-5/2} \sigma_1^3 (b-a)^2  - dj^{2/3} + 2 (\log n)^2 \bigg\} \Bigg] \ge 1- c_1^{-1} e^{-c_1 (\log n)^2}, 
				\end{align*}
			
				\noindent which with the definition of $\sigma_1$ gives \eqref{e:xjupbound}. 
				
				To establish the second part of the lemma, first observe that we may assume $f = -\infty$, by \Cref{monotoneheight}. Next define $\bm{u}', \bm{v}' \in \overline{\mathbb{W}}_n$ by setting $u_j' = u_n$ and $v_j' = v_n$ for each $j \in \llbracket 1, n \rrbracket$. Denote the associated Brownian watermelon $\bm{\mathsf{y}} = (\mathsf{y}_1, \mathsf{y}_2, \ldots , \mathsf{y}_n) \in \llbracket 1, n \rrbracket \times \mathcal{C} \big( [a, b] \big)$, given by $n$ non-intersecting Brownian bridges sampled from the measure $\mathsf{Q}^{\bm{u}'; \bm{v}'}$. Since $\bm{u}' \le \bm{u}$ and $\bm{v}' \le \bm{v}$, we may by \Cref{monotoneheight} couple $\bm{\mathsf{x}}$ and $\bm{\mathsf{y}}$ in such a way that $\mathsf{x}_j (t) \ge \mathsf{y}_j (t)$ for each $(j, t) \in \llbracket 1, n \rrbracket \times [a, b]$. Hence there exists a  constant $C = C(A) > 1$ such that, with probability at least $1 - C e^{-(\log n)^5}$, we have
				\begin{align*}
					\mathsf{y}_n(t)
					&\geq \frac{t -a}{b-a} \cdot v_n + \frac{b -t}{b-a} \cdot u_n - (b-a)^{1/2} n^{1/2} - (A + 1) (\log n)^9 \\
					&\geq -B n^{2/3}-M-(b-a)^{1/2}n^{1/2} - (A+1)(\log n)^9 \geq-(B+ A^{1/2} + 1) n^{2/3}-M,
				\end{align*}
			
				\noindent where in the first inequality we used the first part of \Cref{estimatexj} (with the facts that $b-a \le An^{1/3}$, that $(b-t)(t-a) \le (b-a)^2 / 4$, and that $\gamma_{\semci; n} (n) \ge -2$, by the first part of \Cref{gammaderivative}); in the second we used the fact that $\min \{ u_j, v_j \} \ge -Bj^{2/3} - M$ for each $j \in \llbracket 1, n \rrbracket$; and in the third we used the fact that $b-a \le An^{1/3}$. Together with the coupling $\mathsf{x}_j (t) \ge \mathsf{y}_j (t)$, this implies that
				\begin{flalign}
					\label{xjtba} 
					\mathbb{P} \Bigg[ \bigcap_{t \in [a, b]} \big\{ \mathsf{x}_n (t) \ge -(A^{1/2} + B + 1) n^{2/3} - M \big\} \Bigg] \ge 1 - Ce^{-(\log n)^5}. 
				\end{flalign}
				
				The claim \eqref{e:xjlowbound} follows by using \eqref{xjtba} to compare $\bm{\mathsf{x}}$ to another rescaled parabolic Airy line ensemble $\widehat{\bm{\mathcal{S}}} = \big( \widehat{\mathcal{S}}_1, \widehat{\mathcal{S}}_2, \ldots \big) \in \mathbb{Z}_{\ge 1} \times \mathcal{C}(\mathbb{R})$, defined by for any $(j, t) \in \mathbb{Z}_{\ge 1} \times \mathbb{R}$ setting 
				\begin{align*}
					\widehat{\mathcal{S}}_j(t) =  \mathcal{S}^{(\sigma_2)}_j \Big(t- \displaystyle\frac{a+b}{2} \Big) + 2^{-5/2} \sigma_2^3 (b-a)^2-M-(\log n)^2,
				\end{align*}
			
				\noindent where we have denoted $\sigma_2 = 2^{-7/6} (3\pi)^{2/3} (2A^2 + B + 3)^{-1}$. Using \eqref{anrsigmaprobability} and the facts that for $(j, t) \in \llbracket 1, n-1 \rrbracket \times \{ a, b \}$ we have 
				\begin{flalign*}
					2^{-5/2} \sigma_2^3 (b-a)^2 - M - 2^{-1/2} \sigma_2^3 \Big( t - \displaystyle\frac{a+b}{2} \Big)^2 - \displaystyle\frac{(3 \pi)^{2/3} j^{2/3}}{2^{7/6} \sigma_2} \le -B j^{2/3} - M \le \min \{ u_j, v_j \} 
				\end{flalign*}
				
				\noindent and for $t \in [a, b]$ we have (since $\sigma_2 \le 1$, $b-a \le An^{1/3}$, and $A^2 + B + 3 \ge A^{1/2} + B + 1$)
				\begin{flalign*}
					2^{-5/2}  \sigma_2^3 ( & b-a)^2 -M - 2^{-1/2} \sigma_2^3 \Big( t-\displaystyle\frac{a+b}{2} \Big)^2 - \displaystyle\frac{(3 \pi)^{2/3} n^{2/3}}{2^{7/6} \sigma_2} \\
					& \le  (b-a)^2 - M - (2A^2 + B + 3) n^{2/3} \le -(A^{1/2} + B + 1) n^{2/3} - M,
				\end{flalign*}
			
				\noindent the proof of \eqref{e:xjlowbound} closely follows that of \eqref{e:xjupbound}, so further details are omitted.
			\end{proof}

	\section{Proofs of Results From Chapter \ref{GAPSCALE}}
	
	\label{Proof0}

	\subsection{Convergence of the Alternating Dynamics} 
	
	\label{DynamicConverge}

	Let $\Omega$ be a measurable space with $\sigma$-algebra $\mathcal{F}$; let $\mathscr{P} (\Omega)$ denote the space of probability measures on $(\Omega, \mathcal{F})$.

	\begin{assumption} 
		
		\label{vx} 
		
		Adopting the above notation, let $\mathsf{K} : \Omega \times \mathcal{F} \rightarrow \mathbb{R}_{\ge 0}$ be a Markov transition kernel. For any function $\varphi : \Omega \rightarrow \mathbb{R}_{\ge 0}$ and measure $\mu$ on $\Omega$, define the function $\mathfrak{\mathsf{K}} \varphi : \Omega \rightarrow \mathbb{R}_{\ge 0}$ and measure $\mathsf{K} \mu$ on $\Omega$ by setting
		\begin{flalign*} 
			\mathsf{K} \varphi (x) 	= \displaystyle\int_{\Omega} \varphi (y) \mathsf{K} (x, dy); \qquad \mathsf{K} \mu (A) = \displaystyle\int_{\Omega} \mathsf{K} (x, A) \mu (dx),
		\end{flalign*} 
		
		\noindent for any $x \in \Omega$ and measurable set $A \in \mathcal{F}$. Assume that there exist constants $\alpha \in (0, 1)$, $\gamma \in (0, 1)$, $B \ge 0$, and $R > \frac{2B}{1 - \gamma}$; a potential function $V : \Omega \rightarrow \mathbb{R}_{\ge 0}$; and a probability measure $\nu$ on $\Omega$, such that the following two conditions hold.
		
		\begin{enumerate} 
			\item For each $x \in \Omega$, we have $\mathsf{K}V (x) \le \gamma V (x) + B$, for each $x \in \Omega$.
			\item For each $x \in \Omega$ with $V(x) \le R$, and any measurable set $A \in \mathcal{F}$, we have $\mathsf{K} (x, A) \ge \alpha \nu (A)$. 
		\end{enumerate} 
		
	\end{assumption} 
	
	The following result provides a convergence theorem for Harris chains \cite{ESMP}. It appears in \cite{CSS}, though it is stated as written below in \cite{ET}.

	\begin{lem}[{\cite[Theorem 1.2]{ET}}]
		
		\label{convergeprocess}
		
		Adopt \Cref{vx}, and fix some measure $\mu$ on $\Omega$. Then, the Markov process defined by $\mathsf{K}$ has a unique stationary measure $\mu_0$, and $\lim_{m \rightarrow \infty} \| \mathsf{K}^m \mu - \mu_0 \|_{\TV} = 0$. 
		
	\end{lem} 
	
	We next apply \Cref{convergeprocess} to the alternating dynamics of \Cref{dynamicalternating}. Throughout the remainder of this section, we adopt the notation of that definition. These include the function $f : \llbracket 0, T \rrbracket \rightarrow \overline{\mathbb{R}}$; the family $\bm{\mathsf{y}}$ of $n$ non-intersecting $(T+1)$-step walks $\bm{\mathsf{y}} (t) = \big( \mathsf{y}_1 (t), \mathsf{y}_2 (t), \ldots , \mathsf{y}_n (t) \big) \in \overline{\mathbb{W}}_n$ (over $t \in \llbracket 0, T \rrbracket$), and the associated Markov operator $\mathsf{P}$ (which we also interpret as a kernel) for the alternating dynamics. We also set $\bm{u} = \bm{\mathsf{y}}(0)$ and $\bm{v} = \bm{\mathsf{y}}(T)$, which are fixed throughout the dynamics. Then the state space for the alternating dynamics $\mathsf{P}$ can be viewed as 
	\begin{flalign*}
		\Omega_0 = \Big\{ \big( \bm{\mathsf{y}} (t) \big)_{t \in \llbracket 1, T - 1 \rrbracket} \in \mathbb{W}_n^{T-1} : \displaystyle\min_{t \in \llbracket 1, T-1 \rrbracket} \big( \mathsf{y}_n (t) - f(t) \big) \ge 0 \Big\}.
	\end{flalign*} 
	
	\noindent We define the associated potential function $V_0$ to be
	\begin{flalign} 
		\label{v0} 
		V_0 (\bm{\mathsf{y}}) = \displaystyle\max_{j \in \llbracket 1, n \rrbracket } \displaystyle\max_{t \in \llbracket 1, T-1 \rrbracket } \Big( \big| \mathsf{y}_j (t) \big| + 1 \Big).
	\end{flalign} 
	
	\noindent We then have the following two lemmas verifying \Cref{vx} for the alternating dynamics; the former is proven in \Cref{Proofuvy} below, and the proof of the latter is similar to that of \cite[Lemma B.13]{U}.
	
	\begin{lem}
		
		\label{vy1} 
		
		There exist constants $\gamma = \gamma (f, \bm{u}, \bm{v}) \in (0, 1)$ and $B = B(f, \bm{u}, \bm{v}) \ge 0$ such that, for any family $\bm{\mathsf{y}}$ of $n$ non-intersecting $(T+1)$-step walks with $\bm{\mathsf{y}} (0) = \bm{u}$ and $\bm{\mathsf{y}} (T) = \bm{v}$, we have $\mathsf{P}^2 V_0 (\bm{\mathsf{y}}) \le \gamma V_0 (\bm{\mathsf{y}}) + B$.
		
	\end{lem}

	\begin{lem}
		
		\label{vy2} 
		
		For any real number $R > 1$, there exists a constant $\alpha = \alpha (f, \bm{u}, \bm{v}, R) > 0$ such that the following holds. Letting $\nu_0$ denote the Lebesgue measure on the set 
		\begin{flalign*}
			\Omega_1 = \bigg\{ \bm{\mathsf{x}} \in \Omega_0 : V_0 (\bm{\mathsf{x}}) \le   \max_{t \in \llbracket 1, T-1 \rrbracket} \max \big\{ f(t), 0 \big\}  + R + 1 \bigg\},
		\end{flalign*} 
		
		\noindent we have  $\mathsf{P}^2 (\bm{\mathsf{y}}, A) \ge \alpha \nu_0 (A)$, for each $\bm{\mathsf{y}} \in \Omega_1$ and any measurable subset $A \subseteq \mathbb{W}_n^{T-1}$. 
		
	\end{lem}
	
	\begin{proof}

		For any integer $T' \ge 2$; two $n$-tuples $\bm{u}', \bm{v}' \in \mathbb{W}_n$; and function $f' : \llbracket 0, T' \rrbracket \rightarrow \overline{\mathbb{R}}$, the density of the measure $\mathsf{G}_{f'}^{\bm{u}'; \bm{v}'}$ on sequences $\bm{\mathsf{x}}(t) = \big( \mathsf{x}_1 (t), \mathsf{x}_2 (t), \ldots , \mathsf{x}_n (t) \big)$ is given by
		\begin{flalign}
			\label{g0} 
			C \cdot \textbf{1}_{\bm{\mathsf{x}} \in \Omega_0} \cdot \displaystyle\prod_{j = 1}^n \Bigg( \textbf{1}_{\mathsf{x}_j (0) = u_j'} \textbf{1}_{\mathsf{x}_j (T') = v_j'} \displaystyle\prod_{t=1}^{T'} \exp \bigg( -\displaystyle\frac{1}{2} \big( \mathsf{x}_j (t) - \mathsf{x}_j (t-1) \big)^2 \bigg) \displaystyle\prod_{t = 1}^{T'-1} d \mathsf{x}_j (t) \Bigg),
		\end{flalign}
		
		\noindent for some normalization constant $C = C (f', \bm{u}', \bm{v}') > 0$. Observe that there exist some constant $c_1 = c_1 (f', \bm{u}', \bm{v}') > 0$ such that $C > c_1$, since the interior of $\Omega_0$ is nonempty. Further observe that, for any fixed real number $R_0 \ge \max_{t \in \llbracket 1, T-1 \rrbracket} \max \big\{ f(t), 0 \big\} + 1$, when restricting to the set of $\bm{\mathsf{x}} \in \Omega_0$ such that $V_0 (\bm{\mathsf{x}}) \le R_0$, the density \eqref{g0} is uniformly bounded above and below (in a way dependent on $R_0$). Thus, there exists a constant $c_1 = c_1 (f'; \bm{u}'; \bm{v}', R_0) > 0$ such that, on $\big\{ \bm{\mathsf{x}} \in \Omega_0 : V_0 (\bm{\mathsf{x}}) \le R_0 \big\}$, the measure $\mathsf{G}_{f'}^{\bm{u}'; \bm{v}'}$ is absolutely continuous with respect to the Lebesgue measure on this set, and its Radon--Nikodym derivative is bounded above by $c_1^{-1}$ and below by $c_1$. 
		
		We use this twice, with $(T'; \bm{u}'; \bm{v}') = \big( 2; \bm{u}; \bm{\mathsf{y}}(2) \big)$ and then with $(T'; \bm{u}'; \bm{v}') = \big( T-1; \bm{\mathsf{y}}(1); \bm{v} \big)$; by \Cref{dynamicalternating}, the former corresponds to the first application of $\mathsf{P}$ and the latter to the second application of $\mathsf{P}$. The former yields a constant $c_1 = c_1 (f, \bm{u}, \bm{v}, R) > 0$ such that 
		\begin{flalign}
			\label{y2} 
			\mathbb{P} \Bigg[ \bigcap_{j = 1}^n \big\{ \mathsf{P} \mathsf{y}_j (1) \in S_j \big\} \Bigg] \ge c_1 \displaystyle\prod_{j=1}^n \displaystyle\int_{S_j} dy,
		\end{flalign} 
		
		\noindent for any measurable subsets 
		\begin{flalign}
			\label{siomega0} 
			S_1, S_2, \ldots , S_n  \subseteq \bigg\{ x \in \mathbb{R} : x \ge f(1), |x| \le  \max_{t \in \llbracket 1, n \rrbracket} \max \big\{ f(t), 0 \big\} + R + 1 \bigg\}.
		\end{flalign}
		
		\noindent The second application yields a constant $c_2 = c_2 (f, \bm{u}, \bm{v}, R) > 0$ such that
		\begin{flalign}
			\label{y3}
			\mathbb{P} \Bigg[ \bigcap_{j = 1}^n \bigcap_{t = 2}^T \big\{ \mathsf{P}^2 \mathsf{y}_j (t) \in S_{t, j} \big\} \Bigg] \ge c_2 \displaystyle\prod_{t=2}^{T-1} \displaystyle\prod_{j = 1}^n \displaystyle\int_{S_{t, j}} dy,
		\end{flalign}
		
		\noindent for any measurable subsets $S_{t,j}$ of the right side of \eqref{siomega0}. The lemma then follows from combining \eqref{y2} and \eqref{y3}.
	\end{proof}
	
	Given these lemmas, we can quicky establish \Cref{pkxtconverge}.
	
	\begin{proof}[Proof of \Cref{pkxtconverge}]
		
		The fact that $\mathsf{G}_f^{\bm{u}; \bm{v}}$ is stationary for $\mathsf{P}^2$ follows from \Cref{gmeasurep}. Thus, the lemma follows from \Cref{convergeprocess}, using \Cref{vy1} and \Cref{vy2} (normalizing $\nu_0$ in the latter so that it becomes a probability measure) to verify \Cref{vx}. 
	\end{proof}

	\subsection{Proof of \Cref{vy1}}
	
	\label{Proofuvy}
	
	In this section we establish \Cref{vy1}. To this end, we first require the following tail bound for non-intersecting Gaussian bridges.
	
	\begin{lem} 
		
		\label{uvf2}
		
		For any integer $T' \ge 2$; two $n'$-tuples $\bm{u}', \bm{v}' \in \overline{\mathbb{W}}_n$; and function $f' : \llbracket 0, T' \rrbracket \rightarrow \mathbb{R} \cup \{ - \infty \}$, there exists a constant $c = c(f', n) > 0$ such that the following holds for any real number $r \ge 0$. Sampling non-intersecting Gaussian bridges $\bm{\mathsf{x}}' (t) = \big(\mathsf{x}_1' (t), \mathsf{x}_2' (t), \ldots , \mathsf{x}_n' (t) \big)$ from the measure $\mathsf{G}_{f'}^{\bm{u}'; \bm{v}'}$, we have
		\begin{flalign*}
			\mathbb{P} \Bigg[ \bigcup_{t=1}^{T-1} \bigcup_{j=1}^n \bigg\{ \big| \mathsf{x}_j' (t) \big| \ge \displaystyle\frac{T'-t}{T'} \cdot \displaystyle\max_{j \in \llbracket 1, n \rrbracket} |u_j'| + \displaystyle\frac{t}{T'} \cdot \displaystyle\max_{j \in \llbracket 1, n \rrbracket} |v_j'| + r \bigg\} \Bigg] \le c^{-1} e^{-c r^2}.
		\end{flalign*}
	\end{lem} 
	
	\begin{proof}
		
		First observe that there exists a constant $c_1 = c_1 (f') > 1$ such that the following holds for any real number $r > 0$. Given a $(T'+1)$-step Gaussian bridge $\big( \mathsf{x} (0), \mathsf{x} (1), \ldots , \mathsf{x} (T') \big)$ conditioned to start and end at some points $u' \in \mathbb{R}$ and $v' \in \mathbb{R}$, respectively, and satisfy $\mathsf{x} (t) \ge f(t)$, we have 
		\begin{flalign}
			\label{xt0} 
			\mathbb{P} \Bigg[ \bigcup_{t \in \llbracket 0, T' \rrbracket} \bigg\{ \big| \mathsf{x} (t) \big| \ge  \displaystyle\frac{T'-t}{T'} \cdot |u'| + \displaystyle\frac{t}{T'} \cdot |v'| + r \bigg\} \Bigg] \le c_1^{-1} e^{-c_1 r^2}. 
		\end{flalign}
		
		 Now, let $\mathfrak{G}_{f'}^{\bm{u}'; \bm{v}'}$ denote the law on sequences $\bm{\mathsf{x}}(t) = \big( \mathsf{x}_1 (t), \mathsf{x}_2 (t), \ldots , \mathsf{x}_n (t) \big)$, with $t \in \llbracket 0, T' \rrbracket$, of $n$ independent Gaussian bridges starting at $\bm{u}'$, ending at $\bm{v}'$, and conditioned to remain above $f'$; as it does not impose the non-intersecting condition, we may view it as the law of ``free'' Gaussian bridges. Then, from \eqref{xt0} and a union bound, we deduce that there exists a constant $c_2 = c_2 (f', n) > 0$ such that 
		\begin{flalign}
			\label{xtt0} 
			\mathbb{P} \Bigg[ \bigcup_{t=1}^{T-1} \bigcup_{j=1}^n \bigg\{ \big| x_j (t) \big| \ge \displaystyle\frac{T'-t}{T'} \cdot \displaystyle\max_{j \in \llbracket 1, n \rrbracket} |u_j'| + \displaystyle\frac{t}{T'} \cdot \displaystyle\max_{j \in \llbracket 1, n \rrbracket} |v_j'| + r \bigg\} \Bigg] \le c_2^{-1} e^{-c_2 r^2}.
		\end{flalign}

		Next, observe that there exists a constant $c_3 = c_3 (f', \bm{u}', \bm{v}') > 0$ such that the walks in $\bm{\mathsf{x}}$ sampled under $\mathfrak{G}_{f'}^{\bm{u}'; \bm{v}'}$ do not intersect (that is, $\bm{\mathsf{x}}(t) \in \mathbb{W}_n$ for each $t \in \llbracket 1, T-1 \rrbracket$), with probability at least $c_3$ under $\mathfrak{G}_{f'}^{\bm{u}'; \bm{v}'}$ (as this is an open condition). Hence, the Radon--Nikodym derivative of the non-intersecting measure $\mathsf{G}_{f'}^{\bm{u}'; \bm{v}'}$ with respect to the free one $\mathfrak{G}_{f'}^{\bm{u}'; \bm{v}'}$ is bounded above by $c_3^{-1}$. Together with \eqref{xtt0}, this establishes the lemma. 
		\end{proof} 
	
		Now we can esetablish \Cref{vy1} 
		
		\begin{proof}[Proof of \Cref{vy1}]
		
		This lemma will follow from two applications of \Cref{uvf2}. First taking the $(T'; \bm{\mathsf{x}}'; \bm{u}', \bm{v}')$ there to be $\big(2; \bm{\mathsf{y}} |_{\llbracket 0, 2 \rrbracket}; \bm{u}; \bm{\mathsf{y}}(2) \big)$ here yields (as $\max_{j \in \llbracket 1, n \rrbracket} \big| \mathsf{y}_j (2) \big| \le V_0 (\bm{\mathsf{y}})$) a constant $c_1 = c_1 (f, \bm{u}) > 0$ such that
		\begin{flalign}
			\label{py} 
			\mathbb{P} \Bigg[ \bigcup_{j = 1}^n \bigg\{ \big| \mathsf{P} \mathsf{y}_j (1) \big| \ge \displaystyle\frac{V_0 (\bm{\mathsf{y}})}{2} + r \bigg\} \Bigg] \le c_1^{-1} e^{-c_1 r^2},
		\end{flalign} 
		
		\noindent for any real number $r \ge 0$. Next applying \Cref{uvf2} with the $(T'; \bm{\mathsf{x}}'; \bm{u}', \bm{v}')$ there to be $\big( T-1; \bm{\mathsf{y}} |_{\llbracket 1, T \rrbracket}; \mathsf{P} \mathsf{y}_j (1); \bm{v} \big)$ here, yields a constant $c_2 = c_2 (f, \bm{v}) > 0$ such that 
		\begin{flalign*} 
			\mathbb{P} \Bigg[ \bigcup_{t = 2}^{T-1} \bigcup_{j=1}^n \bigg\{ \big| \mathsf{y}_j (t) \big| \ge \displaystyle\frac{1}{2} \cdot \displaystyle\max_{j \in \llbracket 1, n \rrbracket} \big| \mathsf{P} \mathsf{y}_j (1) \big| + r \bigg\} \Bigg] \le c_2^{-1} e^{-c_2 r^2},
		\end{flalign*} 
	
		\noindent for any real number $r \ge 0$. This, together with \eqref{py}, the definition \eqref{v0} of $V_0$, and a union bound, yields a constant $c_3 = c_3 (f, \bm{u}, \bm{v}) > 0$ such that 
		\begin{flalign*} 
			\mathbb{P} \Big[ V_0 (\mathsf{P}^2 \bm{\mathsf{y}}) \le \displaystyle\frac{1}{4} \cdot V_0 (\bm{\mathsf{y}}) + \displaystyle\frac{3r}{2} \Big] \le c_3^{-1} e^{-c_3 r^2}.
		\end{flalign*} 
	
		\noindent Integrating this bound then establishes the lemma. 
	\end{proof}

\section{Proofs of Results From Chapter \ref{STATISTICSBRIDGES}}
	
\subsection{Proofs of Continuum Monotonicity Results} 

In this section we first establish \Cref{limitheightcompare} and then outline the proof of \Cref{limitdifferencecompare}.

\label{ProofContinuousCompare}

\begin{proof}[Proof of \Cref{limitheightcompare}]
	
	We only establish the second part of the lemma, as the proof of the first is entirely analogous. For each integer $n \ge 1$, define the $\lfloor An \rfloor$-tuples and $\lfloor \widetilde{A} n \rfloor$-tuples\footnote{For notational simplicity, we will omit the floors in what follows, assuming that $An = \lfloor An \rfloor$, $\widetilde{A} n = \lfloor \widetilde{A} n\rfloor$, and $wn = \lfloor wn \rfloor$; this will barely affect the analysis.} $\bm{u}, \bm{v} \in \overline{\mathbb{W}}_{An}$ and $\widetilde{\bm{u}}, \widetilde{\bm{v}} \in \overline{\mathbb{W}}_{\tilde{A}n}$ by for each $j$ setting 
	\begin{flalign}
		\label{ugvgugvg} 
		& u_j = G^{\star} (a, n^{-1} j); \qquad v_j = G^{\star} (b, n^{-1} j); \qquad \widetilde{u}_j = \widetilde{G}^{\star} (a, n^{-1} j); \qquad \widetilde{v}_j = \widetilde{G}^{\star} (b, n^{-1} j). 
	\end{flalign}
	
	\noindent Sample the two families of non-intersecting Brownian bridges $\bm{x}^n = (x_1, x_2, \ldots , x_{An}) \in \llbracket 1, An \rrbracket \times \mathcal{C} \big( [a, b] \big)$ and $\widetilde{\bm{x}}^n = (\widetilde{x}_1, \widetilde{x}_2, \ldots , \widetilde{x}_{\tilde{A} n}) \in \llbracket 1, \widetilde{A} n \rrbracket \times \mathcal{C} \big( [a, b] \big)$ according to the measures $\mathsf{Q}^{\bm{u}; \bm{v}} (n^{-1})$ and $\mathsf{Q}^{\tilde{\bm{u}}; \tilde{\bm{v}}} (n^{-1})$, respectively. Define the functions $f : [a, b] \rightarrow \mathbb{R}$ and $\widetilde{f}: [a, b] \rightarrow \mathbb{R}$ by setting 
	\begin{flalign}
		\label{ftg} 
		& f(t) = G^{\star} (t, w);   \qquad \widetilde{f} (t) = \widetilde{G}^{\star} (t, w),   \quad \qquad \text{for each $t \in [a, b]$}.
	\end{flalign} 
	
	Further fix a real number $\varepsilon > 0$, and let $\mathscr{E}_n$ denote the event on which $\big| f(t) - x_{wn} (t) \big| \le \varepsilon$ and $\big| \widetilde{f} (t) - \widetilde{x}_{wn} (t) \big| \le \varepsilon$ for all $t \in (a, b)$. We have that $\lim_{n \rightarrow \infty} \mathbb{P}[\mathscr{E}_n] = 1$ by the second part of \Cref{convergepathomega}. Now, recalling \Cref{property}, condition on $\mathcal{F}_{\ext}^{\bm{x}} \big( \llbracket 1, wn-1 \rrbracket \times (a, b) \big)$ and $\mathcal{F}_{\ext}^{\tilde{\bm{x}}} \big( \llbracket 1, wn-1 \rrbracket \times (a, b) \big)$, and restrict to $\mathscr{E}_n$. Define the $(wn-1)$-tuples $\bm{u'}, \bm{v'}, \widetilde{\bm{u}}', \widetilde{\bm{v}}' \in \overline{\mathbb{W}}_{wn-1}$ to be the restriction of $\bm{u}, \bm{v}, \widetilde{\bm{u}}, \widetilde{\bm{v}}$ on $\llbracket 1, wn-1\rrbracket$.  Sample non-intersecting Brownian bridges $\bm{y}^n = (y_1, y_2, \ldots , y_{wn-1}) \in \llbracket 1, wn-1 \rrbracket \times \mathcal{C} \big( [a, b] \big)$ and $\widetilde{\bm{y}}^n = (\widetilde{y}_1, \widetilde{y}_2, \ldots , \widetilde{y}_{wn-1}) \in \llbracket 1, wn-1 \rrbracket \times \mathcal{C} \big( [a, b] \big)$ from $\mathsf{Q}_f^{\bm{u}'; \bm{v}'} (n^{-1})$ and $\mathsf{Q}_{\tilde{f}}^{\tilde{\bm{u}}'; \tilde{\bm{v}}'} (n^{-1})$, respectively.

	By our restriction to $\mathscr{E}_n$, with the Brownian Gibbs property for $\bm{x}$ and $\widetilde{\bm{x}}$ and (the $B = \varepsilon$ case) of \Cref{uvv}, we obtain a coupling between $(\bm{x}, \bm{y})$ and $(\widetilde{\bm{x}}, \widetilde{\bm{y}})$ such that, on $\mathscr{E}$, we have 
	\begin{flalign}
		\label{xyxy2} 
		\displaystyle\max_{j \in \llbracket 1, wn-1 \rrbracket} \displaystyle\sup_{t \in [a, b]} \big| y_j (t) - x_j (t) \big| \le \varepsilon; \qquad \displaystyle\max_{j \in \llbracket 1, wn-1 \rrbracket} \displaystyle\sup_{t \in [a, b]} \big| \widetilde{y}_j (t) - \widetilde{x}_j (t) \big| \le \varepsilon.
	\end{flalign}
	
	\noindent Moreover, since $\bm{u} \le \widetilde{\bm{u}}$ and $\bm{v} \le \widetilde{\bm{v}}$ (as $G^{\star} (t, y) \le \widetilde{G}^{\star} (t, y)$ for $(t, y) \in \{ a, b \} \times [0, w]$), and $f \le \widetilde{f}$ (by \eqref{ftg}, as $G^{\star} (t, w) \le \widetilde{G}^{\star} (t, w)$), it follows from \Cref{monotoneheight} that we may couple $\bm{y}$ and $\widetilde{\bm{y}}$ in such a way that $y_j (t) \le \widetilde{y}_j (t)$ for each $(j, t) \in \llbracket 1, wn-1 \rrbracket \times [a, b]$. By the gluing lemma \cite[Theorem 1.1.10]{UCT}, we may exhibit this coupling on the same probability space as those between $(\bm{x}, \bm{y})$ and $(\widetilde{\bm{x}}, \widetilde{\bm{y}})$ satisfying \eqref{xyxy2}. This induces a coupling between $\bm{x}$ and $\widetilde{\bm{x}}$ such that, on $\mathscr{E}$, we have $x_j (t) - 2 \varepsilon \le \widetilde{x}_j (t)$ for all $(j, t) \in \llbracket 1, wn-1 \rrbracket \times (a, b)$. Recalling that $\lim_{n \rightarrow \infty} \mathbb{P} [\mathscr{E}_n] = 1$, it follows that
	\begin{flalign*} 
		\displaystyle\lim_{n \rightarrow \infty} \mathbb{P} \Bigg[ \bigcap_{j \in \llbracket 1, wn-1 \rrbracket} \bigcap_{t \in [a, b]} \big\{ x_j (t) - 2 \varepsilon  \le \widetilde{x}_j (t) \big\} \Bigg] = 1. 
	\end{flalign*} 
	
	\noindent Together with the first statement of \Cref{convergepathomega}, this implies upon letting $\varepsilon$ tend to $0$ that  
	\begin{flalign} 
		\label{gtygty}
		G^{\star} (t, y) \le \widetilde{G}^{\star} (t, y), \qquad \text{whenever $G^{\star}(t,y')$ and $\widetilde{G}^{\star}(t,y')$ are continuous at $y'=y$}.
	\end{flalign}
	
	\noindent It thus remains to show that \eqref{gtygty} continues for more general $(t, y) \in [a, b] \times [0, w]$.

	To this end, for $t \in (a, b)$ observe (since the densities $\varrho^{\star}$ and $\widetilde{\varrho}^{\star}$ associated with $\bm{\mu}^{\star}$ and $\widetilde{\bm{\mu}}^{\star}$, respectively, are bounded by the third part of \Cref{mutrhot}) that \eqref{htxintegral} and \eqref{gty} together yield
	\begin{flalign}
		\label{gcontinuous}
		G^{\star} (t, y) \quad \text{and} \quad \widetilde{G}^{\star} (t, y) \qquad \text{are lower semicontinuous and non-increasing in $y \in (0, A_0]$}.
	\end{flalign}
	
	\noindent Thus, given any point $y \in [0, A_0]$ and real number $\delta > 0$, there is a point $y_1 \in (y - \delta, y)$ such that $G^{\star} (t, y')$ and $\widetilde{G}^{\star} (t, y')$ are continuous in its second variable at $y' = y_1$. Hence, letting $\delta$ tend to $0$, it follows from \eqref{gtygty} that $G^{\star} (t,y^-) \le \widetilde{G}^{\star} (t, y^-)$, and so by \eqref{gcontinuous} we deduce that $G^{\star} (t, y) \le \widetilde{G}^{\star} (t, y)$ for each $(t, y) \in [a, b] \times [0, w]$, establishing the lemma. 
\end{proof}

\begin{proof}[Proof of \Cref{limitdifferencecompare} (Outline)]
	
	We again only establish the second part of the lemma, as the proof of the first is entirely analogous. Its proof will be similar to that of the second part of \Cref{limitheightcompare}, and so we only outline it. In what follows, we adopt the notation of that lemma and its proof, recalling the entrance and exit data $\bm{u}, \bm{v} \in \overline{\mathbb{W}}_{An}$ and $\widetilde{\bm{u}}, \widetilde{\bm{v}} \in \overline{\mathbb{W}}_{\tilde{A} n}$ from \eqref{ugvgugvg}; the boundaries $f, \widetilde{f} : [a, b] \rightarrow \mathbb{R}$ from \eqref{ftg}; and the event $\mathscr{E}_n$. We also again condition on $\mathcal{F}_{\ext}^{\bm{x}} \big( \llbracket 1, wn-1 \rrbracket \times (a, b) \big)$ and $\mathcal{F}_{\ext}^{\bm{\tilde{x}}} \big( \llbracket 1, wn-1 \rrbracket \times (a, b) \big)$ and restrict to $\mathscr{E}_n$. Additionally define the non-intersecting Brownian bridges $\bm{x}^n \in \llbracket 1, An \rrbracket \times \mathcal{C} \big( [a, b] \big)$, $\widetilde{\bm{x}}^n \in \llbracket 1, \widetilde{A} n \rrbracket \times \mathcal{C} \big( [a, b] \big)$,  and $\bm{y}^n, \widetilde{\bm{y}}^n \in \llbracket 1, wn-1 \rrbracket \times \mathcal{C} \big( [a, b] \big)$ as in the proof of \Cref{limitheightcompare}. Following the proof of \Cref{limitheightcompare}, we may couple $(\bm{x}, \bm{y})$ and $(\widetilde{\bm{x}}, \widetilde{\bm{y}})$ such that, on $\mathscr{E}_n$, \eqref{xyxy2} holds.
	
	Next, observe that $u_j - u_{j+1} \le \widetilde{u}_j - \widetilde{u}_{j+1}$ and $v_j - v_{j+1} \le \widetilde{v}_j - \widetilde{v}_{j+1}$, since $\big| G^{\star} (t, y) - G^{\star} (t, y') \big| \le \big| \widetilde{G}^{\star} (t, y) - \widetilde{G}^{\star} (t, y') \big|$ for each $(t, y), (t, y') \in \{ a,  b \} \in [0, w]$. Moreover, by \eqref{ftg} and \eqref{gdifference121}, we have $r \cdot f (s) - f \big( rs + (1-r) t \big) + (1-r) \cdot f(t) \le r \cdot \widetilde{f} (s) - \widetilde{f} \big( rs + (1-r) t \big) + (1-r) \cdot \widetilde{f} (t)$ for any $s, t \in (a, b)$ and $r \in [0, 1]$. Thus, it follows from gap monotonicity \Cref{monotonedifference} that we may couple $\bm{y}$ and $\widetilde{\bm{y}}$ such that $y_j (t) - y_{j+1} (t) \le \widetilde{y}_j (t) - \widetilde{y}_{j+1} (t)$ for each $(j, t) \in \llbracket 1, wn-1 \rrbracket \times [a, b]$. Combining this with \eqref{xyxy2} and using the fact that $\lim_{n \rightarrow \infty} \mathbb{P} [\mathscr{E}_n] = 1$, it follows that 
	\begin{flalign*}
		\displaystyle\lim_{n \rightarrow \infty} \mathbb{P} \Bigg[ \bigcap_{j, j' \in \llbracket 1, wn-1 \rrbracket} \bigcap_{t \in [a, b]} \Big\{ \big| x_j (t) - x_{j'} (t) - 4 \varepsilon \le \big| \widetilde{x}_j (t) - \widetilde{x}_{j'} (t) \big| \Big\} \Bigg] = 1. 
	\end{flalign*}
	
	\noindent Together with the first statement of \eqref{convergepathomega}, this implies upon letting $\varepsilon$ tend to $0$ that 
	\begin{flalign*} 
		\big| G^{\star} (t, y) - G^{\star} (t, y') \big| \le \big| \widetilde{G}^{\star} (t,y) - \widetilde{G}^{\star} (t, y') \big|,
	\end{flalign*} 
	if $G^{\star}(t,y'')$ and $\widetilde{G}^{\star}(t,y'')$ are continuous at $y''=y$ and $y''=y'$.  Extending this bound to all $(t, y), (t, y') \in (a, b) \times [0, w]$ then follows as in the proof of \Cref{limitheightcompare}, and so further details are omitted; this establishes the lemma.
\end{proof}

\begin{proof}[Proof of \Cref{airyheightcompare} (Outline)]
	
	The proof of this lemma is similar that of \Cref{limitheightcompare}, so we only outline the proof of the first statement in its second part. Let $\bm{\mathcal{S}} = (\mathcal{S}_1, \mathcal{S}_2, \ldots )$ denote a scaled parabolic Airy line ensemble (recall \Cref{ensemblewalks}), and set $\widetilde{x}_j (t) = n^{-2/3} \cdot \mathcal{S}_j (tn^{1/3})$, for each $(j, t) \in \mathbb{Z}_{\ge 1} \times \mathbb{R}$.  Fix a real number $\varepsilon > 0$. For each integer $n \ge 1$, define the $ An$-tuples (again omit the floors in what follows) $\bm{u}, \bm{v} \in \overline{\mathbb{W}}_{An}$ and the $(wn-1)$-tuples $\widetilde{\bm{u}}, \widetilde{\bm{v}} \in \overline{\mathbb{W}}_{wn-1}$ by setting 
	\begin{flalign}
		\label{guvg2} 
		& u_j = G^{\star} (a, n^{-1} j) - \varepsilon; \qquad v_j = G^{\star} (b, n^{-1} j) - \varepsilon; \qquad \widetilde{u}_j = \widetilde{x}_j (a); \qquad \widetilde{v}_j = \widetilde{x}_j (b),
	\end{flalign}
	
	\noindent for each $j$. Also define the functions $f : [a, b] \rightarrow \mathbb{R}$ and $\widetilde{f}: [a, b] \rightarrow \mathbb{R}$ by setting 
	\begin{flalign}
		\label{ftg2} 
		& f(t) = G^{\star} (t, w) - \varepsilon;  \qquad \widetilde{f}(t) = \widetilde{x}_{wn} (t); \qquad \quad \text{for each $t \in [a, b]$}.
	\end{flalign} 

	 Sample non-intersecting Brownian bridges $\bm{x}^n = (x_1, x_2, \ldots , x_{An}) \in \llbracket 1, An \rrbracket \times \mathcal{C} \big( [a, b] \big)$ under the measure $\mathsf{Q}^{\bm{u}; \bm{v}} (n^{-1})$. Fixing a real number $\varepsilon > 0$, define the events
	\begin{flalign*} 
		\mathscr{E}_{n;1} & = \bigcap_{t \in [a, b]} \big| f(t) - x_{wn} (t) \big| \le \varepsilon; \\
		\mathscr{E}_{n;2} & = \bigcap_{t \in [a, b]} \Big\{ \big| \widetilde{f} (t) - \widetilde{G}^{\star} (t, w) \big| \le \varepsilon \Big\}  \cap \bigcap_{j=1}^{wn-1} \big\{ |\widetilde{u}_j - \widetilde{G}^{\star} (a, jn^{-1})| \le \varepsilon \big\} \cap \big\{ |\widetilde{v}_j - \widetilde{G}^{\star} (b, jn^{-1})| \le \varepsilon \big\}.
	\end{flalign*} 
	
	\noindent By the second part of \Cref{convergepathomega}, we have that $\lim_{n \rightarrow \infty} \mathbb{P} [\mathscr{E}_{n;1}] = 1$ and, by \Cref{kdeltad}, we have that $\lim_{n \rightarrow \infty} \mathbb{P}[\mathscr{E}_{n;2}] = 1$. Recalling \Cref{property}, we condition on $\mathcal{F}_{\ext}^{\bm{x}} \big( \llbracket 1, wn-1 \rrbracket \times (a, b) \big)$ and $\mathcal{F}_{\ext}^{\tilde{\bm{x}}} \big( \llbracket 1, wn-1 \rrbracket \times (a, b) \big)$; we further restrict to $\mathscr{E}_{n;1}$ and $\mathscr{E}_{n;2}$.  
	
	Define $\bm{u}', \bm{v}' \in \overline{\mathbb{W}}_{wn-1}$ to be the restrictions of $\bm{u}$ and $\bm{v}$ on $\llbracket 1, wn-1 \rrbracket$. Sample $\bm{y}^n = (y_1, y_2, \ldots , y_{wn-1}) \in \llbracket 1, wn - 1 \rrbracket \times \mathcal{C} \big( [a, b] \big)$ and $\widetilde{\bm{y}}^n = (\widetilde{y}_1, \widetilde{y}_2, \ldots , \widetilde{y}_{wn-1}) \in \mathcal{C} \big( \llbracket 1, wn-1 \rrbracket \times (a, b) \big)$ according to the measures $\mathsf{Q}_f^{\bm{u}'; \bm{v}'} (n^{-1})$ and $\mathsf{Q}_{\tilde{f}}^{\tilde{\bm{u}}; \tilde{\bm{v}}} (n^{-1})$, respectively. By (the $B = \varepsilon$ case of) \Cref{uvv} (observing that the law of $\bm{\widetilde{x}}$, conditional on $\mathcal{F}_{\ext}^{\bm{\tilde{x}}} \big( \llbracket 1, wn-1 \rrbracket \times (a, b) \big)$, is $\mathsf{Q}_{\tilde{x}_{wn}}^{\tilde{\bm{u}}; \tilde{\bm{v}}}$, by the definition of $\widetilde{\bm{x}}$ in terms of $\bm{\mathcal{S}}$ and the fact that $\bm{\mathcal{S}}$ satisfies the Brownian Gibbs property by \Cref{propertya}), we may on $\mathscr{E}_{n;1}$ and $\mathscr{E}_{n;2}$ couple $(\bm{x}; \bm{y})$ and $(\widetilde{\bm{x}}; \widetilde{\bm{y}})$ so that 
	\begin{flalign}
		\label{xyxy3}  
		\displaystyle\max_{j \in \llbracket 1, wn-1 \rrbracket} \displaystyle\sup_{t \in [a, b]} \big| x_j (t) - y_j (t) \big| \le \varepsilon; \qquad \displaystyle\max_{j \in \llbracket 1, wn-1 \rrbracket} \displaystyle\sup_{t \in [a, b]} \big| \widetilde{x}_j (t) - \widetilde{y}_j (t) \big| \le \varepsilon.
	\end{flalign}
	
	On the events $\mathscr{E}_{n;1}$ and $\mathscr{E}_{n;2}$, also we have $\bm{u}' \le \widetilde{\bm{u}}$, as $u_j' = u_j = G^{\star} (a, jn^{-1}) - \varepsilon \le \widetilde{G}^{\star} (a, jn^{-1}) - \varepsilon \le \widetilde{u}_j$; by similar reasoning, we have $\bm{v}' \le \widetilde{\bm{v}}$. We also have $f(t) = G^{\star} (t,w) - \varepsilon \le \widetilde{G}^{\star} (t, w) - \varepsilon \le \widetilde{f}$. Therefore, by \Cref{monotoneheight}, on the event $\mathscr{E}_{n;1} \cap \mathscr{E}_{n;2}$, we may couple $\bm{y}$ and $\widetilde{\bm{y}}$ so that $y_j (t) \le \widetilde{y}_j (t)$ for each $(j, t) \in \llbracket 1, wn-1 \rrbracket \times [a, b]$. By the gluing lemma \cite[Theorem 1.1.10]{UCT}, we may exhibit this coupling on the same probability space as those between $(\bm{x}, \bm{y})$ and $(\widetilde{\bm{x}}, \widetilde{\bm{y}})$ satisfying \eqref{xyxy3}. This induces a coupling between $\bm{x}$ and $\widetilde{\bm{x}}$ such that, on $\mathscr{E}$, we have $x_j (t) - 2 \varepsilon \le \widetilde{x}_j (t)$ for all $(j, t) \in \llbracket 1, wn-1 \rrbracket \times (a, b)$. Recalling that $\lim_{n \rightarrow \infty} \mathbb{P} [\mathscr{E}_{n;1}] = \lim_{n \rightarrow \infty} \mathbb{P} [\mathscr{E}_{n;2}] = 1$, it follows that
	\begin{flalign*} 
		\displaystyle\lim_{n \rightarrow \infty} \mathbb{P} \Bigg[ \bigcap_{j \in \llbracket 1, wn-1 \rrbracket} \bigcap_{t \in [a, b]} \big\{ x_j (t) - 2 \varepsilon \le \widetilde{x}_j (t) \big\} \Bigg] = 1. 
	\end{flalign*} 
	
	\noindent Together with the first statement of \Cref{convergepathomega} and \Cref{kdeltad} (and recalling the shift by $\varepsilon$ in the definitions of $\bm{u}$ and $f$ from \eqref{guvg2} and \eqref{ftg2}, respectively), this implies that $G^{\star} (t, y) - 3 \varepsilon \le \widetilde{G}^{\star} (t, y)$,  whenever $G^{\star}(t,y')$ and $\widetilde{G}^{\star}(t,y')$ are continuous at $y'=y$. Letting $\varepsilon$ tend to $0$, we deduce that $G^{\star} (t, y) \le \widetilde{G}^{\star} (t, y)$ whenever $G^{\star} (t,y')$ and $\widetilde{G}^{\star} (t,y')$ are continuous at $y'=y$.  It thus remains to show that \eqref{gtygty} continues for more general $(t, y) \in [a, b] \times [0, w]$; the verification of this is entirely analogous to that in the proof of \Cref{limitheightcompare} and is therefore omitted.
\end{proof}

\subsection{Proof of \Cref{fg1g2}}

\label{ProofGFG} 

In this section we establish \Cref{fg1g2}; throughout, we recall the notation of that lemma. By the scale invariance (\Cref{scale21} of \Cref{invariancesscale}) of solutions to \eqref{equationxtd}, we may assume that $\ell = 1$ in what follows.\footnote{Indeed, given $G_1^{\pm}$ defined at $\ell = 1$, we set $G^{\pm} (t, x) = \ell^{-1} \cdot G^{\pm} (\ell t, \ell x)$, which continues to satisfy \eqref{equationxtd} (by \Cref{invariancesscale}) and the statements of \Cref{fg1g2} (possibly with different constants $c > 0$ and $C > 1$).} We first define the open rectangles
\begin{flalign*}
	\breve{\mathfrak{S}} = \Big( 0, \displaystyle\frac{1}{L} \Big) \times \Big( \displaystyle\frac{1}{2L}, 1 - \displaystyle\frac{1}{2L} \Big); \qquad \breve{\mathfrak{S}}' = \Big( \displaystyle\frac{1}{8L}, \displaystyle\frac{7}{8L} \Big) \times \Big( \displaystyle\frac{1}{L}, 1 - \displaystyle\frac{1}{L} \Big).
\end{flalign*}

\noindent We further fix functions\footnote{For example, fix a nonnegative, smooth function $\psi : [0, 1] \rightarrow \mathbb{R}$ with $\psi (x) = 1$ if $x \in [ 1 / 5, 4 / 5 ]$, with $\psi (x) = 0$ if $x \in [ 0, 1 / 6 ] \cup [5 / 6, 1]$, and with $0 \le \psi (x) \le 1$ for each $x \in [0, 1]$. Then we may set $g_i^{\pm} (x) = \psi (x) \cdot g_i (x) + \big( 1 - \psi (x) \big) \cdot \big( f_i (x) \pm \vartheta^{8/9} \big)$ for any indices $i \in \{ 0, 1 \}$ and $\pm \in \{ +, - \}$, and any real number $x \in [0, 1]$.} $g_0^-, g_0^+, g_1^-, g_1^+ : [0, 1] \rightarrow \mathbb{R}$ so that, for some constant $C_1 = C_1 (\varepsilon, B, m) > 1$ we have 
\begin{flalign}
	\label{gi1gi2} 
	\begin{aligned}
		& g_i^- (x) = g_i (x) = g_i^+ (x), \qquad \qquad \qquad \qquad \qquad \qquad \qquad \qquad   \text{for each $x \in \Big[ \displaystyle\frac{1}{5}, \frac{4}{5} \Big]$}; \\
		&  g_i^- (x) \le g_i (x) \le g_i^+ (x),   \qquad \qquad \qquad \qquad \qquad \qquad \qquad \qquad \text{for each $x \in [0, 1]$}; \\
		& g_i^- (x) = f_i (x) - \vartheta^{8/9}, \quad \text{and} \quad g_i^+ (x) = f_i (x) + \vartheta^{8/9}, \qquad \qquad   \text{for each $x \in \Big[ 0, \displaystyle\frac{1}{6} \Big] \cup \Big[ \displaystyle\frac{5}{6}, 1 \Big]$}; \\
		& \| g_i^- - g_- \|_{\mathcal{C}_0} + \| g_i^+ - g_i \|_{\mathcal{C}^0} \le C_1 \vartheta^{8/9}, \qquad \qquad \qquad \qquad \qquad \text{for each $i \in \{ 0, 1 \}$}; \\ 
		& \| g_i^- \|_{\mathcal{C}^m} + \| g_i^+ \|_{\mathcal{C}^m} \le C_1, \qquad \qquad \qquad \qquad \qquad \qquad \qquad \qquad \text{for each $i \in \{ 0, 1 \}$},
	\end{aligned} 
\end{flalign}
\noindent where the last three properties can be guaranteed since $\| f_i - g_i \|_{\mathcal{C}^0} \le \vartheta \le \vartheta^{8/9}$ (with the facts that $\| f_i \|_{\mathcal{C}^m} \le \| F_i \|_{\mathcal{C}^m (\mathfrak{R})} \le B$ and $\| g_i \|_{\mathcal{C}^m} \le B$), for each $i \in \{ 0, 1 \}$.

Then, for $\vartheta$ sufficiently small, \Cref{f1f2b} yields (upon translating by $( -1 / 2L, 1)$; replacing the $L$ there by $L / 2$ here; and scaling by a factor of $2$, using \Cref{invariancesscale}; and taking $r = 1 / 4$) solutions $G^-, G^+ \in \Adm_{\varepsilon/2} (\breve{\mathfrak{S}})$ to \eqref{equationxtd} on $\breve{\mathfrak{S}}$ and a constant $C_2 = C_2 (\varepsilon, B, m) > 1$ such that
\begin{flalign}
	\label{g1g0} 
	\begin{aligned}
		& G^{\pm} (iL^{-1}, x) = g_i^{\pm} (x), \qquad \text{for each $x \in \Big[ \displaystyle\frac{1}{2L}, 1 - \displaystyle\frac{1}{2L} \Big]$}; \\
		& \| G^{\pm} - F \|_{\mathcal{C}^m (\breve{\mathfrak{S}}')} \le C_2 L^m \cdot \big(  \| g_0^{\pm} - f_0 \|_{\mathcal{C}^0} + \| g_1^{\pm} - f_1 \|_{\mathcal{C}^0} \big) \\
		& \| G^{\pm} - F \|_{\mathcal{C}^{m-5} (\breve{\mathfrak{S}})} \le C_2 L^{m-5} \cdot \big( \| g_0^{\pm} - f_0 \|_{\mathcal{C}^0}^{3/m} + \| g_1^{\pm} - f_1 \|_{\mathcal{C}^0}^{3/m} \big),
	\end{aligned} 	
\end{flalign} 

\noindent for any indices $i \in \{ 0, 1 \}$ and $\pm \in \{ +, - \}$. This, together with the fact that $\| g_i^{\pm} - f_i \|_{\mathcal{C}^0} \le 2C_1 \vartheta^{8/9}$ for each $i \in \{ 0, 1 \}$ (by \eqref{gi1gi2} and the fact that $\| f_i - g_i \|_{\mathcal{C}^0} \le \vartheta$ for each $i$), implies that 
\begin{flalign}
	\label{fg1fg2}
	\begin{aligned}
	& \| G^{\pm} - F \|_{\mathcal{C}^m (\breve{\mathfrak{S}}')} \le 4 C_1 C_2 L^m \vartheta^{8/9} \le 4 C_1 C_2 \vartheta^{4/5}; \\
	&  \| G^{\pm} - F \|_{\mathcal{C}^{m-5} (\breve{\mathfrak{S}})} \le 4 C_1 C_2 L^{m-5} \vartheta^{8/3m} \le 4 C_1 C_2 \vartheta^{1/m},
	\end{aligned}
\end{flalign}

\noindent where in the last inequalities we also used the facts that $L \le \vartheta^{-1/2m^2}$ and that $m \ge 7$. 

\begin{lem}
	
	\label{g120}
	
	There exists a constant $c = c(\varepsilon, B, m) > 0$ such that the following holds if $\vartheta < c$.
	\begin{enumerate} 
		\item For each $(t, x) \in [0, L^{-1}] \times [ 1 / 4, 3 / 4]$, we have $\big| G^+ (t, x)  - G^- (t, x) \big| \le c^{-1} e^{-cL^{1/8}}$. 
		\item For each $(t, x) \in [0, L^{-1}] \times \big( [1 / 10, 3 / 20] \cup [17 / 20, 9 / 10] \big)$, we have 
		\begin{flalign}
			\label{gtxftx1} 
			\big| G^- (t, x) - F(t, x) + \vartheta^{8/9} \big| \le c^{-1} e^{-cL^{1/8}}; \qquad \big| G^+ (t, x) - F(t, x) - \vartheta^{8/9} \big| \le c^{-1} e^{-cL^{1/8}}.
		\end{flalign}
	\end{enumerate}

\end{lem}

\begin{proof}
	
	This will follow from \Cref{equationcompareboundary}, together with an appropriate rescaling. To establish the first statement of the lemma, define the function $\widehat{G}^{\pm} : \big[ 0, 2 / (L-1) \big] \times [-1, 1] \rightarrow \mathbb{R}$ by setting 
	\begin{flalign}
		\label{gtx2} 
		\widehat{G}^{\pm} (t, x) = \displaystyle\frac{2L}{L-1} \cdot G^{\pm} \bigg( \displaystyle\frac{L-1}{2L} \cdot t, \displaystyle\frac{L-1}{2L} \cdot x +  \displaystyle\frac{1}{2} \bigg),
	\end{flalign}
	
	\noindent for any index $\pm \in \{ +, - \}$ and pair  $(t, x) \in \big[ 0, 2 / (L-1) \big] \times [-1, 1]$; these functions satisfy \eqref{equationxtd} by \Cref{scale21} of \Cref{invariancesscale}. We next apply \Cref{equationcompareboundary}, with the $(F_1, F_2; L; r; \varepsilon; B)$ there equal to $\big( \widehat{G}^-, \widehat{G}^+; (L-1) / 2; 1 / 3; \varepsilon / 2; 4B \big)$ here, using the fact (from the second statement of \eqref{fg1fg2} with the bound $L > 4$) that 
	\begin{flalign*} 
		\| \widehat{G}^{\pm} \|_{\mathcal{C}^5 (\breve{\mathfrak{S}})} \le \frac{2L}{L-1} \cdot \| G^{\pm} \|_{\mathcal{C}^5 (\breve{\mathfrak{S}})} \le 3 \| G^{\pm} - F \|_{\mathcal{C}^5 (\breve{\mathfrak{S}})} + 3 \| F \|_{\mathcal{C}^5 (\breve{\mathfrak{S}})} \le 3 B + 12 C_1 C_2 \vartheta^{1/m} \le 4B,
	\end{flalign*} 

	\noindent for sufficiently small $\vartheta$, to verify the bound on the $\| F_i \|_{\mathcal{C}^5}$ assumed there. This yields a constant $c_1 = c_1 (\varepsilon, B, m) > 0$ such that 
	\begin{flalign*}
		\displaystyle\sup_{t \in [0, 2/(L-1)]}	\displaystyle\sup_{|x| \le 2/3} \big| \widehat{G}^+ (t, x) - \widehat{G}^- (t, x) \big| \le c_1^{-1} e^{-c_1 L^{1/8}}.
	\end{flalign*}
	
	\noindent Together with \eqref{gtx2}, this implies for sufficiently small $\vartheta$ (and thus sufficiently large $L$, as $L > |\log \vartheta|^{20}$) that  
	\begin{flalign*}
		\displaystyle\sup_{t \in [0, 1/L]} \displaystyle\sup_{x \in [1/4, 3/4]} \big| G^+ (t, x) - G^- (t, x) \big| \le \displaystyle\sup_{t \in [0,2/(L-1)]} \displaystyle\sup_{|x| \le 2/3} \big| \widehat{G}^+ (t, x) - \widehat{G}^- (t, x) \big| \le c_1^{-1} e^{-c_1 L^{1/8}},
	\end{flalign*}
	
	\noindent and thus the first statement of the lemma.
	
	\begin{figure}
	\center
\includegraphics[scale=.6, trim = 0 1cm 0 1cm]{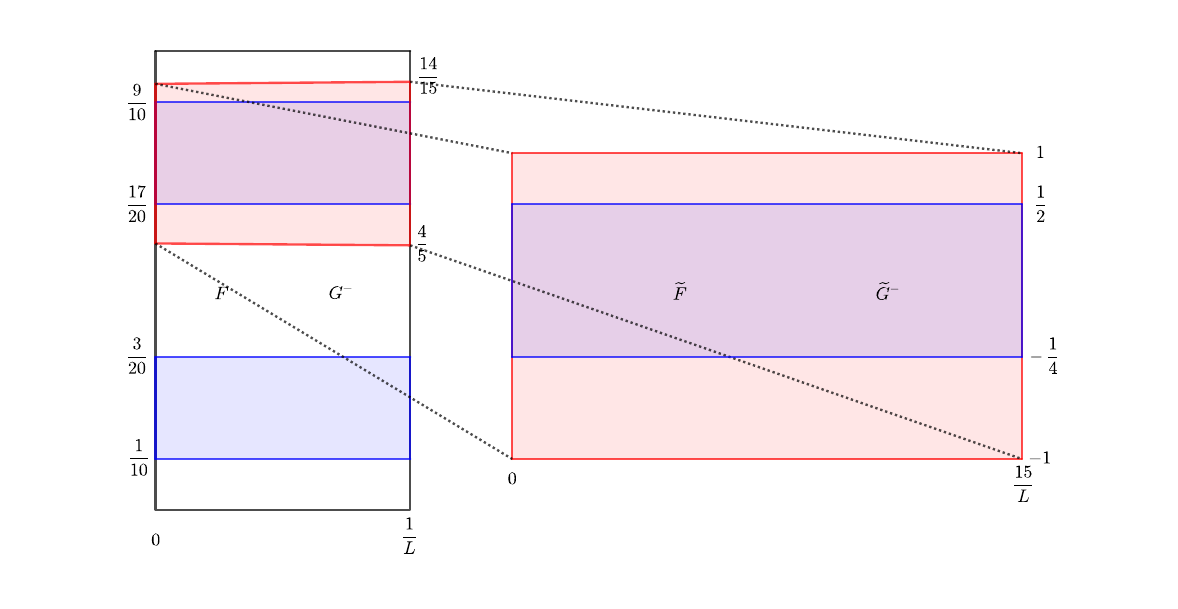}

\caption{Shown above is a depiction of the rescaling used in the proof of the second part of \Cref{g120}.}
\label{f:zoomin}
	\end{figure}

	Next we establish the second part of the lemma. Since the derivation of both statements are entirely analogous, we only detail that of the first, namely, of the bound $\big| G^- (t, x) - F(t, x) + \vartheta^{8/9} \big| \le c^{-1} e^{-cL^{1/8}}$ for $t \in [0, L^{-1}]$ and $x \in [ 1 / 10, 3 / 20 ] \cup [ 17 / 20, 9 / 10]$; we also only address the case when $x \in [ 17 / 20, 9 / 10]$, as the proof in the complementary case is again very similar. As before, we begin by rescaling, namely, we define the functions $\widetilde{F}, \widetilde{G}^- : [0, 15L^{-1}] \times [-1, 1]$ by setting 
	\begin{flalign}
		\label{ftx1} 
		\widetilde{F} (t, x) =  15 F \Big( \displaystyle\frac{t}{15}, \displaystyle\frac{x+13}{15} \Big) - 15 \vartheta^{8/9}; \qquad \widetilde{G}^- (t, x) = 15 G^- \Big( \displaystyle\frac{t}{15}, \displaystyle\frac{x+13}{15} \Big).
	\end{flalign}
	
	\noindent for any pair $(t, x) \in [0, 15L^{-1}] \times [0, 1]$; see \Cref{f:zoomin}. These functions satisfy \eqref{equationxtd} by \Cref{scale21} of \Cref{invariancesscale}. We then apply \Cref{equationcompareboundary}, with the parameters $(F_1, F_2; L; r; \varepsilon; B)$ there equal to $( \widetilde{G}^-, \widetilde{F}; L / 15; 1 / 4; \varepsilon / 2; 20B)$ here (to verify the bounds on the $\| F_i \|_{\mathcal{C}^2}$ assumed there, we used the facts that $\| \widetilde{F} \|_{\mathcal{C}^2} \le 15 \| F \|_{\mathcal{C}^2} \le 15 B$ and $\| \widetilde{G}^- \|_{\mathcal{C}^2} \le 15 \| G^- \|_{\mathcal{C}^2} \le 15 (B + 4C_2 \vartheta^{1/m}) \le 20 B$, by \eqref{fg1g2} and taking $\vartheta$ sufficiently small), which yields a constant $c_2 = c_2 (\varepsilon, B, m) > 0$ such that 
	\begin{flalign*}
		\displaystyle\sup_{t \in [0, 15/L]} \displaystyle\sup_{x \in [-1/4, 1/2]} \big| \widetilde{G}^- (t, x) - \widetilde{F} (t, x) \big| \le \displaystyle\sup_{t \in [0, 15/L]} \displaystyle\sup_{|x| \le 3/4} \big| \widetilde{G}^- (t, x) - \widetilde{F} (t, x) \big| \le c_2^{-1} e^{-c_2 L^{1/8}}.
	\end{flalign*}
	
	\noindent Together with \eqref{ftx1}, this yields
	\begin{flalign*}
		\displaystyle\sup_{t \in [0, 1/L]} \displaystyle\sup_{x \in [17/20, 9/10]} \big| G^- (t, x) - F(t,x) + \vartheta^{8/9} \big| & \le 15 \displaystyle\sup_{t \in [0, 15/L]} \displaystyle\sup_{|x| \le 3/4} \big| \widetilde{G}^- (t, x) - \widetilde{F} (t, x) \big| \\
		& \le 15 c_2^{-1} e^{-c_2 L^{1/8}},
	\end{flalign*}
	
	\noindent which verifies the first bound in \eqref{gtxftx1} when $x \in [ 17 / 20, 9 / 10 ]$. As mentioned previously, the proof of this estimate when $x \in [ 1 / 10, 3 / 20]$ is entirely analogous, as is the proof of the second statement of \eqref{gtxftx1}; this establishes the lemma. 
\end{proof}

Now we can quickly establish \Cref{fg1g2}. 

\begin{proof}[Proof of \Cref{fg1g2}]
	
	As indicated above, we may assume that $\ell = 1$. The $(G^-, G^+)$ of this lemma will be taken to be $(G^-|_{\mathfrak{S}}, G^+ |_{\mathfrak{S}})$ here. Then the first statements of \eqref{g1g0} and \eqref{gi1gi2} together verify that $(G^-, G^+)$ satisfy the first statement of the lemma; moreover, the first statement of \eqref{g1g0} with the second statement of \eqref{gi1gi2} verify the second statement of the lemma. The second bound in \eqref{fg1fg2} (with the fact that $\| F \|_{\mathcal{C}^{m-5} (\mathfrak{S})} \le \| F \|_{\mathcal{C}^m (\mathfrak{R})} \le B$) verifies the third statement of the lemma, and the first bound in \eqref{fg1fg2} verifies the fourth. The first part of \Cref{g120} verifies the fifth part of the lemma, and its second part (together with the fact that $\vartheta^{8/9} - c^{-1} e^{-cL^{1/8}} \ge \vartheta$ for sufficiently small $\vartheta$, since $L > |\log \vartheta|^{20}$) verifies the sixth.
\end{proof}

\section{Proofs of Results From Chapter \ref{EDGESHAPE}}
\label{a:proofdensity}

\subsection{Further Properties of Free Convolutions}

\label{ConvolutionProperties2}

In this section we collect some properties of free convolutions with a rescaled semicircle distribution (which is essentially due to \cite{FCSD}, but stated as below in \cite{ERGIM}) through the following lemma, which will be used repeatedly below. In what follows, we recall the definitions related to free convolutions from \Cref{TransformConvolution} (including the Stieltjes transform $m_0$ of $\mu$ from \eqref{mz0}, the function $M$ and set $\Lambda_t$ from \eqref{mtlambdat}, and the density $\varrho_t \in L^1 (\mathbb{R})$ of $\mu_t = \mu \boxplus \mu_{\semci}^{(t)}$ with respect to Lebesgue measure).

\begin{lem}[{\cite[Lemma 2.3]{ERGIM}}]
	
	\label{yconvolution}
	
	The following statements hold, for any real number $t > 0$.
	
	\begin{enumerate} 
		
		\item \label{i:boundary} Define the function $v_t : \mathbb{R} \rightarrow \mathbb{R}_{\ge 0}$ by setting 
		\begin{flalign}
			\label{vte} 
			v_t (u) = \displaystyle\inf \bigg\{ v \ge 0 :  \displaystyle\int_{-\infty}^{\infty} \displaystyle\frac{\mu(dx)}{(u-x)^2 + v^2} \le t^{-1} \bigg\},
		\end{flalign}
		
		\noindent for each $u \in \mathbb{R}$. Then, $v_t$ is continuous on $\mathbb{R}$. Moreover, the boundary of $\Lambda_t$ is parameterized by $\partial\Lambda_t=\{E+\ri v_t(E): E\in \bR\}$, and the set $\big\{ E\in \bR: v_t(E)>0 \big\}$ consists of countably many open intervals $\bigcup_{i\geq 1}(a_i,b_i)$. 
		\begin{enumerate} 
			\item  \label{vesum} For each $E\in \bigcup_{i\geq 1}(a_i,b_i)$, we have  $\int_{-\infty}^\infty \big|x-E-\ri v_t(E) \big|^{-2} \mu (dx) = t^{-1}$.
			\item For each $E\in\bR\setminus \bigcup_{i\geq 1}(a_i,b_i)$, we have $\int_{-\infty}^\infty \big| x-E-\ri v_t(E) \big|^{-2} \mu(dx) \leq t^{-1}$.
		\end{enumerate}
		\item \label{i:mu} We have $\supp \mu \subseteq \bigcup_{i\geq 1}[a_i,b_i]$, and $\mu \big( \bigcup_{i\geq 1}\{a_i,b_i\} \big)=0$.
		
		\item	\label{i:bijection} The function $M \big(E+\ri v_t(E) \big)$ is (strictly) increasing in $E \in \bR$. Moreover, $v_t (E)$ and $M \big( E + \mathrm{i} v_t (E) \big)$ are smooth for $E$ in the interior of $\mathbb{R} \setminus \bigcup_{i\ge 1} \{ a_i, b_i \}$. 
		
		\item
		\label{wy}A real number $y \in \mathbb{R}$ satisfies $\varrho_t (y) > 0$ if and only if $y =M(w)= w - t m_0 (w)$ for some $w = w(y) \in \partial \Lambda_t \cap \mathbb{H}$. Moreover, the function $w(y)$ is smooth  in $y \in \big\{ y' \in \mathbb{R} : \varrho_t (y') > 0 \big\}$. 
		\begin{enumerate} 
			\item \label{rhoyw} We have $\varrho_t (y) = \pi^{-1} \Imaginary m_0(w)= (\pi t)^{-1} \Imaginary w$.
			\item The Hilbert transform of $\varrho_t$ is given by $H\varrho_t (y)= \pi^{-1} \Real m_0(w) = (\pi t)^{-1} \Real (w-y)$. 
		\end{enumerate}
		
		\item \label{i:edge}  Denote $\mathfrak{e}_+ = \mathfrak{e}_+ (t) = \max \supp \mu_t$ and $\mathfrak{e}_- = \mathfrak{e}_- (t) = \min \supp \mu_t$, and let
		\begin{flalign*} 
			w_+ = \sup \Bigg\{ w \in \mathbb{R} : \displaystyle\int_{-\infty}^{\infty} \displaystyle\frac{\mu (d x)}{(x-w)^2} > t^{-1} \Bigg\}; \qquad w_- = \inf \Bigg\{ w \in \mathbb{R} : \displaystyle\int_{-\infty}^{\infty} \displaystyle\frac{\mu (d x)}{(x-w)^2} > t^{-1} \Bigg\}.
		\end{flalign*}
		Then $\mathfrak{e}_+=w_+-tm_0(w_+)$ and $\mathfrak{e}_- = w_- - tm_0 (w_-)$.
		
	\end{enumerate}

\end{lem}

\subsection{Proof of \Cref{p:densityest}}

\label{Proofrhox} 
	
In what follows, we recall the notation from \Cref{TransformConvolution} and adopt \Cref{blx12}. We denote the Stieltjes transforms of $\nu$ and $\nu_\tau$ as $m = m^{\nu} : \mathbb{H} \rightarrow \mathbb{H}$ and $m_{\tau} = m^{\nu_{\tau}} : \mathbb{H} \rightarrow \mathbb{H}$, respectively. We also recall the function $M = M^{\nu}$ and set $\Lambda_{\tau}$ from \eqref{mtlambdat}, as well as the function $v_{\tau}$ from \eqref{vte}, which is continuous by the first part of \Cref{yconvolution}. In what follows, we abbreviate $v = v_{\tau}$ and further recall from \Cref{i:boundary} of \Cref{yconvolution} that $\partial \Lambda_{\tau} = \big\{ E + \mathrm{i} v(E) : E \in \mathbb{R} \big\}$.

By part (\ref{vesum}) of \Cref{yconvolution}, we have 
\begin{align}\label{e:cha}
    1=\tau\int_{-\infty}^\infty \frac{\nu(dx)}{\big| x-E-\ri  v(E) \big|^2}
    =\tau\int_{-\infty}^\infty \frac{\nu(dx)}{(x-E)^2+  v(E)^2}, \qquad \text{if $v(E) > 0$}.
\end{align}

\noindent Moreover, define the functions $w, y : \mathbb{R} \rightarrow \overline{\mathbb{H}}$ by setting (here, we recall $M$ from \eqref{mtlambdat})
\begin{flalign}
	\label{wy2} 
	w(E) = E + \mathrm{i} v(E), \quad \text{and} \quad y(E) = M \big( w(E) \big) = E + \mathrm{i} v(E) - \tau m \big(E + \mathrm{i} v(E) \big).
\end{flalign} 

\noindent By \Cref{i:boundary} of \Cref{yconvolution}, we have $w(E) \in \partial \Lambda_{\tau}$. Together with \Cref{mz}, this implies that $y(E) \in \mathbb{R}$, so \Cref{wy} of \Cref{yconvolution} gives 
\begin{align}\label{e:yfunc}
     y(E) =  E - \tau \Real m \big(E+\ri  v(E) \big), \qquad \text{so} \quad y(E) = E - \tau \Real m \big( E + \mathrm{i} v(E) \big) \in \supp \varrho_\tau, \quad \text{if $v(E) > 0$},
\end{align}

 \noindent and also that 
\begin{align}\label{e:rhoty}
\varrho_\tau \big( y(E) \big)=\pi^{-1} \Imaginary m \big( w(E) \big) = (\pi \tau)^{-1}  v(E).
\end{align}

\noindent Moreover, for any $E \in \mathbb{R}$, \Cref{i:boundary} in \Cref{yconvolution} gives 
\begin{align}\label{e:cha3}
    \int_{-\infty}^\infty \frac{\nu(dx)}{\big|x-E-\ri  v(E) \big|^2}\leq \tau^{-1}.
\end{align}

Then \eqref{e:cha3} and the fact that $\nu (\mathbb{R}) = L^{3/2}$ together imply for each $E \in \mathbb{R}$ that
\begin{align}\label{e:mvbound}
\Big| m \big(E+\ri  v(E) \big) \Big| = \Bigg| \displaystyle\int_{-\infty}^{\infty} \displaystyle\frac{\nu(dx)}{x-E-\mathrm{i} v(E)} \Bigg| \leq  \Bigg(\int_{-\infty}^{\infty} \frac{ \nu(dx)}{|x-E-\ri  v(E) |^2}\int_{-\infty}^{\infty} \nu(dx) \Bigg)^{1/2} \le \frac{L^{3/4}}{\tau^{1/2}}.
\end{align}

Using \eqref{e:yfunc} and \eqref{e:mvbound}, it follows for any $E \in \mathbb{R}$ that 
\begin{flalign}
	\label{mey}
	\big| y(E) - E \big| = \tau \Big| \Real m \big( E + \mathrm{i} v(E) \big) \Big| \le \tau^{1/2} L^{3/4}.
\end{flalign} 

We next have the following lemma bounding the density $\varrho_{\tau}$ and the endpoints of its support. Below, we define 
\begin{flalign}
	\label{yy} 
	 y_- = \inf (\supp \varrho_{\tau}); \quad y_+ = \sup (\supp \varrho_{\tau}).
\end{flalign}

\begin{lem} 
	
	\label{rhoxyy} 
   
   The following two statements hold.
\begin{enumerate}
	\item For any $x \in \mathbb{R}$, we have $\varrho_\tau (x) \le \pi^{-1} B^{1/2} L^{3/4}$.  
	\item \label{yyestimate} We have $-BL-2 B^{1/2} L^{3/4}\leq y_-\leq y_+\leq 2 B^{1/2} L^{3/4}$.
\end{enumerate}
\end{lem} 

\begin{proof}
From \eqref{e:rhoty}, \eqref{wy2}, and \eqref{e:mvbound}, we have 
\begin{flalign*} 
	\varrho_\tau (x) \le \pi^{-1} \Big| m \big(x + \mathrm{i}v(x) \big) \Big| \le \displaystyle\frac{L^{3/4}}{\pi \tau^{1/2}} \leq \displaystyle\frac{B^{1/2} L^{3/4}}{\pi},
\end{flalign*} 

\noindent where in the last bound we used the fact that $\tau \ge B^{-1}$. This verifies the first part of the lemma. 

From Item \ref{i:edge} in \Cref{yconvolution}, we have $y_- = \min \supp \varrho_\tau = y (E_-)$ and $y_+ = \max \supp \varrho_\tau = y (E_+)$, where $E_-, E_+ \in \bR$ are supremum and infimum, respectively, over all  real numbers $E_0$ satisfying
\begin{align*} 
1 < \tau\int_{-\infty}^{\infty} \frac{\nu(dx)}{|x-E_0 |^2}.
\end{align*}

\noindent Hence, for each index $\pm \in \{ +, - \}$, we can estimate $E_\pm$ through the bound
\begin{align*}
1 \leq \frac{\tau \int_{-\infty}^{\infty}\nu(dx)}{\dist(E_\pm, \supp \nu)^2}=  \frac{\tau L^{3/2}}{\dist(E_\pm, \supp \nu)^2}, \qquad \text{so that} \qquad \dist(E_\pm, \supp \nu) \le \tau^{1/2} L^{3/4},
\end{align*}

\noindent where we used the fact that $\nu (\mathbb{R}) = L^{3/2}$ in the second estimate above. Since by \Cref{blx12} we have $\supp \nu \subseteq  [-BL, 0]$, it follows that 
\begin{flalign*} 
	-BL- \tau^{1/2} L^{3/4}\leq E_\pm \leq \tau^{1/2} L^{3/4}.
\end{flalign*}

\noindent Together with \eqref{mey} and the fact that $y_{\pm} = y(E_{\pm})$, this yields 
\begin{align*}
&y_-\geq E_--L^{3/4} \tau^{1/2} \geq -BL-2 \tau^{1/2} L^{3/4}\geq -BL-2 B^{1/2} L^{3/4}; \\
&y_+\leq E_++L^{3/4} \tau^{1/2} \leq 2 \tau^{1/2} L^{3/4}\leq 2 B^{1/2} L^{3/4},
\end{align*}

\noindent where in the last inequalities of the above bounds we used the fact that $\tau \le B$; this implies the second part of the lemma.
\end{proof}

The below lemma shows that if $v$ is bounded above on $[a, b]$ then, up to a multiplicative factor, $\nu_{\tau} \big( [y(a), y(b)] \big)$ is lower bounded by $\nu \big( [a, b] \big)$; it is established in \Cref{Proofnu2} below.

\begin{lem}\label{l:densitycompare}

Fix real numbers $a<b$. If  $v(E)\le (b-a)/2$ for each real number $E \in [a, b]$, then
\begin{align}\label{e:rhonubound}
\nu_\tau \Big( \big[ y(a), y(b)\big] \Big)\geq \frac{1}{8\pi} \cdot \nu \big([a, b] \big).
\end{align}
\end{lem}

We now fix positive real parameters $\mathfrak{M}, r \in \mathbb{R}$ so that 
\begin{flalign}
	\label{mr} 
	\fM>2; \qquad r= \big\lceil \log_\fM (BL+2B^{1/2} L^{3/4}) \big\rceil.
\end{flalign}

\noindent We further fix a sequence of numbers $y_0>y_1>\cdots>y_{r+1}$ defined by setting
\begin{align}\label{e:difyi}
 y_0 = 0, \qquad \text{and} \qquad   y_{i}-y_{i+1}=\fM^i, \quad \text{for each $i \in \llbracket 0, r \rrbracket$}.
\end{align}

Then, recalling the endpoints $y_-$ and $y_+$ of $\supp \varrho_{\tau}$ from \eqref{yy}, we have
\begin{flalign}
	\label{yr1e} 
	y_{r+1}\leq -\fM^{r}\leq -(BL+2B^{1/2} L^{3/4})\leq y_-,
\end{flalign} 

\noindent where the last inequality follows from \Cref{yyestimate} of \Cref{rhoxyy}; hence, $\supp \varrho_\tau \subseteq  [y_{r+1}, y_+]$. From \eqref{e:difyi} (and \eqref{mr}), we also have $-2\fM^{i-1}\leq y_i\leq 0$ for each integer $i \in \llbracket 0, r+1 \rrbracket$. Recalling the map $y(E)$ from \eqref{wy2} (and \eqref{e:yfunc}), Item \ref{i:bijection} in \Cref{yconvolution} indicates $y : \mathbb{R} \rightarrow \mathbb{R}$ is an increasing bijection. Therefore, there exist real numbers $E_0 > E_1 > \cdots > E_{r+1}$ such that for each $i \in \llbracket 1, r+1 \rrbracket$ we have  
\begin{align}\label{e:defEi}
    y_i=y(E_i)=E_i+\ri v(E_i)-\tau m \big(E_i+\ri v(E_i) \big)=E_i-\tau\Real m \big(E_i+\ri v(E_i) \big).
\end{align}

Defining the real numbers $E_+, E_- \in \mathbb{R}$ by 
\begin{flalign*}
	E_+ = \sup \Bigg\{ E \in \mathbb{R} : \displaystyle\int_{-\infty}^{\infty} \displaystyle\frac{\mu (dx)}{(x-w)^2} > t^{-1} \Bigg\}; \qquad E_- = \sup \Bigg\{ E \in \mathbb{R} : \displaystyle\int_{-\infty}^{\infty} \displaystyle\frac{\mu(dx)}{(x-w)^2} > t^{-1} \Bigg\},
\end{flalign*}

\noindent \Cref{i:edge} of \Cref{yconvolution}, \eqref{yr1e}, and the fact that $y$ is increasing together yield $y(E_-) = y_- \geq y_{r+1}=y(E_{r+1})$, so $E_-\geq E_{r+1}$. Since $\supp \nu \subseteq [E_-, E_+]$ by \Cref{i:mu} of \Cref{yconvolution}, we thus have
\begin{flalign}
	\label{er1e}  
	\supp \nu \subseteq [E_-, E_+]\subseteq [E_{r+1}, E_+], 
\end{flalign} 

\noindent and so by \Cref{blx12} (which also indicates that $B \ge 1$, so $\mathfrak{M}^r \ge BL \ge L$) it follows that
\begin{align}\label{e:ind1}
\nu \big([E_{r+1},E_+] \big)=\nu(\bR)=L^{3/2}\leq (\fM^r)^{3/2}=\fM^{3r/2}.
\end{align}
Using \eqref{mey}, we find for any integer $i \in \llbracket 0, r + 1 \rrbracket$ that $|E_{i}-y_{i}|\leq \tau \big| m(E_i+\ri v(E_i)) \big| \leq \tau^{1/2} L^{3/4}$, and so
\begin{align}\label{e:ind2}
    E_0-E_{r+1}\geq y_0-y_{r+1}-2\tau^{1/2} L^{3/4}\geq \fM^{r}-2\tau^{1/2} L^{3/4}\geq \displaystyle\frac{\fM^r}{3},
\end{align}

\noindent where in the last bound we used that $\mathfrak{M}^r \ge BL + 2 B^{1/2} L^{3/4}$ and $BL \ge  B^{1/2} L^{3/4}\geq \tau^{1/2}L^{3/4}$ (as $B, L > 1$ and $\tau\leq B$).

The following lemma bounds the $E_i$ and $\nu \big( [E_{i+1}, E_i] \big)$ under \Cref{blx122}. 

\begin{lem}\label{c:dyadicbound}
Adopt \Cref{blx122}. Fix two constants $c= (2^{9/2} \pi B)^{-1} > 0$ and $\fK > 1$ with $\fM^{\fK/2-4}\geq 72 c^{-3} B$. For any integers $k \in \llbracket \fK, r+1 \rrbracket$ and $i \in \llbracket k, r + 1 \rrbracket$, we have
\begin{align}\begin{split}\label{e:ind}
    & E_0-E_{k}\geq c\cdot \fM^{k-1}; \qquad \qquad \qquad E_{i}-E_{i+1}\geq c \cdot \fM^{i}; \\
    & \nu \big( [E_{k}, E_+] \big) \leq c^{-1} \cdot \fM^{3(k-1)/2}; \qquad \nu \big([E_{i+1}, E_i] \big) \leq c^{-1} \cdot \fM^{3i/2},
\end{split}\end{align}

\noindent where the second and fourth statements of \eqref{e:ind} are empty if $k=r+1$. Moreover,  
\begin{align}\label{e:etaEbb}
     v(E)\leq \Big( \displaystyle\frac{2\tau}{c}\Big)^{1/2}\fM^{3k/4}, \qquad \text{for each real number $E \ge E_{k}$}.
\end{align}
\end{lem}

\begin{proof}
We prove the lemma by induction on $r-k+1$. The statement \eqref{e:ind} holds for $k=r+1$ by \eqref{e:ind1} and \eqref{e:ind2}. We therefore fix an integer $k \in \llbracket \mathfrak{K}, r+1 \rrbracket$ and assume that the statement \eqref{e:ind} holds for $k+1$. We will then prove that \eqref{e:etaEbb} holds for $k+1$ and that \eqref{e:ind} holds for $k$. We begin with the former.

Fix a real number $E \ge E_{k+1}$. If $v(E) = 0$, then \eqref{e:etaEbb} holds. Otherwise, $v(E) > 0$, so  
\begin{align}
	\label{e:numass}
\begin{aligned}
\frac{1}{\tau}=\int_{-\infty}^{\infty} \frac{\nu(dx)}{(x-E)^2+  v(E)^2}
&\leq \int_{E_{k+2}}^{E_+} \frac{\nu(dx)}{  v(E)^2}
+\int_{E_{r+1}}^{E_{k+2}} \frac{\nu(dx)}{(x-E)^2}\\
&\leq \frac{ \nu \big( [E_{k+2},E_+] \big)}{  v(E)^2}
+\sum_{i = k+2}^r \frac{ \nu \big([E_{i+1}, E_i] \big)}{(E_{k+1}-E_i)^2}\\
&\leq \frac{ \fM^{3(k+1)/2}}{  c v(E)^2}
+\sum_{i = k+2}^r \frac{\fM^{3i/2}}{c (c \fM^{i-1})^2} =\frac{\fM^{3(k+1)/2}}{  c v(E)^2}
+\frac{4}{c^3 \fM^{(k-2)/2}}.
\end{aligned}
\end{align}

\noindent where the first statement follows from \eqref{e:cha}; the second from \eqref{er1e}; the third from the bound $|x-E| \ge E_{k+1} - E_i$ whenever $E \ge E_{k+1}$ and $x \in [E_{i+1}, E_i]$ with $i \ge k+1$; the fourth from the inductive hypothesis (the second, third and fourth statements of \eqref{e:ind}, applied with $k$ there replaced by $i-1 \ge k + 1$, by $k + 2$, and by $i \ge k+2$ here, respectively); and the fifth from performing the sum (and using the fact that $\mathfrak{M} > 2$). It follows from \eqref{e:numass} that, for $E \in [E_{k+1}, E_+]$, we have
\begin{align}\label{e:etaEbb2}
     v(E)\leq \Big( \frac{2\tau}{c} \Big)^{1/2}\fM^{3(k+1)/4}, 
\end{align}

\noindent since $c^{3}\fM^{(k-2)/2}\geq c^{3}\fM^{(\fK-2)/2}\geq 72B \ge 8\tau$. This verifies \eqref{e:etaEbb} with its $k$ replaced by $k+1$.

We next show \eqref{e:ind}, beginning with the first two statements there. From the defining relation \eqref{e:defEi}, we have 
\begin{align*}
    y_0=E_0-\tau \Real m \big( E_0+\ri v(E_0) \big); \qquad y_k=E_k-\tau \Real m \big(E_k+\ri v(E_k) \big).
\end{align*}
By taking the difference and using \eqref{e:difyi}, we get  
\begin{align}\label{e:mdiff0}
    \fM^{k-1}\leq y_0-y_k\leq (E_0-E_k)
    +\tau \Big| m \big(E_0+\ri v(E_0) \big) - m \big(E_k+\ri v(E_k) \big) \Big|
\end{align}

\noindent To estimate the right side of this inequality, we bound $m'$. To this end, for any complex number of the form $z=E+\ri \mathfrak{v} \in \overline{\mathbb{H}}$ with $\mathfrak{v}\geq  v(E)$, we have from \eqref{mz0} that
\begin{align}
	\label{mzderivativeeta}
	\begin{aligned}
    \big| m'(z) \big| = \left|\int_{-\infty}^{\infty} \frac{\nu(dx)}{(x-z)^2}\right|
    \leq\int_{-\infty}^{\infty} \frac{\nu(dx)}{(x-E)^2+ \mathfrak{v}^2}
    \leq \int_{-\infty}^{\infty} \frac{\nu(dx)}{(x-E)^2+  v(E)^2}\leq\frac{1}{\tau},
    \end{aligned}
\end{align}
where we used \eqref{e:cha3} for the last inequality.
 
Thus, to bound $m \big(E_0+\ri v(E_0) \big) - m \big(E_k+\ri v(E_k) \big)$ we introduce the parameter 
\begin{flalign*} 
	\widetilde{\mathfrak{v}}=\Big( \frac{2\tau}{c} \Big)^{1/2}\fM^{3(k+1)/4}\geq  v(E), 
\end{flalign*}

\noindent where the last inequality holds for any  $E \ge E_{k+1}$ by \eqref{e:etaEbb2}. In particular, the vertical segments from $E_0+\ri v(E_0)$ to $E_0+\ri \widetilde{\mathfrak{v}}$ and from  $E_k+\ri v(E_k)$ to $E_k+\ri \widetilde{\mathfrak{v}}$, as well as the horizontal segment from $E_0+\ri \widetilde{ \mathfrak{v}}$ to $E_k+\ri \widetilde{\mathfrak{v}}$ are all in the domain $\big\{ z=E+\ri \mathfrak{v}: \mathfrak{v}\geq  v(E) \big\}$. Thus,
\begin{align}\label{e:mdiff}
	\begin{aligned}
   \Big|m \big( & E_0+\ri v(E_0) \big)-m \big(E_k+\ri v(E_k) \big) \Big| \\
  & \leq \Big|m \big(E_0+\ri v(E_0) \big)-m(E_0+\ri \widetilde{\mathfrak{v}}) \Big| + \big| m(E_0+\ri \widetilde{\mathfrak{v}})-m(E_k+\ri \widetilde{\mathfrak{v}}) \big| \\
  & \qquad  + \Big| m(E_k+\ri \widetilde{\mathfrak{v}})- m \big( E_k+\ri v(E_k) \big) \Big| \\
  &\leq \tau^{-1} \Big( \big( \widetilde{\mathfrak{v}}- v(E_0) \big) + |E_0-E_k|+ \big( \widetilde{\mathfrak{v}} - v(E_k) \big) \Big)\leq \tau^{-1} \bigg( E_0-E_k +2\Big( \frac{2\tau}{c} \Big)^{1/2} \fM^{3(k+1)/4}\bigg),
  \end{aligned}
\end{align}

\noindent where in the second inequality we applied (and integrated) \eqref{mzderivativeeta} and in the third we used the definition of $\widetilde{\mathfrak{v}}$. By plugging \eqref{e:mdiff} into \eqref{e:mdiff0}, we conclude that 
\begin{align}
	\label{e0ek}
   E_0-E_k\geq \frac{1}{2}\bigg( \fM^{k-1}-2 \Big( \displaystyle\frac{2\tau}{c} \Big)^{1/2} \fM^{3(k+1)/4} \bigg) \geq \frac{\fM^{k-1}}{3},
\end{align}

\noindent since $\fM^{(k-7)/4}\geq \fM^{(\fK-7)/4}\geq 2^{3/4} \cdot 6 c^{-3/2} B^{1/2} \ge 6 (2\tau / c)^{1/2}$ (as $c < 1$ and $\tau \le B$). This verifies the first statement of \eqref{e:ind}. 

The proof of the second is similar. In particular, by \eqref{e:difyi} and \eqref{mey}, we have
\begin{align*}
    \fM^{k}= y_k-y_{k+1}\leq (E_k-E_{k+1})
    +\tau \Big| m \big(E_k+\ri v(E_k) \big)-m \big(E_{k+1}+\ri v(E_{k+1}) \big) \Big|,
\end{align*}

\noindent and by following the derivation of \eqref{e:mdiff} we have
\begin{flalign*}
	\Big| m \big( E_k + \mathrm{i} v (E_k) \big) - m \big( E_{k+1} + \mathrm{i} v(E_{k+1}) \big) \Big| \le \tau^{-1} \bigg( E_k - E_{k+1} + 2 \Big( \displaystyle\frac{2\tau}{c} \Big)^{1/2} \mathfrak{M}^{3(k+1)/4} \bigg).
\end{flalign*}

\noindent Together, these two bounds (as in \eqref{e0ek}) yield 
\begin{align}\label{e:distbb}
   E_k-E_{k+1}\geq \frac{\fM^k}{3},
\end{align}

\noindent giving the second bound in \eqref{e:ind}.

To prove the third and fourth bounds in \eqref{e:ind}, beginning with the latter, we use \Cref{l:densitycompare}. First, we have
\begin{align}\label{e:massbound}
    2^{3/2}B  \fM^{3k/2}\geq \nu_\tau \big([y_{k+1}, y_+] \big) \geq \nu_\tau \big([y_{k+1}, y_k] \big),
\end{align}
where the first inequality is from \eqref{e:rhotintbound} and the fact that $|y_{k+1}|\leq 2 \fM^{k}$ (by \eqref{e:difyi} and \eqref{mr}). From \eqref{e:distbb} and \eqref{e:etaEbb2}, we also have for each $E\in[E_{k+1}, E_k]$ that
\begin{align}
	\label{ve0e}
     v(E)\leq \Big( \frac{2\tau}{c} \Big)^{1/2}\fM^{3(k+1)/4}\leq \displaystyle\frac{\fM^k}{6} \leq \displaystyle\frac{E_k-E_{k+1}}{2},
\end{align}

\noindent where we have additionally used the bound $\fM^{(k-3)/4}\geq \fM^{(\fK-3)/4}\geq 2^{5/4} \cdot 6 c^{-3/2} B^{1/2} \ge 6 (2\tau / c)^{1/2}$. By \eqref{e:massbound} and \Cref{l:densitycompare} (whose assumption on the upper bound for $v$ is verified by \eqref{ve0e}), this gives
\begin{align*}
    2^{3/2}B \fM^{3k/2}\geq \nu_\tau \big([y_{k+1}, y_k] \big)\geq\frac{1}{8\pi} \cdot \nu \big([E_{k+1}, E_k] \big), \quad \text{so that} \quad  \nu \big([E_{k+1}, E_k] \big)\leq 2^{9/2}\pi B\fM^{3k/2},
\end{align*}
which gives the fourth estimate in \eqref{e:ind}. The proof of the third estimate there is entirely analogous and thus omitted. This establishes the lemma.
\end{proof}

Now we can quickly establish \Cref{p:densityest}. 

\begin{proof}[Proof of \Cref{p:densityest}]

Throughout this proof, we adopt the notation of \Cref{c:dyadicbound}. The first statement of the proposition follows from the second part of \Cref{rhoxyy}, so it suffices to establish \eqref{e:densitybb0}. This will follow from \eqref{e:etaEbb}. Indeed, fix a real number $x \in \mathbb{R}$; by \Cref{mz}, \Cref{i:boundary} of \Cref{yconvolution}, and \eqref{wy2}, there exists a real number $E = E(x) \in \mathbb{R}$ such that $x = y(E)$. We may assume in what follows that $x \ge y_- \ge -BL - 2B^{1/2} L^{3/4} \ge -\mathfrak{M}^r$ (where the bound follows from \Cref{yyestimate} of \Cref{rhoxyy}), for otherwise $\varrho_{\tau} (x) = 0$; similarly, we may assume that $x \le y_+$.

We first consider the case when $x \in [y_{i+1}, y_i]$ for some integer $i\geq \fK$, so that $E \in [E_{i+1}, E_i]$ (by \eqref{e:defEi} and the fact that $y$ is increasing in $E$, by \Cref{i:bijection} of \Cref{yconvolution}). Then, $-\mathfrak{M}^r \le -BL - 2B^{1/2} L^{3/4} \le y_- \le x \le y_i \le - \mathfrak{M}^{i-1}$, so $i \in \llbracket \mathfrak{K}, r+1 \rrbracket$. In particular,
\begin{align}
   \varrho_\tau(x)= \varrho_\tau \big( y(E) \big) =\frac{ v(E)}{\pi\tau}\leq\Big( \frac{2\tau}{c} \Big)^{1/2} \fM^{3(i+1)/4}\leq \Big( \frac{2\tau}{c} \Big)^{1/2}\fM^{3/2}|x|^{3/4}\leq C|x|^{3/4},
\end{align}

\noindent for $C \ge (2B / c)^{1/2} \mathfrak{M}^{3/2}$. Here, the first statement follows from the fact that $x = y(E)$; the second from \eqref{e:rhoty}; the third from \eqref{e:etaEbb}; and the fourth from the fact that $x \le y_i \le -\mathfrak{M}^{i-1}$. This establishes the proposition if $x \le y_{\mathfrak{K}}$. 

If instead $x=y(E)\in [y_{\fK}, y_+]$, so that $E \ge E_{\mathfrak{K}}$, then by analogous reasoning we have
\begin{align}
   \varrho_\tau(x)= \varrho_\tau(y(E))=\frac{ v(E)}{\pi\tau}\leq \Big( \frac{2\tau}{c} \Big)^{1/2} \fM^{3\fK/4}\leq C,
\end{align}

\noindent for $C\geq ( 2B / c)^{1/2} \fM^{3\fK/4}$, establishing the proposition in this case as well. 
\end{proof}

\subsection{Proof of \Cref{l:densitycompare}} 

\label{Proofnu2} 

In this section we establish \Cref{l:densitycompare}; throughout, we adopt the notation of that lemma, as well as that from \Cref{TransformConvolution} and \Cref{Proofrhox}.

\begin{figure}
	\center
\includegraphics[scale=1, trim=3.5cm 0.5cm 0 0.5cm, clip]{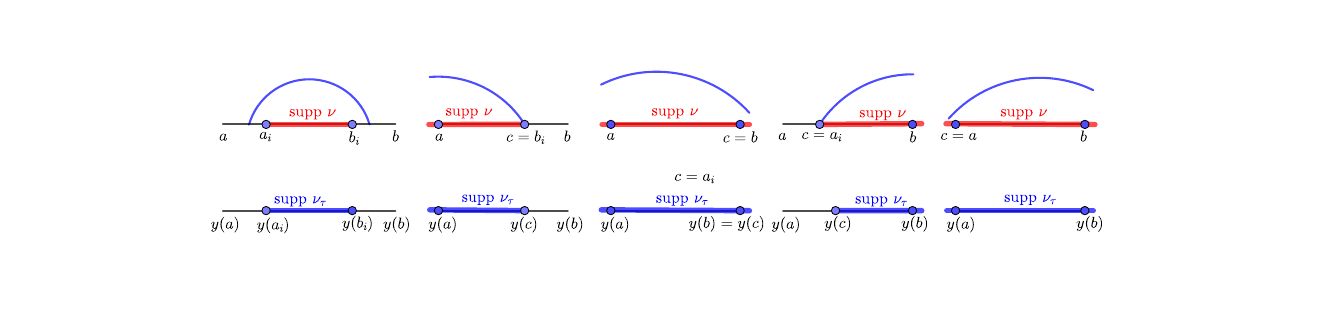}

\caption{Depicted above are the five different cases of \Cref{yab}.}
\label{f:mass_compare}
	\end{figure}

From \Cref{yconvolution}, we have that $v$ is smooth on the set $\big\{E\in \bR: v(E)>0 \big\}$ (by \Cref{wy}), which consists of countably many open intervals $\bigcup_{i\geq 1}(a_i,b_i)$ (by \Cref{i:boundary}); that $\supp \nu \subseteq I$, where $I =\bigcup_{i\geq 1}[a_i,b_i]$, and $\nu \big( \bigcup_{i\geq 1}\{a_i,b_i\} \big)=0$ (by \Cref{i:mu}); and that $\supp \nu_{\tau} \subseteq \bigcup_{i\geq 1} \big[ y(a_i),y(b_i) \big]$ (by \Cref{wy}, \Cref{i:boundary}, and the fact that $y$ is increasing in $E$ by \Cref{i:bijection}). Thus \eqref{e:rhonubound} is equivalent to 
	\begin{align}\label{e:rhonubound2}
		\nu_\tau \Big( \big[ y(a), y(b) \big]\cap y(I) \Big)\geq \frac{1}{8\pi} \cdot \nu \big([a, b]\cap I \big).
	\end{align}
	The intersection $[a,b]\cap I$ consists of two types of intervals:  intervals $[a_i, b_i]$ contained in $[a,b]$, and intervals $[a,c]$ or $[c,b]$ (for some $c \in \bigcup_{i \ge 1} \{ a_i, b_i \}$) containing an endpoint $a$ or $b$ of $I$; see \Cref{f:mass_compare}. The estimate \eqref{e:rhonubound2} follows from summing the two statements of the following lemma.

	\begin{lem} 
		
		\label{yab} 
		
		Adopting the notation and assumptions of \Cref{l:densitycompare}, the following two statements hold.
		
	\begin{enumerate}
		\item  \label{i:wholeinterval} For any integer $i \ge 1$, we have $\nu_\tau \big([y(a_i), y(b_i)] \big)= \nu \big([a_i, b_i] \big)$.
		\item \label{i:partialinterval} 
		\begin{enumerate} 
			\item \label{bci} Assume that $c=b_i$ for some $i \ge 1$ and  $v(E)>0$ for each $E \in [a, b_i)$, or that $c=b$ and $v(E)>0$ for each $E \in [a, b]$. Then, $\nu_\tau \big([y(a), y(c)] \big)\geq (8 \pi)^{-1}  \cdot  \nu \big([a, c] \big)$.
			\item \label{aci} Assume that $c=a_i$ for some $i \ge 1$ and $v(E)>0$ for each $E \in (a_i, b]$, or that $c=a$ and $v(E)>0$ for each $E \in [a, b]$. Then, $\nu_\tau([y(c), y(b)])\geq (8 \pi)^{-1} \cdot \nu \big( [c, b] \big)$. 
		\end{enumerate} 
	\end{enumerate}

	\end{lem}

		\begin{proof}[Proof of \Cref{l:densitycompare}]
			
			Summing the result of \Cref{yab} over all intervals in $\big[ y(a), y(b) \big] \cap y(I)$ yields \eqref{e:rhonubound2}, which as mentioned above, implies the lemma.
	\end{proof} 
	
	We next establish the first part of \Cref{yab}. 
	
	\begin{proof}[Proof of \Cref{i:wholeinterval} of \Cref{yab}]
	Recalling from \eqref{wy2} that $w(E)=E+\ri v(E)$, we have from \Cref{yconvolution} that $w \big( [a_i, b_i] \big)\in \partial\Lambda_\tau$ (by \Cref{i:boundary}), and $w(a_i), w(b_i)\in \bR$ (since $v(a_i) = v(b_i) = 0$ by \Cref{i:boundary}). Then,
	\begin{align}
		\label{nmintegral} 
		\begin{aligned}
			\nu_\tau \Big( \big[y(a_i), y(b_i) \big] \Big)  =\int_{y(a_i)}^{y(b_i)}\varrho_\tau(y)d y &  = \displaystyle\frac{1}{\pi \tau} \int_{a_i}^{b_i}  \Imaginary w(E) \cdot d \big(w(E)-\tau m(w(E) \big)\\
			&  =\Imaginary \frac{1}{\pi \tau}\left( w \big(w-\tau m(w) \big) \Big|_{w(a_i)}^{w(b_i)} -\int_{w([a_i, b_i])} \big(w-\tau m(w) \big)d w\right)\\
			&   =\Imaginary \frac{1}{\pi} \int_{w([a_i, b_i])} m(w)d w.
		\end{aligned} 
	\end{align}
	
	\noindent Here, in the first equality we used the fact that $\nu_{\tau}$ has density $\varrho_{\tau}$; in the second we used \eqref{e:rhoty} and \eqref{wy2}; in the third, we integrated by parts, using the fact that $w(E) - \tau m \big( w(E) \big) = y(E) \in \mathbb{R}$; and in the fourth we used the facts $w(a_i)  \in \mathbb{R}$ and $w(b_i)  \in \bR$ (as $v(a_i) = 0 = v(b_i)$), and that $m \big( w (a_i) \big), m \big( w(b_i) \big) \in \mathbb{R}$ (as $\Imaginary m \big( w(a_i) \big) = \Imaginary w(a_i) = 0$, where the first statement is due to \Cref{rhoyw} of \Cref{yconvolution}, and similarly for $m \big(w(b_i) \big)$). 
	
	Abbreviating the set $\omega= w \big([a_i, b_i] \big)\cup \overline{w \big([a_i, b_i] \big)}$ (which does not intersect the real interval $(a_i, b_i)$, since $\Imaginary w(E) = v(E) > 0$ for $E \in (a_i, b_i)$), the above expression can be written as a contour integral along $\omega$ counterclockwise, by
	\begin{align*}
		\nu_{\tau} \Big( \big[ y(a_i), y(b_i) \big] \Big) &= -\Imaginary \frac{1}{\pi} \int_{w([a_i, b_i])} \Bigg( \displaystyle\int_{-\infty}^{\infty} \displaystyle\frac{\nu(dx)}{w-x} \Bigg) d w \\
		& =\Imaginary \frac{1}{2\pi} \oint_{\omega} \int_{-\infty}^\infty \frac{\nu(d x)}{w-x} dw\\
		&=\Imaginary \frac{1}{2\pi} \int_{-\infty}^\infty  \oint_{\omega}  \frac{dw}{z-x} \nu(dx) =\int_{a_i}^{b_i}\nu(dx) =\nu \big( [a_i,b_i] \big), 
	\end{align*}
	
	\noindent where in the first equality we used \eqref{nmintegral} and the definition \eqref{mz0} of $m$; in the second we used the definition of $\omega$ and the fact that $\overline{(z-x)^{-1}} = (\overline{z} - x)^{-1}$; in the third we interchanged the order of integration between $x$ and $w$; in the fourth we applied the residue theorem; and in the fifth we used the fact that $\nu \big(\{a_i, b_i\} \big)=0$ from Item \ref{i:mu} of \Cref{yconvolution}. This confirms the first statement of the lemma.
	\end{proof}
	
	To establish the second part of \Cref{yab}, we require the below integral estimate.
	
	\begin{lem} 
		
		\label{vintegral}

		Adopt the notation and assumptions of \Cref{i:partialinterval} of \Cref{yab}. Fix a real number $x \in [a, c]$ such that $x \ne c$ if $c = b_i$ for some $i \ge 1$, and $x \ne a$ if $a = a_i$ for some $i \ge 1$. We have
		\begin{flalign*}
			\displaystyle\int_a^c \displaystyle\frac{v(E)^3 \big( 1 + v'(E) \big)^2 dE}{\big( (x-E)^2 + v(E)^2 \big)^2} \ge \displaystyle\frac{1}{16}.
		\end{flalign*}
	\end{lem} 

	\begin{proof}
		
		Throughout this proof, we adopt the notation and assumptions of \Cref{bci} of \Cref{yab}, as the proof is entirely analogous under \Cref{aci} of that lemma. It suffices to show that 
		\begin{align}
			\label{e:lowb}
			\begin{aligned}
				\int_{x}^{c} \frac{v(E)^3 \big(1+v'(E)^2 \big)d E}{\big( (x-E)^2+  v(E)^2 \big)^2}\geq \frac{1}{16}, \qquad \text{if $x \le \displaystyle\frac{a+c}{2}$}; \\
				\displaystyle\int_a^x \displaystyle\frac{v(E)^3 \big( 1 + v'(E) \big) dE}{\big( (x-E)^2 + v(E)^2 \big)^2} \ge \displaystyle\frac{1}{16}, \qquad \text{if $x \ge \displaystyle\frac{a+c}{2}$}.
			\end{aligned}
		\end{align}
		
		\noindent We only show the first statement of \eqref{e:lowb}, as the proof of the second is entirely analogous; so, we assume that $x \le (a+c) / 2$ in what follows. Then, $x \in [a, c)$ with $x \ne a$ if $a = a_i$; by \Cref{i:boundary} of \Cref{yconvolution}, this implies that $v(x) > 0$, so $c = x + \lambda \cdot v(x)$ for some real number $\lambda = \lambda (a, c, x) > 0$. Observe that 
		\begin{flalign}
			\label{lambdacb0} 
			\lambda \ge 1, \qquad \quad \text{or} \qquad \quad c = b_i, \quad \text{so} \quad v(c) = 0.
		\end{flalign}
		
		\noindent Indeed, the fact that $v(c) = 0$ when $c = b_i$ follows from \Cref{i:boundary} and \Cref{yconvolution} (and the continuity of $v$). If instead $c \ne b_i$, then we must have $c = b$, in which case $c -x \ge (b-a) / 2 \ge v(x)$ (where the last bound follows from the assumptions of \Cref{l:densitycompare}), and so $\lambda \ge 1$. It then suffices to show 
		\begin{align}\label{e:toshow}
			\int_{x}^{c} \frac{v(E)^3 \big(1+v'(E)^2 \big)d E}{\big((x-E)^2+  v(E)^2 \big)^2}=\int_{x}^{x+\la v(x)} \frac{v(E)^3 \big(1+v'(E)^2 \big)d E}{\big((x-E)^2+  v(E)^2 \big)^2}\geq \frac{1}{16}
		\end{align}
		
		To prove \eqref{e:toshow}, we first define the function $f : [0, 1] \rightarrow \mathbb{R}$ by setting   
		\begin{flalign*} 
			f(\theta)= \displaystyle\frac{v \big(x+v(x)\theta \big)}{v(x)} > 0, \qquad \text{so that} \qquad f'(\theta) = v' \big( x + v(x) \theta \big).
		\end{flalign*} 
		
		\noindent Then $f(0)=1$; moreover, by \eqref{lambdacb0} we have $v \big( x + \lambda \cdot v(x) \big) = v(c) = 0$ if $\lambda < 1$, meaning that 
		\begin{flalign}
			\label{flambda0lambda1}
			f(\lambda) = 0, \qquad \text{if} \quad \lambda < 1. 
		\end{flalign}
		
		\noindent Changing variables $E = x + \theta \cdot v(x)$, we then find that \eqref{e:toshow} is equivalent to
		\begin{align}\label{e:lowbb}
			\int_0^\la \frac{f(\theta)^3 \big(1+f'(\theta)^2 \big)}{\big(\theta^2+f(\theta)^2 \big)^2}d \theta\geq\frac{1}{16}.
		\end{align}
		
		We now consider several cases. First, if $\lambda \ge 1 / 2$ and $\theta \le f(\theta) \le 2$ for each $\theta \in [ 0, 1 / 2]$, then   
		\begin{align*}
			\int_0^\la \frac{f(\theta)^3 \big(1+f'(\theta)^2 \big)}{\big(\theta^2+f(\theta)^2 \big)^2}d \theta
			&\geq \int_0^{1/2} \frac{f(\theta)^3 \big( 1+f'(\theta)^2 \big)}{\big( 2f^2(\theta) \big)^2}d \theta=\int_0^{1/2} \frac{1+f'(\theta)^2}{4f(\theta)}d \theta
			\geq \int_0^{1/2} \frac{d \theta}{4 \cdot 2}
			=\frac{1}{16},
		\end{align*}
		
		\noindent where in the third statement we used the facts that $f' (\theta)^2 \ge 0$ and that $f(\theta) \le 2$. If otherwise either $\lambda < 1 / 2$ or $\theta \le f(\theta) \le 2$ does not hold for some $\theta \in [ 0, 1 / 2 ]$, then set 
		\begin{flalign*} 
			\theta_0 =\inf \big\{ \theta \ge 0: f(\theta)<\theta \text{ or } f(\theta)>2 \big\}.
		\end{flalign*} 
		
		\noindent Then, we have $\theta_0 \le \lambda$. Indeed, if $\lambda\leq 1 / 2$, then $f(\lambda)=0<\lambda$ (by \eqref{flambda0lambda1}), so  $\theta_0 \le \lambda$. If instead $\lambda \ge 1 / 2$, then either $f(\theta) < \theta$ or $f(\theta) > 2$ for some $\theta \in [ 0, 1 / 2 ]$; in this case, we have $\theta_0 \le 1/2 \le \lambda$. Thus, 
		\begin{align}
			\label{ftheta3}
			\int_0^\la \frac{f(\theta)^3 \big(1+f'(\theta)^2 \big)}{\big(\theta^2+f^2(\theta) \big)^2}d \theta
			&\geq 
			\int_0^{\theta_0} \frac{1+f'(\theta)^2}{4f(\theta)}d \theta \geq \Bigg|  \int_0^{\theta_0} \frac{f'(\theta) d \theta}{2f(\theta)} \Bigg|
			\geq \displaystyle\frac{1}{2} \cdot \bigg| \log \displaystyle\frac{f(\theta_0)}{f(0)} \bigg| \geq \frac{\ln(2)}{2},
		\end{align}
		
		\noindent where in the first inequality we used the facts that $\theta_0 \le \lambda$ and $\theta \le f(\theta)$ for $\theta \in [0, \theta_0]$; in the second we used the fact that $1 + f'(\theta)^2 \ge 2 \big| f'(\theta) \big|$; in the third we performed the integration; and in the fourth we used the fact that either $f(\theta_0) \le 1/2$ or $f(\theta_0) \ge 2$ (and that $f(0) = 1$). Indeed, to verify the latter, observe since $f$ is continuous (as $v$ is continuous and $v(x) \ne 0$) that we either have $f(\theta_0) \le \theta_0$ or $f(\theta_0) \ge 2$. It suffices to address the former case; if $\lambda \le 1/2$, then $f(\theta_0) \le \theta_0 \le \lambda \le 1 / 2$; otherwise, we must have that $\theta_0 \le 1/2$, and so $f(\theta_0) \le \theta_0 \le 1/2$. The bound \eqref{ftheta3} then finishes the proof of \eqref{e:lowbb} and thus of the lemma.
	\end{proof}

	Now we can establish the second part of \Cref{yab}. 
	
	\begin{proof}[Proof of \Cref{i:partialinterval} of \Cref{yab}]
		
		We will only establish \Cref{bci} of the lemma, as the proof of \Cref{aci} is entirely analogous. Then, $v(E)>0$ for each $E \in (a, c)$. By Item \ref{i:bijection} of \Cref{yconvolution}, $y(E)$ is smooth in $E \in (a, c)$. Applying \eqref{e:rhoty}, it follows that
	\begin{align}\label{e:mutau}
		\nu_\tau \Big( \big[y(a), y(c) \big] \Big) 
		=\int_{y(a)}^{y(c)}\varrho_\tau(y)d y
		= \displaystyle\frac{1}{\pi \tau} \int_{a}^{c} v(E)d y(E).
	\end{align}
	
	\noindent To evaluate $y'(E)$, we differentiate both sides of the definition \eqref{wy2} of $y$ to find
	\begin{align}\label{e:yE}
		y'(E)= \big(1+\ri v'(E) \big) \Big(1-\tau m' \big(E+\ri  v(E) \big) \Big) \in \bR,
	\end{align}
	
	\noindent where the last inclusion follows from the fact that $y(E) \in \mathbb{R}$ for each $E \in \mathbb{R}$ (by \Cref{mz} and \Cref{i:boundary} of \Cref{yconvolution}). It follows that there exits some real number $r(E)\in \bR$ such that
	\begin{align}\label{e:defaE}
		1-\tau m' \big( E+\ri  v(E) \big)=r(E) \big( 1 - \ri v'(E) \big),\quad \text{so that} \quad y'(E)=r(E) \big(1+v'(E)^2 \big).
	\end{align}
	By taking real parts on both sides of the first equation in \eqref{e:defaE}, we get 
	\begin{flalign} 
		\label{mre} 
		r(E) =1-\tau \Real m' \big( E+\ri  v(E) \big).
	\end{flalign} 
	
	\noindent To evaluate $\Real m' \big( E+\ri  v(E) \big)$, observe from the definition \eqref{mz0} of $m$ that
	\begin{align}\label{e:redm}
		\Real m' \big(E+\ri  v(E) \big) =\Real \int_{-\infty}^\infty\frac{\nu(dx)}{\big( (x-E)-\ri v(E) \big)^2}=\int_{-\infty}^{\infty} \frac{ \big( (x-E)^2-  v(E)^2 \big) \nu(dx)}{\big( (x-E)^2+  v(E)^2 \big)^2}
	\end{align}
	Thus using \eqref{mre}, \eqref{e:cha}, and \eqref{e:redm}, we can compute $r(E)$ for $E \in (a, c)$ by    
	\begin{align}\begin{split}\label{e:aE}
			r(E) & = \tau\left(\int_{-\infty}^{\infty} \frac{\nu(dx)}{(x-E)^2+  v(E)^2} 
			-\int_{-\infty}^{\infty} \frac{ \big( (x-E)^2-  v(E)^2 \big)\nu(dx)}{\big( (x-E)^2+  v(E)^2 \big)^2}\right)
			\\
			&=2 \tau v(E)^2 \int_{-\infty}^{\infty} \frac{\nu(dx)}{\big( (x-E)^2+  v(E)^2 \big)^2} >0.
	\end{split}\end{align}
	
	By plugging \eqref{e:yE}, \eqref{e:defaE} and \eqref{e:aE} back into \eqref{e:mutau}, we obtain 
	\begin{align*}
			\nu_\tau \Big( \big[y(a), y(c) \big] \Big) = \displaystyle\frac{1}{\pi \tau}\int_{a}^{c} v(E) y'(E) d E
			& = \displaystyle\frac{1}{\pi \tau} \int_{a}^{c} v(E) r(E) \big( 1+v'(E)^2 \big) d E\\
			&=
			\frac{2}{\pi}\int_{-\infty}^{\infty} \nu(dx)\int_{a}^{c} \frac{v(E)^3 \big(1+v'(E)^2 \big)d E}{\big((x-E)^2+  v(E)^2 \big)^2}.
	\end{align*}
	
	\noindent Together with \Cref{vintegral}, this gives 
	\begin{align*}
		&\nu_\tau \Big( \big[y(a), y(c) \big] \Big)\geq \frac{1}{8\pi} \cdot \nu \big([a,c] \big),\qquad \qquad \qquad \qquad \quad \text{if $c = b$, $b \ne b_i$, and $a \ne a_i$}; \\
		& \nu_{\tau} \Big( \big[ y(a), y(c) \big] \Big) \ge \displaystyle\frac{1}{8 \pi} \cdot \nu \big( (a, c] \big) = \displaystyle\frac{1}{8 \pi} \cdot \nu \big( [a, c] \big), \qquad \text{if $a = a_i$ and $c \ne b_i$}; \\ 
		&\nu_\tau \Big( \big[y(a), y(c) \big] \Big)\geq \frac{1}{8\pi} \cdot \nu \big([a,c) \big)=\frac{1}{8\pi} \cdot \nu \big([a,c] \big),\qquad \text{if $c=b_i$ and $a \ne a_i$},
	\end{align*}
	
	\noindent where in the last equalities of the second and third statements we used Item \ref{i:mu} of \Cref{yconvolution}, which indicates that $\nu \big( \{ a_i, b_i\} \big)=0$. Since the case when $[a, c] = [a_i, b_i]$ was addressed in the first part of the lemma, this finishes the proof of the second  part of the lemma.
\end{proof}

	\subsection{Proof of \Cref{c:rhoderbound}}
	
	\label{ProofrhoxDerivative} 
	
	In this section we establish \Cref{c:rhoderbound}. Throughout, we adopt the notation from \Cref{Proofrhox}, recalling in particular the functions $v$, $w$, and $y$ from \eqref{wy2}; the parameters $\mathfrak{M}$ and $r$ from \eqref{mr}, the sequence $0 = y_0 \ge y_1 \ge \ldots \ge y_{r+1}$ from \eqref{e:difyi}; their respective preimages $E_0 \ge E_1 \ge \cdots \ge E_{r+1}$ under $y$ from \eqref{e:defEi}; and the constants $c > 0$ and $\mathfrak{K} > 1$ from \Cref{c:dyadicbound}. 
	
	For any real number $x \in \supp \varrho_\tau$, we let $k = k_x$ denote the minimal integer such that $k\geq \fK$ and $x\geq y_{k+1}$. By \Cref{i:bijection} of \Cref{yconvolution}, there exists a unique real number $E = E(x) \geq E_{k+1}$ with  
	\begin{align}
		\label{xywe}
		x=y(E)=w-\tau m(w), \qquad \text{where} \qquad w = w(E) =E+\ri  v(E).
	\end{align}
	
	\noindent We also have
	\begin{align}\label{e:vEbound}
		\tau\pi \varrho_\tau(x)=  v(E)\leq \Big( \frac{2\tau}{c} \Big)^{1/2}\fM^{3 (k+1)/4},
	\end{align}
	
	\noindent where the first inequality holds by \eqref{e:rhoty}; and the second holds by  \eqref{e:etaEbb} and the fact that $E \ge E_{k+1}$. Moreover, using the fact that $x=y(E)=w-\tau m(w)$, we can interpret $w=w_x$ as a function of $x$. By Item \ref{wy} of \Cref{yconvolution}, the function $w_x$ is smooth in $x \in \big\{ x'\in \mathbb{R} : \varrho_{\tau} (x') > 0 \big\}$.  Thus, differentiating the first equation in \eqref{xywe} with respect to $x$, we find
	\begin{align}\label{e:dw}
		\partial_x w (E) =\frac{1}{1-\tau m'(w)}.
	\end{align}

	The following lemma bounds $\varrho$ and its derivatives. It is established in \Cref{DerivativeProofrho} below.

	\begin{lem}\label{l:densitylowb}
		
		Adopting the notation and assumptions of \Cref{c:rhoderbound}, there exists a constant $C=C(\ell, A, B)>1$ such that for any real number $x\in \big[\gamma_\tau(B/2), \gamma_\tau(2/B) \big]$ we have
		\begin{align}\label{e:densityder}
			\varrho_\tau(x)\geq (2A)^{-3}; \qquad \big| \partial_x^\ell\varrho_\tau(x) \big| \leq C.
		\end{align}
	\end{lem}
	
	Now we can establish \Cref{c:rhoderbound}.
	
	\begin{proof}[Proof of \Cref{c:rhoderbound}]  
		
		By \Cref{l:densitylowb}, for any real number $x\in \big[\gamma_\tau(B/2), \gamma_\tau(2/B) \big]$, we have $\varrho_\tau(x) > 0$. Together with the definition \eqref{e:defgammatau} of $\gamma_{\tau}$, this implies that for any $y\in [2/B, B/2]$ 
		\begin{align}
			\label{yf}
			y = F_\tau \big(\gamma_\tau(y) \big), \qquad \text{where} \quad  F_{\tau} (x) = \int^{\infty}_x \varrho_\tau(x) d x, \quad \text{for any $x \in \mathbb{R}$}.
		\end{align}
		
		\noindent Since $F_\tau'(x)=-\varrho_\tau(x)$, \Cref{l:densitylowb} yields for any integer $\ell \ge 1$ a constant $C_1 = C_1 (\ell, A, B) > 1$ such that $F_{\tau}' (x) \leq -(2A)^{-3}$ and $\big| \partial_x^{\ell} F_\tau(x) \big| \le C_1$ for each $x\in \big[\gamma_\tau(B/2), \gamma_\tau(2/B) \big]$. Together with \eqref{yf} and the Inverse Function Theorem, this yields a constant $C_2 = C_2 (\ell, A, B) > 1$ such that $\big| \partial_y^{\ell} \gamma_\tau(y) \big| \le C_2$, for each $y \in [ 2 / B, B / 2]$, which establishes the proposition.
	\end{proof}	
	
	\subsection{Proof of \Cref{l:densitylowb}}
	
	\label{DerivativeProofrho} 
	
	In this section we establish \Cref{l:densitylowb}, to which end we first show the following lemma bounding  $m$ and its derivatives. Throughout, we recall the notation of \Cref{ProofrhoxDerivative}.

	\begin{lem}\label{l:mzder}
		
		Adopting the notation and assumptions of \Cref{c:rhoderbound}, there exists a constant $C = C(A, B) > 1$ such that the following holds. Fix a real number $x_0 \ge -A$ with $\varrho_{\tau} (x_0) > 0$, and define the associated $w$ as in \eqref{xywe}. We have 
		\begin{align}\label{e:dwbb}
			\big|\partial_w^\ell m(w) \big|\leq \frac{B\ell!}{  (\Imaginary w)^{\ell-1}}; \qquad
			\big| 1-\tau m'(w) \big|\geq \frac{(\Imaginary w)^2}{C}.
		\end{align}
	\end{lem}
	
	\begin{proof}
		
		Throughout this proof, we recall the notation from \Cref{c:dyadicbound}. As in \eqref{xywe}, we let $E = E(x_0)$ be such that $x_0 = y(E)$, so that $w = w(E) = E + \mathrm{i} v(E)$; let $k = k_{x_0} \ge \mathfrak{K}$ be the minimal integer such that $E \ge E_{k+1}$. To deduce the first bound in \eqref{e:dwbb}, observe that 
		\begin{align*}
			\big| \partial_w^\ell m(w) \big| 
			\leq \ell! \displaystyle\int_{-\infty}^{\infty} \frac{\nu(dx)}{|x-w|^{\ell+1}} = \displaystyle\frac{\ell!}{v(E)^{\ell-1}} \displaystyle\int_{-\infty}^{\infty} \displaystyle\frac{\nu (dx)}{|x-w|^2} 
			= \displaystyle\frac{\ell!}{v(E)^{\ell-1}} \cdot \frac{\Imaginary m(w)}{\Imaginary w}
			=\frac{\ell!}{\tau v (E)^{\ell-1}} ,
		\end{align*}
		
		\noindent where the first statement follows from the definition \eqref{mz0} of $m$; the second from the fact that $|x - w| \ge \Imaginary w = v(E)$; the third again from the definition of $m$; and the fourth from the facts that $\Imaginary w = v(E)=\tau \Imaginary m(w)$ from \eqref{e:rhoty}. The first statement in \eqref{e:dwbb} then follows from this, with the facts that $\tau\geq B^{-1}$ and $\Imaginary w = v(E)$.
		
		To deduce the second bound in \eqref{e:dwbb}, first observe from \eqref{e:aE} (and \eqref{mre}) that for any real number $D \ge 0$ we have 
		\begin{flalign} 
			\label{mwv1}
			\begin{aligned}
				\Real \big( 1-\tau m'(w) \big)
				&  = 2\tau v(E)^2 \int_{-\infty}^{\infty} \frac{\nu(dx)}{\big( (x-E)^2+  v(E)^2 \big)^2} \\
				&  \geq 2\tau v(E)^2 \int_{E-D\fM^{3k/4}}^{E+D\fM^{3k/4}} \frac{\nu(dx)}{\big( (x-E)^2+  v(E)^2 \big)^2}\\
				&
				\geq \frac{2\tau v(E)^2}{D^2 \fM^{3k/2}+ v(E)^2}\int_{E-D\fM^{3k/4}}^{E+D\fM^{3k/4}}\frac{\nu(dx)}{(x-E)^2+  v(E)^2},
			\end{aligned}
		\end{flalign} 
		\noindent where in the third line we used the fact that $(x-E)^2 \le D^2 \mathfrak{M}^{3k/2}$ for $x \in [E-D\mathfrak{M}^{3k/4}, E+D\mathfrak{M}^{3k/4}]$. Thus, we must lower bound the integral on the right side of the above inequality. To this end, first observe (following \eqref{e:numass}) that we have the upper bound
		\begin{align}\begin{split}\label{e:numass2}
				&\int_{E_{r+1}}^{E_{k+2}} \frac{\nu(dx)}{(x-E)^2+  v(E)^2}\leq \displaystyle\sum_{i=k+2}^r \displaystyle\frac{\nu \big( [E_{i+1}, E_i] \big)}{(E_{k+1} - E_i)^2} \le \displaystyle\sum_{i=k+2}^r \displaystyle\frac{\mathfrak{M}^{3i/2}}{c(c\mathfrak{M}^{i-1})^2} \le  \frac{4}{c^{3}\fM^{(k-2)/2}}\leq \frac{1}{2\tau},
		\end{split}\end{align}
		
		\noindent where to deduce the first inequality we used the fact that $(x-E)^2 + v(E)^2 \ge (x-E)^2 \ge (E_{k+1} - E_i)^2$ whenever $x \in [E_{i+1}, E_i]$ (as $E \ge E_{k+1}$); to deduce the second we used \eqref{e:ind}; to deduce the third we performed the sum; and to deduce the fourth we used the fact that $\mathfrak{M}^{(k-2)/2} \ge \mathfrak{M}^{(\mathfrak{K}-2)/2} \ge 72c^{-3} B \ge 8 c^{-3} \tau$. 
		
		Further setting the constant $D = ( 8\tau/c)^{1/2} \mathfrak{M}^{3/4}$ and again applying \eqref{e:ind}, we also find (observing that $E + D \mathfrak{M}^{3k/4} \ge E \ge E_{k+1} > E_{k+2}$ and that $\supp \nu \subseteq [E_{r+1}, E_+]$ by \eqref{er1e}) that 
		\begin{align}\begin{split}\label{e:numass3}
				&\int_{x\geq E+D\fM^{3k/4}} \frac{\nu(dx)}{(x-E)^2+  v(E)^2}\leq \frac{\nu\big( [E_{k+2}, E_+] \big)}{D^2 \fM^{3k/2}}\leq  \frac{\fM^{3(k+1)/2}}{c D^2 \fM^{3k/2}}= \frac{1}{8\tau},\\
				&\int_{E_{k+2}\leq x\leq E-D\fM^{3k/4}} \frac{\nu(dx)}{(x-E)^2+  v(E)^2}\leq \frac{\nu \big([E_{k+2}, E_{+}] \big)}{D^2 \fM^{3k/2}}\leq \displaystyle\frac{\mathfrak{M}^{3(k+1)/2}}{c D^2 \mathfrak{M}^{3k/2}} \le \frac{1}{8\tau},
		\end{split}\end{align}
		
		\noindent where in the first inequalities in both estimates we used the fact that $(x-E)^2 + v(E)^2 \ge (x-E)^2 \ge D^2 \mathfrak{M}^{3k/2}$ on the domain of integration. Thus, it follows from combining \eqref{e:numass2} and \eqref{e:numass3}, and using \eqref{e:cha} (and again the fact from \eqref{er1e} that $\supp \nu \subseteq [E_{r+1}, E_+]$) that
		\begin{align*}
			\int_{E-D\fM^{3k/4}}^{E+D\fM^{3k/4}} \frac{\nu(dx)}{(x-E)^2+  v(E)^2}
			=\frac{1}{\tau}-\int_{x\notin [E-D\fM^{3k/4},E+D\fM^{3k/4}] }\frac{\nu(dx)}{(x-E)^2+  v(E)^2}\geq \frac{1}{4\tau}.
		\end{align*}
		
		\noindent Together with \eqref{mwv1}, this yields
		\begin{align}
			\label{mwv}
			\begin{aligned}
				\Real \big( 1-\tau m'(w) \big) \geq \frac{2\tau v^2(E)}{D^2 \fM^{3k/2}+ v(E)^2} \cdot \frac{1}{4\tau}
				\geq \frac{ v(E)^2}{2 \big(D^2 \fM^{3k/2}+ v(E)^2 \big)}, 
			\end{aligned} 
		\end{align}
		
		Using \eqref{e:vEbound} and the definition of $D = ( 8\tau/c)^{1/2} \mathfrak{M}^{3/4}$, we have 
		\begin{flalign*} 
			D^2 \fM^{3k/2}+ v(E)^2\leq 10 c^{-1} \tau\fM^{3(k+1)/2}.
		\end{flalign*} 
		
		\noindent Inserting this into \eqref{mwv} (and using the fact that $\tau \ge B^{-1}$) yields  
		\begin{align}
			\big| 1-\tau m'(w) \big| \ge \Real \big( 1 - \tau m'(w) \big) \geq \frac{c}{20\tau \mathfrak{M}^{3(k+1)/2}} \cdot v(E)^2 \geq \frac{cB}{20 \mathfrak{M}^{3(k+1)/2}} \cdot v(E)^2,
		\end{align}
	
		\noindent which proves the second bound in \eqref{e:dwbb}, since $\Imaginary w = v(E)$ and $\mathfrak{M}^k \le A / 2$ (as $A \ge -x_0 \ge -y_{k+1} \ge 2 \mathfrak{M}^k$, where the last bound holds by \eqref{e:difyi} and \eqref{mr}). 
	\end{proof}

	Now we can establish \Cref{l:densitylowb}.

	\begin{proof}[Proof of \Cref{l:densitylowb}]  
		
		Recall that $\varrho_\tau(x)= (\pi \tau)^{-1} \Imaginary w_x = (\pi \tau)^{-1} v(E)$ from  \eqref{e:rhoty}, where $w = w_x = E + \mathrm{i} v(E)$. By taking derivatives with respect to $x$ on both sides, we get
		\begin{align}\label{e:lthder}
			\big| \partial_x^\ell \varrho_\tau(x) \big|\leq (\tau \pi)^{-1} \cdot \big| \partial_x^\ell w(x) \big|.
		\end{align}
	
		\noindent At $\ell = 1$, this yields for $x \ge -A$ with $\varrho_{\tau} (x) > 0$ that
		\begin{align}\label{e:densitylowerbbb}
			\big| \partial_x \varrho_\tau(x) \big| \leq \frac{1}{\tau\pi \big|1-\tau m'(w) \big|}\leq \frac{C_1}{\tau\pi v(E)^2}=\frac{C_1}{(\tau\pi)^3 \varrho_\tau (x)^2}\leq \frac{C_1 B^3}{\pi^3 \varrho_\tau (x)^2},
		\end{align}
	
		\noindent for some constant $C_1 = C_1 (A, B) > 1$. Here, the first statement uses \eqref{e:lthder} and \eqref{e:dw}; the second uses the estimate \eqref{e:dwbb}; the third uses \eqref{e:rhoty}; and the fourth uses the bound $\tau\geq B^{-1}$. 
		
		By \eqref{e:densitybb0}, for any real number $x\geq \gamma_\tau(B)\geq -A$, we have for some constant $C = C (B) > 1$ that $\varrho_\tau(x)\leq CA^{3/4}$. Together with the definition \eqref{e:defgammatau} of $\gamma_{\tau}$ (and the fact that $\varrho_{\tau}$ is bounded by the third part of \Cref{mutrhot}), this implies for any real numbers $0\leq y\leq y'\leq B$ that
		\begin{align}
			y'-y =\int_{\gamma_{\tau}(y')}^{\gamma_{\tau}(y)}\varrho_\tau(x)d x \leq \int_{\gamma_{\tau}(y')}^{\gamma_{\tau}(y)}CA^{3/4}d x\leq CA^{3/4} \big(\gamma_\tau(y)-\gamma_\tau(y') \big).
		\end{align}
		
		\noindent Hence, $\gamma_\tau(y)-\gamma_\tau(y')\geq (CA^{3/4})^{-1} |y'-y|$ for any $0\leq y\leq y'\leq B$. In particular, by taking $(y,y')=(3B/4, B)$, this yields
		\begin{align}\label{e:gammataul}
			\gamma_{\tau} \Big( \displaystyle\frac{3B}{4} \Big)\geq \gamma_{\tau}(B)+(4CA^{3/4})^{-1} B.
		\end{align}
		
		Now set $\varepsilon=(12 C C_1 A^4 B^3)^{-1}$. Fix any real number $y_0\in [2/B, B/2]$; denote $x_0=\gamma_\tau(y_0)$; and set $y_1 =y_0-\varepsilon/2$ and $y_2 =y_0+\varepsilon/2$. By our choice $\varepsilon\leq B^{-1}$, so $B^{-1}\leq y_1 \leq y_2\leq B$, which gives 
		\begin{align}\label{e:interval}
			\big|\gamma_\tau(y_1)- x_0 \big| = \big|\gamma_\tau(y_1)- \gamma_\tau(y_0) \big|\leq \frac{A\varepsilon}{2};
			\qquad
			\big| \gamma_\tau(y_2)- x_0 \big|= \big|\gamma_\tau(y_2)- \gamma_\tau(y_0) \big|\leq \frac{A\varepsilon}{2},
		\end{align}

		\noindent where we used our hypothesis that, for any $B^{-1}\leq y\leq y'\leq B$ with $y'-y\geq \varepsilon$, we have $|\gamma_\tau(y)-\gamma_\tau(y')|\leq A (y' - y)$. Again using the definition of $\gamma_\tau$ from \eqref{e:defgammatau} (with the fact that $\varrho_{\tau}$ is bounded by the third part of \Cref{mutrhot}), we have
		\begin{align*}
			y_2 - y_1 =\int_{\gamma_{\tau}(y_2)}^{\gamma_{\tau}(y_1)}\varrho_\tau(x)d x & \leq \big( \gamma_{\tau} (y_1) - \gamma_{\tau} (y_2) \big) \cdot \max_{x \in [\gamma_\tau(y_2),  \gamma_\tau(y_1)]}\varrho_\tau(x) \\ 
			&\leq A (y_2 - y_1) \cdot \max_{x \in [\gamma_\tau(y_2), \gamma_\tau(y_1)]}\varrho_\tau(x),
		\end{align*}
		where in the last inequality, we used our assumption that $\big| \gamma_\tau(y)-\gamma_\tau(y') \big|\leq A(y'-y)$ whenever $B^{-1} \le y \le y' \le B$. Hence, there exists some $\widetilde x_0\in \big[ \gamma_\tau(y_1), \gamma_\tau(y_2) \big]$ such that $\varrho_\tau(\widetilde x_0)\geq A^{-1}$. 
		
		By \eqref{e:interval}, we have $|x_0-\widetilde x_0|\leq \max \big\{ |\gamma_\tau(y_1)-x_0|, |\gamma_\tau (y_2)-\widetilde x_0| \big\}\leq A\varepsilon/2$, and so $x_0 \in [ \widetilde{x}_0 - A\varepsilon / 2, \widetilde{x}_0 + A\varepsilon / 2 ] \subseteq [\widetilde x_0-A\varepsilon, \widetilde x_0+A\varepsilon]$. Moreover, using \eqref{e:gammataul}, our assumption $\gamma_\tau(B)\geq -A$, and our choice of $\varepsilon$ (with the facts that $\widetilde{x}_0 \ge \gamma_{\tau} (y_2)$, that $y_2 = y_0 + \varepsilon / 2 \ge (B + \varepsilon) / 2 \le 3B / 4$, and that $\gamma_{\tau}$ is non-increasing), it follows that
		\begin{flalign*} 
			\widetilde x_0-A\varepsilon \geq \gamma_\tau(y_2)-A\varepsilon\geq \gamma_\tau \Big( \displaystyle\frac{3B}{4} \Big) -A\varepsilon\geq \gamma_{\tau}(B)+(4CA^{3/4})^{-1} B-A\varepsilon\geq \gamma_{\tau}(B)\geq -A.
		\end{flalign*} 
	
		\noindent Thus, for any $x\in  [\widetilde x_0-A\varepsilon, \widetilde x_0+A\varepsilon]$ with $\varrho_{\tau} (x) > 0$, \eqref{e:densitylowerbbb} holds. By rearranging that bound, we obtain
		\begin{align}\label{e:dlowerbb3}
			\Big|\partial_x \big( \varrho_\tau(x)^3 \big) \Big|\leq 3 \pi^{-3} C_1 B^3,\qquad \text{for each $x\in [\widetilde x_0-A\varepsilon, \widetilde x_0+A\varepsilon]$ with $\varrho_{\tau} (x) > 0$}.
		\end{align}
		By integrating \eqref{e:dlowerbb3}, we conclude that for any $x\in [\widetilde x_0-A\varepsilon, \widetilde x_0+A\varepsilon]$, we have
		\begin{align}\label{e:dlowbb2}
			\varrho_\tau(x)^3\geq \varrho_\tau(\widetilde x_0)^3- 3 \pi^{-3} C_1B^3 \cdot (2A \varepsilon) \geq A^{-3} - 6 \pi^{-3} C_1 \varepsilon A B^3 \geq (2A)^{-3}.
		\end{align}
		
		\noindent Here, in the second inequality, we used the bound $\varrho_\tau(\widetilde x_0)\geq A^{-1}$; in the last inequality, we used fact fact that $\varepsilon=(12 C C_1A^4 B^3)^{-1} $. Taking $x = x_0 = \gamma_{\tau} (y_0)$, this verifies the first statement in \eqref{e:densityder}.
		
		To establish the second, abbreviate $m^{(k)} (w) = \partial_w^k m(w)$ for each integer $k \ge 0$. To use \eqref{e:lthder}, it is quickly verified using \eqref{e:dw} that 
		\begin{align}\label{e:pset}
			\partial_x^{\ell} w(x) = \displaystyle\sum_{r = 0}^{\ell-1} \displaystyle\sum_{\substack{\bm{\ell} \in \mathbb{Z}_{\ge 2}^r \\ |\bm{\ell}| = \ell + r-1}} D_{\bm{\ell}} \tau^r \cdot \frac{m^{(\ell_1)}(w) m^{(\ell_2)}(w) \cdots m^{(\ell_r)}(w)}{\big( 1-\tau m'(w) \big)^{1 + r}}, 
		\end{align}
		
		\noindent for some constants $\{ D_{\bm{\ell}} \}$, where the second sum is over all $r$-tuples $\bm{\ell} = (\ell_1, \ell_2, \ldots , \ell_r) \in \mathbb{Z}_{\ge 2}^r$ such that $\sum_{i=1}^r \ell_i = \ell + r-1$. 
		Using \eqref{e:rhoty}, the two estimates in \eqref{e:dwbb}, and the first bound in \eqref{e:densityder} (with the facts that $r \le \ell-1$ and $\tau \ge B^{-1}$), the summands on the right side of \eqref{e:pset} are bounded by 
		\begin{align}\begin{split}\label{e:termupbb}
				\Bigg| \displaystyle\frac{m^{(\ell_1)} (w) m^{(\ell_2)} (w) \cdots m^{(\ell_r)} (w)}{\big( 1 - \tau m' (w) \big)^{1+ r}} \Bigg|
				& \le \frac{B^r\ell_1!\ell_2!\cdots \ell_r!}{ v(E)^{\ell-1}} \cdot \left(\frac{C_1}{v(E)^2}\right)^{\ell} \\
				& = \frac{B^r C_1^{\ell}\ell_1!\ell_2!\cdots \ell_r!}{ \big( \tau\pi \varrho_\tau(x) \big)^{3\ell-1}} 
				\le \frac{(2A)^{3\ell-1}B^{4\ell-2} C_1^{\ell}\ell_1!\ell_2!\cdots \ell_r!}{ \pi^{3\ell-1}}.
		\end{split} \end{align}
		Summing over all the terms in the form \eqref{e:pset} using \eqref{e:termupbb} (and the fact that $\tau \le B$), and applying \eqref{e:lthder}, we deduce that there exists a constant $C_2 = C_2 (\ell, A, B) > 1$ such that
		\begin{align}
			\big| \varrho^{(\ell)}_\tau(x) \big|\le (\tau \pi)^{-1} \cdot \big| \partial_x^\ell w(x) \big| \leq C_2.
		\end{align}
		
		\noindent which establishes the second statement in \eqref{e:densityder}.
	\end{proof}

	\subsection{Proof of \Cref{degreed}} 
	
	\label{ProofFunctionG}

	In this section we establish \Cref{degreed} (see \Cref{f:homeo}); throughout, we recall the notation from that proposition. In what follows, for any real number $r > 0$ and point $z \in \mathbb{R}^2$, we let $\mathcal{B}_r (z) = \big\{ z' \in \mathbb{R} : |z' - z| < r \big\}$ denote the open disk of radius $r$ centered at $z$; if $z = (0, 0)$, we abbreviate $\mathcal{B}_r (z) = \mathcal{B}_r$, and if moreover $r = 1$ we abbreviate $\mathcal{B}_1 = \mathcal{B}$. We also denote the winding number of a continuous curve $\breve{\gamma} \subset \mathbb{R}^2$ with respect to a point $z \in \mathbb{R}^2$ by $\wind (\breve{\gamma}; z)$. 
	
	Observe that we may assume that $\mathfrak{R} = \mathcal{B}$, by precomposing $G$ with a strictly positively oriented, real analytic homeomorphism from the unit disk $\mathcal{B}$ to $\mathfrak{R}$ (guaranteed to exist by, for example, the Riemann mapping theorem). Further observe that, since $G$ is real analytic and nonconstant, its set of critical points is discrete; we will use this fact repeatedly in what follows.
	
	We begin with the following two lemmas; the first indicates that $G$ is injective away from its critical points, and the second indicates that $\overline{\mathfrak{W}}$ is in the image of $G$.
	
	\begin{lem} 
		
		\label{g10critical} 
		
		Let $w \in \mathfrak{W}$ be a point in the image of $G$ such that $G^{-1} (w) \subset \mathfrak{R}$ contains no critical point of $G$. Then, $G^{-1} (w)$ consists of one point. 
	\end{lem} 
	
	\begin{proof} 
		
	Assume to the contrary that $G^{-1} (w) = \{ u_1, u_2, \ldots , u_k \} \subset \mathfrak{R}$ for some integer $k > 1$, such that none of the $u_j$ are critcal points of $G$. For each integer $i \in \llbracket 1, k-1 \rrbracket$, let $\ell_j \subset \mathcal{B}$ denote a curve connecting $u_j$ to $u_{j+1}$ that does not pass through a critical point of $G$, such that the interiors of the $\{ \ell_j \}$ are pairwise disjoint.
		
		Let $r > 0$ be a small real number, and let  $\ell_j (r) = \big\{ z \in \mathbb{R}^2 : \dist (z, \ell_j) < r \big\}$. Then, define the sets $\mathfrak{R}', \mathfrak{R}_1, \mathfrak{R}'' \subset \mathcal{B}$ and $\gamma', \gamma_1, \gamma'' \subset \mathcal{B}$ by
		\begin{flalign*} 
			\mathfrak{R}' = \bigcup_{j=1}^k & \mathcal{B}_{100r} (u_j); \qquad \mathfrak{R}_1 = \mathfrak{R}'  \cup \bigcup_{j=1}^k  \ell_j (r); \qquad \mathfrak{R}'' = \mathfrak{R}_1 \setminus \overline{\mathfrak{R}}'; \\ 
			& \qquad \gamma' = \partial \mathfrak{R}'; \qquad \gamma_1 = \partial \mathfrak{R}_1; \qquad \gamma'' = \partial \mathfrak{R}''.
		\end{flalign*} 
		
		\noindent For sufficiently small $r$, the set $\mathfrak{R}_1$ does not contain a critical point of $G$ (as such points are isolated). For sufficiently small $r$, we can also guarantee for $j \ne j'$ that $\mathcal{B}_{200r} (u_j) \subset \mathcal{B}$; that $\mathcal{B}_{100r} (u_j) \cap \mathcal{B}_{100r} (u_{j'})$ is empty; and that $\ell_j (r) \cap \ell_{j'} (r) \subseteq \bigcup_{i=1}^k \mathcal{B}_{2r} (u_i)$. 
		
	\begin{figure}
	\center
\includegraphics[scale=1, trim=0 0.5cm 0 0.5cm, clip]{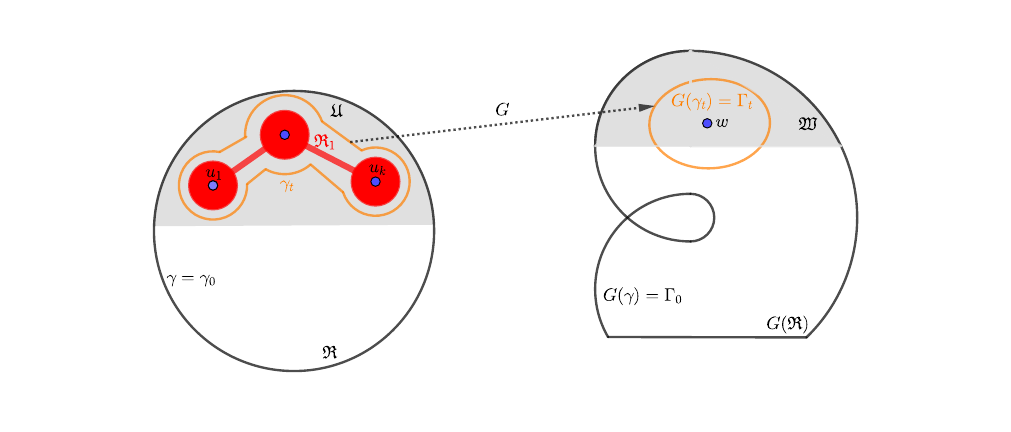}

\caption{Depicted above is an orientation-preserving homotopy $\{\gamma_t\}_{t\in[0,1]}$  from $\gamma_0$ to $\gamma_1=\partial \mathfrak{R}_1$, as in the proof of \Cref{g10critical}.}
\label{f:injection}
	\end{figure}	
		
		It is quickly verified that there exists an orientation-preserving homotopy $\{ \gamma_t \}_{t \in [0, 1]}$ from $\gamma_0 = \gamma$ to $\gamma_1$, such that $\gamma_t \cap \mathfrak{R}_1$ is empty for each $t \in [0, 1]$; see \Cref{f:injection} for a depiction. Since $G$ is positively oriented and continuous, letting $\Gamma_t = G(\gamma_t)$ for each $t \in [0, 1]$, this induces an orientation-preserving homotopy $\{ \Gamma_t \}$ from $\Gamma_0$ to $\Gamma_1$. Since each $\gamma_t$ is disjoint from $G^{-1} (w) \subset \mathfrak{R}_0$, none of the $\Gamma_t$ intersects $w$, meaning that $\wind (\Gamma_t; w)$ is constant in $t \in [0, 1]$. Hence, $\wind (\Gamma_1; w) = \wind (\Gamma_0; w) = 1$, where in the last equality we used the third hypothesis of \Cref{degreed}.
		
		Since $G$ is strictly positively oriented away from its critical points, and none of the $u_j$ are critical points of $G$, the point $w$ is a regular value of $G$. So, we have $\wind  G (\partial \mathcal{B}_{100r} (u_j)); w \big) = 1$ for each $j \in \llbracket 1, k \rrbracket$ and sufficiently small $r > 0$. Thus, $\wind (G(\gamma'); w )= \sum_{j=1}^k \wind ( G(\partial \mathcal{B}_{100r} (u_j)); w ) = k$; we also have $\wind ( G(\gamma''); w ) = 0$ (as $\gamma''$ does not enclose any of the $u_j$). Since $\mathfrak{R}_1 = \mathfrak{R}' \cup \mathfrak{R}''$, it follows that 
		\begin{flalign*} 
			\wind (\Gamma_1; w) = \wind ( G(\partial \mathfrak{R}_1); w ) = \wind ( G(\mathfrak{R}'); w ) + \wind ( G(\mathfrak{R}''); w ) = k.
		\end{flalign*} 
		
		\noindent This is contradicts the fact that $\wind (\Gamma_1; w) = 1$, which confirms the lemma. 
		\end{proof} 
		
		\begin{lem} 
			
			\label{gw} 
			
			The function $G$ surjects onto $\overline{\mathfrak{W}}$. 
		\end{lem} 
		
		\begin{proof}
		
		Since $G$ is continuous and $\overline{\mathfrak{R}}$ is compact, it suffices to show that $G$ surjects onto $\mathfrak{W}$. Suppose to the contrary that this is false, so that there exists some point $w \in \mathfrak{W}$ not in the image of $G$. The first assumption in \Cref{degreed} stipulates the existence of some $w' \in \mathfrak{W}$ that is in the image of $G$. Since $G$ is real analytic and $\mathfrak{W}$ is open, it follows that there are infinitely many points $w_1, w_2, \ldots \in \mathfrak{W}$ not in the image of $G$ and infinitely many points $w_1', w_2', \ldots \in \mathfrak{W}$ in the image of $G$; we may assume that these points are all uniformly bounded away from $\partial \mathfrak{W}$. Since $\mathfrak{W}$ is connected, for each integer $j \ge 1$, there exists a continuous curve $\omega_j : [0, 1] \rightarrow \mathbb{R}$ such that $\omega_j (0) = w_j$ and $\omega_j (1) = w_j'$; we may assume that these curves are pairwise disjoint, in the sense that $\omega_j (r) \ne \omega_{j'} (r')$ unless $(r, j) = (r', j')$. 
		
		Further let $s_j \in [0, 1]$ denote the infimum over all $s \in [0, 1]$ such that $\omega_j (s)$ is in the image of $G$. Then $\omega_j (s_j)$ is in the image of $G$, since $\omega$ and $G$ are continuous; hence, each $s_j \in (0, 1]$, and so the $\omega_j (s_j)$ are mutually distinct over $j \ge 1$. Thus, since the critical points of $G$ are isolated, there exists some integer $j_0 \ge 1$ such that $G^{-1} \big( \omega_{j_0} (s_{j_0}) \big)$ contains no critical points of $G$. Together with the fact that $G$ is strictly positively oriented away from its critical points, this implies that $G^{-1}$ is a diffeomorphism in a neighborhood of $ \omega_{j_0} (s_{j_0})$. Hence, $\omega_{j_0} (s_{j_0} - \varepsilon)$ is in the image of $G$ for sufficiently small $\varepsilon > 0$. This contradicts the minimality of $s_{j_0}$, establishing the lemma. 
		\end{proof} 
	
		Next we show that $G$ is injective, from which we can quickly deduce \Cref{degreed}. 
		
		\begin{lem}
		
		\label{g1}
		
		For any $w \in \overline{\mathfrak{W}}$, there is at most one point $u \in \overline{\mathfrak{R}}$ such that $G(u) = w$. 
		
		\end{lem} 
		
	\begin{proof} 
		
		Assume to the contrary that this is false, so that there exists some $w \in \overline{\mathfrak{W}}$ such that $G^{-1} (w) = \{ u_1, u_2, \ldots , u_k \} \subset \overline{\mathfrak{R}}$ for some integer $k > 1$. By the fourth assumption in \Cref{degreed}, we have $u_j \in \mathfrak{R}$ for each $j \in \llbracket 1, k \rrbracket$ (that is, none of the $u_j$ lie on $\partial \mathfrak{R}$). Denote $\mathfrak{W}_0 = \mathcal{B}_r (w)$ for some small real number $r > 0$ such that $\overline{\mathfrak{W}}_0 \subset \mathfrak{R}$, and set $\mathfrak{U}_0 = G^{-1} (\mathfrak{W}_0) \subset \overline{\mathfrak{R}}$. Since $G$ is continuous and $\overline{\mathfrak{R}}$ is compact, for sufficiently small $r > 0$ we have $\mathfrak{U}_0 = \bigcup_{j=1}^k \mathfrak{U}_j$, for some open sets $\mathfrak{U}_j \subset \mathfrak{R}$ such that $u_j \in \mathfrak{U}_j$ and such that $\overline{\mathfrak{U}}_i$ is disjoint from $\overline{\mathfrak{U}}_j$ for each distinct $i, j \in \llbracket 1, k \rrbracket$.
		
		Set $\mathfrak{U}' = \bigcup_{j=2}^k \mathfrak{U}_j$. We claim that there exists a sequence of distinct points $w_1, w_2, \ldots \in \overline{\mathfrak{W}}$ converging to $w$, such that $\overline{\mathfrak{U}}_1 \cap G^{-1} (w_i)$ and $\overline{\mathfrak{U}}' \cap G^{-1} (w_i)$ are nonempty for each integer $i \ge 1$. We first establish the lemma assuming this claim. Since the set of critical points for $G$ is discrete, we may assume (by taking a subsequence of the $\{ w_i \}$ if necessary) that $G^{-1} (w_i)$ does not contain a critical point of $G$ for each $i \ge 1$. Since $\mathfrak{U}' = \bigcup_{j=2}^k \mathfrak{U}_j$, there exists an integer $j_0 \in \llbracket 2, k \rrbracket$ and an infinite subsequence $w_{i_1}, w_{i_2}, \ldots $ of $(w_1, w_2, \ldots )$ such that $\overline{\mathfrak{U}}_1 \cap G^{-1} (w_{i_m})$ and $\overline{\mathfrak{U}}_{j_0} \cap G^{-1} (w_{i_m})$ are nonempty for each integer $m \ge 1$. However, as $G^{-1} (w_{i_m})$ contains no critical points of $G$, it follows from \Cref{g10critical} that $G^{-1} (w_{i_m}) = v_m$ is one point. Hence, $v_m \in \overline{\mathfrak{U}}_1 \cap \overline{\mathfrak{U}}_{j_0}$, contradicting the disjointness of the $\overline{\mathfrak{U}}_j$. 
		
		It therefore remains to establish the above claim. Set $\mathfrak{W}_1 = G(\overline{\mathfrak{U}_1})$ and $\mathfrak{W}' = G(\overline{\mathfrak{U}}')$; observe that $\mathfrak{W}_1$ and $\mathfrak{W}'$ are closed (since $G$ is continuous, and $\overline{\mathfrak{U}}_1$ and $\overline{\mathfrak{U}}'$ are compact). Since $G$ is real analytic and nonconstant (and $G(u_2) = w$), there exists a sequence of distinct points $u_1', u_2', \ldots \in \mathfrak{U}'$ converging to $u_2$ such that, denoting $w_i' = G(u_i')$ for each $i \ge 1$, we have $w_1', w_2', \ldots \in \mathfrak{W}'$ are mutually distinct and converge to $w$. Thus, if $\mathfrak{W}_1$ contains a neighborhood $\mathfrak{W}_1'$ of $w$, we would be able to take $\{ w_1, w_2, \ldots \} = \{ w_1', w_2', \ldots \} \cap \mathfrak{W}_1'$, confirming the claim.

		Otherwise, there exists a sequence of mutually distinct points $p_1, p_2, \ldots \notin \mathfrak{W}_1$ converging to $w$. Moreover, (again since $G$ is real analytic and nonconstant) there exists a sequence of mutually distinct points $q_1', q_2', \ldots \in \mathfrak{U}_1$ converging to $u_1$ such that, denoting $p_i' = G(q_i')$ for each $i \ge 1$, we have $p_1', p_2', \ldots \in \mathfrak{W}_1$ are mutually distinct and converge to $w$. For each integer $j \ge 1$, let $\omega_j : [0, 1] \rightarrow \mathbb{R}$ denote a continuous curve with $\omega_j (0) = p_j$ and $\omega_j (1) = p_j'$; since $p_j, p_j' \in \overline{\mathfrak{W}}_0 = \overline{\mathcal{B}}_r (w)$, we may assume that $\omega_j (r) \in \overline{\mathfrak{W}}_0$ for each $j \ge 1$ and $r \in [0, 1]$; we may also assume that the $\omega_j$ are mutually disjoint, in that $\omega_j (r) = \omega_{j'} (r')$ if and only if $(j, r) = (j', r')$. 
		
		Then, let $s_j = \inf \big\{ s \in [0, 1] : \omega_j (s) \in \mathfrak{W}_1 \big\}$ and set $w_j = \omega_j (s_j)$, so that the $\{ w_j \}$ are mutually distinct. We have $w_j \in \mathfrak{W}_1$ (as $\mathfrak{W}_1$ is closed and $\omega_j$ is continuous), and $s_j > 0$ (as $\omega_j (0) = p_j \notin \mathfrak{W}_1$). Moreover, observe that $\mathfrak{W}_1 \cup \mathfrak{W}' = \overline{\mathfrak{W}}_0$, since $\overline{\mathfrak{W}}_0 = G (\overline{\mathfrak{U}}_0) = G (\overline{\mathfrak{U}}_1 \cup \overline{\mathfrak{U}}') = \mathfrak{W}_1 \cup \mathfrak{W}'$ (where in the first statement we used the fact that $\mathfrak{U}_0 = G^{-1} (\mathfrak{W}_0)$ and the surjectivity of $G$, provided by \Cref{gw}). Together with the facts that $\omega_j (s) \in \overline{\mathfrak{W}}_0$ for each $s \in [0, 1]$, and that $\omega_j (s_j - \varepsilon) \notin \mathfrak{W}_1$ for each $\varepsilon \in (0, s_j]$, this implies that $\omega_j (s_j - \varepsilon) \in \mathfrak{W}'$ for each $\varepsilon \in (0, s_j]$. Since $\mathfrak{W}'$ is closed, it follows that $w_j = \omega_j (s_j) \in \mathfrak{W}'$. Hence, $w_1, w_2, \ldots \in \overline{\mathfrak{W}}$ is a sequence of mutually distinct points converging to $w$ such that $\overline{\mathfrak{U}}_1 \cap G^{-1} (w_i)$ and $\overline{\mathfrak{U}}' \cap G^{-1} (w_i)$ are nonempty for each integer $i \ge 1$ (as $w_i \in \mathfrak{W}_1 \cap \mathfrak{W}' = G(\overline{\mathfrak{U}}_1) \cap G(\overline{\mathfrak{U}}')$). This establishes the claim and thus the lemma.
		\end{proof}

	\begin{proof}[Proof of \Cref{degreed}]
		
		Let $\mathfrak{V} = G^{-1} (\overline{\mathfrak{W}})$, which is closed and thus compact, as $\mathfrak{V} \subseteq \overline{\mathfrak{R}}$ is bounded. By \Cref{gw} and \Cref{g1}, the map $G : \mathfrak{V} \rightarrow \overline{\mathfrak{W}}$ is bijective. Since $G$ is also continuous, and $\mathfrak{V}$ and $\overline{\mathfrak{W}}$ are compact, it follows that $G : \mathfrak{V} \rightarrow \overline{\mathfrak{W}}$ is a homeomorphism. Let $\mathfrak{U}$ denote the interior of $\mathfrak{V}$; since $\mathfrak{W}$ is the interior of $\overline{\mathfrak{W}}$, it follows that $G : \mathfrak{U} \rightarrow \mathfrak{W}$ is a homeomorphism, establishing the proposition. 
	\end{proof}

\section{Convergence of the KPZ and Log-Gamma Line Ensembles} 

\label{EquationPolymer}

In this section we provide two quick corollaries of \Cref{constantk1}. The former (first proven through a combination of the works \cite{CLE,CEPEFP,HTL,CTLE}) states convergence of the KPZ line ensemble, originally defined in \cite{LE}, to the Airy one. The latter states convergence of the log-gamma line ensemble, originally defined in \cite{TDLEE,SLNP} (and whose top curve provides the free energy distribution for the log-gamma polymer model introduced in \cite{SDPBC}), to the Airy one. 

\begin{cor} 
	
	\label{converge0}
	
	As $S$ tends to $\infty$, the KPZ line ensemble $\bm{\mathfrak{H}}^S = ( \mathfrak{H}_1^S, \mathfrak{H}_2^S, \ldots )$ (defined by \cite[Definition 1.2]{CTLE}) converges to the rescaled parabolic Airy line ensemble $2^{-1/2} \cdot \bm{\mathcal{S}}^{(q)}$ (defined by \eqref{sigmar}) at $q = 2^{-1/6}$, uniformly on compact subsets of $\mathbb{Z}_{\ge 1} \times \mathbb{R}$. 
\end{cor} 

\begin{proof} 
	
	By \cite[Theorem 1.5]{CTLE}, the sequence of line ensembles $\{ \bm{\mathfrak{H}}^S \}_{S > 0}$ is tight as $S$ tends to $\infty$ (under the topology of uniform convergence on compact subsets of $\mathbb{Z}_{\ge 1} \times \mathbb{R}$), and any subsequential limit point satisfies the Brownian Gibbs property. Fixing such a subsequential limit point $\bm{\mathfrak{H}}^{\infty}$, it suffices to show that $\bm{\mathfrak{H}}^{\infty}= \bm{\mathcal{S}}^{(q)}$ at $q = 2^{-1/6}$. To this end, by \cite[Proposition 1.4 and Corollary 1.6]{PFEC}, we have
	\begin{flalign}
		\label{h10} 
		\mathbb{P} \Big[ \mathfrak{H}_1^{\infty} (t) - \displaystyle\frac{t^2}{2} \le a \Big] = F_{\TW} (2^{1/3} a), \qquad \text{for each $a, t \in \mathbb{R}$},
	\end{flalign}
	
	\noindent where $F_{\TW} (s)$ denotes the Tracy--Widom GUE distribution. This verifies \Cref{constantk20} of \Cref{constantk1} at $(\ell, q) = (0, 2^{-1/6})$, so it follows that there exists a rescaled parabolic Airy line ensemble $\bm{\mathcal{S}}^{(q)}$ and an independent random variable $\mathfrak{c} \in \mathbb{R}$ such that $\mathfrak{H}_j^{\infty} (t)= \mathcal{S}_j^{(q)} (t) + \mathfrak{c}$ for each $(j, t) \in \mathbb{Z}_{\ge 1} \times \mathbb{R}$. Together with the $t=0$ case of \eqref{h10} and the fact that 
	\begin{flalign*} 
		\mathbb{P} \big[ \mathcal{S}_1^{(q)} (0) \le a \big] = \mathbb{P} \big[\mathcal{A}_1 (0) \le 2^{1/3} a \big] = F_{\TW} (2^{1/3} a), \qquad \text{for each $a \in \mathbb{R}$},
	\end{flalign*} 

	\noindent where we recall the Airy line ensemble $\bm{\mathcal{A}} = (\mathcal{A}_1, \mathcal{A}_2, \ldots )$ from \Cref{ensemblewalks}, it follows that $\mathfrak{c} = 0$. This confirms that $\bm{\mathfrak{H}}^{\infty} = \bm{\mathcal{S}}^{(q)}$, thus establishing the corollary.
\end{proof}

\begin{cor} 
	
	\label{converge1}
	
	Fix a real number $\theta > 0$; let $\sigma = \Psi' ( \theta / 2)^{1/2}$, where $\Psi : \mathbb{R}_{> 0} \rightarrow \mathbb{R}$ denotes the digamma function; and define the function $d_{\theta} : \mathbb{R}_{> 0} \rightarrow \mathbb{R}$ as in \cite[Equation (1.9)]{STEL}. As $N$ tends to $\infty$, the log-gamma line ensemble $\bm{\mathcal{L}}^N = (\mathcal{L}_1^N, \mathcal{L}_2^N, \ldots , \mathcal{L}_N^N, \ldots )$ with parameter $\theta$ (defined by \cite[Equations (1.12) and (1.3)]{TLE}) converges to the rescaled parabolic Airy line ensemble $\bm{\mathcal{S}}^{(q)}$ (defined by \eqref{sigmar}) at $q = 2^{-5/6} \sigma \cdot d_{\theta} (1)^{-1}$, uniformly on compact subsets of $\mathbb{Z}_{\ge 1} \times \mathbb{R}$. 
	
\end{cor} 

\begin{proof} 
	
	By \cite[Theorem 1.11]{TLE}, the sequence of line ensembles $\{ \bm{\mathcal{L}}^N \}_{N \ge 1}$ is tight as $N$ tends to $\infty$ (under the topology of uniform convergence on compact subsets of $\mathbb{Z}_{\ge 1} \times \mathbb{R}$), and any subsequential limit point satisfies the Brownian Gibbs property. Fixing such a subsequential limit point $\bm{\mathcal{L}}^{\infty}$, it suffices to show that $\bm{\mathcal{L}}^{\infty}= \bm{\mathcal{S}}^{(q)}$ at $q = 2^{-5/6} \sigma \cdot d_{\theta}(1)^{-1}$. To this end, by \cite[Equations (3.4) and (1.3)]{TLE} with \cite[Equation (1.12)]{STEL}, we have
	\begin{flalign}
		\label{h100} 
		\mathbb{P} \big[ \mathcal{L}_1^{\infty} (t) - 2^{-1/2} q^3 t^2 \le a \big] = F_{\TW} (2^{1/2} q a), \qquad \text{for each $a, t \in \mathbb{R}$},
	\end{flalign}
	
	\noindent where $F_{\TW} (s)$ denotes the Tracy--Widom GUE distribution. This verifies part \Cref{constantk20} of \Cref{constantk1} at $\ell = 0$, so it follows that there exists a rescaled parabolic Airy line ensemble $\bm{\mathcal{S}}^{(q)}$ and an independent random variable $\mathfrak{c} \in \mathbb{R}$ such that $\mathcal{L}_j^{\infty} (t)= \mathcal{S}_j^{(q)} (t) + \mathfrak{c}$ for each $(j, t) \in \mathbb{Z}_{\ge 1} \times \mathbb{R}$. Together with the $t=0$ case of \eqref{h100} and the fact that 
	\begin{flalign*} 
		\mathbb{P} \big[ \mathcal{S}_1^{(q)} (0) \le a \big] = \mathbb{P} \big[\mathcal{A}_1 (0) \le 2^{1/2} q a \big] = F_{\TW} (2^{1/2} q a), \qquad \text{for each $a \in \mathbb{R}$,}
	\end{flalign*} 
	
	\noindent where we recall the Airy line ensemble $\bm{\mathcal{A}} = (\mathcal{A}_1, \mathcal{A}_2, \ldots )$ from \Cref{ensemblewalks}, it follows that $\mathfrak{c} = 0$. This confirms that $\bm{\mathcal{L}}^{\infty} = \bm{\mathcal{S}}^{(q)}$, thus establishing the corollary.
\end{proof}

\printindex


\begin{thebibliography}{100}
	
	\bibitem{PWUC}
	M.~Adler, P.~L. Ferrari, and P.~van Moerbeke.
	\newblock Airy processes with wanderers and new universality classes.
	\newblock {\em Ann. Probab.}, 38(2):714--769, 2010.
	
	\bibitem{ULTS}
	A.~Aggarwal.
	\newblock Universality for lozenge tiling local statistics.
	\newblock {\em Ann. of Math. (2)}, 198(3):881--1012, 2023.
	
	\bibitem{UU}
	A.~Aggarwal and A.~Borodin.
	\newblock Colored line ensembles for stochastic vertex models.
	\newblock {\em Selecta Math. (N.S.)}, 30(5):Paper No. 105, 111, 2024.
	
	\bibitem{SSE}
	A.~Aggarwal, A.~Borodin, and A.~Bufetov.
	\newblock Stochasticization of solutions to the {Y}ang-{B}axter equation.
	\newblock {\em Ann. Henri Poincar\'{e}}, 20(8):2495--2554, 2019.
	
	\bibitem{DPG}
	A.~Aggarwal, A.~Borodin, and M.~Wheeler.
	\newblock Deformed polynuclear growth in {$(1+1)$} dimensions.
	\newblock {\em Int. Math. Res. Not.}, (7):5728--5780, 2023.
	
	\bibitem{SLCA}
	A.~Aggarwal, I.~Corwin, and M.~Hegde.
	\newblock Scaling limit of the colored {ASEP} and stochastic six-vertex models.
	\newblock Preprint, ar{X}iv:2403.01341.
	
	\bibitem{ERGIM}
	A.~Aggarwal and J.~Huang.
	\newblock Edge rigidity of {D}yson {B}rownian motion with general initial data.
	\newblock {\em Electron. J. Probab.}, 29:Paper No. 137, 62, 2024.
	
	\bibitem{ESTP}
	A.~Aggarwal and J.~Huang.
	\newblock Edge statistics for lozenge tilings of polygons, {II}: {A}iry line
	ensemble.
	\newblock {\em Forum Math. Pi}, 13:Paper No. e2, 60, 2025.
	
	\bibitem{U}
	A.~Aggarwal and J.~Huang.
	\newblock Local {S}tatistics and {C}oncentration for {N}on-intersecting
	{B}rownian {B}ridges with {S}mooth {B}oundary {D}ata.
	\newblock {\em Comm. Math. Phys.}, 406(3):Paper No. 70, 2025.
	
	\bibitem{IDA}
	C.~D. Aliprantis and K.~C. Border.
	\newblock {\em Infinite dimensional analysis}.
	\newblock Springer, Berlin, third edition, 2006.
	\newblock A hitchhiker's guide.
	
	\bibitem{PFEC}
	G.~Amir, I.~Corwin, and J.~Quastel.
	\newblock Probability distribution of the free energy of the continuum directed
	random polymer in {$1+1$} dimensions.
	\newblock {\em Comm. Pure Appl. Math.}, 64(4):466--537, 2011.
	
	\bibitem{IRM}
	G.~W. Anderson, A.~Guionnet, and O.~Zeitouni.
	\newblock {\em An introduction to random matrices}, volume 118 of {\em
		Cambridge Studies in Advanced Mathematics}.
	\newblock Cambridge University Press, Cambridge, 2010.
	
	\bibitem{DMC}
	K.~Astala, E.~Duse, I.~Prause, and X.~Zhong.
	\newblock Dimer models and conformal structures.
	\newblock Preprint, ar{X}iv:2004.02599.
	
	\bibitem{FMTP}
	A.~Ayyer, S.~Chhita, and K.~Johansson.
	\newblock G{OE} fluctuations for the maximum of the top path in alternating
	sign matrices.
	\newblock {\em Duke Math. J.}, 172(10):1961--2014, 2023.
	
	\bibitem{DLLIS}
	J.~Baik, P.~Deift, and K.~Johansson.
	\newblock On the distribution of the length of the longest increasing
	subsequence of random permutations.
	\newblock {\em J. Amer. Math. Soc.}, 12(4):1119--1178, 1999.
	
	\bibitem{RWDR}
	G.~Barraquand and I.~Corwin.
	\newblock Random-walk in beta-distributed random environment.
	\newblock {\em Probab. Theory Related Fields}, 167(3-4):1057--1116, 2017.
	
	\bibitem{STEL}
	G.~Barraquand, I.~Corwin, and E.~Dimitrov.
	\newblock Spatial tightness at the edge of {G}ibbsian line ensembles.
	\newblock {\em Comm. Math. Phys.}, 397(3):1309--1386, 2023.
	
	\bibitem{GQ}
	Y.~Baryshnikov.
	\newblock G{UE}s and queues.
	\newblock {\em Probab. Theory Related Fields}, 119(2):256--274, 2001.
	
	\bibitem{FGPS}
	R.~Basu, S.~Ganguly, and A.~Hammond.
	\newblock Fractal geometry of {${\rm Airy}_2$} processes coupled via the {A}iry
	sheet.
	\newblock {\em Ann. Probab.}, 49(1):485--505, 2021.
	
	\bibitem{ISEGW}
	R.~Basu, S.~Ganguly, A.~Hammond, and M.~Hegde.
	\newblock Interlacing and scaling exponents for the geodesic watermelon in last
	passage percolation.
	\newblock {\em Comm. Math. Phys.}, 393(3):1241--1309, 2022.
	
	\bibitem{LDSI}
	S.~Belinschi, A.~Guionnet, and J.~Huang.
	\newblock Large deviation principles via spherical integrals.
	\newblock {\em Probab. Math. Phys.}, 3(3):543--625, 2022.
	
	\bibitem{FCSD}
	P.~Biane.
	\newblock On the free convolution with a semi-circular distribution.
	\newblock {\em Indiana University Mathematics Journal}, pages 705--718, 1997.
	
	\bibitem{SVMP}
	A.~Borodin, A.~Bufetov, and M.~Wheeler.
	\newblock Between the stochastic six vertex model and {H}all-{L}ittlewood
	processes.
	\newblock Preprint, ar{X}iv:1611.09486.
	
	\bibitem{P}
	A.~Borodin and I.~Corwin.
	\newblock Macdonald processes.
	\newblock {\em Probab. Theory Related Fields}, 158(1-2):225--400, 2014.
	
	\bibitem{SSVM}
	A.~Borodin, I.~Corwin, and V.~Gorin.
	\newblock Stochastic six-vertex model.
	\newblock {\em Duke Math. J.}, 165(3):563--624, 2016.
	
	\bibitem{AGRS}
	A.~Borodin and P.~L. Ferrari.
	\newblock Anisotropic growth of random surfaces in {$2+1$} dimensions.
	\newblock {\em Comm. Math. Phys.}, 325(2):603--684, 2014.
	
	\bibitem{KTSP}
	A.~Borodin and S.~P\'ech\'e.
	\newblock Airy kernel with two sets of parameters in directed percolation and
	random matrix theory.
	\newblock {\em J. Stat. Phys.}, 132(2):275--290, 2008.
	
	\bibitem{USCG}
	E.~Br\'{e}zin and S.~Hikami.
	\newblock Universal singularity at the closure of a gap in a random matrix
	theory.
	\newblock {\em Phys. Rev. E (3)}, 57(4):4140--4149, 1998.
	
	\bibitem{RFSVM}
	A.~Bufetov, M.~Mucciconi, and L.~Petrov.
	\newblock Yang-{B}axter random fields and stochastic vertex models.
	\newblock {\em Adv. Math.}, 388:Paper No. 107865, 94, 2021.
	
	\bibitem{FSF}
	A.~Bufetov and L.~Petrov.
	\newblock Yang-{B}axter field for spin {H}all-{L}ittlewood symmetric functions.
	\newblock {\em Forum Math. Sigma}, 7:Paper No. e39, 70, 2019.
	
	\bibitem{ITFP}
	J.~Calvert, A.~Hammond, and M.~Hegde.
	\newblock Brownian structure in the {KPZ} fixed point.
	\newblock {\em Ast\'{e}risque}, (441):vi+119, 2023.
	
	\bibitem{FESD}
	M.~Capitaine and S.~P\'{e}ch\'{e}.
	\newblock Fluctuations at the edges of the spectrum of the full rank deformed
	{GUE}.
	\newblock {\em Probab. Theory Related Fields}, 165(1-2):117--161, 2016.
	
	\bibitem{CPM}
	R.~Cerf.
	\newblock {\em The {W}ulff crystal in {I}sing and percolation models}, volume
	1878 of {\em Lecture Notes in Mathematics}.
	\newblock Springer-Verlag, Berlin, 2006.
	\newblock Lectures from the 34th Summer School on Probability Theory held in
	Saint-Flour, July 6--24, 2004, With a foreword by Jean Picard.
	
	\bibitem{LSRT}
	H.~Cohn, N.~Elkies, and J.~Propp.
	\newblock Local statistics for random domino tilings of the {A}ztec diamond.
	\newblock {\em Duke Math. J.}, 85(1):117--166, 1996.
	
	\bibitem{EUC}
	I.~Corwin.
	\newblock The {K}ardar-{P}arisi-{Z}hang equation and universality class.
	\newblock {\em Random Matrices Theory Appl.}, 1(1):1130001, 76, 2012.
	
	\bibitem{TFSVM}
	I.~Corwin and E.~Dimitrov.
	\newblock Transversal fluctuations of the {ASEP}, stochastic six vertex model,
	and {H}all-{L}ittlewood {G}ibbsian line ensembles.
	\newblock {\em Comm. Math. Phys.}, 363(2):435--501, 2018.
	
	\bibitem{ECT}
	I.~Corwin, P.~Ghosal, and A.~Hammond.
	\newblock K{PZ} equation correlations in time.
	\newblock {\em Ann. Probab.}, 49(2):832--876, 2021.
	
	\bibitem{PLE}
	I.~Corwin and A.~Hammond.
	\newblock Brownian {G}ibbs property for {A}iry line ensembles.
	\newblock {\em Invent. Math.}, 195(2):441--508, 2014.
	
	\bibitem{LE}
	I.~Corwin and A.~Hammond.
	\newblock K{PZ} line ensemble.
	\newblock {\em Probab. Theory Related Fields}, 166(1-2):67--185, 2016.
	
	\bibitem{ETFP}
	I.~Corwin, A.~Hammond, M.~Hegde, and K.~Matetski.
	\newblock Exceptional times when the {KPZ} fixed point violates {J}ohansson's
	conjecture on maximizer uniqueness.
	\newblock {\em Electron. J. Probab.}, 28:Paper No. 11, 81, 2023.
	
	\bibitem{TCF}
	I.~Corwin, N.~O'Connell, T.~Sepp\"{a}l\"{a}inen, and N.~Zygouras.
	\newblock Tropical combinatorics and {W}hittaker functions.
	\newblock {\em Duke Math. J.}, 163(3):513--563, 2014.
	
	\bibitem{SLP}
	I.~Corwin, T.~Sepp\"{a}l\"{a}inen, and H.~Shen.
	\newblock The strict-weak lattice polymer.
	\newblock {\em J. Stat. Phys.}, 160(4):1027--1053, 2015.
	
	\bibitem{ELE}
	I.~Corwin and X.~Sun.
	\newblock Ergodicity of the {A}iry line ensemble.
	\newblock {\em Electron. Commun. Probab.}, 19:no. 49, 11, 2014.
	
	\bibitem{LCDP}
	S.~Das and W.~Zhu.
	\newblock Localization of the continuum directed random polymer.
	\newblock {\em Probab. Theory Related Fields}, page To appear.
	\newblock Preprint, ar{X}iv:2203.03607.
	
	\bibitem{GNL}
	D.~Dauvergne.
	\newblock The 27 geodesic networks in the directed landscape.
	\newblock Preprint, ar{X}iv:2302.07802.
	
	\bibitem{NTM}
	D.~Dauvergne.
	\newblock Non-uniqueness times for the maximizer of the {KPZ} fixed point.
	\newblock {\em Adv. Math.}, 442:Paper No. 109550, 41, 2024.
	
	\bibitem{DLE}
	D.~Dauvergne.
	\newblock Wiener densities for the {A}iry line ensemble.
	\newblock {\em Proc. Lond. Math. Soc. (3)}, 129(4):Paper No. e12638, 57, 2024.
	
	\bibitem{TL}
	D.~Dauvergne, J.~Ortmann, and B.~Vir\'{a}g.
	\newblock The directed landscape.
	\newblock {\em Acta Math.}, 229(2):201--285, 2022.
	
	\bibitem{BPLE}
	D.~Dauvergne and B.~Vir\'{a}g.
	\newblock Bulk properties of the {A}iry line ensemble.
	\newblock {\em Ann. Probab.}, 49(4):1738--1777, 2021.
	
	\bibitem{DOL}
	D.~Dauvergne and L.~Zhang.
	\newblock Disjoint optimizers and the directed landscape.
	\newblock {\em Mem. Amer. Math. Soc.}, 303(1524):v+103, 2024.
	
	\bibitem{MCFARS}
	D.~De~Silva and O.~Savin.
	\newblock Minimizers of convex functionals arising in random surfaces.
	\newblock {\em Duke Math. J.}, 151(3):487--532, 2010.
	
	\bibitem{WLE}
	E.~Dimitrov.
	\newblock Airy wanderer line ensembles.
	\newblock Preprint, ar{X}iv:2408.08445.
	
	\bibitem{CHLE}
	E.~Dimitrov.
	\newblock Characterization of {$H$}-{B}rownian {G}ibbsian line ensembles.
	\newblock {\em Probab. Math. Phys.}, 3(3):627--673, 2022.
	
	\bibitem{CLE}
	E.~Dimitrov and K.~Matetski.
	\newblock Characterization of {B}rownian {G}ibbsian line ensembles.
	\newblock {\em Ann. Probab.}, 49(5):2477--2529, 2021.
	
	\bibitem{TLE}
	E.~Dimitrov and X.~Wu.
	\newblock Tightness of {$(H, H^{RW})$}-{G}ibbsian line ensembles.
	\newblock Preprint, ar{X}iv:2108.0784.
	
	\bibitem{CRWB}
	E.~Dimitrov and X.~Wu.
	\newblock K{MT} coupling for random walk bridges.
	\newblock {\em Probab. Theory Related Fields}, 179(3-4):649--732, 2021.
	
	\bibitem{C}
	R.~Dobrushin, R.~Koteck\'{y}, and S.~Shlosman.
	\newblock {\em Wulff construction}, volume 104 of {\em Translations of
		Mathematical Monographs}.
	\newblock American Mathematical Society, Providence, RI, 1992.
	\newblock A global shape from local interaction, Translated from the Russian by
	the authors.
	
	\bibitem{UCT}
	R.~M. Dudley.
	\newblock {\em Uniform central limit theorems}, volume~63 of {\em Cambridge
		Studies in Advanced Mathematics}.
	\newblock Cambridge University Press, Cambridge, 1999.
	
	\bibitem{MERM}
	F.~J. Dyson.
	\newblock A {B}rownian-motion model for the eigenvalues of a random matrix.
	\newblock {\em J. Mathematical Phys.}, 3(6):1191--1198, 1962.
	
	\bibitem{AMT}
	N.~Elkies, G.~Kuperberg, M.~Larsen, and J.~Propp.
	\newblock Alternating-sign matrices and domino tilings. {II}.
	\newblock {\em J. Algebraic Combin.}, 1(3):219--234, 1992.
	
	\bibitem{LSLCD}
	L.~Erd\H{o}s, B.~Schlein, and H.-T. Yau.
	\newblock Local semicircle law and complete delocalization for {W}igner random
	matrices.
	\newblock {\em Comm. Math. Phys.}, 287(2):641--655, 2009.
	
	\bibitem{SLSSD}
	L.~Erd\H{o}s, B.~Schlein, and H.-T. Yau.
	\newblock Semicircle law on short scales and delocalization of eigenvectors for
	{W}igner random matrices.
	\newblock {\em Ann. Probab.}, 37(3):815--852, 2009.
	
	\bibitem{DRM}
	L.~Erd{\H{o}}s and H.-T. Yau.
	\newblock {\em A dynamical approach to random matrix theory}, volume~28.
	\newblock Amer. Math. Soc., 2017.
	
	\bibitem{PM}
	P.~L. Ferrari and S.~Shlosman.
	\newblock The {${\rm Airy_2}$} process and the 3{D} {I}sing model.
	\newblock {\em J. Phys. A}, 56(1):Paper No. 014003, 15, 2023.
	
	\bibitem{RGM}
	P.~L. Ferrari and H.~Spohn.
	\newblock {Random growth models}.
	\newblock In {\em {The Oxford Handbook of Random Matrix Theory}}. Oxford
	University Press, 09 2015.
	
	\bibitem{COUS}
	P.~J. Forrester, T.~Nagao, and G.~Honner.
	\newblock Correlations for the orthogonal-unitary and symplectic-unitary
	transitions at the hard and soft edges.
	\newblock {\em Nuclear Phys. B}, 553(3):601--643, 1999.
	
	\bibitem{SHDLP}
	S.~Ganguly and A.~Hammond.
	\newblock Stability and chaos in dynamical last passage percolation.
	\newblock {\em Commun. Am. Math. Soc.}, 4:387--479, 2024.
	
	\bibitem{UTEL}
	S.~Ganguly and M.~Hegde.
	\newblock Sharp upper tail estimates and limit shapes for the {KPZ} equation
	via the tangent method.
	\newblock Preprint, ar{X}iv:2208.08922.
	
	\bibitem{LGC}
	S.~Ganguly and M.~Hegde.
	\newblock Local and global comparisons of the {A}iry difference profile to
	{B}rownian local time.
	\newblock {\em Ann. Inst. Henri Poincar\'e{} Probab. Stat.}, 59(3):1342--1374,
	2023.
	
	\bibitem{FGSDP}
	S.~Ganguly and L.~Zhang.
	\newblock Fractal geometry of the space-time difference profile in the directed
	landscape via construction of geodesic local times.
	\newblock Preprint, ar{X}iv:2204.01674.
	
	\bibitem{MCNP}
	D.~J. Grabiner.
	\newblock Brownian motion in a {W}eyl chamber, non-colliding particles, and
	random matrices.
	\newblock {\em Ann. Inst. H. Poincar\'{e} Probab. Statist.}, 35(2):177--204,
	1999.
	
	\bibitem{DT}
	V.~Guillemin and A.~Pollack.
	\newblock {\em Differential topology}.
	\newblock AMS Chelsea Publishing, Providence, RI, 2010.
	\newblock Reprint of the 1974 original.
	
	\bibitem{FOAMI}
	A.~Guionnet.
	\newblock First order asymptotics of matrix integrals; a rigorous approach
	towards the understanding of matrix models.
	\newblock {\em Comm. Math. Phys.}, 244(3):527--569, 2004.
	
	\bibitem{LDSCRM}
	A.~Guionnet.
	\newblock Large deviations and stochastic calculus for large random matrices.
	\newblock {\em Probab. Surv.}, 1:72--172, 2004.
	
	\bibitem{LDASI}
	A.~Guionnet and O.~Zeitouni.
	\newblock Large deviations asymptotics for spherical integrals.
	\newblock {\em J. Funct. Anal.}, 188(2):461--515, 2002.
	
	\bibitem{SVMR}
	L.-H. Gwa and H.~Spohn.
	\newblock Six-vertex model, roughened surfaces, and an asymmetric spin
	{H}amiltonian.
	\newblock {\em Phys. Rev. Lett.}, 68(6):725--728, 1992.
	
	\bibitem{ET}
	M.~Hairer and J.~C. Mattingly.
	\newblock Yet another look at {H}arris' ergodic theorem for {M}arkov chains.
	\newblock In {\em Seminar on {S}tochastic {A}nalysis, {R}andom {F}ields and
		{A}pplications {VI}}, volume~63 of {\em Progr. Probab.}, pages 109--117.
	Birkh\"{a}user/Springer Basel AG, Basel, 2011.
	
	\bibitem{MCPWP}
	A.~Hammond.
	\newblock Modulus of continuity of polymer weight profiles in {B}rownian last
	passage percolation.
	\newblock {\em Ann. Probab.}, 47(6):3911--3962, 2019.
	
	\bibitem{PRPW}
	A.~Hammond.
	\newblock A patchwork quilt sewn from {B}rownian fabric: regularity of polymer
	weight profiles in {B}rownian last passage percolation.
	\newblock {\em Forum Math. Pi}, 7:e2, 69, 2019.
	
	\bibitem{ERDP}
	A.~Hammond.
	\newblock Exponents governing the rarity of disjoint polymers in {B}rownian
	last passage percolation.
	\newblock {\em Proc. Lond. Math. Soc.}, 120(3):370--433, 2020.
	
	\bibitem{RLE}
	A.~Hammond.
	\newblock Brownian regularity for the {A}iry line ensemble, and multi-polymer
	watermelons in {B}rownian last passage percolation.
	\newblock {\em Mem. Amer. Math. Soc.}, 277(1363):v+133, 2022.
	
	\bibitem{ESMP}
	T.~E. Harris.
	\newblock The existence of stationary measures for certain {M}arkov processes.
	\newblock In {\em Proceedings of the {T}hird {B}erkeley {S}ymposium on
		{M}athematical {S}tatistics and {P}robability, 1954--1955, vol. {II}}, pages
	113--124. University of California Press, Berkeley-Los Angeles, Calif., 1956.
	
	\bibitem{MCTMGP}
	J.~Huang and B.~Landon.
	\newblock Rigidity and a mesoscopic central limit theorem for {D}yson
	{B}rownian motion for general {$\beta$} and potentials.
	\newblock {\em Probab. Theory Related Fields}, 175(1-2):209--253, 2019.
	
	\bibitem{RTT}
	W.~Jockusch, J.~Propp, and P.~Shor.
	\newblock Random domino tilings and the arctic circle theorem.
	\newblock Preprint, ar{X}iv:9801068, 1998.
	
	\bibitem{SFRM}
	K.~Johansson.
	\newblock Shape fluctuations and random matrices.
	\newblock {\em Comm. Math. Phys.}, 209(2):437--476, 2000.
	
	\bibitem{ULSDM}
	K.~Johansson.
	\newblock Universality of the local spacing distribution in certain ensembles
	of {H}ermitian {W}igner matrices.
	\newblock {\em Comm. Math. Phys.}, 215(3):683--705, 2001.
	
	\bibitem{NPRTM}
	K.~Johansson.
	\newblock Non-intersecting paths, random tilings and random matrices.
	\newblock {\em Probab. Theory Related Fields}, 123(2):225--280, 2002.
	
	\bibitem{DPGDP}
	K.~Johansson.
	\newblock Discrete polynuclear growth and determinantal processes.
	\newblock {\em Comm. Math. Phys.}, 242(1-2):277--329, 2003.
	
	\bibitem{TBP}
	K.~Johansson.
	\newblock The arctic circle boundary and the {A}iry process.
	\newblock {\em Ann. Probab.}, 33(1):1--30, 2005.
	
	\bibitem{RMDP}
	K.~Johansson.
	\newblock Random matrices and determinantal processes.
	\newblock In {\em Mathematical statistical physics}, pages 1--56. Elsevier B.
	V., Amsterdam, 2006.
	
	\bibitem{EFLS}
	K.~Johansson.
	\newblock Edge fluctuations of limit shapes.
	\newblock In {\em Current developments in mathematics 2016}, pages 47--110.
	Int. Press, Somerville, MA, 2018.
	
	\bibitem{SLNP}
	S.~G.~G. Johnston and N.~O'Connell.
	\newblock Scaling limits for non-intersecting polymers and {W}hittaker
	measures.
	\newblock {\em J. Stat. Phys.}, 179(2):354--407, 2020.
	
	\bibitem{MSC}
	I.~Karatzas and S.~E. Shreve.
	\newblock {\em Brownian motion and stochastic calculus}, volume 113 of {\em
		Graduate Texts in Mathematics}.
	\newblock Springer-Verlag, New York, second edition, 1991.
	
	\bibitem{DSGI}
	M.~Kardar, G.~Parisi, and Y.-C. Zhang.
	\newblock Dynamic scaling of growing interfaces.
	\newblock {\em Phys. Rev. Lett.}, 56(9):889, 1986.
	
	\bibitem{LSCE}
	R.~Kenyon and A.~Okounkov.
	\newblock Limit shapes and the complex {B}urgers equation.
	\newblock {\em Acta Math.}, 199(2):263--302, 2007.
	
	\bibitem{landon2017edge}
	B.~Landon and H.-T. Yau.
	\newblock Edge statistics of {D}yson {B}rownian motion.
	\newblock Preprint, ar{X}iv:1712.03881. 
	
	\bibitem{TDMS}
	B.~Laslier and F.~L. Toninelli.
	\newblock Lozenge tilings, {G}lauber dynamics and macroscopic shape.
	\newblock {\em Comm. Math. Phys.}, 338(3):1287--1326, 2015.
	
	\bibitem{EUDM}
	J.~O. Lee and K.~Schnelli.
	\newblock Edge universality for deformed {W}igner matrices.
	\newblock {\em Rev. Math. Phys.}, 27(8):1550018, 94, 2015.
	
	\bibitem{UPCSE}
	A.~M.~S. Mac{\^e}do.
	\newblock Universal parametric correlations at the soft edge of the spectrum of
	random matrix ensembles.
	\newblock {\em Europhys. Lett.}, 26(9):641, 1994.
	
	\bibitem{TFP}
	K.~Matetski, J.~Quastel, and D.~Remenik.
	\newblock The {KPZ} fixed point.
	\newblock {\em Acta Math.}, 227(1):115--203, 2021.
	
	\bibitem{LLI}
	A.~Matytsin.
	\newblock On the large-{$N$} limit of the {I}tzykson-{Z}uber integral.
	\newblock {\em Nuclear Phys. B}, 411(2-3):805--820, 1994.
	
	\bibitem{CSS}
	S.~Meyn and R.~L. Tweedie.
	\newblock {\em Markov chains and stochastic stability}.
	\newblock Cambridge University Press, Cambridge, second edition, 2009.
	\newblock With a prologue by Peter W. Glynn.
	
	\bibitem{FPRM}
	J.~A. Mingo and R.~Speicher.
	\newblock {\em Free probability and random matrices}, volume~35 of {\em Fields
		Institute Monographs}.
	\newblock Springer, New York; Fields Institute for Research in Mathematical
	Sciences, Toronto, ON, 2017.
	
	\bibitem{IDL}
	M.~Nica.
	\newblock Intermediate disorder limits for multi-layer semi-discrete directed
	polymers.
	\newblock {\em Electron. J. Probab.}, 26:Paper No. 62, 50, 2021.
	
	\bibitem{PTRW}
	N.~O'Connell.
	\newblock A path-transformation for random walks and the {R}obinson-{S}chensted
	correspondence.
	\newblock {\em Trans. Amer. Math. Soc.}, 355(9):3669--3697, 2003.
	
	\bibitem{DPQL}
	N.~O'Connell.
	\newblock Directed polymers and the quantum {T}oda lattice.
	\newblock {\em Ann. Probab.}, 40(2):437--458, 2012.
	
	\bibitem{MSE}
	N.~O'Connell and J.~Warren.
	\newblock A multi-layer extension of the stochastic heat equation.
	\newblock {\em Comm. Math. Phys.}, 341(1):1--33, 2016.
	
	\bibitem{AT}
	N.~O'Connell and M.~Yor.
	\newblock Brownian analogues of {B}urke's theorem.
	\newblock {\em Stochastic Process. Appl.}, 96(2):285--304, 2001.
	
	\bibitem{CFPALG}
	A.~Okounkov and N.~Reshetikhin.
	\newblock Correlation function of {S}chur process with application to local
	geometry of a random 3-dimensional {Y}oung diagram.
	\newblock {\em J. Amer. Math. Soc.}, 16(3):581--603, 2003.
	
	\bibitem{okounkov2007random}
	A.~Okounkov and N.~Reshetikhin.
	\newblock Random skew plane partitions and the {P}earcey process.
	\newblock {\em Comm. Math. Phys.}, 269(3):571--609, 2007.
	
	\bibitem{AFP}
	L.~P.~R. Pimentel.
	\newblock Brownian aspects of the {KPZ} fixed point.
	\newblock In {\em In and out of equilibrium 3. {C}elebrating {V}ladas
		{S}idoravicius}, volume~77 of {\em Progr. Probab.}, pages 711--739.
	Birkh\"{a}user/Springer, Cham, [2021] \copyright 2021.
	
	\bibitem{GSSPD}
	V.~L. Pokrovsky and A.~L. Talapov.
	\newblock Ground state, spectrum, and phase diagram of two-dimensional
	incommensurate crystals.
	\newblock {\em Phys. Rev. Lett.}, 42:65--67, 1979.
	
	\bibitem{TTDIC}
	V.~L. Pokrovsky and A.~L. Talapov.
	\newblock The theory of two-dimensional incommensurate crystals.
	\newblock {\em Zh. Eksp.Teo. Fiz.}, 78:269--295, 1980.
	
	\bibitem{SIDP}
	M.~Pr\"ahofer and H.~Spohn.
	\newblock Scale invariance of the {PNG} droplet and the {A}iry process.
	\newblock {\em J. Statist. Phys.}, 108(5-6):1071--1106, 2002.
	
	\bibitem{IT}
	J.~Quastel.
	\newblock Introduction to {KPZ}.
	\newblock In {\em Current developments in mathematics, 2011}, pages 125--194.
	Int. Press, Somerville, MA, 2012.
	
	\bibitem{CEPEFP}
	J.~Quastel and S.~Sarkar.
	\newblock Convergence of exclusion processes and the {KPZ} equation to the
	{KPZ} fixed point.
	\newblock {\em J. Amer. Math. Soc.}, 36(1):251--289, 2023.
	
	\bibitem{NBPP}
	H.~Rost.
	\newblock Non-equilibrium behaviour of a many particle process: Density profile
	and local equilibria.
	\newblock {\em Z. Wahrsch. Verw. Gebiete}, 58(1):41--53, 1981.
	
	\bibitem{ACPIC}
	S.~Sarkar and B.~Vir\'{a}g.
	\newblock Brownian absolute continuity of the {KPZ} fixed point with arbitrary
	initial condition.
	\newblock {\em Ann. Probab.}, 49(4):1718--1737, 2021.
	
	\bibitem{ERTD}
	T.~Sasamoto and M.~Wadati.
	\newblock Exact results for one-dimensional totally asymmetric diffusion
	models.
	\newblock {\em J. Phys. A}, 31(28):6057--6071, 1998.
	
	\bibitem{SDPBC}
	T.~Sepp\"{a}l\"{a}inen.
	\newblock Scaling for a one-dimensional directed polymer with boundary
	conditions.
	\newblock {\em Ann. Probab.}, 40(1):19--73, 2012.
	
	\bibitem{TDLE}
	C.~Serio.
	\newblock Tightness of discrete {G}ibbsian line ensembles.
	\newblock {\em Stochastic Process. Appl.}, 159:225--285, 2023.
	
	\bibitem{ULER}
	T.~Shcherbina.
	\newblock On universality of local edge regime for the deformed gaussian
	unitary ensemble.
	\newblock {\em J. Stat. Phs.}, 143, 2011.
	
	\bibitem{RS}
	S.~Sheffield.
	\newblock Random surfaces.
	\newblock {\em Ast\'{e}risque}, (304):vi+175, 2005.
	
	\bibitem{FNP}
	A.~B. Soshnikov.
	\newblock Gaussian fluctuation for the number of particles in {A}iry, {B}essel,
	sine, and other determinantal random point fields.
	\newblock {\em J. Statist. Phys.}, 100(3-4):491--522, 2000.
	
	\bibitem{EODLE}
	H.~Spohn.
	\newblock Kardar--{P}arisi--{Z}hang equation in one dimension and line
	ensembles.
	\newblock {\em Pramana}, 64(6):847--857, 2005.
	
	\bibitem{LDK}
	C.~A. Tracy and H.~Widom.
	\newblock Level-spacing distributions and the {A}iry kernel.
	\newblock {\em Comm. Math. Phys.}, 159(1):151--174, 1994.
	
	\bibitem{HTL}
	B.~Vir\'{a}g.
	\newblock The heat and the landscape {I}.
	\newblock Preprint, ar{X}iv:2008.07241.
	
	\bibitem{II}
	J.~Warren.
	\newblock Dyson's {B}rownian motions, intertwining and interlacing.
	\newblock {\em Electron. J. Probab.}, 12:no. 19, 573--590, 2007.
	
	\bibitem{EL}
	X.~Wu.
	\newblock The {KPZ} equation and the directed landscape.
	\newblock Preprint, ar{X}iv:2301.00547.
	
	\bibitem{CTLE}
	X.~Wu.
	\newblock Convergence of the {KPZ} line ensemble.
	\newblock {\em International Mathematics Research Notices}, page rnac272, 10
	2022.
	
	\bibitem{TDLEE}
	X.~Wu.
	\newblock Tightness of discrete {G}ibbsian line ensembles with exponential
	interaction {H}amiltonians.
	\newblock {\em Ann. Inst. Henri Poincar\'e{} Probab. Stat.}, 59(4):2106--2150,
	2023.
	
\end{thebibliography}
	\end{document}